%\special{!userdict begin /bop-hook{gsave 
%  /inch{72 mul}def
%  1.25 inch .875 inch translate
%  .2 setgray
%  0 0 moveto
%  0 9.25 inch lineto 
%  6 inch 9.25 inch lineto 
%  6 inch 0 lineto
%  0 0 lineto stroke
%  grestore}def end}

\input psfig.sty
\input epsf.tex
\baselineskip=14truept
\voffset=.75in  
%\voffset=1.40 truein % Setting in Yale
\hoffset=1truein
\hsize=4.5truein
\vsize=7.75truein
\parindent=.166666in
\pretolerance=500 \tolerance=1000 \brokenpenalty=5000
%\abovedisplayskip=12pt
%\abovedisplayshortskip=0pt
%\belowdisplayskip=12pt
%\belowdisplayshortskip=0pt
\medskipamount=7pt

\def\hexdigit#1{\ifnum#1<10 \number#1\else%
  \ifnum#1=10 A\else\ifnum#1=11 B\else\ifnum#1=12 C\else%
  \ifnum#1=13 D\else\ifnum#1=14 E\else\ifnum#1=15 F\fi%
  \fi\fi\fi\fi\fi\fi}

% romain
%\font\elevenrm=cmr12 at 24pt
\font\elevenrm=cmr12 at 11pt%
\font\eightrm=cmr8%
\font\sixrm=cmr6%

% math romain italique
\font\elevenmi=cmmi12 at 11pt%
\font\nineit=cmmi9%
\font\eightmi=cmmi8%
\font\sixmi=cmmi6%

% romain italique
\font\elevenit=cmti12 at 11pt%
%\font\eigthit=cmti8
%\font\sixit=cmti7 at 6pt

% slantered roman
\font\elevensl=cmsl12 at 11pt%
%\font\eightsl=cmsl8%
%\font\sixsl=cmsl8 at 6pt%

% symboles a 11 pt
\font\elevensy=cmsy10 at 11pt%
\font\eigthsy=cmsy8%
\font\sixsy=cmsy6%

% symboles etendus a 11 pt
\font\elevensyex=cmex10 at 11pt%
\font\eightsyex=cmex9 at 8 pt%
\font\sixsyex=cmex9 at 6pt%

% sans serif a 9 pt
\font\niness=cmss9%
%\font\eightss=cmss8%
\font\sixss=cmss8 at 6pt%
\font\fivess=cmss8 at 5pt%

% sans serif et math italiques a 9 pt
\font\ninessit=cmssi9%
\font\sixssit=cmssi8 at 6pt%
\font\fivessit=cmssi8 at 5pt%

% math romain italique
\font\ninemi=cmmi9% sixmi est deja definie
\font\fivemi=cmmi5%

% math romain
\font\ninerm=cmr9%
\font\fiverm=cmr5%

% symboles a 9 pt
\font\ninesy=cmsy9% le sixsy est deja defini
\font\fivesy=cmsy5%

% symboles etendus a 9 pt
\font\ninesyex=cmex9% le sixsyex est deja defini
\font\fivesyex=cmex9 at 5pt%

% ams euler fracture
\font\eleveneufm=eufm10 at 11pt%
\font\nineeufm=eufm9%
\font\eighteufm=eufm8%
\font\sixeufm=eufm6%
\font\fiveeufm=eufm5%

% ams blackboard
\font\elevenbb=msbm10 at 11pt%
\font\ninebb=msbm9%
\font\eigthbb=msbm8%
\font\fivebb=msbm5%
\font\sixbb=msbm6%

% ams pour les signes \leq et \geq
\font\elevenmsam=msam10 at 11pt%
\font\ninemsam=msam9%
\font\eigthmsam=msam8%
\font\sixmsam=msam6%
\font\fivemsam=msam5%

% autres fontes

\font\elevenss=cmss12 at 11pt%
\font\elevenssi=cmssi12 at 11pt%
\font\fifteencmss=cmss12 at 15pt%
\font\twentycmss=cmss12 at 20pt%
\font\twentyfivecmss=cmss12 at 25pt%
\font\thirtycmss=cmss12 at 30pt%
\font\twelvebsy=cmbsy10 at 11pt% bold caligraphique a 11 pt

\font\twelverm=cmr12%
\font\fchapter=cmbx12 at 24pt%    fonte pour le titre du chapitre
\font\fchapterit=cmmib10 at 24pt% fonte pour titre chapitre en itallique
\font\fsection=cmb10 at 12 pt%    fonte pour les sections.
\font\fsectionmath=cmmib10 at 12pt%
\font\ftitre=cmr10 scaled\magstep2%
\font\twelvemib=cmmib10 at 12 pt%
\font\tenmib=cmmib10 at 9pt%

% definition des nouvelles familles
\newfam\itfam%
\newfam\frakfam%
\newfam\bbfam%
\newfam\msamfam%
\newfam\mathit%
\newfam\rmfam%
\newfam\twelvebmathfam%
\newfam\slfam%

\textfont\twelvebmathfam=\twelvemib%
\scriptfont\twelvebmathfam=\tenmib%
\def\twelvebmit{\fam\twelvebmathfam}%

% il faudra redefinir le signe < en 9pt sans serif
%\catcode`<=\active
%\catcode`/=\active
%\catcode`>=\active

% composition en 11pt romain

\def\elevenpointsrm{%
% fam 0 : texte romain en 11pt
\textfont0=\elevenrm%
\scriptfont0=\eightrm%
\scriptscriptfont0=\sixrm%
\global\def\rm{\fam0\elevenrm}
\global\def\sl{\fam\slfam\elevensl}% \sl is family 5
\textfont\slfam=\elevensl
% fam 1 : math romain italique
\textfont1=\elevenmi%
\scriptfont1=\eightmi%
\scriptscriptfont1=\sixmi%
\global\def\mit{\fam1}%
% on redefinit les lettres qu on redefinit en 9 pt sans serif
% sinon, en passant de 9 pt ss a 11 pt rm ces lettres ne seraient
% pas changees
\mathchardef\Gamma"0000%
\mathchardef\Delta"0001%
\mathchardef\Theta"0002%
\mathchardef\Lambda"0003%
\mathchardef\Xi"0004%
\mathchardef\Pi"0005%
\mathchardef\Sigma"0006%
\mathchardef\Upsilon"0007%
\mathchardef\Phi"0008%
\mathchardef\Psi"0009%
\mathchardef\Omega"000A%
\mathchardef\alpha"010B%
\mathchardef\beta"010C%
\mathchardef\gamma"010D%
\mathchardef\delta"010E%
\mathchardef\epsilon"010F%
\mathchardef\zeta"0110%
\mathchardef\eta"0111%
\mathchardef\theta"0112%
\mathchardef\iota"0113%
\mathchardef\kappa"0114%
\mathchardef\lambda"0115%
\mathchardef\mu"0116%
\mathchardef\nu"0117%
\mathchardef\xi"0118%
\mathchardef\pi"0119%
\mathchardef\rho"011A%
\mathchardef\sigma"011B%
\mathchardef\tau"011C%
\mathchardef\upsilon"011D%
\mathchardef\phi"011E%
\mathchardef\chi"011F%
\mathchardef\psi"0120%
\mathchardef\omega"0121%
\mathchardef\varepsilon"0122%
\mathchardef\vartheta"0123%
\mathchardef\varphi"0124%
\mathchardef\varrho"0125%
\mathchardef\varsigma"0126%
\mathchardef\varphi"0127%
\mathcode`<="313C%
\mathcode`/="013D%
\mathcode`>="313E%
%\mathchardef <"313C
%\mathchardef /"313D
%\mathchardef >"313E
\mathchardef\partial"0140%
%
% fam 2 : symboles
\textfont2=\elevensy%
\scriptfont2=\eigthsy%
\scriptscriptfont2=\sixsy%
%
% fam 3 : symboles etendus
\textfont3=\elevensyex%
\scriptfont3=\eightsyex%
\scriptscriptfont3=\sixsyex%
%
% fam itfam : romain italique
\textfont\itfam=\elevenit%
\def\it{\fam\itfam\elevenit}%
%
% fam frakfam : ams euler frakture
\textfont\frakfam=\eleveneufm%
\scriptfont\frakfam=\eighteufm%
\scriptscriptfont\frakfam=\sixeufm%
\def\frak{\fam\frakfam\eleveneufm}%
%
% fam bbfam : ams blackboard
\textfont\bbfam=\elevenbb%
\scriptfont\bbfam=\eigthbb%
\scriptscriptfont\bbfam=\sixbb%
%
%\global\def\bb{\fam\bbfam\elevenbb}
\def\bb{\fam\bbfam\elevenbb}%
%
% fonte ams pour \leq et \geq
\textfont\msamfam=\elevenmsam%
\scriptfont\msamfam=\eigthmsam%
\scriptscriptfont\msamfam=\sixmsam%
\def\msam{\msamfam\elevenmsam}%
% redefinition de \leq et \geq%
\mathchardef\leq"3\hexdigit\msamfam 36%
\mathchardef\geq"3\hexdigit\msamfam 3E%
\rm%
}

% et on compose de suite, pour assurer que tous les caracteres
% speciaux ci apres seront bien en 11pt

\elevenpointsrm%

% Composition en 9 pt sans serif

\def\ninepointsss{%
%
% le romain devient du sans serif a 8 pt
\textfont0=\niness%
\scriptfont0=\sixss%
\scriptscriptfont0=\fivess%
\def\rm{\fam0\niness}%
%
% fam 1 : math romain italique
%\textfont1=\ninemi
%\scriptfont1=\sixmi
%\scriptscriptfont1=\fivemi
%\def\mit{\fam1}
\textfont1=\ninessit%
\scriptfont1=\sixssit%
\scriptscriptfont1=\fivessit%
\def\mit{\fam1}%
%
% on redefini alors les symbole grecs qui ne sont pas dans la fonte itallique
\textfont\mathit=\ninemi%
\scriptfont\mathit=\sixmi%
\scriptscriptfont\mathit=\fivemi%
\textfont\rmfam=\ninerm%
\scriptfont\rmfam=\sixrm%
\scriptscriptfont\rmfam=\fiverm%
\mathchardef\Gamma"0000% % ou \mathchardef\Gamma"0\hexdigit\rmfam 00
                         % pour avoir du greque avec serif - Attention,
                         % les minuscules greques sont avec serif
\mathchardef\Delta"0001%
\mathchardef\Theta"0002%
\mathchardef\Lambda"0003%
\mathchardef\Xi"0004%
\mathchardef\Pi"0005%
\mathchardef\Sigma"0006%
\mathchardef\Upsilon"0007%
\mathchardef\Phi"0008%
\mathchardef\Psi"0009%
\mathchardef\Omega"000A%
\mathchardef\alpha"0\hexdigit\mathit 0B%
\def\alpha{\vbox{\hbox to 6.5pt{\kern .5pt\includegraphics{Graph/alpha_ss.ps}\hfill}}}%
\mathchardef\beta"0\hexdigit\mathit 0C%
\mathchardef\gamma"0\hexdigit\mathit 0D%
\mathchardef\delta"0\hexdigit\mathit 0E%
\mathchardef\epsilon"0\hexdigit\mathit 0F%
\def\epsilon{\vbox{\hbox to 5pt{\includegraphics{Graph/epsilon_ss.ps}\hfill}}}%
\mathchardef\zeta"0\hexdigit\mathit 10%
\mathchardef\eta"0\hexdigit\mathit 11%
\mathchardef\theta"0\hexdigit\mathit 12%
\mathchardef\iota"0\hexdigit\mathit 13%
\mathchardef\kappa"0\hexdigit\mathit 14%
\mathchardef\lambda"0\hexdigit\mathit 15%
\mathchardef\mu"0\hexdigit\mathit 16%
\mathchardef\nu"0\hexdigit\mathit 17%
\mathchardef\xi"0\hexdigit\mathit 18%
\mathchardef\pi"0\hexdigit\mathit 19%
\def\pi{\vbox{\hbox to 6pt{\includegraphics{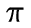}\hfill}}}%
\mathchardef\rho"0\hexdigit\mathit 1A%
\mathchardef\sigma"0\hexdigit\mathit 1B%
\def\sigma{\vbox{\hbox to 6pt{\includegraphics{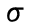}\hfill}}}%
\mathchardef\tau"0\hexdigit\mathit 1C%
\def\tau{\vbox{\hbox to 4.5pt{\includegraphics{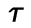}\hfill}}}%
\mathchardef\upsilon"0\hexdigit\mathit 1D%
\mathchardef\phi"0\hexdigit\mathit 1E%
\mathchardef\chi"0\hexdigit\mathit 1F%
\mathchardef\psi"0\hexdigit\mathit 20%
\def\psi{\vbox{\hbox to 6.5pt{\includegraphics{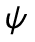}\hfill}}}%
\mathchardef\omega"0\hexdigit\mathit 21%
\def\omega{\vbox{\hbox to 7pt{\includegraphics{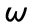}\hfill}}}%
\mathchardef\varepsilon"0\hexdigit\mathit 22%
\mathchardef\vartheta"0\hexdigit\mathit 23%
\mathchardef\varphi"0\hexdigit\mathit 24%
\mathchardef\varrho"0\hexdigit\mathit 25%
\mathchardef\varsigma"0\hexdigit\mathit 26%
\mathchardef\varphi"0\hexdigit\mathit 27%
%\mathchardef <"3\hexdigit\mathit 3C%
%\mathchardef /"3\hexdigit\mathit 3D%
%\mathchardef >"3\hexdigit\mathit 3E%
\mathcode`,="6\hexdigit\mathit 3B%
\mathcode`.="6\hexdigit\mathit 3A%
\mathcode`<="3\hexdigit\mathit 3C%
\mathcode`/="0\hexdigit\mathit 3D%
\mathcode`>="3\hexdigit\mathit 3E%
\mathchardef\partial"0\hexdigit\mathit 40%
\def\partial{\vbox{\hbox to 6pt{\hskip 1pt\includegraphics{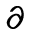}\hfill}}}%
%
% fam 2 : symboles
\textfont2=\ninesy%
\scriptfont2=\sixsy%
\scriptscriptfont2=\fivesy%
%
% fam 3 : symboles etendus
\textfont3=\ninesyex%
\scriptfont3=\sixsyex%
\scriptscriptfont3=\fivesyex%
%
% fam bbfam : ams blackboard
\textfont\bbfam=\ninebb%
\scriptfont\bbfam=\sixbb%
\scriptscriptfont\bbfam=\fivebb%
%
%\global\def\bb{\fam\bbfam\elevenbb}
\def\bb{\fam\bbfam\ninebb}%
%
% fonte ams pour \leq et \geq
\textfont\msamfam=\ninemsam%
\scriptfont\msamfam=\sixmsam%
\scriptscriptfont\msamfam=\fivemsam%
\def\msam{\msamfam\elevenmsam}%
% redefinition de \leq et \geq
\mathchardef\leq"3\hexdigit\msamfam 36%
\mathchardef\geq"3\hexdigit\msamfam 3E%
\def\calD{\vbox{\hbox to 7pt{\includegraphics{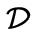}\hfill}}}%
\def\RR{\vbox{\hbox to 7pt{\includegraphics{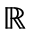}\hfill}}}%
\def\epsilonsmall{\vbox{\hbox to 3.5pt{\includegraphics{Graph/epsilon_6_ss.ps}\hfill}}}%
\def\tausmall{\vbox{\hbox to 3pt{\includegraphics{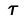}\hfill}}}%
\rm%
}

% composition en 9pt romain

\def\ninepointsrm{%
% fam 0 : texte romain en 9pt
\textfont0=\ninerm%
\scriptfont0=\sixrm%
\scriptscriptfont0=\fiverm%
\global\def\rm{\fam0\ninerm}
%
% fam 1 : math romain italique
\textfont1=\ninemi%
\scriptfont1=\sixmi%
\scriptscriptfont1=\fivemi%
\global\def\mit{\fam1}%
% on redefinit les lettres qu on redefinit en 9 pt sans serif
% sinon, en passant de 9 pt ss a 9 pt rm ces lettres ne seraient
% pas changees
\mathchardef\Gamma"0000%
\mathchardef\Delta"0001%
\mathchardef\Theta"0002%
\mathchardef\Lambda"0003%
\mathchardef\Xi"0004%
\mathchardef\Pi"0005%
\mathchardef\Sigma"0006%
\mathchardef\Upsilon"0007%
\mathchardef\Phi"0008%
\mathchardef\Psi"0009%
\mathchardef\Omega"000A%
\mathchardef\alpha"010B%
\mathchardef\beta"010C%
\mathchardef\gamma"010D%
\mathchardef\delta"010E%
\mathchardef\epsilon"010F%
\mathchardef\zeta"0110%
\mathchardef\eta"0111%
\mathchardef\theta"0112%
\mathchardef\iota"0113%
\mathchardef\kappa"0114%
\mathchardef\lambda"0115%
\mathchardef\mu"0116%
\mathchardef\nu"0117%
\mathchardef\xi"0118%
\mathchardef\pi"0119%
\mathchardef\rho"011A%
\mathchardef\sigma"011B%
\mathchardef\tau"011C%
\mathchardef\upsilon"011D%
\mathchardef\phi"011E%
\mathchardef\chi"011F%
\mathchardef\psi"0120%
\mathchardef\omega"0121%
\mathchardef\varepsilon"0122%
\mathchardef\vartheta"0123%
\mathchardef\varphi"0124%
\mathchardef\varrho"0125%
\mathchardef\varsigma"0126%
\mathchardef\varphi"0127%
\mathcode`<="313C%
\mathcode`/="013D%
\mathcode`>="313E%
%\mathchardef <"313C
%\mathchardef /"313D
%\mathchardef >"313E
\mathchardef\partial"0140%
%
% fam 2 : symboles
\textfont2=\ninesy%
\scriptfont2=\sixsy%
\scriptscriptfont2=\fivesy%
%
% fam 3 : symboles etendus
\textfont3=\ninesyex%
\scriptfont3=\sixsyex%
\scriptscriptfont3=\fivesyex%
%
% fam itfam : romain italique
\textfont\itfam=\nineit%
\def\it{\fam\itfam\sixit}%
%
% fam frakfam : ams euler frakture
\textfont\frakfam=\nineeufm%
\scriptfont\frakfam=\sixeufm%
\scriptscriptfont\frakfam=\fiveeufm%
\def\frak{\fam\frakfam\nineeufm}%
%
% fam bbfam : ams blackboard
\textfont\bbfam=\ninebb%
\scriptfont\bbfam=\sixbb%
\scriptscriptfont\bbfam=\fivebb%
%
%\global\def\bb{\fam\bbfam\ninebb}
\def\bb{\fam\bbfam\ninebb}%
%
% fonte ams pour \leq et \geq
\textfont\msamfam=\ninemsam%
\scriptfont\msamfam=\sixmsam%
\scriptscriptfont\msamfam=\fivemsam%
\def\msam{\msamfam\ninemsam}%
% redefinition de \leq et \geq%
\mathchardef\leq"3\hexdigit\msamfam 36%
\mathchardef\geq"3\hexdigit\msamfam 3E%
\rm%
}

%%%%%%%%%%%%%% Format pour les def, th., prop, cor., etc. %%%%%%%%%

\def\state#1{\noindent{\elevenss#1}\quad}
\def\stateit#1{\noindent{\elevenssi#1}\quad}

%%%%%%%%%%%%%% Notes marginales %%%%%%%%%%%%%%%%%%%%%%%%%%%%%%%%%%%

\def\note#1{%
  \vadjust{%
    \setbox1=\vtop{%
      \hsize 2.5cm\parindent=0pt\sixrm\baselineskip=9pt%
      \rightskip=4mm plus 4mm\raggedright#1%
      }%
    \hbox{\kern-4.4cm\smash{\box1}\hfil}%
    }%
  }

%%%%%%%%%% Creation de la table des matieres %%%%%%%%%%%%%%%%%%%%
%% En debut de chapitre mettre la commande 
%% \chapter{numero_du_ch}{nom_du_ch} et en debut de section,
%% \section{numero_de_section}{nom_de_section}
%% ces macros cree un fichier .cnt qui contient un fichier TeX
%% lui meme contenant les instructions pour la table des matieres
%% Les noms et numeros de chapitre et section sont reutilises
%% pour faire les hauts de pages

\newtoks\chapternumber \chapternumber={}
\newtoks\sectionnumber \sectionnumber={}
\newtoks\chaptername \chaptername={}
\newtoks\sectionname \sectionname={}
%\newwrite\cnt  % table of content
%\immediate\openout\cnt=sommaire.ttex
%\write\cnt{\string\centerline{\bf Sommaire}}
%\write\cnt{\string\bigskip}
\def\chapter#1#2{%
  \ifodd\pageno\relax\else\noheadlinetrue ~ \vfill\eject\noheadlinefalse\fi
  \chapternumber={#1}\chaptername={#2}
%  \sectionnumber={}\sectionname={}
  \sectionnumber=\chapternumber\sectionname=\chaptername
  \mark{\the\sectionnumber .\quad\the\sectionname}
  \titlepagetrue
%  \write\cnt{\string\bigskip} 
%  \write\cnt{\string\noindent \space Chapitre #1}
%  \write\cnt{\string\medskip}
%  \write\cnt{\string\noindent \string{\string\bf \space #2\string} 
%	     \string\quad \string\dotfill \string\quad \space\the\pageno}
%  \write\cnt{}
  }

\def\section#1#2{\sectionnumber={#1}\sectionname={#2}
%  \noindent\mark{\the\chapternumber .\the\sectionnumber\quad\sectionname}
%  \write\cnt{\string\noindent \string\quad \string\S#1. #2 
%             \string\quad \string\dotfill \string\quad \space\the\pageno}
%  \write\cnt{}
  \mark{\the\sectionnumber .\quad\the\sectionname}
  }

%%%%%%%%%% Index des notations %%%%%%%%%%%%%%%%%%%%%%%%%%%%%%%%%

\newif\ifproofmode \proofmodetrue
%%\newwrite\ntntexta
%%\newwrite\ntntextb
%%\newtoks\str
%%\immediate\openout\ntntexta=ntn1_5.text
%%\immediate\openout\ntntextb=ntn6etc.text
%%\def\ntn#1{\toks1={#1}
%%  \ifnum\the\chapternumber<6%
%%    \write\ntntexta{#1 \space \the\sectionnumber , p.\the\pageno}%
%%  \else%
%%    \write\ntntextb{\the\sectionnumber\space #1 , p.\the\pageno}%
%%  \fi%
%%  \ifproofmode\note{#1}\else\relax\fi\toks2=\toks1}%

\def\ntn#1{\relax}

%%%%%%%%%% Haut et bas de pages %%%%%%%%%%%%%%%%%%%%%%%%%%%%%%%%
\newif\iftitlepage \titlepagetrue  % est utilise pour la creation de la table
\newif\ifnoheadline \noheadlinetrue

\def\rightheadline{\ninepointsss\firstmark\hfill\folio}
\def\leftheadline{\ninepointsss\folio\hfill%
                  Chapter \the\chapternumber .\quad\the\chaptername}

\headline{\ifnoheadline\hfill%
	  \else%
	    \iftitlepage\hfill\global\titlepagefalse%
	    \else%
	      \vadjust{\hbox{\vrule height1pt depth2pt width0pt}\hrule}%
	      \ifodd\pageno\rightheadline\else\leftheadline\fi%
	    \fi%
	  \fi}
\footline{\hfill}

%%%%%%%%%%%%%%%%% Macros de mise en page %%%%%%%%%%%%%%%%%%%%%%%%%%

\def\cut{{}\hfill\cr \hfill{}}  % pour couper une formule sur 2 ligne
				% cf. Seroul, p151
\long\def\do_not_print#1{\relax}
\long\def\finetune#1{#1}
\def\color#1{
}

%%%%%%%%%%%%%%%% Macros graphiques %%%%%%%%%%%%%%%%%%%%%%%%%

\def\anote#1#2#3{\smash{\kern#1in{\raise#2in\hbox{#3}}}%
  \nointerlineskip}	% permet de mettre des anotations
			% dans les graphiques inclus avec
			% \special
% macros m(ultiple)ddot, m(ultiple)vdots, m(ultiple)cdots
\newcount\n
\newcount\depth
\def\mddots#1#2{%  met #1 ddots dans une boite elevee de #2pt
  \depth=#2%
  \setbox1=\hbox{$\mathinner{\mkern1mu\raise\depth pt%
  \vbox{\kern\depth pt\hbox{.}}%
  \n=1%
  \loop\ifnum\n<#1 \advance\n by 1%
    \advance\depth by -3%
    \mkern2mu\raise\depth pt\hbox{.}%
  \repeat%
  }$}
  \box1
  }
\def\mvdots#1{% met #1 \cdot verticaux
  \vbox{\baselineskip=4pt\lineskiplimit=0pt
  \kern6pt
  \n=0
  \loop\ifnum\n<#1 \advance\n by 1%
    \hbox{.}
  \repeat }}
\def\mcdots#1{% met #1 \cdot horizontaux
  \mathinner{
  \n=0
  \loop\ifnum\n<#1 \advance\n by 1%
    \cdotp
  \repeat }}
%%%%%%%%%%%%%%%%%%%%%%%%%%%%%%%%%%%%%%%%%%%%%%%%%%%%%%%%%%%%

\catcode`@=11
\def\big#1{{\hbox{$\left#1\vbox to9.5\p@{}\right.\n@space$}}} % 8.5pt in plain
\catcode`@=12

%%%%%%%%%%%%%%%%%%%%%%%%%%%%%%%%%%%%%%%%%%%%%%%%%%%%%%%%%%%%

\def\<{\langle}
\def\>{\rangle}
\def\cjj{C_{j,j}}
\def\cct{(C+C^\T)}
\def\closure{{\rm cl}}
\def\Cov{{\rm Cov}}
\def\cXX{\< CX,X\>}
\def\cxx{\< Cx,x\>}
\def\d{{\rm d}}				% differentielle ATTENTION \d existe en plain
\def\D{{\rm D}}				% differentielle
\def\DA{{\cal D}_A}
\def\da1{{\cal D}_{A_1}}
\def\dat{{\cal D}_{A_t}}
\def\det{{\rm det}}
\def\diag{{\rm diag}}			% diagonal matrix
\def\diam{{\rm diam}}                   %  diameter
\def\doverdt{{\d\over \d t}}
\def\dispar{\displaystyle\partial}
\def\DSLnR{{\rm DSL}(n,\RR)}
\def\Eij{E^{i,j}}
\def\Eik{E^{i,k}}
\def\Ejl{E^{j,l}}
\def\Ekl{E^{k,l}}
\def\Elm{E^{l,m}}
\def\Emi{E^{m,i}}
\def\eIx{e^{-I(x)}}
\def\glnr{{\rm GL}(n,\RR)}
\def\Huv{{H_{u,v}}}
\def\qed{~~~{\vrule height .9ex width .8ex depth -.1ex}}
\def\qed{{\vrule height .9ex width .8ex depth -.1ex}}
\def\Id{{\rm Id}}                       % identity function
\def\imhht{(\Id-hh^\T )}
\def\indic{{\rm I}}			% fct indicatrice.
\def\inj{{\rm inj}}
\def\interior{{\rm int}} 
\def\inv{{\rm inv}}			% inverse
\def\iphth{(\Id+h^\T h)}
\def\iphht{(\Id+hh^\T )}
\def\llst{\log\log\sqrt{t}}
\def\MI{M_\calI}
\def\MnR{{\rm M}(n,\RR )}		% n-n matrices

\def\pinM{{p\in M}}
\def\pisa{\Phi^\leftarrow\!\circ S_\alpha}
\def\Proj{{\rm Proj}}
\def\psixt{\psi (x,t)}
\def\px{\< p,x\>}
\def\Ricc{{\rm Ricc}}
\def\rx{_{\{ x\}^\perp}}		% restriction a {x}^\perp
\def\Saip{S_\alpha^\leftarrow\!\circ\Phi}
\def\salt{\sqrt{\alpha\log t}}
\def\sbullet{{\scriptscriptstyle\bullet}}
\def\frakSigma{{\frak S}}
\def\sign{{\rm sign}}			% sign function
\def\SLnR{{\rm SL}(n,\RR)}
\def\Sn1{S_{n-1}}
\def\SOnR{{\rm SO}(n,\RR)}
\def\span{{\rm span}}
\def\sqrtt{\sqrt{t}}
\def\T{{\rm T}}				% transposition ATTENTION \t existe en plain
\def\trace{{\rm tr}}
\def\ttoi{t\to\infty}
\def\Var{{\rm Var}}                     %variance
\def\vol{{\rm Vol}}
\def\xp{\< x,p\>}
\def\Xp{\< X,p\>}

\def\CC{{\bb C}}
\def\NN{{\bb N}} 	% \def\NN{{\tenmathbb N}}
\def\QQ{{\bb Q}}
\def\RR{{\bb R}}                        % real numbers
\def\ZZ{{\bb Z}}

\def\bfcalA{{\twelvebsy A}}

\def\calA{{\cal A}}
\def\calB{{\cal B}}
\def\calC{{\cal C}}
\def\calD{{\cal D}}
\def\calI{{\cal I}}
\def\calJ{{\cal J}}

\def\calM{{\cal M}}
\def\calN{{\cal N}}
\def\calR{{\cal R}}
\def\calS{{\cal S}}
\def\calT{{\cal T}}
\def\calX{{\cal X}}

\mathchardef\tbftau="0\hexdigit\twelvebmathfam 1C

%\bye

\proofmodefalse
\setbox1=\vbox to \vsize{
  \kern 9.1in
  \includegraphics{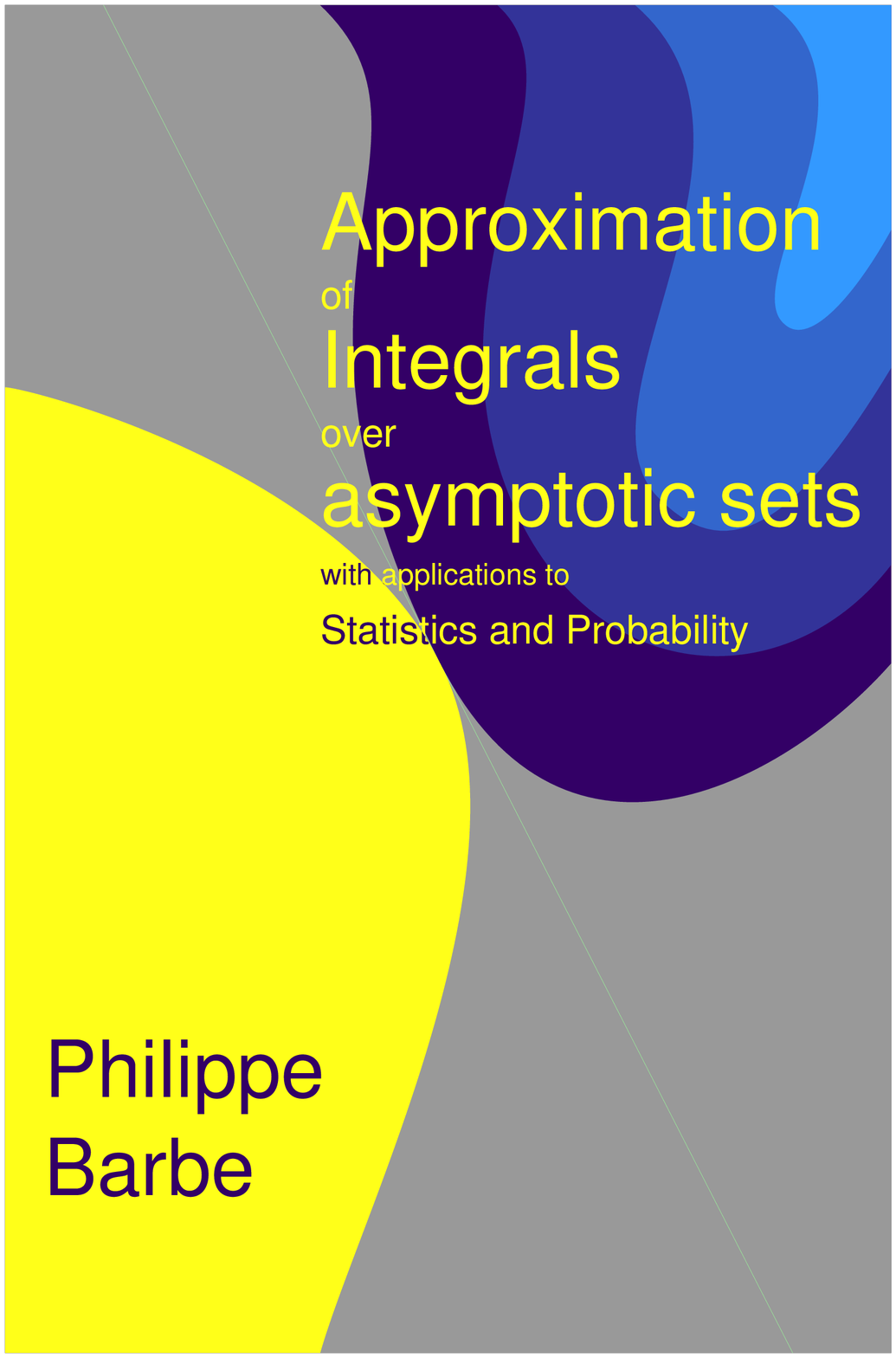}
  \vss}
\noindent
\line{\hskip-1.4in\box1\hfill}\vfill\eject
\pageno=0
\noheadlinetrue
~\vfill\eject  %% A garder si psnup est utilise
\pageno=0
\setbox1=\hbox to 6in{\leavevmode\kern -.75in\vbox to 9.25in{\vfill%
  \hsize=6in%
  \kern 17.6in%
%  \special{psfile=couvbkgd.ps}%
  \vss}%
  \hss}
\wd1=0pt\ht1=0pt\dp1=0pt

\noindent\box1
~
\vfil
%\font\thirtycmss=cmss12 at 30pt
%\font\thwentyfivecmss=cmss12 at 25pt
%\font\fifteencmss=cmss12 at 15pt
\color{0 0 .2}
\line{\hfill\twentyfivecmss Philippe Barbe}
\vskip 25pt
\line{\hfill\thirtycmss APPROXIMATION}
\vskip 10pt
\line{\hfill\twentycmss of}
\vskip 10pt
\line{\hfill\thirtycmss INTEGRALS} 
\vskip 10pt
\line{\hfill\twentycmss over}
\vskip 10pt
\line{\hfill\thirtycmss ASYMPTOTIC SETS}
\vskip 20pt
\line{\hfill\fifteencmss with applications to}
\vskip 10pt
\line{\hfill\twentycmss Statistics and Probability}
\nocolor
\vfil\eject
~\vfil\eject %pour commencer sur une page de droite!
~\vfill
\line{Philippe Barbe\hfill}
\line{CNRS\hfill}
\line{90 Rue de Vaugirard\hfill}
\line{75006 PARIS\hfill}
\line{FRANCE\hfill}
\vskip .5in
\line{\copyright \ Ph. Barbe\hfill}
\vfill\eject
~\vfill\eject
%\bye

%%\input preface.tex
%%~ \vfill\eject
~\vskip 28pt
{\fchapter\hfill Contents}
 
\vskip 54pt  

%\hskip -\parindent {\twelverm Preface}

%\medskip

\hskip -\parindent {\twelverm 1. Introduction \hfill 1}

  Notes \hfill 16

\medskip

\hskip -\parindent {\twelverm 2. The logarithmic estimate \hfill 19}

  Notes \hfill 23

\medskip

\hskip -\parindent {\twelverm 3. The basic bounds \hfill 25}

 1. The normal flow and the normal foliation \hfill 25

 2. Base manifolds and their orthogonal leaves \hfill 36

 Notes\hfill 40

\medskip

\hskip -\parindent {\twelverm 4. Analyzing the leading term for some 
smooth sets \hfill 43}

 1. Quadratic approximation of $\tau_A$ near a dominating manifold\hfill 44

 2. Approximation for $\det \calA (t,v)$\hfill 46

 3. What should the result be? \hfill 50

\medskip

\hskip -\parindent {\twelverm 5. The asymptotic formula \hfill 53}

  Notes \hfill 65

\medskip

\hskip -\parindent {\twelverm 6. Asymptotic for sets translated towards 
  infinity \hfill 69}

  Notes \hfill 77

\medskip

\hskip -\parindent {\twelverm 7. Homothetic sets, homogeneous $I$ and
Laplace's method \hfill 79}

  Notes \hfill 88

\medskip

\hskip -\parindent {\twelverm 8. Quadratic forms of random vectors \hfill 89}

 1. An example with light tail distribution \hfill 89

 2. An example with heavy tail distribution \hfill 93

 3. Heavy tail and degeneracy \hfill 142

 Notes \hfill 162

\medskip

\hskip -\parindent {\twelverm 9. Random linear forms \hfill 165}

 1. Some results on convex sets \hfill 166

 2. Example with light tails \hfill 177

 3. Example with heavy tails \hfill 185

 Notes\hfill 193

\medskip

\hskip -\parindent {\twelverm 10. Random matrices \hfill 197}

 1. Random determinants, light tails \hfill 200

 2. Random determinants, heavy tails \hfill 210

 3. Geometry of the unit ball of $\MnR$ \hfill 218

 4. Norms of random matrices \hfill 229

 Notes \hfill 236

\medskip

\hskip -\parindent {\twelverm 11. Finite sample results for 
autoregressive processes \hfill 239}

 1. Background on autoregressive processes \hfill 239

 2. Autoregressive process of order $1$ \hfill 243

 3. Autoregressive process of arbitrary order \hfill 255

 Notes \hfill 271

\medskip

\hskip -\parindent {\twelverm 12. Suprema of some stochastic processes \hfill 273}

 1. Maxima of processes and maxima of their variances \hfill 273
 
 2. Asymptotic expansions for the tail of the supremum of

 \leavevmode\kern 1pt\phantom{2. }Gaussian processes can be arbitrary bad \hfill 276

 3. Maximum of nonindependent Gaussian random variables \hfill 280

 4. The truncated Brownian bridge \hfill 282

 5. Polar processes on boundary of convex sets \hfill 293

 Notes \hfill 296

\medskip

\hskip -\parindent {\twelverm Appendix 1. Gaussian and Student tails\hfill 297}

\medskip

\hskip -\parindent {\twelverm Appendix 2. Exponential map\hfill 303}

\medskip

\hskip -\parindent {\twelverm References \hfill 305}

\medskip

\hskip -\parindent {\twelverm Notation\hfill 311}

\medskip

\hskip -\parindent {\twelverm Postface\hfill 315}

\vfill\eject
~\vfill\eject~\vfill\eject
\noheadlinefalse
\pageno=1
~\vskip 28pt
\chapter{1}{Introduction}
{\fchapter\hfill 1.\quad Introduction}

\vskip 54pt\elevenrm

\noindent Computing integrals is a 
basic need in many areas
of mathematics and applied sciences. It is a common experience
that some integrals can be calculated ``by hand'', obtaining
a ``closed formula'', while others cannot be reduced to a
simple expression. This fact is at the origin of numerous
approximation techniques for integrals, such as
quadratures, series expansions, Monte-Carlo methods,
etc. The purpose of these notes is to start what seems to be
a new line of investigations in studying the behavior of
integrals of the form\ntn{$A$}
$$
  \int_A f(x) \d x \, , \qquad A\subset \RR^d \, ,
  \eqno{(1.1)}
$$
for sets $A$ far away from the origin. In this integral,
$\d x$ is the Lebesgue measure on $\RR^d$, and the 
function $f$ is integrable over $\RR^d$.

Such integral is in general very difficult 
to compute numerically. 
Indeed, in interesting applications $d$ is between 10 and 50
say, and quadrature methods are essentially out of
considerations. On the other hand, as $A$ is far away from 
the origin, the integral is small --- the interesting range
for some applications is when the integral is of order
$10^{-2}$, $10^{-3}$, or even smaller --- and 
standard Monte-Carlo methods
fail to provide a good evaluation at a cheap
cost.

The motivation for this study comes mainly from applications
in statistics and probability theory. However,
the techniques developed should be useful in other areas
of mathematics as well, wherever such integrals arise.

Before going further, let us show some consequences of our
main approximation result in probability and statistics.
We hope that the reader will find in these examples good
motivation to continue reading.

Consider the following two density functions on $\RR\,$,
$$\eqalign{
  w_\alpha(x)%
  &=K_{w,\alpha}\exp\bigl(-|x|^\alpha/\alpha\bigr) \, ,\cr
  s_\alpha(x)%
  &=K_{s,\alpha}
  \Big(1+{x^2\over\alpha}\Big)^{-(\alpha+1)/2} \, , \cr
  }
$$
where\ntn{$K_{w,\alpha}$, ch1}\ntn{$K_{s,\alpha}$, ch1}
$$
  K_{w,\alpha} 
  = { 1\over\vbox to .2in{} 2\alpha^{{1\over \alpha}-1}
      \Gamma \Big({\displaystyle 1\over\displaystyle\alpha}\Big) }
  \hbox{\quad and \quad} 
  K_{s,\alpha} = {\displaystyle 1\over\displaystyle\sqrt{\pi\alpha}}
  {\vtop to .15in{}\displaystyle\Gamma \Big(
    {\displaystyle\alpha +1\over\displaystyle 2}
    \Big)\over\vbox to .2in{}\displaystyle\Gamma \Big(
    {\displaystyle\alpha\over \displaystyle2}\Big)} \, ;
$$
so $w_\alpha$ and $s_\alpha$ integrate to $1$ on
the real line. The density $w_2$ is the
standard normal distribution. The density
$s_\alpha$ is that of a Student distribution. These two 
functions are very different. The symmetric Weibull-like, $w_\alpha$, decays
exponentially at infinity, while the Student one decays
polynomially.

Let us now consider a real $d\times d$ matrix 
$C={(C_{i,j})}_{1\leq i,j\leq d}$ with at least one positive
term on its diagonal.
We write $\< \cdot \, ,\cdot\>$ the standard inner
product on $\RR^d$, and consider the domain
$$
  A_t=\{\, x\in\RR^d : \<Cx,x\>\geq t \,\} \, .
$$
As $t$ tends to infinity, it is easy to see
that the set $A_t=\sqrt{t}A_1$ pulls away from the origin
since $A_1$ does not contain $0$. We will explain the following in
section 8.1.
If $\alpha\neq 2$, then there exists a function --- rather
explicit, and in any case computable --- $c(\cdot )$ of
$\alpha$, $d$ and $C$, such that
$$
  \int_{A_t} \prod_{1\leq i\leq d} w_\alpha (x_i) \d x_1
  \ldots \d x_d \sim 
  c(\alpha,d,C) e^{-t^{\alpha /2}c_2} 
  t^{-(\alpha -2){d\over 4} -{\alpha\over 2} } \, ,
$$
as $t$ tends to infinity.
When $\alpha =2$, the asymptotic behavior of the integral is 
different. Namely, if $k$ is the dimension of 
the eigensubspace associated
to the largest eigenvalue $\lambda$ of $C^\T +C$,
$$
 \int_{A_t} \prod_{1\leq i\leq d} w_2 (x_i) 
 \d x_1 \ldots \d x_d 
 \sim c(2,d,C) e^{-t/\lambda} t^{(k-1)/2}
 \, , 
$$
as $t$ tends to infinity. Comparing the formula for
$\alpha$ equal or different from $2$, we see that the power of $t$,
namely $(\alpha -2)d/4 -\alpha/2$ for $\alpha$ different from $2$ and $(k-1)/2$
for $\alpha$ equal $2$, has a discontinuity in $\alpha =2$, whenever
$k$ is different than $2$. Moreover, we will see in section 8.2 that
$$
  \int_{A_t}\prod_{1\leq i\leq d}s_\alpha(x_i) \d x_1
  \ldots \d x_d \sim
  K_{s,\alpha} \alpha^{(\alpha+1)/2}2 t^{-\alpha /2}
  \sum_{i: C_{i,i}>0}C_{i,i}^{\alpha/2}
$$
as $t$ tends to infinity.
These estimates give approximations
for the tail probability of the quadratic form $\<CX,X\>$
evaluated at a random vector $X$ having independent components
distributed according to $w_\alpha$ or $s_\alpha$.
We will also see what happens
when all the diagonal elements of $C$ are
negative.

Consider now $d=n^2$ for some positive integer $n$, and
write $\MnR$ for the set of all real $n\times n$
matrices. Let
$$
  A_t=\{\,x\in\MnR : \det\, x\geq t\,\}
$$
be the set of all $n\times n$ real matrices with determinant
at least $t$. The set $A_1$ is closed and does not contains
the origin. One sees that $A_t=t^{1/n}A_1$ moves away
from the origin when $t$ increases. We will prove that,
as $t$ tends to infinity,
$$
  \int_{A_t}\prod_{1\leq i,j\leq n}w_\alpha (x_{i,j}) 
    \d x_{i,j} 
  \sim
  \cases{ c\, e^{-t^\alpha/n}t^{(\alpha(n^2-1)-2n^2)/2n}
          & if $\alpha\neq 2$ \cr
         \noalign{\vskip 0.05 in}
          c\, e^{-nt^2/2} t^{(n^2-n-2)/2} 
          & if $\alpha=2$ \cr}
$$
while
$$
  \int_{A_t}\prod_{1\leq i,j\leq n}s_\alpha (x_{i,j}) 
  \d x_{i,j} 
  \sim
  c\, {(\log t )^{n-1}\over t^\alpha} \, .
$$
This gives a tail estimate of the distribution
of the determinant of a random matrix with independent and
identically distributed entries from
a distribution $w_\alpha$ or $s_\alpha$. The constant $c$
depends on $s_\alpha$ or $w_\alpha$ as well as on the dimension
$n$; but it does not depend on $t$ and will be explicit.

Possible applications are numerous. We will deal with
further interesting examples, such as norms of random matrices,
suprema of linear processes, and other related quantities.

Though some specific examples could be derived with other methods
--- in particular those dealing with the Gaussian density integrated
over rather simple sets --- we hope that the reader will enjoy
having a unified framework to handle all those asymptotic problems.
It has the great advantage of providing a systematic approach. It
breaks seemingly complicated problems into much more manageable
ones. It also brings a much better understanding of how the leading
terms come out of the integrals. Through numerical computations we 
will also see that our approximations are accurate enough
to be of practical use in statistics, when dealing with some small
sample problems --- the type of problem for which there are no
systematic tools as far as the author knows, and which
originally motivated this work.

Another purpose of these notes is to investigate asymptotics for
conditional distributions, and Gibbs conditioning. For instance,
consider $n$ points in $\RR^n$, whose coordinates are
independent and identically distributed
random variables. They form a parallelogram. Given, say, the volume of
this parallelogram, these points are no longer independent. What
is their distribution? In general, this seems to be a very
difficult question, except if one stays at a very theoretical
level and does not say much. We will show that it is possible to
obtain some good approximation of this conditional distribution
for large volumes of the parallelogram --- mathematically, as the
volume tends to infinity.

It is important to realize that conditional distributions are
rather delicate objects. To convince the reader of the truth of
this 
assertion, let us consider four elementary examples where
everything can be calculated explicitly and easily.

\bigskip

\state{Example 1.} Let $X$, $Y$ be two independent
standard normal random variables. We are looking for the
distribution of $(X,Y)$ given $X+Y\geq t$, for large values of
$t$. A simple probabilistic argument to obtain a limiting
distribution is as follows. Let $U=(X+Y)/\sqrt 2$ and
$V=(X-Y)/\sqrt 2$. The distribution of $X$ given $X+Y\geq t$ is
that of $(U+V)/\sqrt 2$ given $U\geq t/\sqrt{2}$. One can
check by hand that the conditional distribution of $U/t$
given $U\geq t/\sqrt 2$ converges weakly* to a point mass
at $1/\sqrt 2$. Since $U$ and
$V$ are independent standard normal, and since $(X,Y)$ is exchangeable
conditionally on $X+Y$, the pair $(X,Y)/t$ converges in
probability to $(1/2,1/2)$ conditioned on $X+Y\geq t$. 

Consequently, if we are given $X+Y\geq t$ for large $t$,
we should expect
to observe $X/t\approx Y/t\approx 1/2$. Hence, $X$ and $Y$ are
both large, and about the same order.

\bigskip

\state{Example 2.} Consider the same problem as in
example 1, but with $X$ and $Y$ independent and both exponential.
For any $\alpha\in (0,1)$,
$$
  P\{\, X\leq \alpha t \, | \, X+Y\geq t\,\} =
  { P\{\, X\leq \alpha t \, ; X+Y\geq t\,\} \over 
    P\{\, X+Y\geq t\,\} } \, .
$$
We can explicitly calculate
$$\eqalign{
  P\{\, X\leq \alpha t\, ; X+Y\geq t \,\} 
  & = \int_{x=0}^{\alpha t}\int_{y=t-x}^\infty
    e^{-x}e^{-y} \d x \, \d y
    = \int_0^{\alpha t} e^{-t} \d x \cr
  & = \alpha t e^{-t}  \, .\cr
  }
$$
On the other hand, $X+Y$ has a gamma distribution with mean 
$2$. Integration by parts yields
$$
  P\{\, X+Y\geq t\,\} =\int_t^\infty xe^{-x} \d x 
  = (t+1)e^{-t} \, .
$$
It follows that
$$
  P\{\, X\leq \alpha t \, | \, X+Y\geq t\,\} 
  = \alpha {t\over t+1} \, .
$$
Thus, given $X+Y\geq t$, the distribution of $X/t$
converges to a uniform distribution over $[\, 0,1\, ]$. Again, by
conditional exchangeability, the same 
result holds for $Y/t$. In conclusion,
the distribution of $(X,Y)/t$ conditioned on
$X+Y\geq t$ converges weakly* to a uniform
distribution over the segment $\alpha (1,0) + (1-\alpha ) (0,1)$
in $\RR^2$. So, we do not have the kind of degeneracy that occurs in
Example 1. Conditioned on $X+Y\geq t$, the random
variables $X$ and $Y$ should still
be of the same order of magnitude, $t$, but not necessarily 
close to each others after
the rescaling by $1/t$.

\bigskip

\state{Example 3.} Consider the same problem, with now
$X$ and $Y$ independent, both having a Cauchy distribution.
Elementary calculation shows that
$$
  P\{\, X\geq s\,\} = \int_s^\infty {\d x\over \pi (1+x^2)}
  \sim \int_s^\infty {\d x\over \pi x^2} \sim {1\over \pi s}
  \qquad\hbox{as } s\to\infty \, .
$$
Moreover, the stability of the Cauchy distribution --- see,
e.g., Feller, 1971, \S II.4 ---
asserts that $(X+Y)/2$ has also
a Cauchy distribution. Consequently,
$$\displaylines{\qquad
    \lim_{t\to\infty} P\{ \, X\geq \epsilon t\, ,
      Y\geq \epsilon t\, |\, X+Y\geq t\, \}
  \cut
    \leq \lim_{t\to\infty} 
    {P\{\, X\geq \epsilon t\, ;\, Y\geq \epsilon t\, \} \over
     P\{\, X+Y\geq t\,\}} 
    = 0 \, .
  \qquad}
$$
This proves that we cannot have both $X$ and $Y$ too large as
$t$ tends to infinity. Moreover, for $s$ positive,
$$
  P\{\, X\geq t(1+s)\,|\, X+Y\geq t\,\}
  \sim {P\{\, X\geq t(1+s)\,\}\over 
        P\{\, X+Y\geq t\,\} }
  \sim {1\over 2(1+s)} 
$$
as $t$ tends to infinity. Consequently, given $X+Y\geq t$, the random
vector $(X,Y)/t$ has a distribution which converges weakly* 
to $(P_x+P_y)/2$ where $P_x$ (resp. $P_y$) is a Pareto 
distribution on the $x$-axis (resp. $y$-axis).

\bigskip

These three examples show that we can obtain very different
behaviors, with various degrees of degeneracy for the limiting
conditional distributions. But they also show examples of
nonconceptual proofs. They provide no insight on what is happening. All
the arguments are rather ad hoc. However, it is intuitively clear
that the limiting conditional distribution
has to do with the tail of the distribution of $(X,Y)$.
Specifically, if we write $f(x,y)$ for the density of $(X,Y)$ and
$A_t=\{\, (x,y) : x+y\geq t\,\}$, we are led to consider
expressions of the form
$$
  { \displaystyle\int_{A_t\cap \{\, x : x\leq ct\,\} } 
    f(x,y)\,\d x \,\d y \over
    \displaystyle\int_{A_t} f(x,y)\, \d x \,\d y } \, .
$$
When looking at $\int_{A_t} f(x,y)\, \d x\, \d y$, examples 1, 2, 3 have
very different features. Intuitively, the behavior of the
integral has to do with how the set $A_t$ lies in $\RR^2$ compared
to the level sets of $f(x,y)$. Examples 1, 2 and 3 correspond to the
following pictures.

\bigskip

%\vskip 2mm
%\hbox{\vbox{
% Ga_tech
%\noindent\vbox to 2in{%
%  %\leavevmode
%  \kern 2in
%  \special{psfile=Graph/ch1_fig1.ps}
%  \leavevmode \kern 2.334in 
%  \special{psfile=Graph/ch1_fig2.ps}
%}

\setbox1=\vbox to 2.5in{\ninepointsss
  \hsize=2in
  \kern 2in
  \includegraphics{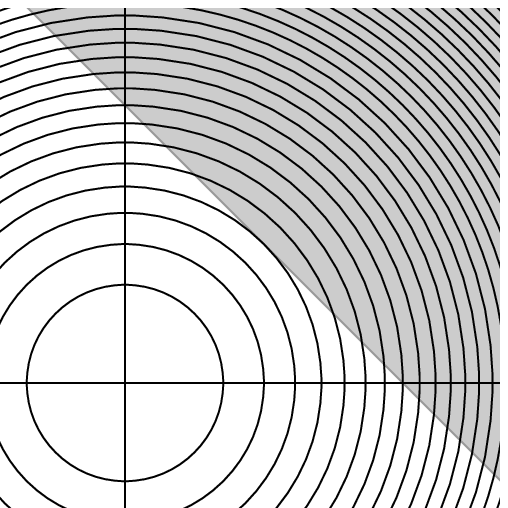}
  \kern .1in
  \centerline{Example 1}\vss}
\setbox2=\vbox to 2.5in{\ninepointsss
  \hsize=2in
  \kern 2in
  \includegraphics{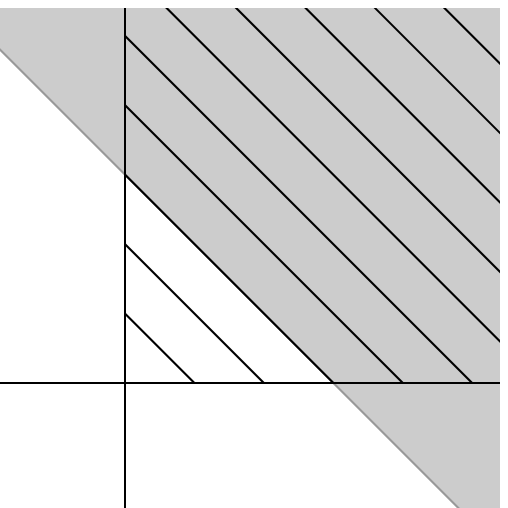}
  \kern .1in
  \centerline{Example 2}\vss}
\setbox3=\vbox to 2.5in{\ninepointsss
  \hsize=2in
  \kern 2in
  \includegraphics{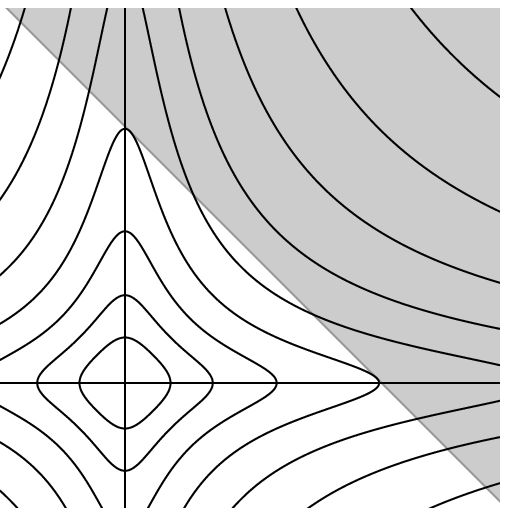}
  \kern .1in
  \centerline{Example 3}\vss}
\line{\box1\hfill\box2}

\finetune{\vfill\eject}

\centerline{\box3}
%}}
% autre
%\hskip 3cm\hbox{\psfig{file=fig0.epsi}}
% Tlse
%\hskip 1.4cm\epsfbox[358 471 549 610]{fig0.epsi}
% puis on met la legend, p, \theta , 0 , et le point (1,0) !
%\vskip -2.3cm \hskip 1.5cm \hskip 1.1cm p
%\vskip 0.7cm \hskip 2.5cm $\theta$
%\vskip 0.3cm\hskip 1.6cm 0 \hskip 1.3cm 1
%\vskip 0.7cm
%}}

In the three pictures, we see the gray shaded set $A_t$, 
its boundary, and the level sets
$$
  \Lambda_c =\{\, (x,y) : f(x,y)=e^{-c}\,\}
$$
for $c=1, 2, 3, \ldots$ The fact that the sets $\Lambda_c$ are
getting closer and closer as $c$ tends to infinity 
in example 1 expresses
the very fast decay of the Gaussian distribution. At the opposite extreme,
in example 3, the level sets are further and further apart because
the distribution has a subexponential decay. These pictures show
why in examples 1 and 3 we obtained some form of degeneracy in
the limiting distribution. In example 1, the density function
is maximal on the boundary of $A_t$ at one point whose two
coordinates are equal. In example 3, the density is maximal
at two points on the boundary of $A_t$. These two points have
one coordinate almost null, while the other coordinate is of
order $t$. But one should be careful, as
rough pictures may not give the right result if they are read
naively in more complicated situations. Our fourth example
is still elementary and illustrates this point.

\bigskip

\state{Example 4.} Consider $(X,Y)$ where $X$ and $Y$ are
independent, both having a Cauchy distribution with density
$s_1$. We are now
interested in the tail distribution of $XY$ and in the
distribution of $(X,Y)$ given $XY\geq t$ for large $t$. By
symmetry, it is enough to consider the same distribution but
adding the conditioning $X\geq 0$ and $Y\geq 0$. 

On $X,Y> 0$, set $U=\log X$ and $V=\log Y$. 
The density of $U$ is 
$$
  f_U(u)=e^us_1(e^u) = {e^u\over \pi (1+e^{2u})}
  \sim {{e^{-u}}\over \pi} \qquad\hbox{ as } u\to\infty \, .
$$
Thus, the level sets of the density of $(U,V)$ far away from the
origin look like those of the exponential distribution of example
2. Hence, when studying the distribution of $U$ given $U+V\geq \log t$
--- i.e., $X$ given $XY\geq t$ --- we could expect to have a
behavior similar to example 2. However, the picture looks as
follows.

\bigskip
\setbox1=\vbox to 2.5in{\ninepointsss
  \hsize=2in
  \kern 2in
  \includegraphics{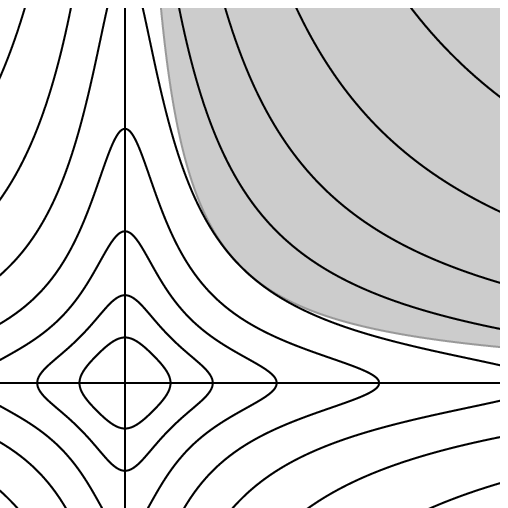}
  \kern .1in
  \centerline{Example 4}\vss}
\centerline{\box1}

\noindent In particular, write $A_t=\{\, (x,y) : xy\geq t\,\}$
for the set on which we integrate the density of $(X,Y)$. 
The boundary $\partial A_t$ has a unique contact point with the
maximal level set of the density intersecting $A_t$, where by
maximal level set we mean
$\Lambda_c$ such that $\Lambda_{c+\epsilon}\cap A_t$ is
nonempty if $\epsilon$ is positive, and is empty if $\epsilon$
is negative. In that aspect, we are close to examples 1 or 3. This
example 4 will show us two more things. First we need to take into
consideration the decay of the density --- which explains why we
will not have degeneracy as in examples 1 or 3.
Second, it shows that
whichever result we can prove will have a lot to do with the
shape of the domain considered and how this domain pulls away from
the level sets of the density.

To see what happens in this fourth example is still rather easy.
We first calculate
$$\eqalignno{
  P\{\, XY\geq t\,\} 
  & = {2\over \pi^2} \int_{xy\geq t\, ,\, x\geq 0\, ,\, y\geq 0}
    {\d y\over 1+y^2}{\d x\over 1+x^2} \cr
  \noalign{\vskip .05in}
  & = {2\over \pi^2} \int_0^\infty 
    \Big( {\pi\over 2} - \arctan {t\over x}\Big) 
    {\d x\over 1+x^2} \, . 
    & (1.2)\cr
  }
$$
It is true that
$$
  {\pi\over 2}-\arctan u
  = \int_u^\infty {\d x\over 1+x^2}
  \sim {1\over u}\qquad\hbox{ as } u\to\infty \, .
$$
However, replacing $(\pi/2) - \arctan (t/x)$ by $x/t$ 
in (1.2) leads to a divergent integral. 
So, we need to be careful. Let $\eta$ be a positive real number.
There exists a positive $\epsilon$ such that for any 
$u> 1/\epsilon$,
$$
  {\pi\over 2}-\arctan u\in 
  \Big[\, {1-\eta\over u} , {1+\eta\over u}\,\Big] \, .
$$
Consequently,
$$\displaylines{
  \qquad
  \int_0^{\epsilon t} 
  \Big({\pi\over 2}-\arctan {t\over x}\Big)
  {\d x\over 1+x^2}
  \hfill\cr\noalign{\vskip .1in}\vbox{\kern 0.4in}\hfill
  \smash{{}
  \cases{%
    \leq (1+\eta )\displaystyle\int_0^{\epsilon t}
       {\displaystyle x \d x\over \displaystyle t(1+x^2)}
      = {\displaystyle 1+\eta\over\displaystyle 2t} 
      \log (1+\epsilon^2t^2 )  &\cr
    \noalign{\vskip .05in}
    \geq (1-\eta) \displaystyle\int_0^{\epsilon t} 
       {\displaystyle x \d x\over \displaystyle t(1+x^2)}
      = {\displaystyle 1-\eta\over\displaystyle 2t} 
      \log (1+\epsilon^2t^2 )  \, . &\cr
        }%
  }\qquad\cr
}$$
\vbox{\kern .5in}
Moreover,
$$
  0\leq \int_{\epsilon t}^\infty 
    \Big( {\pi\over 2}-\arctan {t\over x}\Big) 
    {\d x\over 1+x^2} 
  \leq {\pi\over 2} \int_{\epsilon t}^\infty {\d x\over x^2}
  = {\pi\over 2\epsilon t}  \, .
$$
Therefore, using $\log (1+\epsilon^2t^2 )\sim 2\log t$ as
$t$ tends to infinity, we have
$$
  P\{\, XY\geq t\,\} \sim {2\over \pi^2} {\log t\over t}
  \qquad\hbox{as } \ttoi \, .
$$
Next, for $\alpha\in (0,1)$, the same argument shows that
$$\eqalign{
  P\{\, XY\geq t\, ; X\geq t^\alpha \,\} 
  &= {1\over \pi^2} \int_{t^\alpha}^\infty \int_{t/x}^\infty
    {\d y\over 1+y^2} {\d x\over 1+x^2} \cr
  &= {1\over \pi^2} \int_{t^\alpha}^\infty 
    \Big( {\pi\over 2}-\arctan {t\over x}\Big) 
    {\d x\over 1+x^2} \cr
  &\sim {1\over \pi^2} \int_{t^\alpha}^{\epsilon t}
    {x\over t} {\d x\over 1+x^2} \cr
  &= {1\over t\pi^2} 
    \big( \log \epsilon t - \log t^\alpha +o(1)\big) \cr
  & \sim {\log t\over t\pi^2} (1-\alpha ) \qquad
    \hbox{ as } t\to\infty . \cr
  }       
$$
Consequently, we obtain
$$
  P\{\, X\geq t^\alpha \, | \, XY\geq t\,\} 
  = { P\{\, X\geq t^\alpha \, ; XY\geq t\,\}  \over
      P\{\, XY\geq t\,\}  }
  \sim {1-\alpha \over 2} \, .
$$
It follows that the distribution of $(X,Y)/t$ given ${XY\geq t}$
converges weakly* to $\delta_{(0,0)}$, showing a degeneracy like
in example 1, but of a different nature. A linear
normalization is not suitable, so we must use a logarithmic
scale. Writing $X=\epsilon_1e^U$ and $Y=\epsilon_2 e^V$ with
$\epsilon_1$ and $\epsilon_2$ taking values $+1$ or $-1$ with
probability $1/2$, we proved that the distribution of
$(U,V)/\log t$ given $XY\geq t$ converges weakly* to a uniform
distribution over $\{\, (s,1-s) : s\in [\,0,1\,]\,\}$, which 
is very similar to example 2.

\bigskip

Another interesting feature about examples 1--4 is that they use
some rather specific arguments in each case, since
everything could be calculated quite explicitly. When looking at
more complicated distributions, or at more complicated
conditioning, or in higher dimensions, we cannot 
rely so much on
elementary intuition. One of the main goals of the present work is
to give a systematic procedure for doing the calculation. The key
point is that we will be able to transform the problem into one
in asymptotic differential geometry, i.e., a
conjunction of asymptotic analysis and differential geometry. In
practice, the approximation of the conditional distribution will
boil down to a good understanding of the contact between $\partial A_t$ and
the level sets of the density.  Having to deal with purely geometrical
quantities will be helpful because differential
geometric arguments will not depend on the dimension, 
and also since differential
geometric quantities, such as curvatures, can be computed from
different formulas, using very different parameterizations.
Whichever parameterization is the most convenient will
be used. In some
sense, this is very much like changing variables in integrals, and
we hope to convince the reader that once it has been learned, it
considerably simplifies the analysis. The disadvantage is
that it requires investing some time in learning the basics of
differential geometry and asymptotic analysis --- but can we
really hope to solve every problem with elementary calculus?

At this stage, I urge the reader to look at the conditional
distribution of the parallelogram given its volume mentioned in
the introduction. In terms of linear algebra, the
question amounts to looking at the conditional distribution of a
random matrix with independent and identically distributed
coefficients, given that its determinant
is large. For $2\times 2$ matrices, we are already in $\RR^4$, 
since we have $4$ coefficients. The condition that the
determinant is larger than $t$ determines a subset whose boundary
is of dimension $3$, and it is impossible to visualize
it. Such a simple example already shows the need 
for systematic methods
relying as little as possible on intuition.

\bigskip

Let us now say a few words on how our results are connected to the
classical Laplace method. They can be viewed as a Laplace
method with an infinite dimensional parameter. 

In order to discuss this point, and because some understanding
of Laplace's method may be helpful in reading the next pages,
let us state and prove the following elementary result. Writing
$I=-\log f$, consider an integral of the form 
$\int_A e^{-\lambda I(x)} \d x$.

\bigskip

\state{1.1. THEOREM.} {\it%
  (Laplace's approximation) \qquad
  Let $A$ be a compact set in $\RR^d$ and $I$ be a twice 
  differentiable function on $\RR^d$. Assume that $I$ has a
  minimum on $A$ at a unique interior point $a_*$, and that
  $ \D^2I(a_*)$ is definite. Then,
  $$
    \int_A e^{-\lambda I(x)} \d x \sim
    { (2\pi )^{d/2}\over \big(\det\, \D^2I(a_*)\big)^{1/2} }
    \lambda^{d/2} e^{-\lambda I(a_*)}
    \qquad \hbox{ as } \lambda\to\infty \, .
  $$
  }

\bigskip

\stateit{Proof.} Let $\epsilon$ be a number between $0$ and $1/2$.
Since $a_*$ is in
the interior of $A$ and $I$ is twice differentiable, we can find
an open neighborhood $U$ of $a_*$ on which
$$\displaylines{\qquad
    {1-\epsilon\over 2}\< \D^2 I(a_*)(x-a_*)\, , \, x-a_*\>
    \leq I(x)-I(a_*)
  \cut
    \leq {1+\epsilon\over 2} \< \D^2I(a)(x-a_*),x-a_*\>\, .
  \qquad\cr}
$$
Moreover, there exists a positive $\eta$ such that $I(x)\geq
I(a_*) + \eta$ on $A\cap U^{\rm c}$. Consequently,
$$\displaylines{
  \int_A e^{-\lambda I(x)} \d x
  \leq e^{-\lambda I(a_*)}\int_{A\cap U}\exp \Big( 
  {1-\epsilon\over 2} \< \D^2I(a_*)(x-a_*),x-a_*\>\Big) \d x
  \cut
  + e^{-\lambda (I(a_*)+\eta )} |A\cap U^{\rm c}| \, .
 \quad\cr}
$$
The change of variable $h=\sqrt{\lambda (1-\epsilon)}(x-a_*)$
yields to the upper bound
$$\displaylines{
  e^{-\lambda I(a_*)}
  \int_{\sqrt{\lambda (1-\epsilon)}(A\cap U-a_*)}
  \exp\Big( -{1\over 2} \< \D^2I(a_*)h,h\>\Big) \d h \,
  \lambda^{-d/2} (1-\epsilon )^{-d/2}
  \cut
  + e^{-\lambda (I(a_*)+\eta)} |A\cap U^{\rm c}| \, .
  \quad\cr}
$$
Since $a_*$ is an interior point of $A$, the set $\sqrt{\lambda
(1-\epsilon)}(A\cap U-a_*)$ expands to fill $\RR^d$ as
$\lambda$ tends to infinity. Therefore,
$$
  \int_A e^{-\lambda I(x)} \d x \leq
  {\lambda^{-d/2}\over \big(\det\, \D^2I(a_*)\big)^{1/2}}
  e^{-\lambda I(a_*)} (1-2\epsilon )^{-d/2} 
$$
for $\lambda$ large enough. A similar argument leads to the
lower bound where the term $(1-2\epsilon )^{-d/2}$ is replaced by
$(1+2\epsilon )^{-d/2}$. Since $\epsilon$ is arbitrary, the
result follows.\hfill$\qed$

\bigskip

How is our integral $\int_A e^{-I(x)} \d x$ related to Laplace's
method? If $A=tA_1$ for some fixed set $A_1$, and if $I$ is a
homogeneous function of degree $\alpha$, we see that 
$$
  \int_{tA_1} e^{-I(x)} \d x
  = t^d \int_{A_1} e^{-I(ty)} \d y 
  = t^d \int_{A_1} e^{-t^\alpha I(y)} \d y \, .
  \eqno{(1.3)}
$$
On one hand Laplace's method is directly applicable provided $I$ has a
unique minimum on $A_1$. On the other hand, 
the integral on the left
hand side of (1.3) is like (1.1). Having
a good estimate for (1.1) will yield an estimate of (1.3).

When studying (1.1), one can
try to argue as in Laplace's method. Writing $I(A)$ for
the infimum of $I$ over $A$, we obtain
$$
  \int_A e^{-I(x)} \d x
  = e^{-I(A)} \int_A e^{-( I(x)-I(a_*)) } \d x \, .
$$
We can hope to have a quadratic
approximation for $I(x)-I(a_*)$. 
But if we consider an arbitrary set far away 
from the origin, there is no reason for $I$ to have a minimum at a 
unique point
on $A$, and there is even less reason to have this infimum in
the interior of $A$. In general attains its minimum on a set 
$\DA$ depending on $A$, which can be very messy. 
Roughly, we will limit ourselves to situations
where $\DA$ is a smooth $k$-dimensional manifold, which still
allows some room for wiggling curves and other rather nasty
behaviors. In essence, thinking of $A$ as
parametrized by its boundary $\partial A$, our method consists
in applying Laplace's approximation at every point of $\DA$
along fibers orthogonal to $\DA$, and integrating these
approximations over $\partial A$, keeping only the leading terms. The
difficulty is to obtain a good change of variable formula in
order to extract the leading terms, to keep a good control
over the Jacobian of the transformation, and to control all the
error terms. In doing that, the reader will see that it amounts to a
Laplace method where the parameter $\lambda$ is now the infinite
dimensional quantity $\partial A$.

Since we allow so much freedom on how the set $A$ can
be, we will see that not only should we look at points were $I$ is
minimized, but also at points $x$ in $A$ such that $I(x)-I(A)$ 
stays bounded --- this is of course not a rigorous statement --- or
even is unbounded but not too large compared
to some function of $I(A)$.

Approximating an integral of the 
form $\int_A e^{-I(x)} \d x$
for arbitrary sets $A$ far away from the origin and arbitrary
functions $I$ seems a very difficult task. For the 
applications that we
have in mind, we will only concentrate on sets with smooth
boundary. We will also require that $I$ be convex. This last
assumption may look very restrictive. However, by not putting too
many restrictions besides smoothness on $\partial A$,
one can often make a change of variable in order to get back to the
case where $I$ is convex. This is actually how we obtained the
estimates of the integrals over $s_\alpha$ given at the beginning of
this chapter. As the reader can see, $-\log s_\alpha$ is all but a
convex function. This idea of changing variables will be
systematically illustrated in our examples.

\bigskip

We now outline the content of these notes.

Chapters 2--5 are devoted to the proof of an asymptotic equivalent
for integrals of the form $\int_A e^{-I(x)} \d x$ for smooth sets $A$
far away from the origin and convex functions $I$ --- plus some other
technical restrictions. The goal and culminating point is to prove
Theorem 5.1. The main tools come from differential
geometry and related formulas in integration. A reader with no
interest in the theoretical details can go directly to Theorem
5.1, but will need to read a few definitions in chapters 2--4 in
order to fully understand its statement.

In chapter 6, we consider the special case where $A$ is the
translate of a fixed set.

A second basic situation is studied in chapter 7, where $A=tA_1$ is
obtained by scaling a fixed set $A_1$ which does not contain the
origin, $t$ tends to infinity and $I$ is homogeneous. In this 
case, we will overlap and somewhat extend the classical 
Laplace method.

In chapter 8, we study the tail probability of quadratic forms of
random vectors. If $X$ is a random vector in $\RR^d$ and $C$ is a
$d\times d$ matrix, we seek an approximation of $P\{\, \langle
CX,X\rangle\geq t\,\}$ for large $t$. We focus on the case where $X$
has a symmetric Weibull type distribution or a Student type
distribution. These two cases yield different behavior which we
believe to be representative of the so-called sub-Gaussian and heavy
tails distributions. They also illustrate the use of our main
approximation result.

The next example, about random linear forms, developed in chapter 9, requires
more geometric analysis. The problem is as follows: if $X$ is a
random vector on $\RR^d$ and $M$ is a subset of $\RR^d$, how does
$P\{\, \sup_{p\in M} \langle X,p\rangle\geq t\,\}$ decay as
$t$ tends to infinity? Again, our main theorem provides an answer
and we will also
focus on symmetric Weibull- and Student-like vectors. These two
distributions capture the main features of the applications of our
theoretical result to this situation.

The last example, treated in chapter 10, deals with random matrices.
Specifically, if $X$ is an $n\times n$ matrix with independent and
identically distributed coefficients, we seek for an approximation
of $P\{\, \det\, X\geq t\,\}$ for large $t$. Again, we will deal with
Weibull- and Student-like distributions. 
In the last subsection, we also approximate $P\{\, \| X\|\geq
t\,\}$ for large $t$, the tail distribution of the norm of the
random matrix. This last example yields some interesting geometry
on sets of matrices, and turns to be an application of the
results on random linear forms obtained in chapter 9.

The last chapters address some more applied issues, ranging from
applications in statistics to numerical computation. 

Chapter 11 presents some applications to statistical analysis of
time series. We will mainly obtain tail distribution of 
the empirical covariances of autoregressive models. 
This turns out to be an application of the results of chapter 8. We
will go as far as obtaining numbers useful for some real
applications, doing numerical work.

Chapter 12 deals with the distribution of the suprema of some
stochastic processes. It contains some examples of pedagogical
interest. It concludes with some calculations related to
the supremum of the Brownian bridge and the supremum of an
amusing process defined on the boundary of convex sets.

There are two appendices. One deals with classical estimates on
tail of Gaussian and Student distributions. In the other one, we
prove a technical estimate on the exponential map.

Every chapter ends with some notes, giving information and/or open
problems on related material.

\bigskip

When doing rather explicit calculations, we will focus on two
specific families of distributions. Recall that the Weibull
distribution on the positive half line has cumulative
distribution function $1-e^{-x^\alpha}$, $x\geq 0$. Its density
is $\alpha x^{\alpha-1}e^{-x^\alpha}\indic_{[0,\infty )}(x)$. We
will consider a variant, namely the density
$$
  w_\alpha(x)=K_{w,\alpha}e^{-|x|^\alpha/\alpha} \, ,\qquad
  x\in\RR \, ,
$$
where setting $K_{w,\alpha}=\alpha^{1-(1/\alpha)}/\big( 2\Gamma
(1/\alpha)\big)$ ensures that $w_\alpha$ integrates to $1$ on
the real line. We call this distribution symmetric Weibull-like.
This is the first specific family of distributions that we will
use. 

To introduce the second family, recall that the Student
distribution has density proportional to
$(1+\alpha^{-1}y^2)^{-(\alpha+1)/2}$. Its cumulative
distribution function $\underline S_\alpha$ has tail given by
the asymptotic equivalent --- see Appendix 1 ---
$$
  \underline S_\alpha (-x) \sim 1-\underline S_\alpha (x)
  \sim {K_{s,\alpha}\alpha^{(\alpha-1)/2}\over x^\alpha}
  \qquad \hbox{ as } x\to\infty \, ,
$$
where $K_{s,\alpha}=\int_{-\infty}^\infty
(1+\alpha^{-1}y^2)^{-(\alpha+1)/2}\, \d y$ is the normalizing
constant of the density. Accordingly, we say that a density
$s_\alpha$ is Student like if the corresponding distribution
function $S_\alpha (x)=\int_{-\infty}^x s_\alpha (y) \d y$
satisfies
$$
  S_\alpha (-x)\sim 1-S_\alpha (x) 
  \sim {K_{s,\alpha}\alpha^{(\alpha-1)/2}\over x^\alpha}
  \qquad\hbox{ as } x\to\infty \, ,
$$
for some constant $K_{s,\alpha}$.

Why are we interested in these two specific families? When
$\alpha=2$ the Weibull-like distribution is the standard
Gaussian one. Embedding the normal distribution this way will
allow us to see how specific the normal is. When looking at
product densities, it is only when $\alpha=2$ that the Weibull-like 
distribution is invariant under orthogonal transforms; this
will create discontinuities at $\alpha=2$ in some asymptotic
approximations. 

The Student-like distributions are of interest because they are
rather representative of the so-called heavy tail distributions.
In particular, they include the symmetric stable distributions
with index $\alpha$ in $(0,2)$. But more generally, recall that a
cumulative distribution function $F$ is infinitely divisible if
its characteristic function has the form
$$
  \exp \int_{-\infty}^{+\infty} 
  {e^{i\zeta x}-1-i\zeta\sin x\over x^2} \d M(x)
$$
for some finite measure $M$ --- see, e.g., Feller (1970, \S
XVII). For $x$ positive, define
$$
  M^+(x)=\int_x^{+\infty} {\d M(y)\over y^2}
  \qquad\hbox{ and }\qquad
  M^-(x)=\int_{-\infty}^{-x} {\d M(y)\over y^2} \, .
$$

Recall that a function $f$ on $\RR$ is said to be regularly varying
with index $\rho$ at infinity if, for any positive $\lambda$,
$$
  \lim_{x\to\infty} {f(\lambda x)\over f(x)} = \lambda^\rho
  \, .
$$
It can be proved --- see, e.g., Feller (1970, \S XVII.4) ---
that whenever $M^+$ (resp. $M^-$) is regularly varying with
negative index, then $1-F(x)\sim M^+(x)$ as $x$ tends to
infinity (resp. $F(x)\sim M^-(x)$ as $x$ tends to minus
infinity). Thus, if
$$
  M(-\infty ,-x)\sim M(x,\infty) \sim {c\over x^{\alpha-2}}
  \, ,
$$
with $\alpha>2$, then $F$ is Student-like. More precisely,
$$
  1-F(x)\sim {c(\alpha-2)\over \alpha x^\alpha}
  \qquad\hbox{ as } x\to\infty \, ;
$$
one can take the constant $K_{s,\alpha}$ to be
$c(\alpha-2)\alpha^{-(\alpha+1)/2}$. 

Other distributions are Student-like. For instance, one may start
with a random variable $X$ having a Pareto distribution on
$(a,\infty)$,
$$
  P\{\, X\geq x\,\} 
  = {1\over (x-a+1)^{\alpha}} \,\indic_{[a,\infty)}(x) \, .
$$
We then symmetrize it. Take $\epsilon$ to be a random variable
independent of $X$, with
$$
  P\{\, \epsilon =-1\,\} = P\{\, \epsilon = +1\,\}
  = 1/2 \, .
$$
Then $\epsilon X$ has a symmetric distribution. Its tails are
given by
$$
  P\{\, \epsilon X\geq x\,\} 
  = P\{\, \epsilon X \leq -x \,\}
  = {1\over 2(x-a+1)^\alpha}
$$
for $x$ large enough.
Thus, $\epsilon X$ has a Student-like distribution.

\bigskip

\centerline{\fsection Notes}
\section{1}{Notes}

\bigskip

\noindent There is an extensive literature on calculating integrals.
Concerning the classical numerical analysis, quadratures, 
Monte-Carlo methods, etc., chapter 4 of the {\sl Numerical Recipes} by
Press, Teukolsky, Vetteling and Flannery (1986) is an excellent
starting point. It states clearly the problems, some solutions,
their advantages and drawbacks, and contains useful references.

The statistical literature is full of more or less ad hoc methods
to perform some specific integrations, sometimes related to those
we are interested in here. In the last 20 years or so, Markov
chain Monte-Carlo and important sampling methods have been
blooming, bringing a large body of papers. Unfortunately, I do not
know a reference explaining things simply. Perhaps this area is
just too active and has not yet matured to a stage where classical
textbooks become available.

Concerning the approximation of univariate integrals, Laplace's
method, asymptotic expansions and much more, I find Olver (1974) a
great book. The classical references may be Murray (1984); De
Bruijn (1958) and Bleistein and Handelsman (1975) have been
republished by Dover. Both are inexpensive --- as any book should
be --- and are worth owning. Not so well known is Combet (1982). 
Combet has proofs of
the existence of asymptotic expansions derived from very general
principles, although these expansions are not too explicit. I
tend to believe that work in this direction would bring some
practical results, and this is the first --- loose --- open problem in
these notes. 

Regularly varying functions are described beautifully in 
Bingham, Goldie and Teugels (1987).

\vfil\eject

~\vskip 28pt
\chapter{2}{The logarithmic estimate}
\line{\fchapter\hfill 2.\quad The logarithmic}
\vskip .1in
\line{\fchapter\hfill estimate}

\vskip 54pt

\noindent Throughout this book, we consider a nonnegative function $I$ 
defined on $\RR^d$,\ntn{$I(\cdot )$} convex, such that 
$$
  \lim_{|x|\to\infty} I(x) = +\infty \, .
$$
This assumption holds whenever we use a function $I$, and we 
will not repeat it every time.

For any subset $A$ of $\RR^d$, we consider the infimum of $I$
over $A$,\ntn{$I(A)$}
$$
  I(A)=\inf \{\, I(x) : x\in A\,\} \, .
$$

When studying an integral of the form $\int_A e^{-I(x)}\d x$ with
$A$ far away from the origin, a first natural
question is to investigate if it is close to $e^{-I(A)}$, at least in
logarithmically. If this is the case, we have an analogue of the
logarithmic estimate for exponential integrals
$$
 \log\int_A\exp \big( -\lambda I(x)\big) \d x \sim \lambda I(A)
  \qquad\hbox{as } \lambda\to\infty 
$$
--- $A$ is fixed here! --- which holds, e.g., for sets $A$ with 
smooth boundary.

Our first result in this flavor will be a slightly sharper upper estimate
than a purely logarithmic one, but for some very specific sets $A$, namely,
complements of level sets of $I$. This estimate will be instrumental in the
following chapters.

Let us consider the level sets of $I$,\ntn{$\Gamma_c=I^{-1}\big(
-\infty,c])$}
$$
  \Gamma_c=\big\{\, x : I(x)\leq c\, \} \, , \qquad c\geq 0 \, .
$$
The complement $\Gamma_c^{\rm c}$ of $\Gamma_c$ in $\RR^d$
is the set of all points for which $I$ is strictly larger than $c$. 
Define the function\ntn{$L(c)$}
$$
  L(c)=\int_{\Gamma_c^{\rm c}} \eIx \d x \, , \qquad c\geq 0 \, .
$$

\bigskip

\state{2.1. PROPOSITION.}\quad{\it%
  There exists a constant $C_0$ such that for any positive $c$,
  $$
    L(c)\leq C_0 e^{-c} (1+c)^d \, .
  $$}

\bigskip

\stateit{Proof.} Use Fubini's theorem and the change of
variable $u=v+c$ to obtain
$$\eqalign{
  L(c) & = \int\int \indic_{[c,u]}\big( I(x)\big) \, \d x\,  
         e^{-u} \, \d u 
         = e^{-c} \int_{v\geq 0} e^{-v} \, |\Gamma_{c+v}\setminus \Gamma_c |
         \, \d v \cr
  \noalign{\vskip .05in}
       & \leq e^{-c} \int_{v\geq 0} e^{-v} \, | \Gamma_{c+v} | 
         \, \d v \, . \cr
  }$$
Since the function $I$ is convex, nonnegative and tends 
to infinity with its argument, its graph 
$$
  \big\{\, \big(x,I(x)\big) : x\in \RR^d\, \big\}
  \subset \RR^{d+1}
$$
is included in a cone with vertex the point $(0,\ldots  , 0, -1)$ and
passing trough a ball centered at $0$ in $\RR^d\times \{\, 0\,\}$.

\bigskip

\setbox1=\hbox to 2in{\vbox to 2in{\ninepointsss
  \hsize=2in
  \kern 10.9in
  \includegraphics{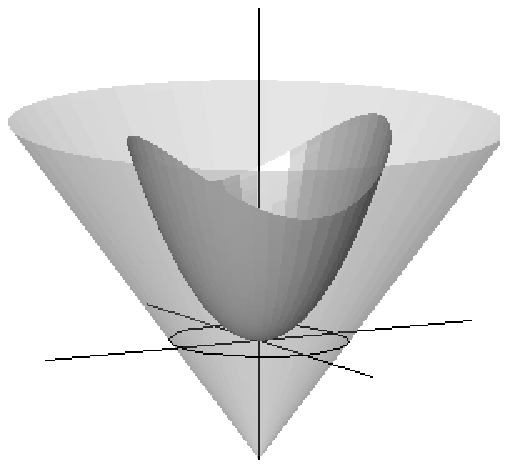}
  \vfill
  \vss
}\hfill}    
\setbox2=\hbox to 2in{\vbox to 2in{\ninepointsss
  \hsize=2in
  \kern 10.9in
  \includegraphics{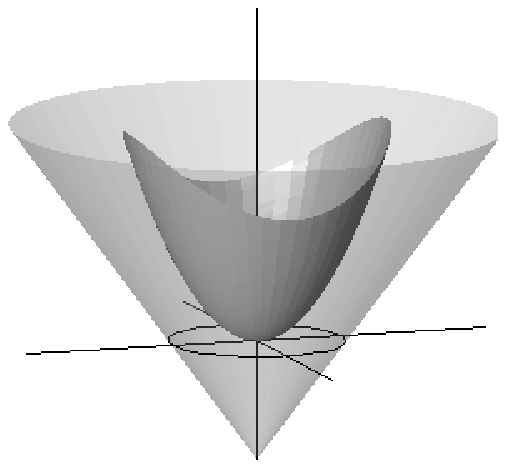}
  \vfill\vss
}\hfill}  

\bigskip

\centerline{\box1\hfill\box2}

\bigskip

Thus, $|\Gamma_t|$ is less than the $d$-dimensional measure 
of the slice of the cone at height $t$, i.e.,
$$
  |\Gamma_t| \leq C_1 (1+t)^d
$$
for some positive constant $C_1$. Thus,
$$\eqalign{
  L(c) & \leq C_1 e^{-c} \int_{v\geq 0} e^{-v} (c+v+1)^d \d v \cr
       & \leq 3^d C_1 e^{-c} \int_{v\geq 0} e^{-v} (c^d+v^d+1) \d v \cr
  }$$
and the result follows.  \hfill$\qed$%

\bigskip

To handle more general sets, it is convenient to introduce the
following notation.

\bigskip

\state{NOTATION.}\quad{\it%
   Let
   $F$, $G$ be two functions defined on the Borel $\sigma$-field
   of $\RR^d$. We write

  \noindent (i) $\lim_{A\to\infty} F(A)= 0$ if and only if
  for any positive $\epsilon$, there exists $c$ such that for all
  Borel set $A$ of $\RR^d$, the inequality $I(A)\geq c$ implies $|F(A)|\leq
  \epsilon$;

  \noindent (ii) $\lim_{A\to\infty} F(A)=0\Rightarrow
  \lim_{A\to\infty} G(A)=0$ if and only if the following holds:

  $\forall\epsilon , \, \exists\delta ,\, \exists c>0 , \,
  (I(A)\geq c\hbox{ and } |F(A)|\leq \delta ) \Rightarrow
  |G(A)|\leq\epsilon$.
  }

\bigskip

Another way to phrase condition (ii) is to say that whenever $I(A)$ is
large enough and $F(A)$ is small enough, then $G(A)$ is small.

Notice that this notion of limit depends on the function $I$. But if
we restrict $I$ to be convex, defined on $\RR^d$ and blowing up at
infinity, then it does not depend on which such specific function we
choose.

The advantage of this notation is that it allows to express approximation
properties related to sets moving away from the origin, but under some
analytical constraints. We will mainly use (ii). We can similarly define
$\liminf_{A\to\infty}F(A)$ and $\limsup_{A\to\infty}F(A)$.

A first example of the use of this notation is to express a condition
which ensures that we have a ``nice'' logarithmic estimate. It
asserts that $\log\int_A \eIx \d x$ is of order $-I(A)$ provided
$A$ is not like a very thin layer attached to $\partial
\Gamma_{I(A)}$.

\bigskip

\state{2.2. PROPOSITION.}\quad{\it%
  The following are equivalent:

  \vskip .05in

  \noindent (i) $\lim_{A\to\infty}I(A)^{-1}
  \log{\displaystyle \int_A}\eIx \d x = -1$,

  \vskip .05in

  \noindent (ii) $\lim_{\epsilon\to 0} \liminf_{A\to\infty}
  I(A)^{-1} \log |A\cap \Gamma_{(1+\epsilon)I(A)}| = 0$.
  }

\bigskip

%\setbox1=\hbox to 1.7in{\vbox to 1.2in{\sevencmss
%  \hsize=1.7in
%  \kern 1in
%  \special{psfile=Graph/ch2_fig1.ps}
%  \anote{.6}{-.1}{$\scriptstyle \Gamma_{\hbox{\fivecmss I(A)}}$}
%  \anote{.9}{-.1}{$\scriptstyle \Gamma_{\hbox{%
%     \fivecmss (1+$\scriptscriptstyle\epsilon$)I(A)}}$}
%  \anote{1.4}{.15}{$\scriptstyle\partial$\sevencmss A}
%  \anote{1.2}{.7}{A}
%  \vfill
%}\hfill}     
\setbox1=\hbox to 1.7in{\vbox to 1.2in{\ninepointsss
  \hsize=1.7in
  \kern 1in
  \includegraphics{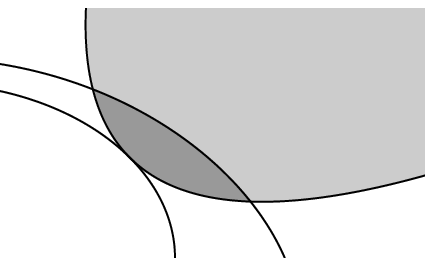}
  \anote{.6}{-.1}{$\Gamma_{I(A)}$}
  \anote{.9}{-.1}{$\Gamma_{(1+\epsilonsmall)I(A)}$}
  \anote{1.4}{.15}{$\partial A$}
  \anote{1.2}{.7}{$A$}
  \vfill
}\hfill}             

\bigskip

\centerline{\box1}

\bigskip

\state{REMARK.} For any set $A$,
$$
  \int_A \eIx \d x \leq \int_{\Gamma_{I(A)}^c} \eIx \d x 
  = L\big( I(A)\big) \, .
$$
We then infer from Proposition 2.1 that
$$
  \limsup_{A\to\infty} I(A)^{-1}\log\int_A \eIx \d x \leq -1 \, .
$$
Thus, statement (i) in Proposition 2.2 is really about the limit inferior 
of the integral as $I(A)$ tends to infinity.

Notice also that the the limit as $\epsilon$ tends to $0$ in (ii)
exists since $\epsilon\mapsto \Gamma_{(1+\epsilon )I(A)}$ is
monotone. The limit is at most $0$ since 
$|A\cap \Gamma_{(1+\epsilon )I(A)}|\leq |\Gamma_{(1+\epsilon )I(A)}|$
is bounded by a polynomial in $I(A)$. Therefore, statement (ii) is
really about the limit inferior as $\epsilon$ tends to $0$.

\bigskip

\stateit{Proof of Proposition 2.2.}
Assume that (ii) holds. For any positive number $\epsilon$,
$$\eqalign{
  I(A)^{-1} \log\int_A \eIx \d x
  & \geq I(A)^{-1} \log\int_{A\cap \Gamma_{(1+\epsilon )I(A)}}\eIx
    \d x \cr
  \noalign{\vskip .2cm}
  & \geq I(A)^{-1} \log\big( e^{-(1+\epsilon )I(A)}
    |A\cap\Gamma_{(1+\epsilon )I(A)}|\big)  \cr
  \noalign{\vskip .2cm}
  &= -1-\epsilon + I(A)^{-1} \log |A\cap\Gamma_{(1+\epsilon )I(A)}| \, .
   \cr
  }$$
Given the above remark, (i) holds.

To prove that (ii) is necessary for (i), let us argue by
contradiction. If (ii) does not hold, then there exists a positive
$\beta$ and a sequence of positive numbers $\epsilon_k$ converging 
to $0$ as $k$ tends
to infinity, such that
$$
  \liminf_{A\to\infty} I(A)^{-1} 
  \log  |A\cap \Gamma_{(1+\epsilon_k)I(A)}| < -\beta \, 
$$
for any $k$ large enough. For $k$ large enough, $\epsilon_k <\beta$.
Thus, using Proposition 2.1,
$$\eqalign{
  \int_A \eIx \,\d x
  & \leq \int_{A\cap \Gamma_{(1+\epsilon_k)I(A)}}\eIx \,\d x
    + \int_{\Gamma_{(1+\epsilon_k)I(A)}^{\rm c}} \eIx \, \d x
    \cr
  & \leq e^{-I(A)}|A\cap\Gamma_{(1+\epsilon_k)I(A)}|
    + L\big( (1+\epsilon_k)I(A)\big) \cr
  & \leq e^{-I(A)(1+\beta)} 
    + C_0 e^{-I(A)(1+\epsilon_k)} \big( 1+(1+\epsilon_k)I(A)\big)^d
    \cr
  & \leq e^{-I(A)(1+\epsilon_k + o(1))} \qquad\hbox{ as } A\to\infty
    \, . \cr
  }
$$
Therefore, 
$$
  \limsup_{A\to\infty} I(A)^{-1}\log \int_A \eIx \, \d x
  \leq -(1+\epsilon_k) < -1
$$
and (i) does not hold.\hfill$\qed$

\bigskip

\centerline{\fsection Notes}
\section{2}{Notes}

\bigskip

\noindent The idea of the proof of Proposition 2.1, writing $e^{-I}$ as 
the integral of $e^{-u}$ over $[\, I,\infty )$ and then using
Fubini's theorem is all but new; at most it has not been used
enough. Lieb and Loss (1997) call a much more general fact the
``layer cake representation''. I have the recollection of
hearing talks using similar tricks and referring to the coarea
formula. See Federer (1969) to know all about it, or Morgan
(1988) to simply know what it is. Proposition 2.2 grew from
Broniatowski and Barbe (200?) which we wrote, I think, in 1996 or
1997. This second chapter has the flavor of large deviation
theory and the notation $I$ for the function in the exponential
is not a coincidence. More will be said in the notes of chapter
5 and in chapter 7. There are by now a few books on large
deviations. Dembo and Zeitouni (1993) and Dupuis and Ellis
(1997) are good starting points to the huge literature. From a
different perspective, and restricted essentially to the
univariate case, Jensen (1995) may be closer to what we are
looking for here.

Introducing the set $\Lambda_c$ suggests that the variations of $-\log
f$ are important. It has to do with the following essential
remark. The negative exponential function is the only one --- up to an
asymptotic equivalence --- for which integrating on an interval of
length of order $1$ produces a relative variation of order $1$ on the
integral. Mathematically, the fact is that
$$
  \int_{t+s}^\infty e^{-x}\,\d x = e^{-s}\int_t^\infty e^{-x}\, \d x \, .
$$
Hence, in approximating $\int_{t+s}^\infty e^{-x}\d x$ by
$\int_t^\infty e^{-x}\, \d x$, we are making a relative error equals
to $e^{-s}-1$. If $s$ is fixed, this relative error stays fixed, event
if $t$ moves. If one writes the analogue formula with a power
function, one obtains
$$
  \int_{t+s}^\infty {\d x\over \alpha x^{\alpha+1}}
  = {1\over (t+s)^\alpha}
  = \Big( {t\over t+s}\Big)^\alpha \int_t^\infty {\d x \over\alpha x^{\alpha+1}}
  \, .
$$
As $t$ tends to infinity, the ratio $\big(t/(t+s)\big)^\alpha-1$ tends
to $0$. Thus the relative variation of the integral tends to $0$ as
$t$ tends to infinity. Finally, in the subexponential case, for $\alpha>1$,
$$
  \int_{t+s}^\infty \alpha x^{\alpha -1} e^{-x^\alpha} \, \d x
  = e^{-(t+s)^\alpha+t^\alpha} \int_t^\infty \alpha x^{\alpha-1}e^{-x^\alpha}
    \, \d x \, .
$$
The relative variation is now driven by
$$
  \exp \big(-st^{\alpha-1}(1+o(1))\big) -1 \, .
$$
If $s$ is positive (resp. negative), it tends to $-1$ (resp. infinity) 
as $t$ tends to infinity.

In other words, the exponential scale that we use on the variations of
$f$ ensures that the variation of $f(x)$ in this scale are of the same
order of magnitude as the variations of $x$ in space.

\vfil\eject

~\vskip 28pt
\chapter{3}{The basic bounds}
{\fchapter\hfill 3.\quad The basic bounds}

\vskip 54pt

\noindent In this chapter, our goal is to obtain lower and upper bounds for 
the integral $\int_A\eIx \d x$. In order to get useful estimates, we
need to decompose $A$ into small pieces on which the integral can be
almost calculated with a closed formula. These pieces are given by the
geometry of the graph of the function $I$, and therefore we will first
devote some time in introducing a few useful quantities related to
this graph.

\bigskip

\noindent{\fsection 3.1. The normal flow and the normal foliation.}
%\chapternumber={3.1}
\section{3.1}{The normal flow and the normal foliation}
\bigskip

\noindent Recall that we assume

\vskip\abovedisplayskip

\hfil\vbox{\hsize 4truein \noindent$I$ is strictly convex on $\RR^d$, nonnegative,
  twice differentiable, with $\lim_{x\to\infty} I(x)=+\infty$.}\hfil

\vskip\belowdisplayskip

\noindent Up to adding a constant to $I$, which amounts to multiplying
the integral $\int_A \eIx \d x$ by a
constant, and up to translating $A$ by a fixed amount, we can assume that
$$
  I(0)=0 \, .
$$
Let $\D I$ denote the differential --- or gradient --- of $I$ and $\D ^2I$ 
its Hessian. The strict convexity assumption implies that $\D ^2I$ is a 
symmetric definite positive matrix.

Recall that we use
$$
  \Gamma_c
  =\big\{\, x\in\RR^d : I(x)\leq c\,\big\}
  = I^{-1}([\, 0,c\, ])
$$
to denote the level set of $I$. We define the level lines --- or level
hypersurfaces --- of $I$ as the
points on which $I$ is constant, that is\ntn{$\Lambda_c$}
$$
  \Lambda_c=\big\{\, x\in\RR^d : I(x)=c\,\big\}
  = I^{-1} (\{ \, c\,\}) \, .
$$
Since $I$ is a convex function, $\Gamma_c$ is a
convex set. Moreover, $I$ being defined on $\RR^d$, we write
$$
  \RR^d=\bigcup_{c\geq 0} \Lambda_c \, .
$$
Let $x$ be a nonzero vector in $\RR^d$ --- equivalently, $I(x)$ is
nonzero.  Set $c=I(x)$. The gradient $\D I(x)$ is an outward normal to
$\Lambda_c$ at $x$ and defines a vector field on $\RR^d\setminus \{\,
0\, \}$. For $t$ nonnegative, we write $\psixt$ the integral curve of
the vector field $\D I$, such that $\psi (x,0)=x$ and $I\big(
\psixt\big)= I(x)+t$. Sometime it will be convenient to use the
notation $\psi_t(x)$ or $\psi_x(t)$ instead of $\psixt$.\ntn{$\psi
(x,s)=\psi_s(x)$}

\bigskip

\state{DEFINITION.} {\it%
  The flow $\psi_t$ is called the normal flow at time $t$.
  }

\setbox1=\vbox to 2.5in{\ninepointsss
  \baselineskip .15in
  \hsize=2.1in
  \kern 2.25in
  \includegraphics{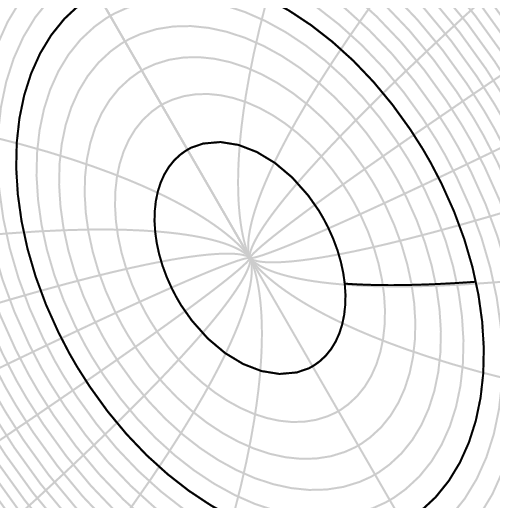}
  \anote{1.3}{0.81}{$x$}
  \anote{1.94}{0.81}{$\psi (x,t)$}
  \kern .1in
  \hsize=2in
  \centerline{$x$, the flow $\psi (x,s)$ for $s$ in $[\, 0\, ,t\, ]$,}
  \centerline{the level sets $\Lambda_{I(x)}$ and $\Lambda_{I(x)+t}$
              hilighted.}
  \vss
}  

\centerline{\box1}

\bigskip\bigskip

It is also convenient to introduce the outward normal unit 
vector to $\Lambda_c$ at $x$,\ntn{$N(x)$}
$$
  N(x)={\D I(x)\over |\D I(x)|} \, .
$$

\bigskip

\state{3.1.1. LEMMA.} {\it%
  The function $t\mapsto\psi (x,t)$ is a solution of the differential 
  equation
  $$
    {\d\over \d t}\psixt 
    ={\D I\big(\psixt\big) \over \big|\D I\big(\psixt\big) \big|^2} 
    \, ,
  $$
  with initial condition $\psi (x,0)=x$.
  }

\bigskip

\stateit{Proof.} By definition of the normal flow,
$$
  1={\d\over \d t} I\big(\psixt \big) 
  = \Big\langle \D I\big(\psixt\big) ,{\d\over \d t} \psixt\Big\rangle
  \, .
$$
The result follows since ${\displaystyle \d\over
\displaystyle \d t}\psi (x,t)$ must be collinear
to $\D I\big(\psi (x,t)\big)$.\hfill$\qed$

\bigskip

The next trivial fact is due to the convexity of $I$. It asserts that the norm of the
gradient $|\D I|$ can only increase along the normal flow.

\bigskip

\state{3.1.2. LEMMA.} {\it%
  The function $t\mapsto \big|\D I\big(\psixt\big)\big|$ is nondecreasing.
  }

\bigskip

\stateit{Proof.} Since $\D ^2I$ is positive definite,
$$
  {\d\over \d t} \big| \D I\big(\psixt\big)\big|^2
  = 2\big\< \D ^2I\big(\psixt\big) N\big(\psixt\big) ,
  N\big(\psixt\big) \big\>
  \geq 0 \, . \eqno{\qed}
$$

\bigskip

The monotonicity of $|\D I|$ along the normal flow yields a Lipschitz property of the map
$t\mapsto \psixt$.

\bigskip

\state{3.1.3. COROLLARY.} {\it%
  For any $x$ in $\RR^d\setminus\{\,0\, \}$ and nonnegative real number $t$,
  $$
    |\psixt-\psi (x,0)| 
    = |\psi (x,t)-x|
    \leq {t\over |\D I(x)|} \, .
  $$}

\bigskip

\stateit{Proof.} The result follows from the fundamental theorem of
calculus, Lemmas 3.1.1 and 3.1.2, since
$$\eqalignno{
  |\psi (x,t) -\psi (x,0)|
  \leq \int_0^t\Big| {\d\over \d s}\psi (x,s) \Big| \d s
  & = \int_0^t {\d s\over \big| \D I\big( \psi (x,s)\big)\big|} \cr
  & \leq {t\over | \D I(x)|} \, . &\qed\cr
  }$$

\bigskip

As a consequence, two points on the same integral curve and 
nearby level lines of $I$
cannot be too far apart in Euclidean distance. This is expressed in the following
result

\bigskip

\state{3.1.4. COROLLARY.} {\it%
  Let $M$ be a positive number. If $y$ is in the forward orbit of $x$ 
  under the flow $\psi_t$, i.e.,
  $$
    y\in\big\{\, \psi_t(x) : t\geq 0\,\big\} \, ;
  $$
  and if $I(y)-I(x)$ is in $[\, 0,M\, ]$, then
  $$
    |y-x|\leq {M\over | \D I(x)|} \, .
  $$}

\bigskip

\stateit{Proof.} Since $y=\psixt$ for some $t$,
$$
  M\geq I(y)-I(x) = I\big( \psixt\big) - I\big( \psi (x,0)\big) = t \, .
$$
Hence, Corollary 3.1.3 implies
$$
  |y-x|\leq {t\over |\D I(x)|}\leq {M\over |\D I(x)|}
  \eqno{\qed}
$$

\bigskip

So far, we have just collected some trivial informations on the normal
flow. Since $I$ is twice differentiable, the level lines $\Lambda_c$
are submanifolds --- hypersurfaces --- of $\RR^d$. The function
$\psi_t$ maps $\Lambda_c$ onto $\Lambda_{c+t}$.  For $p$ in
$\Lambda_c$, we can consider the tangent space $T_p\Lambda_c$ of
$\Lambda_c$. The differential of $\psi_t$ at some point $p$, which we
denote either by $\D \psi_t(p)$ or
$\psi_{t*}(p)$,\ntn{$\psi_{t*}(\cdot )$} maps $T_p\Lambda_{I(p)}$ to
$T_{\psi_t(p)}\Lambda_{I(p)+t}$. We will need some estimates on this
differential, and we first calculate its derivative with respect to
$t$.

\bigskip

\state{3.1.5. LEMMA.} {\it%
  The following holds
  $$
    {\d\over \d t} \psi_{t*} (p)
    = \Big[ \Id - 2 N\otimes N \big(\psi_t (p)\big) \Big]
    {\D ^2I\over |\D I|^2} \big(\psi_t(p)\big) \psi_{t*}(p) \, .
  $$}

\bigskip

\stateit{Proof.} Lemma 3.1.1 yields
$$
  \doverdt \psi_{t*} 
  = \Big( \doverdt \psi_t\Big)_*
  = \Big( {\D I\over |\D I|^2}\circ\psi_t\Big)_*
  = \Big( {\D I\over |\D I|^2}\Big)_* (\psi_t)\circ \psi_{t*}
$$
by the chain rule. Therefore,
$$
  \doverdt\psi_{t*} 
  = {\D ^2I\over |\D I|^2}(\psi_t)\circ \psi_{t*}
    -2{\D I\over |\D I|}\otimes{\D I\over |\D I|}{\D ^2I\over |\D I|^2}(\psi_t)\circ \psi_{t*} \, .
$$
This is the result since $N=\D I/|\D I|$. \hfill$\qed$

\bigskip

Since $I$ is convex, the level sets $\Gamma_c$ are 
convex, and one can see by a
drawing that $\psi_{t*}$ must be a map which expands 
distances --- see the picture before Lemma 3.1.1. 
This is expressed in 
the next result. Before stating
it, notice that for $p$ in $\Lambda_c$, the tangent space
$T_p\Lambda_c$ inherit of the Hilbert space structure of
$\RR^d$. Thus, it can be identified with its dual
$(T_p\Lambda_c)^*$. Since $\psi_{t*}$ maps $T_p\Lambda_c$ to
$T_{\psi_t(p)}\Lambda_{c+t}$, its transpose $\psi_{t*}^\T$ maps
$(T_{\psi_{t*}(p)}\Lambda_{c+t})^*\equiv
T_{\psi_t(p)}\Lambda_{c+t}$ to $T_p\Lambda_c$. Consequently,
$\psi_{t*}^\T\psi^{\phantom T}_{t*}$ is a linear map acting on $T_p\Lambda_c$.

\bigskip

\state{3.1.6. LEMMA.} {\it%
  For any nonnegative $t$, 
  $\det (\psi_{t*}^\T\psi_{t*}^{\phantom T})\geq 1$.
  }

\bigskip

\stateit{Proof.} Since $\psi_0$ is the identity function, we have
$\det (\psi_{0*})=\det (\Id)= 1$. It is then enough to prove that for every $h$
in $T_p\Lambda_c$, with $c=I(p)$, the map
$$
  t\mapsto\big| \psi_{t*}(p)h\big|^2 
  = \big\< h\, ,\psi_{t*}^\T(p)\psi^{\phantom T}_{t*}(p)h\big\>
$$
is nondecreasing. 

Writing $q=\psi_t(p)$, we infer from Lemma 3.1.5 that
$$
  \doverdt \big| \psi_{t*}(p)h\big|^2
  = 2\Big\< \big( \Id - 2N\otimes N(q)\big) {\D ^2I\over |\D I|}(q)
  \psi_{t*}(p)h \, , \psi_{t*}(p)h\Big\> \, .
$$
Since $\psi_{t*}h$ belongs to $T_q\Lambda_{c+t}$ it is 
orthogonal to $N(q)$. Therefore,
$$
  N(q)^\T\psi_{t*}(p)h 
  = \big\< N(q)\, ,\psi_{t*}(p)h\big\> = 0 \, .
$$
Since $\D^2 I$ is nonnegative, it follows that
$$
  \doverdt \big| \psi_{t*}(p)h\big|^2
  = 2\Big\< {\D ^2I\over |\D I|}(q)\psi_{t*}(p)h\, ,\psi_{t*}h\Big\>
  \geq 0 \, .
$$
Consequently, the given map is indeed nondecreasing.\hfill$\qed$

\bigskip

It follows from Lemma 3.1.6 that $\psi_{t*}$ is invertible.
Since $t\mapsto\psi_t$ is a semigroup, $\psi_{t*}$ cannot expand
too fast as $t$ increases.

\bigskip

\state{3.1.7. LEMMA.} {\it%
  For any $p$ in $\RR^d\setminus \{\, 0\,\}$ and any 
  nonnegative $s$, $t$,
  $$
    \| \psi_{t+s,*}(p)\|
    \leq \| \psi_{t*}(p)\| \exp\bigg( \int_0^s
      { \big\| \D ^2I\big(\psi_{t+u}(p)\big)\big\|\over 
        \big| \D I\big( \psi_{t+u}(p)\big) \big|^2} \d u\bigg) \, .
  $$}

\bigskip

\stateit{Proof.} Notice that $\Id-2N\otimes N$ is an inversion.
Its operator norm is $1$. Using Lemma 3.1.5, 
writing $q=\psi_t(p)$ and using the triangle
inequality for the increments of the function 
$t\mapsto \big\| \psi_{t*}(p)\big\|$, we obtain
$$
  \doverdt \big\| \psi_{t*}(p)\big\|
  \leq \Big\| \doverdt \psi_{t*}(p)\Big\|
  \leq {\|\D ^2I\|\over |\D I|^2}(q) \, \|\psi_{t*}(p)\| \, .
$$
Integrate this differential inequality between $t$ and $s$ to obtain the
result.\hfill$\qed$

\bigskip

In the same spirit as in the previous lemma, a little more work gives us a good
control on the growth of $\det (\psi_{t*}^\T
\psi^{\phantom T}_{t*})$, improving upon Lemma 3.1.6.

\bigskip

\state{3.1.8. LEMMA.} {\it%
  For any nonzero $p$ and any nonnegative $t$,
  $$\displaylines{\qquad
      \doverdt \log\det \big(\psi_{t*}^\T\psi_{t*}^{\phantom T}
      (p)\big)
    \cut
      ={2\over \big| \D I\big(\psi_t(p)\big)\big|^2}
      \trace \Big[ \Proj_{T_{\psi_t(p)}\Lambda_{I(p)+t}}
                   \D ^2I\big(\psi_t(p)\big)\Big|_{T_{\psi_t(p)}
                   \Lambda_{I(p)+t}}\Big]
      \, .
    \qquad\cr}
  $$
  In particular, for any nonnegative $s$, 
  $$\displaylines{\qquad
    \det \big( \psi_{t+s,*}^\T\psi_{t+s,*}^{\phantom T}(p)\big)
    \cut
    \leq \det\big( \psi_{t*}^\T\psi_{t*}^{\phantom T}(p)\big)
      \exp\bigg( 2\int_0^s
      { \trace \big[ \D ^2I\big( \psi_{t+u}(p)\big) \big] \over
        \big| \D I\big( \psi_{t+u}(p)\big) \big|^2 } \d u\bigg) \, ,
    \qquad}
  $$
  and the function $s\mapsto \det\big( \psi_{t+s,*}^\T
  \psi_{t+s,*}^{\phantom T}(p)\big)$ is
  nondecreasing.}

\bigskip

\stateit{Proof.} Let $p$ be a nonzero vector. Consider a local 
chart $p(\cdot ) :
U\subset \RR^{d-1}\mapsto\Lambda_{I(p)}$ 
around $p$, such that $p(0)=p$ and the vectors
$\partial_i={\partial\over\partial u_i}p(0)$ form an orthonormal basis of
$T_p\Lambda_{I(p)}$. Lemma 3.1.1 yields
$$\eqalign{
  {\partial\over \partial t}\big( \psi_{t*}(p)\partial_i\big)
  &= {\partial^2\over\partial u_i\partial t}
    \psi \big( p(u_1,\ldots ,u_{d-1}),t)\Big|_{u=0} \cr
  &= {\partial\over\partial u_i}
    \Big( {N\over |\D I|}\big(\psi (p,t)\big)\Big)\Big|_{u=0}  \, . \cr
  }
$$
Since $\psi_{t*}\partial_j$ belongs to $T_{\psi_t(p)}\Lambda_{I(p)+t}$,
it is orthogonal to $N\big(\psi_t(p)\big)$ for all $j$. Moreover,
$\d N$ is self-adjoint when acting on the tangent space of the
level curves $\Lambda_{I(p)+t}$ --- a known fact when studying 
the second fundamental form of immersions of hypersurfaces in $\RR^d$; 
see, e.g., Do Carmo (1992) --- we obtain
$$
  {\doverdt}\big\< \psi_{t*}(p)\partial_i
  \, ,\psi_{t*}\partial_j\big\>
  = {2\over \big| \D I\big(\psi_t(p)\big) \big| }
  \Big\< \d N\big( \psi_t(p)\big) \psi_{t*}(p)\partial_i \, ,
  \psi_{t*}(p)\partial_j\Big\> \, .
$$
Using that $\D \det (x)=\det (x)\trace (x^{-1}\cdot )$ --- a proof 
of this fact is in Lemma 10.0.1 --- and that $\psi_{t*}$ is invertible,
$$\displaylines{%
  \quad%
  \doverdt\det\big(
  \psi_{t*}^\T(p)\psi_{t*}^{\phantom T}(p)\big)
  \cut
  \eqalign{
  =\,&\det\big( \psi_{t*}^\T(p)\psi_{t*}^{\phantom T}(p)\big) 
  \trace\Big[ 
    \big(\psi_{t*}^\T(p)\psi_{t*}^{\phantom T}(p)\big)^{-1} 
    \doverdt \big(\psi_{t*}^\T(p)\psi_{t*}^{\phantom T}(p)\big)\Big] \cr
  =\,&2\, { \det\big( \psi_{t*}^\T(p)\psi_{t*}^{\phantom T}(p)\big) \over 
      \big| \D I\big( \psi_{t*}(p)\big)\big| }
    \, \trace\Big[ \big( \psi_{t*}^\T(p)\psi_{t*}^{\phantom T}(p)
    \big)^{-1} \psi_{t*}^\T(p)
    \d N\big(\psi_t^{\phantom T}(p)\big) \psi_{t*}^{\phantom T}(p)
    \Big] \cr
  =\,&2\, { \det \big( \psi_{t*}^\T(p)\psi_{t*}^{\phantom T}(p)\big)\over 
      \big| \D I\big( \psi_{t*}(p)\big)\big| }
    \, \trace\Big[ \Proj_{T_{\psi_t(p)}\Lambda_{I(p)+t}}
                   \d N\big( \psi_t(p)\big)\Big|_{T_{\psi_t(p)}
                   \Lambda_{I(p)+t}}\Big] \, . \cr
  }\cr}$$
This is the first assertion in the Lemma, since the 
restriction of $\d N$ to the tangent space coincides with
$\D ^2I/|\D I|$. Moreover, since $\D ^2I$ is positive, we have 
$$
  \trace\Big[ \Proj_{T_{\psi_t(p)}\Lambda_{I(p)+t}}
              \d N\big( \psi_t(p)\big) 
              \Big|_{T_{\psi_t(p)}\Lambda_{I(p)+t}}\Big]
  \leq {\trace\big[ \D ^2I\big( \psi_t(p)\big) \big] 
    \over \big| \D I\big(\psi_t(p)\big)\big| } \, ,
$$
and the second statement follows from integration.

The third statement follows from the first one. Indeed, $\D ^2I$ is 
nonnegative, and so
${\displaystyle \d\over\displaystyle \d t} \log\det 
\big(\psi_{t*}^\T\psi_{t*}^{\phantom T}(p)\big)\geq 0$.
\hfill$\qed$

\bigskip

We now have all the informations we need on the normal flow.  Its
integral curves define a foliation of $\RR^d\setminus \{\, 0\,\}$. We
will obtain an expression for $\int_A\eIx \d x$ by integrating first
on $\Lambda_{I(A)}$, and then along the leaves.  The virtue of the
normal flow is that computation of the Jacobian in the change of
variables is easy, because the normal vector field and the level sets
of the function $I$ are orthogonal at {\it every} point of
$\RR^d$. Before rewriting the integral, we need to parameterize the
boundary $\partial A$ in the flow coordinate system.  Specifically,
for any point $p$ in $\RR^d$, define\ntn{$\tau_A(p)$}
$$
  \tau_A (p)=\inf\big\{\, t\geq 0 : \psi_t(p)\in A\,\big\} \, ,
$$
with the convention that $\inf\emptyset=+\infty$. In the sequel, we will always read
$e^{-\infty}$ as $0$.

\finetune{\vfill\eject}

\setbox1=\vbox to 2.5in{\ninepointsss
  \baselineskip .15in
  \hsize=2.1in
  \kern 2.25in
  \includegraphics{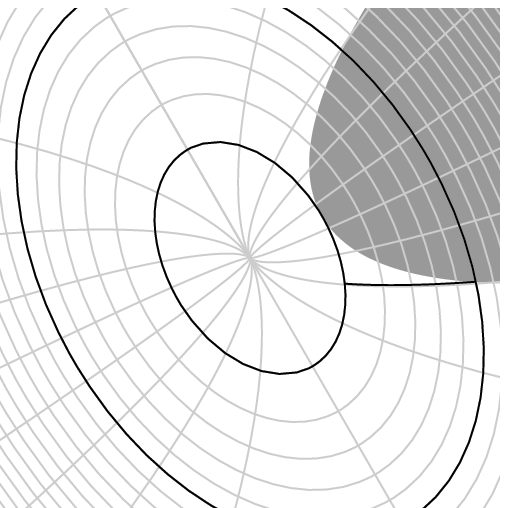}
  \anote{1.8}{1.8}{$A$}
  \anote{1.3}{0.81}{$p$}
  \anote{1.94}{0.81}{$\psi (x,\tau_A(p))$}
  \kern .1in
  \hsize=2.05in
  \centerline{$p$, the flow $\psi (p,s)$ for $s$ in 
              $[\, 0,\,\tau_A(p)\,]$,}
  \centerline{the level sets $\Lambda_{I(A)}$ and 
              $\Lambda_{I(A)+\tausmall_A(p)}$ hilighted.}
  \vfill\vss
}  

\centerline{\box1}

\bigskip

It is convenient to agree on the following convention: Except
if specified otherwise, we will equip submanifolds of $\RR^d$ with
the inner product of $\RR^d$ on their tangent spaces. This defines
their Riemannian measure completely --- not only up to a 
multiplicative constant.

\bigskip

\state{NOTATION.} {\it%
  Whenever a set $M$ is a submanifold of $\RR^d$, we write
  $\calM_M$ for its Riemannian measure.
  }

\bigskip

We can now rewrite our integral.

\bigskip

\state{3.1.9. PROPOSITION.} {\it%
  If $I(A)$ is positive, the following equality holds,
  $$\displaylines{
    \int_A\eIx \d x = e^{-I(A)}\int_{p\in\Lambda_{I(A)}}
    e^{-\tau_A(p)}\int_{s\geq 0} e^{-s} 
    { \indic_A\big( \psi_{\tau_A(p)+s}(p)\big) \over
      \big| \D I\big( \psi_{\tau_A(p)+s}(p) \big) \big| }
    \hfill\cr\noalign{\vskip .05in}\hfill
    \big| \det\big( \psi_{\tau_A(p)+s,*}^\T
          \psi_{\tau_A(p)+s,*}^{\phantom T}(p)\big)\big|^{1/2}
    \d s \, \d\calM_{\Lambda_{I(A)}}(p) \, .\qquad\cr
  }$$
  }

\bigskip

\stateit{Proof.} Use Fubini's theorem to first obtain
$$\eqalignno{
  \int_A\eIx \d x
  & = \int\int \indic_A (x) \indic_{[I(x),\infty)}(c) e^{-c} 
    \, \d c \, \d x \cr
  & =e^{-I(A)} \int e^{-c} |A\cap \Gamma_{I(A)+c}|\, \d c \, .
     & (3.1.1)\cr
  }$$
This brings the leading exponential term out of the integral.

Let $p$ be in $\Lambda_{I(A)}$, and consider a local chart
$$
  p(\cdot ) : U\subset\RR^{d-1}\to \Lambda_{I(A)}
$$
around $p$. If $x$ is in the image of $p(U)$ through the normal flow, we can
parameterize it as $x=\psi\big( p(u_1,\ldots , u_{d-1}),s\big)$. Setting
$\partial_i={\partial\over\partial u_i}p$, the Jacobian of the change of 
variable $x\leftrightarrow (u_1,\ldots ,u_{d-1}, s)$ is
$$\eqalign{
  J
  &= \Big| \det \Big( \psi_{s*}(p)\partial_1,\ldots , \psi_{s*}(p)\partial_{d-1} ,
    {\partial\over \partial s} \psi (p,s)\Big) \Big| \cr
  &= \det\left| \pmatrix{%
    \big\<\psi_{s*}(p)\partial_i\, ,\psi_{s*}(p)\partial_j
    \big\>_{1\leq i,j\leq d-1}
      & 0 \cr
    0 & \Big| {\displaystyle\partial\over\displaystyle\partial s}
        \psi \big( p(u),s\big) \Big|^2 \cr
                        }\right|^{1/2} \cr
  & = \Big|{\d\over \d s} \psi_s\Big| \big( \det (\psi_{s*}^\T
    \psi_{s*}^{\phantom T})\big)^{1/2} 
    \det \big( {\< \partial_i,\partial_j\>}_{1\leq i,j\leq d-1}
    \big)^{1/2} \, ,\cr
  }$$
since $\big\< \psi_{s*}(p)\partial_i\, ,{\displaystyle \d\over 
\displaystyle \d s}\psi (p,s)\big\> = 0$, 
for $i=1,\ldots ,d-1$.

In term of the new parameterization and gluing charts using a
partition of unity, we obtain
$$\eqalign{
  |A\cap \Gamma_{I(A)+c}|
  & = \int \indic_{A\cap \Gamma_{I(A)+c}}(x) \d x \cr
  & = \int_{0\leq s\leq c} \int_{\Lambda_{I(A)}} \indic_A \big( \psi (p,s)\big) 
    \Big| {\d\over \d s} \psi (p,s)\Big|\times \cr
  \noalign{\vskip 2mm}
  & \qquad\qquad\qquad \det (\psi_{s*}^\T\psi_{s*}^{\phantom T})^{1/2} 
    \d\calM_{\Lambda_{I(A)}}(p) \, \d s \, . \qquad\cr
  }$$
Notice that $\indic_A\big(\psi (p,s)\big)$ vanishes
if $s\leq \tau_A(p)$. Using (3.1.1) and Fubini's theorem,
the expression of $\d\psi(p,s)/\d s$ in Lemma 3.1.1 
gives the result.\hfill$\qed$

\bigskip

Just as we introduced $\tau_A(p)$, the first en\underbar{t}rance time
of a point in $A$ through the flow, define\ntn{$\chi_A^L(p)$ : last
entrance time in $A$ of the normal flow}
$$
  \chi_A^L(p)=\sup\big\{\, s : \psi_{\tau_A(p)+s}(p)\in A\,\big\}
$$
the \underbar{l}ast time of e\underbar{x}it of $A$ starting the clock time at
$\tau_A(p)$, and\ntn{$\chi_A^F(p)$ : first exit time in $A$ of the
normal flow}
$$
  \chi_A^F(p)=\inf\big\{\, s>0  : \psi_{\tau_A(p)+s}(p)\not\in A\,\big\} \, ,
$$
the \underbar{f}irst time of e\underbar{x}it of $A$ starting the time at
$\tau_A(p)$. If $A$ has ``holes'', $\chi_A^F$ may be strictly less than 
$\chi_A^L$. Both $\chi_A^L(p)$ and $\chi_A^F(p)$ may be infinite.

Proposition 3.1.9 and our lemmas on the normal flow yield the following basic 
bounds, which will turn out to be surprisingly sharp and useful.

\bigskip

\state{3.1.10. THEOREM.} {\it%
  For any Borel set $A$ of $\RR^d$ with $I(A)$ positive,
  $$\displaylines{
    \hfill%
    \int_A\eIx \d x
    \leq e^{-I(A)} \int_{p\in\Lambda_{I(A)}} e^{-\tau_A(p)}
      { \det \big( \psi_{\tau_A(p)*}^\T(p)
        \psi_{\tau_A(p)*}^{\phantom T}(p) \big)^{1/2}
        \over \big| \D I\big(\psi_{\tau_A(p)}(p)\big)\big| }
      \times
    \hfill\cr\noalign{\vskip .05in}\hfill
    \int_{0\leq s\leq \chi_A^L(p)}\exp \bigg[ 
      \int_{0\leq t\leq s}\Big( d{\|\D ^2I\|\over |\D I|^2} 
      \big(\psi_{\tau_A(p)+t}(p)\big)-1\Big)
      \d t \bigg]
    \cr\noalign{\vskip .05in}\hfill
      \d s \, \d\calM_{\Lambda_{I(A)}}(p) \cr
    }$$
  and
  $$\displaylines{\hfill
      \int_A\eIx \d x
      \geq e^{-I(A)} \int_{\scriptstyle p\in\Lambda_{I(A)}\atop%
                           \scriptstyle\tau_A(p)<\infty}
      e^{-\tau_A(p)}
      {\det (\psi_{\tau_A(p)*}^\T(p)
       \psi_{\tau_A(p)*}^{\phantom T}(p))^{1/2} 
       \over
       \big| \D I\big( \psi_{\tau_A(p)}(p)\big)\big| }
      \times
    \hfill\cr\noalign{\vskip .05in}\hfill
      \int_{0\leq s\leq \chi_A^F(p)}\exp\Bigg[
      \int_{0\leq t\leq s}\bigg( 
        {-\| \D ^2I\|\over |\D I|^2}\big( \psi_{\tau_A(p)+t}(p)\big)
         -1 \bigg)
      \d t\Bigg]
    \cr\noalign{\vskip .05in}\hfill
      \d s \, \d\calM_{\Lambda_{I(A)}}(p) \, . \cr}
    $$
  }

\bigskip

\noindent (Be aware that the letter $d$ in front of $\|\D ^2I\|/|\D I|$ in 
the exponential in the upper bound refers to the dimension $d$ of $\RR^d$, 
and not to some differentiation!)

\bigskip

\stateit{Proof.} The upper bound follows from Proposition 3.1.9, and the
following observations. Clearly, 
$\indic_A\big( \psi_{\tau_A(p)+s}(p)\big)$ is at most $1$, and
vanishes whenever $s\geq \chi_A^L(p)$. Lemma 3.1.8 implies
$$\displaylines{
  \qquad
  \det \big( \psi_{\tau_A(p)+s,*}^\T 
  \psi_{\tau_A(p)+s,*}^{\phantom T}\big)
  \cut
  \leq \det \big( \psi_{\tau_A(p),*}^\T 
    \psi_{\tau_A(p),*}^{\phantom T}\big) 
    \exp \Big( 2d\int_{0\leq t\leq s} {\|\D ^2I\|\over |\D I|^2}
    (\psi_{\tau_A(p)+t}) \d t \Big) \, , \qquad\cr
  }$$
and Lemma 3.1.2 gives
$$
  { 1\over \big| \D I\big( \psi_{\tau_A(p)+s}(p)\big)\big| }
  \leq { 1\over \big| \D I\big( \psi_{\tau_A(p)}(p)\big)\big| } \, .
$$
To prove the lower bound, notice first that $\indic_A\big(\psi_{\tau_A(p)+s}(p)\big)
=1$ for all $s$ in $\big[\,0,\chi_A^F(p)\big)$. Using Lemma 3.1.8
and Proposition 3.1.9,
$$\displaylines{
  \qquad%
  \int_A\eIx \d x\geq
  e^{-I(A)}\int_{p\in\Lambda_{I(A)}} e^{-\tau_A(p)}
  { \det\big(\psi_{\tau_A(p)*}^\T
    \psi_{\tau_A(p)}^{\phantom T}(p)\big) \over
    \big| \D I\big( \psi_{\tau_A(p)}(p)\big) \big| } \times
  \cut
  \int_{0\leq s\leq \chi_A^F(p)} \exp \bigg(
    -s-\log\big| \D I\big( \psi_{\tau_A(p)+s}(p)\big)\big|
    + \log\big| \D I\big( \psi_{\tau_A(p)}(p)\big) \big|\bigg) 
  \cut
  \hfill
  \d s \, \d\calM_{\Lambda_{I(A)}}(p) \, . \qquad \cr
  }$$
Using Lemma 3.1.1, we have
$$
  {\d\over \d s}\log \Big| \D I\big( \psi (p,s)\big)\Big|
  =\Big\< {\D ^2I\cdot N\over |\D I|^2} , N\Big\> 
   \big( \psi(p,s)\big) \, .
$$
Therefore,
$$\displaylines{\qquad
    \log \big| \D I\big(\psi_{\tau_A(p)}(p)\big)\big| 
    -\log \big| \D I\big(\psi_{\tau_A(p)+s}(p)\big)\big|
  \cut
    \geq \int_{0\leq t\leq s} -{\| \D ^2I\|\over |\D I|^2} 
    \big( \psi_{\tau_A(p) +t}(p)\big) \d t \, ,
  \qquad\cr}
$$
and this brings the lower estimate. \hfill$\qed$

\bigskip

\state{REMARK.} The gap between the upper and the lower bounds comes essentially from
the term $\|\D ^2I\|/|\D I|^2$ in the exponential. 
One should expect this ratio to be small for large
arguments and in interesting situations. For instance, when $d=1$ and
$I(x)=|x|^\alpha$, we have 
$$
  {|\D ^2I(x)| \over |\D I(x)|^2}
  = {\alpha-1\over\alpha} {1\over |x|^\alpha } \, .
$$
In the same vein, if $\tau_A(p)$ is not too large, we should be able to replace
$\psi_{\tau_A(p)}(p)$ by $\psi_0(p)=p$, while if $\tau_A(p)$ is large, the term
$e^{-\tau_A(p)}$ will make the contribution of those $p$'s negligible. Therefore, we
hope to obtain the approximation
$$\displaylines{\quad
  \int_A\eIx \d x
  \hfill\cr\qquad\quad
  \approx e^{-I(A)}\int_{p\in\Lambda_{I(A)}} {e^{-\tau_A(p)}\over \big| \D I(p)\big|}
   \int_{s\geq 0}\exp\Big( \int_0^s -1 \d t\Big) 
   \d s \, \d\calM_{\Lambda_{I(A)}}(p) 
  \hfill\cr\qquad\quad
  = e^{-I(A)}\int_{p\in\Lambda_{I(A)}} {e^{-\tau_A(p)}\over \big| \D I(p)\big|}
   \d\calM_{\Lambda_{I(A)}}(p) \, , 
  \hfill (3.1.2)\cr}
$$
which is a quite manageable expression. To prove that such approximation is valid
requires some ideas on the order of magnitude of the final expression
in (3.1.2), and the next
section is devoted to the study of this term.

\bigskip

\noindent{\fsection 3.2. Base manifolds and their orthogonal leaves.}
%\chapternumber={3.2}
\section{3.2}{The normal flow and the normal foliation}

\bigskip

\noindent Following the remark concluding the previous subsection, 
we want to rewrite the integral as $e^{-I(A)}$ times 
$$
  \int_{\Lambda_{I(A)}} {e^{-\tau_A}\over \big| \D I\big|}
  \d\calM_{\Lambda_{I(A)}} \, , 
  \eqno{(3.2.1)}
$$
so that we can isolate the leading terms. We first need 
to recall some notation and introduce some definitions.

For $p$ in $\Lambda_c$, we denote by $\exp_p(\cdot )$ the exponential 
map at $p$ in the manifold $\Lambda_c$. That is, if $u$ belongs 
to $T_p\Lambda_c$, the value of $\exp_p(u)$ is the point on $\Lambda_c$ 
at a distance $|u|$ to $p$ on the geodesic starting at $p$ in the
direction $u$ --- see for instance Do Carmo (1992) or Chavel (1996). 
The exponential map is always defined for a small value of the argument 
in the tangent plane.

If $M$ is a submanifold of $\Lambda_{I(A)}$ and $p$ is a point in $M$, 
the tangent space $T_p\Lambda_{I(A)}$ splits as $T_pM\oplus (T_pM)^\perp$ 
where we denote by $(T_pM)^\perp$ the orthocomplement of $T_pM$ 
in $T_p\Lambda_{I(A)}$. 

\bigskip

\state{DEFINITION.} {\it%
  The projection of a set $A\subset\RR^d$ on $\Lambda_{I(A)}$ --- through
  the flow $\psi$ --- is the set of all points
  $p$ in $\Lambda_{I(A)}$ such that $\psi_t(p)$ belongs to $A$ for 
  some positive $t$. A submanifold $\calD_A\subset\Lambda_{I(A)}$ is 
  called a base manifold\ntn{$\calD_A$}
  --- for $A$ --- if the projection of $A$ through $\psi$ is contained in 
  $\bigcup_{p\in\calD_A}\exp_p\big( (T_p\calD_A)^\perp\big)$.
  }

\bigskip

Often, we will loosely speak of the projection of $A$ on
$\Lambda_{I(A)}$, forgetting to add that it is through the flow $\psi$.

In what follows, if $\calD_A$ is a base manifold, we denote by
$k=k(\calD_A)$ its dimension. In order to replace the integration over
$\Lambda_{I(A)}$ by an integration over a base manifold $\calD_A$ in
(3.2.1), we need to attach orthogonal leaves to $\calD_A$. For this
construction, we borrow some definitions and notations used to
parameterize tubes as in Weyl's (1939) formula.

Consider a base manifold $\calD_A$ of $A$, and define a normal bundle by
$$
  N_p\calD_A=T_p\Lambda_{I(A)}\ominus T_p\calD_A 
  \, , \qquad p\in\calD_A \, .
$$
Let ${\rm dist}(\cdot \, ,\cdot )$ denote the Riemannian distance on 
$\Lambda_{I(A)}$. For every $u$ in $N_p\calD_A$, consider the 
radius of injectivity of $\Lambda_{I(A)}$ in the direction $u$,
$$
  e_A(p,u)
  =\sup\big\{\, t\geq 0  : {\rm dist}\big(
                \exp_p(tu),p\big)=t\,\big\} 
  \, ,
$$
the exponential map being that on $\Lambda_{I(A)}$ as before. 
Define
$$
  \Omega_{A,p}=\big\{\, (t,u) : u\in N_p\calD_A \, , \, |u|=1 \, , \,
  0\leq t< e_A(p,u) \,\big\} \, , \qquad p\in\calD_A \, ,
$$
so that $\bigcup_{p\in\calD_A}\partial \Omega_{A,p}$ is the set of all 
focal points of $\calD_A$ immersed in $\Lambda_{I(A)}$. From its definition, 
we infer that $\exp_p(\Omega_{A,p})$ coincide with the set\ntn{$\omega_A$}
$$\displaylines{
  \qquad%
  \omega_{A,p}=\{\, q\in\Lambda_{I(A)} : \hbox{ there exists a unique minimizing}
  \cut
  \hbox{ geodesic through $q$ which meets $\calD_A$ 
         orthogonally at $p\, \}$.}
  \qquad\cr}$$
Set\ntn{$\omega_A$}
$$
  \omega_A=\bigcup_{p\in\calD_A}\omega_{A,p} \, .
$$
On $\omega_A$, we can define a projection $\pi_A$ onto $\calD_A$ as follows.
Any point $q$ in $\omega_A$ can be written in a unique way as $q=\exp_p(u)$ 
for $(p,u)\in\calD_A\times N_p\calD_A$. We set $\pi_A(q)=p$. In other 
words, $q$ is on a unique geodesic
starting from $\calD_A$ and orthogonal to $\calD_A$; the projection of $q$ on
$\calD_A$ is the starting point of this geodesic 
in $\calD_A$. We call the sets 
$\{\, \pi^{-1}_A(p) : p\in\calD_A\, \}$, the orthogonal 
leaves to $\calD_A$. By construction,
$\omega_{A,p}$ is in $\Lambda_{I(A)}=\Lambda_{I(p)}$, and
it contains $\pi^{-1}_A(p)$.

\bigskip

\setbox1=\hbox to 2in{\vbox to 1.5in{\ninepointsss
  \hsize=2in
  \kern 10.7in
  \includegraphics{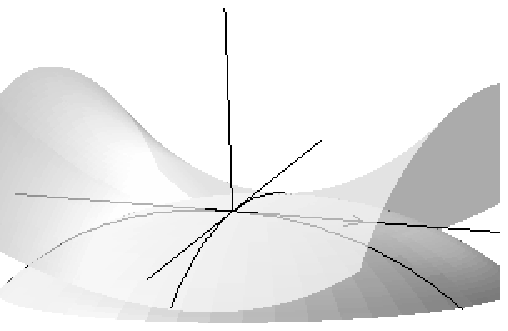}
  \vfill
  \vss
  \anote{.95}{.6}{$p$}
  \anote{1.6}{.25}{$\omega_{A,p}$}
  \anote{1.45}{.1}{or $\pi^{-1}_A(p)$}
  \anote{.3}{1.3}{$\partial A$}
  \anote{0}{.17}{$\Lambda_{I(A)}$}
  \anote{.6}{.17}{$\calD_A$}
}\hfill}    
\setbox2=\hbox to 2in{\vbox to 1.5in{\ninepointsss
  \hsize=2in
  \kern 10.7in
  \includegraphics{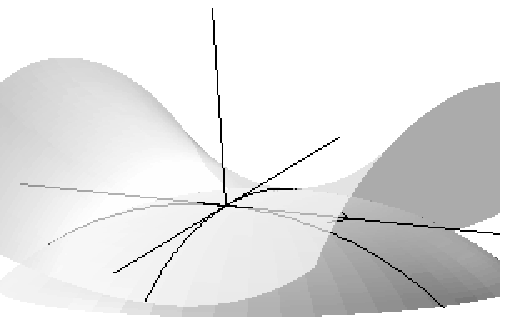}
  \vfill\vss
  \anote{.95}{.6}{$p$}
  \anote{1.53}{.25}{$\omega_{A,p}$}
  \anote{1.38}{.1}{or $\pi^{-1}_A(p)$}
  \anote{.43}{1.31}{$\partial A$}
  \anote{0}{.17}{$\Lambda_{I(A)}$}
  \anote{.5}{.19}{$\calD_A$}
}\hfill}  

\bigskip

\centerline{\box1\hfill\box2}

\bigskip

In order to rewrite the integral (3.2.1) as an integral over a base
manifold and its orthogonal leaves, we need to calculate the Jacobian
of the change of variable
$q\in\omega_A\leftrightarrow \exp_p(u)$, $p\in\calD_A$, $u\in
N_p\calD_A$. For this purpose, notice that the differential
$\pi_{A*}(q)$ maps $T_q\Lambda_{I(A)}$ onto $T_{\pi(q)}\calD_A$. Since
$\pi_A$ is constant on the orthogonal leaves, the orthocomplement of
$\ker \pi_{A*}(q)$ is a vector space of dimension $k$. Let
$b_1(q),\ldots , b_k(q)$ be an orthonormal basis of this
orthocomplement, and define the Jacobian\ntn{$J\pi_A(\cdot )$}
$$
  J\pi_A(q) = \Big|\det 
  \Big( \big( \pi_{A*}(q)\wedge \ldots \wedge \pi_{A*}(q)\big)
  (b_1,\ldots ,b_k)\Big)\Big|
$$
--- when $k=0$, read $J\pi_A(q)=1$.
Federer's (1959) co-area formula --- see the appendix of Howard
(1993) for a simple proof in the smooth case we are using here --- yields
$$\displaylines{\quad
  \int_{\omega_A}{e^{-\tau_A(p)}\over |\D I(p)|}
    \d\calM_{\Lambda_{I(A)}}(p)
    = \int_{u\geq 0}\int_{p\in\calD_A}\int_{q\in\pi^{-1}_A(p)}
    e^{-u}
    { \indic_{(-\infty ,u]}\big(\tau_A(q)\big) \over 
      \big| \D I(q)\big| }
  \cut
    (J\pi_A)^{-1} (q) \d\calM_{\pi_A^{-1}(p)} (q)  \d\calM_{\calD_A}(p)
    \d u
    \, .\qquad\phantom{(3.2.2)}\llap{(3.2.2)}\cr
  }
$$
We can now express the Riemannian measure $\calM_{\pi^{-1}_A(p)}$ over
the leaf $\pi^{-1}_A(p)$ in normal coordinates. For this aim,
following Chavel (1993, chapter 3), let $\calT_{p,t,v}$ denote the
parallel transport along the geodesic from $p$ to $\exp_p(tv)$, and
furthermore, let $R_q$ denote the Riemannian curvature tensor of the leaf
$\pi^{-1}_A\circ\pi_A (q)$ at the point $q$. Denote
$\gamma_{p,v}(t)=\exp_p(tv)$ a point of the geodesic starting from $p$
in the direction $v$ with constant velocity $|v|$. Then, consider the
matrix
$$
  \calR_v(t)=\calT_{p,t,v}^{-1}\circ R_{\exp_p(tv)}
  \big( \gamma'_{p,v}(t),\cdot \big)
  \gamma'_{p,v}(t)\circ\calT_{p,t,v}
$$
defined on $T_p\pi_A^{-1}(p)$. We define a matrix-valued 
function $\calA_p(t,v)$\ntn{$\calA_p(t,v)$} as solving
the differential equation --- in the set of matrices over the 
$(d-k-2)$-dimensional vector space
$T_p\pi^{-1}_A(p)\ominus v\RR=N_p\calD_A\ominus v\RR$ ---
$$
  {\partial^2\over\partial t^2}\calA (t,v)+\calR_v(t)\calA(t,v) = 0 \, ,
$$
subject to the boundary condition
$$
  \calA (0,v)=0 \, , \qquad 
  {\partial\over\partial t}\calA (t,v)\Big|_{t=0}=\Id \, .
$$
The boundary conditions imply
$$
  \det \calA (t,v)\sim \det (t\, \Id_{\RR^{d-k-2}})\sim t^{d-k-2}
  \qquad\hbox{ as } t\to 0 \,
$$
uniformly in $v$ in the unit sphere
$S_{T_p\pi^{-1}_A(p)}(0,1)$ of $N_p\calD_A$. The
expression of the Riemannian measure over a leaf, 
$\calM_{\pi^{-1}_A(p)}$, in normal coordinates yields
$$\eqalignno{
  &\kern -10pt\int_{\omega_A}{e^{-\tau_A(q)}\over |\D I(q)|} 
   \d\calM_{\Lambda_{I(A)}}(q) \cr
  & = \int_{u\geq 0} e^{-u} \int_{p\in\calD_A} 
    \int_{v\in S_{T_p\pi^{-1}_A(p)}}\int_{t\in [0,e_A(p,v)]} 
    { \indic_{[0,u]}\Big(\tau_A\big(\exp_p(tv)\big)\Big)
    \over \big| \D I\big( \exp_p(tv)\big) \big| } \cr
  & \kern 25pt{1\over \big| J\pi_A\big(\exp_p(tv)\big)\big|}
    \det\big( \calA(t,v)\big)
    \, \d t \, \d\mu_p(v) \, \d\calM_{\calD_A}(p) \, \d u \, ,
  &(3.2.3)\cr
  }$$
where $\mu_p$ is the Riemannian measure over the unit
sphere of $T_p\pi^{-1}_A(p)$, centered at the origin.

This last expression looks quite complicated. However, we are almost
done, and the intuition goes as follows.  Roughly, we want to choose
$\calD_A$ as $\Lambda_{I(A)}\cap\partial A$, so that $I$ is minimal in
$A$ over $\calD_A$ --- but we will actually need to have a little bit
of freedom for some applications and make a slightly more subtle
choice.  Due to the term $e^{-u}$, let us concentrate on the range
$u=O(1)$. The assumption $|\D I(p)|$ tends to infinity with $|p|$ will
imply that in good situations $\tau_A\big(\exp_p(tv)\big)$ grows very
fast as a function of $t$, since $\exp (tv)$ is transverse to
$\calD_A$ where $\tau_A$ is minimal. Hence, if the indicator function
of $\tau_A\big(\exp_p(tv)\big) \leq u =O(1)$ is not zero, we must 
have $t$ small. But for small $t$'s and $p\in\calD_A$,
$$
 (J\pi_A)^{-1}\big( \exp_p (tv)\big) \approx (J\pi_A)^{-1}(p) = 1
$$
since $\exp_p(tv)\approx \exp_p(0)=p$ and $\pi_A(p)=p$; 
furthermore, as we mentioned earlier, for small $t$'s
$$
  \det \big( \calA (t,v)\big) \approx t^{d-k-2}
$$
and again, since $\exp_p(tv)\approx p$,
$$
  \big| \D I\big( \exp_p (tv)\big) \big| 
  \approx \big| \D I(p)\big| \, .
$$
Thus, we should expect the right hand side of (3.2.3) to be approximately
$$\displaylines{
  \quad
  \int_{u\geq 0} e^{-u} 
    \int_{p\in\calD_A}{e^{-\tau_A(p)}\over \big| \D I(p)\big|}
    \int_{v\in\calS_{T_p\pi^{-1}_A(p)}(0,1)} 
    \int_{t\in [0,e_A(p,v)]} 
  \cut
  \indic_{(-\infty,u]}\Big( \tau_A\big(\exp_p(tv)\big)-\tau_A(p)\Big)
         t^{d-k-2}
  \, \d t \, \d\mu_p(v) \, \d\calM_{\calD_A}(p) \, \d u 
  \, .\hfill (3.2.4)
  }$$
Assuming that $\partial A$ is smooth, for 
a unit vector $v$ in $T_p\pi_A^{-1}(p)$, the function 
$t\mapsto\tau_A\big(\exp (tv)\big) -\tau_A(p)$,
is minimum for $t=0$. This function should be
approximately quadratic near $0$. The integration in $t$ and $v$
then gives the volume of an ellipsoid, which is related to 
the curvatures of $\Lambda_{I(A)}$ and $\partial A$ near $p$. 
We will also need to prove that in good situations (3.2.1) 
is equivalent to (3.2.2), in which we restricted the 
integration to $\omega_A$.

At this stage it should be noticed that pinching the curvature $R$
yields differential inequalities for $\calA$ and therefore bounds for
(3.2.3). We shall not pursue this line, but the reader may
notice that there are situations where $R$ is easy to calculate ---
for instance if $I(x)=|x|^2$ as in the Gaussian case, then $R=\Id/|x|$
--- or to pinch. Thus, more precise estimates could be obtained, and
even a control of the error terms. This could be useful in some
applications, but it is not clear that it is worth investigating in a
general setting.

\bigskip

\centerline{\fsection Notes}
\section{3}{Notes}

\bigskip

\noindent This chapter builds upon the classical theory of surfaces 
and integration on manifolds. If you don't know any differential
geometry, don't give up!  It took me a long time to find a good
starting point, that is a book that I could read and understand. I
found it when I was visiting the Universit\'e Laval at Quebec! Buy Do
Carmo's (1976, 1992) two books, and start reading the one on curves
and surfaces. If you are as bad learner as I am, do what I did, that
is, all the exercises. Once you read about two dimensional surfaces and
understand that curvature is a geometric name for a second order
Taylor formula, you will have enough intuition to digest the abstract
Riemannian manifolds --- which really copy the classical theory of
surfaces in $\RR^3$. After reading Do Carmo's books, I found Chavel
(1996) and some parts of Spivak (1970) most valuable. Some colleagues
liked Morgan's (1992) book very much as a starting point, others
McCleary's (1994).

The change of variable comes from that in Weyl's (1939) tube
formula. Weyl's paper is very nice to read, and a good part does not
require much knowledge of differential geometry.

Notice that the formula derived in Proposition 3.1.9 does not rely upon
convexity of $I$. It could be of some use in other places.

It is certainly possible to extend the estimates obtained in this
chapters, the previous ones and the next ones, to some integral over
noncompact manifolds. In this case, one should replace the Lebesgue
measure by the Riemannian one. I do not know if this could be of any use.

\vfil\eject

~\vskip 28pt
\chapter{4}{Analyzing the leading term for some smooth sets}
\line{\hfill{\fchapter 4.\quad Analyzing the leading}}
\vskip .1in
\line{\hfill{\fchapter term for some smooth sets}}

\vskip 54pt

\noindent In this chapter we analyze further the quantities
involved in the integral (3.2.3). We will assume that
$$
  \hbox{$\partial A$ and $\Lambda_{I(A)}$ are smooth 
  --- i.e., C$^2$ --- submanifolds of $\RR^d$,}
  \eqno{(4.0.1)}
$$
and furthermore that
\setbox1=\vtop{\hsize=3.5in\noindent%
  the base manifold $\calD_A$ is a smooth --- i.e., C$^2$ ---
  submanifold of $\RR^d$, possibly with boundary.}
$$
  \raise 7pt\box1\eqno{(4.0.2)}
$$
Notice that so far we have a lot of freedom to choose the base
manifold, so that (4.0.2) is not much of a restriction.
Borrowing from the theory of large deviations, it is convenient
to introduce the following notion.

\bigskip

\state{DEFINITION.} {\it%
  A base manifold $\calD_A$ for $A$ is called a dominating
  manifold (for $A$) if, for all $p$ in $\calD_A$ and all
  unit vectors $v$ in $T_p\pi^{-1}_A(p)$, the function
  $$
    t\in \big(-e_A(p,-v),e_A(p,v)\big)\mapsto 
    \tau_A\big( \exp_p(tv)\big) \in [0,\infty )
  $$
  has a local finite minimum at $t=0$. A point in a 
  dominating manifold is called a dominating point.}

\bigskip

Thus, if $\calD_A$ is a dominating manifold,
the set $\partial A$ is pulling away from $\Lambda_{I(p)}$
along the lift through the normal flow of the geodesics 
orthogonal to $\calD_A$. Notice the important fact that if a point $p$ 
is in a dominating manifold $\calD_A$, it may not be in $\partial A$, i.e.,
$\tau_A(p)$ may still be positive. This is an essential difference
from the large deviation theory, or analogously to what would be
considered in Laplace's method. This distinction will turn out to
be crucial for some applications --- see sections 8.2, 8.3 or 10.2. The
downside of allowing this extra freedom is that the main
result has a slightly more involved statement.
But it is worth the extra power provided. Another feature is
that only the points for which $\tau_A\big(\exp_p(tv)\big)$ is
finite matter.

\bigskip

\noindent{\fsection 4.1. Quadratic approximation of
{$\twelvebmit\tbftau_A$} near a dominating manifold.}
%\chapternumber={4.1}
%\section{4.1}{Quadratic approximation of $\tau_{\hbox{\fivecmssit A}}$ near 
%              a dominating manifold}
\section{4.1}{Quadratic approximation of $\tau_A$ near 
              a dominating manifold}
\bigskip

\noindent Under the assumptions (4.0.1)--(4.0.2), if $\calD_A$ is a
dominating manifold and $p$ is one of its points, the function
$t\mapsto\tau_A\big(\exp_p (tv)\big)$ is minimal at $0$ for all unit
vectors $v$ tangent to $\pi^{-1}_A(p)$ at $p$.  It admits a quadratic
approximation. To obtain it, let us denote by $\Pi_{\Lambda_c,q}$
(resp. $\Pi_{\partial A,q}$) the second fundamental form of
$\Lambda_c$ (resp. $\partial A$) at $q\in\Lambda_c$ (resp.
$q\in\partial A$). Those are defined using the unit outward normal
vector $N$ to $\Lambda_c$. At points of $\Lambda_{I(A)}\cap\partial
A$, this unit normal $N$ is also a unit normal for $\partial A$, and
so $\partial A$ is oriented by an extension of $N$. The submanifold
$\pi^{-1}_A(p)\subset \Lambda_{I(A)}\subset \RR^d$ admits a second
fundamental form also associated with the normal field $N$. For $q$
belonging to $\pi^{-1}_A(p)$, we denote it by
$\Pi_{\Lambda_{I(A)},q}^\pi$.\ntn{$\Pi_{\Lambda_{I(A)},q}^\pi$} It is
nothing but the restriction of $\Pi_{\Lambda_{I(A)},q}$ to
$T_q\pi^{-1}_A(p)$.  The leaf $\pi^{-1}_A(p)$ can be lifted to
$\partial A$ through the normal flow in considering the
set
$$
  \psi_A\big(\pi^{-1}_A(p)\big)
  = \big\{\, \psi \big( q,\tau_A(q)\big)\, , \,
          q\in\pi^{-1}_A(p) \,\big\} \, .
$$

\bigskip

\setbox1=\hbox to 2in{\vbox to 1.5in{\ninepointsss
  \hsize=2in
  \kern 10.7in
  \includegraphics{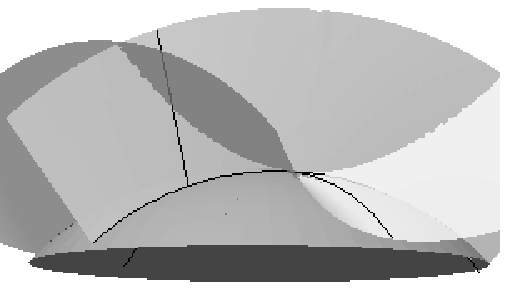}
  \vfill
  \vss
  \anote{1.2}{.75}{$p$}
  \anote{.2}{.25}{$\omega_{A,p}$}
  \anote{.2}{.1}{or $\pi^{-1}_A(p)$}
  \anote{0}{1.2}{$\partial A$}
  \anote{1.2}{.17}{$\Lambda_{I(A)}$}
  \anote{1.5}{.45}{$\calD_A$}
  \anote{1}{1.4}{$\psi_A\left(\pi_A^{-1}(p),\RR\right)$}
  \anote{.7}{.55}{$q$}
  \anote{.65}{.95}{$\circ$}
}\hfill}    
\setbox2=\hbox to 2in{\vbox to 1.5in{\ninepointsss
  \hsize=2in
  \kern 10.7in
  \includegraphics{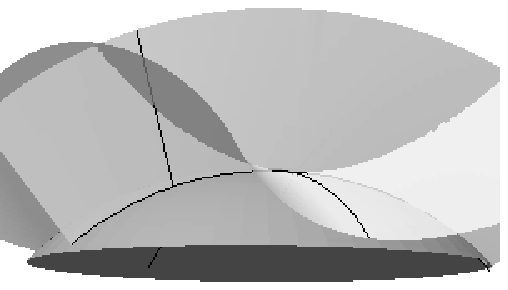}
  \vfill\vss
  \anote{1.12}{.75}{$p$}
  \anote{.15}{.25}{$\omega_{A,p}$}
  \anote{.15}{.1}{or $\pi^{-1}_A(p)$}
  \anote{0}{1.2}{$\partial A$}
  \anote{1.2}{.17}{$\Lambda_{I(A)}$}
  \anote{1.42}{.45}{$\calD_A$}
  \anote{1}{1.4}{$\psi_A\left(\pi_A^{-1}(p),\RR\right)$}
  \anote{.65}{.55}{$q$}
  \anote{.58}{.94}{$\circ$}
}\hfill}  

%\bigskip

\centerline{\box1\hfill\box2}
\kern.1in
\centerline{
\vbox{\hsize=4in\ninepointsss\noindent On this picture, the line leaving the point
  $q$ is the normal flow $\psi (q,t)$ for $t\geq 0$. It crosses
  the boundary $\partial A$ at the circled point 
  $\psi \left(q,\tau_A(q)\right)$.
  When $q$ moves on the geodesic $\pi^{-1}_A(p)$, the crossing 
  point describes $\psi_A(\pi^{-1}(p))$.}
}

\bigskip

This lifted leaf admits
a second fundamental form relative to the extension of $N$ on
$\partial A$. At any point $p$ in $\calD_A$ for which
$\tau_A(p)$ vanishes, this second fundamental form is the restriction
of $\Pi_{\partial A,p}$ to the tangent space of the lifted
manifold, which coincides with the tangent space
$T_p\pi^{-1}_A(p)=N_p\calD_A$. We denote by $\Pi_{\partial A,p}^\pi$
this restriction. The difference of the fundamental forms,
$\Pi_{\Lambda_{I(A)},p}^\pi-\Pi_{\partial A,p}^\pi$ can be
interpreted as follows. Assume we live on the leaf
$\pi^{-1}_A(p)\subset\Lambda_{I(A)}$, and look 
at $\partial A$ along the ``vertical''
direction given by the normal flow. From $\pi^{-1}_A(p)$, the
boundary $\partial A$ pulls away as we move away from
$p\in\pi^{-1}_A(p)\cap \partial A$. The difference of the
fundamental forms is a measure of the curvature of $\partial A$
viewed from $\pi^{-1}_A(p)$, with vertical distance measured in
Euclidean distance along the normal flow. However, as far as the
integration goes, the right measure of vertical distance is on the
increments of the function $I$. At an infinitesimal level,
the increment is $1 /|\D I(p)|$ near $p$ --- see Lemma
3.1.1. Thus, in the geometry of the level sets of $I$, the
bending of $\partial A$ away from $\Lambda_{I(A)}$ is measured
by $|\D I(p)| (\Pi_{\Lambda_{I(A)}}^\pi -\Pi_{\partial
A,p}^\pi)$.

Equipped with this interpretation, the following result is very
natural.

\bigskip

\state{4.1.1. PROPOSITION.} {\it%
  Let $p(s)$ be a curve on $\pi^{-1}_A(p)$, such that $p(0)=p$ 
  belongs to $\partial A\cap \Lambda_{I(A)}$. Under
  (4.0.1)--(4.0.2),
  $$
    \tau_A \big( p(s)\big) ={s^2\over 2} |\D I(p)|
    \big\< (\Pi_{\Lambda_{I(A)},p}^\pi - \Pi_{\partial A,p}^\pi 
    )p'(0)\, , p'(0)\big\> + o(s^2)
  $$
  as $s$ tends to zero.}

\bigskip

\stateit{Proof.} We will make use of the
following elementary fact. If $q(s)$ is a curve on a surface
$M$ and $\nu$ is a unit normal vector field on $M$, then
$\big\< \nu(q(s)) ,q'(s)\big\>=0$. Differentiating with
respect to $s$,
$$
  \< \d\nu (q) q',q'\> + \< \nu (q) ,q''\> = 0 \, .
$$

Since $p$ belongs to $\partial A\cap \Lambda_{I(A)}$, the function
$q\in\Lambda_{I(A)}\mapsto\tau_A(q)\in\RR$ is minimal at $p$.
Hence, $\D\tau_A(p)$ is collinear to $\D I(p)$. Differentiating at
$s=0$ the relation $\tau_A(q)-s=\tau_A\big(\psi(q,s)\big)$
yields $\D\tau_A(p)=-\D I(p)$. 

Next, let $p(s)$ be a curve on $\Lambda_{I(A)}$ such that
$p(0)=p$. Let $t=p'(0)$ be its tangent vector at $p$, and set
$v=p''(0)$. Define
$$
  \theta_1=-\< \D I(p)\, ,v\> + 
  \< \D^2\tau_A(p)t\, ,t\> \, .
$$
A Taylor expansion yields $\tau_A\big( p(s)\big)
=s^2\theta_1/2+o(s^2)$. To prove Proposition 4.1.1, we 
just need to find an expression for $\theta_1$. Denote 
$$
  g(s)=\psi\big[ p(s),\tau_A\big( p(s)\big)\big]
  = \psi \big( p(s),s^2\theta_1/2 + o(s^2)\big) \, ,
$$
the curve $p(s)$ lifted to $\partial A$ through the normal flow.
A Taylor ex\-pan\-sion and Lemma 3.1.1 gives
$$
  g(s)
  = p + st +{s^2\over 2} \Big( v+\theta_1{\D I\over |\D I|^2}(p)
  \Big) + o(s^2) \, .
$$
In particular, the curves 
$\psi \big[ p(s),\tau_A\big(p(s)\big)\big]$ and $p(s)$ have 
the same tangent vector $t$ at $s=0$. Writing $\nu$
for a unit normal vector field extending $N(p)$ on $\partial
A$, using the Weingarten map, and the elementary fact at the
beginning of this proof,
$$\eqalign{
  \< \Pi_{\partial A,p}t\, ,t\>
  & = \< \d\nu(p)g'(0)\, , t\>
    = -\< \nu (0)\, , g''(0)\> \cr
  & = \Big\< N(p),v+\theta_1{N(p)\over |\D I(p)|}
    \Big\> \, . \cr
  }
$$
This gives an expression for $\theta_1$ involving
$$
  \< N(p)\, ,v\>
  = \big\< N(p)\, ,{\d^2\over \d s^2}p(s)\big|_{s=0}\big\>
  = -\< \d N(p)t\, ,t\> 
  = \< \Pi_{\Lambda_{I(A)},p}t\, ,t\>
$$
--- using the elementary fact again --- and ultimately the result.
\hfill$\qed$

%%%%%%%%%%%%%%%%%%%%%%%%%%%%%%%%%%%%%%%%%%%%%%%%%%%%%%%%%
% Section 4.2.
%%%%%%%%%%%%%%%%%%%%%%%%%%%%%%%%%%%%%%%%%%%%%%%%%%%%%%%%%

\bigskip

\noindent{\fsection 4.2. Approximation for det\bfcalA ($\twelvebmit t,v$).}
%\chapternumber={4.2}
\section{4.2}{Approximation for $\det\calA (t,v)$}

\bigskip

\noindent In this section we obtain some bounds on the
determinant of the matrix $\calA (t,v)$ involved in formula
(3.2.3). They will be instrumental in approximating
(3.2.3) further. These bounds are given by the inequalities of Gunther
(1960) and Bishop (1977) --- see, e.g., Chavel (1993, pp.
118-121).
They involve the curvature tensor of the surface
$\pi^{-1}_A(p)\subset\Lambda_{I(p)}$ and so we will first
compute that of the level set $\Lambda_c$.

In order to avoid any ambiguity, recall that if $E$ is a
subspace of $\RR^d$, the compression of a $d\times d$ matrix
$M$ to $E$ is the linear operator from $E$ into $E$ obtained in
restricting $M$ to $E$ and then projecting the image of $M$
into $E$. Thus, writing $\Proj_E$ for the projection from
$\RR^d$ to $E$ and $|_E$ for the restriction to $E$, the
compression of $M$ to $E$ is $\Proj_E \circ M|_E$.

\bigskip

\state{4.2.1. LEMMA.} {\it%
  Let $S(p)$ be the compression of $\D^2I(p)/\big| \D I(p)\big|$
  to $T_p\Lambda_{I(p)}$. The curvature tensor of the surface
  $\Lambda_{I(p)}$ at $p$ is given by
  $$\displaylines{
    \qquad
    R : x,y,z\in T_p\Lambda_{I(p)}\mapsto
    R(x,y)z
    \cut
    =\< S(p)x,y\> S(p) y
    - \< S(p)y,z\> S(p) x \in T_p\Lambda_{I(p)} \, .
    \qquad \cr}
  $$}

\bigskip

\stateit{Proof.} Let $\Proj_{T_p\Lambda_{I(p)}}$
denote the orthogonal projection from $\RR^d$ onto
$T_p\Lambda_{I(p)}$. Furthermore, we write $\nabla$ for the standard
Riemannian connection on $\RR^d$. For any $x$ in $T_p\Lambda_{I(p)}$,
$$
  -\Proj_{T_p\Lambda_{I(p)}}\nabla_xN(p)
  = {1\over |\D I(p)|}\Proj_{T_p\Lambda_{I(p)}}
  \nabla_x \D I(p) = S(p)x \, .
$$
It follows that the second fundamental form of
$\Lambda_{I(p)}\subset \RR^d$ is the bilinear map
$$
  B: x,y\in T_p\Lambda_{I(p)}\mapsto
  B(x,y)=\< S(p)x,y\> N(p)\in\RR^d\ominus
  T_p\Lambda_{I(p)} \, .
$$
The result follows from the immersion of $\Lambda_{I(p)}$ in
the space $\RR^d$ of null curvature, and say, Theorem 2.2 in
Chavel (1993). \hfill$\qed$

\bigskip

We now obtain an expression for the sectional and the Ricci
curvature of $\Lambda_{I(p)}$. 

\bigskip

\state{4.2.2. LEMMA} {\it%
  The sectional curvature of the surface $\Lambda_{I(p)}$ at $p$
  is 
  $$\displaylines{\qquad
    K: x,y\in T_p\Lambda_{I(p)}\mapsto K(x,y)
    \cut
    ={\< S(p)x,x\> \< S(p)y,y\> - \< S(p)x,y\>^2\over 
      |x|^2|y|^2-\< x,y\>^2}\in\RR \, .
    \qquad\cr}
  $$
  Its Ricci curvature is
  $$\displaylines{\qquad
    \Ricc : x,y\in T_p\Lambda_{I(A)}\mapsto \Ricc (x,y)
    \cut
    =\sum_{1\leq j\leq d-1} \big(
    \< S(p)x,y\> \< S(p)e_j,e_j\>
    -\< S(p)x,e_j\> \< S(p)e_j,y\>\big)
    \in\RR \, , \qquad\cr}
  $$
  where $e_1,\ldots , e_{d-1}$ is any orthogonal basis of
  $T_p\Lambda_{I(p)}$.}

\bigskip

\stateit{Proof.} The expression of the sectional
curvature follows from Theorem 2.2 in Chavel (1993) say, and
the calculation of the Weingarten map
$\Proj_{T_p\Lambda_{I(p)}}\nabla_xN$ in the proof of Lemma
4.2.1. The expression for the Ricci curvature follows 
from Lemma 4.2.1. \hfill$\qed$

\bigskip

We can now pinch --- i.e., bound above and bellow --- the curvatures,
and for this purpose, let us denote $\lambda_{\min}\big(\D^2I(p)\big)$
and $\lambda_{\max}\big( \D^2I(p)\big)$ the smallest and the largest
eigenvalues of the symmetric definite positive $d\times d$ matrix
$\D^2I(p)$.

\bigskip

\state{4.2.3. PROPOSITION.} {\it%
  For any $x$ in $T_p\Lambda_{I(p)}$,
  $$
    \Ricc (x,x)\geq (d-2)
    {\lambda_{\rm min}\big( \D^2I(p)\big)\over |\D I(p)|^2}
    |x|^2 \, .
  $$
  If $x,y$ are in $T_p\Lambda_{I(p)}$ and are orthogonal --- i.e., $\< x,y\>
  = 0$ --- then
  $$
    K(x,y)\leq 
    {\lambda_{\rm max}\big( \D^2I(p)\big)^2\over |\D I(p)|^2}
    |x| |y| \, .
  $$}

\bigskip

\stateit{Proof.} Since $\D^2I$ is a symmetric
definite positive matrix, so is $S(p)$. Let $0\leq
s_1\leq\cdots\leq s_{d-1}$ be the eigenvalues and $e_1,\ldots
,e_{d-1}$ be the corresponding orthonormal basis of
eigenvectors of $S(p)$. Let us denote by $x=\sum_{1\leq i\leq
d-1}x_ie_i$ and $y=\sum_{1\leq j\leq d-1}y_ie_i$ two vectors in
$T_p\Lambda_{I(p)}$. We then have
$$\eqalign{
  \Ricc (x,x)
  & = \sum_{1\leq j\leq d-1}\Big( \sum_{1\leq i\leq d-1}
    x_i^2s_is_j-x_j^2s_j^2\Big) \cr
  & = \sum_{1\leq i\leq d-1} x_i^2 s_i
   \Big( \trace \big( S(p)\big)-s_i\Big) \, . \cr
  }$$
If $s_{d-1}\geq \trace \big( S(p)\big) /2$, the function $s\in
[\, s_1,s_{d-1}\, ]\mapsto s\big( \trace S(p)-s\big)$ 
is increasing on $[\, s_1,\trace S(p)/2\, ]$, decreasing 
on $[\, \trace S(p)/2,s_{d-1}\, ]$. 

Otherwise, if $s_{d-1}\leq \trace S(p)/2$, this function is
increasing. Thus, in any case, for any $i=1,\ldots , d-1$,
$$\eqalign{
  s_i\big( \trace S(p)-s_i\big) 
  &\geq s_1\big( \trace S(p)-s_1\big)\wedge
    s_{d-1} \big( \trace S(p)-s_{d-1}\big)  \cr
  \noalign{\vskip 1pt}
  &\geq (d-2) s_1^2 \, . \cr
  }$$
Since $S(p)$ is a compression of $\D^2I(p)/|\D I(p)|$, the
smallest eigenvalue of $\D^2I(p)/|\D I(p)|$ is larger than
$s_1$. The bound on the Ricci curvature follows.

To obtain an upper bound for the sectional curvature, notice that if
$x$ and $y$ are of unit norm and orthogonal, then
$$\eqalign{
  K(x,y)
  &= \sum_{1\leq i\leq d-1} s_ix_i^2\sum_{1\leq i\leq d-1}
    s_iy_i^2 - \Big( \sum_{1\leq i\leq d-1} x_iy_is_i\Big)^2
    \leq s_{d-1}s_{d-1} \cr
  & \leq \lambda_{\rm max}\big( \D^2I(p)\big)^2 \big/
    |\D I(p)|^2 \, ,\cr
  }$$
where the last inequality comes from the fact that $S(p)$ is a
compression of $\D^2I(p)/\D I(p)$. \hfill$\qed$

\bigskip

Let us denote by $\inj_{\pi^{-1}_A(p)}(p)$ the radius of 
injectivity of
$p$ in the manifold $\pi^{-1}_A(p)$. Since $\pi^{-1}_A(p)$ 
is built of geodesics of $\Lambda_{I(p)}$, we have
$$
  \inj_{\pi^{-1}_A(p)}(p)
  =\inf\big\{\, c(p,u)\wedge e_A(p,u) : u\in T_p\pi^{-1}_A(p)
  \, , \, |u|=1\,\big\}  \, ,
$$
where $c(p,u)$ is the Riemannian distance from $p$ to its cut
point on $\Lambda_{I(p)}$ along the geodesic leaving $p$ in the
direction $u$ --- notice that this property on 
$\inj_{\pi^{-1}_A(p)}(p)$ is very specific to the way we
constructed the leaves, and to the fact that we only consider
the point $p$.

From the previous lemmas, we can deduce some bounds on
$\det\calA (t,v)$, using some classical volume comparison
theorems.

\bigskip

\state{4.2.4. PROPOSITION.} {\it%
  (i) (Bishop-Gunther) Assume that for any unit vector
  $v$ in $T_p\pi^{-1}_A(p)$, the Ricci curvature along the
  geodesic $\exp_p (tv)$ in $\pi^{-1}_A(p)$ is nonnegative.
  Then, for all $t$ in $\big[\, 0,\inj_{\pi^{-1}_A(p)}(p)\big)$, 
  the inequality $\det\calA (t,v)\leq t^{d-k-2}$ holds.

  \noindent (ii) For any positive $t$, let
  $$
  K_{\rm max}(p,t) = \sup\Big\{\, 
    { \lambda_{\rm max}\big(\D^2I(q)\big)^2 \over |\D I(q)|^2 } 
    : q\in\pi^{-1}_A(p)\, , \,  d(p,q)\leq t\,\Big\} \, .
  $$
  Then, for any positive $t_0$ and any $t$ in $\Big[ \, 0,
  {\displaystyle\pi \over\displaystyle\sqrt{K_{\rm max}(p,t_0)}}
  \wedge e_A(p,v)\Big)$,
  $$
    \det\calA (t,v) \geq
    \bigg[ { \sin\big( \sqrt{K_{\rm max}(p,t_0)}t\big)\over
             \sqrt{K_{\rm max}(p,t_0)}\big) }
    \bigg]^{d-k-2} \, .
  $$}

\bigskip

\stateit{Proof.} (i) is Bishop's (1977) or
G\"unther's (1960) theorem. Assertion (ii) follows from
Bishop's (1977) theorem --- see, e.g., Chavel, 1996, 
pp.118-120 ---
provided we have the proper upper bound on the Ricci curvature
along the geodesic $\exp_p(tv)$ in $\pi^{-1}_A(p)$. 
For $x,y$ in $T_q\pi^{-1}_A(p)$, denote $K_{\pi^{-1}_A(p),q}(x,y)$
the sectional curvature at $q$ of the manifold $\pi^{-1}_A(p)$.
Then (ii) follows from the bound
$$
  K_{\pi^{-1}_A(p),q}(x,y)
  \leq K_{\Lambda_{I(p)},q}(x,y) 
  \leq {\lambda_{\rm max}\big( \D^2I(q)\big)^2\over%
    \| \D I(q)\|^2} \, .
$$
To prove this inequality, extend $y$ to a vector field in
the tangent bundle
$T\pi^{-1}_A(p)$. The second fundamental form of the immersion
$\pi^{-1}_A(p)\subset \Lambda_{I(p)}$ is given by
$B(x,y)=\Proj_{T_q\pi^{-1}_A(p)}( \nabla_xy)$ where
$\nabla$ is the connection on $\Lambda_{I(p)}$. In particular,
if $x$ is parallel --- which is the case for $x={\d\over
\d t}\exp_p(tv)$ --- then $B(x,x)=0$. Hence --- see, e.g., 
Theorem 2.2 in Chavel, 1996 --- if $x$ and $y$ are of unit norm,
$$
  K_{\pi^{-1}_A(p),q}(x,y) 
  = K_{\Lambda_{I(p)},q}(x,y) - \big| B(x,y)\big|^2
  \leq K_{\Lambda_{I(q)},q} \, .
$$
We conclude in using Proposition 4.2.3. \hfill$\qed$

%%%%%%%%%%%%%%%%%%%%%%%%%%%%%%%%%%%%%%%%%%%%%%%%%%%%%%%%%
% Section 4.3.
%%%%%%%%%%%%%%%%%%%%%%%%%%%%%%%%%%%%%%%%%%%%%%%%%%%%%%%%%

\bigskip

\noindent{\fsection 4.3. What should the result be?}
%\chapternumber{4.3}
\section{4.3}{What should the result be?}

\bigskip

\noindent Combining the results of sections 4.1 and 4.2,
the heuristic argument at the end of section 3 and (3.1.2), 
we see that we should expect
$$\displaylines{
  \int_A\eIx \d x \approx
  e^{-I(A)} \int_{u\geq 0} e^{-u} \int_{p\in\calD_A}
  {e^{-\tau_A(p)}\over\big| \D I(p)\big|}
  \int_{v\in S_{T_p\pi^{-1}_A(p)}(0,1)} 
  \cut
  \int_{t\in [0,e_A(p,v)]} 
  \indic_{[0,u]}\Big( {1\over 2}\big| \D I(p)\big|
  \big\< G_A(p)tv,tv\big\> \Big)
  t^{d-k-2} 
  \cr\hfill
  \d t \, \d\mu_p(v) \, \d\calM_{\calD_A}(p) \, 
  \d u \, ,
  \cr
  }$$
where $G_A(p)=\Pi_{\Lambda_{I(A)},p}-\Pi_{\partial A,p}$.
Assuming for the time being that everything works as expected,
we can complete the calculation of the asymptotic equivalent.
Recall that the volume of the unit ball of $\RR^n$ is
$\omega_n=\pi^{n/2}/\Gamma({n\over 2}+1)$.

\bigskip

\state{4.3.1. PROPOSITION.} {\it%
  For any $a,b$ both positive, the following holds,
  $$\eqalign{
    & \lim_{u_0\to\infty} \int_{u\in [0,u_0]} e^{-au}
      \int_{p\in\calD_A} {e^{-\tau_A(p)}\over |\D I(p)|}
      \int_{v\in S_{T_p\pi^{-1}_A(p)}(0,1)}\int_{t\geq 0}
    \cr
    & \phantom{\lim_{u_0\to\infty} }
      \indic_{[0,bu]}\Big(\, {1\over 2} |\D I(p)|
      \< G_A(p)tv,tv\> \,\Big) t^{d-k-2} \, 
      \d t \, \d\mu_p(v) \, \d\calM_{\calD_A}(p) \, \d u
    \cr
    \noalign{\vskip .1in}
    =\,&{b^{(d-k-1)/2}\over a^{(d-k+1)/2}}
      2^{(d-k-1)/2}\Gamma\Big( {d-k+1\over 2}\Big) 
      \omega_{d-k-1}\times
    \cr
    & \phantom{\lim_{u_0\to\infty} }
      \int_{p\in\calD_A} 
      { e^{-\tau_A(p)}\over 
        |\D I(p)|^{(d-k+1)/2}\big( \det\, G_A(p)\big)^{1/2} } \, 
      \d\calM_{\calD_A}(p)\, .
    \qquad\cr
  }$$}

\bigskip

\state{REMARK.} In Proposition 4.3.1, the integral in $t$ is in
the range $[\, 0,\infty )$, while it is in the range 
$[\, 0,e_A(p,v))$ in the integral at the beginning of 
this section. But this will not make any difference ultimately.

\bigskip

\stateit{Proof.} Notice that $t^{d-k-2} \d t
\, \d\mu_p(v)$ is the Lebesgue measure in
the tangent space $T_p\pi^{-1}_A(p)\equiv
\RR^{d-k-1}$. Therefore,
$$\displaylines{\quad
    \int_{ v\in S_{T_p\pi^{-1}_A(p)}(0,1) } \int_{t\geq 0}
    \indic_{[0,bu]}\Big( {1\over 2}|\D I(p)|
    \big\< G_A(p)tv,tv\big\>\Big) t^{d-k-2} \d t \, \d\mu_p(v)
  \cut
  \eqalign{
    = \, & \bigg| \Big\{\, x\in\RR^{d-k-2}:
           \big\< G_A(p)x,x\big\> \leq 
           { 2bu\over |\D I(p)| } \Big\}\bigg|_{\RR^{d-k-1}} \cr
    = \, & \Big( { 2bu\over |\D I(p)|} \Big)^{(d-k-1)/2}
           { \omega_{d-k-1}\over \big( \det\, G_A(p)\big)^{1/2} } \, . 
    \cr}
  \qquad\cr
  }$$
Hence, the integral for which we want to take the limit
as $u_0$ tends to infinity is
$$\displaylines{\qquad
    \int_{u\in [0,u_0]} e^{-au} \int_{p\in\calD_A}
    {e^{-\tau_A(p)}\over \big| \D I(p)\big|^{(d-k+1)/2} }
    (2b)^{(d-k-1)/2}\omega_{d-k-1} \times
  \cut
    u^{(d-k-1)/2}
    {\d\calM_{\calD_A}(p)\over \big( \det\, G_A(p)\big)^{1/2} }
    \, \d u \, .
  \qquad\cr}
$$
We perform the integration in $u$ and let $u_0$ tends to
infinity to obtain the term $\Gamma\big( (d-k+1)/2\big)$ in the
Proposition. \hfill$\qed$
\vfil\eject

~\vskip 28pt
\chapter{5}{The asymptotic formula}
\line{\hfill{\fchapter 5.\hfill The asymptotic formula}}
\vskip 54pt

\noindent We now have everything that we need to derive our
asymptotic approximation for $\int_A \eIx \d x$ as the set $A$
tends to infinity nicely. The assumptions that we require may
look quite bad at first glance. However, the reader will see in
the next chapters that they are in fact quite well tailored for
applications.

We will assume that for any fixed positive $M$, there 
exists some positive number
$c_{A,M}$, depending on $A$ and meeting the following
requirements. 

We first assume that there exists a manifold $\calD_A$ such that
\setbox1=\vtop{\hsize=3.5in\noindent%
  $\calD_A$ is a dominating manifold for the set 
  $A\cap\Gamma_{I(A)+c_{A,M}}$, of fixed dimension $k$.
  }
$$
  \raise 7pt\box1\eqno{(5.1)}
$$
Hypothesis, (5.1) contains two key requirements. First, $k$ does not
depend on $A$, restricting the class of sets $A$ that we
consider. Second, if $\calD_A$ is a dominating manifold, it has
to be a base manifold. Therefore, $\partial A$ can pull away
from $\calD_A$ only in the orthogonal directions. This is a
restriction for instance if $\calD_A$ is a closed curve with two
boundary points; we may want to allow $\partial A$ to pull in the
outward tangent directions at the boundary points. However, by
breaking $A$ into smaller pieces, the restriction can be
overcome in practice.

Let us denote by $\underline{A}_M$\ntn{$\underline{A}_M$} the projection of
$A\cap\Gamma_{I(A)+c_{A,M}}$ on $\Lambda_{I(A)}$ through the
normal flow $\psi$, that is
$$
  \underline{A}_M=\big\{\, p\in\Lambda_{I(A)} :
  \tau_A(p)\leq c_{A,M}\,\big\} \, .
$$
Under (5.1), it makes sense to assume that
$$
  \underline{A}_M\subset \bigcup_{p\in\calD_A}\omega_{A,p} \, ,
  \eqno{(5.2)}
$$
where $\omega_{A,p}$ is defined in section 3.2.

When $\partial A$ is smooth, Proposition 4.1.1 gives a quadratic
approximation of $\tau_A$. Of course, such an approximation is
local and does not have any kind of uniformity with respect to
$A$. Moreover, looking at (3.2.4), we also would like to have a
quadratic approximation for $\tau_A(q)-\tau_A(p)$ for $q$ near
$\psi \big(p,\tau_A(p)\big)$. How near? Well, we need to 
cover $\underline{A}_M$. But Proposition 4.1.1 suggests that 
$$
  {1\over 2}\Big| \D I\Big( \psi \big(p,\tau_A(p)\big)\Big)\Big| \,
  \big\< \big( \Pi_{\Lambda_{I(p)+\tau_A(p)}}
              -\Pi_{\partial A\, ,\psi (p,\tau_A(p))} \big)
         (q-p) , q-p \big\>
$$
should be a good approximation for $\tau_A(q)-\tau_A(p)$.
Since $c_{A,M}$ is expected to be small compare to $I(p)=I(A)$,
we should have $\psi \big( p,\tau_A(p)\big)$ quite close to $p$,
and so maybe the quadratic approximation given in Proposition
4.1.1 is just fine. This is what happens in many interesting
examples. But for the
time being, there is no other way than to force it, and assume
that there exists a linear map $G_A(p)$\ntn{$G_A(p)$} on $T_p\Lambda_{I(A)}$
such that
$$\displaylines{
  \lim_{A\to\infty} \sup_{q\in\underline{A}_M}
  \bigg| {\tau_A(q)-\tau_A\big(\pi_A(q)\big) \over
          {1\over 2} \big| \D I\big( \pi_A(p)\big) \big|
	  \big\< G_A\big( \pi_A(q)\big) \exp^{-1}_{\pi_A(q)}(q)
		 \, , \exp^{-1}_{\pi_A(q)}(q)\big\> }
	 -1 \bigg| 
  \hfill\cr\noalign{\vskip .1in}\hfill
  = 0 \, .\quad(5.3)\cr
  }
$$
This can be rewritten as
$$
  \lim_{A\to\infty} 
  \sup_{\scriptstyle v\in T_p\Lambda_{I(A)}\ominus 
    T_p\calD_{A} \atop \scriptstyle \tau_A(\exp_p(|v|))
    \leq c_{A,t}}
  \bigg| {\tau_A\big( \exp_p(v)\big)-\tau_A(p)\over {1\over 2}
	  |\D I(p)| \< G_A(p)v,v\> } -1 \bigg|
  = 0 \, .\,
$$
Now that we have a dominating manifold and a quadratic
approximation, it makes sense to proceed as the heuristic
argument at the end of section 3 suggests. Then, we can guess the
asymptotic equivalent of the integral in using the result of
section 4.3. Define\ntn{$c_{A,M}^*$}
$$
  c_{A,M}^*
  =\inf\bigg\{\, c : L\big( I(A)+c\big) \leq {1\over M}
               e^{-I(A)}\int_{\calD_A}
               { e^{-\tau_A}d\calM_{\calD_A}\over
                 |\D I|^{d-k+1\over 2} (\det\, G_A)^{1\over 2} } 
       \,\bigg\} \, .
$$
Since all our asymptotic analysis is driven by the
desire to have the integral influenced mainly by the behavior
of $\partial A$ near $\calD_A$, and to reduce the integral to
$A\cap\Gamma_{I(A)+c_{A,M}}$\ntn{$c_{A,M}$} where things go well, 
we assume that
$$
  c_{A,M}\geq c_{A,M}^* \, .
  \eqno{(5.4)}
$$
Any $c_{A,M}$ larger than $c_{A,M}^*$ will do, and the reader will
easily see that the larger $c_{A,M}$ is, the more stringent our
assumptions are. Hence, $c_{A,M}^*$ is the best choice, but can seldom
be calculated exactly. In practice, picking for $c_{A,M}$ an
asymptotic equivalent of a multiple of $c_{A,M}^*$ will do. To fix the
ideas, a typical order of magnitude of $c_{A,M}$ is $\log I(A)$ for
the applications that we will study.

Though we want to be able to localize the study of the integral to
points of $A$ near $\calD_A$, we still want $A$ to have some
thickness! In particular, we do not want the main contribution in
the integral to come from the thinness of $A$ --- think for instance of
taking $A=\Gamma_{c+\epsilon}\setminus\Gamma_c$ for $\epsilon= \exp
(-e^c)$ or even smaller, and looking for asymptotics as $c$ tends to
infinity. This can be ruled out by assuming that the first exit time
of the normal flow after a time $\tau_A$ is large enough, namely, that
for all positive $M$,
$$
  \lim_{A\to\infty}\inf_{p\in\underline{A}_M}\chi_A^F(p) =
  +\infty \, .
  \eqno{(5.5)}
$$
Since $\chi_A^F$ is less or equal to $\chi_A^L$, assumption (5.5) implies
$$
  \lim_{A\to\infty} \inf_{p\in\underline{A}_M}\chi_A^L(p)
  = +\infty \, .
$$

In the same spirit, we see that it does not make any
difference in the asymptotic analysis if we replace $A$ by
$A\cap\Gamma_{I(A)+c_{A,M}}$. Thus, if $c_{A,M}$ stays bounded,
the set $A$ is quite small in the geometry of the level sets of
$I$. In this case, we would need to deal with an analogue of
Proposition 4.3.1, but keeping $u_0$ fixed. This would introduce
an incomplete gamma function. In order to simplify the result, 
we assume that $c_{A,M}$ is chosen such that
$$
  \lim_{A\to\infty} c_{A,M}=+\infty \, .
  \eqno{(5.6)}
$$
This assumption is satisfied in all the applications that
follow; but the reader will see that up to changing some
constant in the asymptotics, the proof still goes through
without it.

Our next condition can be explained by first thinking of two difficult
situations. Imagine that around the dominating manifold $\calD_A$, the
set $\partial A$ is pulling away very slowly from
$\Lambda_{I(A)}$. Thus, going in the normal direction to $\calD_A$ on
$\Lambda_{I(A)}$, we need to take $q$ very far away from $\calD_A$ in
order to have $\tau_A(q)$ not too close to $0$. In such circumstances,
we should need an extra rescaling so that in the new scale $\tau_A$
grows faster.

Another difficult situation would be to have the level sets $\Gamma_c$
concentrated along a proper subspace of $\RR^d$ --- think for
instance, if $d=2$, of some ellipsoid with one axis growing like $e^c$
and the other like $c$. On a large scale, the problem would be
essentially lower dimensional. Notice however that in such
circumstances, the curvature of the level set should be of different
orders of magnitude at different points.  Near points of high
curvature, the normal flow could pull away very slowly in the
Euclidean geometry. And so there is not much hope to localize the
problem near the boundary of $A$.

A way to take care of these different situations is to relate
the curvature of the level sets with the rate at which $\partial
A$ pulls away from $\calD_A$ in the normal directions. 
Define\ntn{$t_{0,M}(p)$}
$$
  t_{0,M}(p)=\sup\Big\{\, t : 
  \inf_{v\in S_{T_p\pi^{-1}_A(p)} (0,1)} \tau_A \big( \exp_p (tv)
  \big) \leq c_{A,M}\,\Big\} \, .
$$
Whenever $q$ in $\pi_A^{-1}(p)$ is at distance $t_{0,M}(p)$ or more
from $p$, it satisfies $\tau_A(q)\geq c_{A,M}$. We then assume that
for all positive $M$,
$$
  \lim_{A\to\infty}\sup_{p\in\calD_A}
  \sqrt{K_{\rm max}\big(p,t_{0,M}(p)\big)} \, t_{0,M}(p) = 0 \, .
  \eqno{(5.7)}
$$

In order to be able to use Proposition 4.2.4.i, we impose that
for any $p$ in $\calD_A$ and any unit vector $v$ in 
$T_p\pi^{-1}_A(p)$,
\setbox1=\vtop{\hsize=3.5in\noindent%
  the Ricci curvature along the geodesic $\exp_p(tv)$
  in $\pi^{-1}_A(p)\cap \underline{A}_M$ is nonnegative.
  }
$$
  \raise 7pt\box1\eqno{(5.8)}
$$
%$$\eqalign{
%  &\hbox{the Ricci curvature along the geodesic $\exp_p(tv)$}\cr
%  &\hbox{in $\pi^{-1}_A(p)\cap \underline{A}_M$ is nonnegative.}\cr
%  }\eqno{(5.8)}
%$$
--- I am inclined to believe that convexity of $I$ is
enough to guarantee (5.8) but could not prove it.

It is also possible that the set $A$ and the function $I$ are
such that $|\D I|$ varies widely in a small neighborhood of some
point in the dominating manifold. In such situation, $I$ would
increase very fast in some specific directions, and much slower
in others. Then, the rescaling needed, even in $\pi_A^{-1}(p)$
--- when $\pi_A^{-1}(p)$ is of dimension at least $2$ --- could not be
homogeneous in different directions. This can be ruled out by
assuming that for any positive $M$,
$$
  \lim_{A\to\infty}\sup\Big\{\, \Big| 
  {|\D I(p)|\over |\D I(q)|} -1\Big|   :
  p\in\calD_A\, ,\, q\in\underline{A}_M\cap\pi^{-1}_A(p)\,\Big\}
  = 0 \, .
  \eqno{(5.9)}
$$
In order to proceed along the lines of the final remark in
section 3.1, we need that
$$
  \lim_{A\to\infty}\sup\bigg\{\, \sup_{t\geq 0}
  { \big\| \D^2I\big( \psi_t(p)\big)\big\|\over
    \big| \D I\big( \psi_t(p)\big)\big|^2 }  :
  p\in \underline{A}_M\, , \, \tau_A(p) <\infty\,\bigg\}
  = 0 \, .
  \eqno{(5.10)}
$$
This assumption mainly controls the growth of $I$. A sufficient
condition to guarantee (5.10) is of course to 
have $\|\D^2I(p)\| / |\D I(p)|$ tends to $0$ 
as $|p|$ tends to infinity, which means essentially that
$I$ grows slower than any exponential function.

Actually, we also need a rate of convergence in (5.10), but in a
rather weak sense, namely that
$$
  \lim_{A\to\infty} \sup\bigg\{\, \int_0^{\tau_A(p)}
  {\|\D^2I\|\over |\D I|^2}\big( \psi_u(p)\big)\, \d u  :
  p\in\underline{A}_M\,\bigg\} = 0 \, .
  \eqno{(5.11)}
$$

Finally, we need two technical assumptions in order to carry
over the intuitive argument at the end of section 3.2, namely
that for any positive $M$,
$$
  \lim_{A\to\infty} \sup_{q\in \underline{A}_M}
  \big| J\pi_A(q)-1\big|  = 0 \, ,
  \eqno{(5.12)}
$$
and that for any fixed positive $w$, and, as $A$ moves to
infinity,
$$\displaylines{\qquad
    \int_{\scriptstyle \calD_A\hfill\atop\scriptstyle \tau_A
          \geq  c_{A,M}-w}
    { e^{-\tau_A} \over |\D I|^{d-k+1\over 2} (\det\, G_A)^{1\over 2} }
    \,\d\calM_{\calD_A}
  \cut
    = o\bigg( \int_{\calD_A}
    { e^{-\tau_A} \over |\D I|^{d-k+1\over 2} (\det\, G_A)^{1\over 2} }
    \,\d\calM_{\calD_A} \bigg) \, .
  \qquad (5.13)\cr}
$$

Often $\tau_A(\cdot )$ vanishes on $\calD_A$. Then (5.13) is satisfied
whenever (5.6) is. In particular, this is always the case when
$\calD_A$ is a single point or the union of a finite number of points.
Similarly, (5.11) and (5.12) always hold if $\calD_A$ reduces to a
point.

The reader may legitimately be suspicious about these assumptions,
and how they can be checked in applications. We will show in
nontrivial examples that they are not too difficult
to verify. They reduce the problem of approximating the integral
to much more manageable small problems, which can
be handled by systematic methods. The key point to understand
is perhaps that as $c_{A,M}$ is usually much smaller than
$I(A)$, the set $\underline{A}_M$ is actually quite close to
$\calD_A$ in the scale given by $I$. 

We can now state our main result.

\bigskip

\state{5.1. THEOREM.} {\it%
  Assume that $\partial A\cap \Gamma_{I(A)+c_{A,M}}$ is
  a smooth --- twice differentiable --- manifold for every
  positive $M$. Then, under (5.1)--(5.13),
  $$
    \int_A \eIx \d x \sim
    e^{-I(A)} (2\pi )^{d-k-1\over 2} 
    \int_{\calD_A}
    {e^{-\tau_A}\over |\D I|^{d-k+1\over 2} 
    (\det\, G_A)^{1\over 2}}
    \,\d\calM_{\calD_A} 
  $$
  as $A$ moves to infinity, and with $k=\dim\calD_A$.}

\bigskip

\state{REMARK.} If $\calD_A$ is an open subset of $\Lambda_{I(A)}$, 
then $k=d-1$, and the formula should be read with the determinant of $G_A$
to be $1$.

\bigskip

\stateit{Proof.} We obtain an upper and a lower
bound for the integral.

Let us start with the upper bound. Let $\epsilon$ be a
positive number. Assume first that, for $M$ large enough,
$$
  A=A\cap \Gamma_{I(A)+c_{A,M}}
  \eqno{(5.14)}
$$
provided $I(A)$ is large enough. Thus, for $p$ in
$\Gamma_{I(A)}$, either $\tau_A(p)$ is less than 
$c_{A,M}$ or is infinite.

Recall that $\psi_{0*}$ is the identity since $\psi_0$ is the
identity as well. Consequently,
assumption (5.11), Lemmas 3.1.6 and 3.1.8 imply that 
for $\tau_A(p)\leq c_{A,M}$,
$$\displaylines{\qquad
    1\leq \det(\psi_{\tau_A(p),*}^\T
    \psi_{\tau_A(p),*}^{\phantom T})^{1/2}
  \cut
    \leq \exp\bigg( d \int_0^{\tau_A(p)}
         {\|\D^2I\|\over |\D I|^2} \big( \psi_u(p)\big) \d u \bigg) 
    \leq 1+\epsilon
  \qquad\cr}
$$
provided $I(A)$ is large enough. It then follows from (5.5), 
(5.9), (5.11), Lemma 3.1.2 and Theorem 3.1.10 that
$$
  \int_A\eIx \d x
  \leq e^{-I(A)} \int_{\Lambda_{I(A)}} 
       {e^{-\tau_A(p)}\over |\D I(p)|} \d\calM_{\calD_A}(p)
  \big( 1+o(1)\big) 
$$
as $I(A)$ tends to infinity. Under (5.2), the set $\{\,\tau_A
<\infty\,\}$ is included in $\omega_A$. Consequently, (3.2.3) 
yields, as $A$ tends to infinity,
$$\displaylines{
    \int_A\eIx \d x
    \leq e^{-I(A)}\int_{u\geq 0} e^{-u} \int_{p\in\calD_A}
    {e^{-\tau_A(p)}\over |\D I(p)|}
    \int_{v\in S_{T_p\pi^{-1}_A(p)}(0,1)}
  \cut
    \int_{t\in [0,e_A(p,v))} 
    { \indic_{[0,u]}\Big(\tau_A\big(\exp_p(tv)\big)-\tau_A(p)
		    \Big) \over 
      \big| J\pi_A\big( \exp_p(tv)\big) \big| }
    { |\D I(p)| \over \big| \D I\big( \exp_p(tv)\big)\big| }
    \times
  \cr\noalign{\vskip .1in}\hfill
    \det\calA (t,v)\,
    \d t\, \d\mu_p(v) \, \d\calM_{\calD_A} (p)  \, \d u \,
    \big(1+o(1)\big) \, . 
  \cr
  }
$$
Combine (5.8) and Proposition 4.2.4.i to upper $\det \calA(t,v)$
from above. Use (5.9) to get an upper bound $|\D I(p)| / \big|
\D I\big(\exp_p(tv)\big)\big|$ by $1+\epsilon$, and (5.12) 
bound $\big|J\pi_A \big( \exp_p(tv)\big)\big|$ by $1+\epsilon$. 
Then, assumption (5.3) yields, for $I(A)$ large enough,
\finetune{\hfuzz=4pt}
$$\eqalign{
  \int_A\eIx \d x
  & \leq e^{-I(A)}\int_{u\geq 0} \int_{p\in\calD_A}
    {e^{-\tau_A(p)}\over \big| \D I(p)\big|}
    \int_{v\in S_{T_p\pi^{-1}(p)}(0,1)}
    \int_{t\in [0,e_A(p,v))} \cr
  &\kern 1.1in \indic_{[0,u(1+\epsilon)]} \bigg( {1\over 2}
    \big| \D I(p)\big| \big\< G(p)tv,tv\big\> \bigg)
    t^{d-k-2} \cr
  &\kern 1.1in \, \d t \, \d\mu_p(v)\, \d\calM_{\calD_A}(p) \, 
    \d u \, 
    (1+\epsilon ) \cr
  }$$
\finetune{\hfuzz=0pt}
Extend the integration in $t$ over the domain $[\, 0,+\infty )$
and use Proposition 4.3.1 to conclude
$$\displaylines{\qquad
    \int_A\eIx \d x
    \leq e^{-I(A)} 2^{d-k-1\over 2} 
    \Gamma \Big( {d-k+1\over 2}\Big) \omega_{d-k-1}\times
  \cut
    \int_{p\in\calD_A} {e^{-\tau_A(p)}\over |\D I|^{d-k+1\over 2}}
    {\d\calM_{\calD_A}\over (\det\, G_A)^{1\over 2}}
    (1+\epsilon )^{d-k+7\over 2} \, . \qquad\cr
  }$$
Since $\epsilon$ is arbitrary, this yields the proper upper bound
since
$$
  2^{d-k-1\over 2} \Gamma\bigg( {d-k-1\over 2} \bigg)
  \omega_{d-k-1} = (2\pi )^{d-k-1\over 2} \, .
$$

Before we can drop assumption (5.14), we need to prove the lower
bound. To do so, apply Theorem 3.1.10, Lemma 3.1.6 with
assumptions (5.5) and (5.10) to obtain, as $A$ moves to
infinity,
$$
  \int_A\eIx \d x \geq
  e^{-I(A)} \int_{p\in\underline{A}_M}
  { e^{-\tau_A(p)}\over 
   \big| \D I\big( \psi_{\tau_A(p)}(p)\big)\big| }
  \d\calM_{\Lambda_{I(A)}}(p) \big( 1+o(1)\big)^{-2} \, .
  \eqno{(5.15)}
$$
Notice that
$$\eqalign{
  {\big| \D I\big( \psi_{\tau_A(p)}(p)\big)\big|
   \over |\D I(p)|}
  & = \exp \Big( {1\over 2} \int_0^{\tau_A(p)}
    {\d\over \d s} \log \big| \D I\big( \psi (p,s)\big)\big|^2 
    \d s \Big) \cr
  & = \exp \Big( {1\over 2} \int_0^{\tau_A(p)} 2 \,
    {\< \D^2I\cdot N,N\>\over |\D I|^2}\big(\psi (p,s)\big) \d s
    \Big) \cr
  & \leq \exp \Big( \int_0^{\tau_A(p)}{\| \D^2I\|\over |\D I|^2}
    \big( \psi_s (p)\big) \d s \Big) \, . \cr
  }
$$
This last inequality and assumption (5.11) imply
$$
  \lim_{A\to\infty}\sup\Big\{ \, \Big| 
  { \big| \D I\big( \psi_{\tau_A(p)}(p)\big)\big|
    \over |\D I(p)| } -1 \Big| \, : p\in \underline{A}_M
  \, \Big\} = 0 \, .
$$
Consequently, up to a multiplicative factor of $\big(1+o(1)\big)$, we
can replace $\D I\big( \psi_{\tau_A(p)}(p)\big)$ by $\D I(p)$ in
(5.15). We then use equality (3.2.3) and proceed as follows.  First,
we change the variable $u$ into $\tau_A(p)+w$, and restrict the
integration to $w$ between $0$ and $w_0$ for some positive
$w_0$. Second, we restrict further the domain by integrating only over
the points $p$ in $\calD_A$ with $\tau_A (p)$ less than
$c_{A,M}-w_0$. On this range, we can use assumptions (5.9), (5.7),
Proposition 4.2.4-ii, assumptions (5.12) and (5.3) to obtain, for
$I(A)$ tending to infinity,
$$\displaylines{
    \int_A\eIx \d x 
    \geq
    e^{-I(A)}\int_{w\in [0,w_0]} e^{-w}
    \int_{\scriptstyle p\in\calD_A\hfill\atop 
          \scriptstyle \tau_A(p)\leq c_{A,M}-w_0}
    { e^{-\tau_A(p)}\over \big| \D I(p)\big| } 
  \cut
    \int_{v\in S_{T_p\pi^{-1}(p)}(0,1)} 
    \int_{t\in [0,e_a(p,v))}\indic_{[0,w(1-\epsilon)]}
    \Big( {1\over 2} \big| \D I(p)\big| 
      \big\< G_A(p)tv, tv\big\> \Big) \cr
  \cut
    \d t \, \d\mu_p(v) \, \d\calM_{\calD_A}(p) \, \d w
    \big(1+o(1)\big)^{-2} \, . \cr
  }
$$
Arguing as in the proof of Proposition 4.3.1, we obtain that for
$I(A)$ large enough,
$$\displaylines{
    \int_A\eIx \d x
    \geq e^{-I(A)} 2^{d-k-1\over 2} 
    \Gamma\Big( {d-k+1\over 2}\Big) \omega_{d-k-1}\times
  \cut
    \int_{\scriptstyle p\in\calD_A\hfill\atop 
          \scriptstyle \tau_A(p)\leq c_{A,M}-w_0}
    { e^{-\tau_A(p)}\over |\D I(p)|^{d-k+1\over 2}
      \big( \det\, G_A(p)\big)^{1\over 2} }
    \d\calM_{\calD_A}(p)(1+\epsilon )^{-1} \, .\cr
  }$$
Since $\epsilon$ is arbitrary, assumption (5.13) gives then
$$\displaylines{\quad
  \int_A\eIx \d x 
  \cut
  \geq e^{-I(A)} 2^{d-k-1\over 2}
  \int_{p\in\calD_A} 
  { e^{-\tau_A(p)} \over 
    |\D I(p)|^{d-k+1\over 2} \big(\det\, G_A(p)\big)^{1\over 2} }
  \, \d\calM_{\calD_A}(p)
  \, .\quad}
$$
It remains for us to drop the assumption that $A=A\cap
\Gamma_{I(A)+c_{A,M}}$ is the upper bound. This is immediate; for
a general set $A$, write
$$
  \int_A \eIx \d x = \int_{A\cap \Gamma_{I(A)+c_{A,M}}} \eIx \d x
  + \int_{A\cap\Gamma_{I(A)+c_{A,M}}^{\rm c}} 
  \eIx \d x \, .
$$
For the first integral in the right hand side of the above
equality, the theorem --- proved in this case! --- gives the asymptotic
equivalent. For the second one, it is less than the integral
over $\Gamma_{I(A)+c_{A,M}}^{\rm c}$, that is 
$L\big( I(A)+c_{A,M}\big)$. Assumption (5.4) shows that it has
a negligible contribution to the asymptotics. 

When $k=d-1$ and $\calD_A$ is an open subset of
$\Lambda_{I(A)}$, the asymptotic equivalent in
Theorem 5.1 is nothing but (3.2.1). \hfill$\qed$

\bigskip

\state{5.2. REMARK.} Assumption (5.8) turns to be difficult to
check in practice. The main reasons are that geodesics can seldom be
explicitly calculated, and that the curvature tensor may be difficult
to calculate. However, we only used it to apply Proposition 4.2.4.i
when deriving the upper bound in the proof of Theorem 5.1. It would be
enough to have
$$\displaylines{\qquad
    \lim_{A\to\infty}\sup\bigg\{\, \bigg|
    {\det \calA (t,v)\over t^{d-k-2}}-1\bigg|  :
    v\in S_{T_p\pi^{-1}(p)}(0,1) \, , \,
  \cut\qquad
    t\in [\, 0,e_A(p,v)) \, , \, 
    \exp_p(tv)\in\underline{A}_M \,\bigg\} = 0 \, .
    \qquad\qquad\llap{\hfill(5.16)}\cr
  }
$$
This condition will turn out to be easier to check in many cases.
This could also replace (5.7) as well.

In a similar spirit, (5.3) may be tedious to verify. 
Often $\tau_A$ cannot be easily calculated, but is only
known via an asymptotic expansion as $A$ moves to infinity. 
Due to the
error term in the asymptotic expansion, the uniformity in (5.3),
for small $q$ very close to $p$, may be difficult to check. 
Therefore, we will make a
rather systematic use of the following weaker hypothesis. Assume
that there exists a function $\tilde\tau_A$ defined on
$\underline{A}_M$ such that
$$
  \lim_{A\to\infty}\sup_{q\in\underline{A}_M}
  |\tau_A(q)-\tilde\tau_A(q)| = 0
$$
and
$$
  \lim_{A\to\infty} \sup_{q\in\underline{A}_M}
  \bigg| { \tilde\tau_A(q)-\tilde\tau_A\big(\pi (q)\big)
           \over {1\over 2} \big| \D I\big( \pi_A (q)\big)\big|
           \big\< G_A\big( \pi_A(q)\big) 
                  \exp^{-1}_{\pi_A(q)}(q)
                  \, ,\exp^{-1}_{\pi_A(q)} \big\> 
          } -1 \bigg| = 0
  \eqno{(5.17)}
$$
for some $G(p)$ on $T_p\Lambda_{I(A)}$. Then, Theorem 5.1 holds
when (5.3) is replaced by (5.17). Indeed, let $\epsilon$ be a
positive number. In the proof of Theorem 5.1, we now use 
the bound
$$\displaylines{\qquad
  \indic_{[0,u(1-\epsilon)-\epsilon]} 
    \Big({1\over 2} |\D I(p)|\< G_A(p)tv,tv\>\Big)
  \hfill\cr\vadjust{\vskip 0.05in}\hfill
  \eqalign{
    \leq &\, \indic_{[0,u]} \Big( \tau_A\big( \exp_p(tv)\big) 
         -\tau_A(p)\Big) \cr
    \leq &\, \indic_{[0,u(1+\epsilon )+\epsilon]}
	 \Big( {1\over 2} |\D I(p)| 
	 \< G(p)tv,tv\> \Big) \, . \cr}
  \qquad}
$$
By the dominated convergence theorem,
$$
  \lim_{\epsilon\to 0} \int_{u\geq 0}
  e^{-u} \big( u+(1\pm \epsilon )\pm\epsilon\big)^{d-k-1\over 2}
  \d u 
  = \int e^{-u} u^{d-k-1\over 2} \d u \, ,
$$
and one easily sees that the proof of Theorem 5.1 goes through.

It is sometimes convenient to weaken even a tiny bit this
assumption, only assuming that $\tau_A(q)-\tau_A\big(\pi_A(q)
\big)$ can be approximated by some
$\tilde\tau_A(q)-\tilde\tau_A\big(\pi_A(q)\big)$ in the sense that
$$
  \lim_{A\to\infty} \sup_{q\in\underline{A}_M}
  \Big| \tau_A (q)-\tau_A\big(\pi_A(q)\big)
        -\Big( \tilde\tau_A(q)-\tilde\tau_A\big(\pi_A (q)\big)
 	\Big)\Big|= 0
$$
and of course keeping requirement (5.17).

\bigskip

Let us now explain how Theorem 5.1 can be used to obtain
information on limiting conditional distributions. 

Assume that we consider a log-concave density function
proportional to $e^{-I}$ on $\RR^d$. As $A$ moves away
to infinity, there is
not much hope for the conditional distribution 
$$
  \int_{A\cap B}\eIx \d x \;\Big/ \int_A \eIx \d x
$$
to converge to a nontrivial limit. Indeed, a fixed 
bounded set $B$ does not intersect
$A$ if $I(A)$ is large enough. So, we need to rescale
$B$. For this, consider a normalizing function
$A\mapsto \lambda_A\in (0,\infty )$. We will require that the dominating
manifold $\calD_A/\lambda_A$ converges in a weak sense. But for
the time being, consider the rescaled conditional 
distribution\ntn{$\mu_A(B)$}
$$
  \mu_A(B)=\int_{A\cap\lambda_A B} \eIx \d x
  \;\Big/ \int_A \eIx \d x \, .
$$
It converges weakly* if for any continuous and bounded function
$f$ on $\RR^d$, the integral
$$
  \int f\, \d\mu_A
  = \int_A f(x/\lambda_A)\eIx \d x \;\Big/ \int_A \eIx \d x
$$
converges as $A$ moves to infinity. Its limit is a linear form in $f$,
associated to a measure, the weak* limit of $\mu_A$.

Let us assume that
 
\vskip\abovedisplayskip
{ \line{\kern .3truein{}as $A$ moves to infinity, the 
        family of rescaled measures\hfill}
  \kern .1in
  \line{\kern .3truein
    ${ \displaystyle\lambda_A^k e^{-\tau_A(\lambda_A q)}
      \d\calM_{\calD_A/\lambda_A} (q)
      \over \displaystyle |\D I(\lambda_Aq)|^{d-k+1\over 2} 
      \big( \det\, G_A(\lambda_Aq)\big)^{1\over 2} }
    \; \bigg/
      {\displaystyle\int}_A e^{-\tau_A}{\displaystyle \d\calM_{\calD_A} \over
      \displaystyle |\D I|^{d-k+1\over 2}(\det\, G_A)^{1\over 2} }
    $
    %\hfill(5.18)%
  }
  %\vskip\belowdisplayskip
  \kern .1in
  \line{\kern.3truein{}converges weakly* to a probability 
        measure $\nu$.\hfill (5.18)}
}
\vskip\belowdisplayskip

\noindent Assume furthermore that
$$
  \lim_{A\to\infty} {c_{A,M}\over\lambda_A}
  \sup\big\{ \, |\D I(p)|^{-1}\, :\, p\in \underline{A}_M\, \big\}
  = 0 
  \eqno{(5.19)}
$$
and
$$
  \lim_{A\to\infty}\sup\big\{\, |q-p|/\lambda_A \, :\, 
  q\in\pi^{-1}_A(p) \, , \, \tau_A(p)\leq c_{A,M}\, \}
  = 0 \, .
  \eqno{(5.20)}
$$
We then have the following convergence.

\bigskip

\state{5.3. COROLLARY.} {\it
  Under the assumptions of Theorem 5.1 and (5.15)--(5.20), the
  conditional distribution $\mu_A$ converges weakly* to
  $\nu$ as $I(A)$ tends to infinity.
  }

\bigskip

\stateit{Proof.} Argue as in the proof of Proposition 3.1.9 to
obtain
\finetune{\hfuzz=3.5pt}
$$\hskip -2pt\eqalign{
  \int_A \eIx f(x/\lambda_A) \, \d x
  & = \int\int \indic_A (x)\indic_{[I(x),\infty )}(c)
    e^{-c} f(x/\lambda_A) \, \d c \, \d x \cr
  & = e^{-I(A)} \int\int e^{-c}\, \indic_{A\cap\Gamma_{I(A)+c}}
    (x) f(x/\lambda_A) \, \d c \, \d x \cr
%  \vadjust{\vfill\eject} 
  & = e^{-I(A)}\int_{c\geq 0} e^{-c}\int_{0\leq s\leq c}
    \int_{p\in\Lambda_{I(A)}}\indic_A\big( \psi (p,s)\big) \cr
  & \qquad\qquad
    \Big| {\d\over \d s} \psi (p,s)\Big| 
    \det \big( \psi_{s*}^\T \psi_{s*}^{\phantom T}\big)^{1/2} 
    f\big( \psi (p,s)/\lambda_A\big) \cr
  \noalign{\vskip 2pt}
  & \qquad\qquad
    \d\calM_{\Lambda_{I(A)}}(p) \, \d s \, \d c \, . \cr
  }
$$
\finetune{\hfuzz=0pt}
Consider a function $f$, bounded and continuous. After adding a
constant to $f$, we may assume that $f$ is larger than some
positive number. Then, up to
introducing a term $f\big( \psi (p,s)/\lambda_A\big)$, the
bounds in Theorem 3.1.10 remain valid.

Under (5.14) and arguing as in the proof of Theorem 5.1, as
$A$ tends to infinity, we have
$$\displaylines{\quad
    \int_A\eIx f(x/\lambda_A)\, \d x
    \sim 
    e^{-I(x)}
    \int_{p\in\Lambda_{I(A)}\hfill\atop\tau_A (p)<\infty}
    { e^{-\tau_A(p)}\over |\D I(p)|} \times
  \cut
    \int_{0\leq s\leq \chi_A(p)} e^{-s}
    f\big( \psi_{\tau_A(p)+s}(p)/\lambda_A\big)
    \, \d s\, \d\calM_{\Lambda_{I(A)}}(p) 
  \cr
  }
$$
where $\chi_A^F(p)\leq \chi_A(p)\leq \chi_A^L(p)$. Assumption
(5.14), (5.19) and Corollary 3.1.4 imply
$$
  {1\over \lambda_A} |\psi_{\tau_A(p)+s}(p)-p|
  \leq {\chi_A(p)\over \lambda_A |\D I(p)|}
  \leq {c_{A,M}\over \lambda_A |\D I(p)| }
  = o(1)
  \eqno{(5.21)}
$$
as $I(A)$ tends to infinity, uniformly in $p\in\Lambda_{I(A)}$ with
$\tau_A(p)<\infty$.

Assume that $f$ is uniformly continuous. Since $f$ is
bounded and larger than some positive number, (5.21) implies
$$
  \lim_{A\to\infty}\sup\Big\{ \, \Big|
    { f\big(\psi_{\tau_A(p)+s}(p)/\lambda_A\big)
      \over f(p/\lambda_A)}
    -1 \Big| 
  \, : \, p\in\underline{A}_M\, , \, 0\leq s\leq c_{A,M}
  \,\Big\} = 0 \, .
$$
Thus, using (5.5),
$$
  \int_A \eIx f(x/\lambda_A) \d x
  \sim e^{-I(A)}
  \int_{p\in\Lambda_{I(A)}\hfill\atop\tau_A(p)<\infty}
  {e^{-\tau_A(p)}\over |\D I(p)|}
  f(p/\lambda_A) \, \d\calM_{\Lambda_{I(A)}}(p) \, .
$$
We make a change of variable as we did in (3.2.2). Noting that
(5.20) implies
$$\displaylines{\qquad
  \lim_{A\to\infty}\sup\Big\{ \, \Big|
  { f\big(\exp_p(v)/\lambda_A\big)\over f(p/\lambda_A) }
  - 1\Big| \, : \, 
  v\in T_p\pi_A^{-1}(p)\, , \, 
  \cut
  \tau_A\big( \exp_p(v)\big) 
  \leq c_{A,M} \, , \, p\in\calD_A\, \Big\} = 0 \, ,
  \qquad}
$$
we obtain 
$$\displaylines{\qquad
  \int_A \eIx f(x/\lambda_A) \d x
  \sim e^{-I(A)} (2\pi )^{d-k-1\over 2} \times
  \cut
  \int_{\calD_A}
  { e^{\tau_A(p)}f(p/\lambda_A)\over |\D I(p)|^{d-k+1\over 2}
    \big( \det\, G_A(p)\big)^{1\over 2} }
  \, \d\calM_{\calD_A}(p) \, .
  \qquad}
$$
Make a change of variable $q=p/\lambda_A$ and use (5.18) to obtain
$$
  \lim_{A\to\infty} 
  {\int_A \eIx f(x/\lambda_A) \, \d x \over \int_A \eIx \d x}
  = \int f \d\nu \, ,
  \eqno{(5.22)}
$$
for all uniformly continuous, positive functions $f$.

If $f$ is bounded, we drop the restriction (5.14), as we did in
the proof of Theorem 5.1. Using Theorem 5.1, we can consider
arbitrary bounded uniformly continuous function $f$ in (5.22).
This implies --- see, e.g., Pollard (1984) --- that the 
conditional distribution converges weakly*.\hfill$\qed$

\bigskip

\centerline{\fsection Notes}
\section{5}{Notes}

\bigskip

\noindent There are many things related to Theorem 5.1 that I 
wanted to do but could not.

A first one is to understand to what extent an exponentially
integrable density may be approximated by a log-concave one at
infinity. Here is the beginning of what could be a proof. Let $f$
be a density on $\RR^d$, such that the moment generating
function 
$$
  \phi (t)=\int e^{\< t,x\>} f(x) \d x 
$$
is finite in a neighborhood of the origin. Under some classical
steepness conditions --- see, e.g., Barndorff-Nielsen (1978) or
Brown (1986) ---  the differential $m=\D\phi$ is a diffeomorphism.
Denote by $m^\inv$ its inverse. We now follow word for word the
construction of Barbe and Broniatowski (200?), but in a different
setting. 

The function $\log\phi$ is convex. Let $I$ be its convex
conjugate, that is
$$
  I(x)=\sup\{\, \log \phi (t) -\< t,x\>\,\} \, .
$$
Using the change of variable $a=s-x$ and Fubini's
theorem, we obtain
$$
  \int\int e^{\<m^\inv (s-x),x\>} \indic_A(s-x) f(s) \, 
  \d s \,\d x 
  = \int_A e^{-I(a)} \d a \, .
$$
This allows us to define a new density
$$
  g_A(x) = 
  { \int e^{\<m^\inv (s-x),x\>} \indic_A(s-x) f(s) \d s \over
    \int_A e^{-I(a)}\d a } \, .
$$
The interesting fact is that $g_A(0)=\int_A f(s) \d s \Big/
\int_A e^{-I(a)} \d a$. Consider the rescaled density
$r^dg_A(rx)$ with $r$ possibly depending on $A$. If we can prove
that $r^dg_A(r\,\cdot )$ converges to a limit, say $g(\cdot )$, 
as $A$ moves to infinity, in such a way that pointwise convergence 
at $0$ holds, then
$$
  \int_A f(s) ds \sim r^{-d}g(0)\int_A e^{-I(a)} \d a \, .
$$
Thus, when integrating over A, we can approximate the density 
$f$ by a multiple of
$e^{-I(A)}$, and then use Theorem 5.1.

To achieve this approximation, we can calculate the Fourier
transform of $r^dg_A(r\,\cdot )$. It is
$$
  \hat g_{r,A}(\lambda )
  = r^d\int e^{-\<\lambda ,x\>} g_A(rx) \, \d x
  = \int e^{-k(a,\lambda/r)}\, \d \mu_A(a) \, ,
$$
with
$$
  k(a,\lambda)=i\<a,\lambda\>+\log\phi\big( m^\inv (a)\big)
  -\log\phi\big( m^\inv(a)+i\lambda\big) \, ,
$$
and
$$
  \mu_A(B)=\int_{A\cap B}e^{-I(a)}\d a\Big/\int_A e^{-I(a)}\d a
  \, .
$$
In particular, as $A$ moves to infinity, the support of $\mu_A$
--- that is $A$ --- moves to infinity. A Taylor expansion of
$k(a,\lambda/r)$ near $\lambda /r =0$ gives
$$
  k(a,\lambda /r)
  = {1\over 2} \Big\< { (\D^2\log\phi )\circ m^\inv (a)\over r^2}
  \lambda ,\lambda \Big\> + o(\lambda^2) \, .
$$
If this can be done as $a$ tends to infinity --- which 
is in the spirit of what we did using Proposition 4.1.1 to 
prove Theorem 5.1 --- we can hope to approximate
$$
  \hat g_{r,A}(\lambda )\approx
  \int \exp\Big( {1\over 2}\Big\< 
  {(\D^2\log\phi)\circ m^\inv(a)\over r^2}\lambda,\lambda\Big\>
  \Big) \d\mu_A(a) \, .
$$
Let $V(a)=(\D^2\log\phi)\circ m^\inv(a)$ be the so-called variance
function of $f$. Inverting the Fourier transform of the
approximation, we should obtain
$$
  r^d g(0)\approx \int {1\over (2\pi )^{d/2}\det (V(a)/r)^{1/2} }
  \d\mu_A (a) \, .
$$
If we can find $r$ depending on $A$ such that the right hand side
has a limit as $A$ tends to infinity, and use Corollary 5.3, we are done.

Unfortunately, I could not come up with useful conditions for
this idea to work.

When $d=1$, it is possible to prove that if $g_A$ converges, then its
limit is given by a mixture of either normal densities --- as we
outlined here --- or gamma ones; this follows from Balkema,
Kl\"uppelberg and Resnick (1999). In higher dimensions, a related
approximation is in Barndorff-Nielsen and Kl\"uppelberg (1999). I
somewhat believe that the whole virtue of saddlepoint approximations
used in statistics is to provide some form of log-concave
approximation in the spirit of what is outlined here.  But most of the
time, in the multivariate setting, it relies on assumptions similar
to the convergence of $g_A$, which I don't find too appealing. I have
been searching unsuccessfully for a decent condition on $f$ itself.

A second project, which perhaps would be desirable to carry out,
is to obtain higher order expansions. Theorem 5.1 provides a one
term asymptotic expansion. Starting from the equality in
Proposition 3.1.9, one could do the change of variable using the
dominating manifold and the orthogonal leaves; then one would use
asymptotic expansions for whatever function is involved, and
obtain the desired approximation. Higher order differential
geometry is involved. The difficulty is to come up with a set of
usable conditions to perform all the approximations. Another
route would be to mimic the practice of Edgeworth expansions in
statistics. There are essentially two types of them: those that
are proved rigorously, and a vast majority that are called
``formal''. To do the formal ones, the argument is pretty much to
neglect what one believes to be negligible under some quite
unknown conditions, and proceed. One could obtain formal asymptotic
expansion in the same way. This may be of some value in a few
applications. Indeed, sometimes one may not look
for a theorem but maybe more for a guideline.

A third path to explore would be to derive explicit upper
bounds, starting either from Proposition 3.1.9 or Theorem 3.1.10.
In particular, I wonder if the technique developed here
could be of any use to investigate ``asymptotic'' isoperimetric
problems.

As pointed out in the notes to chapter 3, Proposition 3.1.9
does not use the convexity of $I$. This proposition is true for
any smooth function $I$ for which the sets $\Lambda_c$ are
smooth hypersurfaces. There may be some examples where the
normal flow can be calculated explicitly and other arguments
used in order to derive an estimate similar to that of Theorem
5.1.

To conclude these notes, Corollary 5.3 is inspired by the Gibbs
conditioning principle in large deviations. In the large
deviation context, the reader may consult Csisz\'ar (1984)
and Bolthausen (1993)

\vfil\eject

~\vskip 28pt
\chapter{6}{Asymptotics for sets translated towards infinity}
{\fchapter
  \line{\hfil 6.\quad Asymptotics for sets}
  \vskip .1in
  \line{\hfil translated towards infinity}}

\vskip 54pt

\noindent In this section, we study integrals of the form
$\int_{A+t}\eIx \d x$ as $|t|$ tends to infinity. To avoid
any ambiguity, recall that if $A$ is a set and $t$ is a vector,
both in $\RR^d$, the translation of $A$ by $t$ is
$$
  A+t=\{\, x+t \, :\, x\in A\,\} \, .
$$
We will assume that
\setbox1=\vtop{\hsize=3.5in\noindent%
  $A\subset\RR^d$ is a closed bounded convex 
  neighborhood of the origin, with smooth boundary and positive 
  curvature.}
$$
  \raise 14pt\box1\eqno{(6.1)}
$$
The only restriction here is convexity --- and smoothness, but we want
to be able to use differential geometric methods! It could be dropped
at the cost of a more sophisticated discussion on how $A+t$ and
$\Lambda_{I(A+t)}$ intersect. Up to changing $t$ by a fixed amount, we
can always assume that $A$ contains the origin.

We will control the growth of $I$ at infinity, assuming that
$$
  \lim_{p\to\infty} {\log I(p)\over |\D I(p)|}=0
  \qquad\hbox{ and }\qquad
  \lim_{p\to\infty} {\|\D^2I(p)\|\over |\D I(p)|}=0 \, .
  \eqno{(6.2)}
$$
The second condition forces $I$ not to increase too fast.
Indeed, for $d=1$, it reads $|I''(p)/I'(p)|$ tends to $0$ as
$|p|$ tends to infinity. Hence, for any small positive 
$\epsilon$ and any $p$, $q$ large enough,
$$
  \Big| {I'(p)\over I'(q)}\Big|
  = \Big|\exp\int_q^p {I''(t)\over I'(t)} \d t\Big|
  \leq \exp\big( \epsilon |p-q|\big) \, .
$$
For instance, the function $I(p)=\exp (|p|^\alpha)$ satisfies
$(I''/I')(p)$ tends to $0$ as $p$ tends to infinity if and only if
$\alpha <1$.  So, roughly, $I$ should have a subexponential growth. A
polynomial growth, like $I(p)=|p|^\alpha$ with a positive $\alpha$, is
admissible.

The first condition forces $I$ to increase fast enough. Indeed,
when $d=1$, it implies that for any positive $\epsilon$ 
and $x<y$ large enough,
$$
  I(y)-I(x)
  =\int_x^y I'(t)\d t
  \geq {1\over\epsilon} \int_x^y \log I(t) \d t
  \geq (y-x){\log I(y)\over \epsilon} \, .
$$
Hence, ultimately, $I$ has to grow faster that any linear
function. For instance, the function $I(t)=t(\log t)^\alpha$
satisfies $\big(\log I(t)\big)/I'(t)$ tends to $0$ as 
$t$ tends to infinity if and only if $\alpha > 1$. 

We will need to strengthen the second condition in (6.2) by
assuming
$$  
  \lim_{\ttoi}\log I(A+t)\sup\Big\{\, {\| \D^2I(q)\|\over |\D I(q)|}  :
  I(q)\geq I(A+t)\Big\} = 0 \, .
  \eqno{(6.3)}
$$

For a strictly convex function $I$ the level sets $\Gamma_c$ are
strictly convex. There is a unique point $x$ in $\RR^d$ at which
$I$ is minimal. Assumption (6.1) implies that for any $t$
with $|t|> 2\, \diam (A)+|x|$, 
$$
  (A+t)\cap \Lambda_{I(A+t)} = \{\, p_t\} 
$$
for a unique point $p_t$. In particular, $p_t-t$ is in $\partial A$.

\bigskip

\state{6.1. THEOREM.} {\it%
  If $I$ is convex and (6.1)--(6.3) hold, then
  $$
    \int_{A+t} \eIx \d x \sim
    {(2\pi )^{d-1\over 2}\over |\D I(p_t)|^{d+1\over 2} K_t}
    \, e^{-I(A+t)}
    \qquad\hbox{ as } t\to\infty \, ,
  $$
  where $K_t$\ntn{$K_t$} is the Gauss-Kronecker curvature of $\partial A$
  at $p_t-t$.}

\bigskip

\stateit{Proof.} Write $A_t=A+t$.
In order to apply Theorem 5.1, we
need to have a candidate for the dominating manifold
$\calD_{A_t}$.
Clearly $\{\, p_t\}$ should do, and we set $\dat=\{\,
p_t\}$. Since $\tau_{A+t}(p_t)=0$ by definition of $p_t$, the
result of Theorem 5.1 reads
$$
  \int_{A+t} \eIx \d x \sim
  { e^{-I(A+t)} (2\pi )^{d-1\over 2} \over
    |\D I(p_t)|^{d+1\over 2} \big( \det\, G_{A+t}(p_t)\big)^{1\over 2} }
  \qquad\hbox{ as } t\to\infty \, .
$$
Now, recall that we should expect $G_{A+t}(p_t)$ to be the difference
of the second fundamental forms of $\Lambda_{I(A+t)}$ and $A+t$
at $p_t$ --- not restricted to anything here, since $\dat$ is a
point and so $\pi_{A_t}^{-1}(p_t)$ is $\Lambda_{I(A+t)}$, up
to what is in the cut locus of $p_t$. However, the
second part of assumption (6.2) asserts that asymptotically, the
second fundamental form of $\Lambda_{I(A+t)}$ degenerates, and
so, locally, $\Lambda_{I(A+t)}$ is almost flat. Thus, $G_{A+t}(p_t)$
should be the second fundamental form of $\partial (A+t)$ at
$p_t$, which is equal to that of $\partial A$ at $p_t-t$. Its
determinant is exactly $K_t$. This explains how to guess the
result. It is hoped that this twelve line argument convinces
the reader that Theorem 5.1 can be useful.

Now that the result is guessed, let us find a candidate for
$c_{A+t,M}$. Define
$$
  c(t)=d\log I(A+t) + {d+2\over 2} \log |\D I(p_t)|
$$
--- in this formula, $d$ refers to the dimension of $\RR^d$.
Since $I(p_t)$ tends to infinity with $t$, the first part 
of (6.2) implies $\lim_{t\to\infty}|\D I(p_t)|=\infty$.
Proposition 2.1 yields
$$
  L\big( I(A+t)+c(t)\big) 
  = {e^{-I(p_t)} \over |\D I(p_t)|^{d+1\over 2}} o(1) 
  \qquad \hbox{ as }t\to\infty \, .
$$
Given (5.4) and our twelve line argument, 
$c_{A+t,M}=c(t)$ is a good candidate, no matter what $M$
is. It guarantees (5.4) as well as (5.6).

We now check all the assumptions of Theorem 5.1.

As noted in chapter 5, since $\dat=\{\, p_t\}$ is a point, (5.1) is
trivial. Notice that $\underline{A}_{t,M}$ is included in the
projection of $A+t$ on $\Lambda_{I(A+t)}$. So it is enough to check
(5.2) with $\underline{A}_{t,M}$ replaced by the projection of $A+t$.
Assumption (6.2) asserts that the second fundamental form of
$\Lambda_{I(A+t)}$ tends to $0$ uniformly over this surface.  Thus,
its curvature tensor vanishes asymptotically and the radius of
injectivity of any point in $\Lambda_{I(A+t)}$ tends to infinity
uniformly over the surface --- this follows from Rauch's (1951)
theorem or Klingenberg's (1959) lemma; see, e.g., Do Carmo (1992) or
Chavel (1996).  As $A+t$ stays of finite diameter, (5.2) follows.

To prove (5.3) and find $G_{A+t}$, we need to have some more
information on $\underline{A}_{t,M}$ and on the normal flow.  The idea
is that $\underline{A}_{t,M}$ should be very close to $\dat=\{\,
p_t\}$. To prove this fact, we first define a family of local
parameterizations of $\partial A$. For $p$ belonging to $\partial A$,
we denote by $\nu (\cdot )$\ntn{$\nu (p)$} the inward unit normal vector
to $\partial A$ at $p$.  By compactness of $A$, there exists a
positive $\epsilon_1$, independent of $p$, such that $\partial
A\cap\calB (p,\epsilon_1)$ can be parametrized as all points of the
form $p+u+f_p(u)\nu (p)$ for $u$ in $T_p\partial A$ and some
nonnegative function $f_p$.\ntn{$f_p(\cdot )$} Notice that the
curvature assumption (6.1) ensures that there exists a positive matrix
$Q_p$\ntn{$Q_p$} such that $f_p(u)={1\over 2} \< Q_pu,u\>\big(
1+o(1)\big)$ as $|u|$ tends to $0$. Moreover, since $\partial A$ is
smooth and compact, the term $o(1)$ is uniform when $p$ varies in
$\partial A$, and the matrices $Q_p$ are bounded bellow by a fixed
positive one.

To prove that $\underline{A}_{t,M}$ shrinks around $p_t$,
notice that $\nu (p_t-t)$ and $\D I(p_t)$ are collinear since
$\partial (A+t)$ and $\Lambda_{I(A+t)}$ are tangent at $p_t$.
Convexity of $I$ implies
$$\eqalign{
  & I\big(p_t+u+f_{p_t-t}(u)\nu(p_t-t)\big) - I(p_t) \cr
  & \qquad \geq {\d\over \d s} I\Big( p_t+s\big( u+ f_{p_t-t}(u)
    \nu (p_t-t)\big)\Big) \Big|_{s=0} \cr
  \noalign{\vskip 2pt}
  & \qquad = f_{p_t-t}(u) |\D I(p_t)| \, . \cr
  }$$
Consequently, the points $q$ in $\partial A+t$ such 
that $I(q)\leq I(A+t)+c_{A+t,M}$ can be parametrized 
as $p_t+u+f_{p_t}(u)\nu (p_t-t)$ with 
$|u|^2\leq O\big( c(t)/|\D I(p_t)|\big) = o(1)$. They can 
also be written as
$$
  p_t+u+{1\over 2} \<Q_{p_t-t}u ,u\>\nu (p_t-t)\big( 1+o(1)\big)
  \,
  \eqno{(6.4)}
$$
where the $o(1)$ is uniform in $t$ and $|u|^2\leq
O\big(c(t)/|\D I(p_t)|\big)$.

Since the curvature of the level set of $I$ tends to $0$, the
normal flow should be almost like straight lines on sizeable
intervals. In order to make this statement rigorous, and seeking
a linear approximation of the normal flow with good error
bounds, an elementary calculation shows that
$$
  \D\Big({\D I\over |\D I|^2}\Big)
  = {1\over |\D I|}\Big( {\D^2I\over |\D I|}- 2N\otimes N 
    {\D^2I\over |\D I|} \Big)
  = {1\over |\D I|}(\Id - 2N\otimes N){\D^2I\over |\D I|} \, .
$$
In particular, this implies the inequality
$$
  \Big\| \D \Big( {\D I\over |\D I|^2}\Big) \Big\|
  \leq {1\over |\D I|} {\| \D^2I\|\over |\D I|} \, .
  \eqno{(6.5)}
$$
Notice also that assumption (6.2) insures that
$$
  \eta (t)=\sup\Big\{\, { \| \D^2I(q)\|\over |\D I(q)|}
   : I(q)\geq I(p_t)\Big\} = o(1)
  \qquad \hbox{ as } t\to\infty \, .
$$
Using Lemma 3.1.1 and (6.5), it follows that with
$q=\exp_{p_t}(u)$,
$$\eqalignno{
  \Big| \psi (q,s)-q-s {\D I\over |\D I|^2}(q)\Big|
  & = \Big| \int_0^s\int_0^r {\d\over \d v} 
    {\D I\over |\D I|^2}
    \big( \psi (q,v)\big) \d v \, \d r\Big| \cr
  & \leq {\eta(t)\over |\D I(q)|} {s^2\over 2} \, .
    &(6.6)\cr
  }$$
Next, let us prove that $|\D I(q)|\sim |\D I(p_t)|$ as
$t$ tends to infinity provided $u$ stays bounded. 
Writing $\gamma (s)=\exp_{p_t}\big( su/|u|\big)$,
$$\eqalign{
  \log |\D I(q)|-\log |\D I(p_t)|
  & = \int_0^{|u|}{\d\over \d s} 
    \log \big| \D I\big( \gamma (s)\big)\big|\, \d s \cr
  & = \int_0^{|u|}{1\over 2} {\big\< \D^2I\big(\gamma (s)\big) 
    \gamma'(s) \, ,\D I\big( \gamma (s)\big) \big\> \over
    \big| \D I\big(\gamma (s)\big)\big|^2 } \,\d s \, , \cr
  }$$
from which we obtain the bound
$$\eqalignno{
  \big| \log |\D I(q)|-\log |\D I(p_t)|\big|
  & \leq {|u|\over 2} \eta (t)
    = {|u|\over 2}o(1) \qquad\hbox{ as } \ttoi \, . 
  &(6.7)\cr
  }$$

Furthermore, we have a good control on the
oscillations of $\D I/|\D I|^2$ in using (6.5); namely, 
for $t$ large enough and say $|u|\leq 1$,
$$\eqalign{
  \Big| {\D I\over |\D I|^2}(q) - {\D I\over |\D I|^2}(p_t)\Big|
  & = \int_0^{|u|} {\d\over \d s} {\D I\over |\D I|^2}
    \big( \gamma (s)\big) \d s \cr
  & \leq {2|u|\over |\D I(p_t)|} \eta (t) \, . \cr 
  }$$
Consequently, for $|u|\leq 1$, the inequality (6.6) gives the bound
$$
  \Big| \psi (q,s)-q-s{\D I\over |\D I|^2}(p_t)\Big|
  \leq {2\eta (t)\over |\D I(p_t)|} \big( s^2+s|u|\big) \, .
$$
This is the linear approximation of the normal flow that we were
looking for. Considering $s=\tau_{A+t}(q)$ and using the linear
approximation for the exponential map in Proposition A.2.1 ---
remember that $q=\exp_{p_t}(u)$ --- we then obtain
$$\displaylines{\qquad
    \psi\big( q, \tau_{A+t}(q)\big) 
    = p_t + u +\tau_{A+t}(q){N(p_t)\over |\D I(p_t)|}
  \cut
    + {\eta (t)\over |\D I(p_t)|}
    \big(\tau_{A+t}(q)^2+\tau_{A+t}(q)|u|\big)
    O(1) + \eta(t)|u|^2O(1)
    \qquad\cr}
$$
where the $O(1)$-terms are uniform in $|u|\leq 1$ as $|t|$ tends to
infinity. Since $\psi \big( q, \tau_{A+t}(q)\big)$ is in the boundary
of $A+t$ by the very definition of $\tau_{A+t}(q)$, and since $u$
belongs to $T_{p_t}\Lambda_{I(p_t)}$, (6.4) forces us to have
$$\displaylines{
    \tau_{A+t}(q){N(p_t)\over |\D I(p_t)|}
    + {\eta (t)\over |\D I(p_t)|}
    \big( \tau_{A+t}(q)^2+\tau_{A+t}(q)|u|\big)O(1)
    + \eta (t) |u|^2 O(1)
  \cut
    = {1\over 2} \< Q_{p_t-t}u,u\> \nu (p_t-t)
  \big( 1+o(1)\big) \quad\cr
  }$$
as $t$ tends to infinity, and uniformly in $|u|\leq 1$. 
Therefore, provided $\eta (t)\tau_{A+t}(q) = o(1)$, we obtain
$$
  \tau_{A+t}(q)={|\D I(p_t)|\over 2} \< Q_{p_t-t}u,u\>
  \big( 1+o(1)\big) 
$$
uniformly in $|u|\leq 1$, as $t$ tends to infinity.
Since $\eta (t)\tau_{A+t}(q)\leq \eta(t)c(t)$, assumption (6.3)
ensures that $\eta(t)\tau_{A+t}(q)=o(1)$, and we proved that
(5.3) holds with $G_{A+t}$ being the
second fundamental form of $\partial A$ at $p_t-t$.

Given our checking of (5.3), (5.5) is obvious since
$\underline{A}_{t,M}$ shrinks around $p_t$ --- see the proof that
$|u|^2=o(1)$ before equation (6.4) --- and (5.5) holds.

Assumption (5.7) is trivially satisfied. The shrinking of
$\underline{A}_{t,M}$ to $\{\, p_t\}$ 
and assumption (6.2) --- which
implies that the curvature tends to $0$; see also Proposition
4.2.3 --- imply that $t_{0,M}(p)= o(1)$ and $K_{\rm max}\big( p,
t_{0,M}(p)\big) = o(1)$ uniformly over $\underline{A}_{t,M}$ as
$t$ tends to infinity.

Assumption (5.8) is satisfied since $I$ is convex, $\calD_{A+t}$ is a
point, and the first part of Proposition 4.2.3 holds.

Assumption (5.9) follows from (6.7) and the shrinking of
$\underline{A}_{t,M}$ around $p_t$.

Clearly, (6.2) implies (5.10).

Assumption (5.11) is implied by (6.3) and Lemma 3.1.2, while (5.12)
holds systematically for $k=0$.

Since $\tau_{A+t}$ vanishes on $\calD_{A+t}$, (5.13) holds as well,
and this concludes the proof of Theorem 6.1. \hfill$\qed$

\bigskip

We obtained the conclusion of Theorem 6.1 by a brute application
of Theorem 5.1. A little extra work makes the asymptotic formula
nicer, replacing the term $|\D I(p_t)|$ by $|\D I(t)|$.

\bigskip

\state{6.2. COROLLARY.} {\it%
  Under the assumptions of Theorem 6.1,
  $$
    \int_{A+t}\eIx \d x \sim
  {(2\pi )^{(d-1)/2} \over |\D I(t)|^{(d+1)/2} K_t}
  e^{-I(A+t)} \qquad\hbox{ as } t\to\infty \, .
  $$
  }

\bigskip

\stateit{Proof.} Since 
$$
  {\d\over \d s} \log |\D I(t+su)|^2 
  = 2 {\< \D^2I(t+su)u\, , \D I(t+su)\> \over 
       | \D I(t+su)|^2 } \, , 
$$
the inequality
$$
  \big|\log |\D I(t+u)|^2-\log |\D I(t)|^2\big|
  \leq 2 \int_0^1 {\|\D^2I(t+su)\|\over |\D I(t+su)|}
  {|u|\over |\D I(t+su)|} \, \d s
$$
holds. Compactness of $A$, convexity of $I$ and (6.2) imply
that the right hand side of the above inequality is $o(1)$,
uniformly in $u$ belonging to $A$ as $t$ tends to infinity. Thus, 
$|\D I(t+u)|\sim |\D I(t)|$ uniformly over $u$ in $A$ as $t$
tends to infinity; and we can replace $|\D I(p_t)|$
by $|\D I(t)|$ in the statement of Theorem 6.1.\hfill$\qed$

\bigskip

In general, we do not have $I(A+t)-I(t)=o(1)$ as $t$ tends to
infinity.  This is easily seen when $d=1$ and $I(x)=x^2$ for
instance. Thus we cannot replace $I(A+t)$ by $I(t)$ in Theorem 6.1 or
Corollary 6.2. In some instances, it is possible to obtain an
asymptotic expansion for $I(A+t)$. We illustrate this fact in the
important case when $I$ is $\alpha$-positively homogeneous. To state
the result, recall the notation $N=\D I/|\D I|$ for the outward unit
normal vector field to the level lines of $I$, and set $\Pi=\D^2I/|\D
I|$.\ntn{$\Pi (\cdot )$} The compression of $\Pi$ to the tangent space
of a level line $\Lambda_c$ is its second fundamental form.

\bigskip

\state{6.3. PROPOSITION.} {\it%
  Assume that $I$ is $\alpha$-positively homogeneous and smooth.
  Assume also that $A$ is a neighborhood of the origin with
  a smooth boundary. Let $e$ be a unit vector in $\RR^d$.
  Define $r$\ntn{$r$} by the condition $-rN(e)\in\partial A$.
  Then, as $\lambda$ tends to infinity, $I(\lambda e+A)$ 
  admits an asymptotic expansion over the
  powers $\lambda^{\alpha -i}$, $i\in\NN$, and
  $$\displaylines{\qquad
      I(A+\lambda e) = \lambda^\alpha I(e) 
      - \lambda^{\alpha-1} r|\D I(e)| 
      + {\lambda^{\alpha-2}\over 2} r^2\<\Pi N(e)\, ,N(e)\>
    \cut
    \eqalign{
      &{} - {\lambda^{\alpha-3}\over 6} r^3\D^3I(e)
       \big( N(e)\, ,N(e)\, ,N(e)\big)\cr
      &{} -{\lambda^{\alpha-3}\over 2} 
       r^2 \< \Pi (e)Q^{-1}_{-rN(e)} \Pi (e)^\T N(e)\, ,N(e)\>
       +O(\lambda^{\alpha-4}) \, . \cr}
    \cr}
  $$}

\finetune{%\bigskip
}

\stateit{Proof.} Since $A$ is compact and $I$ is smooth 
and $\alpha$-positively homogeneous, we have, uniformly in
$u$ belonging to $A$ and as $\lambda$ tends to infinity,
$$\eqalign{
  I(\lambda e + u)
  &= \lambda^\alpha I(e+u/\lambda )\cr
  &= \lambda^\alpha I(e) +\sum_{1\leq i\leq k} 
    {\lambda^{\alpha-i}\over i!} \D^iI(e)
    \underbrace{(u\, ,\ldots ,u)}_{i\rm\;  times}
    + O(\lambda^{\alpha-k-1}) \, . \cr
  }$$
The expansion of $I(\lambda e+A)$ follows by induction. The
computation of the first terms can be done by introducing the point
$u_\lambda = u(\lambda ,e)$ in $A$, such that $I(e+A/\lambda )=
I(e+u_\lambda/\lambda)$. Since $A$ is convex, compact, and $I$
is convex, $u_\lambda$ in $\partial A$ for
$\lambda$ large enough.

Taylor's expansion gives, uniformly in $u$ belonging to $A$,
$$\displaylines{\qquad%
  I\Big( e+ {u\over\lambda} \Big) 
  = I(e) + {|\D I(e)|\over \lambda} \< N(e)\, ,u\> 
    + {|\D I(e)|\over 2\lambda^2} \< \Pi (e)u\, ,u\> 
  \cut
    {} + {1\over 6\lambda^3}\D^3I(e)(u\, ,u\, ,u) + O(\lambda^{-4})
  \qquad\cr
  }$$
as $\lambda$ tends to infinity. Consequently, $u_\lambda =
-rN(e)+u_{1,\lambda}/\lambda$ where $u_{1,\lambda}= O(1)$ as $\lambda$
tends to infinity --- because $u_\lambda$ has to minimize
$I(e+u/\lambda)$, and so should minize $\< N(e),u\>$ as well, up to a
term of order $\lambda^{-2}$.  Since $e+u_\lambda/\lambda$ belongs to
$\partial A$, using the notation of the proof of Theorem 6.1, there
exists a vector $v=v_\lambda$ in $T_{-rN(e)}\partial A=N(e)^\perp$
such that
$$\eqalign{
  {u_{1,\lambda}\over \lambda}
  & = {v\over\lambda} + f_{-rN(e)}\Big( {v\over \lambda}\Big)
    N(e) \cr
  & = {v\over\lambda} + {1\over 2\lambda^2}
    \< Q_{-rN(e)}v,v\> N(e) + O(\lambda^{-3}) \cr
  }$$
as $\lambda$ tends to infinity. It follows that
$$\displaylines{\qquad
    I\Big( e+{u_\lambda\over\lambda}\Big)
    = I(e) - {r|\D I(e)|\over \lambda} 
    + {r^2\over 2\lambda^2} |\D I(e)|\< \Pi (e)N(e)\, ,N(e)\>
  \cut
  \eqalign{
    &{} + {|\D I(e)|\over 2\lambda^3} \< Q_{-rN(e)}v\, ,v\>
      - {r\over \lambda^3}|\D I(e)|\< \Pi (e)v\, ,N(e)\> \cr
    &{} - {r^3\over 6\lambda^3} \D^3I(e) 
      \big( N(e)\, ,N(e)\, ,N(e)\big) + O(\lambda^{-4}) \, . \cr}
  \qquad\cr
  }$$
The term in $1/\lambda^3$ is smallest when $v= r
Q_{-rN(e)}^{-1}\Pi (e)^T N(e)$ --- we used that $Q$ is 
symmetric --- and its minimum value is
$$\displaylines{\qquad
    -{|\D I(e)|\over 2\lambda^3} r^2
    \< \Pi (e) Q^{-1}_{-rN(e)}\Pi(e)^TN(e)\, ,N(e)\>
  \cut
    {}- {r^3\over 6\lambda^2} \D^3I(e)
    \big( N(e)\, ,N(e)\, ,N(e)\big) \, .
  \qquad\cr}
$$
This completes the proof. \hfill$\qed$

\bigskip

In particular, for $\alpha=2$, we can replace $I(\lambda e +A)$
in the exponential term of the asymptotic equivalent by
$$
  \lambda^2I(e)-\lambda r |\D I(e)|+ {r^2\over 2}
  \< \Pi (e)N(e)\, ,N(e)\> \, .
$$

In the Gaussian setting, $I(x)=-|x|^2/2$. Thus, $|DI(e)|=1$
and $\Pi (c)=\Id$. The exponential term simplifies to
$$
  -{\lambda^2\over 2}-\lambda r + {r^2\over 2} \, .
$$
A neat expression, but very specific to the Gaussian
distribution$\dots$

\bigskip

Similarly to what we did in Corollary 5.2, we can obtain a result
on conditional distribution. It is easy to prove
that if $X$ is a random variable with density $e^{-I}$, then,
the conditional distribution of $X/|t|$ given $X\in A+t$ can be
approximated by a point mass at $t/|t|$, under the assumptions
of Theorem 6.1. However, Theorem 6.1 itself leads to a more
precise result.

\bigskip

\state{6.4. COROLLARY.} {\it%
  Let $X$ be a random variable with density proportional to
  $e^{-I}$. Let $e$ be a unit vector in $\RR^d$. Under the 
  assumptions of Theorem 6.1, the conditional distribution
  of $X-\lambda e$ given $X\in A+\lambda e$ converges
  weakly* to a point mass at $-rN(e)\in\partial A$ as $\lambda$
  tends to infinity.
  }

\bigskip

\stateit{Proof.} Let $U$ be a neighborhood of $-rN(e)$. We can
find a closed convex set $B$ in $U$ with smooth boundary and
positive curvature, such that $\partial B$ and $\partial A$
coincide in a neighborhood of $-rN(e)$. Applying Theorem 6.1
twice, we see that
$$
  \lim_{\lambda\to\infty}
  { \displaystyle\int_{B+\lambda e} \eIx \d x\over 
    \displaystyle\int_{A+\lambda e}\eIx \d x}
  = 1 \, .
$$
Consequently, the conditional distribution of $X-\lambda e$
given $X\in A+\lambda e$ is asymptotically concentrated on
$B\subset U$. Since $U$ is an arbitrary small neighborhood of
$-rN(e)$, the result follows.\hfill$\qed$

\bigskip

A slightly more involved proof would show that the conditional
distribution of $X-t$ given $X\in A+t$ can be approximated by a
point mass at $p_t-t$ as $|t|$ tends to infinity under the 
assumption of Theorem 6.1.

\bigskip

\centerline{\fsection Notes}
\section{6}{Notes}

\bigskip

\noindent This chapter has three motivations. First it provides a simple
example of applying of Theorem 5.1, and I hope it is of
pedagogical interest. Second, translating a set away from the
origin may be one of the most intuitive and natural ways to make
it moving to infinity. Third, and this is more important, the
Gaussian case has received some attention, due to statistical
applications. In LeCam's theory of local asymptotic normality
--- see e.g., LeCam (1986) and LeCam and Yang (1990) --- the
asymptotic power of a test is given by the probability that a
noncentered Gaussian vector lies in a given domain. Thus, one
issue is to calculate the probability that a centered Gaussian
vector hits a translated set, typically an ellipsoid. The work
of Breitung (1994), Breitung and Hohenbichler (1989),
Breitung and Richter (1996) are most relevant here. The remarks
following Proposition 6.3 somewhat enlighten the Gaussian case.

\vfil\eject

~\vskip 28pt
%\chapter{7}{Homothetic sets, {\sevencmssit I} homogeneous, Laplace's method}
\chapter{7}{Homothetic sets, homogeneous $I$ and Laplace's method}
{\fchapter
  \line{\hfil 7.\quad Homothetic sets,}
  \vskip .1in
  \line{\hfil homogeneous {\fchapterit I}}
  \vskip .1in
  \line{\hfil and Laplace's method}}

\vskip 54pt

\noindent In this chapter we consider a set $A_1$\ntn{$A_1$} such that
$$
  \hbox{there exists a neighborhood of $0$ not 
        intersecting $A_1$}
  \eqno{(7.1)}
$$

\noindent Equivalently, we could say that the complement of $A_1$ is a
neighborhood of the origin.
This assumption ensures that the sets $A_t=tA_1$ are moving
to infinity as $t$ tends to infinity. Assume furthermore that\ntn{$I(\cdot )$}
$$
  I \hbox{ is a strictly convex, $\alpha$-positively
  homogeneous function,}
  \eqno{(7.2)}
$$
that is $I(tx)=t^\alpha I(x)$ for all nonnegative $t$, all $x$ in
$\RR^d$, and some positive $\alpha$. Under such assumptions, $\alpha$
must be strictly larger than $1$ to ensure strict convexity. Setting
$x=ty$, we see that
$$
  \int_{A_t}\eIx \d x 
  = t^d\int_{A_1}e^{-I(ty)} \d y
  = t^d\int_{A_1} e^{-t^\alpha I(y)} \d y \, .
$$
The asymptotic decay of the last term in the equality is related to
the Laplace method. When $I$ has a unique minimum in $A_1$, not on the
boundary of $A_1$, this type of integral has been well studied.
However, here, $I(A_1)$ is achieved on the boundary of $A_1$,
eventually on a $k$-dimensional submanifold of $\RR^d$.  A direct
proof of an asymptotic equivalent of the right hand side, working out
a multivariate Laplace method, is quite tractable. However, for purely
pedagogical reasons, we will obtain an asymptotic equivalent of the
right hand side of the above equality by using Theorem 5.1. This proof
does not require more work than a direct one. The equality with the
right hand side makes it easy to understand how Theorem 5.1 works. It
also shows that Theorem 5.1 can be thought as a generalization of
Laplace's method.

\bigskip

Consider $\calD_{A_1}=A_1\cap\Lambda_{I(A_1)}$\ntn{$\calD_{A_1}$}
and assume that
\setbox1=\vtop{\hsize=3.5in\noindent%
  $\partial A_1$, $\calD_{A_1}$ and $\Lambda_c$ 
  are smooth --- twice continuously differentiable 
  --- manifolds.}
$$
  \raise 7pt\box1\eqno{(7.3)}
$$
Let $k$ be the dimension of $\da1$. We assume that $A_1$ separates
from $\Lambda_{I(A_1)}$ with contact of order $1$ exactly, and
therefore,
$$
  \det\, G_{A_1}(p)\neq 0 \qquad \hbox{ for any } p\in\da1 \, .
  \eqno{(7.4)}
$$
We also need to make sure that $A_1$ is not a thin $(d-1)$-dimensional
layer against $\Lambda_{A_1}$. For instance we could assume that it is
equal to the closure of its interior. Such an assumption is global.
We can work with a much weaker local one. Roughly speaking, for $p$ in
$\calD_{A_1}$, we need to be able to squeeze a ball in the
intersection of $A_1$ with the forward image of a leaf
$\pi^{-1}_{A_1}(p)$ through the normal flow. This guarantees some
thickness near $\calD_{A_1}$ along the section of $A_1$ orthogonal to
$\calD_{A_1}$. The exact assumption is that
\setbox1=\vtop{\hsize=3.5in\noindent%
there exists a positive $\epsilon$ such that for all 
$p$ in $\da1$ and any unit vector $v$ in $T_p\pi^{-1}_{A_1}$, 
any $s,\eta$ in $[\, 0,\epsilon \, ]$ the set $A_1$ contains
$\psi\big[\,\exp_p(\eta v)\, ,\tau_{A_1}\big( \exp_p (\eta v)\big)+s\, 
\big]$}
$$
  \raise 14pt\box1\eqno{(7.5)}
$$

The following is then a consequence of Theorem 5.1 and is
a multivariate Laplace type approximation.

\bigskip

\state{7.1. THEOREM.} {\it%
  Under (7.1)--(7.5), and if $\da1$ is
  a base manifold for $A_1$, then
  $$
    \int_{A_t} \eIx \d x \sim c_1 e^{-t^\alpha I(A_1)} 
    t^{k-(\alpha -2){\scriptstyle d-k\over\scriptstyle 2}
      -{\scriptstyle \alpha\over\scriptstyle 2}}
    \qquad\hbox{ as } t\to\infty \, ,
  $$
  where
  $$
    c_1 = (2\pi )^{d-k-1\over 2}
    \int_{\da1}
    {\d\calM_{\da1}\over 
     |\D I|^{d-k+1\over 2} 
    (\det\, G_{A_1})^{1\over 2}}
    \, .
  $$}

\bigskip

As in Theorem 5.1, the asymptotic equivalent in Theorem 7.1
must be read with $\det\, G_{A_1}=1$ if $\calD_{A_1}$ is an open
subset of $\Lambda_{I(A_1)}$ and $k=d-1$.

Before proving Theorem 7.1, notice first that for $\alpha =2$,
the polynomial term in $t$ in the approximation has exponent
$k-1$; it does not depend on the dimension $d$ of the
ambient space. More importantly, no matter what $d$ is, this
exponent can be written as $k{\alpha \over 2}+ d {2-\alpha\over
2} -{\alpha\over 2}$; it is an increasing function of $k$, as
one should expect.

The proof of Proposition 4.1.1 shows that
$G_{A_1}(p)=
\Pi_{\Lambda_{I(A_1),p}}^\pi-\Pi_{\partial A_1,p}^\pi$ is the
difference of the fundamental forms of $\Lambda_{I(A_1)}$ and
$\partial A_1$ compressed to the direction orthogonal to
$\calD_{A_1}$.

During the proof of Theorem 7.1, we will make use of the
following result, relating the large scale analysis of $tA_1$
to that of $A_1$ as far as the normal flow is concerned.

\bigskip

\state{7.2. LEMMA} {\it%
  If $I$ is $\alpha$-homogeneous, then

  \noindent (i) $\psi (p,s)=t^{-1} \psi (tp, t^\alpha s)$, and

  \noindent (ii) $\tau_{tA_1}(tp)=t^\alpha \tau_{A_1}(p)$.
  }

\bigskip

\stateit{Proof.} To prove (i), write $\tilde\psi
(p,s)=t^{-1}\psi (tp, t^\alpha s)$. Since $\D I$ is $(\alpha
-1)$-homogeneous, Lemma 3.1.1 yields
$$
  {\d\over \d s} \tilde\psi (p,s)
  = t^{\alpha -1} {\D I\over |\D I|^2} 
    \big( \psi (tp,t^\alpha s)\big) 
  ={\D I\over |\D I|^2} \big( \tilde\psi (p,s)\big) \, .
$$
Thus, $\tilde\psi$ obeys the differential equation of Lemma
3.1.1. It equals $\psi$ since $\tilde\psi (p,0)=p=\psi (p,0)$.

Assertion (ii) follows since $\tau_{A_1}(p)$ is the first positive
time $s$ such that $t\psi (p,s)$ is in $tA_1$, and 
$$
  \psi \big(tp,t^\alpha \tau_{A_1}(p)\big)
  = t\psi \big(p,\tau_{A_1}(p)\big) 
$$
thanks to assertion (i). \hfill$\qed$

\bigskip

\stateit{Proof of Theorem 7.1.} Since all the assumptions used
in Theorem 5.1 depend on $c_{A_t,M}$, we first need to guess
its value, and then proceed. To this aim, we first need an
estimate on the integral itself. We obtain it in evaluating the
asymptotic equivalent given by Theorem 5.1.

It is natural to consider
$$
  \calD_{A_t}=t\da1
  = t\, \{\, x\in\partial A_1  : I(x)=I(A_1)\,\} \, .
$$
Using the homogeneity of $I$,
$$
  \D I(tp) = t^{\alpha-1} \D I(p) \, .
$$
Let $p$ be in $\da1$, or equivalently, $tp$ be in $\calD_{A_t}$. 
To estimate $G_{A_t}(tp)$, the equality 
$(A_t,\Lambda_{I(A_t)})=(tA_1,t\Lambda_{I(A_1)})$ gives
$$
   G_{A_t}(tp)  = t^{-1}G_{A_1}(p) 
$$
--- the curvature tensors are rescaled by $1/t$; think of the sphere of
radius $t$ whose curvature is $1/t$ times that of a sphere of
radius $1$. We also have
$$
  \int f(p) \d\calM_{\calD_{A_t}}(p) 
  = t^k \int f(tq) \d\calM_{\da1}(q) \, .
$$
Consequently,
$$\displaylines{\qquad
  e^{-I(A_t)} (2\pi )^{d-k-1\over 2}
  \int_{\calD_{A_t}} 
  {\d\calM_{\calD_{A_t}}\over |\D I|^{d-k+1\over 2} 
   (\det \, G_{A_t})^{1\over 2}}
  \cut
  = c_1 e^{-t^\alpha I(A_1)} 
    t^{k-(\alpha-2) {\scriptstyle d-k\over\scriptstyle 2}
       -{\scriptstyle\alpha\over\scriptstyle 2}}
  \qquad\cr}
$$
where $c_1$ is given in the statement of Theorem 7.1. Notice
again that in a very few lines, Theorem 5.1 allows us to guess
the result.

As pointed out in chapter 5, the larger $c_{A_t,M}$ is, the
stronger the assumptions are. However, it is important to
remember that all that we need is to find $c_{A_t,M}$ larger than
$c_{A_t,M}^*$. So, consider a positive $\epsilon$ and let
$$
  c(t)=\Big( d\alpha -k+(\alpha -2) {d-k\over 2} 
  + {\alpha\over 2} +\epsilon \Big) \log t \, .
$$
From Proposition 2.1, we infer that
$$
  L\big( I(A_t)+c(t)\big) 
  = L\big(t^\alpha I(A_1)+ c(t)\big) 
  = o(e^{t^\alpha I(A_1)} t^{k-(\alpha -2)
      {\scriptstyle d-k\over\scriptstyle 2} - 
      {\scriptstyle\alpha\over\scriptstyle 2}}) 
  \eqno{(7.6)}
$$
as $t$ tends to infinity. Thus, $c_{A_t,M}^*$ is less than $c(t)$ for
$t$ large enough and any positive $M$. We can try to choose
$c_{A_t,M}$ to be $c(t)$.  In this case, we just proved that (5.4) is
satisfied.

It should be noticed that our choice of $c(t)$ is very naive. We
inverted asymptotically the function $L$ and evaluated the inverse at
the guessed asymptotic equivalent for $\int_{A_t}\eIx \d x$. The
addition of the term $\epsilon\log t$ in $c(t)$ is only to obtain
(7.6).

We now proceed in checking all the assumptions needed to apply Theorem
5.1. We already chose a candidate for the dominating manifold,
$$
  \calD_{A_t}=t\, \da1 = t\, \big\{\, x\in\partial A_1  :
  I(x)=I(A_1)\,\big\} \, .
$$
We postpone the check of (5.1)--(5.2) to the end of the
proof since it requires some discussion.

To describe $\underline{A}_{t,M}$, let
$$
  \tilde A_{1,t}=\big\{\, p\in\Lambda_{I(A_1)}  :
  \tau_{A_1}(p)\leq c(t)/t^\alpha \,\big\} \, .
$$
This choice ensures that
$$
  \underline{A}_{t,M} = 
  \big\{\, p\in\Lambda_{I(A_t)}  : \tau_{A_t}(p)\leq 
           c(t)\,\big\}
  = t\tilde A_{1,t} \, .
$$
Thus, we can look at $\underline{A}_{t,M}$ through its rescaled
version $\tilde A_{1,t}$.

We check (5.3) through a rescaling. Let $tq$ be 
in $\underline{A}_{t,M}$, and define $tp=\pi_{A_t}(tq)$. The
point $p$ is in $\da1$. Lemma 7.2 implies
$$
  \tau_{A_t}(tq)-\tau_{A_t}(tp)
  = t^\alpha \big( \tau_{A_1}(q)-\tau_{A_1}(p)\big) \, .
$$
Since $q$ belongs to $\tilde A_{1,t}$, we have $\tau_{A_1}(q)\leq
c(t)/t^\alpha$; in particular, $\tau_{A_1}(q)$ converges to $0$ as $t$
tends to infinity. Since $\tau_{A_1}$ is continuous on
$\pi_{A_1}^{-1}(p)$ with a strict minimum on $\da1$, we have
$|q-p|=o(1)$ uniformly in $q$ belonging to $\tilde{A}_{1,t}$ ---
recall $p=\pi_{A_1}(q)$. It then follows from Proposition 4.1.1 that
$$
  \tau_{A_1}(q)-\tau_{A_1}(p)
  = {1\over 2} \, |\D I(p)|\,  
    \< G_{A_1}(p)\exp_p^{-1} (q)\, ,\exp_p^{-1} (q)\>
    \big( 1+o(1)\big) 
$$
uniformly in $q\in\tilde{A}_{1,t}$. Consequently,
$$
  \tau_{A_t}(tq)-\tau_{A_t}(tp)
  = {t^\alpha\over 2\,} \, |\D I(p)| \, 
  \big\< G_{A_1}(p)\exp_p^{-1} (q)\, ,\exp_p^{-1} (q)\big\>
    \big( 1+o(1)\big)
$$
as $t$ tends to infinity, uniformly for $tq$
in $\underline{A}_{t,M}$. This proves (5.3) here.

Notice that since $(\exp_p)_*(0)=\Id$, we also have the approximation
$$
  \tau_{A_1}(q)-\tau_{A_1}(p)
  = {1\over 2} \, |\D I(p)| \,
  \< G_{A_1}(p)(q-p) \, , (q-p)\> \big( 1+o(1)\big) \, ,
$$
as $q$ converges to $p$, which may look more familiar.

To verify (5.5), Lemma 7.1.2 shows that
\smash{$\chi_{A_t}^F(tp)=t^\alpha\chi_{A_1}^F(p)$}. Thus, it suffices to
prove that \smash{$\inf_{p\in\tilde{A}_{1,t}}\chi_{A_1}^F(p)$} is
uniformly bounded below by some positive number for $t$ large
enough. If $q$ belongs to $\tilde{A}_{1,t}$, write
\smash{$q=\exp_{\pi_{A_1}(q)}(\eta v)$} for
$\eta={\rm dist}\big(q,\pi_{A_1}(q)\big)$ and a unit vector $v$.  As
we have seen, $\eta$ converges to $0$ uniformly over $q$ in
\smash{$\tilde{A}_{1,t}$} as $t$ tends to infinity. Therefore, (7.5)
implies that for $t$ large enough,
\smash{$\chi_{A_1}^F(\cdot )\geq\epsilon$} over \smash{$\tilde{A}_{1,t}$}. 
Thus, (5.5) holds.

Our choice of $c_{A_t,M}= c(t)$ ensures that (5.6) holds as well.

To check (5.7) is not much more complicated. Define
$$
  c= \sup\Big\{\,  { \lambda_{\rm max}\big( \D^2I(q)\big)^2\over 
	       \big| \D I(q)\big|^2 }  : 
  q\in \Lambda_{I(A_1)}\,\Big\} \, .
$$
We first notice
that for any $p$ in $\Lambda_{I(A_1)}$ and $s$ positive, the definition of
$K_{\rm max}$ in Proposition 4.2.4 and homogeneity of $I$ imply
$$
  K_{\rm max}(tp,s) \leq
  \sup\Big\{\, { \lambda_{\rm max}\big( \D^2I(tq)\big)^2\over 
	       |\D I(tq)|^2 }  : 
  q\in \Lambda_{I(A_1)}\,\Big\}
  = {c\over t^2} \, .
$$
Identifying $T_p\Lambda_{I(p)}$ and $T_{tp}\Lambda_{I(tp)}$, the
equality $\exp_{tp}(w)=t\exp_p (w/t)$ holds. It gives,
$$\eqalign{
  t_{0,M}(tp)
  & = \sup\Big\{\, s : 
    \inf_{v\in S_{T_{tp}\pi_{A_t}^{-1}(tp)}(0,1)}
    \tau_{tA_1}\big( \exp_{tp}(sv)\big) \leq c(t)\,\Big\}  \cr
  & = t\sup\Big\{\, s : 
    \inf_{v\in S_{T_{tp}\pi_{A_t}^{-1}(tp)}(0,1)}
    \tau_{A_1}\big( \exp_p(sv)\big) \leq c(t)/t^\alpha\,\Big\}  
    \, . \cr
  }$$
Again, as $c(t)/t^\alpha$ converges to $0$ as $t$ tends to
infinity, the requirement
$\tau_{A_1}\big( \exp_p(sv)\big) \leq c(t)/t^\alpha$ forces $sv$
to be $o(1)$ as $t$ tends to infinity. Proposition 4.1.1 yields
\finetune{\hfuzz=3pt}
$$
  t_{0,M}(tp) \leq
  t\sup\Big\{\, s : \inf_{v\in S_{T_p\pi_{A_1}^{-1}(p)}(0,1)}
  {1\over 2} |\D I(p)| \< G_{A_1}(p)sv\, , sv\> 
  \leq {2c(t)\over t^\alpha} \,\Big\}
$$
\finetune{\hfuzz=0pt}
for $t$ large enough.
Since $G_{A_1}(p)$ does not have null eigenvalues on
$T_p\pi^{-1}_{A_1}(p)$ thanks to (7.2), $t_{0,M}(tp)$
is at most
$$\displaylines{\qquad
    t\sup\Big\{\, s : s^2
    \inf_{v\in S_{T_p\pi^{-1}_{A_1}(p)}(0,1)}
    {1\over 2} |\D I(p)| \< G_{A_1}(p)v\, , v\> 
    \leq {2c(t)\over t^\alpha} \,\Big\}
  \cut
    = t^{1-\alpha/2}\sqrt{c(t)}\, O(1)
  \qquad\cr}
$$
as $t$ tends to infinity. Consequently,
$$
  \sup_{tp\in\underline{A}_{t,M}}
  \sqrt{ K_{\rm max}\big(tp,t_{0,M}(tp)\big) }
  t_{0,M}(tp) 
  = t^{-\alpha/2} \sqrt{ c(t) } \, O(1) 
  = o(1)
$$
as $t$ tends to infinity, and (5.7) holds.

Following Remark 5.2, we will not check (5.8) but (5.16)
instead. Adding subscripts to distinguish on which manifold we
are working, we have
$$
  \calA_{\pi^{-1}_{A_t}(tp)}(s,v)
  =t\, \calA_{\pi^{-1}_{A_1}(p)}(s/t,v) \, ,
$$
for all $s$ in $\big[\,0,e_A(tp,v)\big)$ and all
$v$ in the unit sphere of \smash{$T_{tp}\pi^{-1}_{A_t}(tp)$}. This
unit sphere can be identified with that of \smash{$T_p\pi^{-1}_{A_1}(p)$}. 
As we have seen, if
$\exp_{tp}(sv)$ belongs to $\underline{A}_{t,M}$, then
$t^{1-\alpha/2}s/t^\eta= o(1)$ for any positive $\eta$, as 
$t$ tends to infinity, and uniformly in $p$ belonging 
to $\da1$. It follows from compactness of $\da1$
and the classical expansion for $\calA_{\pi^{-1}_{A_1}}(s,v)$ 
as $s$ tends to $0$ --- see, e.g., Chavel (1996) ----
that --- recall $\alpha>1$! ---
$$
  \calA_{\pi^{-1}_{A_t}(tp)}(s,v)
  = s\, \Id_{T_{tp}\pi^{-1}_{A_t}(tp)} + O(s^3/t^2)
  = s\, \Id_{T_{tp}\pi^{-1}_{A_t}(tp)} + o(1)
  \eqno{(7.7)}
$$
as $t$ tends to infinity, and uniformly for $p$ in $\da1$. 
This implies (5.16).

Assumption (5.9) is easy to check. Since $p$ belongs to $\da1$ and
$q$ to $\tilde{A}_{1,t}\cap \pi^{-1}_{A_1}(p)$, we have
$$
  { |\D I(tp)| \over |\D I(tq)| }
  = { |\D I(p)| \over |\D I(q)| }
  = { |\D I(tp)| \over \big|\D I\big(p+o(1)\big)\big| }
$$
where the $o(1)$-term is uniform in $q$ in the given range and
$p$ in $\da1$.

Assumption (5.10) is trivial since 
$$
  { \|\D^2I(tp)\| \over |\D I(tp)|^2 }
  = t^{-\alpha} { \|\D^2I(p)\| \over |\D I(p)|^2 }
$$
by homogeneity of $I$.

Since $\tau_{A_1}(\cdot )$ vanishes on $\da1$ here, (5.11) holds as
well as (5.13).

We check (5.12) by rescaling. Indeed, for $q$ in $\tilde{A}_{1,t}$,
or equivalently, $tq$ in $\underline{A}_{t,M}$, we have
$$
  \pi_{A_t}(tq)=t\pi_{A_1}(q) \, .
$$
Identifying $T_{tq}\underline{A}_{t,M}$ and
$T_q\tilde{A}_{1,t}$, we also have 
$$
  \pi_{A_t,*}(tq)(v)= t\pi_{A_1,*}(q)(v/t) \, .
$$ 
But $\big|q-\pi_{A_1}(q)\big|$ tends to $0$ uniformly 
in $\tilde{A}_{1,t}$ as $t$ tends to infinity. Since
the differential $\pi_{A_1,*}\big( \pi_{A_1}(q)\big)$ 
is the orthogonal projection onto ${T_p}\da1$
--- this comes form the fact that $\pi_{A_1}(p)=p$ for all
$p$ in $\da1$, and that $\pi^{-1}_{A_1}(p)$ is orthogonal to
$\da1$ --- (5.12) follows.

We are left to check that in (5.1), $\calD_{A_t}$ is indeed a
dominating manifold for the set
$A_t\cap\Gamma_{I(A)+c_{A_t},M}$ and that (5.2) holds. 
If $k=0$, i.e., $\da1$ is made
of a finite number of points, (5.1) is clear.
Assumption (5.2) is then checked by rescaling, using the fact
that $\tilde{A}_{1,t}$ shrinks around $\calD_{A_1}$. 
Assume $k\geq 1$. Define
\finetune{\hfuzz=2.1pt}
$$
  A_t^{(k)} = A_t\cap 
  \big\{\, \psi \big(\exp_x(u),s\big) : x\in\calD_{A_t} \, ,\, 
  u\in T_p\Lambda_{I(A_t)}\ominus T_p\calD_{A_t}\, ,\, s\geq 0
  \,\big\} \, .
$$
\finetune{\hfuzz=0pt}
Then, $\calD_{A_t}$ is a dominating manifold for $A_t^{(k)}=
tA_1^{(k)}$. From what we have done, we can apply Theorem
5.1 to obtain an asymptotic approximation for
$\int_{A_t^{(k)}}\eIx \d x$, and the result is nothing but the
statement of Theorem 7.1. Instead of checking (5.1), it now
suffices to prove that
$$
  \int_{A_t\setminus A_t^{(k)}} \eIx \d x
  = o\Big( \int_{A_{t^{(k)}}}\eIx \d x\Big)
  \qquad \hbox{ as } t\to\infty \, .
$$
To apply Theorem 5.1 to evaluate the integral in the left hand
side of the above inequality, we need to find a
dominating manifold for $A_t\setminus A_t^{(k)}$. Pick
$\calD_{A_1\setminus A_1^{(k)}}$ to be $\partial \calD_{A_1}$. 
It is a dominating manifold, of dimension at most $k-1$ for 
the set
$$
  A_t^{(k-1)}= \big( A_t\setminus A_t^{(k)} \big)
  \bigcap \Big\{\, \psi (x,s) : 
  x\in t\calD_{A_1\setminus A_1^{(k)}} \,\Big\} \, .
$$
Now, it is possible that $A_t\setminus A_t^{(k)}$ or
$A_t^{(k-1)}$ do not have smooth boundaries. But since $\partial A_1$
is smooth, the set \smash{$A_t^{(k-1)}$} can be included in a set
$A_t^{(k-1)\prime}$ say, with smooth boundary, and for which
$t\calD_{A_1\setminus A_1^{(k)}}$ is again a dominating manifold.
This is done in parameterizing $\tau_{A_t^{(k-1)}}(p)$ for 
$p$ in $t\calD_{A_1\setminus A_1^{(k)}}$ so that
$\tau_{A_t^{(k-1)\prime}}$ and $\tau_{A_t^{(k-1)}}$ coincide on the
projection of $A_t^{(k-1)}$ on $\Lambda_{I(A_t)}$ for instance, and extending
$\tau_{A_t^{(k-1)}}$ into a differentiable function on
$\Lambda_{I(A_t)}$. This way, we need to estimate the integral
over \smash{$A_t^{(k-1)\prime}$}, for which we have a $(k-1)$-dimensional
candidate for a dominating manifold. Iterating this process, we go
down to a $0$ dimensional dominating manifold, apply Theorem 5.1
in this case, and obtain the order of all the terms with
dominating manifold of dimension between $0$ and $k-1$. In
particular, this implies
$$\eqalign{
  \int_{A_t\setminus A_t^{(k)}}\eIx \d x
  & = O(1)e^{-I(A_t)} t^{(k-1)-(\alpha -2)
    {\scriptstyle d-(k-1)\over\scriptstyle 2} -
    {\scriptstyle \alpha\over\scriptstyle 2} } \cr
  & = o\Big( \int_{A_t^{(k)}} \eIx \d x \Big) \, ,
  }$$
and this proves Theorem 7.1. \hfill\qed

\bigskip

\state{7.3. REMARK.} It is essential to notice the following. In
checking (5.3), (5.9), (5.12), the only feature we used 
besides homogeneity of $I$ is that
a point $q$ in $\tilde{A}_{1,t}$ converges to $\pi_{A_1}(q)$ as
$t$ tends to infinity, uniformly in $\tilde{A}_{1,t}$. For (5.7) and
(5.8) --- actually (5.16) --- we used slightly more, namely that
$\big| q-\pi_{A_1}(q)\big|= o(t^{-1/3})$ in order to obtain
(7.7). The conclusion is that if $I$ is homogeneous, all that
we need to do to check the assumptions is to check that
$\underline{A}_{t,M}/t$ concentrates to $\calD_{A_t}/t$ at the
rate $o(t^{-1/3})$. So, the work done in the proof of Theorem
7.1 may save us some effort in other applications.

\bigskip

\state{7.4. REMARK.} In many applications, the set $A_1$ is of the
form $A_1=\big\{\, x : g(x)\geq 0 \,\big\}$ for some smooth --- twice
continuously dif\-fer\-en\-ti\-able --- function $g$, and $I(A)$ is
achieved at boundary points $x$ such that $g(x)=0$. In such cases,
assumption (7.5) can be simplified into an analytical condition. To
see this, define $q=\exp_p(\eta v)$ and $r=\psi\big(q,\tau_{A_1}(q)\big)$. 
Since
$$
  \psi\big( q,\tau_{A_1}(q)+s\big)
  = r + s {\D I\over |\D I|^2}(r) + O(s^2)
$$
thanks to Lemma 3.1.1, we have
$$\eqalign{
  g\Big(\psi\big( q,\tau_{A_1}(q)+s\big)\Big) 
  & = g(r) + s \Big\< \D g(r),{\D I\over |\D I|^2}(r)\Big\>
    + O(s^2) \cr
  & = s \Big\< \D g,{\D I\over |\D I|^2}\Big\>(r) 
    + O(s^2) \, . \cr
  }$$
For this expression to be nonnegative for $s$ nonnegative, it is
enough to have
$$
  \< \D g(r),\D I(r)\> > 0 \, .
$$
Since $I$ and $g$ are smooth and $q$ is close to $p$ for small
$\eta$, the condition
$$
  \inf_{p\in\calD_{A_1}} \< \D g(p),\D I(p)\> > 0
$$
is sufficient for (7.5) to hold.  Since a point $p$ in $\da1$
minimizes $I$ on $A_1$, the vectors $\D g(p)$ and $\D I(p)$ are
positively proportional. The condition $\< \D g(p),\D I(p)\>$ positive
holds provided $\D g(p)$ is nonzero --- notice that $\D I$ does not
vanish for $I$ is strictly convex.  In other words, (7.5) is fulfilled
as soon as $g$ has no critical points on $\calD_{A_1}$.

\bigskip

Applying Corollary 5.2, we can obtain results on conditional
distributions.  The following statement asserts that the probability
measure with density propositional to $\indic_{A_1}(x) e^{-I(tx)}$
converges weakly* as $t$ tends to infinity to the probability measure
absolutely continuous with respect to $\calM_{\da1}$, and density
proportional to $|\D I|^{-(d-k+1)/2}(\det\, G_{A_1})^{-1/2}$.

\bigskip

\state{7.5. THEOREM.} {\it%
  Under the assumptions of Theorem 7.1, if $B$ is a 
  set of continuity of $\calM_{\da1}$, then
  $$
    \lim_{t\to\infty} 
    { \displaystyle\int_{t(B\cap A_1)}\eIx \d x \over 
      \displaystyle\int_{tA_1}\eIx \d x }
    =
    { \displaystyle \int_{B\cap\da1} 
      { \displaystyle \d\calM_{\da1}\over
        \displaystyle |\D I|^{(d-k+1)/2} 
        (\det\, G_{A_1})^{1/2} } \over
      \displaystyle \int_{\da1} 
      { \displaystyle \d\calM_{\da1}\over 
        \displaystyle |\D I|^{(d-k+1)/2} 
        (\det\, G_{A_1})^{1/2} } } \, .
  $$}

\bigskip

\stateit{Proof.} We just need to check (5.18)--(5.20).
Here, we consider $\lambda_{A_t}=t$. The measure in (5.18)
rcan be rewritten as
$$
  {\d\calM_{\da1}(q)\over 
   |\D I|^{d-k+1\over 2}\big( \det\, G_{A_1}(q)\big)^{1\over 2} }
  \; \Big/
  \int_{A_1} {\d\calM_{\da1}(q)\over 
   |\D I|^{d-k+1\over 2}\big( \det\, G_{A_1}(q)\big)^{1\over 2} }
$$
and does not depend on $t$.

Since
$$
  {c_{A_t,M}\over \lambda_{A_t}|\D I(tp)|}
  \leq 2 d \alpha {\log t\over t^\alpha |\D I(p)|} \, ,
$$
assumption (5.19) is satisfied.

To check (5.20), use rescaling to obtain
$$\displaylines{\qquad
    \sup\Big\{ \, {|tq-tp|\over \lambda_{A_t}} \, : \, 
    tq\in\pi^{-1}_{A_t}(tp) \, , \, \tau_{A_t}(tq)\leq c(t) 
    \, \Big\}
  \cut 
    = \sup\big\{ |q-p|\, :\, q\in\pi_{A_1}^{-1}(p)\, , \,
    \tau_{A_1}(q)\leq c(t)/t^\alpha \, \} \, . 
  \qquad\cr
  }
$$
Since $\tilde{A}_{1,t}$ shrinks around $\da1$, assumption
(5.20) holds true. Applying Corollary 5.2 yields the conclusion. 
\hfill$\qed$

\bigskip

Another way to formulate Theorem 7.5 is to say that the probability
measures with density proportional to $\indic_{A_1}(x)e^{-I(tx)}$
converge weakly* to the one with density proportional to
$|DI|^{(d-k+1)/2}(\det\, G_{A_1})^{-1/2}$ with respect to the
Riemannian measure on $\calD_{A_1}$.

\bigskip

\centerline{\fsection Notes}
\section{7}{Notes}

\bigskip

\noindent The notes of chapter 1 contain references on Laplace's method.
Also very much related to this chapter is the work of
Breitung (1994) in a Gaussian setting. Theorem 5.1 is related to
Hwang (1980). The Laplace method in dimension larger than one
with a dominating manifold of minimizing points is developed in
Barbe and Broniatowski (200?), motivated by large deviation theory.

\vfil\eject

~\vskip 28pt
\chapter{8}{Quadratic forms of random vectors}
{\fchapter
  \line{\hfil 8.\quad Quadratic forms}
  \vskip .1in
  \line{\hfil of random vectors}}

\vskip 54pt

\noindent In this chapter, we illustrate the use of Theorems 5.1
and 7.1 to deal with the following question. Consider a random
vector $X=(X_1,\ldots , X_d)$ in $\RR^d$, and a real $d\times d$ 
matrix\ntn{$C$} $C={(C_{i,j})}_{1\leq i,j\leq d}$. What is the decay of 
$P\{\, \< CX,X\>\geq t\,\}$ as $t$ tends to infinity?

Of course, this decay depends on the distribution of $X$ as well
as on the matrix $C$. We will deal with two different types of
distributions: symmetric Weibull- and Student-like. The Weibull-like
tail will be handled through application of Theorem 7.1, while
the Student-like one will be handled by a change of variable
technique and Theorem 5.1.

\bigskip

\noindent{\fsection 8.1. An example with light tail
distribution.}
%\chapternumber={8.1}
\section{8.1}{An example with light tail}

\bigskip

\noindent Consider a random vector $X=(X_1,\ldots
,X_d)$ in $\RR^d$, having density $f=e^{-I}$ for some convex
function $I$ on $\RR^d$. Assume moreover that $I$ is
$\alpha$-positively homogeneous. Writing\ntn{$A_t$}
$$
  A_t=\big\{\, x\in\RR^d  : \<Cx,x\>\geq t\,\big\}
  = \sqrt t A_1 \, ,
$$
we see that
$$
  P\big\{\, \< CX,X\>\geq t\,\big\} 
  = \int_{A_t}\eIx \d x \, ,
$$
and Theorem 7.1 is relevant here. 

To be more specific, assume that the components $X_i$ of $X$ are
independent and identically distributed, all with 
density\ntn{$w_\alpha (\cdot )$}
$$
  w_\alpha (x) = {\alpha^{1-(1/\alpha)}\over 2\Gamma (1/\alpha)}
  \exp \Big( {-|x|^\alpha\over \alpha}\Big) \, ,
  \qquad x\in\RR \, , \quad \alpha > 1 \, . 
$$
To apply Theorem 7.1 we need to describe the points of
$\partial A_1$ at which\ntn{$I(\cdot )$}
$$
  I(x)={1\over\alpha}\sum_{1\leq i\leq d}|x_i|^\alpha
$$
is minimum. Surprisingly, this problem seems quite difficult,
and I have not been able to solve it in general. The result will
rely on the following conjecture.

\bigskip

\state{8.1.1. CONJECTURE.} {\it If $\alpha\neq 2$ and
  $C+C^\T$ has no vanishing eigenvalue, then 
  $\sum_{1\leq i\leq d}|x_i|^\alpha$ admits a finite number
  of minima in $A_1$; moreover $\det\, G_{A_1}$ is not null
  at these minima.
  }

\bigskip

Hence, if this conjecture is indeed true, $\calD_{A_1}$ is of
dimension $k=0$ when $\alpha\neq 2$ and $C+C^\T$ has no
degeneracy.

The application of Theorem 7.1 is then trivial. For 
$\alpha\neq 2$, we obtain
$$
  P\{\, \< CX,X\> \geq t\,\} 
  \sim {\alpha^{d-(d/\alpha)}\over 2^d\Gamma (1/\alpha)^d}
  c_1e^{-t^{\alpha/2}I(A_1)} 
  t^{-(\alpha-2){d\over 4} -{\alpha\over 4}} 
$$
as $t$ tends to infinity --- recall that $A_t=\sqrt tA_1$ here and not
$tA_1$ as in chapter 6 --- where
$$
  c_1= (2\pi )^{d-1\over 2} \sum_{x\in\calD_{A_1}}
  | \D I(x)|^{-(d+1)/2}(\det\, G_{A_1})^{-1/2} \, .
$$
The term in $c_1$ can be made more explicit. Indeed,
$$
  \D I(x)={\big( \sign (x_i)|x_i|^{\alpha -1}\big)}_{1\leq i\leq d}
$$
for $\alpha >1$. Moreover, the remark following the statement of
Theorem 7.1 asserts that $G_{A_1}$ is obtained by the
difference of the fundamental form of $\Lambda_{I(A_1)}$ and
$\partial A_1$. Since $\D^2I(x)=(\alpha -1)\diag \big(
|x_i|^{\alpha-2}\big)$, we have
$$
  \Pi_{\Lambda_{I(A_1)}}(x)
  = \Proj_{T_x\Lambda_{I(A_1)}} {\alpha -1\over |\D I(x)|} 
    \diag \big( |x_i|^{\alpha -2} \big)
    \Big|_{T_x\Lambda_{I(A_1)}}
$$
while
$$
  \Pi_{\partial A_1}(x) 
  = \Proj_{T_x\partial A_1} {C+C^\T\over |(C+C^\T)x|}
    \bigg|_{T_x\partial A_1} \, .
$$

For $\alpha =2$, the calculation can be done explicitly.
Let $\lambda$ be the largest eigenvalue of $C+C^\T$, and assume
that $\lambda$ is positive --- otherwise $\< Cx,x\>$ is nonpositive for any $x$
and $P\{\, \< CX,X\>\geq t\,\}$ is null for any
positive $t$. Let\ntn{$H$}
$$
  H=\{\, x  : (C+C^\T)x=\lambda x\,\}
$$
be the eigenspace associated to the largest eigenvalue $\lambda$.

\bigskip

\state{8.1.2. THEOREM.} {\it%
  Let $k$ be $\dim H-1$ and $M$ be 
  the compression of $\,\Id - {\lambda^{-1}(C+C^\T)}$ to $H^\perp$.
  For $\alpha=2$ and $\lambda$ positive,
  $$
    P\{\, \cXX\geq t\,\} \sim 
    {1\over\lambda^{(k-1)/2} \Gamma\big( {k+1\over 2}\big)
     (\det M)^{1/2}} \, e^{-t/\lambda}\, t^{(k-1)/2}
  $$
  as $t$ tends to infinity. 
  }

\bigskip

\stateit{Proof.} Assumptions (7.1), (7.2) and
(7.3) hold. To check (7.4), we only need to calculate
$\calD_{A_1}$. We claim that $\calD_{A_1}$
is the sphere of radius $\sqrt{2/\lambda}$ centered at the origin 
in $H$. Indeed, if $\cxx = 1$ then
$$
  2= \< \cct x,x\> \leq \lambda |x|^2 \, .
$$
So $|x|\geq \sqrt{2/\lambda}$. On the other hand, if $\cct
x=\lambda x$ and $|x|^2=2/\lambda$, then $\< \cct x,x\>=\lambda
|x|^2 = 2$, and then $\cxx =1$. 

Applying Theorem 7.1, we obtain,
$$
  \int_{A_t} e^{-|x|^2/2} \d x \sim
  c_1 e^{-tI(A_1)}t^{k-1\over 2} 
  \qquad \hbox{ as } t\to\infty \, ,
$$
and $I(A_1)=1/\lambda$ from the preceding argument. To calculate the
constant $c_1$, notice that in our case, $\D I=\Id$.
Thus, $|\D I(x)|=\sqrt{2/\lambda}$ on $\da1$.

To calculate $G_{A_1}$, observe that $T_x\Lambda_{I(x)}=\{\,
x\,\}^\perp$, for the level lines of $I$ are spheres. Thus, for
$x$ in $\da1$, the second fundamental form of $\Lambda_{I(A_1)}$ at
$x$ is
$$
  \Pi_{\Lambda_{I(A_1)}}(x) 
  = {\Proj_{\{ x\}^\perp}\Id\big|\rx\over |\D I(x)|}
  = \sqrt{\lambda\over 2}\, \Proj_{\{ x\}^\perp}\Id\big|\rx \, .
$$
On the other hand, the second fundamental form of $\partial A_1$
at some $x$ in $\calD_{A_1}$ is
$$
  \Pi_{\partial A_1}(x)
  = {\Proj_{\{ x\}^\perp}\cct\big|\rx\over |\cct x|}
  = {1\over\sqrt{2\lambda}}\Proj_{\{ x\}^\perp}\cct\big|\rx \, .
$$
Clearly, since $\da1$ is a sphere, its tangent space at $x$ is
$\{\, x\,\}^\perp\cap H$. It follows that for $x$ belonging to $\da1$,
$$
  T_x\Lambda_{I(A_1)}\ominus T_x\da1
  = \{\, x\,\}^\perp \ominus \big( \{\, x\,\}^\perp\cap H\big)
  = \{\, x\,\}^\perp\ominus H
  = H^\perp
$$
since $x$ is in $H$ as well. Thus, for $x$ in $\da1$,
$$
  G_{A_1}(x) = \Proj_{H^\perp}\sqrt{\lambda\over 2} 
  \Big( \Id - {\cct\over\lambda} \Big)\Big|_{H^\perp} \, .
$$
This matrix does not depend on $x$. Therefore, putting all the 
pieces together,
$$
  c_1={ (2\pi)^{d-k-1\over 2} 
  \vol \big(S_H(0,\sqrt{2/\lambda})\big)
  \over (2/\lambda)^{d-k+1\over 4}(\lambda/2)^{d-k-1\over 4}
  \det\Big(\Proj_{H^\perp}\big( 
           \Id-{\cct\over\lambda}\big)\Big|_{H^\perp}
  \Big)^{1\over 2}
  } \, .
$$
Using the classical fact
$$
  \vol (S_{n-1})=
  { 2\pi^{n/2}\over 
    \Gamma\big( {\displaystyle n\over\displaystyle 2}\big)} \, ,
$$
we obtain 
$$
  c_1={ (2\pi)^{d/2} \over \lambda^{(k-1)/2}\Gamma
  \Big({\displaystyle k+1\over\displaystyle 2}\Big)
  \det\Big(\big(\Id-{\displaystyle\cct\over\displaystyle\lambda}
           \big)\Big|_{H^\perp}
  \Big)^{1/2}} \, .
$$
The result follows after dividing $c_1$ by $(2\pi)^{d/2}$, the
normalizing factor of the $d$-dimensional standard normal
density.\hfill$\qed$

\bigskip

\state{REMARK.} It is interesting to notice the discontinuity in
the polynomial term in $t$ in Theorem 8.1.2. In (8.1.1), this
term has degree $-(\alpha-2){d\over 4}-{\alpha\over 4}$, which
is strictly less than $-1/2$ for $\alpha\geq 2$, and equals
$-1/2$ for $\alpha =2$. For $\alpha=2$, Theorem 8.1.2 gives a
polynomial term of degree at most $0$ since $k$ is at least $1$. 
Notice also that the map $C\mapsto k$ is not continuous for any
standard topology on the set of matrices.

The determinant of $M$ involved in Theorem 8.1.2 can be given more
explicitly as a function of the matrix $C$. Write $\lambda_1\leq
\cdots \leq \lambda_d$ for the spectrum of $C+C^\T$.  Since the dimension
of $H$ is $k+1$, we have $\lambda_{d-k}=\lambda_{d-k+1}=\ldots
=\lambda_d$.  Diagonalizing $C+C^\T$ and noticing that $H^\perp$ is
invariant under $\Id + \lambda_d^{-1}\cct$, we have
$$
  \det M =\prod_{1\leq i\leq d-k-1} 
  \Big( 1-{\lambda_i\over\lambda_d}\Big) \, .
$$

\bigskip

Applying Theorem 7.5, we obtain also the following.

\bigskip

\state{8.1.3. PROPOSITION.} {\it%
  For $\alpha=2$ and $\lambda$ positive, the conditional distribution
  of $X/\sqrt t$ given $\cXX\geq t$ converges weakly* to the
  uniform distribution over the sphere centered at the origin of radius
  $\sqrt{2/\lambda}$ of $H$.}

\bigskip

\stateit{Proof.} Notice that $\calM_{\da1}$ is proportional to the
uniform distribution on $S_H(0,\sqrt{2/\lambda})$ and that $|\D I|$ 
as well as $\det\, G_{A_1}$ are constant on $\da1$. The result follows
from the proof of Theorem 8.1.2. and Theorem 7.5. \hfill$\qed$

\bigskip

\chapternumber={8}
\chaptername={Application to quadratic forms of random vectors}
\noindent{\fsection 8.2. An example with heavy tail
distribution.}
%\chapternumber{8.2}
\section{8.2}{An example with heavy tail distribution}

\bigskip

\noindent In this section, we consider a random
vector $X=(X_1,\ldots , X_d)$ in $\RR^d$, with independent and
identically distributed components, all having a Student-like
distribution with parameter $\alpha$.  Thus, $X_i$ has a density, and
there exists a constant $K_{s,\alpha}$\ntn{$K_{s,\alpha}$} such that
$$
  P\{\, X_i\leq -x\,\} \sim P\{\, X_i\geq x\,\} 
  \sim {K_{s,\alpha}\alpha^{(\alpha-1)/2}\over x^\alpha}
$$
as $x$ tends to infinity. Define\ntn{$A_t$}
$$
  A_t=\big\{\, x\in\RR^d  : \cxx\geq t\,\big\} 
  = \sqrt{t} A_1 \, .
$$
Writing $s_\alpha(\cdot )$ for the density of a single $X_i$, the
density of the vector $X$, given 
by $\prod_{1\leq i\leq d}s_\alpha(x_i)$. It is not log-concave. It is
not even specified at all, except by an asymptotic equivalent! Thus,
we cannot use Theorem 5.1 in a straightforward way to approximate
$$
  P\big\{\, \cXX\geq t\,\big\} = \int_{A_t}\prod_{1\leq i\leq d}
  s_\alpha (x_i)\d x_i \, .
$$
However, as pointed out in the introduction, chapter 1, we can make a
change of variables, and then try to use Theorem 5.1. This will
require all the power of Theorem 5.1, in particular the freedom on the
set $A$ that is allowed.

To state our first result, define\ntn{$J_1$}
$$
  J_1 = \{\, j : C_{j,j}>0\,\} \, .
$$

\bigskip

\state{8.2.1. THEOREM.} {\it%
  Let $X$ be a $d$-dimensional random vector with independent
  and identically distributed components having a Student-like
  distribution. Let $C$ be a $d\times d$ matrix.
  If $J_1$ is not empty, then
  $$
    P\big\{\, \cXX\geq t\,\big\} \sim
    K_{s,\alpha} \alpha^{(\alpha-1)/2}
    {2\over t^{\alpha /2}}
    \sum_{j\in J_1} C_{j,j}^{\alpha/2}
    \qquad \hbox{ as } t\to\infty \, .
  $$}
\finetune{
%\bigskip
}
\stateit{Proof.} Let us first make a change of
variable so that we will be able to apply Theorem 5.1.

Let $Y=(Y_1,\ldots ,Y_d)$ be a random vector with centered
normal distribution, with identity covariance matrix. Its
density,
$$
  {1\over (2\pi )^{d/2}}\exp \Big( -{|y|^2\over 2}\Big)
  \, , \qquad y\in\RR^d \, ,
$$
is log-concave. Let us write\ntn{$\Phi (\cdot )$}
$$
  \Phi(y)=\int_{-\infty}^y 
  {e^{-u^2}\over \sqrt{2\pi}} \d u \, ,
  \qquad y\in\RR \, ,
$$
the normal cumulative distribution function.
Similarly, denote by\ntn{$S_\alpha (\cdot )$}
$$
  S_\alpha (y) 
  = \int_{-\infty}^y s_\alpha (u) \d u
$$
the Student-like cumulative distribution function of each
individual $X_i$. Writing
$$
  S_\alpha^\leftarrow (u)=\inf\big\{\, y  : S_\alpha(y)\geq u\,\big\}
$$
for the inverse function of $S_\alpha$ and analogously
$\Phi^\leftarrow$ for the inverse function of $\Phi$, we see that
$\Saip (Y_i)$ has the same distribution as $X_i$.

\bigskip

\state{NOTATION.} {\it 
  Let us agree that a function $g$
  defined on $\RR$ is extended componentwise to $\RR^d$, 
  writing $g(x_1,\ldots ,x_d)$ for $\big( g(x_1),\ldots , 
  g(x_d)\big)$. }

\bigskip

It follows that $X$ has the same distribution as $\< C\Saip
(Y)\, ,\Saip (Y)\>$. In other words, defining\ntn{$B_t$}
$$
  B_t=\big\{\, y\in\RR^d  : \< C\Saip (y)\, ,\Saip (y)\>
  \geq t\,\big\} 
  = \pisa (A_t) \, ,
$$
we have
$$
  P\big\{\, \cXX \geq t\,\big\}
  = \int_{A_t} \prod_{1\leq i\leq d} s_\alpha(x_i)\d x_i
  = \int_{B_t} { e^{-|y|^2/2}\over (2\pi )^{d/2} } \d y \, .
  \eqno{(8.2.1)}
$$
Since $\Saip$ is continuous and defined on the whole real line, we see
that for any positive $M$ and any $t$ large enough, the ball centered
at the origin and of radius $M$ does not intersect $B_t$. Thus, $B_t$
moves to infinity as $t$ tends to infinity, and the right hand side of
(8.2.1) is the integral of a log-concave function over a set moving to
infinity as $t$ tends to infinity. We can try to apply Theorem 5.1.

It should be noticed that we could make a change of variable
leading to a different distribution than the standard Gaussian one.
However, this one is rather convenient since its level sets and
their geodesics are known explicitly.

The disadvantage of the change of variable is of course that the
set $B_t$ is more complicated than $A_t$. Nevertheless, whatever
information is needed on $B_t$ can be first read on $A_t$, and
then pulled back to $B_t$. This fact is illustrated by
Proposition 8.2.4 bellow, where we will calculate $\calD_{B_t}$.
This change of variable technique works mainly because $\pisa$ has
an explicit and simple asymptotic equivalent.

To apply Theorem 5.1, let us define\ntn{$I(\cdot )$}
$$
  I(y)={|y|^2\over 2} + \log (2\pi )^{d/2} \, ,
$$
that is minus the logarithm of the Gaussian density. The
function $I$ is convex.

We will make use of the following elementary result in
asymptotic analysis, whose proof can be found in appendix 1,
\finetune{\hfuzz=1.2pt}
$$
  \big(\pisa(x)\big)^2 = 2\alpha\log x -\log\log x 
  - 2\log (K_{s,\alpha}\alpha^{\alpha /2}) 
  - 2\log (2\sqrt\pi ) + o(1)
$$
as $x$ tends to infinity. It implies
$$
\finetune{\hfuzz=0pt}
  \pisa (x)=\sqrt{2\alpha\log x}+ o(1) \qquad 
  \hbox{ as } x\to\infty \, .
$$
It is also convenient to introduce the canonical basis
$e_1,\ldots ,e_d$ of $\RR^d$. For any $j$ in $J_1$ and $\epsilon$
in $\{\, -1,1\,\}$, the point $p_{\epsilon ,j,t}=\epsilon \sqrt{
t/C_{\smash{j,j}}}\, e_j$ belongs to $\partial A_t$. Thus, $q_{\epsilon ,
j,t} = \pisa (p_{\epsilon ,j,t})$ belongs to $\partial B_t$. 
The following lemma gives a
parameterization of $\partial A_t$ and $\partial B_t$ near
$p_{\epsilon ,j,t}$ and $q_{\epsilon ,j,t}$. This
describes these boundaries locally.

\bigskip

\state{8.2.2. LEMMA.} {\it%
  The tangent space of the boundary $\partial A_t$ at 
  $p=p_{\epsilon,j,t}$
  is $\big\{\, \cct p\,\big\}^\perp$.  Near $p$, the boundary
  of $A_t$ can be parametrized as\ntn{$p(v)$}\ntn{$p_{\epsilon ,j,t}(v)$}
  $$
    p(v)=p_{\epsilon ,j,t} (v) =
    \epsilon \sqrt{t\over C_{j,j}}\Big[ 1-{1\over 2t} \< Cv,v\>
    + O\Big({\< Cv,v\>\over t}\Big)^2\Big] e_j + v \ ,
  $$
  for $v$ in $T_p\partial A_t$, and $|v|=o(\sqrt{t})$ 
  as $t$ tends to infinity.

  The boundary $\partial B_t$ near $q=q_{\epsilon ,j,t}
  =\pisa (p)$ can be parametrized as\ntn{$q(v)$}\ntn{$q_{\epsilon,j,t}(v)$}
  $$\displaylines{\qquad
    q(v)=q_{\epsilon ,j,t}(v) = \epsilon \bigg( 
    \sqrt{ \alpha\log{t\over C_{j,j}} }
    - {\llst\over 2\sqrt{\alpha\log t}}
    - { \log (K_{s,\alpha} \alpha^{\alpha/2}2\sqrt{\pi})
       \over \sqrt{\alpha\log t} }
    \cut
    {}+ o\Big( {1\over\sqrt{\log t}} \Big)\bigg) e_j
      + \sum_{\scriptstyle 1\leq i\leq d\atop\scriptstyle i\ne j}
      \pisa (v_i) e_i
    \qquad\cr
  }$$
  for $v$ in $T_p\partial A_t$ and $|v|=o(\sqrt{t})$
  as $t$ tends to infinity.
  }

\bigskip

\stateit{Proof.} The assertion on the tangent space of $\partial
A_t$ at $p$ is plain since the differential of the map $x\mapsto
\< Cx,x\>$ at $p$ is $(C+C^\T )p$. Near the point
$\sqrt{1/C_{i,i}}$, we can parameterize $\partial A_1$ by its
tangent plane. This leads to the following parameterization of
$\partial A_t$. Let $h(v)$ be such that
$$
  \tilde p(v) 
  = \sqrt{t/C_{j,j}} \big(1 +h(v)\big) e_j + v\in\partial A_t
$$
for all $v$ in $T_p\partial A_t$ with $|v|$ not too large. This
inclusion becomes
$$
  t= \< C\tilde p(v),\tilde p(v) \>
  = t\big( 1+h(v)\big)^2 + \< Cv,v\> \, .
  \eqno{(8.2.2)}
$$
An approximation of $h$ follows either by working out an asymptotic
expansion for $h(\cdot )$ or using the following argument. For
$|v|=o(\sqrt{t})$, (8.2.2) implies $h(v)=o(1)$ as $t$ tends to
infinity. Rewriting (8.2.2) as the quadratic equation in $h$,
$$
  0= th(v)^2 + 2th(v) + \< Cv,v\> \, ,
$$
we obtain 
$$
  h(v)=-1+ \Big(1-{\< Cv,v\>\over t}\Big)^{1/2}
  \qquad \hbox{ as } t\to\infty \, .
$$
This gives the asymptotic expansion for $p(v)$.

We then pull back the expression of $p(v)$ to parameterize
$\partial B_t$ by $q(v)=\pisa \big( p(v)\big)$. Notice first
that
$$
  \log\bigg( \sqrt{t\over C_{j,j}}\big( 1+ O( t^{-1}|v|^2)\big)
      \bigg)
  = {1\over 2} \log{t\over C_{j,j}} + O( t^{-1}|v|^2) \, .
$$
In the range $|v|=o(\sqrt{t})$, the asymptotic expansion for
$\pisa$ in Lemma A.1.5 gives 
$$\eqalign{
  \epsilon \< q(v),e_j\>\,
  & { = \pisa \bigg( \sqrt{t\over C_{j,j}} 
    \big(  1+ O( t^{-1}|v|^2)\big)\bigg) } \cr
  & { = \sqrt{\alpha\log{t\over C_{j,j}}} 
    - {\llst \over 2\sqrt{\alpha\log t}} }
    -{ \log (K_{s,\alpha} \alpha^{\alpha/2} 2\sqrt{\pi})\over
       \sqrt{\alpha\log t} } \cr
  & \phantom{\smash{ = \sqrt{\alpha\log{t\over C_{j,j}}}} 
    - {\llst \over 2\sqrt{\alpha\log t}} }  
    + o(\log t)^{-1/2} \, . \cr
  }
$$
On the other hand, for $i\neq j$,
$$
  \epsilon \< q(v),e_i\> = \pisa (v_i)\ .
$$
This proves Lemma 8.2.2. \hfill$\qed$

\bigskip

It is somewhat important for what follows to have some intuition
on the shape of $\partial B_t$ near $q_{\epsilon ,j,t}$. This is
precisely what the last assertion of Lemma 8.2.2 gives us.
Recall that $p_{\epsilon ,j,t}$ is collinear to $e_j$. As
$v$ varies in $T_{p_{\epsilon ,j,t}}\partial A_t=\{\,
(C+C^\T )e_j\,\}^\perp$, the term $\sum_{1\leq i\leq d ; i\neq
j}\pisa (v_i)e_i$ in the expression of $q_{\epsilon ,j,t}(v)$
varies too. If the $v_i$'s were allowed to vary independently,
then $\sum_{1\leq i\leq d ; i\neq j}\pisa (v_i)e_i$ would
describe the hyperplane $\{\, e_j\,\}^\perp$, and $\partial B_t$
would be a hyperplane perpendicular to $e_j$, passing through
$q_{\epsilon,j,t}$. This is not quite the case of course, but
almost, provided we look at the right scale. This is the meaning
of the next claim.

\bigskip

\state{8.2.3. CLAIM.} {\it%
  For $t$ large enough and $j$ in $J_1$, the set
  $$
    \Big\{\, \sum_{\scriptstyle 1\leq i\leq d\atop 
    \scriptstyle i\neq j} \pisa (v_i)e_i \, : \,
    v\perp (C+C^\T )e_j \, , \, |v|=o(\sqrt{t})\,\Big\}
  $$
  is contained in 
  $$
    \Big\{\, \sum_{\scriptstyle 1\leq i\leq d\atop 
    \scriptstyle i\neq j} w_ie_i \, :\, |w|\leq 
    {1\over 2} \sqrt{\alpha\log t}\, \Big\} \, .
  $$
  }

\bigskip

\stateit{Proof.} One may argue as follows. Notice first that
$\{\, (C+C^\T )e_j\,\}^\perp$ does not contains $e_j$. Indeed, if this
were the case, we would have $C_{j,j}=\<(C+C^\T )e_j\, , e_j\>=0$,
contradicting the fact that $j$ belongs to $J_1$. Consequently, as $v$
varies in $\{\, (C+C^\T )e_j\,\}^\perp$, the vector
$\sum_{\scriptstyle 1\leq i\leq d ; i\neq j}\pisa (v_i)e_i$ describes
the space spanned by the $e_i$'s for $1\leq i\leq d$ and $i\neq
j$. Finally, if $w$ is orthogonal to $e_j$ and of norm less than
$(1/2)\sqrt{\alpha\log t}$, then $w=\sum_{\scriptstyle 1\leq i\leq d ;
i\neq j}\pisa (v_i)e_i$ for some $v$ in $(C+C^\T )e_j$. Furthermore,
$\pisa (v_i)^2\leq {\alpha\over 2} \log t$. From Lemma A.1.5, we then
infer $|v_i|\leq t^{3/8}$ for $t$ large enough. The relation $v\perp
(C+C^\T )e_j$ forces then $|v_j|\leq O(t^{3/8})$, and so
$|v|=o(\sqrt{t})$.  This proves our claim. \hfill$\qed$

\bigskip

We can now locate the interesting minima of $I$ over $\partial
B_t$. They will provide a good guess for a dominating
manifold, as well as an estimation of $I(B_t)$.

\bigskip

\state{8.2.4. PROPOSITION.} {\it%
  Assume that $J_1$ is nonempty. If $y$ belongs to $\partial B_t$
  and $I(y)\leq I(B_t)+ O(1)$ as $t$ tends to infinity, 
  then $y$ is in
  a $O(1)$-neighborhood of a points $q_{\epsilon ,j,t}$ for
  some $j$ in $J_1$ and some $\epsilon$ in $\{\, -1,1\,\}$. 
  Moreover, as $t$ tends to infinity,
  $$\eqalign{
      I(q_{\epsilon ,j,t})
      = {\alpha\over 2}\log t -{1\over2}\llst
    & - \log (K_{s,\alpha}\alpha^{\alpha/2}2\sqrt{\pi}) \cr
    & + \log (2\pi)^{d/2} - {\alpha\over 2} \log C_{j,j}
      + o(1) \, . \cr
  }$$}

\bigskip

\stateit{Proof.} By the very definition of $I$
and $B_t$, 
$$\displaylines{\qquad
  I(B_t)-\log (2\pi )^{d/2}
  \hfill\cr\noalign{\vskip 1pt}\hfill
  \eqalign{
  \noalign{\vskip 2mm}
    =& \inf\big\{\, |y|^2/2  : \< C\Saip (y)\, ,\Saip (y)\> 
       \geq t\,\big\} \cr
  \noalign{\vskip 1mm}
    =& {1\over 2} \inf\Big\{\, \sum_{1\leq i\leq d} \pisa (
      \lambda\sqrt{t} u_i)^2  : \< Cu\, ,u\> = 1 \,  , \,
      \lambda \geq 1\,\Big\}  \cr
    =& {1\over 2} \inf\Big\{\, \sum_{1\leq i\leq d} \pisa (
      \sqrt{t} u_i)^2  : \< Cu\, ,u\> = 1
      \,\Big\}  \, , \cr }
  \qquad
  }$$
the second equality coming from the change of variable
$y_i=\pisa (\lambda \sqrt{t}u_i)$; the last one comes from the
fact that the function $\lambda\in [\,0,\infty
)\mapsto \pisa (\lambda\sqrt{t}u_i)^2$ is increasing for $t$
large enough and $u_i$ fixed.

If $s$ is positive and such that $\sqrtt s$ tends to
infinity and $\log s/\log t$ tends to $0$ as $t$ tends 
to infinity, the asymptotic expansion for $(\pisa )^2$ in
Lemma A.1.5. shows that
$$
  \pisa (\sqrt{t}s)^2 = \alpha\log t \big(1+o(1)\big)
  \qquad\hbox{ as } t\to\infty \, .
  \eqno{(8.2.3)}
$$
If $\< Cu\, ,u\>= 1$, and $r$ of the $u_i$'s, say $u_1,\ldots ,
u_r$, are of order larger than $1/\log t$, i.e. $\min_{1\leq
i\leq r}|u_i|\gg 1/\log t$ as $t$ tends to infinity, 
then (8.2.3) yields
$$
  \sum_{1\leq i\leq d} \pisa (\sqrtt u_i)^2
  \geq \sum_{1\leq i\leq r} \pisa (\sqrtt u_i)^2
  \sim r\alpha \log t
  \eqno{(8.2.4)}
$$
as $t$ tends to infinity. Hence, to minimize the left hand side of
(8.2.4), we should have $r$ as small as possible. But $r$ must be at
least $1$, for $\< Cu,u\>=1$. Moreover $r=1$ can be achieved by
considering $j$ in $J_1$ and $u=\epsilon e_j/\sqrt{C_{\smash{j,j}}}$
for some $\epsilon$ in $\{\, -1,1\,\}$. This leads us to look at the
function $I$ near $q_{\epsilon ,j,t}=\pisa (\epsilon\sqrtt
e_j/\sqrt{C_{\smash{j,j}}})$.  Furthermore, if $I\big(\pisa
(\sqrt{t}u)\big)$ is minimal, $u$ must be on the boundary of a $O(\log
t)^{-1}$-neighborhood of $\epsilon e_j/\sqrt{C_{\smash{j,j}}}$ for
some $\epsilon$ in $\{\, -1,1\,\}$ and $j$ in $J_1$. Consequently,
$\sqrt{t}u$ is in $\partial A_t$ and in an $O(\sqrt{t}/\log t)=
o(\sqrt{t})$-neighborhood of $p_{\epsilon ,j,t}$. Therefore, when
studying such a point, we can use the parameterization given in Lemma
8.2.2. This also leads us to look at the function $I\big( q_{\epsilon
,j,t}(v)\big)$ for $v$ in $T_{p_{\epsilon,j,t}}\partial A_t$ and
$|v|=o(\sqrt{t})$.

Write $v_i=\< v,e_i\>$ for the components of the vector
$v$ belonging to $T_{p_{\epsilon ,j,t}}\partial A_t$. 
Using Lemma 8.2.2, we obtain
$$\displaylines{
  \qquad
  I\big( q_{\epsilon ,j,t}(v)\big)
  = {\alpha\over 2} \log {t\over C_{j,j}} 
  - {1\over 2}\log\log \sqrt{t}
  - \log (K_{s,\alpha}\alpha^{\alpha/2}2\sqrt{\pi})
  \cut
  {}+{1\over 2}
  \sum_{\scriptstyle 1\leq i\leq d\atop\scriptstyle i\neq j}
  \pisa (v_i)^2 + \log (2\pi )^{d/2} + o(1)
  \qquad\qquad\llap{(8.2.5)}\cr}
$$
as $t$ tends to infinity, and uniformly in $|v|=o(\sqrtt )$. 
Therefore, up to $o(1)$ as $t$ tends to infinity, the function 
$v\mapsto I\big(q_{\epsilon ,j,t}(v)\big)$ is
minimum at $0$, and its minimum value is $I(q_{\epsilon ,j,t})$ 
as claimed. \hfill$\qed$

\bigskip

Notice that the proof of Proposition 8.2.4 gives actually a
little bit more, and this will be useful. Indeed, if
$x$ is in $\partial B_t$ and $I(x)=I(B_t)+ o(\log\log t)$, then 
(8.2.3)--(8.2.4) and Lemma A.1.5 show that $x$ is in a 
$o(\log\log t)$-neighborhood of some $q_{\epsilon ,j,t}$. Indeed,
we must have 
$$
  \max_{\scriptstyle 1\leq i\leq d\atop \scriptstyle i\neq j}
  |\pisa (v_i)|=O(\log\log t)^{1/2} \, .
$$

In view of Proposition 8.2.4 and its proof, we can start 
to apply Theorem 5.1
in calculating a few terms of the asymptotic formula. Indeed,
define\ntn{$\gamma_1$}
$$
  \gamma_1=\max_{1\leq j\leq d}C_{j,j} \, .
$$
We have immediately
$$\eqalign{
    I(B_t)= {\alpha\over 2} \log t -{1\over 2} \llst
    &- \log (K_{s,\alpha}\alpha^{\alpha/2}2\sqrt{\pi}) \cr
    &+ \log (2\pi )^{d/2} - {\alpha\over 2}\log \gamma_1 + o(1)
     \cr}
$$
as $t$ tends to infinity. Moreover, a candidate for the dominating
manifold is
$$
  \big\{\, q_{\epsilon ,j,t} : j\in J_1 \, ,\,
  \epsilon \in \{\, -1,1\,\} \,\big\} \, .
$$
Unfortunately, this choice does not match with our definition of
a dominating manifold. It is indeed required in the definition
that it is also a base manifold, and as such belongs to
$\smash{\Lambda_{I(B_t)}}$. The expression $I(y)=(|y|^2/2)+\log
(2\pi )^{d/2}$ and Proposition 8.2.4 shows that the points
$q_{\epsilon ,j,t}$, for $j$ in $J_1$, cannot lie on the same
sphere centered at the origin. But they almost do!

Let us denote by $\rho_t$ the radius of the sphere
$\Lambda_{I(B_t)}$. The expression for $I$ shows that
$\rho_t=\sqrt{I(B_t)+d\log (2\pi )}$. Define $r_{\epsilon
,j,t}=\rho_tq_{\epsilon ,j,t}/|q_{\epsilon ,j,t}|$. We consider
the candidate\ntn{$\calD_{B_t}$}
$$
  \calD_{B_t} =\big\{\, r_{\epsilon ,j,t} \, :\, j\in J_1 \, , \,
  \epsilon\in\{\, -1,1\,\} \,\big\}
$$
for a dominating manifold. It will be clear after Lemma 8.2.7 that
what we are really doing here is moving the points
$q_{\epsilon,j,t}$ through the normal flow, until they reach the
level line $\Lambda_{I(B_t)}$; this gives $r_{\epsilon ,j,t}$
--- somehow unfortunately for the clarity of the argument, but
luckily for the calculation, this move along the normal flow 
and the Euclidean projection on the sphere coincide
when working with the normal distribution.

Since $\calD_{B_t}$ is of dimension $0$, its Riemannian measure
is a sum of point masses,
$$
  \calM_{\calD_{B_t}} = 
  \sum_{\scriptstyle j\in J_1\atop 
        \scriptstyle\epsilon\in \{ -1,1\}} 
  \delta_{r_{\epsilon ,j,t}} \, .
$$
From Proposition 8.2.4 and the above expression for $I(B_t)$, 
we infer that for $j$ in $J_1$,
$$
  \tau_{B_t}(r_{\epsilon ,j,t}) = {\alpha\over 2}
  \log{\gamma_1\over \cjj} + o(1) \qquad
  \hbox{ as }t\to\infty \, .
$$
Since $\D I$ is the identity function
and $I(y)={\displaystyle |y|^2\over\displaystyle 2}+\log (2\pi )^{d/2}$, 
we also have
$$
  |\D I(r_{\epsilon ,j,t})|
  =\rho_t
  =\sqrt{2I(B_t)+d\log (2\pi )}
  \sim \salt 
  \qquad \hbox{ as } t\to\infty \, .
$$
If we can apply Theorem 5.1, we obtain

\setbox100=\hbox{$\displaystyle\int_{B_t}\eIx \d x$}
\wd100=0pt
$$\eqalignno{
  \box100\qquad
  & \cr
  & \hskip -.2in\sim e^{-I(B_t)} (2\pi )^{(d-1)/2}
    \sum_{\scriptstyle j\in J_1\atop
          \scriptstyle\epsilon\in\{-1,1\}}
    { \exp \big( -\tau_{B_t}(q_{\epsilon ,j})\big) \over
	|\D I(r_{\epsilon ,j,t)}|^{(d+1)/2} 
	\big( \det\, G_{B_t}(r_{\epsilon ,j,t})\big)^{1/2} } \cr
  & \hskip -.2in\sim { (\log t)^{(1-d)/4}\over t^{\alpha /2} }
    K_{s,\alpha} \alpha^{(2\alpha -1-d)/4}
    \sum_{\scriptstyle j\in J_1\atop 
          \scriptstyle \epsilon\in\{ 1,1\}}
    { C_{j,j}^{\alpha /2} \over 
      \big( \det\, G_{B_t}(r_{\epsilon ,j,t})\big)^{1/2} }
    &(8.2.6)\cr
  }
$$
as $t$ tends to infinity. We are left with calculating
$G_{B_t}(r_{\epsilon ,j,t})$ and checking the assumptions of Theorem
5.1. In order to calculate $G_{B_t}(r_{\epsilon ,j,t})$, we need to
calculate $\tau_{B_t}$, and ultimately the normal flow. This turns out
to be particularly easy for the normal distribution.

\bigskip

\state{8.2.5. LEMMA.} {\it%
  For the Gaussian distribution $\calN (0,\Id )$ on $\RR^d$,
  the normal flow is given by $\psi (q,s)=
  \sqrt{ 1+{\displaystyle 2s\over\displaystyle |q|^2}}\, q$. 
  }

\bigskip

\stateit{Proof.} The level lines $\Lambda_c$ are
spheres centered at the origin since $I$ is a spherical function.
Hence, $\psi (q,s)$ moves on a straight line through the origin
as $s$ varies, and $\psi (q,s)=a(s) q$ for some function
$a(\cdot )$. We obtain $a$ from the equation
$$
  I(q)+s
  = {|q|^2\over 2} + \log (2\pi )^{d/2} + s  
  = I\big( \psi (q,s)\big) 
  = {a(s)^2\over 2} |q|^2 + \log (2\pi )^{d/2} \, .
$$
That is,
$a(s) = \sqrt{1+{\displaystyle 2s\over \displaystyle |q|^2}}\,$. 
\hfill$\qed$

\bigskip

Since the exponential map on the level line $\Lambda_{I(B_t)}$ 
is involved in the definition of the curvature term $G$, 
we first recall its expression in the Gaussian case.

\bigskip

\state{8.2.6. LEMMA.} {\it%
  If $I$ is a spherical function, then
  for $q$ in $\Lambda_c$, we have
  
  \noindent (i) $T_q\Lambda_c= \{\, q\,\}^\perp$,

  \noindent (ii) $\exp_q (w) = \cos\Big(
    {\displaystyle |w|\over\displaystyle |q|}\Big)q
  + \sin\Big({\displaystyle |w|\over\displaystyle |q|}\Big)|q| 
    {\displaystyle w\over\displaystyle |w|}$, for all
  $w$ in $T_q\Lambda_c$ with $|w|\leq \pi |q|$.
  }

\bigskip

\stateit{Proof.} Since $\Lambda_c$ is a sphere centered at the
origin, (i) follows. The geodesics on $\Lambda_c$ are circles of
maximal diameter. By cutting $\Lambda_c$ along the plane determined by
$q$ and $w$, the expression of the maximal circle leaving $q$ in the
direction $v$, that is (ii), follows.
\hfill$\qed$

\bigskip

Since we calculated $I(q_{\epsilon,j,t})$ in Proposition 8.2.4, 
it is easier to calculate
$\tau_{B_t}\big(\exp_{q_{\epsilon ,j,t}}(v)\big)$ than
$\tau_{B_t}\big( \exp_{r_{\epsilon ,j,t}}(v)\big)$. The
following lemma will be instrumental in relating these
quantities. It is specific to the Gaussian situation. Since 
$T_q\Lambda_{I(q)}=\{\, q\,\}^\perp$, we can identify 
$T_q\Lambda_{I(q)}$ and 
$T_{\lambda q}\Lambda_{I(\lambda q)}$ for any nonzero $\lambda$.

\bigskip

\state{8.2.7. LEMMA.} {\it%
  For the standard Gaussian distribution, i.e., $I(x)=|x|^2/2
  +\log (2\pi)^{d/2}$,

  \noindent (i) for any nonzero $q$, any positive $\lambda$ 
  and any $w$ in $T_q\Lambda_{I(q)}\equiv T_{\lambda q}
  \Lambda_{I(\lambda q)}$, we have
  $$
    \exp_{\lambda q}(w)=\lambda \exp_q (w/\lambda ) \, ;
  $$

  \noindent (ii) moreover, for any set $B$, any positive 
  $\lambda$, and any $q$ in $\RR^d$ such
  that the line segment between $\lambda q$ and $q$ does 
  not intersect $B$,
  $$
    \tau_B(\lambda q) = {|q|^2\over 2} (1-\lambda^2) + \tau_B(q)\, .
  $$}

\bigskip

\stateit{Proof.} (i) follows from Lemma 8.2.6 since
$$
  \exp_{\lambda q}(w) = \cos\Big( {|w|\over\lambda |q|} \Big) \lambda q
  + \sin \Big( {|w|\over\lambda |q|} \Big) |\lambda q|{w\over |w|}
  = \lambda\exp_q\Big({w\over\lambda}\Big) \, .
$$
To prove (ii), the condition that the segment
between $\lambda q$ and $q$ does not intersect $B$,
the fact that the normal flow moves along straight lines through
the origin, and the definition of $\tau_B$ imply
$\psi\big(q,\tau_B(q)\big) 
= \psi \big( \lambda q,\tau_B(\lambda q)\big)$. The
expression of the normal flow in Lemma 8.2.5 gives 
$$
  \sqrt{1+2{\tau_B(q)\over |q|^2}}\, q
  =\sqrt{1+2 {\tau_B(\lambda q)\over |\lambda q|^2}}\, \lambda q
  \, .
$$
The formula for $\tau_B(\lambda q)$ follows.\hfill$\qed$

\bigskip

\state{REMARK.} The essence of Lemma 8.2.7 is to relate
$\exp_{\psi(q,s)}\circ\psi_{s,*}$ and $\psi_s\circ\exp_q$.
Both maps act on $T_q\Lambda_{I(q)}$. We can write one as
the other one composed with some transform of $T_q\Lambda_{I(q)}$.
Lemma 8.2.7 makes this explicit in the Gaussian case.

\bigskip

We can obtain an approximation of $\tau_{B_t}$ provided its argument
is not too far away from $q_{\epsilon ,j,t}$ for some $j$ in $J_1$ and
some $\epsilon$ equal to $-1$ or $+1$. Since we made everything
explicit up to $o(1)$, our approximation will not be good enough to
check (5.3), but perfectly fine to check (5.17) --- one may try to
check (5.3) and hopefully will agree that (5.17) is a useful
refinement. Hence, we are ready to calculate the curvature term
$G_{B_t}$.

\bigskip

\state{8.2.8. LEMMA.} {\it%
  In the range $|w|= o(\sqrt{\log t})$, we have
  $$
    \tau_{B_t}\big( \exp_{r_{\epsilon ,j,t}}(w)\big)
    = \tilde\tau_{B_t}\big( \exp_{r_{\epsilon ,j,t}}(w)\big)
    + o(1) \, , 
  $$
  with
  $$
    \tilde\tau_{B_t}\big( \exp_{r_{\epsilon ,j,t}}(w)\big)
    -\tilde\tau_{B_t}(r_{\epsilon ,j,t})
    = {1\over 2}\< w,w\> \, .
  $$
  Consequently,
  $$
    G_{B_t}(r_{\epsilon ,j,t})
    \sim {\Id_{\RR^{d-1}}\over |\D I(q_{\epsilon ,j})|}
    \sim {\Id_{\RR^{d-1}}\over \salt}
    \qquad \hbox{ as } \ttoi \, .
  $$}

\bigskip

\stateit{Proof.} Let us first obtain an approximation for
$\smash{\tau_{B_t}\big(\exp_{q_{\epsilon ,j,t}}(w)\big)}$.
Lemma 8.2.2 shows that in the
range $|v|=o(\smash{\sqrtt} )$,  near $q_{\epsilon ,j,t}$, 
the surface $\partial B_t$ parametrized
by $v\mapsto q_{\epsilon ,j,t}(v)$, is given by the equation of
the hyperplane $\< x,e_j\> = |q_{\epsilon ,j,t}|$, up
to $o(\log t)^{-1/2}$.
Consequently, for $q=q_{\epsilon ,j,t}$, for $w$ in
$T_q\Lambda_{I(B_t)}$ and $s=\tau_{B_t}\big( \exp_q (w)\big)$,
using Lemmas 8.2.5, 8.2.6, we obtain
$$
  |q|+o(\log t)^{-1/2}
  = \big\< \psi\big(\exp_q(w) ,s\big) , e_j\big\>
  = \sqrt{1+{2s\over |q|^2}} \cos \Big({|w|\over |q|}\Big)|q| 
  \, .
$$
It follows that
$$\eqalign{
  \tau_{B_t}\big(\exp_q(w)\big)
  & = {|q|^2\over 2} \bigg[ \bigg( 
    { |q|+o(\log t)^{-1/2}\over |q|\cos (|w|/|q|)}\bigg)^2-1
    \bigg] \cr
  & = {|w|^2\over 2} + o(1) \qquad \hbox{ as } \ttoi \, ,\cr
  }$$
in the range $|w|=o(|q|)= o(\log t)^{1/2}$. From Lemma 8.2.6, we
deduce
$$\eqalign{
  \tau_{B_t}\big(\exp_{r_{\epsilon ,j,t}}(w)\big)
  & = \tau_{B_t}\Big[ {\rho_t\over |q|}\exp_q\Big(w 
    {|q|\over \rho}\Big)\Big] \cr
  & = {|q|^2\over 2}\Big( 1-{\rho_t^2\over |q|^2}\Big)
    + \tau_{B_t}\Big[ \exp_q\Big( w{|q|\over \rho_t}\Big)\Big]
    \, . \cr
  }
$$
Therefore, since $\tau_{B_t}(q)=o(1)$,
$$\eqalign{
  \tau_{B_t}\big(\exp_{r_{\epsilon ,j,t}}(w)\big)
  - \tau_{B_t}(r_{\epsilon ,j,t})
  & = \tau_{B_t}\Big[ \exp_q\Big(w{|q|\over \rho_t}\Big)\Big]
    +o(1) \cr
  & = {|w|^2\over 2} +o(1) \, . \cr
  }
$$
This is the result, setting
$$
  \tilde{\tau}_{B_t}\big(\exp_{r_{\epsilon,j,t}}(w)\big) 
  = {|q_{\epsilon ,j,t}|^2-\rho_t^2\over 2}+{|w|^2\over 2} \, .
$$

The second statement follows since we proved the asymptotic
equivalence $|\D I(r_{\epsilon ,j,t})|\sim\salt$. \hfill$\qed$

\bigskip

Combining Lemma 8.2.8 and result (8.2.6) yields the asymptotic
equivalence given in Theorem 8.2.1. It remains to check
the assumptions of Theorem 5.1.

Our choice of $\calD_{B_t}$ as a discrete set ensures that (5.1)
holds.

We now need a candidate for $c_{B_t,M}$. Let $c(t)=2d\log\log
t$. From Proposition 2.1 and our calculation of $I(B_t)$, 
we infer that
$$
  L\big( I(B_t)+c(t)\big)
  \leq {1\over t^{\alpha /2}(\log t)^{d-{1\over 2}}}\, O(1)
  = o(t^{-\alpha/2})
  \qquad\hbox{ as } \ttoi \, .
$$
Thus, $c(t)$ is a good candidate for $c_{B_t,M}$, no matter what
$M$ is. 

From the proof of Proposition 8.2.4 and Lemma 8.2.8, we infer
that for any $t$ large enough,
$$\displaylines{
  \qquad 
  \underline{B}_{t,M}\subset \big\{\, 
  \exp_{r_{\epsilon ,j,t}}(w) : |w|\leq \sqrt{5d\log\log t}
  \, ;\,
  \cut
  w\in T_{r_{\epsilon ,j,t}}\Lambda_{I(B_t)} \, ,
  \epsilon\in\{\, -1,1\,\}\, , \, j\in J_1 \,\big\} \, .
  \hfil\hskip 60pt\llap{(8.2.7)}\cr}$$
Since the level set $\Lambda_{I(B_t)}$ is a sphere of radius
$\sqrt{2I(B_t)}\big( 1+o(1)\big)\sim\salt$ as $t$ tends to infinity,
its radius of injectivity is of order $\sqrt{\log t}\gg
\sqrt{\log\log t}$, and assumption (5.2) holds.

Assumption (5.4) holds thanks to our choice of $c(t)$. 

To verify (5.5), let $r$ be a point in $\underline{B}_{t,M}$. 
The point $q=\psi\big( r,\tau_{B_t}(r)\big)$ is in the boundary 
$\partial B_t$. Since $\tau_{B_t}(r)$ is less than 
$c(t)=2d\log\log t$, the proof of Proposition
8.2.4 and (8.2.5) show that $q=q_{\epsilon ,j,t}(v)$ for some
$\epsilon$ in $\{\, -1,1\,\}$, some $j$ in $J_1$ and
$$
  \sum_{\scriptstyle 1\leq i\leq d\atop \scriptstyle i\neq j}
  \pisa (v_i)^2\leq 7d\log\log t \, .
$$
Consequently, Lemma A.1.5 shows that
$|v_i|\leq (\log t)^{2d/\alpha}$, for $1\leq i\leq d$
with $i\neq j$.
Since $v$ is in $\{\, (C+C^\T )e_i\,\}^\perp$, and this 
hyperplane does not contain $e_j$, we have 
$|v|=O(\log t)^{2d/\alpha}$. Notice that
$$
  \psi \big( r,\tau_{B_t}(r)+s\big)
  = \psi (q,s)
  = \sqrt{ 1+{2s\over |q|^2}} q\, ,
$$
and that $|q|\geq |r|=\rho_t$ tends to infinity with $t$. To
prove that $\chi_{B_t}^F\geq \rho_t^2$ for instance --- which is
more than enough to guarantee (5.5) --- it suffices to prove that
$\lambda q$ is in $B_t$ for any $1\leq \lambda\leq 2$. This is plain
from Lemma 8.2.2 and Claim 8.2.3.

Assumption (5.6) is plain.

To check (5.7), notice that for $q$ in $\calD_{B_t}$, the inequality
$t_{0,M}(q)\leq\sqrt{5d\log\log t}$ holds thanks to (8.2.7).
Furthermore, as $\D^2I(q)/|\D I(q)|$ equals $\Id/|q|$,
$$
  K_{\rm max}(q,t_0)
  \leq \sup\big\{\, |q|^{-2}  : q\in\Lambda_{I(B_t)}\,\big\}
  \sim {1\over 2I(B_t)}
  \sim {1\over \alpha\log t}
$$
as $t$ tends to infinity. Therefore, (5.7) holds.

Assumption (5.8) holds as well since $\pi^{-1}_{B_t}(p)$ is
essentially a finite union of spherical caps, and the Ricci 
curvature of a sphere is positive.

Assumption (5.9) is trivially satisfied since $\Lambda_c$ is a
sphere. Thus, two points in $\Lambda_c$ have equal norms.

It is no harder to verify (5.10), since
$$
  {\|\D^2I(p)\|\over |\D I(p)|^2}
  = { \|\Id\|\over |p|^2}
  = {1\over |p|^2} \, .
$$

To check assumption (5.11), again, we have 
$\| \D^2I\|=\| \Id\|=1$. Moreover, if $q$ belongs to
$\underline{B}_{t,M}\subset\Lambda_{I(B_t)}$ and $u$ is
nonnegative,
$$
  \big| \D I\big( \psi_u(q)\big) \big|^2
  \geq |\D I(q)|^2
  \sim 2 I(B_t)
  \sim \alpha\log t \, .
$$
Consequently, for $q$ in $\underline{B}_{t,M}$, that is 
$\tau_{B_t}(q)$ is less than $c(t)$, and for $t$ large enough,
$$
  \int_0^{\tau_{B_t}(q)} {\| \D^2I\|\over |\D I|^2}
  \big( \psi_u(q)\big) \d u 
  \leq {2c(t)\over \alpha\log t}
  = {4d\log\log t \over \alpha\log t} \, .
$$
Since $\calD_{B_t}$ is discrete, (5.12) holds automatically.

In conclusion, all the assumptions of Theorem 5.1 are
satisfied. This proves Theorem
8.2.1.\hfill$\qed$

\bigskip

From the work done, we can easily infer the following
conditional result.

\bigskip

\state{8.2.9. THEOREM.} {\it%
  Under the assumption of Theorem 8.2.1, if $J_1$ is nonempty, 
  the conditional distribution of the vector 
  $$
    (\log\sqrtt )^{-1} 
    {\big( \sign (X_i)\log |X_i|\big)}_{1\leq i\leq d}
  $$
  given $\cXX\geq t$ converges weakly* to 
  $$
    \sum_{j\in J_1}
    { C_{j,j}^{\alpha/2}
      \over \sum_{i\in J_1}{C_{i,i}^{\alpha /2}} }
    {\delta_{-e_j}+\delta_{e_j}\over 2} \, .
  $$}

\bigskip

\stateit{Proof.} In order to apply Corollary 5.3, let us check
its assumptions. Set $\lambda_{B_t}=\rho_t$. The numerator of the
measure involved in (5.18) is
$$
  \sum_{\epsilon\in\{ -1,1\}}\sum_{j\in J_1}
  { \exp \big( -\tau_{B_t}(r_{\epsilon ,j,t})\big)
    \delta_{r_{\epsilon ,j,t}/\rho_t} \over
    |\D I(r_{\epsilon,j,t})|^{(d+1)/2}
    \big( \det\, G_{B_t}(r_{\epsilon ,j,t})\big)^{1/2} } \, .
$$
We already calculated
$$
  \tau_{B_t}(r_{\epsilon ,j,t})
  = I(q_{\epsilon ,j,t})-I(B_t)
  = {\alpha\over 2}\log {\gamma_1\over C_{j,j}} + o(1)
  \qquad \hbox{ as } \ttoi \, .
$$
Moreover, as $r_{\epsilon ,j,t}$ is in the sphere
$\Lambda_{I(B_t)}$ and $\D I=\Id$, we have $|\D I(r_{\epsilon
,j,t})|=\rho_t$. Lemma 8.2.8 gives the value of $\det\,
G_{B_t}(r_{\epsilon ,j,t})$. Moreover, $r_{\epsilon ,j,t}/\rho_t
= q_{\epsilon ,j,t}/|q_{\epsilon ,j,t}| =\epsilon e_j$.
Consequently, the measure in (5.18) is
$$
  { \sum_{\epsilon\in\{ -1,1\}}\sum_{j\in J_1}
    C_{j,j}^{\alpha/2}\big(1+o(1)\big) \delta_{\epsilon e_j} \over
    2\sum_{j\in J_1}  C_{j,j}^{\alpha/2} \big( 1+o(1)\big) }
  \, .
$$
It certainly converges weakly* to
$$
  \nu =
  { \sum_{j\in J_1} C_{j,j}^{\alpha /2}
    (\delta_{-e_j}+\delta_{e_j})\over
    2\sum_{j\in J_1}C_{j,j}^{\alpha/2} } \, .
$$
Assumption (5.19) is trivial to verify since
$\underline{B}_{t,M}$ is on the sphere $\Lambda_{I(B_t)}$ of
radius $\rho_t$. Consequently, (5.19) becomes
$$
  \lim_{\ttoi} {2d\log\log t\over \rho_t^2}
  = \lim_{\ttoi} {2d\log\log t\over \alpha \log t} 
  = 0 \,  .
$$
To check (5.20) is as simple. The inclusion (8.2.7) shows that
if $q$ belongs to $\underline{B}_{t,M}$; then the Riemannian 
distance on
$\Lambda_{I(B_t)}$ between $q$ and $\pi_{B_t}$ is at most
$\sqrt{5d\log\log t}$. Since $\Lambda_{I(B_t)}$ is a sphere of
radius $\rho_t$, simple trigonometry shows that
$|q-\pi_{B_t}(q)|\leq \rho_t\sin (\sqrt{5d\log\log t}/\rho_t)$.
Since $\rho_t$ is of order $\sqrt{\alpha\log t}$ as $t$ tends to
infinity, assumption (5.20) is fulfilled.

Applying Corollary 5.2, the distribution of $Y/\sqrt{\alpha\log t}$
given $Y\in B_t$ converges weakly* to $\nu$. In other words, the
distribution of $\pisa (X)/\sqrt{\alpha\log t}$ given $\< CX,X\>\geq
t$ converges weakly* to $\nu$.

To rephrase this conclusion directly on $X$, we can use the Skorokhod
(1956) representation theorem. It implies the existence of a random
variable $Y_t$ having the same distribution as $Y$ given $Y\in B_t$,
and a random variable $Y_\infty$ having distribution $\nu$ such that
$Y_t/\sqrt{\alpha\log t}$ converges almost surely to $Y_\infty$.
Thus, $X$ given $\< CX,X\>\geq t$ has the same distribution as
$$
  \Saip (Y_t\sqrt{\alpha\log t})
  = \Saip \Big( Y_\infty \sqrt{\alpha\log t}\big( 1+o(1)\big)
  \Big) \, .
$$
Since $\Saip$ is ultimately sign preserving on $\RR$ and 
Lemma A.1.6 yields
$$\displaylines{\qquad
  \log \Big| \Saip \Big( \epsilon e_j\sqrt{\alpha\log t}
  \big( 1+o(1)\big)\Big)\Big|
  \hfill\cr\noalign{\vskip .08in}\hfill
  \eqalign{
   ={} & e_j\log \Saip\Big(\sqrt{\alpha \log t}
     \big( 1+o(1)\big)\Big) \cr
   ={} & e_j{\log t\over 2}\big( 1+o(1)\big) \, ,\cr
  }\qquad
  \cr}
$$
the result follows\hfill$\qed$

\bigskip

A careful sharpening of all the estimates could certainly lead
to more precise information on the conditional distribution of
$X$ given $\< CX,X\>\geq t$, and even an asymptotic expansion of
this conditional distribution. We will not pursue in that
direction for mainly two reasons: First, such calculation would
be quite specific to this example. Second, we will see hereafter
in this section that the kind of degeneracy at the limit --- the
limiting distribution is concentrated on a finite number of
points --- is not due to a bad rescaling but only to the fact that
$J_1$ is nonempty.

\bigskip

It may happen that all the diagonal coefficients of the matrix $C$ are
nonpositive, that is $\max_{1\leq i\leq d} C_{i,i}\leq 0$. What is the
analogue of Theorem 8.2.1 then? When $\gamma_1$ was positive, we could
essentially set all the $u_i$'s but one equal to $0$ in order to
optimize $I\big( \pisa (\sqrtt u)\big)$ --- see the proof of
Proposition 8.2.4 and inequality (8.2.4). When $\gamma_1$ is negative,
we need to take at least two components $u_i$, $u_j$ to be
nonzero. Setting $p_{i,j}= \sqrtt\, u_ie_i+\sqrtt\, u_je_j$, the
equation $p_{i,j}\in\partial A_t$ becomes
$$
  1 = u_i^2C_{i,i}+ u_iu_j (C_{i,j}+C_{j,i}) + u_j^2 C_{j,j} 
  \, .
$$
This equation admits a solution in $u_i$, $u_j$ if and only if
$$
  (C_{i,j}+C_{j,i})^2 - 4C_{i,i}C_{j,j} > 0 \, .
$$
Consequently, if
$$
  J_2=\big\{\, (i,j) : i\neq j\, ,\, 
  (C_{i,j}+C_{j,i})^2-4C_{i,j}C_{j,i} >0 \,\big\} 
$$
is nonempty, we can indeed consider only two nonzero 
components $u_i$, $u_j$ with $(i,j)$ in $J_2$.

How many components do we need to consider in general?
To answer this question, it is more convenient to change 
the notation. 

Let\ntn{$T$} $T$ be the set of all subsets of $\{\, 1,2\ldots ,d\,\}$.
For a set\ntn{$\calI$} $\calI=\{\, i_1,\ldots ,i_k\,\}$ with distinct
elements, denote by $|\calI |=k$ its cardinality. To $\calI$, we
associate the subspace\ntn{$V_\calI$} 
$V_\calI = \span \{\, e_{i_1},\ldots , e_{i_k}\,\}$ 
of dimension $|\calI|$. To the matrix $C$, we associate\ntn{$N(C)$}
$$
  N(C)= \min\big\{\, |\calI | : \calI\in T\, , \, 
  \exists u\in V_\calI \, , \, \< Cu,u\> > 0\,\big\}
$$
and\ntn{$J(C)$} 
$$
  J(C)=\big\{\, \calI \in T : \exists u\in V_\calI \, ,\,
  \< Cu,u\> > 0 \, , \, |\calI |= N(C)\,\big\}  \, .
$$
So, if $N(C)=1$, the set $J(C)$ is $J_1$. The integer $N(C)$ 
is the smallest cardinal of a set $\calI$
such that the inequation $\< Cu,u\>\geq 0$ has a solution in
$V_\calI\setminus\{\, 0\,\}$.
We exclude some degeneracy, assuming that
\setbox1=\vtop{\hsize=3.5in\noindent%
  for any $\calI$ in $T$ of cardinal $|\cal I|<N(C)$,
  the matrix $C$ is negative on $V_\calI$,}
\finetune{\kern 5pt}
$$
  \raise 7pt\box1\eqno{(8.2.8)}
$$
that is $\< Cu,u\>$ is negative for any nonzero vector $u$ in $V_\calI$.
This typically prevents having $\gamma_1$ null, and the analogue when
more components need to be considered. Notice that this
nondegeneracy is typical with respect to the matrix $C$. 
In particular, the equation $\<Cu,u\>= 1$ has a solution in 
any subspace $V_\calI$ for $\calI$ in $J(C)$, and has no 
solution in any subspace of the form $\span \{\,
e_{\alpha_1},\ldots , e_{\alpha_k}\,\big\}$ for any $k<N(C)$. 

We keep the notation\ntn{$A_t$}
$$
  A_t=\{\, x\in\RR^d : \< Cx,x\>=t\,\}
  = \sqrt{t}A_1 \, .
$$
The sets\ntn{$M_\calI$}
$$
  M_\calI = \big\{\, m\in V_\calI  : \< Cm\, ,m\> = 1\,\big\}
  \, , \qquad \calI\in J(C)\, ,
$$
are $( |\calI |-1)$-dimensional submanifolds
of $\RR^d$ and $\partial A_1$ as well. Indeed, it suffices to
prove that $1$ is a regular value of the map $x\in
V_\calI\mapsto \< Cx,x\>$. The differential of this map at $m$
is $\cct m$. If $m$ belongs to $M_\calI$,
$$
  \big\< \cct m,m\big\> = 2\< Cm,m\> = 2\neq 0 \, .
$$
Consequently, $\cct m$ does not vanish, or, equivalently, the
differential $\cct m$ is of full rank, and $1$ is a regular
value.

The result is then as follows.

\bigskip

\state{8.2.10. THEOREM.} {\it%
  Let $X$ be a $d$-dimensional random vector with
  independent and identically distributed components with
  Student-like distribution with parameter $\alpha$. Let $C$ 
  be a $d\times d$ matrix and $N=N(C)$. For $\calI$ in $J(C)$
  and $m$ in $M_\calI$, denote\ntn{$\tilde{G}(m)$}
  $\tilde{G}(m)$ the compression of the diagonal matrix
  $\sum_{i\in\calI}e_i\otimes e_i/|m_i|$ to 
  $\big\{\, \cct m \,\big\}^\perp \cap V_\calI$.
  Under (8.2.8), for $\alpha >2/N$,
  $$
    P\{\, \cXX \geq t\,\} \sim
    { K_{s,\alpha}^N\alpha^{\alpha N/2}\over t^{\alpha N/2}
      \sqrt{\alpha N} }
    \sum_{\calI\in J(C)}\int_{M_\calI} 
    { \det\,\tilde{G}(m) \over
      \prod_{i\in\calI} |m_i|^\alpha }
    \d\calM_{M_\calI} (m)
  $$
  as $t$ tends to infinity.
  }

\bigskip

\state{REMARK.} The assumption $\alpha >2/N$ guarantees that 
the integral over $M_\calI$ in the equivalence is finite. But we will
see that the result is true whenever the integral over $M_\calI$ is
finite. It is not clear whether $\alpha >2/N$ is required, though
$\alpha$ too small makes the integral diverge. This can be seen in
Lemma 8.2.18 below.  If $N(C)=1$, Theorem 8.2.10 is exactly Theorem
8.2.1. Clearly, in this case, the term
$K_{s,\alpha}^N\alpha^{\alpha N/2} t^{-\alpha N/2}/\sqrt{\alpha
N}$ in Theorem 8.2.10 gives the term $K_{s,\alpha}
\alpha^{(\alpha-1)/2} t^{-\alpha /2}$ in Theorem 8.2.1. For $N(C)=1$,
we have $J(C)=J_1$, provided we identify $\{\, i\,\}$ and $i$. If
$\calI= \{\, i\,\}$ belongs to $J(C)$, then $V_\calI= \RR\, e_i$, and
$$
  M_\calI
  = \big\{\, x\in \RR\, e_i : \cxx = 1 \,\big\}
  = \{\, -C_{i,i}^{-1/2}e_i,C_{i,i}^{-1/2}e_i\,\} \, .
$$
Thus, the Riemannian measure on $M_\calI$ is 
$$
  \calM_{M_\calI}=\delta_{-C_{i,i}^{-1/2}e_i}
  +\delta_{C_{i,i}^{-1/2}e_i} \, . 
$$
Moreover, for $m$ in $M_\calI$, the
matrix $\tilde{G}(m)$ is the compression of $e_i\otimes
e_i/m_i^2$ to $\big\{\, \cct m\,\big\}^\perp\cap V_\calI$. 
But, in our case, 
$$
  \big\{\, \cct m\,\big\}^\perp\cap V_\calI = \emptyset\, ,
$$
because, if $x$ is in $V_\calI$, we have $x=se_i$ for some
real number $s$, and
$$
  \big\< \cct m\, ,e_i\big\>
  = \pm C_{i,i}^{-1/2}\big\<\cct e_i\, ,e_i\big\>
  = \pm 2 C_{i,i}^{1/2}
  \neq 0
$$
--- the last inequality holds since we assume $N(C)=1$ here 
and $i$ in $J_1$, i.e., $\{\, i\,\}$ is in $J(C)$; we actually 
did the same work in the proof of claim 8.2.3. So, the term
$\det\,\tilde{G}(m)$ has to be omitted, and we obtain
$$
  \sum_{\calI\in J(C)}\int_{M_\calI} 
  {\d\calM_{M_\calI}\over \prod_{i\in\calI}|m_i|^\alpha}
  = 2\sum_{i\in J_1}C_{i,i}^{\alpha/2} \, ,
$$
as in Theorem 8.2.1.

\bigskip

\stateit{Proof of Theorem 8.2.10.} The proof is actually very
similar to that of Theorem 8.2.1, except that the dominating
manifold will no longer be a discrete set, and the
parameterizations will be slightly more sophisticated.

As in the proof of Theorem 8.2.1, we denote by $S_\alpha$ the
cumulative distribution function of $X_i$. We consider the
two sets
$$
  A_t=\{\, x\in\RR^d : \cxx \geq t\,\} 
  = \sqrtt A_1 \, , 
$$
and
$$
  B_t=\{\, y\in\RR^d : \< C\Saip (y)\, , \Saip (y)\> 
  \geq t\,\}  \, .
$$

Notice that if $\calI$ and $\calI'$ are distinct and in $J(C)$, then
$M_\calI$ does not intersect $M_{\calI'}$. Indeed, if $x$ is in both
$M_\calI$ and $M_{\calI'}$, then it is in 
$V_\calI\cap V_{\calI'}=V_{\calI\cap\calI'}$, 
and moreover $\< Cx,x\>=1$.  This contradicts the minimality of
$N(C)$, since $\calI$ different than $\calI'$ implies
$|\calI\cap\calI'|<N(C)$.

We first need a technical lemma, saying that whenever $\calI$
is in $J(C)$, the manifold $M_\calI$ stays away from 
any subspace $\{\, e_i\,\}^\perp$ with $i$ in $\calI$. 
We denote by $m_i=\< m,e_i\>$ the components of the vector 
$m$ belonging to $\RR^d$.

\bigskip

\state{8.2.11. LEMMA.} {\it%
  Under the assumption (8.2.8), there exists a positive 
  $\epsilon_0$ such that for all $\calI$ in $J(C)$,
  all $m$ in $M_\calI$, and all $i$ in $\calI$, 
  the inequality $|m_i|\geq\epsilon_0$ holds.  }

\bigskip

\stateit{Proof.} Searching for a contradiction, assume that there 
exists $\calI$ in $J(C)$, some index $i$ in $\calI$ 
and a sequence $m(n)$ in $M_\calI$ such that
$\lim_{n\to\infty}m_i(n)=0$. Write
$$
  m(n)=m_i(n)e_i+s(n)v(n)
$$
where $v(n)$ is a unit vector in 
$V_{\calI\setminus \{\, i\,\}} = V_\calI\ominus e_i\RR\,$, 
and $s(n)$ is a real number. Dropping the index $n$
for notational simplicity, the condition $m=m(n)$
belonging to $M_\calI$ becomes
$$
  1 = \< Cm\, ,m\>
  = m_i^2 \< Ce_i\, ,e_i\> + m_is\<\cct e_i\, ,v\> 
    + s^2\< Cv\, ,v\> \, .
$$
Assumption (8.2.8) guarantees 
$$
  \sup\big\{\, \< Cv,v\>  : v\in V_{\calI\setminus\{\, i\,\}}
  \, , \, |v|=1 \,\big\}
  <0 \, .
$$
Thus, the above quadratic equation in $s$ does not have any
solution, since its discriminant is $m_i^2\< \cct e_i,v\> -
4(m_i^2C_{i,i}-1)\< Cv,v\>= 4\< Cv,v\> + o(1)$, which is negative as
$n$ tends to infinity; this is a contradiction. \hfill$\qed$

\bigskip

In order to parameterize the boundaries $\partial A_t$ 
and $\partial B_t$, we consider the normal bundle of the 
immersion $M_\calI\subset \partial A_1$, namely,\ntn{$\calN_\calI$}
$$
  \calN_\calI = \big\{\, (m,v) : m\in M_\calI \, ,\, 
  v\in T_m\partial A_1\ominus T_mM_\calI\,\big\} \, .
$$
The analogue of the parameterization $p_{\epsilon ,j,t}(v)$ of
$\partial A_t$ in the proof of Theorem 8.2.1 is now a map
defined on $o(\sqrtt )$-sections of the normal bundle
$\calN_\calI$. 

It is convenient to introduce\ntn{$Q(t)$}
$$
  Q(t)=\salt - {\log\log \sqrtt\over 2\salt} - 
  {\log (K_{s,\alpha}\alpha^{\alpha/2}2\sqrt{\pi})\over
   \sqrt{\alpha\log t}} \, .
$$

\bigskip

\state{8.2.12. LEMMA.} {\it 
  Let $m$ be in $\MI$. The boundary $\partial A_t$ near $\sqrtt m$ can
  be parameterized as\ntn{$p_{\calI,t}(m,v)$}
  $$
    p_{\calI,t}(m,v)
    = \sqrtt m \Big( 1-{1\over 2t}\< Cv\, ,v\> + o(1/t)\Big) 
    + v \, ,
  $$
  $(m,v)\in\calN_\calI$, $|v|=o(\sqrtt )$ as $t$ tends to 
  infinity. The
  boundary $\partial B_t$ near $\pisa (\sqrt t m)$ can be 
  pa\-ram\-e\-terized as $q_{\calI,t} (m,v)+ o(\log t)^{-1/2}$ 
  where\ntn{$q_{\calI,t}(m,v)$}
  $$
    q_{\calI,t} (m,v) 
    = \sum_{i\in\calI} \sign (m_i)
    \bigg( Q(t)+{\sqrt{\alpha}\log |m_i|\over \sqrt{\log t}}
    \bigg) e_i
    + \sum_{\scriptstyle 1\leq i\leq d\atop\scriptstyle
            i\not\in\calI} \pisa (v_i) e_i
  $$
  and in the range $|v|=o(\sqrtt )$, 
  $\log |m_i|= o(\log t)^{1/2}$.
  }

\bigskip

\stateit{Proof.}  Since $\partial A_t=\sqrtt
\partial A_1$ and $\partial A_1$ is a manifold, there exists a
function $h$ and some small positive $\epsilon$ such that for any
$|v|\leq \epsilon\sqrtt$, any $m$ in $\MI$ and $v$ in
$T_m\partial A_1\ominus T_m\MI$,
$$
  p_{\calI,t}(m,v)
  =\sqrtt m \big( 1+h(v)\big) + v \in \partial A_t \, .
$$
This equation can be rewritten as
$$
  t\< Cm,m\> \big( 1+ h(v)\big)^2 + \sqrtt \big( 1+h(v)\big)
  \< \cct m\, ,v\> + \< Cv,v\> = t \, .
$$
Since $\< Cm,m\>=1$ and $\<\cct m,v\>= 0$ --- recall that $v$
belongs to $T_p\partial A_1=\{\, \cct m\,\}^\perp$ --- 
we obtain
$$
  0 = th(v)^2+ 2th(v) + \< Cv\, ,v\>
$$
as in the proof of Lemma 8.2.2. Thus, $h(v)\sim
-\< Cv\, ,v\>/(2t)$ as $t$ tends to infinity, uniformly 
in $|v|=o(\sqrtt )$. This gives the asymptotics for 
$p_{\calI,t}(m,v)$ in Lemma 8.2.12.

We pull back this parameterization to $\partial B_t$ by
introducing
$$
  \tilde{q}_{\calI,t} (m,v)
  = \pisa \big( p_{\calI,t} (m,v)\big)\, .
$$ 
Lemma A.1.5 implies
$$
  \pisa \Big(\sqrtt m_i\big( 1+|v|^2O(t^{-1})\big)+v_i\Big)
  = Q(t)+\sqrt\alpha\, {\log |m_i|\over \sqrt{\log t}}
  + o(\log t)^{-1/2}
$$
as $t$ tends to infinity, uniformly in the range $|v|=o(\sqrtt )$, $\log
|m_i|= o(\log t)^{1/2}$ and $|m_i|\geq \epsilon_0>0$. Since Lemma 8.2.11
ensures that $|m_i|$ stays away from $0$, the last condition,
$|m_i|\geq\epsilon_0$, may be omitted in the statement of Lemma
8.2.12.\hfill$\qed$

\bigskip

Define\ntn{$R(t)$}
$$\eqalign{
  R(t)={N(C)\over 2} \alpha\log t 
  & - {N(C)\over 2}\llst \cr
  & -N(C)\log (K_{s,\alpha}\alpha^{\alpha/2}2\sqrt{\pi})
    + \log (2\pi )^{d/2} \, .\cr
  }
$$
Recall that $I(x)=(|x|^2/2) +\log (2\pi )^{d/2}$. In the range
$|v|=o(\sqrtt )$ and $\log |m_i|=o(\log t)^{1/2}$,
it follows from Lemma 8.2.12 that
$$\eqalignno{
  I\big( q_{\calI,t}(m,v)\big)
  & = {1\over 2} \sum_{i\in\calI}\Big( Q(t)^2+ 2Q(t)
    {\sqrt{\alpha}\log |m_i|\over \sqrt{\log t}} + o(1)\Big) \cr
  & \qquad\qquad + {1\over 2} 
    \sum_{ \scriptstyle 1\leq i\leq d\atop\scriptstyle
           i\not\in\calI}\pisa (v_i)^2
    +\log (2\pi )^{d/2} \cr
  & \eqalign{ &= R(t) +\alpha\sum_{i\in\calI}\log |m_i|
                      +{1\over 2} 
                      \sum_{\scriptstyle 1\leq i\leq d\atop
                      \scriptstyle i\not\in I} \pisa (v_i)^2 \cr
                   \noalign{\vskip -.1in}
              &\qquad\qquad + o(1) \, \cr}
  &(8.2.9)\cr
  }$$
as $t$ tends to infinity, for all $i$ in $\calI$.

For a set $\calI$ belonging to $J(C)$, define\ntn{$\gamma_\calI$}
$$
  \gamma_\calI
  =\min\Big\{\, \prod_{i\in\calI}|m_i| : m\in\MI\,\Big\}
  \, .
$$
Furthermore, set\ntn{$\gamma$}
$$
  \gamma=\inf\big\{\, \gamma_\calI  : \calI\in J(C)\,\big\} \, .
$$
The following result locates the points on $\partial B_t$ where
the function $I$ is nearly minimal.

\bigskip

\state{8.2.13. LEMMA.} {\it%
  We have
  $$
    I(B_t)=R(t) + \alpha \log \gamma + o(1)
    \qquad \hbox{ as } \ttoi \, .
  $$
  Moreover, there exists $\eta$ in $(0,1/2)$ such that for any 
  number positive $M_1$, any $M_2>M_1/\alpha$, and any $t$ 
  large enough, the set
  $$
    \partial B_t\cap \Gamma_{I(B_t)+M_1\log\log t}
  $$
  is included in 
  $$\displaylines{\qquad
    \bigcup_{\smash{\calI\in J(C)}} \big\{\, q_{\calI,t}(m,v)  
    : m\in M_\calI\, ,\, |m|\leq (\log t)^{M_2} \, ,
    \cut
    |v|\leq t^{(1/2)-\eta} \, ,\, 
    |\Proj_{V_\calI^\perp}v|\leq (\log t)^{M_2} \,\big\}
    \qquad\cr}
  $$
  where $\Proj_{V_\calI^\perp}$ is the projection onto
  $V_\calI^\perp$.
  }

\bigskip

\stateit{Proof.} We first prove the second assertion
of the lemma. By construction, $\pisa (\partial\sqrtt A_1) =\partial
B_t$. Let $u$ be a point in $\partial A_1$. Lemma A.1.5 implies that
for any $\eta$ in $(0,1/2)$ and provided $t$ is large enough,
$$
  I\big( \pisa (\sqrtt u)\big)
  \geq \alpha (1-\eta )\sum_{1\leq i\leq d} 
  \indic_{[t^{-\eta},\infty )}\big( |u_i|\big) 
  \log \big( \sqrtt |u_i|\big) \, .
$$
Consequently, if 
$$
  I\big( \pisa (\sqrtt u)\big) \leq R(t)+M_1\log\log t + O(1) \, ,
$$
we must have, for $t$ large enough,
$$
  (1+\eta){N(C)\over 2} \alpha \log t
  \geq \alpha\, \sharp 
  \big\{\, 1\leq i\leq d  : |u_i|\geq t^{-\eta}\,\big\} 
  \log (t^{(1/2)-\eta}) \, .
$$
Taking $\eta$ positive and small enough so that the integer part of
${1+\eta\over 1-2\eta}N(C)$ is $N(C)$, the previous inequality yields
$$
  N(C)\geq \sharp \big\{\, 1\leq i\leq d : |u_i|\geq t^{-\eta}
  \,\} \, .
$$
Then, since $u$ is in $\partial A_1$, the minimality of $N(C)$ and Lemma
8.2.11 implies that we must have
$$
  N(C)=\sharp \big\{\, 1\leq i\leq d : |u_i|\geq t^{-\eta}
  \,\big\} \, .
$$
Thus, $\sqrtt u$ is in a $t^{(1/2)-\eta}= o(\sqrtt
)$-neighborhood of $\sqrtt\bigcup_{\calI\in J(C)}\MI$.
Therefore, it can be written as $p_{\calI,t}(m,v)$ for some 
$\calI$ in $J(C)$, some $m$ in $\MI$ and $|v|=o(t^{(1/2)-\eta})$. 

Lemma 8.2.11 ensures that $m_i\geq \epsilon_0$ wherever $i$ is in some 
set $\calI$ of $J(C)$. Thus, for any positive $\epsilon$ and $t$ 
large enough, Lemma A.1.5 implies
\setbox100=\hbox{$\displaystyle I\big( \pisa ( p_{\calI,t}) (m,v)\big)$}
\wd100=0pt
$$\eqalignno{
  \box100\hskip .5in & \cr
  \noalign{\vskip .1in}
  & \geq {1\over 2} \sum_{i\in\calI} \pisa 
    \big(\sqrtt |m_i| (1-\epsilon)\big)^2 \cr
  & \geq \alpha \sum_{i\in\calI} \log \big( \sqrtt |m_i|\big)
    -\log\log \big( \sqrtt |m_i|\big) + o(1) \cr
  & = \Big( {N(C)\over 2} \alpha\log t
    + \alpha \sum_{i\in\calI}\log |m_i|\Big) \big( 1+o(1)\big)
    \qquad\qquad
    &(8.2.10)\cr
  }$$
Therefore, the inequality
$$
  I\Big( \pisa \big(p_{\calI,t}(m,v)\big)\Big) 
  \leq R(t) + M_1\log\log t + O(1)
$$
implies $\alpha\max_{i\in\calI}\log |m_i|\leq M_1\log\log t
\big( 1+o(1)\big)$ for $t$ large enough. So, for $t$ large 
enough, $\log |m|\leq M_2\log\log t$ for $M_2$ larger than
$M_1/\alpha$.  Notice that $M_2$ can be chosen independently of $m$ in
$M_\calI$ since all the bounds are uniform in such $m$'s.  This proves
the second assertion of Lemma 8.2.13, except for the restriction
$|\Proj_{V_\calI^\perp}v|\leq (\log t)^{M_2}$.

In the range obtained so far, $I\big( q_{\calI,t}(m,v)\big)$
has minimal value $R(t)$ up to $o(1)$ as $t$ tends 
to infinity. Hence, $I(B_t)$ is as claimed. Finally,
on this range, if $|\Proj_{V_\calI^\perp}v|\geq (\log t)^{M_2}$,
then at least one component $v_j$, for some $j$ not in $\calI$,
is larger than $(\log t)^{M_2}/\sqrt{d}$. Consequently, as $t$
tends to infinity,
\finetune{\hfuzz=4pt}
$$\eqalign{
  I\big( \pisa (p_{\calI,t} (m,v)\big)
  & \geq R(t) + \alpha \log \gamma 
    + \sum_{\scriptstyle 1\leq j\leq d\atop\scriptstyle  
      j\not\in I}{1\over 2} \pisa (v_i)^2 + o(1) \cr
  & \geq R(t) + \alpha \log \gamma + \alpha M_2
    \log\log t \big( 1+o(1)\big) \, ,\cr
  }$$
\finetune{\hfuzz=0pt}
thanks to Lemma A.1.5. So, provided $M_2$ is large enough, the 
condition $q_{\calI,t} (m,v)$ belonging to $\Gamma_{I(B_t)+M_1\log\log
t}$ imposes to have $|\Proj_{V_\calI^\perp}v|$ less 
than $(\log t)^{M_2}$.\hfill$\qed$

\bigskip

As we did in the proof of Theorem 8.2.1, we can start to apply
Theorem 5.1 in calculating the terms of the asymptotic formula.

How do we choose the dominating manifold?

Looking at the parameterization of $\partial B_t$ by 
$q_{\calI,t}(m,v)$
in Lemma 8.2.12, we see that variations of $\prod_{i\in\calI}|m_i|$
of order $1$ yield small variations of order $1/\log t$ on
$q_\calI(m,v)$, i.e., in term of the Euclidean distance in $\RR^d$, 
while, according to (8.2.9), they give variations of order
$1$ in term of $I$. This suggests that on the
dominating manifold we should have $\prod_{i\in\calI}|m_i|$ constant.
Once $m$ is fixed, (8.2.9) shows that fluctuations in space of
order $1$ in $v$ yields fluctuations of order $1$ on $I$ as well.
Keeping in mind that $\partial B_t$ is of order $\pisa (\sqrtt
)\sim \sqrt{\alpha\log t}$ as far as its size is concerned, we see that
small fluctuations in $v$ --- on the scale of $\log t$ --- brings
sizable fluctuations of $I$. So, we should have $v$ constant in
the dominating manifold; and (8.2.9) suggests $v=0$. This leads us to
consider the set
$$
  \big\{\, m\in M_\calI  : \prod_{i\in \calI}|m_i|=\gamma_\calI
  \,\big\}
$$
of all points in $\MI$ which minimize the product of their
nonvanishing components. And then, one can try to choose the image
of $\sqrtt$ times this set by $\pisa$ as the dominating manifold, 
i.e.,
$$
  \big\{\, \pisa (\sqrtt m)  : m\in \MI \, , \,
  \prod_{i\in \calI}|m_i|=\gamma_\calI \,\big\} \, .
$$
This is not quite right since it does not belong to
$\Lambda_{I(B_t)}$, but a projection would fix this
detail. For some reason that the author does not quite understand
--- can someone give an explanation? --- this does not work and breaks
down when looking for a quadratic approximation of $\tau_{B_t}$.
The right manifold to consider seems to be a projection of\ntn{$\calD'_{B_t}$}
$$
  \calD'_{B_t}=\big\{\, q_\calI (m,0) : m\in \MI \, ,\, 
  \max_{i\in\calI}\log |m_i|\leq (\log t)^{1/4}\, , \,
  I\in J(C)\,\big\}
$$
on $\Lambda_{I(B_t)}$. The condition $\log |m_i|\leq (\log
t)^{1/4}$ in the definition of $\calD'_{B_t}$ guarantees that
$\log |m_i|=o(\log t)^{1/2}$ and will allow us to use Lemma
8.2.12 and equality (8.2.9).
Since $\Lambda_{I(B_t)}$ is a ball of radius\ntn{$\rho_t$}
$\rho_t=\sqrt{2\big( I(B_t)-\log (2\pi )^{d/2}\big)}$, 
this leads us to introduce\ntn{$r_{\calI,t}(m,v)$}
$$
  r_{\calI,t}(m,v)
  =\rho_t{q_{\calI,t}(m,v)\over |q_{\calI,t}(m,v)|} \, ,
$$
and\ntn{$\calD(B,t)$}
$$
  \calD_{B_t}=\Big\{\, r_{\calI,t}(m,0) : m\in\MI \, , \, 
  \max_{i\in\calI}\log |m_i|\leq (\log t)^{1/4}\, , \,
  \calI\in J(C)\,\Big\}  \, .
$$
Again, the projection on the sphere that we are doing is
actually a mapping through the normal flow to the level set
$\Lambda_{I(B_t)}$. In more general situations, we would define
$r_{\calI,t}(m,v)=\psi_{B_t}\big(q_{\calI,t}(m,v),s\big)$ with
$s=I(B_t)-I\big(q_{\calI,t}(m,0)\big)$.

Let us agree on the notation\ntn{$q_{\calI,t}(m)$}\ntn{$r_{\calI,t}(m)$}
$$
  q_{\calI,t}(m)=q_{\calI,t}(m,0)
  \qquad\hbox{ and } \qquad
  r_{\calI,t}(m)=r_{\calI,t}(m,0) \, .
$$
Define\ntn{$p(m)$}
$$
  p(m)=\sqrt{\alpha}\sum_{i\in\calI}\sign (m_i)\log |m_i|e_i
  \, , \qquad m\in \MI \, , \, \calI\in J(C) \, .
$$
For any $\calI$ in $J(C)$ and any $m$ in $\MI$ with $\log |m|=o(\log
t)^{1/2}$, we have
$$
  q_{\calI,t} (m)=Q(t) \sum_{i\in\calI}
  \sign (m_i)e_i + {p(m)\over \sqrt{\log t}}
$$
as can be seen from Lemma 8.2.12. From Lemma 8.2.11, we infer
that the vector $\sum_{i\in\calI} \sign (m_i)e_i$ is constant on
each connected component of $\bigcup_{\calI\in J(C)}\MI$. On each
of these connected components, we can think of the set of all
$q_{\calI,t} (m)$ --- i.e., the connected components of 
$\calD_{B_t}'$ --- as a translation by a fixed vector of 
length $\sqrt{|\calI |}Q(t)$ and a rescaling by 
$1/\sqrt{\log t}$ of the corresponding connected component 
of the set
$$
  \calD_\calI = \big\{\, p(m) : m\in \MI\,\big\}\, ,
  \qquad \calI\in J(C) \, .
$$

Following what we did in the proof of Theorem 8.2.1, we can start
to apply Theorem 5.1. The value of $I(B_t)$ is available from Lemma
8.2.13.

For $m$ in $M_\calI$, the differential  $p_*(m)$ is 
the restriction to $T_mM_\calI
=\big\{\, \cct m\,\big\}^\perp\cap V_\calI$ of the matrix
$$
  \sqrt{\alpha}\,\diag {(1/|m_i|)}_{i\in\calI}
  =\sqrt{\alpha}\sum_{i\in\calI}e_i\otimes e_i/|m_i| \, .
$$ 
The change of variable $m\leftrightarrow q_{\calI,t}(m)$ gives
$$
  \d\calM_{\calD_{B_t}'} = \det \big( p_*^\T p_*^{\phantom\T}\big)^{1/2}
  \d\calM_{\cup_{\calI\in J(C)}M_I} (\log t)^{(1-N(C))/2}\, ,
$$
provided we integrate on the restricted range $\log |m|=o(\log
t)^{1/2}$. Since $r_{\calI,t}(m)$ is reasonably close to
$q_{\calI,t}(m)$, we should be able to approximate the measure
$\d\calM_{\calD_{B_t}}$ by $\d\calM_{\calD'_{B_t}}$
when applying Theorem 5.1.

From equation (8.2.9), we infer that
$$
  \tau_{B_t}\big( q_{\calI,t} (m)\big)
  = \alpha \sum_{i\in\calI}\log |m_i|-\alpha\log\gamma + o(1)
  \quad \hbox{ as } \ttoi \, .
$$

To obtain the needed quadratic approximation for $\tau_{B_t}$ along
the geodesic leaves orthogonal to $T_r\calD_{B_t}$, we first need to
determine $T_r\Lambda_{I(B_t)}\ominus T_r\calD_{B_t}$. What should
it be? Let us argue informally. Later, we will make the argument
rigorous. First, we should be able to replace
$T_r\calD_{B_t}$ by $T_{q_{\calI,t}(m)}\calD_{B_t}'$ for some 
$\calI$ in $J(C)$ and $m$ in $\MI$. Since
$T_{q_{\calI,t}(m)}\Lambda_{I(B_t)}$ is the orthocomplement
of $q_{\calI,t}(m)\RR\,$, we
should expect $T_{q_{\calI,t}(m)}\Lambda_{I(B_t)}\ominus
T_{q_{\calI,t}(m)}\calD_{B_t}'$ to be roughly 
$\{\, q_{\calI,t} (m)\,\}^\perp\ominus 
T_{q_{\calI,t}(m)}\calD_{B_t}'$. Given the 
expression for $q_{\calI,t}(m)$, we have
$$\eqalign{
  T_{q_{\calI,t} (m)}\calD_{B_t}'
  & = T_{p(m)}\calD_\calI = p_*(m)T_m\MI \cr
  \noalign{\vskip 1mm}
  & = \Big\{\, \sum_{i\in\calI} {v_i\over |m_i|} e_i  : 
    v\in V_\calI\ominus \cct m\RR\,\Big\}
    \subset V_\calI \, . \cr
  }$$
Intuitively, our idea of projecting $\calD_{B_t}'$ on
$\Lambda_{I(B_t)}$ to obtain $\calD_{B_t}$ will work well if the
projection does not create a singularity or reduce the dimension;
in other words, if $\calD_{B_t}'$ is transverse to the
direction of $q_{\calI,t}(m)=Q(t)\sum_{i\in\calI}\sign (m_i)e_i
\big(1+o(1)\big)$, or, roughly, if $\sum_{i\in\calI}\sign
(m_i)e_i$ does not belong to $T_{p(m)}\calD_\calI$ for all 
$m$ in $\MI$. Since this latter vector is in $\MI$, our
calculation of $T_{q_{\calI,t}(m)}\calD'_{B_t}$ allows
to rewrite this condition as
$$
  \big\< \sum_{i\in\calI}\sign (m_i)|m_i|e_i\, ,\cct m\big\> 
  \neq 0 \, .
$$
It holds since $\sum_{i\in\calI}\sign(m_i)|m_i|e_i=m$ and
$\<(C+C^\T )m,m\>\geq 1$.
In this case, since $q_{\calI,t}(m)\RR$ is approximately
$\sum_{i\in\calI}\sign (m_i)e_i\RR\subset V_\calI$, 
the subspace $V_\calI$ is
almost spanned by $q_{\calI,t}(m)\RR$ and $T_{p(m)}\calD_\calI$. 
So, approximately,
$$
  \big\{\, q_{\calI,t} (m)\,\big\}^\perp\ominus 
  T_{q_{\calI,t}(m)}\calD_{B_t}'\approx V_\calI^\perp \, .
$$
Then, the expression of $q_{\calI,t}$ in Lemma 8.2.12 shows 
that in $V_\calI^\perp+q_{\calI,t}(m)\RR\,$, the boundary 
$\partial B_t$ is approximately the 
$\big(d-N(C)\big)$-dimensional affine space
$V_\calI^\perp+q_{\calI,t}(m,0)$ --- consider $\pisa (v_i)$ 
as a new coordinate, say, $w_i$. But then, working on
$V_\calI^\perp+q_{\calI,t} (m)\RR\,$, the very same argument as in the
proof of Theorem 8.2.1 shows that
$\tau_{B_t}\big(\exp_{q_{\calI,t}(m)}(w)\big)$
is approximately $|w|^2/2$ for $w$ in the normal bundle. So, we
should have
$$
  G_{B_t}\big(q_{\calI,t}(m)\big)
  \approx {\Id_{V_\calI^\perp}\over |q_{\calI,t}(m)|}
  = {\Id_{\RR^{d-N(C)}}\over Q(t)\sqrt{N(C)}} \, .
$$
Thus, as $t$ tends to infinity,
$$\eqalign{
  \det\, G_{B_t}\big(q_{\calI,t}(m)\big)^{-1/2}
  &\sim Q(t)^{(d-N(C))/2} N(C)^{(d-N(C))/4} \cr
  &\sim \big( \alpha N(C)\log t\big)^{(d-N(C))/4}
  \, .\cr
  }
$$
We can then put all the pieces together, find a way to drop
the restriction $\log |m_i|\leq (\log t)^{1/4}$ in the range
of integration on $\calD_{B_t}$, and use the formula given
in Theorem 5.1 to obtain the --- hypothetical --- approximation
$$\displaylines{
  P\big\{\, \cXX\geq t\,\big\}
  \sim { K_{s,\alpha}^{N(C)}\alpha^{\alpha N(C)/2} \over
         \sqrt{N\alpha}\, t^{\alpha N(C)/2} }
  \times
  \cut
  \int_{\cup_{\calI\in J(C)}\MI} 
  { \det \big( p_*(m)^\T p_*(m)\big)^{1/2} \over 
    \prod_{i\in\calI} |m_i|^\alpha } \, 
  \d\calM_{\cup_{\calI\in J(C)}M_\calI}(m) \, ,\cr}
$$
which is the result.

Let us now work out the proper arguments.

We first determine our candidate for $c_{B_t}$. Let 
$$
  c(t)= \Big(d+{N(C)\over 2}+1\Big)\log\log t \, . 
$$
From Proposition 2.1 and Lemma 8.2.13, we infer that
$$\eqalignno{
  L\big( I(B_t)+c(t)\big)
  & = O(1)e^{-I(B_t)}e^{-c(t)} I(B_t)^d \cr
  \noalign{\vskip 2pt}
  & = t^{-\alpha N(C)/2} (\log t)^{d+(N(C)/2)} 
    (\log t)^{-d- (N(C)/2)-1} O(1) \cr
  \noalign{\vskip 2pt}
  & = o(t^{-\alpha N(C)/2}) \, .
      &(8.2.11)\cr
  }$$
Thus, $c(t)$ is a good candidate for an upper bound of $c_{B_t}$
--- notice again that $c(t)$ is of order $\log I(B_t)$. Many assumptions
in Theorem 5.1 deal with the behavior of $\tau_{B_t}$ near the
dominating manifold. Our goal is now to more or less calculate the
value of $\tau_{B_t}$ where we need it.

The following lemma will give us the normal bundle
of $\calD_{B_t}$ immersed in $\Lambda_{I(B_t)}$.

\bigskip

\state{8.2.14. LEMMA.} {\it%
  Let $\calI$ be in $J(C)$.
  Let $m$ be a point in $\MI$ with
  $\max_{i\in\calI}\log |m_i|=o(\log t)^{1/2}$.
  Then, $T_{r_{\calI,t}(m)}\calD_{B_t}+r_{\calI,t}(m)\RR
  = V_\calI$. Consequently,
  $$
    T_{r_{\calI,t}(m)}\Lambda_{I(B_t)}\ominus T_{r_{\calI,t}(m)}
    \calD_{B_t} = V_\calI^\perp \, .
  $$}
\finetune{}%\bigskip

\stateit{Proof.} The map $u\in\RR^d\setminus\{\, 0\,\}
\mapsto u/|u|\in S_{d-1}$ has differential 
$|u|^{-1}\Proj_{u^\perp}$ at $u$. Consequently,
$$\displaylines{\quad
  T_{r_{\calI,t}(m)}\calD_{B_t}
  \cut
  \eqalign{
  ={} & {\rho_t\over |q_{\calI,t}(m)|}\Proj_{q_{\calI,t}(m)^\perp}
    T_{q_{\calI,t}(m)}\calD_{B_t}' \cr
  ={} & \Proj_{q_{\calI,t}(m)^\perp}\big\{\, q_{\calI,t,*}(m)v  : 
    v\in T_m\MI \,\big\} \cr
  ={} & \Proj_{q_{\calI,t}(m)^\perp} \Big\{\, \sum_{i\in\calI}
    {v_i\over |m_i|}e_i : \< \cct m\, ,v\>=0\,\Big\} 
    \subset V_\calI \, . \cr
  }
  \quad\cr}$$
Consequently, $T_{r_{\calI,t}(m)}\calD_{B_t}$ is of
dimension 
$$
  \dim \big( \{\, \cct m\,\}^\perp \cap V_\calI\big)
  = \dim V_\calI -1 = N(C)-1
$$ 
and $T_{r_{\calI,t}(m)}\calD_{B_t}+q_{\calI,t}(m)\RR = V_\calI$.
Since $r_{\calI,t}(m)$ is collinear to $q_{\calI,t}(m)$, this
gives the first statement of the lemma.

The last statement of Lemma 8.2.14 follows from the fact that
the tangent space of $\Lambda_{I(B_t)}$ at $r_{\calI,t}(m)$ is
orthogonal to $q_{\calI,t} (m)$ since the level sets are
spheres.\hfill$\qed$

\bigskip

As we did in the proof of Theorem 8.2.1 with Lemma 8.2.7,
we now relate $\tau_{B_t}$ on $\calD_{B_t}$ to
$\tau_{B_t}$ on $\calD_{B_t}'$. Since $\calD_{B_t}'$ is somewhat
more explicit than $\calD_{B_t}$, and better parametrized, this
will make further calculations easier. Using
Lemma 8.2.7, we see that for any $w$ orthogonal to $r_{t,I}(m)$,
$$\eqalign{
  \tau_{B_t}\big( \exp_{r_{\calI,t}(m)}(w)\big)
  & = \tau_{B_t}
   \Big[ { \rho_t\over |q_{\calI,t} (m)| }
         \exp_{q_{\calI,t}(m)}
         \Big( w{|q_{\calI,t}(m)|\over\rho_t}\Big)
   \Big] \cr
  & = {|q_{\calI,t}(m)|^2\over 2} 
   \Big( 1-{\rho_t^2\over |q_{\calI,t}(m)|^2} \Big) \cr
  & \qquad\qquad
    + \tau_{B_t}\Big[ \exp_{q_{\calI,t} (m)}
    \Big( w { |q_{\calI,t} (m)|\over\rho_t } \Big)\Big]
    \cr
  }\eqno{(8.2.12)}$$
In particular, since $\exp_q(0)=q$,
$$\displaylines{\qquad
     \tau_{B_t}\big( \exp_{r_{\calI,t}(m)}(w)\big)
      - \tau_{B_t}\big( r_{\calI,t}(m)\big)
  \cut
  \eqalign{
     =\,& \tau_{B_t}\Big[\exp_{q_{\calI,t}(m)}\Big( 
        w {|q_{\calI,t}(m)|\over \rho_t}\Big)\Big]
        - \tau_{B_t}\big( q_{\calI,t}(m)\big) \cr
     =\,& \tau_{B_t}\Big[\exp_{q_{\calI,t}(m)}\Big( 
        w {|q_{\calI,t}(m)|\over \rho_t}\Big)\Big] + o(1) \, ,
     \cr}
   \qquad\cr}
$$
the last equality coming from the parameterization of $\partial
B_t$ in Lemma 8.2.12. We can now prove the analogue of 
Lemma 8.2.8.

\bigskip

\state{8.2.15. LEMMA.} {\it%
  For $m$ in $\calD_{B_t}$ and 
  $w$ in $V_\calI^\perp = T_{r_{\calI,t}(m)}
  \Lambda_{I(B_t)}\ominus T_{r_{\calI,t}(m)}\calD_{B_t}$
  in the range $|w|=o(\sqrt{\log t})$, we have
  $$
    \tau_{B_t}\big( \exp_{r_{\calI,t}(m)}(w)\big)
    - \tau_{B_t}\big( r_{\calI,t}(m)\big)
    = {|w|^2\over 2}+o(1) \qquad \hbox{ as } \ttoi \, .
  $$
  Consequently, 
  $$
    G_{B_t}\big( r_{\calI,t}(m)\big) 
    = { \Id_{\RR^{d-N(C)}}
        \over \big| \D I\big( r_{\calI,t}(m)\big)\big| }
    \sim {\Id_{\RR^{d-N(C)}}\over \sqrt{N(C)\alpha\log t} }
    \qquad\hbox{ as }\ttoi \, .
  $$}

\bigskip

\stateit{Proof.} Let $w$ be a vector orthogonal
to $V_\calI$. From the calculation preceding Lemma 8.2.15, 
we see that it is enough to evaluate $s=
\tau_{B_t}(z)$ for $z=\exp_{q_{\calI,t}(m)}
\big( {w|q_{\calI,t}(m)|/\rho_t}\big)$. By definition, 
$\psi(z,s)$ is in the boundary of $B_t$. Lemma
8.2.12 shows that on directions orthogonal to $V_\calI$, 
the boundary $\partial B_t$ behaves like a
$\big(d-N(C)\big)$-dimensional 
linear subspace --- at least in the range $w=\pisa (v)$ with
$v=o(\sqrt{t})$, that is $|w|=o(\log t)^{1/2}$. So, 
looking on the 
components on $V_\calI$, we must have for $(m,v)$ in $\calN$,
$$
  \big\< \psi (z,s)\, ,e_i\big\>
  = \big\< q_{\calI,t} (m,v),e_i\big\>+o(\log t)^{-1/2} 
  \, , \qquad i\in\calI \, .
$$
To evaluate the left hand side term of the above equality, we
use
$$\displaylines{
    \psi (z,s)
  \cut
    =\sqrt{ 1+{2s\over |q_{\calI,t}(m)|^2} }
    \Big[ \cos \Big( {|w|\over\rho_t}\Big) q_{\calI,t}(m)
    + \sin \Big({|w|\over \rho_t}\Big) \big| q_{\calI,t}(m)\big|
    {w\over |w|}\Big]  \, .
  \cr}
$$
To evaluate the right hand side, we have
$$
  \< q_{\calI,t}(m,v)\, ,e_i\> 
  = \< q_{\calI,t} (m,0)\, ,e_i\> + o(\log t)^{-1/2} \, .
$$
Then, we obtain
$$
  s= {|q_{\calI,t}(m,s)|^2\over 2} 
  \Big( {1+o(\log t)^{-1}\over \cos^2\big( |w|/\rho_t\big) }
        -1\Big)
  = {|w|^2\over 2} +o(1) 
$$
in the range $|w|=o(\rho_t)$, i.e., $|w|=o(\log t)^{1/2}$ as
announced. The matrix $G_{B_t}\big( r_{\calI,t}(m)\big)$ in
the statement of Lemma 8.2.15 is that corresponding to the
assumption (5.17) as weakened at the end of 
Remark 5.2.\hfill$\qed$

\bigskip

We are now in position to verify that the assumptions of 
Theorem 5.1 hold.

It will be helpful to keep in mind two orders of magnitudes.
Since $\calD_{B_t}$ is contained in 
$\Lambda_{I(B_t)}=S_{d-1}(0,\rho_t)$,
points in $\calD_{B_t}$ are of order $\rho_t$, that is of
order $\sqrt{\log t}$. On the
other hand, Lemmas 8.2.12, 8.2.13 and A.1.5 show that any point of
$\partial B_t\cap \Gamma_{I(B_t)+c(t)}$ is at a distance at most
$O(\log\log t)$ of a point $q_{\calI,t}(m)$. Indeed, if
$|\Proj_{V_\calI^\perp}v|\leq (\log t)^{M_2}$, the component of
$q_{\calI,t}(m,v)$ on $V_\calI^\perp$ is $\sum_{1\leq i\leq d
\, ; \, i\not\in\calI}\pisa(v_i)e_i$, which is of order $O(\log\log t)$
thanks to Lemma A.1.5. Since the projection onto the sphere
$\Lambda_{I(B_t)}$ is a Lipschitz function when acting on
$\Gamma_{I(B_t)}^{\rm c}$,
the points in $\underline{B}_{t,M}$ are also
at a distance $O(\log\log t)$ of $\calD_{B_t}$.

Assumption (5.1) holds trivially since $k=\dim M_\calI= N-1$. 

To check (5.2) amounts to proving that any point
$s$ on $\Lambda_{I(B_t)}$, in a $O(\log\log t)$-neighborhood of
$\calD_{B_t}$, can be
written in a unique way as $\exp_{r_{\calI,t}(m)}(w)$ for some
$\calI$ in $J(C)$, some $m$ in $M_\calI$ and
$$
  w\in T_{r_{\calI,t}(m)}\Lambda_{I(B_t)}\ominus
  T_{r_{\calI,t}(m)}\calD_{B_t}\equiv V_\calI^\perp \, . 
$$
Lemma 8.2.11 implies that such a point $s$ is actually 
in an $O(\log\log t)$-neighborhood of a unique 
set $\calD_{B_t}\cap V_\calI = \big\{\,r_{\calI,t}(m) : m\in\MI\,\big\}$ 
for some $\calI$ in $J(C)$.  For $w$ belonging 
to $T_{r_{\calI,t}(m)}\Lambda_{I(B_t)}\ominus 
T_{r_{\calI,t}(m)}\calD_{B_t}$, we have
$$
  \Proj_{V_\calI^\perp}\exp_{r_{\calI,t}(m)}(w)
  = \cos\Big( {|w|\over \rho_t}\Big) r_{\calI,t}(m) \, .
$$
Consequently, the component of $s$ on $V_\calI$ is in a one-to-one
relation with $\cos\big( |w|/\rho_t\big)r_{\calI,t}(m)$. This
last point being a positive multiple of $r_{\calI,t}(m)\in
S_{d-1}(0,\rho_t)$, it identifies $r_{\calI,t}(m)$ and consequently
$m$. Looking at the component of $s$ in $V_\calI^\perp$, we can
then calculate $w$ in a unique way, and (5.2) holds.

Our choice of $c(t)$ satisfying (8.2.11) will imply (5.4) ultimately, while (5.6) is
plain.

Assumption (5.5) can be verified exactly as in the proof of
Theorem 8.2.1. Indeed, Lemmas 8.2.12 and 8.2.13 show that
$\partial B_t\cap\Gamma_{I(B_t)+c(t)}$ can be approximated by a
ruled hypersurface based on $\calD_{B_t}'$ where the generators
are $(d-N)$-dimensional Euclidean balls of radius $O(\log\log
t)$.

Assumption (5.7) is verified exactly in the same way as in the
proof of Theorem 8.2.1. The bound $t_{0,M}(p)= O(\log\log t)$
follows from Lemma 8.2.15 ---  or from the discussion at the
beginning of this assumptions checklist, after the proof of Lemma
8.2.15. Assumptions (5.8)--(5.11) are obtained in the very same
way as we did in the proof of Theorem 8.2.1.

Checking (5.12) requires some more work. Let $m$ be in $\MI$. 
Consider a curve $m(s)$ on $\MI$, such that $m(0)=m$. The curve
$r_{\calI,t}\big(m(s)\big)$ lies on $\calD_{B_t}$. In
Lemma 8.2.14 we proved that  $V_\calI^\perp$ can be 
identified with
$T_{r_{\calI,t}m(s)}\Lambda_{I(B_t)}\ominus
T_{r_{\calI,t}(m(s))}\calD_{B_t}$ for all small $s$.
Hence, for any $w$ in the unit sphere of $V_{\calI}^\perp$ and 
$\lambda$ positive, the curve $\exp_{r_{\calI,t}(
m(s))}(\lambda w)$ on $\Lambda_{I(B_t)}$ is well defined. Its
tangent vector field at $s=0$ is given by
$$\displaylines{\qquad
  {\d\over\d s}\exp_{r_{\calI,t}(m(s))}(\lambda w)\Big|_{s=0}
  \hfill\cr\noalign{\vskip .05in}\kern 1in
    = {\d\over \d s}\bigg( 
    \cos\Big( {\lambda |w|\over\rho_t}\Big)
              r_{\calI,t}\big( m(s)\big) 
    + \sin\Big( {\lambda |w|\over\rho_t}\Big) \rho_t{w\over |w|}
    \bigg)\bigg|_{s=0}
  \hfill\cr\noalign{\vskip .05in}\kern 1in
    = \cos\Big( {\lambda |w|\over\rho_t}\Big) r_{\calI,t,*}(m)
    m'(0) \, .
    \hfill(8.2.13)\cr
  }$$
Let $\tilde p=\exp_{r_{\calI,t}(m(s))}(\lambda w)$. 
Since $\pi_{B_t}(\tilde p)=r_{\calI,t}\big(m(s)\big)$
provided $\tilde p$ stays in $\omega_{B_t,r_{\calI,t}(m)}$,
and since the orthocomplement of $\ker\pi_{B_t,*}(\tilde p)$ 
has dimension $k$, this orthocomplement can be identified as
$$
  \big(\ker\pi_{B_t,*}(\tilde p)\big)^\perp
  \equiv r_{\calI,t,*}(m)T_mM_\calI
  =T_{r_{\calI,t}(m)}\calD_{B_t} \, .
$$ 
Moreover, (8.2.13) implies that
for $u$ in $T_{r_{\calI,t}(m)}\calD_{B_t}$,
$$
  \pi_{B_t,*}(p)u= u\Big/\cos\Big( {\lambda |w|\over \rho}\Big) \, .
$$
In other words, the restriction of $\pi_{B_t,*}(p)$ to
$\big(\ker \pi_{B_t,*}(p)\big)^\perp$ is the map
$$
  \Id_{\RR^k}\Big / \cos\Big({\lambda |w|\over\rho}\Big) \, .
$$
Consequently, $J\pi_{B_t}(p)=\cos\big( \lambda |w|/\rho)^{-k}
= 1+o(1)$ uniformly in the range $\lambda \leq O(\log\log
t)=o(\rho_t)$ and $m$ in $\calD_{B_t}$. This proves (5.12).

Before checking (5.13), we need to evaluate the candidate for
the limiting integral, namely
$$
  \int_{\calD_{B_t}} 
  { e^{-\tau_{B_t}}\indic_{[0,c(t)]}(\tau_{B_t})
    \over |\D I|^{d-k+1\over 2} (\det\, G_{B_t})^{1/2} } 
  \d\calM_{\calD_{B_t}}
$$
--- notice that we use the Riemannian measure on $\calD_{B_t}$ 
and not that on $\calD_{B_t}'$. Since $r_{\calI,t}(m)$ is in the
image of $q_{\calI,t}(m)$ by the map $u\in\RR^d\setminus \{\,
0\,\}\mapsto \rho_t u/|u|$ whose differential at $q_{\calI,t}(m)$
is $\big(\rho_t/|q_{\calI,t}(m)|\big)
\Proj_{q_{\calI,t}(m)^\perp}$, we have
$$
  \d\calM_{\calD_{B_t}}\big(r_{\calI,t}(m)\big)
  = \d\calM_{\calD_{B_t}'}\big( q_{\calI,t}(m)\big)
  \big( 1+o(1)\big) \, ,
  \eqno{(8.2.14)}
$$
in the range $|q_{\calI,t}(m)|\sim\rho_t$. This range includes
that for which $\tau_{B_t}\big( r_{\calI,t}(m)\big)$ is
less than $c(t)$, thanks to Lemmas 8.2.12 and 8.2.13. 
Furthermore,
$$\eqalign{
  \tau_{B_t}\big( r_{\calI,t}(m)\big)
  & = \tau_{B_t}\big( q_{\calI,t}(m)\big)
    + {1\over 2}\Big( \big| q_{\calI,t} (m)\big|^2-\rho_t^2\Big)
    \cr
  & = I\big( q_{\calI,t}(m)\big) - I(B_t) \cr
  & =\alpha\sum_{i\in\calI}\log |m_i| -\alpha \log \gamma + o(1)
    \qquad \hbox{ as } \ttoi \cr
  }$$
uniformly over $\tau_{B_t}\big( r_{t,\calI}(m)\big) \leq c(t)$ --- the
first equality comes from (8.2.12) with $w=0$; the second from the
fact that $q_{\calI,t}(m)$ belongs to $\partial B_t$, and so
$\tau_{B_t}\big( q_{\calI,t}(m)\big)=0$; the third from (8.2.9) and
Lemma 8.2.13.

Uniformly in $\tau_{B_t}\big( r_{\calI,t}(m)\big) \leq c(t)$, we
have 
$$
  \big| \D I\big( r_{\calI,t}(m)\big)
  \sim \big| \D I\big( q_{\calI,t}(m)\big)\big|
  \sim \sqrt{N\alpha\log t} \, .
$$
Consequently,
\finetune{\hfuzz=1pt}
$$\displaylines{
    \int_{\calD_{B_t}} {e^{-\tau_{B_t}}\indic_{[0,c(t)]}
    (\tau_{B_t}) \over |\D I|^{(d-k+1)/2}
    (\det\, G_{B_t})^{1/2} } \d\calM_{\calD_{B_t}}
  \cut
  \sim  \int_{\calD_{B_t}'}
         { \gamma^\alpha \, \indic_{[0,c(t)]}(\tau_{B_t})
           \over \prod_{i\in\calI}|m_i|^\alpha 
           (N\alpha\log t)^{(d-(N-1)+1)/4} 
           (N\alpha\log t)^{-{(d-N)/4}} }
         \d\calM_{D_{B_t}'} 
  \cut
  \sim  \sum_{\calI\in J(C)}\int_{\MI} 
         {\gamma^\alpha\over\prod_{i\in\calI}|m_i|^\alpha}
	 \det \big( p_*(m)^\T p_*(m)\big)^{1/2}
	 \indic_{[0,c(t)]}\Big( \tau_{B_t}\big(r_{t,\calI}(m)
	  \big)\Big)
   \cut
        \qquad\qquad\qquad \d\calM_{\MI}(m)
        {1\over (N\alpha\log t)^{1/2}}
        {1\over (\log t)^{(N-1)/2}} 
	\qquad\hbox{ as } \ttoi \, . 
  \cr
  }$$
\finetune{\hfuzz=0pt}
Thus, the formula in Theorem 5.1 becomes
$$\displaylines{\qquad
    e^{-I(B_t)} (2\pi )^{(d-k-1)/2} \int_{\calD_{B_t}}
    { e^{-\tau_{B_t}}\over |\D I|^{(d-k+1)/2} 
      (\det\, G_{B_t})^{1/2} }
    \d\calM_{\calD_{B_t}}
  \cut
  \sim  {(K_{s,\alpha}\alpha^{\alpha/2})^N \over \sqrt{N\alpha}}
	 {1\over t^{N\alpha/2}}
	 \sum_{\calI\in J(C)}\int_{M_\calI} 
	 {\det \big( p_*(m)^\T p_*(m)\big)^{1/2}\over
	  \prod_{i\in\calI}|m_i|^\alpha}
   \cut\qquad\qquad\qquad
	 \indic_{[0,c(t)]}\Big(\tau_{B_t}\big( r_{t,\calI}(m)
	 \big)\Big) \d\calM_{\MI}(m) \, . 
  \qquad\cr }
$$
This formula is valid whenever $\alpha$ is positive.
To remove the term 
$\indic_{[0,c(t)]}\Big(\tau_{B_t}\big(r_{\calI,t}(m)\big)\Big)$
from the formula, it is enough to prove that for all $\calI$ 
in $J(C)$,
$$
  \int_{\MI} 
  { \det\big( p_*(m)^\T p_*(m)\big)^{1/2}\over
    \prod_{i\in\calI}|m_i|^\alpha }
  \d\calM_{\MI}(m) < \infty \, ,
  \eqno{(8.2.15)}
$$
because $\lim_{\ttoi} c(t)=\infty$. This is where we need
the assumption $\alpha >N/2$. For this purpose, we
first derive a bound for $\det\big( p_*(m)^\T p_*(m)\big)$. The
following result will be useful. Since $M_\calI$ is in the
complement of a neighborhood of the origin, it extends Lemma
8.2.11.

\bigskip

\state{8.2.16. LEMMA.} {\it%
  There exists a constant $K$ depending only on the matrix
  $C$ such that for any $\calI$ in $ J(C)$, any $m$ in $\MI$ 
  and $i$ in $\calI$,
  $$
    {|m|\over K} \leq |m_i| \leq K|m| \, .
  $$}

\bigskip

\stateit{Proof.} Assume that the lower bound were false, and
that, for instance, $\inf_{m\in\MI}|m_i|/|m|=0$ for some $i$ in
$\calI$ and $\calI$ in $J(C)$. Then, there exists a sequence $m(n)$ in
$\MI$ such that $\lim_{n\to\infty}|m_i(n)|/|m(n)|=0$. Given Lemma
8.2.11, this forces $\lim_{n\to\infty}|m(n)|=\infty$. Set
$m^i(n)=m(n)-m(n) e_i=\Proj_{e_i^\perp}m(n)$. Define $s(n)=|m^i(n)|$
and $v(n)=m^i(n)/|m^i(n)|$.  Since $|m(n)|$ tends to infinity with
$n$, the condition $m_i(n)/|m(n)|=o(1)$, ensures that $|m^i(n)|$ does
not vanish for $n$ large enough.  Thus $v(n)$ is well defined for $n$
large enough. Moreover, $|m_i(n)|=o\big( s(n)\big)$ as $n$ tends to
infinity. We then have
\finetune{\hfuzz=1.4pt}
$$\eqalign{
  1 & = \< Cm(n)\, ,m(n)\> \cr
    & = m_i(n)^2C_{i,i} + m_i(n)s(n)\big\<\cct e_i\, ,v(n)\big\>
      + s(n)^2\< Cv(n)\, ,v(n)\> \cr
    & = o\big( s(n)\big)^2 + s(n)^2\< Cv(n)\, ,v(n)\> \, .
      \cr
  }$$
\finetune{\hfuzz=0pt}
Since $v(n)$ is in the compact sphere $S_{N-1}\subset
V_{\calI\setminus \{\, i\,\}}$, we can assume, up to extracting a
subsequence, that $v=\lim_{n\to\infty}v(n)$ exists. Then $v$
belongs to $V_{\calI\setminus \{\, i\,\}}$ and we must 
have $\< Cv,v\>\sim s(n)^{-2}=o(1)$, i.e., $\< Cv,v\>=0$. 
This contradicts assumption (8.2.8).

The upper bound is trivial since $|m_i|\leq |m|$ 
anyway.\hfill$\qed$

\bigskip

\state{8.2.17. LEMMA.} {\it%
  If $\calI$ is in $J(C)$ and $m$ belongs to $\MI$, then
  $$
    0\leq \det\big( p_*(m)^\T p_*(m)\big)
    \leq \alpha^{N-1}\Big( {K\over |m|}\Big)^{2(N-1)} \, .
  $$}

\bigskip

\stateit{Proof.}  From Lemma 8.2.16, we infer
$$
  \sqrt{\alpha}\sum_{i\in\calI}{e_i\otimes e_i\over |m_i|}
  \leq {\sqrt{\alpha}K\over |m|} \Id_{V_\calI} 
$$
in the sense that the difference between the right and
left hand sides is a nonnegative matrix.
Consequently, since restriction preserves matrix ordering,
$$
  p_*(m)^\T p_*(m) \leq {K^2\alpha\over |m|^2}
  \Id_{\{\, \cct m\,\}^\perp\cap V_\calI} \, .
$$
Taking the determinant preserves the ordering as well, and the result
follows. \hfill$\qed$

\bigskip

We can now prove (8.2.15).

\bigskip

\state{8.2.18. LEMMA.} {\it%
  If $\alpha > 2/N$, then
  $$
    \int_{\MI} 
    { \det\big( p_*(m)^\T p_*(m)\big)^{1/2}\over \prod_{i\in\calI}
      |m_i|^\alpha} \d\calM_{\MI}(m) < \infty \, .
  $$}

\bigskip

\stateit{Proof.} Using Lemmas 8.2.16 and 8.2.17,
it is enough to prove that 
$$
  \int_{\MI}{\d\calM_{\MI}(m)\over |m|^{N-1+N\alpha}} <\infty 
  \, .
$$
First, using a change of variable in polar coordinate, we see
that for any $\epsilon$ positive and $s>N$,
$$
  \int_{\scriptstyle |x|\geq\epsilon\atop\scriptstyle x\in\RR^N}
  {\d x\over |x|^s}
  < \infty \,  .
$$
Take $\epsilon <\epsilon_0$, where $\epsilon_0$ is as in Lemma
8.2.11. In $V_\calI$ identified with $\RR^N$, make the change of 
variable $x=\theta m$, with $\theta$ positive and $m$ in $M_\calI$. 
We claim that
$$
  \int_{\scriptstyle |x|\geq\epsilon\atop
        \scriptstyle x\in V_\calI} {\d x\over |x|^s}
  \geq \int_{\theta\geq 1}\int_{m\in\MI}
  {\big| \Proj_{(T_m\MI)^\perp} (m)\big| \over \theta^s|m|^s}
  \d\calM(m)\theta^{N-1}\d\theta \, .
$$
To obtain the right hand side and prove the claim, we argue as
follows. Let $m(t_1,\ldots , t_{N-1})$ be a parameterization of
$\MI$. The Jacobian of the 
transformation $x=\theta m(t_1,\ldots ,t_{N-1})$ is
$$
  \det\Big( {\partial x\over\partial t_1}\, , \ldots \, , 
  {\partial x\over\partial t_{N-1}} \, , 
  {\partial x\over\partial\theta}\Big)
  = \theta^{N-1}\det\Big( {\partial m\over \partial t_1} \, ,
  \ldots \, , {\partial m\over \partial t_{N-1}}\, , m\Big) \, .
$$
Using the multilinearity of the determinant, writing $m$ as a
component on $\span \Big\{\, 
{\displaystyle\partial m\over\displaystyle\partial t_1}
\, ,\ldots \, ,
{\displaystyle\partial m\over\displaystyle\partial t_{N-1}}
\,\Big\} = T_m\MI$ and a
component on $(T_m\MI)^\perp$, we obtain
$$
  \det \Big( {\partial m\over \partial t_1}\, ,
  \ldots \, , {\partial m\over \partial t_{N-1}} , m\Big)
  = \det\Big({\partial m\over \partial t_1}\, , \ldots \, ,
    {\partial m\over\partial t_{N-1}},
    \Proj_{(T_m\MI)^\perp}m \Big) \, .
$$
Up to a sign, this last determinant is
$$
  |\Proj_{(T_m\MI)^\perp}m| 
    \bigg( \det \Big( {\partial m\over \partial t_1}\, 
                      , \ldots \, ,
                      {\partial m\over\partial t_{N-1}}
                \Big)^\T 
                \Big({\partial m\over \partial t_1}\, 
                      , \ldots \, ,
                      {\partial m\over\partial t_{N-1}}
                \Big)\bigg)^{1/2} \, .
$$
This proves our claim. If $s>N$, we obtain, after
performing the integration in $\theta$,
$$
  \infty > \int_{m\in\MI} { \big|\Proj_{(T_m\MI)^\perp}(m)\big|
  \over |m|^s}\d\calM_{M_\calI} (m) \, .
$$
Recall that $(T_m\MI)^\perp\cap V_\calI=\Proj_{V_\calI}
\cct m\RR\,$. Since $\MI$ is included in $V_\calI$, it follows 
that for any $m$ belonging to $\MI$,
$$\eqalign{
  \Proj_{(T_m\MI)^\perp}(m)
  &= {\< m,\cct m\> \over |\Proj_{V_\calI}\cct m|^2} 
    \,\Proj_{V_\calI}\cct m \cr
  \vadjust{\vskip .05in}
  &= {2\over |\Proj_{V_\calI}\cct m|^2 } 
     \,\Proj_{V_\calI}\cct m \, . \cr
  }
$$
Thus,
$$
  |\Proj_{(T_m\MI)^\perp}(m)|
  \geq {2\over \|\Proj_{V_\calI}\cct\| \, |m|}
  \geq {1\over \| C\| \, |m|} \, . 
$$
Consequently, for $s>N$, the integral
$\displaystyle\int_{m\in\MI}{\displaystyle \d\calM(m)\over
\displaystyle |m|^{s+1}}$ is finite. This proves Lemma 8.2.18.
\hfill$\qed$

\bigskip

Putting all the pieces together, we are left with only
assumption (5.13) to check. But this is plain from the
calculation we did in Lemmas 8.2.16 and 8.2.18. 
This concludes the proof of Theorem 8.2.10.\hfill$\qed$

\bigskip

We can now state an analogue of Theorem 8.2.9, that is a result
on the limiting behavior of $X$ given $\cXX\geq t$ as $t$ tends
to infinity. For $\calI$ in $J(C)$, Lemma 8.2.11 
implies that the vector $\big(\sign (m_i)\big)_{i\in \calI}$ 
is constant on the connected components of $\MI$. Hence, 
to a connected component $\calC$
of $M_\calI$, we can associate a unique unit vector\ntn{$\epsilon_\calC$}
$$
  \epsilon_\calC=N^{-1/2}\sum_{i\in\calI}\sign (m_i)e_i\, ,
$$
where $m$ is any point in $\calC$. Moreover, for $\alpha >2/N$ 
and $\tilde G(m)$ as defined in Theorem 8.2.10, the number\ntn{$\mu_\calC$}
$$
  \mu_\calC =\int_\calC {\big( \det\, \tilde G(m)\big)^{1/2}
  \over \prod_{i\in\calI}|m_i|^\alpha } \d\calM_{M_\calI}(m)
$$
is finite.

\bigskip

\state{8.2.19. THEOREM.} {\it%
  Under the assumptions of Theorem 8.2.10, the distribution
  of the random vector 
  $\big( \sign (X_i)\log |X_i|\big)_{1\leq i\leq d}\big/
  \log\sqrtt$ given $\cXX\geq t$ converges weakly* to
  $$
    \nu =\sum_{\calC}\mu_\calC\delta_{\epsilon_\calC}
    \Big/ \sum_\calC \mu_\calC \, ,
  $$
  where the sums are over all connected components $\calC$ of
  $\bigcup_{\calI\in J(C)}M_\calI$.
  }

\bigskip

\stateit{Proof.} We proceed in a similar way as we did for
proving Theorem 8.2.9. Set $\lambda_{B_t}=\rho_t$. The numerator
of the measure in (5.18) when looking for the conditional
distribution of $Y/\rho_t$ given $Y\in B_t$ is
$$\displaylines{\qquad
  { \rho_t^{N-1}\exp \big( -\tau_{B_t}(\rho_tr)\big) \over
    \rho_t^{(d-(N-1)+1)/2}(N\alpha\log t)^{-(d-N)/4} }
    \,\d\calM_{\calD_{B_t}/\rho_t}(r)
  \cut
  \sim (N\alpha\log t)^{N/2} \exp\big( -\tau_{B_t}(\rho_tr)\big)
  \,\d\calM_{\calD_{B_t}/\rho_t}(r) \, .\qquad\cr
  }
$$
Let $f$ be a positive continuous function. Then (8.2.14) implies
$$\displaylines{\qquad
  \int f(r) \exp\big( -\tau_{B_t}(\rho_tr)\big)
  \,\d\calM_{\calD_{B_t}/\rho_t}(r)
  \cut
  \sim
  \int f(r)\exp\big( -\tau_{B_t}(\rho_tr)\big)
  \,\d\calM_{\calD'_{B_t}/\rho_t}(r)
  \qquad}
$$
as $t$ tends to infinity. Using the parameterization 
of $\calD_{B_t}$ and the expression for
$\tau_{B_t(r_{\calI,t}(m))}$ obtained from equation (8.2.14), 
this last integral is also asymptotically equivalent to
$$\displaylines{\qquad
  \sum_{\calI\in J(C)}\sum_{\calC}\int_\calC
  f\Big( {r_{\calI,t}(m)\over\rho_t}\Big)
  {\gamma^\alpha\over\prod_{i\in\calI}|m_i|^\alpha}\times
  \cut
  \det\big( p_*(m)^\T p_*(m)\big)^{1/2}
  \d\calM_{M_\calI}(m) \, ,
  \qquad (8.2.16)\cr}
$$
where the second summation is over all connected components $\calC$ of
$M_\calI$. But if $m$ is in a connected component $\calC$, then
$$
  \lim_{\ttoi} {r_{\calI,t}(m)\over \rho_t} = \epsilon_\calC \, ,
$$
uniformly in the range of $m$'s such that
$q_{\calI,t}(m)$ belongs to $\calD'_{B_t}$ --- i.e.,
$\max_{i\in\calI}\log |m_i|$ is less than $(\log t)^{1/4}$. 
Consequently, as $t$ tends to infinity, (8.2.16) converges to
$$
  \sum_{\calI\in J(C)}\sum_{\calC\in M_\calI}
  f(\epsilon_\calC) \gamma^\alpha\mu_\calC
  = \gamma^\alpha\sum_\calC \mu_\calC \int f \d\nu \, .
$$
Since the measure involved in (5.18) is normalized, we proved
that the probability measure with density proportional to
$$
  { \gamma^k\exp\big( -\tau_{B_t}(\rho_tq)\big) \over
    |\D I(\rho_tq)|^{d-k+1\over 2}
    \big(\det\, G_{B_t}(\rho_tq)\big)^{1\over 2} }
  \,\d\calM_{\calD_{B_t}/\rho_t}(q)
$$
converges weakly* to $\nu$. Thus, (5.18) holds.

Assumption (5.19) holds trivially. Assumption (5.20) follows
from the fact that $\rho_t$ is of order $\sqrt{\log t}$, while
for $p$ in $\calD_{B_t}$, the set $\pi_{B_t}^{-1}(p)\cap
\underline{B}_{t,M}$ is in an $O(\log\log t)$-neighborhood of
$\calD_{B_t}$ thanks to Lemma 8.2.13.

Applying Corollary 5.3, we deduce that the conditional
distribution of $Y/\sqrt{N\alpha\log t}$ given $Y\in B_t$
converges weakly* to $\nu$.

Using the Skorokhod representation Theorem as we did in
proving Theorem 8.2.9, Theorem 8.2.19 follows from
$$\displaylines{\qquad
  \Saip \Big( \epsilon_\calC \sqrt{N\alpha\log t}\big(1+o(1)\big)
    \Big)
  \cut
  \eqalign{
    ={} & \sqrt{N} \epsilon_\calC \Saip\Big( \sqrt{\alpha\log t}
      \big( 1+o(1)\big)\Big) \cr
    ={} & \sqrt{N}\epsilon_\calC \exp\Big(\log \sqrt{t}
      \big( 1+o(1)\big)\Big) \cr
    }
  \qquad\cr}
$$
as $t$ tends to infinity, thanks to A.1.6.\hfill$\qed$

\bigskip

One way to interpret Theorem 8.2.19 is to say that $X$ given
$\cXX\geq t$ is distributed as $t^{(1/2)+o(1)}S$ where
$S$ is distributed according to $\nu$. In other words, we
obtained something like the first term of an expansion of $X$ given
$\cXX\geq t$. Given the work done, a little extra effort will
give us a second order term. That is, given $\cXX\geq t$, and
given that $\big( \sign (X_i)\log |X_i|\big)_{1\leq i\leq
d}/\log\sqrtt$ is close to some $\epsilon_\calC$, we can derive
the limiting distribution of $X$. Refining the asymptotic
analysis, we could obtain an asymptotic expansion of the
distribution of $X$ given $\cXX\geq t$, as $t$ tends to infinity, 
in term of successive conditional distributions. Such type of result is
easier to phrase in term of random variables than in term of
distributions.

\bigskip

\state{8.2.20. THEOREM.} {\it%
  Under the assumptions of Theorem 8.2.10, the random variable
  $X$ given $\cXX\geq t$ can be represented as
  $\sqrt{t}m_gH\big( 1+o(1)\big) + T_g +o(1)$, where 
  the random variables $g$, $m_g$, $H$ and $t_g$ are 
  as follows.

  The discrete random variable $g$ has distribution $\nu$.

  Let $\calC$ be a connected component of some $\MI$.
  The conditional density of $m_g$ given 
  $g=\epsilon_\calC$ is proportional to 
  $$
    { |\Proj_{(T_{p(m)}p(\calC ){)}^\perp}\epsilon_\calC |
      \over \prod_{i\in\calI}|m_i|^\alpha} 
    \,\det\,\tilde G(m)\,\d\calM_\calC (m) \, .
  $$ 
  The random variable $H$ has a Pareto distribution 
  $$
    P\{ H\geq \lambda\}=\lambda^{-\alpha N} \, , 
    \qquad \lambda\geq 1\, ,
  $$
  and is independent of $g$ and $m_g$. Finally,
  $T_g$ given $g$ is a random vector in $V_\calI^\perp$,
  with independent components, all having the original
  Student-like $S_\alpha$ distribution.
  }

\bigskip

One way to read Theorem 8.2.20 is in terms of simulating $X$ from
its conditional distribution given
$\cXX\geq t$ for large $t$. We pick a connected component
$\calC$ with probability proportional to $\mu_\calC$. Once
$\calC$ is picked, it lies in a unique subspace $V_\calI$. In
$V_\calI$, we simulate $m$ with distribution proportional to 
$$
  { |\Proj_{(T_{p(m)}p(\calC ){)}^\perp}\epsilon_\calC |
    \over \prod_{i\in\calI}|m_i|^\alpha} \,\det\,\tilde G(m)\,\d\calM_\calC (m) \, .
$$ 
Next, we simulate $H$ with the Pareto distribution. Then the part of
$X$ in $V_\calI$ is $\sqrtt mH\big(1+o(1)\big)$ as $t$ tends to
infinity. The part of $X$ in $V_\calI^\perp$ is a random vector with
independent and identically distributed components from the initial
Student-like distribution, up to an additive term of order $o(1)$ as
$t$ tends to infinity. We will explicitly calculate the norm of
$\Proj_{T_{p(m)}}p(\calC )\epsilon_\calC$ at the end of the proof of
Theorem 8.2.20.

\bigskip

\stateit{Proof of Theorem 8.2.20.}
The intuition behind the proof is extremely
simple given all that we did. Looking at $Y=\pisa (X)$, we want
to obtain an approximation of $Y$ given $Y\in B_t$, and invert
it to obtain one of $X=\Saip (Y)$ given $\cXX\geq t$. Lemma
8.2.12 asserts that the points in $B_t$ near
$Q(t)\sqrt{N}\epsilon_\calC$ are of the form $\psi\big(
q_{\calI,t}(m,v),s\big)$. Set $w=\pisa (v)$, and define
$$
  \tilde q_{\calI,t}(m,w)=\sum_{i\in\cal I}\sign (m_i)
  \Big( Q(t)+{\sqrt\alpha\log |m_i|\over\sqrt{\log t}}\Big)
  e_i + w \, , \qquad w\in V_\calI^\perp \, .
$$
Formula (8.2.9), the definition of the normal flow, and Lemma
8.2.13 show that the normal density at $\psi\big(\tilde
q_{\calI,t}(m,w),s\big)$ is
$$\eqalignno{
  \exp\Big[-I\Big( \psi_s\big( \tilde q_{\calI,t}(m,w)\big)
  \Big)\Big] 
  & = \exp\Big[ -I\big(\tilde q_{\calI,t}(m,w)\big)-s\Big]
    \cr
  & = { e^{-R(t)}e^{-s}e^{-|w|^2/2}\over \prod_{i\in\calI}
        |m_i|^\alpha }\big( 1+o(1)\big) \kern 25pt
  &(8.2.17)\cr
  }
$$
In this expression, the term $e^{-s}$ is an exponential density,
the term $e^{-|w|^2/2}$ is a Gaussian one, and the term
$\prod_{i\in \calI}|m_i|^\alpha$ will give us another density. 
We can parameterize points $y$ of $B_t$ in term of $s,w,m$. Then, we
interpret these parameters as random variables. This will give
a representation of the random variable $Y$ given $Y\in B_t$,
and we will pull back this representation to $X$.

However, we need to be careful with the different scales. For
$m$ in a connected component $\calC$ of $M_\calI$,
$$
  \sqrt{\log t}\,\Proj_{V_\calI^\perp}\big( \tilde q_{\calI,t}(m,w)
  - Q(t)\sqrt{N}\epsilon_\calC\big) = p(m) \, ,
$$
while
$$
  \Proj_{V_\calI}\tilde q_{\calI,t}(m,w) = w \, .
$$
Thus, the part of $y$ in $V_\calI$ should be centered and
rescaled by $\sqrt{\log t}$, while that in $V_\calI^\perp$ 
is already of order $1$.

To proceed rigorously, let $f$ be a nonnegative smooth function
defined on $\RR^N\times\RR^{d-N}$, with compact support.

For some $\calI$ in $J(C)$ and a connected component $\calC$ of
$M_\calI$, let us evaluate the integral
$$
  \int_{B_t}f\Big( \sqrt{\log t}\,\Proj_{V_\calI} 
  \big( y-Q(t)\sqrt{N}\epsilon_\calC\big),\Proj_{V_\calI^\perp}
  y\Big) e^{-I(y)} \d y \, .
  \eqno{(8.2.18)}
$$
Divided by $P\{\, Y\in B_t\,\}$, this integral will give
us the limiting con\-di\-tion\-al behavior of $Y$ given $Y\in B_t$,
after proper nor\-mal\-iza\-tion.

If the projection of $y$ onto $V_\calI$ is not in a neighborhood of
$Q(t)\sqrt{N}\epsilon_\calC$, then
$\sqrt{\log t}\, \Proj_{V_\calI}\big( y-Q(t)\sqrt{N}\epsilon_\calC\big)$
diverges as $t$ tends to infinity. Since $f$ is compactly supported, such
points $y$ do not contribute to the integral for large $t$. From
the preceding, since $f$ is bounded, we can also restrict the
range of integration for those $y$'s such that $I(y)\leq I(B_t)+c(t)$. 
Let us make the change of variable
$$
  y= \psi\big( \tilde{q}_{\calI,t}(m,w),s\big)
  = \sqrt{1+{2s\over |\tilde{q}_{\calI,t}(m,w)|^2}}\,
  \tilde q_{\calI,t}(m,w) \, .
$$
Notice that
$$
  |\tilde q_{\calI,t}(m,w)|^2\geq |q_{\calI,t}(m)|^2
  \sim N\alpha\log t
$$
as $t$ tends to infinity, and uniformly in the 
range $\tau_{B_t}\big( \tilde q_{\calI,t}(m)\big)\leq c(t)$.
Therefore,
$$
  y=\tilde q_{\calI,t}(m,w) +{s\over \sqrt{N\alpha\log t}}
  \epsilon_\calC + o\Big({1\over \sqrt{\log t}}\Big)
$$
as $t$ tends to infinity, uniformly in the 
range $I(y)\leq I(B_t)+c(t)$. Consequently, in that range 
of $y$'s,
$$
  \sqrt{\log t}\,\Proj_{V_\calI}
  \big( y-Q(t)\sqrt{N}\epsilon_\calC\big)
  = p(m) + {s\over\sqrt{N\alpha}}\epsilon_\calC + o(1)
  \qquad\hbox{ as } \ttoi \, .
$$
Using (8.2.17), we can rewrite (8.2.18) as
$$\displaylines{
  \int_{s\geq 0}\int_{m\in M_\calI}\int_{w\in V_\calI^\perp}
  \indic_{[0,I(B_t)+c(t)]}\big( I(y)\big)
  f\Big( p(m)+{s\over\sqrt{N\alpha}}\epsilon_\calC,w\Big)\times
  \cut
  {e^{-R(t)}\over\prod_{i\in\calI}|m_i|^\alpha}
  e^{-s}e^{-|w|^2/2} J(m,w,s)
  \d\calM_{M_\calI}(m)\, \d w\, \d s \, \big( 1+o(1)\big) \, ,
  \qquad\cr
  }
$$
where $J(m,w,s)$ is a Jacobian term. To calculate it, we first
have
$$
  {\partial y\over\partial s}
  = { \tilde q_{\calI,t}(m,w)\over |\tilde q_{\calI,t}(m,w)|^2
      \sqrt{ 1+{ \displaystyle 2s\over\displaystyle 
                 |\tilde q_{\calI,t}(m,w)|^2} } }
  = { \epsilon_\calC \over \sqrt{N\alpha\log t} }
    \big( 1+o(1)\big)
$$
as $t$ tends to infinity, uniformly in $s$ in any compact 
set of $\RR$ and $m,w$ such that $I(y)\leq I(B_t)+c(t)$.

The explicit expression of $\psi_s(q)$ for the Gaussian
distribution gives
$$
  \psi_{s,*}(q)
  = \sqrt{ 1+{2s\over |q|^2}}\,\Id - {2s\over |q|^2}
  {q\otimes q\over |q|^2} 
  {1\over \sqrt{1+{\displaystyle 2s\over\displaystyle |q|^2}}}
  \, .
$$
Let $m=m(u_1,\ldots , u_{N-1})$ be a local parameterization of
$\calC$. Then
$$
  {\partial\tilde q_{\calI,t}\over \partial u_j}(m,w)
  = {1\over\sqrt{\log t}}\, p_*(m){\partial m\over \partial u_j}
  \, .
$$
Consequently
$$
  {\partial y\over\partial u_j}
  = \psi_{s,*}\big( \tilde q_{\calI,t}(m,w)\big)
    {\partial \tilde q_{\calI,t}(m,w)\over\partial u_j}
  = {1\over\sqrt{\log t}} p_*(m) 
    {\partial m\over \partial u_j} \big( 1+o(1)\big) \, ,
$$
uniformly in $y$ such that $I(y)\leq I(B_t)+c(t)$. Finally, since
${\displaystyle\partial \over \displaystyle\partial w_j}
\tilde q_{\calI,t}(m,w)=e_j$, we have
$$
  {\partial y\over\partial w_j}
  = \psi_{s,*}\big( \tilde q_{\calI,t}(m,w)\big) e_j
  = e_j + o(1) \qquad \hbox{ as } \ttoi \, .
$$
Define the $(d\times d)$-matrix with columns indexed by $i$ and
rows indexed by $j$,
$$
  F=\bordermatrix{
  &{\scriptstyle N-1} & {\scriptstyle 1} & {\scriptstyle d-N}
    \cr
  \noalign{\vskip 0.05in}
  {\scriptstyle N-1} &
    \Big\< {\dispar y\over\dispar u_i} ,{\dispar y\over\dispar u_j}
    \Big\>_{i,j}
  & \Big\< {\dispar y\over\dispar u_i} ,{\dispar y\over\dispar s}
    \Big\>_{i}
  & \Big\< {\dispar y\over\dispar u_i} ,{\dispar y\over\dispar w_j}
    \Big\>_{i,j}
  \cr\noalign{\vskip 0.05in}
  {\hbox{\kern 8pt}\scriptstyle 1} &
    \Big\< {\dispar y\over\dispar s} ,
         {\dispar y\over\dispar u_j}\Big\>_{j}
  & \Big|{\dispar y\over\dispar s}\Big|^2
  & \Big\< {\dispar y\over\dispar s} ,{\dispar y\over\dispar w_j}
    \Big\>_{j}
  \cr\noalign{\vskip 0.05in}
  {\scriptstyle d-N} &
    \Big\< {\dispar y\over\dispar w_i} ,{\dispar y\over\dispar u_j}
    \Big\>_{i,j}
  & \Big\< {\dispar y\over\dispar w_j} ,{\dispar y\over\dispar s}
    \Big\>_i
  & \Big\< {\dispar y\over\dispar w_i} ,{\dispar y\over\dispar w_j}
    \Big\>_{i,j} 
  \cr} \, .
$$
\vskip .1in
This ensures that $J(m,w,s)=(\det\, F)^{1/2}$.
Since $\partial y/\partial w_i$ is orthogonal to $V_\calI$, 
while $\partial y/\partial u_i$ and
$\partial y/\partial s$ are in $V_\calI$, and since
$\partial y/\partial w_j$ is roughly $e_j$, 
the determinant of
$F$ is equal to $o(\log t)^{-N}$ plus the determinant of
the upper left $N\times N$ block of $F$, that is
$$\det\pmatrix{
      {\displaystyle 1\over\displaystyle \log t}  
      \Big\< p_*(m){\dispar m\over\dispar u_i} , 
             p_*(m){\dispar m\over\dispar u_j}
      \Big\>_{i,j}
    & 
      {\displaystyle 1\over\displaystyle \sqrt{N\alpha}\log t}  
      \Big\< p_*(m){\dispar m\over\dispar u_i} ,
             \epsilon_\calC \Big\>_i 
    \cr\noalign{\vskip 0.05in}
      {\displaystyle 1\over \displaystyle\sqrt{N\alpha}\log t} 
      \Big\< \epsilon_\calC, p_*(m){\dispar m\over\dispar u_j}
      \Big\>_j
    &
      {\displaystyle 1\over\displaystyle N\alpha\log t} 
    \cr}
$$
--- in this determinant, $i,j$ run over $1,\ldots, N-1$.
This determinant is equal to
$$\displaylines{
     {1\over N\alpha(\log t)^N}\, \det
     \pmatrix{
        \Big\< p_*(m){\dispar m\over\dispar u_i},
               p_*(m){\dispar m\over\dispar u_j}
        \Big\>_{i,j} 
      & \Big\< p_*(m){\dispar m\over \dispar u_i},\epsilon_\calC
        \Big\>_i
      \cr\noalign{\vskip 0.05in}
        \Big\< \epsilon_\calC,p_*(m){\dispar m\over \dispar u_j}
        \Big\>_j
      & 1\cr }
  \hfill\cr\noalign{\vskip .1in}\hfill
   =  {|\Proj_{(T_{p(m)} p(\calC))^\perp}
           \epsilon_\calC|^{\kern .8pt 2}
       \over N\alpha (\log t)^N }
     \det\Big(\Big\< p_*(m){\dispar m\over \dispar u_i} ,
                     p_*(m){\dispar m\over \dispar u_j}
              \Big\>_{i,j}\Big) \, .
  \cr}
$$
Consequently, up to $o(\log t)^{-N}$,
$$
  \d y = 
  {|\Proj_{(T_{p(m)}p(\calC))^\perp}\epsilon_\calC|
      \over \sqrt{N\alpha}(\log t)^{N/2}} \,
  \det \big(p_*(m)^\T p_*(m)\big)^{1/2}\,
  \d\calM_\calC (m) \, \d s\, \d w \, .
$$
Let us write
$$
  p(\calC )=\{\, p(m) \, :\, m\in \calC\,\} \, ,
$$
the image of a connected component $\calC$ by $p(\cdot )$.
If follows from our calculation that that the integral in (8.2.18) is
$$\displaylines{
  \int f\Big( p(m)+{s\over \sqrt{N\alpha}}\epsilon_\calC,w\Big)
  {e^{-R(t)}\over\prod_{i\in\calI}|m_i|^\alpha}
  e^{-s} e^{-|w|^2/2}
  \, \d s \, |\Proj_{(T_{p(m)}p(\calC))^\perp}\epsilon_\calC| 
  \cut
  \det\big(p_*^\T (m)p_*(m)\big)^{1/2} 
  \, \d \calM_\calC (m) \, \d w \,
  {1\over\sqrt{N\alpha}(\log t)^{N/2}} \big( 1+o(1)\big) \, .
  \cr}
$$
Elementary algebra shows that
$$
  {e^{-R(t)}\over\sqrt{N\alpha}(\log t)^{N/2}}
  = { (K_{s,\alpha}\alpha^{\alpha/2})^N\over
      t^{N\alpha/2}\sqrt{N\alpha} }
  {1\over (2\pi )^{(d-N)/2}} \, .
$$
Consequently, as $t$ tends to infinity, (8.2.18) is equivalent
to
$$\displaylines{\qquad
  \int f\Big( p(m)+{s\over\sqrt{N\alpha}}\epsilon_\calC,w\Big)
  \prod_{i\in\calI}|m_i|^{-\alpha}
  e^{-s} {e^{-|w|^2/2}\over (2\pi )^{(d-N)/2}}\,
  \times
  \hfill\cr\vadjust{\vskip .1in}\hfill
  |\Proj_{(T_{p(m)} p(\calC ))^\perp}\epsilon_\calC| \,
  \det\big(p_*(m)^\T p_*(m)\big)^{1/2}
  \d s \, 
  \d\calM_\calC (m) \, \d w \, .
  \qquad\cr}
$$
Combining this estimate with that for $P\{\, X\in B_t\,\}$
given by Theorem 8.2.10, the conditional distribution of
$$
  \Big( \sqrt{\log t}\, \Proj_{V_\calI}
  \big( Y-Q(t)\sqrt{N}\epsilon_\calC\big) ,
  \Proj_{V_\calI^\perp}(Y) \Big)
  \qquad\hbox{ given } Y\in B_t
$$
and $Y$ in a neighborhood of $\epsilon_\calC$ converges weakly*
to that of $\big(p(M)+(S\epsilon_\calC/\sqrt{N\alpha}),W\big)$ where $M$, $S$ and
$W$ are independent with respective densities proportional to
$$\displaylines{\qquad
  |\Proj_{T_{p(m)}(p(\calC ))^\perp}\epsilon_\calC| 
  {\det\big( p_*(m)^\T p_*(m)\big)^{1/2}\over \prod_{i\in\calI}|m_i|^\alpha }\,
  \d\calM_{p(\calC )}(p) \, ,
  \hfill\cr\qquad
  e^{-s}\d s \, ,
  \qquad\hbox{ and }\quad
  \calN (0,\Id_{V_\calI^\perp}) \, .
  \hfill\cr}
$$
In particular, $S$ has an exponential density. 

Using Skorokhod's representation theorem and up to changing the
versions of the random variables, given $Y\in B_t$ and $Y$ in
the neighborhood of $Q(t)\sqrt{N}\epsilon_\calC$, we have
$$\eqalign{
  \Proj_{V_\calI}Y
  & = Q(t)\sqrt{N}\epsilon_\calC
    + {1\over\sqrt{\log t}} \Big( p(M)+{S\over\sqrt{N\alpha}}\epsilon_\calC\Big)
    + o(\log t)^{-1/2} \cr
  \Proj_{V_\calI^\perp}Y 
  & = W+o(1) \, . \cr
  }
$$
Then, given $X\in A_t$, we have
$$
  X=\Saip (\Proj_{V_\calI}Y) + \Saip \big(W+o(1)\big) \, .
$$
The term $T_\calI =\Saip\big( W+o(1)\big)$ is asymptotically a
random vector in $V_\calI^\perp$ with independent coefficients
having a Student-like $S_\alpha$ distribution. Then, for $i$
in $\calI$ and $t$ large enough, Lemma A.1.6 yields
\finetune{\hfuzz=5.5pt}
$$\displaylines{
  \Saip (\< \Proj_{V_\calI}Y,e_i\>)
  \hfill\cr\noalign{\vskip .07in}\hfill
  \eqalign{
    =\, & \Saip \Big( Q(t)\sqrt N \<\epsilon_\calC,e_i\>
          +{1\over\sqrt{\log t}} 
           \Big\< p(M)+{S\over\sqrt{N\alpha}}\epsilon_\calC,e_i\Big\>\cr
        & \kern 3in +o(\log t)^{-1/2} \Big) \cr
    =\, & \sign(m_i)\Saip\Big( Q(t)+{\sign (m_i)\over\sqrt{\log t}}
          \< p(M),e_i\>+{S\over N\sqrt\alpha}
          + o(\log t)^{-1/2} \Big) \cr
    =\, & \sqrt{t}\, \sign (m_i)\exp
          \Big( {S\over\alpha N} + {\sign (m_i)\over\sqrt\alpha}
                \< p(M),e_i\> + o(1)\Big)  \, . \cr
  }
  \cr}
$$
\finetune{\hfuzz=0pt}
If both $M$ and $m$ belong to $\calC$, then 
$$
  \alpha^{-1/2}\sign (m_i)\< p(M),e_i\> = \log |\<M,e_i\>| \, .
$$
Therefore, 
$$\eqalign{
  X & = \sqrt t \sum_{i\in\calI} e^{S/(\alpha N)}\<M,e_i\>e_i 
       \big(1+o(1)\big)
       + T_\calI + o(1) \cr
    & = \sqrt te^{S/(\alpha N)} M \big(1+o(1)\big) + T_\calI +o(1) \cr
  }
$$
To conclude the proof, notice first that $e^{S/(\alpha N)}$ has a Pareto
distribution, since 
$$
  P\Big\{\, \exp\Big( {w\over\alpha N}\Big) \geq x\,\Big\}
  = P\{\, w\geq \alpha N\log x\,\}
  = e^{-\alpha N\log x} = x^{-\alpha N} \, .
$$
Finally, as announced, let us calculate the norm of the projection
of $\epsilon_\calC$ onto the orthocomplement of the tangent space
of $p(\calC )$ at $p(m)$. Since 
$$
  p_*(m) = \sum_{i\in\calI} {e_i\otimes e_i\over |m_i|} \, ,
$$
the tangent space of $p(\calC )$ can be identified as
$$\eqalign{
  T_{p(m)}p(\calC )
  & = \sqrt\alpha \sum_{i\in\calI} {e_i\otimes e_i\over |m_i|} T_m\calC \cr
  & = \sqrt\alpha \sum_{i\in\calI} 
    {e_i\otimes e_i\over |m_i|}\{ (C+C^\T)m\}^\perp \, . \cr
  }
$$
Defining $v(m)=\sum_{i\in\calI} |m_i|\< (C+C^\T)m,e_i\>e_i$, we have
$$
  T_{p(m)}p(\calC ) = V_\calI\cap \{ v\}^\perp \, .
$$
Consequently, using that $\< Cm,m\>=1$, we have
$$\eqalign{
  |\Proj_{T_{p(m)}p(\calC )^\perp}\epsilon_\calC |
  & = { \< v(m),\epsilon_\calC\> \over |v(m)| } \cr  
  & = { 1\over\sqrt N} {\sum_{i\in\calI} mi\< (C+C^\T)m,e_i\> \over
        \big( \sum_{i\in\calI} m_i^2 \< (C+C^\T)m,e_i\>^2\big)^{1/2} } \cr
  & = \Big( N\sum_{i\in\calI}m_i^2\< (C+C^\T)m,e_i\>^2\Big)^{1/2} \, .\cr
  }
$$
This concludes the proof of Theorem 8.2.20.\hfill$\qed$
\bigskip

Theorems 8.2.1 and 8.2.10 do not settle completely the tail
behavior of $\cXX$. For instance, it may happen that
the largest diagonal term of $C$ vanishes, or that $N(C)>1$ 
but (8.2.8) does not hold. The situation not covered 
assumes that
\setbox1=\vtop{\hsize=3.5in\noindent%
  $\<Cu,u\>=0$ for some nonzero $u$ in $V_\calI$ 
  and some $\calI$ in $T$, of cardinality strictly less than $N(C)$.}
$$
  \raise 7pt\box1\eqno{(8.2.19)}
$$

The technique developed in this section may work in this situation;
but we will see in the next section that, under (8.2.19), some extra
complication is added. For the time being, we will only prove a rather
weak result. It will be useful for the statistical applications
developed in chapter 9.

Denote\ntn{$N_0(C)$} $N_0(C)$ the smallest cardinal of a set $\calI$
such that (8.2.19) holds. Under (8.2.19), the next result asserts that
the tail probability of $\cXX$ is much lighter 
than $t^{-\alpha N_0(C)/2}$.

\bigskip

\state{8.2.21. THEOREM.} {\it%
  Let $X$ be a $d$-dimensional random vector with independent
  and identically distributed components having a Student
  like distribution with parameter $\alpha$. Let $C$ be a $d\times d$ 
  matrix, and write $N_0=N_0(C)$. Under (8.2.19), and if 
  $\alpha>2/N_0$,
  $$
    \lim_{\ttoi} t^{\alpha N_0/2} P\{\, \cXX\geq t\,\} = 0 
    \, .
  $$
  }

%\bigskip
\finetune{\kern -.1in}

\stateit{Proof.} The statement is obvious if $C$ is negative.
Thus, we assume that $C$ is not negative. For $\epsilon$
positive, denote $C_\epsilon = C+\epsilon \Id$. The proof relies on
two very basic observations. The first one is that for any
positive $\epsilon$ the matrix $C_\epsilon-C$ is positive.
Consequently, for any $t$,
$$
  P\{\, \cXX\geq t\,\} \leq 
  P\{\,\< C_\epsilon X,X\> \geq t\,\} \, .
$$

Assumption (8.2.19) states that $N_0< N(C)$. Let\ntn{$J_0(C)$}
$$
  J_0(C)=\{\, \calI\in T \, :\, |\calI |=N_0(C) \, ,\, 
  \exists u\in V_\calI\setminus \{\, 0\,\} \, ,\, 
  \< Cu,u\>=0\,\} \, .
$$

The second observation is stated in the following lemma.

\bigskip

\state{8.2.22. LEMMA.} {\it%
  If $\epsilon$ is positive and small enough, 
  then $N(C_\epsilon)=N_0(C)$ and $J(C_\epsilon )=J_0(C)$. 
  Moreover, $C_\epsilon$ satisfies assumption (8.2.8).
  }

\bigskip

\stateit{Proof.} Denote by $T_{N_0}$ the set of all subsets of
$\{\, 1,\ldots ,d\,\}$ of cardinal
at most $N_0$. Assumption (8.2.19) is equivalent to the
following. Whenever $\calI$ is in $T_{N_0}\setminus J_0(C)$, the
compression of $C$ to $V_\calI$ is negative; moreover,
for any $\calI$ in $J_0(C)$, this compression is nonpositive 
and there exists a nonzero $z_\calI$ in $V_\calI$ 
for which
$\<Cz_\calI,z_\calI\>=0$. Consequently, (8.2.19) implies the
existence of a positive $\epsilon_0$ such that 
$\<Cu,u\> <-\epsilon_0$ for any $\calI$ in
$T_{N_0}\setminus J_0(C)$ and any unit vector $u$ of $V_\calI$.
If $\epsilon$ is positive, less than $\epsilon_0$, 
and $\calI$ is in $T_{N_0}\setminus J_0(C)$, 
the compression of $C_\epsilon$ to $V_\calI$ is negative. 
Furthermore, if $\calI$ is in $J_0(C)$, then $\<C_\epsilon
z_\calI,z_\calI\> = \epsilon |z_\calI|^2$ is positive. 
This proves
$N(C_\epsilon )=N_0$ and $J(C_\epsilon)=J_0(C)$ as well as
$C_\epsilon$ satisfies (8.2.8). \hfill$\qed$

\bigskip

Let $\epsilon$ be positive and small enough so that the conclusions of
Lemma 8.2.22 hold. For $\calI$ in $J_0(C)=J(C_\epsilon)$,
denote\ntn{$M_{\epsilon,\calI}$}
$$
  M_{\epsilon,\calI} =\{\, m\in  V_\calI : \<C_\epsilon m,m\>
  = 1 \,\} \, .
$$
Using our two observations and Theorem 8.2.10, we infer that
for any $\epsilon$ positive and small enough,
$$\displaylines{\qquad
    \limsup_{\ttoi} t^{\alpha N_0/2}
    P\{\, \cXX\geq t\,\} 
  \cut
    \leq
    {K_{s,\alpha}\alpha^{\alpha N_0/2}\over \sqrt{\alpha N_0}}
    \sum_{\calI\in J_0(C)}\int_{M_{\epsilon,\calI}}
    { \det\,\tilde G(m)\over 
      \prod_{i\in\calI} |m_i|^\alpha }
    \d\calM_{M_{\epsilon ,\calI}}(m) \, .
  \qquad}
$$
To conclude the proof, we show that the above upper bound
tends to $0$ as $\epsilon$ tends to $0$. For this purpose, we need a
good description of $M_{\epsilon ,\calI}$. It will be helpful
to have some understanding of the sets\ntn{$Z_\calI (C)$}
$$
  Z_\calI (C)=\{\, u\in V_\calI \, :\, \<Cu,u\>=0\,\}\, ,
  \qquad \calI\in J_0(C) \, .
$$
Despite the fact that these sets are defined by a quadratic
equation, minimality of $N_0$ implies that they are linear
spaces.

\bigskip

\state{8.2.23. LEMMA.} {\it%
  Under (8.2.19), the set $Z_\calI(C)$ is a one dimensional
  vector space. It coincides with the null eigensubspace in
  $V_\calI$ of the compression of $C+C^\T$ to $V_\calI$. 
  }

\bigskip

\stateit{Proof.} Since the compression of $C$ to $V_\calI$ is
nonpositive,
$$
  Z_\calI (C)=\big\{\, u\in V_\calI \, : \, \< Cu,u\>
  = \sup_{v\in V_\calI}\<Cv,v\>\,\big\} \, .
$$
The equality $2\< Cu,u\>=\< (C+C^\T)u,u\>$, shows that $Z_\calI(C)$ is
the eigensubspace associated to the largest eigenvalue of the
compression of $(C+C^\T)$ to $V_\calI$. Thus, $Z_\calI(C)$ is indeed a
linear space.

Let $i$ be in $\calI\in J_0(C)$ and $u$ be a nonzero vector in
$Z_\calI(C)$. If $\<u,e_i\>$ vanishes, $u$ belongs to the
$(N_0-1)$-dimensional subspace $V_{\calI \setminus\{ i\}}$,
contradicting the minimality of $N(C)$.  Consequently, no component of
$u$ in $V_\calI$ vanishes.

Assume now that we can find two linearly independent vectors
$u$ and $v$ in $Z_\calI(C)$. Since none of their components
in $V_\calI$ vanish, there exists a linear combination of
these two vector with at least one component in $V_\calI$
vanishing. This linear combination is a vector which 
contradicts what we just showed. Consequently, there
are not two linearly independent vectors in $Z_\calI(C)$.
Since $Z_\calI(C)$ in nonempty by definition of $N_0$, this
concludes the proof. \hfill$\qed$

\bigskip

It follows from Lemma 8.2.23 that the unit sphere of the space
$Z_\calI(C)$ is actually made of two points, $z$ and $-z$.

To describe $M_{\epsilon ,\calI}$, we now introduce the sphere
corresponding to the compression of $-C$ to the
orthocomplement of $Z_\calI(C)$ in $V_\calI$, that is\ntn{$S_\calI(C)$}
$$
  S_\calI(C)=\{\, v\in V_\calI \, :\,
  \< Cv,v\>=-1 \, ;\, v\perp Z_\calI (C)\, \} \, .
$$
Since $Z_\calI(C)$ is the eigenspace associated to the
simple null
eigenvalue of $\Proj_{V_\calI}(C+C^\T)\big|_{V_{\calI}}$, the
set $S_\calI(C)$ is a compact ellipsoid. We can now make
$M_{\epsilon,\calI}$ explicit.

\bigskip

\state{8.2.24. LEMMA.} {\it%
  For any $\calI$ in $J_0(C)$, for any $\epsilon$ positive 
  and small enough,
  $$\displaylines{\qquad
      M_{\epsilon,\calI}=\Big\{\, {\eta\over \sqrt\epsilon}
      \sqrt{1+s(1-\epsilon)}\, z+ \sqrt{s}\, v\, :\, 
      v\in S_\calI(C)\, ,\,
    \cut
      s\geq 0\, ,\,
      \eta\in\{\, -1,1\,\}\,\Big\} \, ,
    \qquad\cr}
  $$
  where $z$ is in the unit sphere of $Z_\calI(C)$.
  }

\bigskip

\stateit{Proof.} Let $\calI$ be in $J_0(C)$. A point $m$ in $V_\calI$
belongs to $M_{\epsilon ,\calI}$ if and only if 
$$
  1=\< Cm,m\> +\epsilon |m|^2 \, .
  \eqno{(8.2.20)}
$$
The nonpositivity of $C$ on $V_\calI$ implies $\< Cm,m\>\leq
0$. Therefore $|m|^2\geq 1/\epsilon$. In particular, $m$
is nonzero. Write $\lambda=|m|^2$ and $p=m/|m|$. 
Equality (8.2.20) and
nonpositivity of $C$ on $V_\calI$ imply
$$
  0\geq \< Cp,p\> = -\epsilon +{1\over\lambda}\geq -\epsilon 
  \, .
  \eqno{(8.2.21)}
$$
The function $u\mapsto \<Cu,u\>$ is continuous. Its restriction
to the --- compact --- unit sphere of $V_\calI$ is maximum exactly
on $Z_\calI(C)$. Consequently, if $U$ is an arbitrary
neighborhood of the unit sphere of $Z_\calI(C)$ on 
the unit sphere of $V_\calI$, the inclusion
$$
  \Big\{ {m\over |m|}\, :\, m\in M_{\epsilon ,\calI}\Big\}
  \subset U
$$
holds for any positive $\epsilon$ small enough.  Then, we can write
$p=\cos\theta z + \sin \theta q$ for some $q$ in $V_\calI$ of norm $1$
and orthogonal to $z$, and some $\theta$ close to $0$ mod $\pi$ ---
the coefficients $\cos\theta$ and $\sin\theta$ are imposed by $|p|=1$;
speaking geometrically, we parameterize the sphere in normal
coordinates. Since $\Proj_{V_\calI}(C+C^\T)z=0$ thanks to Lemma
8.2.23, equation (8.2.21) forces
$$
  \sin^2\theta\, \< Cq,q\> = -\epsilon +{1\over\lambda} \, .
  \eqno{(8.2.22)}
$$
If $\theta=0$ mod $\pi$, this forces $\lambda=1/\epsilon$ and
$p=z$. Thus $m=z/\sqrt\epsilon$. Assume that $\theta\neq 0$
mod $\pi$. Since $S_\calI(C)$ is invariant under $\{\,
-\Id,\Id\,\}$, equality (8.2.22) imposes
$$
  q= {1\over\sin \theta}\sqrt{\epsilon-{1\over \lambda}}\, v
  \qquad\hbox{ for some } v\in S_\calI(C) \, .
$$
We can then write
$$
  m=\sqrt{\lambda} p
  = \sqrt{\lambda}\cos\theta \, z 
  + \sqrt{\epsilon\lambda -1} \, v
  \, .
$$
But equation (8.2.20) yields
$$
  1= -\epsilon\lambda + 1 + \epsilon (\lambda\cos^2\theta
  + \epsilon\lambda-1)\, ,
$$
from which we can obtain $\cos^2\theta$. Thus,
$$
  m=\eta\sqrt{1+\lambda (1-\epsilon)} \, z 
  +\sqrt{\epsilon\lambda -1} \, v
  \qquad \hbox{ for } \eta\in \{\, -1,1\,\} \, .
$$
Set $s=\epsilon\lambda-1$ to obtain the parameterization given
in the lemma. Notice that for $s=-1$, we obtain
$m=\eta\sqrt{\lambda} z$. \hfill$\qed$

\bigskip

We can now obtain some good bounds on the components of
points $m$ belonging to $M_{\epsilon ,\calI}$.
We write $m_i=\< m,e_i\>$ the $i$-th component of $m$, as 
we did in the proof of Theorem
8.2.10. The following statement is the analogue of Lemma
8.2.16.

\bigskip

\state{8.2.25. LEMMA.} {\it%
  There exists a positive $M$ such that for any
  positive $\epsilon$ small enough, any point
  $m={\displaystyle\eta\over\displaystyle\sqrt\epsilon}
  \sqrt{1+s(1-\epsilon)}\, z +
  \sqrt{s}\, v$ in $M_{\epsilon ,\calI}$ satisfies for all
  $i$ in $\calI$ the inequality
  $$
    {1\over M}\leq \sqrt{\epsilon\over 1+s} |m_i|
    \leq M \, .
  $$}

\bigskip

\stateit{Proof.} For $\epsilon$ small enough and any
positive $s$, we
have $s/(2\epsilon)\leq s(1-\epsilon)/\epsilon\leq s/\epsilon$.
If $\eta\< z,e_i\>$ is positive, this implies
$$
  \sqrt{\epsilon\over 1+s}m_i
  \leq \eta\< z,e_i\> +\sqrt{s\epsilon\over1+s}\< v,e_i\>
$$
and
$$
  \sqrt{\epsilon\over 1+s}m_i\geq \sqrt{1+(s/2)\over 1+s}
  \eta\< z,e_i\> +\sqrt{s\epsilon\over 1+s}\< v,e_i\> \, .
$$
Since $S_\calI(C)$ is compact and the functions 
$s\mapsto s/(1+s)$ and $s\mapsto (2+s)/(1+s)$ are bounded
on $[0,\infty )$, the result follows.

If $\eta\< z,e_i\>$ is negative, we proceed similarly.
\hfill$\qed$

\bigskip

We can conclude the proof of Theorem 8.2.21 with
the next result.

\bigskip

\state{8.2.26. LEMMA.} {\it%
  For any $\calI$ in $J_0(C)$,
  $$
    \lim_{\epsilon\to 0}\int_{M_{\epsilon,\calI}}
    {\det\, \tilde G(m)\over \prod_{i\in\calI}
    |m_i|^\alpha} \d\calM_{M_{\epsilon,\calI}}(m) = 0 \, .
  $$
  }

\bigskip

\stateit{Proof.} Arguing as in Lemma 8.2.17, we obtain
$$
  \det\, \tilde G(m)\leq M/|m|^{N_0-1}
$$
for any $m$ in $M_{\epsilon,\calI}$ and some fixed
positive $M$, not depending on $\epsilon$ or $m$. 
Using Lemma 8.2.25, we have
$$\displaylines{\qquad
  \int_{M_{\epsilon,\calI}} {\det\,\tilde G(m)
  \over \prod_{i\in\calI}|m_i|^\alpha}
  \d\calM_{M_{\epsilon,\calI}}(m)
  \cut
  \leq \int_{M_{\epsilon,\calI}}
  \bigg({\epsilon\over 1+s}\bigg)^{N_0-1+N_0\alpha\over2}
  \d\calM_{M_{\epsilon,\calI}}(m)  O(1)\, ,
  \qquad (8.2.22)\cr
  }
$$
as $\epsilon$ tends to $0$. To express the Riemannian 
measure on $M_{\epsilon,\calI}$,
notice that $S_\calI(C)$ is a manifold of dimension $N_0-2$.
Write $v(x_1,\ldots ,x_{N_0-2})$ for a local parameterization of
$S_\calI(C)$. It induces a parameterization
$$
  m(s,x_1,\ldots ,x_{N_0-2})
  = {\eta\over \sqrt\epsilon}\sqrt{1+s(1-\epsilon)}\, z
  + \sqrt{s}\, v(x_1,\ldots, x_{N_0-2})
$$
of $M_{\epsilon,\calI}$. In this local chart
$$\displaylines{
    \d\calM_{M_{\epsilon,\calI}}(m)
    = \bigg| { \sqrt{1-\epsilon} \, s^{(N_0-2)/2} \over
               \sqrt{\epsilon}\sqrt{1+s(1-\epsilon)} }
    \det \Big( z,{\partial v\over \partial x_1},\ldots ,
    {\partial v\over\partial x_{N_0-2}}\Big) 
  \cut
    + { s^{(N_0-2)/2} \over 2\sqrt s }
    \det \Big( v,{\partial v\over \partial x_1},\ldots ,
    {\partial v\over\partial x_{N_0-2}}\Big)\bigg| \, 
    \d s\wedge \d x_1\wedge \ldots \wedge \d x_{N_0-2} \, .
  \cr
  }
$$
Since $z$ is orthogonal to $S_\calI(C)$ and of unit norm, 
we have
$$
  \det \Big( z,{\partial v\over \partial x_1},\ldots ,
  {\partial v\over\partial x_{N_0-2}}\Big) 
  \d x_1\wedge \ldots \wedge \d x_{N_0-2}
  = \d\calM_{S_\calI(C)} \, .
$$
Therefore,
$$\displaylines{
    \d\calM_{M_{\epsilon,\calI}}(m)
    \leq \sqrt{ 1-\epsilon\over 1+s(1-\epsilon)}
    {1\over \sqrt\epsilon} \, s^{(N_0-2)/2} \, \d s\,
    \d\calM_{S_\calI(C)}(v)
  \cut
    + {s^{(N_0-3)/2}\over 2} 
      |\Proj_{(T_vS_\calI(C))^\perp}v| \,
    \d s\, \d\calM_{S_\calI(C)}(v) \, .
  \qquad\cr}
$$
Consequently, (8.2.21) is less than
$$\displaylines{
    \int_{s>0}\int_{S_\calI(C)} 
    { \epsilon^{(N_0-2+N_0\alpha)/2}\over 
      (1+s)^{(N_0-1+N_0\alpha )/2} }
    {  s^{(N_0-2)/ 2}\over \sqrt{1+s(1-\epsilon)} }
  \cut
    + { \epsilon^{(N_0-1+N_0)/ 2}\over
        (1+s)^{(N_0-1+N_0\alpha )/ 2 } }
    { s^{(N_0-3)/ 2}\over 2 } |v|
    \, \d s\, \d\calM_{S_\calI(C)}(v) 
   \, .
  \qquad\cr}
$$
We can bound the integral in $v$ by the constant 
$\diam S_\calI(C)\vol S_\calI (C)$, for instance. 

Then, the integral in $s$ is less than
$$\displaylines{\qquad
    2+2\int_1^\infty 
    s^{ {N_0-2\over 2}-{N_0-1+N_0\alpha\over 2}
        -{1\over 2}} + s^{{N_0-3\over 2} 
        - {N_0-1+N_0\alpha\over 2} }
    \, \d s
  \cut
    = 2+ 4\int_1^\infty s^{-1-{N_0\alpha\over 2}} \d s 
    < \infty \, .
  \qquad\cr}
$$
Therefore, the right hand side of (8.2.21) is less than a
constant times
$$
  \epsilon^{{N_0(1+\alpha)\over 2}-1}
  + \epsilon^{N_0-1+N_0\alpha\over 2} \, .
$$
This tends to $0$ as $\epsilon$ tends to $0$, provided 
$\alpha$ is strictly larger than $(2/N_0)-1$. This concludes the proof 
of Lemma 8.2.26 as well
as that of Theorem 8.2.21.\hfill$\qed$

\bigskip

\noindent{\fsection 8.3. Heavy tail and degeneracy.}
%\chapternumber={8.3}
\section{8.3}{Heavy tail and degeneracy}

\bigskip

\noindent The results of section 8.2 are incomplete. The tail
behavior of $\cXX$ is not described when the largest diagonal
coefficient of $C$ vanishes, or when (8.2.8) does not hold. 
The proof of Theorem 8.2.10 breaks down in these cases.

In example 4 of chapter 1, we dealt with the tail of the product
of two independent Cauchy random variables. In the language of the
current chapter, we considered the $2\times 2$-matrix
$$
  C={1\over 2} \pmatrix{ 0 & 1 \cr 1 & 0 \cr} \, .
$$
We saw that the tail of $\cXX$ is in $t^{-1}\log t$. For this
specific matrix $N(C)=2$. Comparing with the result obtained in the
previous section, we have an extra logarithmic factor in $t$.
Going back to chapter 1, the factor $t^{-1}$ is explained by the
fact that for most points $(x,y)$ on the boundary of
$$
  A_t=\{\, (x,y) \, :\, xy\geq t\,\} = \sqrt{t} A_1\, ,
$$
the minimum of $x$ and $y$ tends to infinity, sharing some similarity
with the case where $N(C)=2$ and (8.2.8) holds.  The factor $\log t$
is explained by the closeness of the boundary $\partial A_t$ to the
level sets of the Cauchy density, along a sizeable part of the
hyperbola $xy=t$.

Quite amazingly, Theorem 5.1 can still be used in this
situation, but the analysis is a bit more involved than needed
to prove Theorems 8.2.1 or 8.2.10. Our goal 
in this section is rather modest.
We do not intend to obtain the tail behavior of the quadratic
form whenever (8.2.8) does not hold. We will concentrate only on
matrices $C$ with vanishing largest diagonal coefficient. The
reason is twofold. First, this simple degeneracy is sufficient
to understand what we should do when (8.2.8) does not hold.
Second, we will make use of this specific case in chapter 11,
when studying a statistical application.

\bigskip

\state{8.3.1. THEOREM.} {\it%
  Let $C$ be a $d\times d$-real matrix, with vanishing largest
  diagonal coefficient. Let $X$ be a random vector
  in $\RR^d$, with independent components, having a Student-like
  distribution $S_\alpha$. The following asymptotic equivalence
  holds as $t$ tends to infinity,
  $$
    P\{\, \cXX\geq t\,\} \sim
    {\log t\over t^\alpha} K_{s,\alpha}^2\alpha^\alpha
    \sum_{i:C_{i,i}=0}\sum_{1\leq j\leq d}
    |C_{i,j}+C_{j,i}|^\alpha \, .
  $$}

\bigskip

\stateit{Proof.} Since $\cXX= \< (C+C^\T )X,X\>/2$, we will assume
throughout the proof that $C$ is symmetric. The proof builds
upon that of Theorem 8.2.10 in the case $N(C)=2$. However, Lemma
8.2.11 cannot be used anymore since (8.2.8) does not hold. 

As before, we denote
$$
  A_t=\{\, p\in \RR^d \,:\, \<Cp,p\>\geq t\,\} = \sqrt{t}A_1
$$
and
$$
  B_t= \pisa (A_t) = \{\, \pisa (p) \, :\, p\in A_t\,\} \, .
$$
This guarantees
$$
  P\{\, \cXX \geq t\,\} 
  = {1\over (2\pi )^{d/2}}\int_{B_t} e^{-|y|^2/2} \d y \, .
$$
The main part of the proof is to precisely locate and describe
the points in $B_t$ with nearly minimal norm. This will be done
through the next five lemmas.

Since all the diagonal coefficients of $C$ are nonpositive, $N(C)\geq
2$. Our first lemma shows that $N(C)=2$, an easy fact. It also
introduces points which will be essential to determine a
dominating manifold. Define\ntn{$J_*(C)$}
$$
  J_*(C)=\big\{\, \{\,i,j\,\} \, :\, 1\leq i,j\leq d \, ,\,
  C_{i,i}C_{j,j}=0 \hbox{ and } C_{i,j}\neq 0 \, \big\} \, .
$$
If $\{\, i,j\,\}\in J_*(C)$ then either $C_{i,i}$ or
$C_{j,j}$ is null. 

\bigskip

\state{8.3.2. LEMMA.} {\it%
  Let $\calI = \{\, i,j\,\}$ be in $J_*(C)$, and assume that
  $C_{i,i}=0$. The solutions of the equation $\<Cm,m\>=1$ 
  in $V_\calI$ are
  $$
    {1\over 2C_{i,j}}\Big({1\over m_j}-m_jC_{j,j}\Big) e_i
    + m_j e_j \, , \qquad m_j\in \RR\setminus \{\, 0\,\} \, .
  $$}

\finetune{}%\bigskip

\stateit{Proof.} Let $m=m_ie_i+m_je_j$ be a point in $V_\calI$. 
The equation $1=\< Cm,m\>$ can be rewritten as
$$
  1=2C_{i,j}m_im_j + C_{j,j}m_j^2 \, .
$$
Since $(m_i,0)$ is not solution, we obtain the result
in expressing $m_i$ as a function of $m_j$. \hfill$\qed$

\bigskip

It follows from Lemma 8.3.2 that $J(C)$ contains $J_*(C)$. As we
did in the proof of Theorem 8.2.10, define\ntn{$M_\calI$}
$$
  M_\calI = \{\, m\in V_\calI \, : \, \< Cm,m\>=1 \,\} \, , 
  \qquad \calI\in J(C) \, ,
$$
and\ntn{$\gamma_\calI$}
$$
  \gamma_\calI = \inf\Big\{\, \prod_{i\in \calI} |m_i| \, :\, 
  m\in M_\calI \,\Big\} \, .
$$
It will be convenient later to extend slightly the notation
introduced so far. For $\calI=\{\, i,j\,\}$ in $J(C)$, we will write
$M_{i,j}$ or $M_{j,i}$ for $M_\calI$.

If $\calI$ belongs to $J(C)\setminus J_*(C)$, Lemma 8.2.11 
shows that $\gamma_\calI$ is po\-si\-tive. If $\calI=\{\,i,j\,\}$
is in $J_*(C)$, Lemma 8.3.2 gives
$$
  \gamma_\calI = \inf\bigg\{\, 
  {|1-m_j^2C_{j,j}|\over 2|C_{i,j}|}
  \, :\, m_j\in \RR\setminus\{\, 0\,\} \,\bigg\}
  = {1\over 2|C_{i,j}|} \, ,
$$
because $C_{j,j}$ is nonpositive. Consequently,\ntn{$\gamma$}
$$
  \gamma=\min\{\, \gamma_\calI \, : \, \calI\in J(C)\,\}
$$
is positive. As in Theorem 8.2.10, we hope\ntn{$I(\cdot )$}
$$
  I(y)={|y|^2\over 2}-\log (2\pi)^{d/2}
$$
will be minimum on $B_t$ at points lying in $\pisa
(\sqrt{t}\bigcup_{\calI\in J(C)} M_\calI )$. Let\ntn{$R(t)$}
$$\displaylines{\qquad
    R(t)=2\alpha\log\sqrt{t} -\log\log \sqrt{t} 
    - 2\log (K_{s,\alpha}\alpha^{\alpha/2}2\sqrt{\pi} )
  \hfill\cr\noalign{\vskip 1pt}\hfill
    + \log (2\pi )^{d/2} + \alpha\log\gamma \, .
  \qquad\cr}
$$
If $\calI$ and $m$ are fixed, respectively in $J(C)$
and $M_\calI$, the proof of Theorem 8.2.10 --- see the 
calculation of $I\big(q_{\calI,t}(m,v)\big)$ before Lemma 
8.2.13 --- shows that
$$
  I\big(\pisa (\sqrt t m)\big) = R(t) - \alpha\log\gamma
  + \alpha\log \Big(\prod_{i\in\calI} |m_i|\Big) + o(1)
$$
as $t$ tends to infinity.
This suggests that $I(B_t)= R(t)+o(1)$ as $t$ tends to infinity. 
To locate the
minima of $I$ on $B_t$, we try to imitate the second assertion
of Lemma 8.2.13. The coming lemma gives a coarser estimate.
After its proof, we will be able to explain more precisely the
difficulty created by the nonemptyness of $J_*(C)$.

\bigskip

\state{8.3.3. LEMMA.} {\it%
  Let $M$ be a number strictly larger than $2$, and $\epsilon$
  be positive. The set of all points $p$ in $\partial A_1$ 
  such that
  $$
    I\big( \pisa (\sqrt{t}p)\big) \leq R(t) + M\log\log t
  $$
  is included in an 
  $\epsilon t^{-1/6}(\log t)^{M/\alpha}$-neighborhood of
  $\bigcup_{\calI\in J(C)} M_\calI$.
  }

\bigskip

\stateit{Proof.} For any point $p$, we use Lemma A.1.5 to obtain
$$\displaylines{
  I\big(\pisa (\sqrt{t}p)\big)\vtop to .15in{}
  \cut
  \eqalign{
    & \geq {1\over 2} \sum_{1\leq i\leq d} \pisa\big( \epsilon
      t^{{1\over 2}-{1\over 6}}(\log t)^{M/\alpha}\big)
      \indic_{[\epsilon t^{-1/6}(\log t)^{M/\alpha},\infty )} 
      (|p_i|) \cr
    & = {1\over 2} \bigg( 2\alpha\log\big( 
      t^{1/3}(\log t)^{M\over\alpha}\big)^2 
      -\log\log \big( t^{1/3}
      (\log t)^{M/\alpha}\big) \big( 1+ o(1)\big) \bigg)
      \times\cr
    & \kern 1.7in \sharp\big\{\, 1\leq i\leq d \, :\, 
      |p_i|\geq \epsilon t^{-1/6}(\log t)^{M/\alpha}\,\big\} \, .%
    \cr}
  \cr}
$$
Consequently, for all the points under consideration in 
the lemma, and for $t$ large enough,
$$\displaylines{\qquad
    \sharp\big\{\, 1\leq i\leq d \, :\, 
    |p_i|\geq \epsilon t^{-1/6}(\log t)^{M/\alpha}\,\big\}
  \hfill\cr\noalign{\vskip .05in}\hfill
  \eqalign{
    \leq\;& {\alpha\log t + (M-1)\log\log t + O(1)
            \over
            {\alpha\over 3} \log t + 
             \big(M-{1\over 2}\big)\log\log t
             \big( 1+o(1)\big) } \cr
    \leq\;& 3 - {3\over\alpha}(2M-4)
             {\log\log t\over \log t}\big( 1+o(1)\big) \cr
       <\;& 3 \, . \cr}
  \qquad\cr
  }
$$
Thus, at most two of the $p_i$'s have their absolute value larger
than $\epsilon t^{-1/6}(\log t)^{M\alpha}$. Since $N(C)=2$, this
concludes the proof. \hfill$\qed$

\bigskip

Now, let us see precisely why the proof of Theorem 8.2.10 breaks
down. Lemma 8.2.11 still holds for $\calI$ in $J(C)\setminus
J_*(C)$. But it fails if $\calI=\{\, i,j\,\}$ is in $J_*(C)$, since
Lemma 8.3.2 shows that we can have $|m_i|$ as small as we like;
and whenever $C_{j,j}$ vanishes, $|m_j|$ can also be as close 
to $0$ as desired. Lemma 8.2.11 was used in deriving the expression
of $q_{\calI,t}(m,v)$ in Lemma 8.2.12. Then all the proof was
more or less based on this approximation of $q_{\calI,t}(m,v)$.
Notice that we can still use this expression whenever we can
make the expansion which was used in its proof. Thus, for
$\calI=\{\, i,j\,\}$ in $J(C)$, we still parameterize 
$\partial A_1$ near a point $m$ of $M_\calI$, as\ntn{$p_\calI(m,v)$}
$$
  p_\calI(m,v)= m\big( 1+h(v)\big) + v \, , \qquad
  v\in T_m\partial A_1\ominus T_mM_\calI \, .
$$
We need to be able to prove that only those $v$'s such that
$v\ll m$ componentwise are of interest for us. Thanks to Lemma
8.3.3, this can be done right away on the range where
$$
  |m_i|\wedge |m_j|\geq t^{-1/6}(\log t)^{M/\alpha}
  \eqno{(8.3.1)}
$$
say. But it can be seen easily that we need a larger range when
$\calI$ belongs to $J_*(C)$.

Another difference with the situation of Theorem 8.2.10 is that
the sets $M_\calI$ are no longer far apart. It is true that the
distance between $\bigcup_{\calI\in J_*(C)}M_\calI$ and 
$\bigcup_{\calI\in J(C)\setminus J_*(C)}M_\calI$ is strictly
positive. Also the distance between $M_\calI$ and $M_\calJ$ for
$\calI$ and $\calJ$ distinct in $J(C)\setminus J_*(C)$ is
positive. This follows from Lemma 8.2.11. However, if $\{\,i,j\,\}$
and $\{\,i,k\,\}$ are in $J_*(C)$ with $j\neq k$ then $M_{i,j}$ and
$M_{i,k}$ are at zero distance. They connect at infinity on the axis
$e_i\RR\,$. This can be seen as follows. For $m_k=m_jC_{i,j}/C_{i,k}$,
the corresponding point on $M_{i,k}$ given by Lemma 8.3.2 is
$$\displaylines{\qquad
    {1\over 2 C_{i,k}}\Big( {1\over m_k}-m_kC_{k,k}\Big) e_i +
    m_k e_k
  \cut
    = \Big( {1\over 2C_{i,j}m_j} - 
    {C_{i,j}C_{k,k}\over 2C_{i,k}^2} m_j\Big) e_i
    + {m_jC_{i,j}\over C_{i,k}}e_k \, . 
  \qquad\cr
  }
$$
Its distance to
$$
  {1\over 2 C_{i,j}}\Big( {1\over m_j}-m_j C_{j,j}\Big) e_i
  + m_j e_j \in M_{i,j}
$$
is
$$
  \bigg( \Big( {C_{i,j}C_{k,k}\over 2C_{i,k}^2}\Big)^2
  + \Big( {C_{i,j}\over C_{i,k}}-1\Big)^2\bigg)^{1/2}
  |m_j| \, .
$$
It tends to $0$ as $m_j$ tends to $0$. This has a dramatic
consequence. For $m_j$ small, we can have $p_{\{ i,j\}}(m,v)
= p_{\{ i,k\}}(m',v')$ for some $m$ in $M_{i,j}$, some $m'$ in
$M_{i,k}$ and $v$, $v'$ satisfying the a priori 
estimate given in Lemma 8.3.2.
In other words some components of $v$, even small, may cancel
with the corresponding components of $m$. Thus, the naive
approach used in Lemma 8.2.12 cannot succeed.

It is essential to understand that what goes wrong here is the
pa\-ram\-e\-teri\-za\-tion of the set $\partial A_1$. 
The naive parameterization with $(m,v)$ is onto, but is not
one-to-one anymore, at least in the interesting range.
The trick is then to introduce a new pa\-ram\-e\-teri\-za\-tion, 
well tailored to handle small
components of the points $m$ in $\bigcup_{\calI\in J_*(C)}M_\calI$.
Lemma 8.3.2 shows that when (8.3.1) fails, then one component of
$m$ has to be of order larger than $t^{1/6}(\log t)^{M/\alpha}$.
The corresponding component of $v$ in $p_\calI (m,v)$ can be
neglected thanks to Lemma 8.3.3. Writing $p$ for the corresponding
point in $\partial A_1$, we have $p_i\sim m_i$. Since $p_i$ is
very large and Lemma 8.3.3 tells us to look in a neighborhood of
$m$, all the other components of $p$ must be small, going to $0$
as $t$ tends to infinity. We can then single out another
component of $p$, say $p_k$, such that we can locate $p$ near
$M_{i,k}$. To do so, write
$$
  p= p_ie_i + p_k e_k + w
$$
where $w$ is orthogonal to $V_{\{i,k\}}$. Since $p$ is in
$\partial A_1$ and $C_{i,i}$ vanishes,
$$
  1= 2 C_{i,k} p_ip_k + p_k^2 C_{k,k} + 2 p_i\< Ce_i,w\>
  + 2p_k\< Ce_j,w\> + \< Cw,w\> \, .
$$
Using Lemma 8.3.2, we focus on a zone where $p_k$ and $w$ are
$o(1)$ as $t$ tends to infinity, and obtain
$$
  1= 2C_{i,k}p_ip_k + 2p_i\< Ce_i,w\> +o(1) \, .
$$
Notice that
$$
  Ce_i \in \span \big\{\, e_l \, :\, \{\,i,l\,\}\in J_*(C)
  \,\big\} \, .
$$
Thus, whenever $\{\,i,l\,\}$ is in $J_*(C)$, if we can
guarantee $|w_l|\leq |p_k|$, we obtain
$$
  1\leq ( 2|C_{i,k}|+\| C\|\sqrt{d}) |p_ip_k| + o(1) \, .
$$
This gives us a lower bound on $|p_k|$. It will be good enough
to bound $\sqrt{t} p_k$ away from $0$, and
to do an asymptotic expansion of $\pisa\big(
\sqrt{t}p_k)\big)^2$ in the calculation of $I\big(\pisa
(\sqrt{t}p)\big)$. We will then be able to improve the a priori
bound of Lemma 8.3.3. Only then we will use the parameterization
$p_\calI(m,v)$.

We have not said how to guarantee that $|w_l|\leq |p_k|$. This is
easy. Take $p_k$ to be the largest of the $|p_j|$'s for
$\{\,i,j\,\}$ in $J_*(C)$. This index $k$ depends on $p\in\partial
A_1$, and this parameterization is one-to-one up to a set of
Lebesgue measure zero.

To proceed along these lines, define the function\ntn{$Q(\cdot )$}
$$
  Q(x) = \sqrt{2\alpha\log\sqrt x} 
  - { \log\log \sqrt x \over 2\sqrt{2\alpha\log\sqrt x} }
  - { \log (K_{s,\alpha}\alpha^{\alpha /2}2\sqrt{\pi} )\over
      \sqrt{2\alpha\log\sqrt x} } \, , \qquad x> 0 \, .
$$

\bigskip

\state{8.3.4. LEMMA.} {\it%
  Let $\calI=\{\,i,j\,\}$ be in $J(C)$. Uniformly in the range 
  $m$ in $M_\calI$ and $v$ in $T_m\partial A_1\setminus T_mM_\calI$ such
  that $\sqrt t (|m_i|\wedge |m_j|)$ tends to infinity, 
  $v_i=o(m_i)$ and $v_j=o(m_j)$,
  $$\displaylines{\quad
      \pisa \big( \sqrt t p_\calI (m,v)\big)
      = \sum_{k\in\calI} \sign (m_k)\Big( Q(tm_k^2) 
      + o\big(\log (tm_k^2)\big)^{-1/2}\Big) e_k 
    \cut
      + \pisa (\Proj_{V_\calI^\perp} v) \, .
    \quad\cr}%
  $$}%
\finetune{
%\bigskip
\vfill\eject
}
\stateit{Proof.} In the given range, for $k=i,j$,
$$
  \< \sqrt t p_\calI (m,v),e_k\>
  = \sqrt t m_k \big( 1+o(1)\big)
$$
tends to infinity. We can apply Lemma A.1.5 to obtain an
asymptotic expansion of $\pisa \big(\<\sqrt t p_\calI
(m,v),e_k\>\big)$. This gives the result. \hfill$\qed$

\bigskip

In particular, for the points considered in Lemma 8.3.4, as
$t$ tends to infinity,
\setbox100=\hbox{$I\Big(\pisa \big( \sqrt t p_\calI(m,v)\big)\Big)$}
\wd100=0pt
$$\eqalign{
  \box100\hskip .8in & \cr
  ={} & \sum_{k\in \calI} 
    \Big( \alpha\log\big( \sqrt t |m_k|\big)
    -{1\over 2} \log\log\big(\sqrt t |m_k|\big)\Big) \cr
  & \qquad - 2 \log (K_{s,\alpha}\alpha^{\alpha/2}2\sqrt\pi )
    + \log (2\pi )^{d/2} \cr
  \noalign{\vskip .05in}
  & \qquad + {1\over 2} |\pisa (\Proj_{V_\calI^\perp} v)|^2
    + o(1) \, .\cr
  }
  \eqno{(8.3.2)}
$$

As we explained above, Lemmas 8.3.3 and 8.3.4 give us a suitable 
approximation of $\partial B_t$ near points
parametrized by $m=m_ie_i+m_je_j$ in $M_{i,j}$ with 
$|m_i|\wedge |m_j|\geq t^{-1/6}(\log t)^{M/\alpha}$.
Those points have two components greater than
$t^{-1/6}(\log t)^{M/\alpha}\big( 1+o(1)\big)$. If now a point 
$p$ in $\partial A_1$ has only one component larger than 
$t^{-1/6}(\log t)^{M/\alpha}$ and is in a $o(1)$-neighborhood of
$\bigcup_{\calI\in J(C)}M_\calI$, then one component, say $p_i$,
has to blow up with $t$. Then, $p$ has to be in a 
$o(1)$-neighborhood of $\bigcup_{k:C_{i,k}\neq 0} M_{i,k}$. 
The following lemma gives some control on the largest component 
of $p_k$ for $C_{i,k}$ nonzero.

\bigskip

\state{8.3.5. LEMMA.} {\it%
  Let $\epsilon$ be a positive number. Let $i$ be such that 
  $C_{i,i}$ is null. Furthermore, let $p$ be a point 
  in $\partial A_1$ and in a $o(1)$-neighborhood of 
  $\bigcup_{k: C_{i,k}\neq 0}M_{i,k}$, as $t$ tends to
  infinity, with $|p_i|\geq\epsilon$.
  Choose $j$ such that
  $$
    |p_j|=\max\{\, |p_k|\, :\, C_{i,k}\neq 0\,\} \, .
  $$
  Then, for $t$ large enough,
  $$
    |p_ip_j|\geq {1\over 8(d\| C\|+\epsilon^{-1})} \, .
  $$
  In particular, for $t$ large enough, the inequality 
  $|p_i|\leq r$ forces
  $$
    |p_j|\geq {1\over 8r(d\| C\|+\epsilon^{-1})} \, .
  $$  }
\finetune{
%\bigskip
\vfill\eject
}
\stateit{Proof.} Since $\bigcup_{k:C_{i,k}\neq 0} M_{i,k}$ is
included in the union of planes generated by $(e_i,e_k)$ for
$C_{i,k}$ nonzero, at most two components of $p$ are not $o(1)$.
Moreover, one component is $p_i$ and the other one is a $p_k$
for $C_{i,k}$ nonzero, that is $p_j$. Write $p=p_ie_i+p_je_j+w$ 
with $w$ orthogonal to both $e_i$ and $e_j$. By definition of $p_j$,
each component of $w$ is smaller than $|p_j|$. From what precedes, 
the norm of $w$ has to be less than $\sqrt d\, |p_j|$, as 
well as $o(1)$. Since $p$ belongs to $\partial A_1$ and $Ce_i$ 
is in $\span\{\, e_k\,:\, C_{i,k}\neq 0\,\}$,
$$\eqalign{
  1= \< Cp,p\>
  & = 2 p_ip_j C_{i,j} + 2p_i\<Ce_i,w\> + 2p_j\< Ce_j,w\>
    + \<Cw,w\> \cr
  \vadjust{\vskip .15in}
  & \leq 2 |p_ip_j|\| C\| + 2\sqrt d \, |p_ip_j|\| C\| + |p_j|o(1)
    + o(1) \cr
  }
$$
as $t$ tends to infinity. Since $|p_i|\geq\epsilon$, we obtain
$$
  1\leq |p_ip_j|\big( 2(1+\sqrt d)\| C\| + 4\epsilon^{-1}\big) 
  +o(1) \, .
$$
This implies the first statement of the lemma. The second follows
trivially.\hfill$\qed$

\bigskip

The lower bounds in Lemma 8.3.5 are useful only if $p_i$ and
$p_j$ are not too large, so that we can have $\sqrt
t(|p_i|\wedge |p_j|)$ going to infinity. Thus we need to shrink
the domain on which we need to perform the integration. This
relies on the simple observation that for Student-like
distributions, a Bonferoni type inequality gives
$$
  P\Big\{\, \max_{1\leq i\leq d} |X_i|\geq 
          {tM\over (\log t)^{1/\alpha}} \,\Big\}
  \leq 2d {\log t\over t^\alpha M^\alpha} K_{s,\alpha}
  \alpha^{(\alpha -1)/2}
  \eqno{(8.3.3)}
$$
for $t$ large enough. Thus, we define the subset\ntn{$A'_t$}
$$
  A'_t=\sqrt t A_1' 
  = \sqrt t \,\Big\{\, p\in\RR^d \, :\, \< Cp,p\>\geq 1 \, , \,
  \max_{1\leq i\leq d} |p_i| \leq 
  {\sqrt t \log\log t\over (\log t)^{1/\alpha}} \,\Big\} \, .
$$
We denote\ntn{$B'_t$} $B'_t=\pisa (\sqrt t A'_t)$.

We can now improve dramatically upon Lemma 8.3.3.

\bigskip

\state{8.3.6. LEMMA.} {\it%
  The set of all $p$'s in $\partial A_1'$ such that
  $$
    I\big( \pisa (\sqrt t p)\big) \leq R(t)+M\log\log t
  $$ 
  can be parameterized as $p(m,v)$ with $m$ in some $\MI$,
  some $\calI$ in  $J(C)$ and
  $$
    \Proj_{V_\calI^\perp} v 
    = O\big( t^{-1/2}(\log t)^{M/(2\alpha)}
             (\log\log t)^{1/(2\alpha)}\big) \, .
  $$
  Moreover, whenever $|m_i|\wedge |m_j|\gg t^{-1/2}
  (\log t)^{M/(2\alpha)}(\log\log t)^{1/(2\alpha)}$, then
  $v_i=o(m_i)$ and $v_j=o(m_j)$. In addition,
  $$
    I\big( \pisa (\sqrt t A'_1)\big) \geq R(t) 
    - \alpha\log\gamma -\alpha \log (16\, d\,\| C\| )
    + o(1)
  $$
  as $t$ tends to infinity.
  }

\bigskip

\stateit{Proof.} Let $p$ be a point in $\partial A_1'$ as in the
statement of the lemma. Lemma 8.3.3 guarantees that it belongs to
an $o(t^{-1/6}(\log )^{M/\alpha})$-neighborhood of
$\bigcup_{\calI\in J(C)}M_\calI$. If it is in an $o(t^{-1/6}(\log
t)^{M/\alpha})$-neigh\-bor\-hood of $\bigcup_{\calI\in J(C)\setminus
J_*(C)}M_\calI$, then we are done thanks to Lemmas 8.3.4 and
8.2.11. Thus, assume that it is not. Let $i$ be such that
$$
  |p_i|=\max\{\, |p_k|\, :\, 1\leq k\leq d\,\} \, .
$$
Since $1=\< Cp,p\>$ is less than $\| C\| |p|^2$, we have 
$$
  |p_i|\geq 1/d\| C\|\gg t^{-1/6}(\log t)^{M/\alpha} \, .
$$ 
Therefore, $p$ is in a
$o\big(t^{-1/6}(\log t)^{M/\alpha}\big)$-neighborhood of
$\bigcup_{k:C_{i,k}\neq 0} M_{i,k}$. Since $p$ belongs to $A_1'$,
$$
  |p_i|\leq {\sqrt t\log\log t\over (\log t)^{1/\alpha}} \, .
$$
Let $p_j$ be as in Lemma 8.3.5. The second inequality in 
Lemma 8.3.5 gives
$$
  |p_j|\geq {(\log t)^{1/\alpha}\over \sqrt t \log\log t}
  {1\over 16 \, d\, \| C\| } \, .
$$
In particular, both $\sqrt t |p_i|$ and $\sqrt t |p_j|$ tend to
infinity. Write $p=p_ie_i + p_je_j+w$ with $w$ orthogonal to
both $e_i$ and $e_j$. We calculate
$$\eqalign{
  I\big( \pisa (\sqrt t p)\big)
  & = {1\over 2}|\pisa (\sqrt tp_i)|^2 
    +{1\over 2} |\pisa (\sqrt t p_j)|^2 \cr
  & \qquad\qquad  + {1\over 2}|\pisa (\sqrt t w)|^2 
    + \log (2\pi )^{d/2} \cr
  & =\alpha \log (\sqrt t |p_i|) 
    - {1\over 2} \log\log (\sqrt t|p_i|)
    + \alpha\log (\sqrt t |p_j|) \cr
  & \qquad\qquad  - {1\over 2} \log\log (\sqrt t |p_j|) 
    - 2\log (K_{s,\alpha}\alpha^{\alpha/2}2\sqrt\pi ) \cr
  \noalign{\vskip 1mm}
  & \qquad\qquad  + {1\over 2} |\pisa (\sqrt t w)|^2 
    +\log (2\pi )^{d/2} + o(1) \cr
  }
$$
as $t$ tends to infinity. From the first estimate in Lemma 8.3.5, we
deduce
$$
  \alpha \log (\sqrt t |p_i|) +\alpha\log (\sqrt t |p_j|)
  \geq 2\alpha \log \sqrt t - \alpha \log ( 16\, d\,\| C\| )
  \, .
$$
Moreover, since $\log (1+h)\leq h$ for any $h> -1$, we also have
$$
  \log\log (\sqrt t |p_i|)
  = \log (\log \sqrt t +\log |p_i|)
  \leq \log\log \sqrt t + {\log |p_i|\over \log \sqrt t} \, .
$$
Thus,
$$\eqalign{
  I\big(\pisa (\sqrt t p)\big)\geq\,
  & 2\alpha\log\sqrt t -\log\log\sqrt t 
    - 2\log (K_{s,\alpha}\alpha^{\alpha/2}2\sqrt\pi ) \cr
  & \qquad + \log (2\pi )^{d/2} 
    -\alpha \log (16d\| C\| ) 
    +o(1) \cr
  & \qquad + {1\over 2}|\pisa (\sqrt t w)|^2 \, , \cr
  }
$$
where the $o(1)$-term does not depend on $p$ in $A_1'$.  This gives
the lower bound on $I\big(\pisa (\sqrt t A_1')\big)$ stated in the
Lemma. Moreover, the inequality
$$
  I\big(\pisa (\sqrt tp)\big)\leq R(t) + M\log\log t
$$
implies
$$
  |\pisa (\sqrt t w)|^2 \leq M\log \log t + c
$$
for some constant $c$. This forces, for any $k$ not in 
$\{\, i,j\,\}$,
$$
  \sqrt t |w_k|\leq S_\alpha^\leftarrow\!\circ\phi 
  (\sqrt{M\log\log t + c})
  = (\log t)^{M/(2\alpha)}(\log\log t)^{1/(2\alpha)} O(1)
$$
as $t$ tends to infinity, thanks to Lemma A.1.6.

Set $\calI=\{\, i,j\,\}$. We just showed that the projection of
$p$ on $V_\calI^\perp$ is at most of order $t^{-1/2}(\log
t)^{M/(2\alpha)}(\log\log t)^{1/(2\alpha )}$ as $t$ tends to
infinity. From Lemma 8.3.3, we then know that $p=p(m,v)$ for
some $m$ in $\MI$ and $v$ in $T_m\partial A_1\setminus T_m\MI$
with $|\Proj_{V_\calI} v|=o\big(t^{-1/6}(\log
t)^{M/\alpha}\big)$. Therefore, if $|m_i|$ and $|m_j|$ are at
least $t^{-1/6}(\log t)^{M/\alpha}$, the lemma is proved. Our
choice of $i$ implies that the only case left to investigate is
$|m_j|\leq t^{-1/6}(\log t)^{M/\alpha}$. In this case, Lemma
8.3.2 shows that $C_{i,i}$ must vanish --- because $m_i$ is lower
bounded if $C_{j,j}$ is not null in Lemma 8.3.2; this is what we are
claiming up to a permutation of $i$ and $j$. Then, Lemma 8.3.2
shows that $m_i\sim 1/(2C_{i,j}m_j)$ as $t$ tends to infinity.
Since $v$ is in $T_m\partial A_1\ominus T_m\MI$, it is
orthogonal to the vector spanning $T_m\MI$, that is, to
$$
  {1\over 2C_{i,j}} \Big( -{1\over m_j^2}-C_{j,j}\Big) e_i
  + e_j \, .
$$
Thus, there exists a real $r$ such that
$$
  \Proj_{V_\calI} v = re_i + {r\over 2C_{i,j}}
  \Big({1\over m_j^2}+C_{j,j}\Big) e_j \, .
  \eqno{(8.3.4)}
$$
Since $v$ belongs to $T_m\partial A_1$, we must have 
$\< Cm,v\>=0$, or, equivalently,
$$
  \< \Proj_{V_\calI}Cm,\Proj_{V_\calI}v\>
  = \< Cm,\Proj_{V_\calI}v\>
  = -\< Cm,\Proj_{V_\calI^\perp} v\> \, .
  \eqno{(8.3.5)}
$$
Since $C_{i,i}$ is null,
\finetune{\hfuzz=3pt}
\setbox100=\hbox to -2.5pt{}
$$\eqalign{\box100
  \< \Proj_{V_\calI}Cm,\Proj_{V_\calI}v\>
  & = \Big\< {1\over 2} \Big( {1\over m_j}-m_jC_{j,j}\Big) e_j
    + m_j C_{i,j} e_i + m_j C_{j,j}e_j \, ,\cr
  & \kern 1.35in  re_i + {r\over 2C_{i,j}}
    \Big({1\over m_j^2}+C_{j,j}\Big) e_j\Big\> \cr
  & = r_jC_{i,j} + {r\over 4C_{i,j}}
    \Big( {1\over m_j}+m_jC_{j,j}\Big) 
    \Big( {1\over m_j^2}+C_{j,j}\Big) \, . \cr
  }
$$
\finetune{\hfuzz=0pt}
Thus, if $m_j$ tends to $0$ as $t$ tends to infinity,
$$
  \< \Proj_{V_\calI}Cm,\Proj_{V_\calI}v\>
  \sim {r\over 4C_{i,j}m_j^3} \, .
$$
On the other hand, using the bound on the projection of $p$
on $V_\calI^\perp$ that we obtained earlier in the proof
of this lemma,
$$\eqalign{
    |\< Cm,\Proj_{V_\calI^\perp}v\>|
  & \leq \| C\| \,  |m| \, |\Proj_{V_\calI^\perp}(v)| \cr
  & = {1\over |m_j|}O\big( t^{-1/2}(\log t)^{M/(2\alpha)}
    (\log\log t)^{1/(2\alpha )}\big) \, . \cr
  }
$$
Therefore, (8.3.5) yields
$$
  r = m_j^2  O\big( t^{-1/2}(\log t)^{M/(2\alpha)} 
  (\log\log t)^{1/(2\alpha )}\big) \, .
$$
as $t$ tends to infinity. Therefore, going back to (8.3.4), 
$$
  \Proj_{V_\calI} v
  = o(m_i)e_i 
  + \Big( O\big( t^{-1/2}(\log t)^{M/(2\alpha)} 
          (\log\log t)^{1/(2\alpha )}\big) + o(m_j)\Big)
  e_j .
$$
This concludes the proof of Lemma 8.3.6.\hfill$\qed$

\bigskip

Though we do not need this right now, Lemmas 8.3.3 and 8.3.6 imply
$$
  I(B'_t)\sim R(t) \sim \alpha\log t \qquad\hbox{ as } \ttoi 
  \, .
$$
More importantly, if $i$ belongs to $\calI=\{\,i,j\,\}\in J_*(C)$, and 
$m$ is in $M_\calI$, Lemma 8.3.2 shows that $m_j\sim 1/2C_{i,j}m_i$ as
$|m_i|$ tends to infinity. Therefore, Lemma 8.3.6 allows 
us to use the expansion of Lemma 8.3.4 in the range 
$|m_i|\ll \sqrt t /(\log t)^{M/(2\alpha)}(\log\log t)^{1/(2\alpha)}$. 

To conclude what we have done so far, we can use Lemma 8.3.4 over all
points $q=\pisa (\sqrt t p)$ in $B_t$, with $p$ belonging to $A'_t$,
except for those in an $O\big( t^{-1/6}(\log t)^{M/\alpha}\big)$-neighborhood 
of a point $m=m_ie_i+m_je_j$ with $\{\,i,j\,\}$ in $J_*(C)$ and, say,
$$
  \sqrt t/(\log t)^{M/(2\alpha)}(\log\log t)^{1/\alpha}
  \leq |m_i| \vee |m_j|\leq M\sqrt t / (\log t)^{1/\alpha}
  \, . 
$$

We now discard this missing
range by an ad hoc argument. The key observation is that Lemma
8.3.5 guarantees that this set of $p$'s is included in\ntn{$\Omega_t$}
$$\displaylines{
    \Omega_t = \bigcup_{\{i,j\}\in J_*(C)}\Big\{\,
    p\in \RR^d \, :\, 
    { \sqrt t \over (\log t)^{M/\alpha} 
      (\log\log t)^{1/(2\alpha)} }
    \leq |p_i|\vee |p_j|
  \cut
    \leq {\sqrt t \log\log t\over (\log t)^{1/\alpha}} \, ; \,
    |p_ip_j|\geq {1\over 16\, d\,\| C\|}
    \,\Big\} \, . \cr
  }
$$

\bigskip

We now show that ultimately we will be able to discard the set
$\Omega_t$ in our computation.

\bigskip

\state{8.3.7. LEMMA.} {\it%
  We have
  $$
    P\{\, X\in \sqrt t \Omega_t\,\}
    = o\Big( {\log t\over t^\alpha}\Big) \qquad \hbox{ as }
    \ttoi \, .
  $$
  }

\bigskip

\stateit{Proof.} The upper and lower tail of the Student-like
distributions are asymptotically equivalent. The result is then
a consequence of the independence of the $X_i$'s and the
following calculation,
$$\displaylines{
    P\Big\{\, 
    {t\over (\log t)^{M/(2\alpha )}(\log\log t)^{1/\alpha}}
    \leq X_i\leq {t\log\log t\over (\log t)^{1/\alpha}} \, ;
    \, X_iX_j\geq {t\over 16 d \| C\|} \,\Big\}
  \hfill\cr\vadjust{\vskip .15in}\hfill
  \eqalign{
    ={} & \int_{ t(\log t)^{-M/(2\alpha)}(\log\log t)^{-1/\alpha}
             }^{ t\log\log t (\log t)^{-1/\alpha} }
        1-S_\alpha (t/16 d\| C\| x) \,\d S_\alpha (x) \cr
    ={} & O(1)\int_{ t(\log t)^{-M/(2\alpha)}(\log\log t)^{-1/\alpha}
             }^{ t\log\log t (\log t)^{-1/\alpha} }
        {x^\alpha\over t^\alpha} \,\d S_\alpha (x) \, . \cr
    }
  \cr
  }
$$
An integration by parts shows that the last integral is 
$$\displaylines{\qquad
    {1\over t^\alpha} \big[ x^\alpha \big( S_\alpha(x)-1)\big)
    \big]_{ t(\log t)^{-M/(2\alpha)}(\log\log t)^{-1/\alpha}
          }^{ t\log\log t (\log t)^{-1/\alpha} }
  \cut
    + {\alpha\over t^\alpha} 
    \int_{ t(\log t)^{-M/(2\alpha)}(\log\log t)^{-1/\alpha}
               }^{ t\log\log t (\log t)^{-1/\alpha} }
    x^{\alpha-1}\big( 1-S_\alpha (x)\big)\, \d x \, .
  \qquad\cr}
$$
Since $x^\alpha \big( 1-S_\alpha (x)\big)$ tends to a constant
as $x$ tends to infinity, the above sum is
$$\displaylines{\qquad
    o(t^{-\alpha}) + O(1)t^{-\alpha} 
    \int_{ t(\log t)^{-M/(2\alpha)}(\log\log t)^{-1/\alpha}
         }^{ t\log\log t (\log t)^{-1/\alpha} }
    x^{-1} \d x\vtop to.2in{\relax}
  \cut
    = o(t^{-\alpha}) + O(1)t^{-\alpha} \log\log t
    = o( t^{-\alpha}\log t)
    \qquad\cr
  }
$$
as $t$ tends to infinity.\hfill$\qed$

\bigskip

Define\ntn{$B''_t$}
$$
  B''_t =B'_t\setminus \pisa (\sqrt t \Omega_t) \, .
$$
Combining (8.3.3) and Lemma 8.3.7, we obtain
$$
  P\{\, \cXX \geq t\,\} = \int_{B''_t}e^{-I(y)} \d y
  + o\Big( {\log t\over t^\alpha}\Big) \, .
$$

From now on, we concentrate on the integral of $e^{-I}$ on
$B''_t$. Combining Lemmas 8.3.4, 8.3.5 and 8.3.6, the set
$B''_t\bigcap\Gamma_{R(t)+M\log\log t}$ can be identified, for
$t$ large enough, with some points 
$\pisa \big(\sqrt t p_\calI(m,v)\big)$, with
\vskip\abovedisplayskip
%
%\halign{\qquad\qquad $#$\hfill&\qquad$#$\hfill\qquad \cr
%  \calI\in J(C) \, , & m\in M_\calI \, , \cr
%  \noalign{\vskip 1mm}
%  v\in T_m\partial A_1\setminus T_mM_\calI \, ,
%  & |\Proj_{V_\calI^\perp}|\leq t^{-1/2}(\log t)^{2M/\alpha}
%    (\log\log t)^{1/(2\alpha)} \, , \cr
% }
%\vskip .1in
%\line{\qquad\qquad $\min_{k\in \calI}|m_k|\geq t^{-1/2}
%      (\log\log t)^{1/\alpha} (\log t)^{M/(2\alpha )}$\hfill}
%\vskip .1in
%\line{\rlap{and}\qquad\qquad 
%      $|\Proj_{V_\calI}v| = o(m)$ componentwise.\hfill}
\line{\qquad\qquad$\calI\in J(C)$, \qquad $m\in M_\calI$,\qquad
      $v\in T_m\partial A_1\setminus T_mM_\calI$ \hfill}
\vskip .05in
\line{\qquad\qquad$|\Proj_{V_\calI^\perp}|\leq 
      t^{-1/2}(\log t)^{2M/\alpha}(\log\log t)^{1/(2\alpha)}$,
      \hfill}
\vskip .05in
\line{\qquad\qquad $\min_{k\in \calI}|m_k|\geq t^{-1/2}
      (\log\log t)^{1/\alpha} (\log t)^{M/(2\alpha )}$,\hfill}
\vskip .1in
\line{\rlap{and}\qquad\qquad 
      $\Proj_{V_\calI}v = o(m)$ componentwise.\hfill}
\vskip\belowdisplayskip

For those points, we can
use Lemma 8.3.4 and (8.3.2). This allows us to obtain the value
of $I(B''_t)$ up to an $o(1)$-term as $t$ tends to infinity.

\bigskip

\state{8.3.8. LEMMA.} {\it%
  The equality $I(B''_t)=R(t)+o(1)$ holds as $t$ tends to 
  infinity.
  }

\bigskip

\stateit{Proof.} Making use of formula (8.3.2) we first need to
evaluate the minimum of the function
$$
  \sum_{k\in\cal I}\Big(\alpha\log (\sqrt t|m_k|)-{1\over 2}
  \log\log (\sqrt t |m_k|)\Big)
$$
on $M_\calI$. We also need an approximate location of the
minimum, to check that $\pisa (\sqrt t\,\cdot)$ maps 
it to a point in $B''_t$ --- and not in $B_t$. When $\calI$ 
is in $J(C)\setminus J_*(C)$, we can argue as in (8.2.9). 
We concentrate on the case where $\calI =\{\,i,j\,\}$ 
belongs to $J_*(C)$. Considering Lemma 8.3.2, 
this leads us to define the function
$$\displaylines{
    f_{t,C_{j,j}}(m_j) =
    \alpha\log\Big( {\sqrt t\over 2|C_{i,j}|}
    \Big| {1\over m_j}-m_jC_{j,j}\Big|\Big)
    +\alpha \log (\sqrt t |m_j|)
  \hfill\cr\noalign{\vskip .05in}\hfill
    - {1\over 2} \log\log \Big( {\sqrt t\over 2|C_{i,j}|}
    \Big| {1\over m_j}-m_jC_{j,j}\Big|\Big)
    - {1\over 2}\log\log (\sqrt t |m_j|) \, . 
  \cr
  }
$$
When $C_{j,j}$ vanishes, this function is
$$\displaylines{\qquad
  f_{t,0}(m_j)=\alpha\log t - \alpha\log (2|C_{i,j}|)
  - {1\over 2}\log\log {\sqrt t\over |C_{i,j}m_j|}
  \cut
  - {1\over 2} \log\log (\sqrt t |m_j|) \, .
  \qquad\cr}
$$
The function $f_{t,0}(\cdot )$ is minimum when $r=\log |m_j|$ 
maximizes
$$
  (\log \sqrt t - \log |C_{i,j}|-r)(\log \sqrt t + r) \, ,
$$
that is
$$
  r= -{1\over 2}\log | C_{i,j}| \, .
$$
Moreover, as $t$ tends to infinity,
$$
  f_{t,0}\Big( {1\over \sqrt{|C_{\smash{i,j}}|}}\Big)
  = \alpha\log t -\alpha\log (2|C_{i,j}|) 
  -\log\log \sqrt t + o(1) \, .
$$
When $m_j= 1/\sqrt {|C_{\smash{i,j}}|}$ and $C_{j,j}$ vanishes,
Lemma 8.3.2 gives  $m_i=\sign (C_{i,j})/(2\sqrt{|C_{\smash{i,j}}|})$. 
Clearly $m_ie_i+m_je_j$ is mapped into $B''_t$ by 
$\pisa (\sqrt t\,\cdot )$.

If $C_{j,j}$ is nonzero and $m_j$ is positive and fixed, then
$$\displaylines{\quad
  f_{t,C_{j,j}}(m_j) = \alpha\log t -\alpha\log ( 2|C_{i,j}|)
  + \alpha \log (1-m_j^2 C_{j,j}) 
  \cut 
  - \log\log \sqrt t + o(1)
  \quad (8.3.6)\cr}
$$
--- recall that $C_{j,j}$ is negative. Thus, for
$f_{t,C_{j,j}}(m_j)$ to be minimum, we must have $m_j^2$ tends
to $0$ as $t$ tends to infinity. But if 
$\sqrt t p_\calI (m,v)$ belongs to $B''_t$, expansion
(8.3.6) always holds since $|m_i|\wedge |m_j|\gg 1/\sqrt t$.
Taking $|m_j|=1/\log t$ for instance, we have
$$
  f_{t,C_{j,j}}(1/\log t) = \alpha\log t 
  - \alpha \log (2|C_{i,j}|) - \log\log \sqrt t + o(1) 
  \, .
$$
as $t$ tends to infinity. Given the uniformity over 
$B''_t$ in (8.3.2), this gives the
asymptotic minimum of $f_{t,C_{j,j}}$, and ultimately
$I(B''_t)$. \hfill$\qed$

\bigskip

We can now apply Theorem 5.1. Denote by\ntn{$\rho_t$}
$$
  \rho_t=\sqrt{2I(B''_t)-\log(2\pi)^{d/2}}
$$
the radius of the ball $\Lambda_{I(B_t)}$. As we have seen after
the proof of Lemma 8.3.6,
$$
  \rho_t\sim \sqrt{2\alpha\log t} \qquad\hbox{ as } \ttoi \, .
$$
In view of Lemma 8.3.4, for $\calI$ in $J(C)$ and $m$
belonging to $M_\calI$, define\ntn{$q_{\calI,t}(m)$}
$$
  q_{\calI,t}(m)=\sum_{k\in \calI}\sign (m_k)Q(tm_k^2)e_k
$$
and its projection onto $\Lambda_{I(B_t)}$ through the normal
flow,\ntn{$r_{\calI,t}(m)$}
$$
  r_{\calI,t}(m)=\rho_t
  { q_{\calI,t}(m)\over |q_{\calI,t}(m)| } \, .
$$
We consider the dominating manifold of dimension $k=1$,\ntn{$\cal_{B''_t}$}
$$\displaylines{\qquad
    \calD_{B''_t} = \bigcup_{\calI\in J(C)}\Big\{\,
    r_{\calI,t}(m) \, :\, m\in \bigcup_{\calI\in J(C)}M_\calI
    \, ; \, m\not\in\Omega_t \, ; \,
  \cut
    \max_{1\leq i\leq d} |m_i|\leq 
               {\sqrt t \log\log t\over (\log t)^{1/\alpha}} \, ,
    \tau_{B''_t}\big( r_{\calI,t}(m)\big)\leq
               \log\log t \Big\} \, . 
  \qquad\cr 
  }
$$
It is fairly clear that $\calD_{B''_t}$ is a dominating manifold
for the set $B''_t\bigcap \Gamma_{I(B''_t)+M\log\log t}$. We can
take $M$ large enough such that
$$
  \int_{\Gamma_{R(t)+M\log\log t}^{\rm c}} e^{-I(y)} \d y
  = o\Big( {\log t\over t^\alpha}\Big) \qquad\hbox{ as }
  \ttoi \, .
$$
When applying the formula given in Theorem 5.1, the integral over
$\calD_{B''_t}$ splits into two parts. The first part comes from the
contribution of points $r_{\calI,t}(m)$ with $\calI$ belonging to
$J(C)\setminus J_*(C)$. This part is exactly like the integral we
dealt with in the proof of Theorem 8.2.10. It is of order
$1/t^\alpha$. Thus, we concentrate on the second part, coming from
points $r_{\calI,t}(m)$ with $\calI$ in $J_*(C)$.

Let $\calI=\{\,i,j\,\}$ be in $J_*(C)$ and
$r_{\calI,t}(m)$ be a point of $\calD_{B''_t}$. From (8.3.2) 
we infer
$$\eqalign{
  \tau_{B''_t}\big( r_{\calI,t}(m)\big)\,
  & {= {1\over 2} \big( |q_{\calI,t}|^2-\rho_t^2\big)} \cr
  & {= \alpha\log |m_im_j| -\alpha\log\gamma}%
    - {1\over 2}\log{\log (\sqrt t|m_i|)\over \log\sqrt t} \cr
  & \phantom{ = \alpha\log |m_im_j| -\alpha\log\gamma } {}
    - {1\over 2} \log {\log (\sqrt t|m_j|)\over \log\sqrt t}
    + o(1) \cr
  }
$$
Consequently,
$$\displaylines{\qquad
    \exp\Big(-\tau_{B''_t}\big( r_{\calI,t}(m)\big)\Big)
  \cut
    = {\gamma^\alpha\over |m_im_j|^\alpha}
    \sqrt{\Big( 1 + {\log |m_i|\over\log\sqrt t}\Big)
          \Big( 1 + {\log |m_j|\over\log\sqrt t}\Big)}
    \big( 1+ o(1)\big)
  \qquad\cr}
$$
uniformly over the part of $\calD_{B''_t}$ corresponding to
points $r_{\calI,t}(m)$ with $m$ in $M_\calI$. Notice that if
$C_{j,j}$ is nonzero, the inequality
$$
  \tau_{B''_t}\big( r_{\calI,t}(m)\big)
  \leq M\log\log t
$$
required in $\calD_{B''_t}$ imposes $|m_j|\leq (\log
t)^{2M/\alpha}$ for $t$ large enough; this can be seen from the
above expression of $\tau_{B''_t}\big( r_{\calI,t}(m)\big)$ and
using the same lower bound argument as in the proof of Lemma
8.3.6 to handle the term in $\log\big( \log (\sqrt t |m_i|)/\log
\sqrt t \big)$. Therefore, when applying Theorem 5.1 to a
component coming from $M_{i,j}$ with $C_{j,j}$ nonzero, we need
only to integrate for $|m_j|\leq (\log t)^{2M/\alpha}$. 

Still in order to apply the formula in Theorem 5.1, on
$\calD_{B''_t}$,
$$
  \big|\D I\big( r_{\calI,t}(m)\big)\big|
  = |r_{\calI,t}(m)|
  \sim \rho_t \sim \sqrt{2\alpha\log t}
  \qquad \hbox{ as } \ttoi \, .
$$

We can now calculate the matrix $G_{\smash{B''_t}}\big( r_{\calI,t}(m)\big)$. In
the range of $m$'s that we are considering in $\calD_{B''_t}$,
Lemma 8.3.4 shows that $\partial B''_t$ behaves like a ruled
surface made of $(d-2)$-dimensional flat subspaces in directions
orthogonal to $V_\calI$. Hence, the very same argument as in
Lemma 8.2.8 shows that
$$
  G_{B''_t}\big( r_{\calI,t}(m)\big)
  \sim {\Id_{\RR^{d-2}}\over \sqrt {2\alpha\log t} }
  \qquad \hbox{ as } \ttoi \, ,
$$
uniformly over the part of $M_\calI$ we are considering. We can
then write explicitly the contribution of the part related to
$M_\calI$ in the formula given in Theorem 5.1. This contribution
is
$$
  e^{-I(B''_t)} (2\pi )^{d-2\over 2}
  \int {e^{-\tau_{B''_t}}\over |DI|^{d/2} (\det\, G_{B''_t})^{1/2}}
  \,\d\calM_{r_{\calI,t}(M_\calI)} \, ,
$$
where we integrate over all points $r_{\calI,t}(m)$, for $m$ in
$M_\calI$, with $|m_i|\vee |m_j|\leq \sqrt t \log\log t / (\log
t)^{1/\alpha}$; and $|m_j|\leq (\log t)^{2M/\alpha}$ if
$C_{j,j}$ is nonzero. As in Theorem 8.2.10, we can replace the
Riemannian measure on $r_{\calI,t}(M_\calI)$ by that on
$q_{\calI,t}(M_\calI)$. Its expression in the local
parameterization given by $m_j$ in Lemma 8.3.2 is
$$
  \bigg(\Big( {\d\over \d m_j}Q(tm_i^2)\Big)^2
  + \Big( {\d\over \d m_j}Q(tm_j^2)\Big)^2\bigg)^{1/2} \d m_j \, .
$$
Since $Q'(x)\sim \sqrt\alpha /(2x\sqrt{\log x})$ as $x$ tends to
infinity, this expression is equivalent to
$$
  \bigg( {\alpha\over m_i^2\log (tm_i^2)}
  \Big({dm_i\over dm_j}\Big)^2
  + {\alpha\over m_j^2\log (tm_j^2)}\bigg)^{1/2} \d m_j\, .
$$
Putting all the pieces together, the contribution to the
integral in Theorem 5.1 coming from $M_\calI$ is
$$\displaylines{\qquad
    {\log\sqrt t\over t^\alpha} 
    {K_{s,\alpha}^2\alpha^\alpha 4\pi\over (2\pi )^{d/2}
     \gamma^\alpha}
    (2\pi )^{(d-2)/ 2} \times
  \hfill\cr\hskip .5in
    \int{\gamma^\alpha\over |m_im_j|^\alpha}
    {\sqrt{ \Big( 1+{\displaystyle\log |m_i|\over
                     \displaystyle\log\sqrt t}\Big)
            \Big( 1+{\displaystyle\log |m_j|\over
                     \displaystyle\log\sqrt t}\Big) } \over
     (2\alpha\log t)^{d/4}(2\alpha\log t)^{-(d-2)/4} }\times
  \hfill\cr\noalign{\vskip .05in}\hfill
    \bigg( {\alpha\over m_i^2\log (tm_i^2)}
    \Big( {\d m_i\over \d m_j}\Big)^2 
    + {\alpha\over m_j^2\log (tm_j^2)}\bigg)^{1/2}
    \d m_j \, ,\qquad (8.3.7) \cr
  }
$$
where we integrate over $m_j$ such that
$$
  { (\log t)^{M/(2\alpha)}(\log\log t)^{1/\alpha}\over 
    2|C_{i,j}|\sqrt t } \big( 1+o(1)\big)
  \leq |m_j| \leq (\log t)^{2M/\alpha}
$$
if $C_{j,j}$ is nonzero, and
$$\displaylines{\qquad
    { (\log t)^{M/(2\alpha)}(\log\log t)^{1/\alpha}\over 
      2|C_{i,j}|\sqrt t }
    \big( 1+o(1)\big)
    \leq |m_j| 
  \cut
    \leq {\sqrt t \over 
      (\log t)^{M/(2\alpha)}(\log\log t)^{1/\alpha}}
  \qquad\cr}
$$
if $C_{j,j}$ is null.
The finale of the proof consists in showing how simple this
integral is, at least asymptotically! Up to multiplying it 
by $2$, we can restrict the 
range of integration to $m_j$ positive.

We first integrate in the range
$$
  {(\log t)^{M/(2\alpha )}(\log\log t)^{1/\alpha}
  \big( 1+o(1)\big) \over 2|C_{i,j}|} 
  \leq m_j\leq {1\over \log t} \, .
$$
In this range, $m_i\sim 1/(2C_{i,j}m_j)$, and
$$
  {1\over m_i^2}\Big( {\d m_i\over \d m_j}\Big)^2
  \sim {1\over m_j^2} \, .
$$
Using the change of variable $m_j=t^s$, the integral becomes 
--- on that range ---
$$\displaylines{\qquad
      2 {\log\sqrt t\over t^\alpha} 
      {K_{s,\alpha}^2\alpha^\alpha 2 \over \sqrt{2\alpha\log t} }
      \int |2C_{i,j}|^\alpha \big( 1+o(1)\big)
      \sqrt{ \big( 1-2s+o(1)\big)(1+2s)}\times
  \cut
      \qquad \sqrt\alpha \, {1\over t^s\sqrt{\log t}}
      \Big( {1\over 1-2s+o(1)} + {1\over 1+2s}\Big)^{1/2}
      \log t \; t^s \d s
  \cr\noalign{\vskip 0.05in}\qquad
      = 4 {\log\sqrt t\over t^\alpha}
      {K_{s,\alpha}^2\alpha^\alpha\over \sqrt 2}
      |2C_{i,j}|^\alpha
      \int \sqrt{2+o(1)}\, \d s \, , \hfill\cr
  }
$$
where we integrate for 
$$
  -{1\over 2}-{M\log\log t\over 2\alpha\log t}\big( 1+o(1)\big) 
  \leq s\leq -{\log\log t\over \log t} \, .
$$ 
This gives a term
$$
  2{\log\sqrt t\over t^\alpha} K_{s,\alpha}^2 \alpha^\alpha
  |2C_{i,j}|^\alpha \, .
$$

When $1/\log t\leq m_j\leq (\log t)^{2M/\alpha}$, then $\log
|m_i|$ and $\log |m_j|$ are $O(\log\log t)$. This part of the
integral contributes a term less than
$$\displaylines{\qquad
    {\log\sqrt t\over t^\alpha} O(1)
    \int {2^\alpha |C_{i,j}|^\alpha\over (1-m_j^2C_{j,j})^\alpha}
    {(\log\log t)^{1/2}\over (\log t)^{1/2} }
    {1\over \sqrt{\log t}} \times
  \cut
    \bigg( {1\over m_i^2}\Big( {\d m_i\over \d m_j}\Big)^2
    + {1\over m_j^2}\Big)^{1/2} \d m_j \, .
  \qquad\cr}
$$
Since $C_{j,j}$ is nonpositive, $1-m_j^2C_{j,j}\geq 1$. Therefore
$$
  \Big|{1\over m_i} {\d m_i\over \d m_j}\Big| 
  = \Big|{1\over m_j}{1+m_j^2C_{j,j}\over 1-m_j^2C_{j,j}}\Big|
  \leq {1\over m_j} \, .
$$
The contribution in the integral is of order at most
$$
  O\Big( {\log\log t\over t^\alpha}\Big)
  \int_{1/\log t}^{(\log t)^{2M/\alpha}} {\d m_j\over m_j}
  = O\Big( {(\log\log t)^2\over t^\alpha}\Big)
  = o\Big( {\log t\over t^\alpha}\Big) \, .
$$
Consequently, if $C_{j,j}$ is not zero, (8.3.7) is equivalent to
$$
  2 {\log \sqrt t\over t^\alpha} K_{s,\alpha}^2
  \alpha^\alpha |2C_{i,j}|^\alpha \, ,
$$
as $t$ tends to infinity.

If $C_{j,j}$ vanishes, we need to add a contribution for the part where
$$
  (\log t)^{2M/\alpha}\leq m_j\leq 
  {\sqrt t\over (\log t)^{M/(2\alpha)}(\log\log t)^{1/\alpha}} 
  \, .
$$
We argue as in the case $m_j\leq 1/\log t$ above --- this amounts to
exploit the symmetry between $m_i$ and $m_j$ when $m_j$ is
large. Therefore, if $C_{j,j}$ is null, (8.3.7) is equivalent to
$$
  2{\log \sqrt t\over t^\alpha} K_{s,\alpha}^2\alpha^\alpha
  |2C_{i,j}|^\alpha \, ,
$$
as $t$ tends to infinity.
Putting all the estimates together, we obtain
$$\displaylines{
    \int_{B''_t}e^{-I(y)} dy
    \sim 2{\log \sqrt t\over t^\alpha} K_{s,\alpha}^2 \alpha^\alpha
    \times
  \cut
    \sum_{(i,j)\in J_*(C)} |2C_{i,j}|^\alpha
    \big( \indic_{\RR\setminus \{ 0\}}(C_{j,j})
      + 2\indic_{\{ 0\} }(C_{j,j}) \big)
  \quad(8.3.8)\cr}
$$
Noticing that
$$
  \sum_{\{ i,j\}\in J_*(C)} |2C_{i,j}|^\alpha
  \big( \indic_{\RR\setminus \{ 0\} }(C_{j,j})
        + 2\indic_{\{ 0\} }(C_{j,j}) \big)
  = \sum_{i:C_{i,i}=0}\sum_{1\leq j\leq d} |2C_{i,j}|^\alpha 
  \, , 
$$
we see that (8.3.8) is the expression given in the statement of
Theorem 8.3.1 since we replaced $C$ by $(C+C^\T )/2$ in this proof.

To conclude the proof, we need to check the assumptions of
Theorem 5.1. This part of the proof of Theorem 8.2.10 can be
copied almost word for word, and this concludes the proof of
Theorem 8.3.1.\finetune{~}\hfill$\qed$

\bigskip

Combining Theorem 8.2.21 and Lemmas 8.3.2, we obtain
the following result.

\bigskip

\state{8.3.9. COROLLARY.} {\it%
  Let $X$ be a $d$-dimensional random vector with independent
  and identically distributed components, all having a 
  Student-like distribution with parameter $\alpha$. Let $C$ be 
  a $d\times d$ matrix, with $N(C)>2$. Then
  $$
    \lim_{\ttoi} t^\alpha P\{\, \cXX\geq t\,\} = 0 \, .
  $$
  }

\finetune{}%\bigskip

\stateit{Proof.} Notice that if the largest diagonal element of
$C$ is positive, then $N(C)=1$. Moreover, if this largest diagonal
term vanishes, Lemma 8.3.2 gives $N(C)=2$.  Then, $N(C)>2$ implies
that all the diagonal coefficients of $C$ are negative. Consequently,
$N_0(C)$ --- defined after (8.2.19) --- is at least $2$. Apply Theorem
8.2.21 to conclude. \hfill$\qed$

\bigskip

At this point, the reader who doubts of the usefulness of Theorem
5.1 in providing a systematic technique should try to obtain the
result of this chapter by other methods. Maybe once the results
are known, Theorems 8.2.10, 8.2.21 and 8.3.1 can be proved more
simply. It is also hoped that though Theorem 8.2.1 can be
derived by easier methods, the path taken makes the current
proof rather didactic. As it may have been noticed,
Theorems 8.2.10 and 8.3.1 add extra arguments to the basic ones
developed to prove Theorem 8.2.1.

\bigskip

\centerline{\fsection Notes}
\section{8}{Notes}

\bigskip

\noindent This chapter is motivated by the statistical 
applications in time series developed in chapter 11.

Concerning section 8.1, there is a classical argument for the
Gaussian case. Replacing $C$ by $(C+C^\T )/2$, there is no loss
of generality in assuming that $C$ is symmetric. Thus we can
diagonalize the matrix, writing $C=QDQ^\T$ for a diagonal
matrix $D$ and an orthogonal one $Q$. Define $Y=Q^\T X$. Since
the standard Gaussian distribution is invariant under orthogonal
transformation, $Y$ has again a standard normal distribution.
Thus $\cXX =\sum_{1\leq i\leq d}Y_i^2D_{i,i}$ is a weighted sum of
independent chi-square random variables. This can be generalized
to noncentered Gaussian distributions --- see Imhoff (1961) ---
and opens the door for saddlepoint approximations ---
Barndorff-Nielsen (1990). 

The orthogonal invariance argument can be used for other ad hoc
distributions. However, the Gaussian one is the only
orthogonally invariant distribution with independent marginals.

For quadratic forms with heavy tail distribution, not much seems
to be known. Davis and Resnick (1986) contains some
asymptotic results as $d$ tends to  infinity in a time series
context, based on the point process technique exposed in 
chapters 3--4 of Resnick (1987).

This chapter 8 is full of open questions. To state a few, what
happens in the degenerate cases when $N(C)>2$? Can one find a
closed formula for the integral involved in Theorem 8.2.10? What
is a good numerical scheme to compute such integral? Can we find
more terms and obtain an asymptotic expansion? Can one obtain
good upper bounds instead of asymptotic equivalents? Answers to
these last two questions would be useful in the applications
developed in chapter 11. Can one prove conjecture 8.1.1?

When dealing with heavy tails, there is a fashionable extension
which consists in replacing any power function by itself times a
slowly varying function. This is done mainly for linear
functions of random variables. For quadratic functions, things
turn to be much more complicated, and the classical guess,
consisting of putting the same slowly varying function in the
tail equivalent, is plain wrong.  This can be seen already when
multiplying two heavy tail random variables.

\vfill\eject

~\vskip 28pt
\chapter{9}{Random linear forms}
{\fchapter\hfill 9.\quad Random linear forms}

\vskip 54pt

\noindent In this chapter, we investigate the following problem.
A $d$-dimensional random vector $X=(X_1,\ldots ,X_d)$
defines a random linear form $X(p)=\< X,p\>$ on $\RR^d$. Given
a subset $M$ of $\RR^d$ --- eventually $M$ could be of dimension much
smaller than $d$ --- we can look at the restriction of the linear
form on $M$, and at the distribution of its supremum,\ntn{$X(M)$}
$$
  X(M)=\sup\{\, \< X,p\>  : p\in M\,\} \, .
$$
Writing\ntn{$A_t$}
$$
  A_t=\big\{\, x\in\RR^d  : \sup_\pinM \xp\geq t\,\big\}
  = tA_1 \, ,
$$
we see that
$$
  P\{\, X(M)\geq t\,\} = \int_{tA_1}\d P
$$
where $P$ is the probability measure of $X$. Provided
$0$ is not in the closure of $A_1$, the sets $tA_1=A_t$ are 
moving to infinity with $t$. Theorem 5.1 may provide 
tail approximations of the distribution of $X(M)$.

There has been a tremendous amount of work on this problem, but from a
somewhat different perspective. Traditionally, the set $M$ is
parametrized as $M=\big\{\, p(s) : s\in S\big\}$ for some set
$S$ and some function $p : S\subset \RR^k\to\RR^d$. When $k=1$,
the random variable $X(M)$ is the supremum of a linear
stochastic process; for $k>1$, it is the supremum of a
linear random field. Though it is widely used, parameterizing the
set $M$ has the disadvantage of hiding the fact that $X(M)$ is
completely parameterization free. Any surjective --- even not
smooth, not injective, not measurable, etc. --- change of
parameterization of $M$ leaves $X(M)$ invariant.

Often, one is interested in the supremum norm
$$
  |X|(M) = \sup\big\{\, |\Xp|  : p\in M\,\big\} \, ,
$$
more than in $X(M)$. Up to changing $M$ into $M\cup (-M)$, it is
enough to consider $X(M)$.

The key point to understand is that obtaining an 
approximation of $P\big\{\, X(M)\geq t\,\big\}$ 
for large $t$ covers several distinct
questions. First, notice that $A_1$ is the complement in $\RR^d$
of the convex set\ntn{$C_M$}
$$
  C_M=\big\{\, x\in\RR^d  : x(M)\leq 1\,\big\}
  = \bigcap_{p\in M}\big\{\, x  : \xp \leq 1\,\big\}
$$
--- an intersection of half spaces. So, a first question is to
understand how the probability that $X(M)$ is large is related
to the geometry of $\partial C_M$. A second question is to study
how the geometry of $\partial C_M$ is related to 
that of $M$. When dealing with a parameterization of $M$, 
a third question is how to
read in the parameterization the geometric information we need
about $M$.

Among some amusing features of our point of view, we cannot 
resist mentionning the degenerate case $M=\{\, (1,\ldots
,1)\,\}\subset \RR^d$. Then, one has $X(M)=X_1+\cdots + X_d$.
From our point of view,
approximating the tail distribution of a sum of random
variables, approximating the tail distribution of the supremum of
a linear process, approximating the distribution of the supremum
of a linear random field are all the same problem.

As in the previous chapters, we will investigate light and heavy
tail distributions. Our purpose is again to illustrate the use
of Theorem 5.1 and to show that very different geometric
features govern the tail behavior of $X(M)$ according to the
distribution of $X$.

\bigskip

\noindent{\fsection 9.1. Some results on convex sets.}
%\chapternumber{9.1}
\section{9.1}{Some results on convex sets}

\bigskip

\noindent In this section we relate the set $M\subset\RR^d$ to 
the convex set $C_M$.
Notice that $C_{M\cup M/2}=C_M$. Hence, $C_M$ does not characterize
$M$. To study which part of $M$ is characterized by
$C_M$, we first prove that we can also assume $M$ is closed.

\bigskip

\state{9.1.1. LEMMA.} {\it%
  For any subset $M$ of $\RR^d$, the equality $C_{\closure M}=
  C_M$ holds. Moreover, $C_M$ contains the origin if and only 
  if $M$ is unbounded.}

\bigskip

\stateit{Proof.} The inclusion of $M$ in its closure implies
that of $C_{\closure M}$ in $C_M$. On the other hand, any point
$p$ in the closure of $M$ is the limit of a sequence of points
$p_n$ belonging to $M$. Any point $x$ in $C_M$ satisfies $\<x,p_n\>\leq 1$
and consequently $\xp\leq 1$. This proves the
inclusion of $C_M$ in $C_{\closure M}$.

The second statement in the lemma follows from the equivalence
between the inclusion of the ball of radius $\epsilon$ centered
at the origin in $C_M$ and that of $M$ in the ball of 
radius $1/\epsilon$ around the origin.\hfill$\qed$

\bigskip

In what follows, we assume that $M$ is closed and bounded,
i.e.,
$$
  \hbox{$M$ is compact in $\RR^d$.}
$$
Another way to think of this assumption is that we can take $M$
to be closed since $X(M)=X(\closure M)$. Moreover,
$X(M)$ is infinite almost surely if and only if
$M$ is unbounded and $X$ is nondegenerate. In short, the 
behavior of $X(M)$ is trivial if and only if $M$ is unbounded.
So, we may as well assume $M$ to be compact.

Recall that for a convex set $C$, there is a dense set in
its boundary containing points for which the tangent space to
$\partial C$ is well defined --- see, e.g., Schneider (1993).
Consequently, it makes sense to define\ntn{$M_0$}
$$\displaylines{
  \qquad
    M_0=\closure \big\{\, p\in\RR^d  : \hbox{there exists }
    x\in\partial C_M\, ,\, \xp=1 \, ,
  \cut
    T_x\partial C_M \hbox{ exists, and } 
    p\perp T_x\partial C\,\big\}
    \qquad\cr
  }
$$
The next proposition shows that $M_0$ is the smallest closed set
in $M$ necessary to describe $\partial C_M$.

\bigskip

\state{9.1.2. PROPOSITION.} {\it%
  Assume $M$ is compact.

  \noindent (i) The inclusion $M_0\subset M$ holds.

  \noindent (ii) Moreover, $C_M=C_{M_0}$.

  \noindent (iii) If $M_1$ is closed in $M$ and 
  $C_{M_1}=C_M$, then $M_0\subset M_1$.
  }

\bigskip

\stateit{Proof.} We will use the following claim: If
$x$ is in $\partial C_M$ and $T_x\partial C_M$ exists, then there
exists $p$ in $M$, orthogonal to $T_x\partial C_M$, 
such that $\< x,p\>=1$. Indeed, for such $x$, 
there exists $p$ in $M$ with $\<
x,p\>=1$. Let $u$ be a tangent vector to $C_M$ at $x$. 
We can find a curve
$x (\epsilon )$ in $\partial C_M$ such that $x(\epsilon
)=x+\epsilon u + o(\epsilon )$. Since $x(\epsilon )$ is in 
$C_M$, we
must have $1\geq \< x(\epsilon ),p\>=1+\epsilon \< u,p\> +
o(\epsilon )$ as $\epsilon$ tends to $0$. Thus, $\< u,p\>$
vanishes, and indeed, $p$ is orthogonal to $T_x\partial C_M$.

\noindent (i) Let $q$ be a point in $\RR^d$, such that there 
exists $x$ in $\partial C_M$ with $\< q,x\>=1$, the tangent space
$T_x\partial C_M$ exists, and $q$ is normal to $\partial C_M$
at $x$. By the above claim, there exists $p$ in $M$ and 
orthogonal to $T_x\partial C_M$ such that $\< x,p\>=1$.
Since $M$ is compact, $0$ is in the interior of $C_M$ 
and $T_x\partial C_M$ is of dimension $d-1$. Consequently, 
$p$ and $q$ must be collinear. They are equal 
since $\<x,p\>=\<x,q\>=1$. Thus $q$ belongs to $M$. 
Since $M$ is closed, the inclusion $M_0\subset M$ follows.

\noindent (ii) The inclusion of $M_0$ in $M$ implies that of 
$C_M$ in $C_{M_0}$. To obtain the reverse inclusion, convexity 
of $C_{M_0}$ and $C_M$ shows that we just need to prove $\partial
C_{M_0}\subset \partial C_M$. Assume that 
$\partial C_{M_0}\setminus \partial C_M$ contains a point $x$. 
Then
$$
  \< x,p\> \leq 1 \hbox{ for all } p\in M_0 \, .
  \eqno{(9.1.1)}
$$
Since $C_M$ is included in $C_{M_0}$, such an $x$ cannot belong to
$C_M$. Moreover, since $C_M$ contains the origin and $x$ is
not in $C_M$, there exists $\lambda$ in $(0,1)$
such that $y=\lambda x$ belongs to $\partial C_M$. For any 
positive $\epsilon$, there exists $y_\epsilon$ in 
$\partial C_M$ such that
$|y-y_\epsilon |\leq\epsilon$ and $T_{y_\epsilon}\partial C_M$
exists. Using our claim, there exists $p_\epsilon$ in $M$,
orthogonal to $T_{y_\epsilon}\partial C_M$ and such that
$\< p_\epsilon ,y_\epsilon\> = 1$, i.e., 
$p_\epsilon$ belongs to $M_0$. Since $M_0$ is
closed and is in the compact set $M$, it is compact.
As $\epsilon$ tends to $0$, the points $p_\epsilon$ admit a
cluster point $p$ belonging to $M_0$, thus belonging to $M$. 
Then, as $\lim_{\epsilon\to 0}y_\epsilon = y$,
$$
  1= \< p,y\> = \lambda \< p,x\> \, .
$$
Consequently, $\< p,x\> = 1/\lambda > 1$, which contradicts
(9.1.1). It follows that $\partial C_{M_0}\subset \partial C_M$,
and therefore $C_{M_0}\subset C_M$. 

\noindent (iii) is clear from (i) and (ii), and shows the
minimality of $M_0$. \hfill$\qed$

\bigskip

Proposition 9.1.2 motivates the following definition.

\bigskip

\state{DEFINITION.} {\it%
  We say that a set $M$ is reduced if $M=M_0$.
  }

\bigskip

If $\partial C_M$ has a well defined tangent space at $x$, it admits a
unit normal vector\ntn{$N(\cdot )$} $N (x)$ pointing outward from
$C$. The claim in the proof of Proposition 9.1.2 shows that there
exists $p$ in $M$, orthogonal to $T_x\partial C_M$, and such that
$\<p,x\>=1$.  Such a $p$ is collinear to $N(x)$ and satisfies
$|p|\<x,N(x)\>=1$.  In conclusion,
$$
  M_0=\closure \{\, N(x)/\< x,N(x)\> :
  x\in\partial C_M \, , \, T_x\partial C_M \hbox{ exists}\,\}
  \, . \eqno{(9.1.2)}
$$
This set is called --- traditionally assuming that $\partial C_M$
is smooth --- the polar reciprocal of $\partial C_M$ --- 
see, e.g., Schneider (1993).

Whenever $\partial C_M$ is locally a C$^2$-manifold,
the following lemma shows that $M_0$ is also locally a
C$^2$-manifold. Moreover, 
the second fundamental form of $M_0$ is related to that 
of $\partial C_M$.

\bigskip

\state{9.1.3. LEMMA.} {\it%
  Let $x : U\subset\RR^{d-1}\mapsto\partial C_M$ be a local
  parameterization of $\partial C_M$. 
  Then, $x_0=N\circ x/\< x,N\circ x\>$ defines a local
  parameterization of $M_0$. If ${(b_{i,j})}_{1\leq i,j\leq d-1}$
  is the the matrix of the second fundamental form of $\partial
  C_M$ at $x$ in this parameterization, then $b_{0,i,j}= b_{i,j}
  \big/ \big(\< N(x)\, ,x\> |x|\big)$ is the second
  fundamental form of $M_0$ in the corresponding parameterization
  $x_0$.
  }

\bigskip

\stateit{Proof.} We follow closely Hasani and
Koutroufiotis (1985). Let $u=(u_1,\ldots , u_{d-1})\in U$.
Define $x_i=\partial x/\partial u_i$ and $N_i=\partial N\circ
x/\partial u_i$. The function $x_0=N\circ x/\< x,N\circ x\>$ defines
a local parameterization of $M_0$ according to (9.1.2). 
Define $f=\< x,N\circ x\>$.
Since $\< x_i,N\circ x\>=0$, we have
$$
  x_{0,i}={\partial x_0\over \partial u_i} 
  = {N_i\over f} - {\< x,N_i\>\over f^2} N \, .
$$
Moreover, $x$ is normal to $M_0$ at $x_0$ since
$$
  \< x_{0,i},x\> = {1\over f} \< N_i,x\> - \< x,N_i\> {1\over f}
  = 0 \, .
$$
for all $i=1,\ldots , d-1$. Thus, the components of the second
fundamental form of $M_0$ are
$$\eqalignno{
  b_{0,i,j}
  & =-\Big\< {\partial x_0\over\partial u_i} ,
	    {1\over |x|}{\partial x\over\partial u_j}\Big\>
    =-\Big\< {N_i\over f}-{\< x,N_i\>\over f^2} N ,
	     {x_j\over |x|}\Big\>
    =-{1\over f|x|}\< N_i,x_j\> \cr
  & = {b_{i,j}\over f|x|} \, . &$\qed$\cr
  }$$

\bigskip 

In Lemma 9.1.3, we derived properties of $M_0$ from knowledge on
$\partial C_M$. We will also need to go the other way,
that is obtain some information on the boundary of $C_M$ from the
knowledge of $M$ or $M_0$.

Our next result shows that whenever $C_M$ is bounded,
an half line starting at the origin
can cut $M_0$ in at most one point.

\bigskip

\state{9.1.4. LEMMA.} {\it%
  Assume that $C_M$ is compact in $\RR^d$. If $q$ is in $M_0$, 
  then $q\RR^+\bigcap M_0=\{\, q\,\}$.
  }

\bigskip

\stateit{Proof.} Consider a point $q$ in $M_0$ and $\lambda$
positive such that $\lambda q$ belongs to $M_0$ as well.
Representation
(9.1.2) of $M_0$ implies that $q=\lim_{n\to\infty} N(x_n)/\<
x_n,N(x_n)\>$ and $\lambda q=\lim_{n\to\infty}N(y_n)/\<
y_n,N(y_n)\>$ for points $x_n,y_n$ in $\partial C_M$ at which
$\partial C_M$ has a well defined tangent plane.

Both $\partial C_M$ and $S_{d-1}$ are compact. Up to extracting a
subsequence, we can assume that $x_n$, $y_n$, $N(x_n)$ and
$N(y_n)$ converge respectively to $x$, $y$, $N_1$, $N_2$.
Then $q=N_1/\< x,N_1\>$ and $\lambda q=N_2/\< y,N_2\>$.
Since $q$ and $\lambda q$ are collinear, we have $N_1=N_2$ or
$N_1=-N_2$. Since $C_M$ is convex and contains the origin,
$\< x, N_1\>$ and $\< x,N_2\>$ are nonnegative. The positivity
of $\lambda$ implies $N_1=N_2$.

Next, $C_M$ being convex, it lies on one side of its tangent
spaces. Thus $\< x_n-y_n,N(x_n)\>\leq 0$.  Taking the limit as $n$ tends
to infinity, we obtain $\< x-y,N_1\>\leq 0$.  Permuting $x$ and $y$
yields $\< y-x,N_2\>\leq 0$. Since $N_1=N_2$, the vector $x-y$ is
orthogonal to $N_1$ and $\< x,N_1\>=\< y,N_2\>$.  Thus, $\lambda$
equals $1$.\hfill$\qed$

\bigskip

As a consequence of Lemma 9.1.4, the next result asserts that
whenever $C_M$ is compact and $M_0$ is a manifold, 
the space $p\RR$ is transverse to
$T_pM_0$ for a typical point $p$ in $M_0$.

\bigskip

\state{9.1.5. LEMMA.} {\it%
  If $M$ is compact in $\RR^d$ and $M_0$ is a manifold, then
  $\{\, p\in M_0 : p\RR\subset T_pM_0\,\}$ is nowhere
  dense in $M_0$.
  }

\bigskip

\stateit{Proof.} Assume that the set under consideration is
dense in an open set $U$ of $M_0$. Since $M_0$ is smooth,
$p\RR$ is included in $T_pM_0$ for every $p$ in $U$. Up 
to considering an
open subset of $U$, we can assume that $U$ does not
contain the origin. Thus, $p/|p|$ is a unit vector field in 
the tangent bundle of $M_0$. Let $\gamma$ be an integral 
curve of this field, with $\gamma (0)$ in $U$. It satisfies
the equation $\gamma'=\gamma/|\gamma |$. Consequently, $\gamma =
|\gamma |\gamma'$. Differentiating, we obtain
$\gamma'=\gamma' +|\gamma |\gamma''$. Hence $\gamma''$
vanishes since
$\gamma$ does not. Consequently, $\gamma'$ is constant and
$\gamma (s)=\gamma (0)+ s{\gamma (0)\over |\gamma (0)|}
= \gamma (0)\big( 1+{s\over |\gamma (0)|}\big)$. Hence the curve
$\gamma$ is in the ray $\gamma (0)\RR$, which contradicts 
Lemma 9.1.4.\hfill$\qed$

\bigskip

Making use of Lemmas 9.1.1--9.1.3 and Proposition 9.1.2, we can
reduce $M$ to $M_0$. Moreover, $M_0$ is smooth if $\partial
C_M$ is. In the following, we assume that
$$
  \hbox{$M$ is reduced and is a smooth $m$-dimensional 
  submanifold of $\RR^d$.}
  \eqno{(9.1.3)}
$$

As announced, our next task is to relate the differential geometric properties
of $\partial C_M$ to those of $M$. For this purpose, we build a local
parameterization of $\partial C_M$ starting from one for $M$.
Lemma 9.1.5 asserts that under (9.1.3), if $C_M$ is compact, then,
at a typical point $p$ of $M_0$, the direction $p\RR$ is not contained
in the tangent space $T_pM_0$. Actually, the proof of 
Lemma 9.1.5 shows slightly more; namely, that we cannot have
$p\RR$ included in $T_pM_0$ along a submanifold of $M_0$.
We will assume more, even when $C$ is not compact; namely 
that given $p_0$ in $M$, 
$$
  \hbox{ $p\RR$ is transverse to $T_pM$ for all $p$ in a
  neighborhood of $p_0$.}
  \eqno{(9.1.4)}
$$
Consider a local parameterization
$$
  \underline{p}:(u_1,\ldots ,u_m)\in U\subset \RR^m \mapsto 
  \underline{p}(u_1,\ldots ,u_m)\in M
$$ 
of $M$ around $p_0$. We can assume without
any loss of generality that $U$ is an open neighborhood of the
origin and that $\underline{p}(0)=p_0$. For what follows, it
is convenient to extend $\underline{p}$ to a map defined on a
neighborhood of the origin of $\RR^{d-1}$. Thus, let $V$ be a
neighborhood of $0$ in $\RR^{d-1-m}$ and consider\ntn{$p(\cdot )$}
$$
  p : u\in U\times V\subset \RR^m\times \RR^{d-1-m}
  \mapsto p(u)=\underline{p}(u_1,\ldots ,u_m) \, .
$$
Under (9.1.4), to each point $p$ near $p_0$, we can associate a unit
normal vector\ntn{$\nu (p)$}
$$
  \nu(p)\in (T_pM+p\RR )\ominus T_pM \, ;
$$
that is, $\nu(p)$ is normal to $T_pM$ in $T_pM+p\RR$.

Let\ntn{$X_i(p)$} $X_{m+1}(p),\ldots ,X_{d-1}(p)$ be an orthonormal moving
frame in $(T_pM+p\RR)^\perp = \RR^d\ominus (T_pM+p\RR )\equiv
\RR^{d-m-1}$. For $u$ in $U\times V$, let\ntn{$X(u)$}
$$
  X(u) = {\nu\circ p(u)\over \big\< \nu\circ p(u),p(u)\big\>}
  + \sum_{m+1\leq j\leq d-1}\big( \phi_j\circ p(u)+u_j\big)
  X_j\circ p(u)
$$
where the function $\phi_j$'s\ntn{$\phi_j (p)$} are C$^2$ and such that
$$
  (u_{m+1},\ldots ,u_{d-1})\in V\mapsto 
  I\big( X(u_1,\ldots ,u_m, u_{m+1}, \ldots , u_{d-1})\big)
$$
is minimum at $0$ --- the existence
of these functions and their smoothness comes from the implicit
function theorem and the smoothness of $I$. By construction,
$X(u)$ is normal to $T_{p(u)}M$ for all $u$ in $U\times V$.

Our next lemma asserts that if $M$ is positively curved, then
$X$ defines a local parameterization of $\partial C_M$. Combined
with Lemma 9.1.3, it allows us to parameterize $M$ or $\partial
C_M$, whichever is the most convenient.

\bigskip

\state{9.1.6. LEMMA.} {\it%
  Under (9.1.3)--(9.1.4), if the second fundamental form 
  of $M$ relative
  to the normal field $X$ is definite positive at every point,
  then $X(u)$ defines a local parameterization of $\partial C_M$.
  Moreover, $p(u)$ is outward normal at $\partial C_M$ at
  $X(u)$.}

\bigskip

\stateit{Proof.} Denote by $\Pi_M^X(p)$ the second fundamental form
of $M$ at $p$ along the normal field $X$.
Notice first that $X(u)$ is orthogonal to $T_{p(u)}M$. Also,
$$
  \< p(u),X(u)\> = 1 \, , 
$$
and $p$ is orthogonal to 
$\span \{\, X_{m+1},\ldots , X_{d-1}\,\}$.
Consequently, if $q(s)$ is a curve on $M$, parametrized by arc
length, such that $q(0)=p(u)$ and $q'(0)=t$, then,
as $s$ tends to $0$,
$$\eqalign{
  \< q(s),X(u)\> 
  & = 1 + {s^2\over 2} \< q''(0),X(u)\> + o(s^2) \cr
  & = 1 - {s^2\over 2} \Pi_M^X\big( p(u)\big)(t,t) + o(s^2) 
    \, .\cr
  }$$
Thus, $s\mapsto \< q(s),X(u)\>$ is maximal and equals $1$
at $s=0$. This proves that
$X(u)$ is in $\partial C_M$ and that $X(u)$ indeed defines 
a local parameterization of
$\partial C_M$.

To prove that $p(u)$ is normal to $\partial C_M$, notice first
that
$$
  {\partial \over\partial u_i}X(u)
  = X_i(u) \qquad \hbox{ for } i=m+1,\ldots , d-1 \, .
$$
Furthermore, writing $\nu_i={\dispar\over \vbox to 8pt{}\dispar u_i} 
\nu\circ p(u)$
and $p_i={\dispar\over\vbox to 8pt{}\dispar u_i} p(u)$, for $i=1,\ldots ,m$,
$$\eqalign{
    X_i={\partial \over \partial u_i}X 
    = {\nu_i\over \< \nu,p\>} 
    & - {\< p,\nu_i\>\over \< p,\nu\>^2} \nu
      + \sum_{m+1\leq j\leq d-1}\d\phi_j(p)\cdot p_i \, X_j \cr
    & +\sum_{m+1\leq j\leq d-1} (\phi_j\circ p+u_j)
      \d X_j(p)\cdot p_i \, .
  \cr}
  \eqno{(9.1.5)}
$$
From the very definition of $X_j$, for $j=m+1,\ldots , d-1$, we
infer
$$
  \< X_j\circ p(u),p(u)\>=0 \qquad \hbox{ for } j=m+1,\ldots ,
  d-1 \,  .
  \eqno{(9.1.6)}
$$
Differentiating this equality yields
$$
  \< \d X_j(p)\cdot p_i , p\> + \< X_j\circ p,p_i\> = 0 \, ,
  \quad i=1,\ldots m \, ,\quad j=m+1,\ldots , d-1 \, .
$$
Therefore, since $X_j$ is orthogonal to $T_pM$, we have
$$
  \< \d X_j(p)\cdot p_i,p\> = 0 \qquad \hbox{ for }
  j=m+1,\ldots , d-1 \, .
  \eqno{(9.1.7)}
$$
Combining (9.1.5)--(9.1.7) yields $\< X_i,p\>=0$ for $i=1,\ldots
, m$. Thus, $p(u)$ is normal to $\partial C_M$ at $X(u)$. 

Since $C_M$ contains the origin and $\< p,x\>=1$ is positive, 
the vector $p(u)$ must be
pointing outward $\partial C_M$.\hfill$\qed$

\bigskip

We can now explain how to compute the second fundamental form of
$\partial C_M$ at $X(u)$. It is determined by its value on the
basis $X_i$ of the tangent space. Let us denote by $N = p/|p|$
the outward unit normal vector field to $\partial C_M$. The
proof of Lemma 9.1.6 shows that
$$\displaylines{
    \< \d N\cdot X_i,X_j\>
     = {1\over |p|}\< p_i,X_j\>
  \hfill\cr\noalign{\vskip 0.05in}\hfill
    = \left\{\,
    \vcenter{\hbox{$0$ \quad if $i=m+1,\ldots ,d-1$ or $j=m+1,\ldots ,d-1$ \hfill}
             \break
             \hbox{ ${ \displaystyle\vtop to 5pt{}\< \nu_j,p_i\>\over
		       \vbox to 9pt{}\displaystyle|p|\< \nu,p\> }
                     + { \displaystyle\vtop to 5pt{}1
                         \over \displaystyle\vbox to 9pt{}|p| }
                     \sum_{m+1\leq r\leq d-1}(\phi_r\circ p+u_r)
	             \< p_i,{ \displaystyle\partial \over 
                              \displaystyle\vbox to 9pt{}\partial u_j }X_r\>$
                      otherwise} 
            }\right.
  \cr}
$$
Consequently, if $\Pi_M^Y$ denotes the second fundamental form of
$M$ relative to a unit normal vector field $Y$, we have
$$\displaylines{
  -\< \d N\cdot X_i,X_j\> \hfill (9.1.8)%
  \cr\noalign{\vskip .05in}\hfill
  {}= \cases{ {\displaystyle\vtop to 5pt{}\Pi_M^\nu(p_i,p_j)\over
	       \vbox to 9pt{}
               \displaystyle |p|\< \nu ,p\>}
		+ {\displaystyle\vtop to 3pt{}1
                   \over\displaystyle\vbox to 9pt{}|p|}
                \sum_{m+1\leq r\leq d-1}
		(\phi_r\circ p+u_r) \Pi_M^{X_r}(p_i,p_j) 
		& \cr
                \hbox{\kern 2.5in if $i,j=1,\ldots ,m$.} &\cr
             \noalign{\vskip .05in}
	     0 \hbox{\qquad otherwise.} & \cr}
  \cr}
$$
This second fundamental form vanishes whenever $i$ or
$j$ is not in $\{\, 1,\ldots , m\,\}$, expressing the 
fact that $C_M$ is a ruled
surface with nontrivial generators if $m<d-1$. Along a
generator of dimension $d-1-m$, the set $C_M$ is flat and has
vanishing curvature.

\bigskip

Now, consider a convex function $I$ on $\RR^d$. Since we are
able to relate points and geometry of $\partial C_M$ to points
and geometry of $M$, we should be able to relate points in
$\partial C_M$ which minimize $I$ to some specific points in
$M$. For this purpose, to a function $f$ defined on
$\RR^d$ we associate the functions\ntn{$f_\sbullet$}\ntn{$f^\sbullet$}
\ntn{$I_\sbullet$}
$$\eqalign{
  f^\sbullet (x)&=\sup\big\{\, f(y) : \< x,y\> = 1\,\big\} \, ,\cr
  f_\sbullet (x)&=\inf\big\{\, f(y) : \< x,y\> = 1\,\big\} \, .\cr
  }$$
The basic properties of these transforms will be of some use and
are stated in the following proposition. Notice that the
statement (iv) in the Proposition does not require any smoothness.

\bigskip

\state{9.1.7. PROPOSITION.} {\it%
  Let $I$ be a convex function on $\RR^d$,
  such that $\lim_{|x|\to\infty}I(x)=\infty$. Then,

  \noindent (i) $I_\sbullet$ is continuous on $\RR^d\setminus
  \{\, 0\,\}$ with $\lim_{|x|\to 0}I_\sbullet (x)=\infty$.
  Moreover, if $I$ is minimal at $0$, then 
  $\lim_{|x|\to\infty}I_\sbullet (x)\allowbreak =I(0)$;

  \noindent (ii) if $I(C_M^{\rm c})=I(\partial C_M)$, then
  $$
    I(C_M^{\rm c})
    =\inf\big\{\, I_\sbullet (p) : p\in M\,\big\} 
    = I_\sbullet (M)= I_\sbullet (M_0) \, ;
  $$ 
  (even if $M$ is not reduced)

  \noindent (iii) points in $C_M^{\rm c}$ minimizing $I$ 
  correspond naturally to points in $M_0$ minimizing $I$ in
  the sense that
  $$\displaylines{\qquad
      \bigcup_{ \{ x\in C_M^{\rm c} \, :\, I(x)=I(C_M^{\rm c})\} }
      \{\, m\in M_0 : \<m,x\>=1\,\}
    \cut
      = \{\, m\in M_0 : I_\sbullet (m)=I_\sbullet (M_0)\,\} \, ;
    \qquad\cr}
  $$
  \noindent (iv) for any convex function $f$ with its minimum at $0$, 
  the equality $(f_\sbullet)^\sbullet = f$ holds.
  }

\bigskip

\stateit{Proof.} (i) If $\xp
=1$, then $1\leq |x| |p|$. Consequently,
$$
  I_\sbullet (x)\geq \inf\big\{\, I(p) : |p|\geq 1/|x|\,\big\} \, ,
$$
and 
$$
  \lim_{|x|\to 0}I_\sbullet (x)=\lim_{|p|\to\infty}I(p)=\infty
  \, .
$$
Moreover, if $\<x,p\>=1$, then 
$p={\displaystyle x\over\displaystyle |x|^2}+\Proj_{x^\perp}p$. 
Therefore,
$$\eqalign{
  I(0)=\inf\big\{\, I(p) : |p|\in \RR^d\,\big\} 
  \leq I_\sbullet (x)
  & \leq \inf\Big\{\, I\Big( {x\over |x|^2}+q \Big)
     : q\perp x\,\Big\} \cr
  & \leq I\Big( {x\over|x|^2}\Big) \, ,\cr
  }$$
The equality
$\lim_{|x|\to\infty}I_\sbullet (x)=I(0)$ follows.

Let us now prove that $I_\sbullet$ is continuous. Let
$p$ be a nonzero vector in $\RR^d$. Let $\epsilon$ be positive and less
than $|p|$. Consider a point $q$ at distance $\epsilon$ from $p$.
Since $I$ is continuous and blows up at infinity,
there exists $x$ in $\RR^d$ such that $\xp=1$ and $I(x)=I_\sbullet
(p)$. Then,
$$
  |\< x,q\> -1|
  = |\< x,q-p\>|
  \leq |x|\epsilon \, .
$$
Hence, there exists $y$ in $\RR^d$ such that $\< y,q\>=1$ and
$|x-y|\leq \epsilon |x|/|q|$ --- take $y=x-\big( \< x,q\> -1\big)
q/|q|^2$. Therefore,
$$
  I_\sbullet (q)\leq I(y)\leq 
  \sup\Big\{\, I(v)  : |v-x|
               \leq {|x|\over |q|}\epsilon\,\Big\}
  \, .
$$
Since $I$ is continuous and $p$ is nonzero, it follows that
$$
  \limsup_{q\to p}I_\sbullet (q) \leq I(x)=I_\sbullet (p) \, .
$$
Next, consider a sequence $p_k$ in the ball $B(p,\epsilon )$ of radius
$\epsilon$ centered at $p$, converging to $p$, and such that
$\lim_{k\to\infty}I_\sbullet (p_k)=\liminf_{q\to p}I_\sbullet
(q)$. Let $x_k$ be such that $\< p_k,x_k\>=1$ 
and $I_\sbullet (p_k)=I(x_k)$.  The function $I$ is bounded 
on the set $q/|q|^2$ for $q$ in $B(p,\epsilon )$ provided $\epsilon$
is strictly less than $|p|$. Then, the 
sequence $I_\sbullet (p_k)\leq I\big(p_k/|p_k|^2\big)$ is bounded, 
and so is the sequence $I(x_k)$.  Therefore, $x_k$ is in a compact set
$I^{-1}([0,c])$ for some positive $c$, and admits a clustering point
$x$. After taking a subsequence, we can assume that $x_k$ converges to
$x$ as $k$ tends to infinity. Since $\< x_k,p_k\> =1$, 
we have $\< x,p\> =1$ and thus $I_\sbullet (p)\leq I(x)$. Moreover, 
continuity of $I$ implies
$$
  I_\sbullet (p)\leq I(x) =\lim_{k\to\infty}I(x_k)
  = \lim_{k\to\infty} I_\sbullet (p_k)
  = \liminf_{q\to p}I_\sbullet (q) \, .
$$
Overall, $\lim_{q\to p}I_\sbullet (q)=I_\sbullet(p)$ 
and $I_\sbullet$ is continuous on $\RR^d\setminus \{\, 0\,\}$.

\noindent (ii) Let $\epsilon$ be positive and $p$ in $M$ 
such that $I_\sbullet (p)\leq I_\sbullet (M)+\epsilon$. 
There exists $x$ in $\RR^d$ such that $\xp=1$
and $I_\sbullet (p)=I(x)$. Since $\xp =1$,
the point $x$ is not in the interior of $C_M$.
Since $I$ is continuous,
$$
  I(C_M^{\rm c})=I\big( (\interior C_M)^{\rm c}\big)
  \leq I(x) = I_\sbullet (p)
  \leq I_\sbullet (M)+\epsilon \, .
$$
Since $\epsilon$ is arbitrary, we established the inequality 
$I(C_M^{\rm c})\leq I_\sbullet (M)$.

Next, we use the assumption $I(C_M^{\rm c})=I(\partial C_M)$.
Continuity of $I$ ensures that there exists
$x$ in $\partial C_M$ such that $I(x)=I(\partial C_M)$. Hence, 
there exists $p$ in  the closure of $M$ with $\px =1$ 
--- otherwise, $x$ would not be on
the boundary of $C_M$. Therefore, $I_\sbullet (p)\leq I(x)$, and
since $I_\sbullet$ is continuous,
$$
  I_\sbullet (\closure M) = 
  I_\sbullet (M)\leq I_\sbullet (p)\leq I(x) 
  = I(C_M^{\rm c}) \, .
$$
This proves $I(C_M^{\rm c})=I_\sbullet (M)$. Since $C_M=C_{M_0}$, we
also obtain $I_\sbullet (M)=I_\sbullet (M_0)$. This proves assertion
(ii). Notice also that we proved $I_\sbullet (p)=I(x)$.

\noindent (iii) Let $x$ be in $\partial C_M$ minimizing $I$
over $C_M^{\rm c}$. Since $M_0$ is closed, $\<x,m\>=1$ for some $m$ in
$M_0$. For this $m$, if $I_\sbullet (m)<I(x)$, then there exists $y$
such that $\< m,y\>=1$ and $I(y)< I(x)$. Given how we choose $x$, the
point $y$ is not in $C_M^{\rm c}$.  Since $\< m,y\>=1$ and $M_0$ is
reduced, $y$ is in $\partial C_M$, and $x$ does not minimize $I$ over
$C_M^{\rm c}$. Consequently, $I_\sbullet (m)=I(x)$. Then, assertion
(ii) implies that $I_\sbullet (m)=I_\sbullet (M_0)$. This proves that
in statement (iii), the set in the left hand side is included in that
in the right hand side.

To prove the reverse inclusion, let $m$ in $M_0$ minimizing
$I_\sbullet$. Since $\lim_{|x|\to\infty}I(x)=\infty$ and $I$ is
continuous, $I_\sbullet (m)=I(x)$ for some $x$ such 
that $\< m,x\>=1$. Then, $I(x)=I_\sbullet (M_0)$ and $x$ is not in the 
interior of $C_M$. Then, assertion (ii) implies 
$I(x)=I(C_M^{\rm c})$. Consequently, the set in the right hand side 
of statement (iii) is included in that in the left hand side.

\noindent (iv) If $\<x,p\>=1$, then $f_\sbullet (p)\leq f(x)$.
Thus, $(f_\sbullet)^\sbullet (x)\leq f(x)$. Now seeking a
contradiction, let $x$ be such that $f(x)$ is positive and  assume that
$(f_\sbullet)^\sbullet (x)<f(x)$. Let $H$ be a supporting
hyperplane at $x$ of the level set $\{\, y: f(y)\leq f(x)\, \}$.
We can write $H=\{\, e\,\}^\perp$ for some $e$ in $\RR^d$. If
$e$ is orthogonal to $x$, then $x$ is in $H$. Since $f$ is 
convex, $f(x)\leq f(x+h)$ for any $h$ in $H$. Taking $h=-x$ leads
$f(0)\geq f(x)>0$, a contradiction. Thus, $e$ is not
orthogonal to $x$.

Notice that
$$
  \{\, p :\<p,x\>=1\,\}
  = \Big\{\, {x\over |x|^2}+y : y\perp x\,\Big\} \, .
$$
Moreover, if $y$ is orthogonal to $x$,
$$\displaylines{\qquad
  \Big\{\, q : \Big\< {x\over |x|^2}+y,q\Big\>=1\,\Big\}
  \hfill\cr\noalign{\vskip .05in}\hfill
  = \Big\{\, \alpha x +(1-\alpha ){y\over |y|^2} +z :
  \alpha\in \RR \, , \, z\perp \span (x,y)\,\Big\} \, .
  \qquad\cr}
$$
Going back to the definition of $f_\sbullet$ and
$(f_\sbullet)^\sbullet$, the assumption $(f_\sbullet)^\sbullet
(x)<f(x)$ can be rewritten: There exists $\epsilon$ positive such that
whenever $y$ is orthogonal to $x$, the inequality 
$f\big( \alpha x+(1-\alpha )y|y|^{-2}+z\big) \leq f(x)-\epsilon$ 
holds for some real $\alpha$ and $z$ orthogonal to $\span (x,y)$.

Consider the following $y$. If $e$ is collinear to $x$, let $y$
be any nonzero vector in $\{\, e\,\}^\perp$. Otherwise, let
$y=a\Proj_{x^\perp}e+bx$, where $a$, $b$ are such that
$$\eqalign{
  0 & = \< y-|y|^2x,e\> \cr
    & = -a^2\,|\Proj_{x^\perp}e|^2 \<x,e\>
      +a\,\<\Proj_{x^\perp}e,e\>
      + (b^2|x|^2+b)\<x,e\> \, . \cr
  }
$$
This quadratic equation in $a$ always has a solution for $|b|$
large enough as well as $|b|$ small enough, since $\< x,e\>$ is nonzero. 
Moreover, for $|b|$ small enough but nonzero, the solution is not $0$, 
and $y$ is not null.

This choice of $y$ ensures that $e$ is in the space spanned by
$x$ and $y$. Consequently,
the orthocomplement of $\span (x,y)$ is in $H$. Moreover, 
for any real $\alpha$,
$$
  \alpha x+(1-\alpha ){y\over |y|^2} 
  = x +{1-\alpha\over |y|^2} (y-|y|^2x)\in x+H \, .
$$
Consequently, for any real $\alpha$ and any $z$ orthogonal
to $\span (x,y)$, the point $\alpha x+ (1-\alpha )y|y|^{-2} + z$
belongs to $x+H$. Since $f$ is convex,
$$
  f\Big( \alpha x+(1-\alpha ){y\over |y|^2}+z\Big)
  \geq f(x)
$$
contradicting $(f_\sbullet )^\sbullet(x)<f(x)$.\hfill$\qed$

\bigskip

Given Lemma 9.1.1 and Proposition 9.1.2, we can replace $M_0$ by
$M$ in statement (iii) of Proposition 9.1.5.

The explicit calculation of $I_\sbullet$ depends of course on
$I$. Notice however that homogeneity is preserved, as indicated
in the following result.

\bigskip

\state{9.1.8. LEMMA.} {\it%
  If $I$ is positively $\alpha$-homogeneous, 
  then $I_\sbullet$ is $-\alpha$-homogeneous.}

\bigskip

\stateit{Proof.} The result is straightforward since
$$\eqalignno{
    I_\sbullet (tx)
    = \inf\big\{\, I(y)  : \< tx,y\>=1\,\big\}
  & = \inf\big\{\, I(y/t) : \< x,y\> =1\,\big\} \cr
  & = t^{-\alpha}I_\sbullet (x) \,  .
    &\qed\cr}
$$

\bigskip

Finally, we calculate $I_\sbullet$ in an important case for
applications. For $r$ positive , let $|x|_r=(\sum_{1\leq i\leq
d}|x_i|^r)^{1/r}$. If $r$ is larger than $1$ this is the $\ell_r$-norm 
of $x$.

\bigskip

\state{9.1.9. LEMMA.} {\it%
  If $I(x)=c\big( |x_1|^\alpha+\cdots + |x_d|^\alpha \big)$,
  then $I_\sbullet (p)=c/|p|_\beta^\alpha$ where $\alpha$ and 
  $\beta$ are conjugate --- i.e., $\alpha^{-1}+\beta^{-1}=1$.}

\bigskip

\stateit{Proof.} Without any loss of generality, we
can assume that $c=1$. If $\px =1$, H\"older's inequality yields
$1\leq |x|_\alpha |p|_\beta$, and so $I(x)=|x|_\alpha^\alpha
\geq 1/|p|_\beta^\alpha$. 

On the other hand, for $x_i=\sign (p_i)|p_i|^{1\over \alpha -1}/
|p|_\beta^\beta$, we have $\xp =1$ and $I(x)=1/|p|_\beta^\alpha$.
\hfill$\qed$

\bigskip

In the situation described in Lemma 9.1.9, one sees
that $I_\sbullet (M)$ is
related to $\sup\big\{\, |p|_\beta  : p\in M\,\big\}$, that is
to the radius of the smallest ball in $\ell_\beta$ which
contains $M$.

\bigskip

\noindent{\fsection 9.2. Example with a light tail.}
%\chapternumber={9.2}
\section{9.2}{Example with a light tail}

\bigskip

\noindent In this section, we consider a random 
vector $X=(X_1,\ldots ,X_d)$ in $\RR^d$, having a 
log-concave density $\exp (-I)$. Our first
result is elementary. It is inspired by its Gaussian analogue,
where $I(x)={1\over 2}|x|^2+\log (2\pi )^{d/2}$.
It shows that under a growth control on $I$, we can estimate the
exponential decay of $P\{\, X(M)\geq t\,\}$ as $t$ tends to
infinity.

\bigskip

\state{9.2.1. PROPOSITION.} {\it%
  If $M$ is compact and
  $$
    \lim_{|v|\to 0}\limsup_{|x|\to\infty} {I(x+v)\over I(x)}
    \leq 1 \, ,
    \eqno{(9.2.1)}
  $$
  then
  $$
    \lim_{\ttoi}{1\over I_\sbullet (M/t)}
    \log P\big\{\, X(M)\geq t\,\big\} = -1 \, .
  $$}

\bigskip

\stateit{Proof.} Recall that the events $\{\, X(M)\geq
t\,\big\}$ and  $\{\, X\not\in tC_M\,\}$ are equal. 
Let $a_t=I(tC_M^{\rm c})$. Proposition 2.2 yields
$$
  \lim_{\ttoi}{1\over a_t}\log \big\{\, X(M)\geq t\,\big\} = -1
$$
if and only if
$$
  \lim_{\epsilon\to 0}\liminf_{\ttoi} {1\over a_t}
  \log |tC_M^{\rm c}\cap\Gamma_{(1+\epsilon)a_t}| = 0 \, .
  \eqno{(9.2.2)}
$$
Let $\epsilon$ be positive. Since $I$ is convex and $C_M$ is a
neighborhood of the origin  the points minimizing $I$ 
over $\RR^d$ are included in $tC_M$ for $t$ large enough. 
Consequently, for $t$ large enough, 
there exists $x_t$ in $\partial (tC_M)$ such
that $I(tC_M^{\rm c})= I(x_t)$. 

Let $\eta$ be a positive number. Consider a set of orthogonal
vectors $v_1,\ldots , v_d$ in a supporting hyperplane of
$\partial C_M$ at $x_t$, such that $\eta /2\leq |v_i|\leq \eta$
for all $i=1,\ldots , d-1$. The supporting hyperplane is the
orthocomplement of a unit vector $e$, pointing outward from
$\partial C_M$. We define $v_d=\eta e$. Since $M$ is compact,
the interior of $C_M$ contains the origin. Therefore, $x_t$
tends to infinity with $t$. Consequently, if $\eta$ is small
enough, (9.2.1) implies that for $t$ large enough
$$
  I(x_t+v_i)\leq (1+\epsilon )I(x_t) = (1+\epsilon ) a_t \, .
$$
Since $\Gamma_{(1+\epsilon )a_t}$ is convex, the simplex with
vertices $x_t$, $x_t+v_i$, $1\leq i\leq d$, is in $tC_M^{\rm
c}\bigcap \Gamma_{(1+\epsilon )a_t}$. Its volume does not
depend on $t$; it bounds the volume of 
$tC_M^{\rm c}\bigcap\Gamma_{(1+\epsilon )a_t}$ below. 
Consequently, (9.2.2) holds.

It remains to prove that $I(tC_M^{\rm c})=I_\sbullet (M/t)$. This is
clear since $tC_M=C_{M/t}$, and Proposition 9.1.7 holds.
\hfill$\qed$

\bigskip

Condition (9.2.1) looks good, but is far from being the best
that we can obtain. In particular, it does not cover the function
$I(x)=\exp (|x|^\alpha )$ with $\alpha >1$. The following will
do, but assumes that $I$ is differentiable.

\bigskip

\state{9.2.2. PROPOSITION.} {\it%
  In Proposition 9.2.1, one can replace assumption (9.2.1) by
  $$
    \lim_{\delta \downarrow 0}\limsup_{x\to\infty}
    { e^{-\delta I(x)}\sup_{|u|\leq 1} |\D I(x+ue^{-\delta I(x)})|
     \over I(x) } = 0\, .
  $$
  }

\bigskip

Notice that since $I$ blows up at infinity, $e^{-\delta I(x)}$
is very tiny for large $x$. In essence, when $d=1$, the new
condition asserts that $I'e^{-\delta I}/I$ tends to $0$. 
When $I$ is large, 
$$
  I'e^{-\delta I}/I\leq I'e^{-\delta I/2}
  = 2(e^{-\delta I/2})'/\delta \, .
$$
Therefore, the only possible limit for $I'e^{-\delta I}/I$ as
its argument tends to infinity is $0$ --- but in this discussion,
nothing guarantees that the limit exists. This does not show that
the new condition holds for any convex function; but it suggests
that those which do not satisfy this condition are rather 
pathological.

\bigskip

\stateit{Proof of Proposition 9.2.2.} We follow the proof of
Proposition 9.2.1. All what we need to do is to specify how to
pick the vectors $v_i$, $1\leq i\leq d$. The new assumption
implies that there exists a function $\delta (x)$ tending to
$0$ at infinity, such that
$$
  \lim_{x\to\infty} 
  { e^{-\delta (x)I(x)}\sup_{|x-y|\leq \exp(-\delta (x)I(x))}
    |\D I (y)| \over I(x) } 
  = 0 \, .
$$
Let $\eta=\exp \big(-\delta (x_t)I(x_t)\big)$, and let
$\epsilon$ be an arbitrary positive number. The previous limit
shows that
$$
  \eta \sup_{|y-x_t|\leq \eta}|\D I(y)| \leq \epsilon I(x_t)
$$
for $t$ large enough. With this new $\eta$, take the $v_i$'s 
exactly as in the proof of Proposition 9.2.1. Since $|v_i|$ is
at most $\eta$, we have, for $i=1,\ldots , d$,
$$
  I(x_t+v_i) 
  \leq I(x_t) + |v_i|\sup_{y : |x-y|\leq |v_i|} |\D I(y)|
  \leq I(x_t)(1+\epsilon )\, .
$$
Consequently, the simplex with vertices $x_t$, $X_t+v_i$, $1\leq
i\leq d$ lies in $tC_M\bigcap \Gamma_{(1+\epsilon )I(x_t)}$.
Its volume is of order $\eta^d$. Thus, to check (9.2.2), it
suffices to show that
$$
  \lim_{\ttoi} {\log \eta\over I(x_t)} 
  = \lim_{\ttoi} -\delta (x_t) 
  = 0 \, .
$$
This is plain from the definition of $\delta (\cdot )$.\hfill$\qed$

\bigskip

Clearly, one can do many variations on the theme, and get
different conditions for the conclusion of Proposition 9.2.1 to
hold.

To obtain a sharper result than in Propositions 9.2.1 or
9.2.2, that
is to estimate $P\big\{\, X(M)\geq t\,\big\}$ and not its logarithm,
we need further assumptions. Many results could be obtained
under various hypotheses. We will suppose that the
random vector $X$ has density $ae^{-I}$, where
$$
  \hbox{$I$ is convex, $\alpha$-positively homogeneous for some
  $\alpha >1$.}
  \eqno{(9.2.3)}
$$
Under (9.2.3), Theorem 7.1 settles more or less the question
of approximating
$$
  P(tC_M^{\rm c})= a\int_{tC_M^{\rm c}}e^{-I(x)} \d x \, .
$$
Of course, we need to verify the assumptions of Theorem 7.1.
Since $C_M$ is convex, assumption (7.5) is always satisfied,
while (7.1) is guaranteed by the boundedness of $M$. Thus, only
(7.3) and (7.4) are left to check.

It does not seem that we can work out a general theory much further.
But let us show how we can obtain a tail equivalent for
$P\{\, X(M)\geq t\,\}$ from Theorem 7.1.

Proposition 9.2.1 and Lemma 9.1.8 imply
$$
  \lim_{\ttoi} {1\over t^\alpha}\log P\big\{\, X(M)\geq t\,\big\}
  = -I_\sbullet (M)
$$
and so we are done as far as the exponential term is concerned.
Reading the formula in Theorem 7.1, we have
$I(A_1)=I_\sbullet (M)$.

We then need to relate the differential geometric quantities
involved in Theorem 7.1 to those of $M$. From section 9.1, we
infer that we can reduce $M$ to $M_0$. So, let us assume that
$$
  \hbox{$M_0$ is a closed, connected, $m$-dimensional 
  submanifold of $\RR^d$.}
$$
Since we can replace $M$ by $M_0$, we will drop the subscript
and write $M$ instead of $M_0$ until the end of this section.
Consider a local parameterization $p(u)$ of $M$. We have seen
in section 9.1 that it induces a parameterization $X(u)$ on
$\partial C_M$. Assume that 
$$\displaylines{\hfill
    \calD_{\sbullet , M} 
    =\big\{\, p\in M  : I_\sbullet (p) = I_\sbullet (M)\,\big\}
  \cut
     \hbox{is a $k$-dimensional submanifold of $M$.}
  \hfill\cr}
$$
We can choose the parameterization $p(\cdot )$ such that 
$$
  (u_1,\ldots ,u_k)\in U'\subset \RR^k\mapsto 
  p(u_1,\ldots ,u_k, 0, \ldots , 0)\in\calD_{\sbullet ,M}
$$
is a local parameterization of $\calD_{\sbullet, M}$. With the
notation of section 9.1, it follows
from Lemma 9.1.6 and Proposition 9.1.7 that
$$
  (u_1,\ldots , u_k)\in U'\mapsto X(u_1,\ldots , u_k,0, 
  \ldots , 0)\in\calD_{C_M^{\rm c}}
$$
is a parameterization of $\calD_{C_M^{\rm c}}$. To compute the
Riemannian measure $\calM_{\calD_{C_M^{\rm c}}}$, we first 
compute the first fundamental form of the surface 
$\calD_{C_M^{\rm c}}$. We obtain
$$
  \d\calM_{\calD_{C_M^{\rm c}}}\big( X(u)\big)
  = \Big( \det {\big( \< X_i(u)\, ,X_j(u)\>\big)}_{1\leq i,j
  \leq k}\Big)^{1/2} \d u_1\wedge \cdots \wedge \d u_k \, ,
$$
where $X_i(u)$ is given in (9.1.5). The first fundamental form
of $\partial C_M$ involves inner products of $\nu$, $\nu_i$, $X_i$,
$\d X_j\cdot p_i$, which can be expressed in term of the third
fundamental form of $M$ --- or analogously, in term of the
connection form on its normal frame principal bundle. In
general, such an expression is rather involved; we will see how
it simplifies in some cases.

To compute the curvature term $G_{\smash{C_M^{\rm c}}}$ 
in Theorem 7.1, the comment
after the statement of Theorem 7.1 shows that it equals
$$
  G(x)=\Pi_{\Lambda_{I(C_M^{\rm c})},x}^\pi-
       \Pi_{\partial C_M,x}^\pi
$$
for all $x$ in $\calD_{C_M^{\rm c}}$. Since
\smash{$\overline{G}(x)=\Pi_{\Lambda_{I(C_M^{\rm c})},x}-\Pi_{\partial
C_M,x}$} vanishes on directions tangent to \smash{$\calD_{C_M^{\rm c}}$}, 
one way to compute $\det\, G(x)$ is actually to diagonalize
$\overline{G}(x)$, and take the product of its positive
eigenvalues --- by construction, the eigenvalues of
$\overline{G}(x)$ are nonnegative; they are null only on the
eigensubspace 
\smash{$T_x\partial C_M^{\rm c}\ominus T_x\calD_{C_M^{\rm c}}$}.
Ultimately, we need to calculate $\det\big( G(x)-\lambda
\Id\big)$. For this purpose, consider an orthonormal basis
$e_1,\ldots , e_{d-1}$ of 
$T_x\partial (C_M^{\rm c})=T_x\Lambda_{I(C_M^{\rm c})}$. 

Denote by $\calX$ the matrix obtained in writing the 
vectors $X_1,\ldots,\allowbreak X_{d-1}$ in the basis $e_1,\ldots , e_{d-1}$, 
that is $\calX ={(\< e_i,X_j\>)}_{1\leq i,j\leq n}$. This 
matrix is nonsingular and

$$\eqalign{
  \det (\overline{G}(x)-\lambda \Id ) = 0
  & \Longleftrightarrow \det (\calX^\T\overline{G}(x)\calX 
    -\lambda\calX^\T\calX )= 0 \cr
  & \Longleftrightarrow \det \big( (\calX^\T\calX)^{-1}\calX^\T
    \overline{G}(x)\calX-\lambda\Id )= 0 \, . \cr
  }$$
How do we compute $\calX^\T\overline{G}\calX$ and
$(\calX^\T\calX )^{-1}$?

Write $N_{\partial C_M}$ for the outward unit normal to
$\partial C_M$. Since
$$
  \calX^\T \overline{G}(x)\calX
  = {\bigg( \Big\< \Big( {\D^2I(x)\over |\D I(x)|}
                   -\d N_{\partial C_M}(x)\Big) X_i \, , X_j
            \Big\>\bigg)}_{1\leq i,j\leq n-1} \, ,
$$
we first compute $\< \D^2I(x)X_i,X_j\>$, $1\leq i,j\leq d-1$. The
terms $\< \d N_{\partial C_M}X_i,X_j\>$, $1\leq i,j\leq d-1$ can
be computed with formula (9.1.8) and depend on the curvature of
$M$ via its second fundamental forms $\Pi_M^\nu$ and
$\Pi_M^{X_r}$, $m+1\leq r\leq d-1$. 

The first fundamental form $\calX^\T\calX
= (\< X_i,X_j\> )_{1\leq i,j\leq n-1}$ can be computed in 
the same way. One may notice however that 
for $i=1,\ldots , m$ and $j=m+1,\ldots , d-1$,
$$
  \< X_i,X_j\>={\< \nu_i,X_j\>\over \< \nu,p\>}
  + \d\phi_j(p)\cdot p_i + \sum_{m+1\leq r\leq d-1}\phi_r
  \< \d X_r\cdot p_i\, ,X_j\> \, .
$$
Moreover, if $i,j=m+1,\ldots , d-1$, then 
$$
  \< X_i,X_j\> =\delta_{i,j}\, .
  \qquad \hbox{ (Kronecker symbol)}
$$
Finally, if $i,j=1,2,\ldots , m$, the expression of $\<
X_i,X_j\>$ involves again the third fundamental form of $M$. 

At this point, it does not seem possible to push the abstract
calculation much further. The author hopes that it is clear that
the tail behavior of $X(M)$ is governed by the differential
geometry of $M$ immersed in $\RR^d$, a somewhat known fact. All
the calculations can be implemented on a computer.

Some simplifications may occur in
some specific cases. We now discuss some important ones.

Assume for instance that $I$ is a radial function, namely that
\setbox1=\vtop{\hsize=3.5in\noindent%
  $I(x)=J\big(|x|^2\big)$ for a function 
  $r\in [0,\infty )\mapsto J(\sqrt{r})$ which is convex on $\RR^+$, 
  and increasing.}
$$
 \raise 7pt\box1\eqno{(9.2.4)}
$$
Under (9.2.4), $I$ is minimal if and only if $|x|$ is. Thus, the
normalization $I\big( X(u)\big)$ minimum at $u_{m+1}=\cdots
=u_{d-1}=0$ forces $\phi_j=0$. Moreover, at the minimum, i.e., on
$\calD_{C_M^{\rm c}}$, the norm $|X(u)|$ is minimum, and
therefore,
$$
  X\perp X_i \, , \qquad i=1,2,\ldots , d-1 \hbox{ on } 
  \calD_{C_M^c} \, .
$$
Since $\nu_i$ is orthogonal to $X_i$ for $i=1,2,\ldots , d-1$, the vectors $X$
and $\nu$ are collinear on $\calD_{C_M^{\rm c}}$. 
Consequently, $X$, $\nu$,
$p$ are collinear on $\calD_{C_M^{\rm c}}$, and all point outward
$\Lambda_{I(C_M^{\rm c})}$.

Notice also that $|\nu |=1$ forces $\< \nu,\nu_i\> = 
\<p\, ,\nu_i\>=0$.

Equation (9.1.5) becomes very simple then, namely
$$
  X_i={\nu_i\over \< \nu,p\>} = {\nu_i\over |p|}
  \qquad\hbox{ on } \calD_{C_M^{\rm c}}\hbox{ and for } 
  i=1,2,\ldots , m \, .
$$
It follows that along $\calD_{C_M^{\rm c}}$, the first fundamental
form of $\partial C_M$ is given by
$$
  \< X_i,X_j\> =
  \cases{ |p|^{-2}\< \nu_i\, ,\nu_j\> \, , 
		& $i,j= 1,2,\ldots  , m$. \cr
  \noalign{\vskip 2mm}
	  |p|^{-1} \< \nu_i\, ,X_j\> \, , 
		& $i=1,2,\ldots , m$, $j=m+1, \ldots ,d-1$.\cr
  \noalign{\vskip 2mm}
	  \delta_{i,j}\, , 
		& $i,j=m+1,\ldots , d-1$.\cr }
$$
The second fundamental form of $\partial C_M$ also undergoes
some simplifications. Indeed equation (9.1.8) becomes
$$
  \< \d N_{\partial C_M}X_i,X_j\>
  =\cases{ 0 & if $i\vee j=m+1,\ldots , d-1$. \cr
  \noalign{\vskip 2mm}
	   -|p|^{-2}\<p_i,\nu_j\> & if $i,j=1,2,\ldots , m$.\cr}
  \eqno{(9.2.5)}
$$
--- remember that the normalization $u_{m+1}=\ldots = u_{d-1}=0$ at
the minimum forces $\phi_j=0$ on $\calD_{C_M^c}$.

\bigskip

\state{REMARK.} Be careful when using (9.2.5). The $X_j$'s, for 
$j=m+1,\ldots ,d-1$ have unit norm. But for $j=1,\ldots , m$, the
norm of $X_j$ may not be $1$. If one wants an expression of the fundamental
form in a basis of unit vectors, one should divide $X_j$ by its norm in 
(9.2.5). Nothing guarantees that the $X_j$'s for $j=1,\ldots ,m$ are
orthogonal to the $X_j$'s for $j=m+1,\ldots , d-1$. Thus to use
(9.2.5) in computations, one needs to use some form of orthogonalization
technique to obtain the matrix $\d N_{\partial C_M}$ in an orthonormal
basis.

\bigskip

\state{REMARK.} 
If we rescale $M$ to $\lambda M$, then $C_M$
becomes $C_{\lambda M}=C_M/\lambda$. So, the second fundamental
form of $C_{\lambda M}$ should be proportional to $\lambda$.
Consequently, the expression on the right hand side of (9.2.5)
is homogeneous in $\lambda$. With our choice, the parameterization
of $C_{\lambda M}$ is $X_\lambda(u)=X(u)/\lambda$. Therefore
$X_{\lambda ,i}=X_i/\lambda$. For $i,j=1,2,\ldots , m$, the
left hand side of $(9.2.5)$ written for $C_{\lambda M}$ reads
$$
  \< \d N_{\partial C_{\lambda M}}X_{\lambda,i}\, ,
     X_{\lambda,j}\>
  = \Big\< \d N_{\partial C_{\lambda M}}{x_i\over \lambda}\, ,
           {x_j\over \lambda}\Big\>
  \, .
$$
Thus, for $C_{tM}$, formula (9.2.5) is
$$
  \< \d N_{\partial C_{\lambda M}}X_i,X_j\>
  = \lambda\< p_i\, ,N_j\> /|p|^2 \, .
$$
This is indeed homogeneous in $\lambda$ since $X_i$ and $X_j$ do
not depend on $\lambda$.

\bigskip

Similarly to what happens for the second fundamental form of
$\partial C_M$ along $\calD_{C_M^{\rm c}}$, the second fundamental
form of the level sets $\Lambda_{I(C_M^{\rm c})}$ undergoes a great
simplification when $I$ is a radial function. Since the level
sets are spheres, their second fundamental form at $X$ is
$$
  \Pi_{\Lambda_{I(C_M^{\rm c})}}(X)(X_i,X_j)
  = |X|^{-1}\< X_i,X_j\> 
  = |p|\< X_i,X_j\>
  \qquad\hbox{ on } \calD_{C_M^{\rm c}}
$$
--- recall that for a radial $I$, we proved that $X(u)$ is
collinear to $p(u)$, and 
so $1=\< X(u),p(u)\> = |X(u)|\, |p(u)|$.

Further simplifications may occur by a good choice of the
parameterizations. For example, it may happen that $\calD_{C_M^{\rm c}}$ is
is parametrized by
$(u_1,\ldots , u_k,0,\ldots ,0)$, 
where $k=\dim\calD_{C_M^{\rm c}}$. In the case $\dim M=1$,
we can take $X_3$ to be collinear to the torsion vector of the
curve $M$. This ensures that $\nu_i=\nu'$ is orthogonal to
$X_4,\ldots , X_{d-1}$. The case $m=d-1$ is also rather
specific. There is a vast number of possible specializations
where more or less remarkable formulas can be obtained. However,
it is not obvious that such extensive developments would
bring more insight. They may be worthwhile for some specific
applications and we will see some in chapters 10--12.

\bigskip

To conclude, we mention that when Theorem 7.5 applies, it
shows that the distribution of $X/t$ given $X(M)$ larger than
$t$ converges to a distribution supported by 
$\calD_{C_M^{\rm c}}=\{\, x\in\partial C_M :I(x)=I(C_M^{\rm c})\,\}$
as $t$ tends to infinity. To a point $x$ in
$\calD_{C_M^{\rm c}}$ correspond points $m$ in $M_0$ such that $\<
x,m\>=1$. Proposition 9.1.7.iii shows that as $x$ varies in
$\calD_{C_M^{\rm c}}$, these $m$'s vary among the points in $M_0$
minimizing $I_\sbullet$. This implies that when $M$ is closed,
$\{\, m\in M_0 : X(m)=X(M)\,\}$ given $X(M)\geq t$ is a random
closed set whose distribution given $X(M)\geq t$ tends to be
concentrated on points in $M_0$ that minimize $I_\sbullet$.
When $X(M)$ is achieved at a unique point $m(X)$ in
$M_0$, Theorem 7.5 even gives the limiting conditional 
distribution of $\arg\max_{m\in M_0}X(m)$ given $X(M)\geq t$, 
as $t$ tends to infinity. It is the image measure by $m$ 
of the limiting conditional distribution of $X/t$ given 
$X\in tC_M^{\rm }$; this latter limiting distribution is 
given in Theorem 7.5. In particular, if $I_\sbullet$ is 
minimum at a unique point $m_*$ of $M_0$, then 
$\arg\max_{p\in M} X(p)$ given $X(M)\geq t$ converges in
probability to $m_*$ as $t$ tends to infinity..

\bigskip

%\chapternumber={9}
%\chaptername={Application to random linear forms}
\noindent{\fsection 9.3. Example with heavy tail
distribution.}
%\chapternumber={9.3}
\section{9.3}{Example with heavy tail}

\bigskip

\noindent In this section we consider a random vector
$X=(X_1,\ldots , X_d)$ in $\RR^d$ with independent coefficients, 
all having a Student-like distribution with parameter $\alpha$.
Thus, the distribution of $X_i$ is absolutely continuous with
respect to the Lebesgue measure, and satisfies
$$
  P\{\, X_i\leq -x\,\} \sim P\{\, X_i\geq x\,\}
  \sim K_{s,\alpha}\alpha^{(\alpha-1)/2}x^{-\alpha}
$$
as $x$ tends to infinity, for some constant $K_{s,\alpha}$.

The estimation of $P\{\, X\not\in C\,\}$ for an arbitrary convex set
$C$ containing a neighborhood of the origin turns out to be amazingly
simple. Recall that $(e_1,\ldots , e_d)$ denotes the canonical
basis of $\RR^d$. It is convenient to introduce the following
terminology.

\bigskip

\state{DEFINITION.} {\it%
  Let $S$ be a set in $\RR^d$, containing a neighborhood
  of the origin. A point of the form $\lambda e_i$ on 
  $\partial S$ is called an axial point of $\partial S$.}

\bigskip

As their name suggests, axial points of $\partial C$ are points on
$\partial C$ which lie on a canonical axis. Notice that if $C$ is not
the whole space $\RR^d$, then $\partial C$ has at least one axial
point. Moreover, $\partial C$ has at most $2d$ axial points. These
lower and upper bound can be achieved. For instance, the half space
$\{\, x\in\RR^d : x_1\leq 1\,\}$ has a unique axial point, $e_1$; and
the centered unit ball has $2d$ of them, namely plus or minus the
vectors of the canonical basis.

\bigskip

\state{9.3.1. THEOREM.} {\it%
  Let $C$ be a convex neighborhood of $0$ in $\RR^d$.
  Assume that $C$ is not the whole space $\RR^d$. Let
  $X$ be a random vector with independent
  and identically distributed components having a 
  Student-like distribution $S_\alpha$, with parameter $\alpha$, 
  and such that $S_\alpha (0)=1/2$. Then
  $$
    P\{\, X\not\in tC\,\} \sim t^{-\alpha} K_{s,\alpha} 
    \alpha^{(\alpha-1)/2}
    \sum_a |a|^{-\alpha}
    \qquad\hbox{ as } \ttoi \, , 
  $$
  where the sum $\sum_a$ is taken over all the axial points
  of $\partial C$.}

\bigskip

\stateit{Proof.} Notice first that if $\partial
C$ is not smooth at some of its axial points, say $a_1,\ldots,
a_k$, then we can sandwich $C$ between two smooth convex sets
with axial points $(1-\epsilon )a_1,\ldots , (1-\epsilon )a_k$
and $(1+\epsilon )a_1,\ldots , (1+\epsilon)a_k$ respectively.
For these approximating convex sets, the asymptotic formula in
Theorem 9.3.1 can be proved assuming that $\partial C$ is
smooth. We then let $\epsilon$ tend to $0$. So, there is no loss of
generality in assuming that $\partial C$ is smooth, which we do
from now until the end of the proof.

The proof of Theorem 9.3.1 is then essentially the same as that
of Theorems 8.2.1 and 8.2.10. It is hoped that the reader 
will be convinced that our unifying formalism is actually quite
convenient, even though, one more time, some specific examples
could be treated more easily with ad hoc methods.

Let $s_\alpha$ denote the density of a single $X_i$.
To prove Theorem 9.3.1, we need to approximate the integral
$$
  P\{\, X\not\in tC\,\}
  = \int_{\RR^d\setminus tC}s_\alpha (x_1)\ldots 
  s_\alpha (x_d)\,\d x_1\ldots \d x_d
$$
for large $t$. As in the proof of Theorem 8.2.1, extend the
function $\pisa$ to $\RR^d$ by considering it acting
componentwise on each coordinate. Making the change of variable
$Y=\pisa (X)$ leaves us to approximate
$$
  \int_{\pisa (t{C}^{\rm c})}{e^{-|y|^2/2}\over (2\pi)^{d/2}} 
  \,\d y \, .
$$
This leads us to consider the convex function
$$
  I(y)={|y|^2\over 2} -\log (2\pi)^{d/2}
$$
and the sets\ntn{$A_t$}\ntn{$B_t$}
$$
  A_t=tC^{\rm c} \, , \qquad \hbox{ and } \qquad 
  B_t=\pisa (tC^{\rm c}) \, .
$$

Our first proposition hereafter evaluates $I(B_t)$ and locates
the points of interests in $B_t$ as far as minimizing the
function $I$ is concerned. To this aim, define\ntn{$\gamma$}
$$
  \gamma = \min\big\{\, |a|  : a\hbox{ axial point of } 
  \partial C\,\big\} \, .
$$

\bigskip

\state{9.3.2. PROPOSITION.} {\it%
  As $t$ tends to infinity, we have
  $$\displaylines{\qquad
    I(B_t)=\alpha \log t-{1\over 2}\log\log t 
    + \alpha \log\gamma 
    \cut
    - \log (K_{s,\alpha}\alpha^{\alpha/2} 2\sqrt{\pi})
    +\log (2\pi )^{d/2} + o(1) \, .
    \qquad}
  $$
  If $a$ is an axial point of $\partial C$, then
  $$
    \tau_{B_t}\big( \pisa (ta)\big)
    =\alpha \log\big( |a|/\gamma \big) + o(1)
    \qquad \hbox{ as } \ttoi \, .
  $$
 
  For any positive number $M_1$, the set $\partial B_t\cap 
  \Gamma_{I(B_t)+M_1\log\log t}$ lies in an
  $O(\sqrt{\log\log t})$-neighborhood of the points
  $\pisa (ta)$, where $a$ is an axial point of $\partial C$.
  Consequently, there exists a positive $M_2$ such that this set
  lies in the image through $\pisa$ of a 
  $(\log t)^{M_2}$-neighborhood of the axial points of $tC$.}

\bigskip

\stateit{Proof.} An easy application of the
expansion for $(\pisa )^2$ given in Lemma A.1.5 gives for any axial
point $a$ of $\partial C$,
\finetune{\hfuzz=1pt}
$$\eqalign{
    I\big( \pisa (ta)\big) 
    =\alpha \log t -{1\over 2}\log\log t +\alpha\log |a|
    & - \log (K_{s,\alpha}\alpha^{\alpha /2}2\sqrt{\pi}) \cr
    & +\log (2\pi )^{d/2} + o(1) \cr
  }
$$
\finetune{\hfuzz=0pt}
as $t$ tends to infinity.

Let $M_1$ be a positive real number. Consider a point $u$ in
$C^{\rm c}$ such that
$$
  I\big( \pisa (tu)\big) \leq \alpha \log t + M_1\log\log t
  \, . \eqno{(9.3.1)}
$$
Consider $\eta <1-(1/\sqrt{2})$. For $t$ large enough, the
inequality
$$\eqalign{
  I\big( \pisa (tu)\big)
  &\geq \alpha (1-\eta)\sum_{1\leq i\leq d} 
    \indic_{[t^{-\eta},\infty )}\big( |u_i|\big) 
    \log\big( t|u_i|\big) \cr
  &\geq \alpha (1-\eta )^2 \sharp \big\{\, 1\leq i\leq d :
    |u_i|\geq t^{-\eta}\,\big\}\log t \cr
  }$$
implies that for $u$ as considered,
$$
  \sharp\big\{\, 1\leq i\leq d : |u_i|\geq t^{-\eta}\,\big\}
  \leq (1-\eta )^{-2}\big( 1+o(1)\big) <2 \, .
$$
However, since $u$ is not in $C$, we must have
$$
  \sharp \big\{\, 1\leq i\leq d : |u_i|\geq
                  t^{-\eta}\,\big\}
  \geq 1
$$
for $t$ large enough. Consequently, $u$ has exactly one coordinate
larger than $t^{-\eta}$. Thus, it must be in a
$t^{-\eta}$-neighborhood of the canonical axes in $\RR^d$. Such a
point is of the form $u=\lambda a + v$ for some axial point $a$
of $\partial C$, some $\lambda \geq 1-t^{-\eta}|a|^{-1}$ and 
$v$ orthogonal to $a$ with $|v|\leq t^{-\eta}$. Consequently,
$$
  I\big( \pisa (tu)\big) 
  ={1\over 2}\pisa \big( t\lambda |a|\big)^2
  + {1\over 2}|\pisa (tv)|^2 +\log (2\pi )^{d/2} \, .
$$
This expression is asymptotically minimum for $v=0$ and $\lambda =1$ and $|a|$
minimum, i.e., $|a|=\gamma$. Combined with Lemma A.1.5, this 
gives the value of $I(B_t)$ up to $o(1)$ as $t$ tends to
infinity. Then,
the value for $\tau_{B_t}\big( \pisa (ta)\big)$ follows from
Lemma A.1.5 as well.

Still assuming (9.3.1), we must have
$$\displaylines{\qquad
    \alpha\log t+M_1\log\log t
    \geq {1\over 2} \pisa \Big( t|a|\big( 1+o(1)\big)\Big)^2 
  \cut
      + {1\over 2}\max_{1\leq i\leq d} \pisa (tv_i)^2 
      +\log (2\pi )^{d/2} \, .
  \qquad\cr}
$$
Hence, using Lemma A.1.5 to approximate 
$\pisa\Big( t|a|\big(1 +o(1)\big)\Big)^2$,
$$
  \max_{1\leq i\leq d} \pisa (tv_i)^2 \leq (M_1+1)\log\log t
$$
for $t$ large enough, that is $\pisa (tu)=O(\log\log t)^{1/2}$. 
Therefore, $\pisa (tv)$ is indeed in an 
$O(\log\log t)^{1/2}$-neighborhood of $\pisa (ta)$.

Given the expression for $\pisa$ in Lemma A.1.5, we also must
have
$$
  \max_{1\leq i\leq d} \alpha \log \big( t|v_i|\big)
  \leq (M_1+1)\log\log t \, ,
$$
that is, $\max_{1\leq i\leq d}|tv_i|
\leq (\log t)^{(M_1+1)/\alpha}$ for $t$ large enough, 
which is the last statement of the proposition. \hfill$\qed$

\bigskip

We can now try to calculate the asymptotic equivalent given by
Theorem 5.1. Given the proof of Proposition 9.3.2, it is natural
to try the projection of the axial points of $\partial B_t$ onto
$\Lambda_{I(B_t)}$ as a dominating manifold. So, let $\rho_t$ be
the radius of the ball $\Lambda_{I(B_t)}$, and set\ntn{$\calD_{B_t}$}
$$
  \calD_{B_t}=\Big\{\, \rho_t
  { \pisa (ta)\over |\pisa (ta)|} : a
  \hbox{ axial point of }\partial C\,\Big\} \, .
$$
Notice that because we assumed $S_\alpha (0)=1/2$, the equality
$\pisa (se_i)=\pisa (s)e_i$ holds for any canonical vector
$e_i$ of $\RR^d$ and any real number $s$. Equivalently, we have
$$
  \calD_{B_t}=\Big\{\, \rho_t {a\over |a|} : 
  a \hbox{ an axial point of }\partial C\,\Big\} \, .
$$
The dimension of $\calD_{B_t}$ is $k=0$.

From the values of $I(B_t)$ in Proposition 9.3.2, we infer that
$$\eqalign{
  \rho_t
  & =\big[ 2\alpha\log t-\log\log t+2\alpha\log \gamma
    - 2\log (K_\alpha\alpha^{\alpha /2}2\sqrt{\pi})
     +o(1)\big]^{1/2}
    \cr
  & \sim \sqrt{2\alpha\log t} \qquad\hbox{ as }\ttoi \, . \cr
  }$$
Putting all the pieces together, and assuming that we can verify
its assumptions, the approximation formula in Theorem 5.1 yields,
$$\eqalign{
  P(A_t)
  & \sim e^{-I(B_t)} (2\pi )^{(d-1)/2}
    \sum_a { \exp \Big( -\alpha \log \big( |a|/\gamma\big)\Big) 
	     \over \rho_t^{(d+1)/2} \det 
             \big( G_{B_t}(\rho_ta/|a|) \big)^{1/2} } \cr
  & \sim { \sqrt{\log t}\over t^\alpha} K_{s,\alpha}
    { \alpha^{\alpha /2}\sqrt{2}
      \over (2\alpha\log t)^{(d+1)/4} }
    \sum_a {1\over |a|^\alpha \big( \det\, G_{B_t}(\rho_t a/|a|)
     \big)^{1/2}} \, , \cr
  }$$
as $t$ tends to infinity, with the sum taken over all the axial points $a$ of
$\partial C$. So, it remains for us to calculate 
$\det\, G_{\smash{B_t}}(\rho_t a/|a|)$ for all the axial points 
of $\partial C$, and check the assumptions of Theorem 5.1. Our 
next lemma does half of the task.

\bigskip

\state{9.3.3. LEMMA.} {\it%
  For any axial point $a$ of $\partial C$,
  $$
    G_{B_t}(\rho_ta/|a|)\sim 
    {\Id_{\RR^{d-1}}\over\sqrt{2\alpha\log t}}\qquad\hbox{ as }
    \ttoi \, .
  $$}

\finetune{}%\bigskip

\stateit{Proof.} Let $a$ be an axial point of $\partial C$ 
and let $u$ be orthogonal to $a$.  Since the origin is in
the interior of the convex set $C$, the line $a\RR$ intersects
$\partial C$ transversally.  Therefore, we can parameterize 
$\partial C$ around $a$ by a ball in $\{\, a\,\}^\perp$, centered at $0$. 
In other words, there exists a smooth 
function $h_a : \{\, a\,\}^\perp\to\RR$ such that
$$
  p_a : u\in \{\, a\,\}^\perp\mapsto p_a(u)
  =\big( 1+h_a(u)\big)a + u\in\partial C
$$
defines a parameterization of $\partial C$ around $a$ --- i.e., 
for $|u|$ small enough, it is a parameterization. Since $u$ is
orthogonal to $a$ and $\pisa$ acts componentwise,
$$
  \pisa \big( tp_a(u)\big)
  = \pisa \Big( t\big(1+h(u)\big)\Big){a\over |a|} 
    + \pisa (tu) \, .
$$
If $I\big(\pisa (tp_a)\big)\leq I(B_t)+M_1\sqrt{\log\log t}$ 
--- the domain which will interest us after we choose $c_{B_t}$
--- Proposition 9.3.2 asserts that
$|u|\leq (\log t)^{M_2}/t$ for some $M_2$.
Since $\partial C$ is smooth, we have $|h(u)| \leq M_3
|u|\leq M_3(\log t)^{M_2}/t$ for some positive $M_3$ and in 
this range of $u$. Consequently
$$
  \pisa \Big( t\big( 1+h(u)\big)a\Big)
  = \pisa (ta) +o(1)
  \qquad \hbox{ as } \ttoi \, ;
$$
and in the range $|u|\leq (\log t)^{M_2}/t$,
$$
  \pisa \big( tp_a(u)\big)
  = \pisa (ta)+\pisa (tu)+o(1)
  \eqno{(9.3.2)}
$$
as $t$ tends to infinity.  Up to the term in $o(1)$, this last
equation defines a plane orthogonal to $\pisa (ta)$. Following the
proof of Lemma 8.2.8, it follows that
$$
  G_{B_t}\big( \pisa (ta)\big)
  \sim {\Id_{\RR^{d-1}}\over\rho_t}
  \sim {\Id_{\RR^{d-1}}\over \sqrt{2\alpha\log t}} \, .
  \eqno{\qed}
$$

\bigskip

Lemma 9.3.3 implies 
$$
  \det\, G_{B_t}\big( \rho_ta/|a|\big)^{1/2}
  \sim (2\alpha\log t)^{-(d-1)/4} \, ,
  \qquad\hbox{ as } \ttoi \, .
$$
With the estimate of $P(A_t)$ obtained before the statement 
of Lemma 9.3.3, we obtain
$$
  P(A_t)\sim{ K_{s,\alpha}\alpha^{(\alpha-1)/2}\over t^\alpha }
  \sum_a |a|^{-\alpha}\qquad \hbox{ as } \ttoi \, .
  \eqno{(9.3.3)}
$$
This is the asymptotic equivalent given in Theorem 9.3.1. Thus, it
remains for us to check the assumptions of Theorem 5.1. Taking the
risk of making the rest of the proof boring, we will do it in a
systematic way, showing that this is a rather easy task.

Our candidate for $c_{B_t}$ is $c_t=(d+1)\log\log t$. Indeed,
combining Propositions 2.1 and 9.3.2, we obtain
$$\eqalign{
  L\big( I(B_t)+c_t\big)
  &\leq c_0 e^{-I(B_t)-c_t}\big( 1+I(B_t)+c_t\big)^d \cr
  & = o(t^{-\alpha}) \qquad\hbox{ as } \ttoi \, .
    \cr
  }$$
Given (9.3.3), this choice of $c_t$ guarantees that (5.4) holds.

Given Proposition 9.3.2 and the way we constructed $\calD_{B_t}$ ---
using a projection on the sphere of radius $\rho_t$ --- the set
$\underline{B}_{t,M}$ lies on $O(\log\log t)$-neighborhood of
$\calD_{B_t}$ on $\Lambda_{I(B_t)}$. Since the radius of injectivity
of $\Lambda_{I(B_t)}$ --- a sphere of radius $\rho_t\sim
\sqrt{2\alpha\log t}$ --- is $\pi\rho_t/2$, assumption (5.2) holds for
$t$ large enough.

Assumption (5.5) is almost plain. Equation (9.3.2) is the
analogue of Lemma 8.2.2 or 8.2.12. It shows that the boundary
$\partial B_t$ near an axis point $\pisa (ta)$ is a plane
orthogonal to $a$ up to an $o(1)$-term. The very same argument as
that used in the proof of Theorems 8.2.1 and 8.2.10 ensures that
(5.5) is verified here.

(5.6) is clear.

(5.7) follows in the very same way as in the proof of Theorem
8.2.1. We still have $t_{0,M}(p)=O(\log\log t)$ while
$K_{\rm max}(q,t_0)=O(\log t)^{-1}$ for $q\in\calD_{B_t}$.

(5.8) is plain, for a sphere has positive Ricci curvature in
$\RR^d$.

(5.9), (5.10) and (5.11) follow exactly as in the proof of
Theorem 8.2.1.

(5.12) and (5.13) are plain as well given the proof of Theorem
8.2.1 or 8.2.10, and this concludes the
proof of Theorem 9.3.1.\hfill$\qed$

\bigskip

\state{REMARK.} Note that if we drop the assumption $S_\alpha (0)=1/2$
in Theorem 9.3.1, that is $S_\alpha$ of median zero, axial points are
not mapped anymore to axial points by $\pisa$. The argument we
developed would still be valid though, since the asymptotic expansion
for $\pisa$ shows that the component of $\pisa (te_i)$ 
on $\{ e_i\}^\perp$ is asymptotically negligible compared to that on
$e_i\RR$. Thus, the conclusion of Theorem 9.3.1 still holds true
without assuming $S_\alpha (0)=1/2$.

\bigskip

From Theorem 9.3.1, we deduce the following limiting behavior
of the conditional distribution of $X$ given $X\not\in tC$.

\bigskip

\state{9.3.4. COROLLARY.} {\it%
  Let $C$ be any convex neighborhood of $0$ in $\RR^d$, such that
  $C\neq\RR^d$. If $X$ is a random vector with independent
  and identically distributed components having a Student-like 
  distribution with parameter $\alpha$, then the distribution of 
  $X/t$ given $X\not\in tC$ converges weakly* to
  $$
    \sum_a |a|^{-\alpha}P_a\big/ \sum_a|a|^{-\alpha} \, ,
  $$
  where the sums are taken over all axial points $a$ of
  $\partial C$ and $P_a$ is a Pareto distribution concentrated on
  $a\RR^+$, whose cumulative distribution function is
  given by
  $$
    P_a\big\{\, \< X, a/|a|\>\geq |a|+\lambda\,\big\}
    = {|a|^\alpha\over \big( |a|+\lambda\big)^\alpha} \, ,
    \qquad \lambda\geq 0 \, .
  $$}%
\finetune{
%\bigskip
}
\stateit{Proof.} Let $\lambda$ be a positive number. Consider 
the convex set
$$
  D=\big\{\, x : \< x,a/|a|\>\leq |a|+\lambda \,\big\} \, .
$$
The set of all axial points of $\partial (D^{\rm c}\cap
C^{\rm c})$ is just
$\big\{\, a {|a|+\lambda\over |a|}\,\big\}$,
and $a$ is an axial point of the convex set $C$. Consequently,
applying Theorem 9.3.1 with the convex set $C$ and 
$(D^{\rm c}\cap C^{\rm c})^{\rm c}$
yields
$$\eqalign{
  P\{\, X\not\in tD|X\not\in tC\,\}
  & = P\big\{\, X\in t(D^{\rm c}\cap C^{\rm c})\,\big\} 
    \big/ P\{\, tC^{\rm c}\,\} \cr
  & \sim \big( |a|+\lambda \big)^{-\alpha}\big/ 
    \sum_a |a|^{-\alpha} \,  . \cr
  }$$
This is the result, for the conditional distribution of $X/t$
given $X\not\in tD$ converges trivially to $P_a$. \hfill$\qed$

\bigskip

Let us now go back to the study of processes of the 
form $\< X,p\>$, for $p$ in a set $M$ of $\RR^d$. Theorem 9.3.1 
gives us the asymptotic 
behavior of $P\big\{\, X(M)\geq t\,\big\}$ as $t$ tends to
infinity. To obtain a more readable
statement, we need to express the axial points of $\partial C_M$
in term of $M$. This is done in the next result.

\bigskip

\state{9.3.5. PROPOSITION.} {\it%
  The set of all axial points of \smash{$\partial C_M$} coincide with the
  set of vectors $\epsilon a_{\epsilon , i}e_i$ where $\epsilon$
  is in \smash{$\{\, -1,+1\,\}$},
  and $i$ is such that
  $\epsilon \< p,e_i\> >0$ for some $p$ belonging to $M$, and 
  $$
    1/a_{\epsilon , i}=\sup\big\{\, \epsilon \,\< p,e_i\> : 
    p\in M\,\big\} \, .
  $$
  }

\bigskip

\stateit{Proof.} Let $a$ be an axial point of
$\partial C_M$. Necessarily $a=\epsilon |a|e_i$ for some
$\epsilon$ in $\{\, -1,1\,\}$, and $e_i$ a vector of the canonical
basis of $\RR^n$. Since $a$ is in $\partial C_M$, we have $\<
a,p\>\leq 1$ for all $p$ in $M$, and $\sup\big\{\, \< a,p\>  : 
p\in M\,\big\} =1$. Thus, $\sup\big\{\, \epsilon |a|\< e_i,p\> :
p\in M\,\big\} =1$ and the result follows. \hfill$\qed$

\bigskip

If we have a parameterization $f(t)=\big( f_1(s), \ldots ,
f_d(s)\big)$ of $M$ indexed by $s$ in some set $S$, it is 
particularly easy to relate
the behavior of $f(s)$ to the geometry of the set $M$ captured
in Proposition 9.3.5. This yields immediately the following
result, where the reader will notice that the function $f$ is
completely arbitrary --- no need for measurability, or any kind of
regularity whatsoever!

\bigskip

\state{9.3.6. THEOREM.} {\it%
  Let $S$ be a set, and $f(s)=\big( f_1(s),\ldots , f_d(s)\big)$
  be a bounded function defined on $S$. Let $X$ be a random 
  vector in $\RR^d$, with independent and identically distributed
  components having a Student-like distribution with
  parameter $\alpha$. Let
  $$
    X(S)=\sup_{s\in S}X_1f_1(s)+\cdots + X_df_d(s) \, .
  $$
  Then,
  $$
    P\big\{\, X(S)\geq t\,\big\} 
    \sim K_{s,\alpha} {\alpha^{\alpha-1/2}\over t^\alpha}
    \sum_{\scriptstyle 1\leq i\leq d
          \atop\scriptstyle\epsilon\in\{-1,1\}}
    c_{\epsilon ,i}^{\alpha} \qquad \hbox{ as }
    \ttoi \, ,
  $$
  \finetune{\vskip -.1in\noindent} 
  where
  $$
    c_{\epsilon ,i} 
    =\sup\big\{\, \big(\epsilon f_i(s)\big)_+  : s\in
               S\,\big\}
    \, ,\qquad \epsilon\in\{\, -1,1\,\}\, , \, i=1,\ldots ,d\, .
  $$}

\finetune{}%\bigskip

Another way to interpret $c_{\epsilon ,i}$ is in looking at the
projection of $M$ on the $i$-th canonical axis $e_i\RR$. The
value $c_{\epsilon ,i}$ is $0$ if this projection is
concentrated on the set $-\epsilon e_i\RR\,$; otherwise
$c_{\epsilon , i}$ is the coordinate of the largest point of
this projection.

It is quite amusing to notice the following. Set
$f(s)=(1,\ldots , 1)$ for all $s$. Then $X(S)=X_1+\cdots + X_d$.
Theorem 9.3.6 implies that
$$
  P\{\, X_1+\cdots + X_d\geq t\,\} 
  \sim {K_\alpha\alpha^{\alpha -1\over 2}\over\lambda^\alpha}d
  \, .
$$
Using again Theorem 9.3.6 for $d=1$, we then infer
$$
  P\{\, X_1+\cdots + X_d\geq t \,\} 
  \sim d\big( 1-S_\alpha (t)\big)
  = dP\{\, X_1\geq t \,\} \, .
$$
This is a known asymptotic identity showing that the Student-like
distributions are subexponential!

\bigskip

\centerline{\fsection Notes}
\section{9}{Notes}

\bigskip

\noindent This chapter is connected with a huge literature. To
proceed in order, I first cannot quite believe that the
transform $I\mapsto I_\sbullet$ and its use are new. It
allows one to read the minimization of a convex function on the
complement of a convex set over its polar reciprocal.
However, I have not found it in the literature. Similarly to 
what is called the infimum convolution in convex analysis,
it would be natural to call $I_\sbullet$ the infimum Radon
transform of $I$. The relation $(I_\sbullet )^\sbullet =I$ in
Proposition 9.1.7.iv is an inversion formula for this infimum
Radon transform. 
 
Lemma 9.1.3 is essentially contained in Hassanis and Koutroufiotis
(1985). Though I refer to Schneider (1993) for convexity
theory, the differential viewpoint in section 9.1 is closer
to Bruce and Giblin (1992). The polar reciprocal is 
sometimes called the pedal surface. It would be 
desirable to connect further the global properties of
$M_0$ and $C_M$.

Proposition 9.2.1 is a generalization of the finite
dimensional version of Fernique (1970) and Landau and Shepp
(1970). There are many proofs in the Gaussian setting, and
the lectures by Ledoux (1996) are most illuminating. The
current literature on related problems is connected with
notions such as concentration of measure and a set of
inequalities: isoperimetric, Sobolev logarithmic,
Poincar\'e. A couple of pointers to this literature are
Talagrand (1995) and Bobkov and Ledoux (1997, 2000). However, it
is not quite clear that Proposition 9.2.1 can be recovered
from the existing results in the literature. An open question is
if the conclusion of Proposition 9.2.1 holds for any convex
functions.

The result of Proposition 9.2.1 also makes sense in infinite dimension,
using of course the dual space to define $I_\sbullet$. But the proof 
given here breaks down in infinite dimensions.

The part of section 9.2 following Proposition 9.2.1 is
connected with a flourishing literature on the Gaussian
case and a few related distributions such as the chi-square. 
The point of view given here is close to an abstraction
of Diebolt and Posse (1996). A radically different line of
investigation is in Piterbarg (1996). A most interesting survey 
of the literature on supremum of Gaussian processes is in
Adler (2000).

In the heavy tail case, Theorem 9.3.6 seems to be part of
the folklore; but I have not found it in the literature. Both
Theorems 9.3.1 and 9.3.6 can be proved directly by ad hoc methods.

Concerning sums of heavy tailed random variables, Bingham,
Goldie and Teugels (1987) contains invaluable material. I
learned about those things in part in Broniatowski and Fuchs
(1995). Much nicer results than the one presented in this
section exist, including second and higher order formulas.
But, unfortunately, the accuracy of these expansions is 
incredibly poor, especially when the tail parameter 
$\alpha$ is large.

It is interesting to rework the proofs of this specific
chapter assuming that the cumulative distribution function 
of $X_i$ is $1-F(x)=x^{-\alpha}\ell (x)$ for some slowly
varying function $\ell$ at infinity. Linearity of the
functional considered here leads to neat results --- but they
depend heavily on the linearity!

\vfill\eject

~\vskip 28pt
\chapter{10}{Random matrices}
{\fchapter\hfill 10.\quad Random matrices}

\vskip 54pt

\noindent In this chapter we consider a matrix
$X={(X_{i,j})}_{1\leq i,j\leq d}$ with random coefficients
that are independent and have the same distribution. Many quantities
associated to $X$ are of interest. For instance, its trace
$\trace (X)$, its determinant $\det (X)$, or its operator norm
$\| X\|$. All these quantities have in general complicated
distributions which cannot be calculated very explicitly.
Hence, it makes sense to investigate their tail behavior.

Before going further, let us mention that the trace of $X$ 
is nothing but a sum of independent and identically distributed 
random variables. Results from section 9 give the tail
approximations for $P \big\{\, \trace (X)\geq t\,\big\}$ when the
coefficients of $M$ have a Weibull or a Student-like distributions.
Therefore, in this section, we will concentrate on the
determinant and on the norm of $X$. We will see that their 
tail behavior often turns to be quite interesting, if not 
fascinating.

Throughout this chapter, it will be convenient to think of
matrices as vectors in $\RR^{n^2}$ as well as linear operators
acting on $\RR^n$. In particular, we denote by $\Eij$ the
canonical orthonormal basis of $\RR^{n^2}$ viewed as matrices.
Thus, $\Eij$ denotes the matrix with $1$ on the $(i,j)$-entry,
and $0$ elsewhere. In other words.
$$
  \Eij={ (\delta_{(i,j),(k,l)})}_{1\leq k,l\leq n}
$$
where $\delta_{u,v}$ is the Kronecker symbol.

Let $\MnR$\ntn{$\MnR$} denote the set of all $n\times n$ matrices
with real coefficients. Also, we write $\glnr$\ntn{$glnr$} for the group
of all invertible matrices in $\MnR$, that is the linear
group.

Since our method is differential geometric, we will need the
differential and Hessian of the determinant as a map
from $\glnr$ to $\RR\,$. For the sake of completeness, 
we recall them.

\bigskip

\state{10.0.1. LEMMA.} {\it%
  For all $x$ in $\glnr$ and $h,k$ in $\MnR$,
  $$
    \D\det (x)h = \det (x)\trace (x^{-1}h) \, , 
  $$
  and
  $$
    \D^2\det (x)(h,k)=\det (x)\trace(x^{-1}h)\trace (x^{-1}k)
    -\det (x)\trace (x^{-1}hx^{-1}k) \, .
  $$}
\finetune{
%\bigskip
}
\stateit{Proof.} Let $x$ be an invertible matrix. Since $\det
(x+h)=\det (x)\det (\Id+x^{-1}h)$, it is enough to calculate
$\D\det (\Id )$ and $\D^2\det (\Id )$.

Let\ntn{$\frakSigma_n$} $\frakSigma_n$ denote the group of 
permutations of $n$ elements. Define
$\frakSigma_n^0=\{\, \Id\,\}\subset\frakSigma_n$ 
and let $\frakSigma_n^1$ be the subset of $\frakSigma_n$ made of all
transpositions.  The signature of a permutation $\sigma$ is $\epsilon
(\sigma )=+1$ (resp. $-1$) if $\sigma$ is the composition of an even
(resp. odd) number of transpositions. We have
$$
  \det (\Id +sh)
  = \sum_{\sigma\in\frakSigma_n}\epsilon (\sigma )\prod_{1\leq i\leq n}
  (\delta_{i,\sigma (i)}+sh_{i,\sigma (i)})\, .
$$
This sum over $\frakSigma_n$ can be decomposed as a sum 
over $\frakSigma_n^0$, plus one over $\frakSigma_n^1$, 
plus a remainder term. The sum over $\frakSigma_n^0$ has a
unique term,
$$\displaylines{\qquad
  \epsilon (\Id )\prod_{1\leq i\leq n} (\delta_{i,i}+sh_{i,i})
  \cut
  \eqalign{%
  =\,& 1+s\sum_{1\leq i\leq n}h_{i,i}
       +s^2\sum_{1\leq i<j\leq n}h_{i,i}h_{j,j} + O(s^3) \cr
  =\,& 1+s\, \trace (h)+{s^2\over 2}\Big( \sum_{1\leq i,j\leq n}
       h_{i,i}h_{j,j}-\sum_{1\leq i\leq n}h_{i,i}^2\Big) 
       + O(s^3) \cr
  =\,& 1+s\, \trace (h) +{s^2\over 2}
       \Big( \trace (h)^2-\sum_{1\leq i\leq n}h_{i,i}^2\Big)
       + O(s^3) \, .\cr}
  \qquad}$$
Next, we also obtain
$$\displaylines{\qquad
    \sum_{\sigma\in\frakSigma_n^1} \epsilon (\sigma )
    \prod_{1\leq i\leq n} (\delta_{i,\sigma (i)}
    +sh_{i,\sigma (i)})
  \cut
  \eqalign{
  =&-\sum_{1\leq i<j\leq n}
    \prod_{1\leq k\leq n\atop k\not\in \{ i,j\}} (1+sh_{k,k})
    s^2h_{i,j}h_{j,i} \cr
  =& -s^2\sum_{1\leq i<j\leq n} h_{i,j}h_{j,i} + O(s^3) \cr
  =& -{s^2\over 2} \Big(\sum_{1\leq i,j\leq n} h_{i,j}h_{j,i}
    -\sum_{1\leq i\leq n}h_{i,i}^2\Big) + O(s^3) \cr
  =& -{s^2\over 2}
    \Big( \trace (h^2)-\sum_{1\leq i\leq n}h_{i,i}^2\Big) 
    + O(s^3) \, .\cr}
  \qquad}$$
If $\sigma$ is in $\frakSigma_n\setminus (\frakSigma_n^0\cup \frakSigma_n^1)$, 
at least 3 integers in $\{\, 1, 2, \ldots , n\,\}$ are not invariant under
$\sigma$. For such permutation
$$
  \prod_{1\leq i\leq n} (\delta_{i,\sigma (i)}+\epsilon 
  \delta_{i,\sigma (i)})=O(s^3) \, .
$$
It follows that
$$
  \det (\Id +sh)= 1 + s\, \trace (h) + {s^2\over 2} 
  \big( \trace (h)^2-\trace (h^2)\big) + O(s^3)
$$
as $s$ tends to $0$. Consequently, for $x$ in $\glnr$,
$$
  \D\det (x) (h) 
  = \det (x)\D\det (\Id )(x^{-1}h)
  = \det (x) \trace (x^{-1}h) \, .
$$
Also, we have
$$
  \D^2\det (\Id ) (h,h) = \trace (h)^2-\trace (h^2) 
$$
and by polarization
$$
  \D^2\det (\Id ) (h,k) = \trace (h) \trace (k)-\trace (hk) \, .
$$
Consequently,
$$\eqalign{
  \D^2\det (x) (h,k) 
  & = \det (x) \D^2\det (\Id ) (x^{-1}h,x^{-1}k) \cr
  & = \det (x) \big( \trace (x^{-1}h)\trace (x^{-1}k)
    -\trace (x^{-1}hx^{-1}k)\big) \cr
  }$$
as claimed. \hfill$\qed$

\bigskip

We conclude this section by a trivial but useful formula. When
needed, we write $\<\cdot ,\cdot \>_{\RR^{n^2}}$ for the inner
product in the Euclidean space $\RR^{n^2}$. We use
the tensor product notation, $E^{i,j}\otimes E^{l,m}$ to denote
the bilinear form $x,y\in\MnR\mapsto E^{i,j}\otimes
E^{l,m}(x,y)=\< E^{l,m},y\>_{\RR^{n^2}}
\<E^{i,j},x\>_{\RR^{n^2}}$.

\bigskip

\state{10.0.2. LEMMA.} {\it%
  On the basis $\Eij$, the bilinear map $(h,k)\in \MnR\mapsto 
  \trace (x^{-1}hx^{-1}k)$ has the form
  $$
    \sum_{1\leq i,j,l,m\leq n}
    (x^{-1})_{m,i}(x^{-1})_{j,l} \Eij\otimes\Elm \, .
  $$}

\bigskip

\stateit{Proof.} It is straightforward,
$$\eqalignno{
  \trace (x^{-1}hx^{-1}k)
  & = \sum_{1\leq i,j,l,n\leq n} (x^{-1})_{i,j} h_{j,l}
    (x^{-1})_{l,m}k_{m,i} \cr
  & = \sum_{1\leq i,j,l,m\leq n} (x^{-1})_{i,j} (x^{-1})_{l,m}
    \Ejl\otimes \Emi (h,k) \, . &$\qed$\cr
  }$$

\bigskip

Viewing $\MnR$ as $\RR^{n^2}$, we have a natural inner product

$$
  x,y\in \MnR\mapsto
  \< x,y\> = \sum_{1\leq i,j\leq n} x_{i,j}y_{i,j}
  = \trace (x^\T y) = \trace (xy^\T ) \, .
$$
Consequently, viewing $\D\det (x)$ in $\RR^{n^2}$, Lemma 10.0.1
implies that
$$
  \D\det (x)=\det (x)(x^{-1})^\T \in\RR^{n^2} \, .
  \eqno{(10.0.1)}
$$

\bigskip

%%%%%%%%%%%%%%%%%%%%%%%%%%%%%%%%%%%%%%%%%%%%%%%%%%
%% Section 10.1
%%%%%%%%%%%%%%%%%%%%%%%%%%%%%%%%%%%%%%%%%%%%%%%%%%

\noindent{\fsection 10.1. Random determinants, light tails.}
%\chapternumber={10.1}
\section{10.1}{Random determinants, light tails}

\bigskip

\noindent In this section, we consider a random matrix
$X={(X_{i,j})}_{1\leq i,j\leq n}$ where the $X_{i,j}$ are
independent and identically distributed, each having a
symmetric Weibull-like density\ntn{$w_\alpha (\cdot )$}
$$
  w_\alpha (u) 
  = {\alpha^{1-(1/\alpha )}\over 2\Gamma (1/\alpha)}
  \exp\Big({-|u|^\alpha\over \alpha}\Big)\, , \qquad u\in\RR 
  \, ,\, \alpha >1 \, .
$$
One of the aim of this section is to show that the Gaussian 
case, obtained for $\alpha=2$, is rather specific.
The main reason is of course the invariance of the Gaussian
distribution under the special orthogonal group.

Given the densities $w_\alpha$ of interest, let us define\ntn{$I(\cdot )$}
$$
  I(x)={1\over \alpha}\sum_{1\leq i,j\leq n}|x_{i,j}|^\alpha 
  \, , \qquad x\in\MnR\equiv \RR^{n^2} \, .
$$
Furthermore, define\ntn{$A_t$}
$$
  A_t=\big\{\, x\in\MnR  : \det X\geq t\,\big\} 
  = t^{1/n}A_1 \, .
$$

Since $I$ is homogeneous, we can use Theorem 7.1 in order to
approximate
$$
  P\{\, \det X\geq t\,\} 
  = \bigg( { \alpha^{1-(1/\alpha)}\over 2\Gamma (1/\alpha) }
  \bigg)^{n^2} \int_{t^{1/n}A_1} e^{-I(x)} \d x \, .
$$
This requires us to compute the dominating manifold $\calD_{A_1}$.
Unfortunately, I have not been able to do so in general. 
The following result will rely on an explicit calculation in some 
special cases, and a conjecture in general.

Our first lemma provides a necessary condition for a matrix to
be in $\calD_{A_1}$. We denote\ntn{$\SLnR$} by $\SLnR$ the special 
linear group on $\RR^n$, that is the group of all matrices of 
determinant $1$.

\bigskip

\state{10.1.1. LEMMA.} {\it%
  If $x$ is an $n\times n$ real matrix minimizing $I(x)$ 
  subject to the constraint $\det x\geq 1$, then 
  $x$ is of determinant $1$, that is belong to $\SLnR$.
  Moreover, for such a matrix,
  $$
    (x^{-1})_{i,j} 
    = {1\over \lambda}\sign (x_{j,i}) |x_{j,i}|^{\alpha -1}
    \, , \qquad 1\leq i,j\leq n \, ,
    \eqno{(10.1.1)}
  $$
  where $n\lambda =\min \big\{\, I(x) : \det\, x=1\,\big\}$.
  }

\bigskip

\stateit{Proof.} Since $\det (\lambda x)=\lambda^n\det x$ 
and $I(\lambda x)=|\lambda |^\alpha I(x)$,
we clearly have $\det x=1$ at the constrained minimum. So, we
need to find $\inf \big\{\, I(x) : \det\, x=1\,\big\}$. At the
minimum, the normal vector of the level set of $I$ and $\{\, x\, :
\, \det\, x= 1\,\}$ are collinear. Using (10.0.1), this condition
writes
$$
  \Big( {1\over\lambda} \sign (x_{i,j})|x_{i,j}|^{\alpha -1}
  \Big)_{1\leq i,j\leq n}
  = (x^{-1})^\T
$$
for some nonzero $\lambda$.

Since $1=(x^{-1}x)_{i,i}$ for all $i=1,2,\ldots , n$, we obtain
$$\eqalign{
  n & = \sum_{1\leq i,k\leq n} (x^{-1})_{i,k} x_{k,i}
      = {1\over \lambda} \sum_{1\leq i,k\leq n} \sign (x_{k,i})
      |x_{k,i}|^{\alpha -1} x_{k,i} \cr
    & = {I(x)\over\lambda} \, , \cr
  }$$
and the result follows. \hfill$\qed$

\bigskip

In general, I have been unable to solve (10.1.1) explicitely.
But a very partial solution can be given, suggesting that the
general one may be quite involved. As customary, we denote by\ntn{$\SOnR$}
$\SOnR$ the special orthogonal group, that is the subgroup of
$\SLnR$ of all matrices $X$ such that $X^\T X=\Id$.

\bigskip

\state{10.1.2. LEMMA.} {\it%
  (i) If $\alpha =2$, then $\da1 = \SOnR$.

  \noindent (ii) If $n=2$ and $\alpha <2$, then
  $$
    \da1=\bigg\{\, \pmatrix{ 0&-\epsilon \cr \epsilon & \hfill 0 \cr}\,
    ,\, \pmatrix{ \epsilon &0\cr 0&\epsilon \cr} : 
    \epsilon\in\{\, -1,1\,\}\,\bigg\} \, .
  $$
  
  \noindent (iii) If $n=2$ and $\alpha >2$, then
  $$
    \da1=\bigg\{\, {1\over \sqrt{2}}
	\pmatrix{ \epsilon_1 &-\epsilon_2 \cr 
		  \epsilon_2 & \hfill\epsilon_1 \cr }  : 
    \epsilon_1,\epsilon_2\in\{\, -1,1\,\}\,\bigg\} \, .
  $$}

\bigskip

\stateit{Proof.} (i) for $\alpha=2$, equation
(10.1.1) becomes 
$$
  (x^{-1})_{i,j}={x_{j,i}\over\lambda} \, \qquad
  1\leq i,j\leq n\, .
$$
Hence, $x^{-1}=x^\T /\lambda$ and $\Id = x^{-1}x=x^\T
x/\lambda$. Moreover, $\det\, x=1$ thanks to Lemma 10.1.1.

When $n$ is odd, we deduce that $1=\det\, \Id = \lambda^{-n}$.
Consequently, $\lambda =1$ and $x^{-1}=x^\T $. This 
proves $\da1\subset \SOnR$. Since $I$ is invariant under 
the action of $\SOnR$, we have $\da1=\SOnR$ for $n$ odd.

When $n$ is even, we can also have $\lambda= -1$. But this implies
$-\Id = x^\T x$. Since  $x^\T x$
is symmetric nonnegative, this is impossible for real matrices.

\noindent (ii)--(iii): Set $x=\pmatrix{ a&b\cr c&d\cr}$.
Using equation (10.1.1), we rewrite the equality $x^{-1}x=\Id$
as
$$\eqalign{
  \lambda &= \sign (a)|a|^{\alpha -1} a 
	    + \sign (b) |b|^{\alpha -1} b 
           = |a|^\alpha + |b|^\alpha\cr
  0       &= \sign (a)|a|^{\alpha -1} c 
	    + \sign (b) |b|^{\alpha -1} d \cr
  0	  &= \sign (c)|c|^{\alpha -1} a
	    + \sign (d) |d|^{\alpha -1} b \cr
  \lambda &= \sign (c)|c|^{\alpha -1} c 
	    + \sign (d) |d|^{\alpha -1} d 
           = |c|^\alpha + |d|^\alpha \, .\cr
  }\eqno{(10.1.2)}
$$
Multiplying the second equality by $ab$, the third by $cd$ and
subtracting yields
$$
  \big( |a|^\alpha - |d|^\alpha \big) bc
  + \big( |b|^\alpha -|c|^\alpha \big) ad = 0 \, .
$$
But $1=\det\, x=ad-bc$ implies then
$$
  \big( |a|^\alpha -|d|^\alpha + |b|^\alpha - |c|^\alpha\big)
  bc + |b|^\alpha -|c|^\alpha = 0\, .
$$
At this stage, the first and last equations in (10.1.2) 
yields $|b|^\alpha
= |c|^\alpha$, i.e., $|b|=|c|$, which then implies $|a|=|d|$. Set
$c=\epsilon_1b$ and $d=\epsilon_2a$ for $\epsilon_1,\epsilon_2$
in $\{\, -1,1\,\}$.

If $bd$ is nonzero, the second equation in (10.1.2) gives, after
multiplication by $|bd|$,
$$
  \sign (a)|d|^\alpha \epsilon_1\sign (b)|b|^2
  + \sign (b)|b|^\alpha \sign (d)|d|^2 = 0\, .
$$
Thus, $\sign (a)\epsilon_1=-\sign (d)$ and $|d|=|b|=|a|=|c|$.
Hence, the matrix is of the form
$$
  x=|a|
  \pmatrix{ \epsilon_1 &\epsilon_3\cr\epsilon_2 &\epsilon_4\cr}
  \qquad\hbox{ with } \epsilon_1,\epsilon_2,\epsilon_3,\epsilon_4
  \in\{\, -1,1\,\} \, .
$$
The condition $1=\det
x=a^2(\epsilon_1\epsilon_4-\epsilon_2\epsilon_3)$ forces
$\epsilon_1\epsilon_4-\epsilon_2\epsilon_3$ to be positive. 
Consequently, $\epsilon_1\epsilon_4=+1$, and $1=2a^2$, i.e.,
$a=\pm1/\sqrt{2}$. Thus,
$$
  x_1={1\over\sqrt{2}}\pmatrix{ \epsilon_1 &-\epsilon_2 \cr
				\epsilon_2 &\hfill\epsilon_1\cr}
$$
solves (10.1.2).

Next, if $bd=0$, let us assume that, say, $d=0$. Equation
(10.1.2) reads
$$\eqalign{
  \lambda & = |a|^\alpha + |b|^\alpha \cr
  0	  & = \sign (a)|a|^{\alpha -1} c\cr
  \lambda & = |c|^\alpha \, .\cr
  }$$
Since $x$ is of determinant $1$, the matrix $x$ is not zero. 
Therefore, the relation
$2\lambda = |a|^\alpha + |b|^\alpha + |c|^\alpha + |d|^\alpha
>0$ forces $c\neq 0$. Hence, $a=0$ and $|b|^\alpha=\lambda$. The
matrix is of the form
$$
  x_2=|c|\pmatrix{ 0 & \epsilon_2 \cr \epsilon_1 & 0\cr } \, .
$$
The condition $\det\, x_2=1$ forces $|c|=1$ and
$\epsilon_1=-\epsilon_2\in\{\, -1,1\,\}$.

Finally, if $b=0$, similar arguments yields a solution
$x_3=\pmatrix{ \epsilon & 0\cr 0 & \epsilon \cr}$, with
$\epsilon$ in $\{\, -1,1\,\}$. 

We then have
$$\eqalign{
  \sum_{1\leq i,j\leq 2}|(x_1)_{i,j}|^\alpha
  & = {4\over 2^{\alpha /2}} = 2^{2-{\alpha\over 2}} \, , \cr
  \sum_{1\leq i,j\leq 2} |(x_k)_{i,j}|^\alpha
  & = 2\, , \qquad k=2,3 \, . \cr
  }$$
If $\alpha <2$, then $2<2^{2-(\alpha /2)}$, while if $\alpha
>2$, we have the reverse inequality $2> 2^{2-(\alpha /2)}$.
The result follows.  \hfill$\qed$

\bigskip

In general, I conjecture the following.

\bigskip

\state{10.1.3. CONJECTURE.} {\it%
  If $\alpha \neq 2$, then 
  $\inf\big\{\, I(x) : \det\, x=1\,\big\}$\
  is achieved at a finite number of matrices. Moreover, the
  difference of the two fundamental forms of 
  $\Lambda_{I(\SLnR)}$ and $\SLnR$ at those matrices is positive.
  Finally, if
  $\alpha <2$, these matrices have unique nonzero elements on
  each row and each column, whose absolute value is $1$; hence, 
  up to signs, they are permutation matrices.}

\bigskip

Some numerical computations support Conjecture 10.1.3. Also,
permutation matrices satisfy equation (10.1.1).

We can now state our approximation of the tail probability for
$\det\, X$. 

\bigskip

\state{10.1.4. THEOREM.} {\it%
  Let $X={(X_{i,j})}_{1\leq i,j\leq n}$ be a random matrix with
  independent and identically distributed coefficients,
  each having a symmetric Weibull like density $w_\alpha$.

  \noindent (i) If $\alpha = 2$ then, as $t$ tends to infinity,
  $$
    P\{\, \det\, X\geq t\,\} 
    \sim {\pi^{(n-1)(n+2)/4}\over (2\pi )^{n^2/2}\sqrt{n}}
    \vol \big( \SOnR\big) e^{-{nt^{2/n}/2}} t^{(n^2-n-2)/2}
    \, .
  $$
  
  \noindent (ii) If $\alpha >2$, under conjecture 10.1.3,
  $$
    P\{\, \det\, X\geq t\,\} \sim c_1e^{-I(A_1)t^{\alpha /n}}
    t^{(\alpha (n^2-1)-2n^2)/2n} \qquad\hbox{ as } \ttoi \, , 
  $$
  where $c_1>0$ is a constant.}

\bigskip

\state{REMARK.} The constant $c_1$ in (ii) can be numerically
computed as will be clear from the proof and Lemma 10.1.1. The
volume of $\SOnR$ in (i) is the volume when $\SOnR$ is viewed as
a submanifold of $\RR^{n^2}$. This volume is the
$n(n-1)/2$-dimensional Hausdorff-Lebesgue measure of the special
orthogonal group.

When plugging $\alpha=2$ in the exponent of $t$ in (ii), we
obtain $-1/n$, which is clearly different from the exponent of
$t$ in (i). Hence, the exponent of $t$ has a discontinuity at
$t=2$.

When $\alpha =2$ and $n=1$, then $\det\, x=x$; the formula reads
$$
  P\{\, X\geq t\,\} \sim {1\over \sqrt {2\pi}}e^{-t^2/2}t^{-1}
  \, \qquad\hbox{ as } \ttoi \, , 
$$
a well known fact! ${\rm SO}(1,R)=\{\, 1\,\}$, and its volume
measure is obtained by putting a Dirac mass at $1$.

\bigskip

\stateit{Proof of Theorem 10.1.4.} The result is an application of
Theorem 7.1.

Let us first determine the exponential term of the asymptotic
equivalent. Since $A_t=t^{1/n}A_1$ and $I$ is
$\alpha$-homogeneous, this term is $e^{-t^{\alpha /n}I(A_1)}$.
The calculation of $I(A_1)$ relies upon Lemma 10.1.2 and
Conjecture 10.1.3. For $\alpha=2$, Lemma 10.1.2 implies that
$I(A_1)=I(\Id )=n/2$. For $n=2$, Lemma 10.1.2 gives
$I(A_1)=2/\alpha$ if $\alpha <2$, while $I(A_1)=4 /(\alpha
2^{\alpha /2})= 2^{2-\alpha/2}/\alpha$ if $\alpha >2$.

In general, $I(A_1)$ can be computed numerically.

To obtain the polynomial term in $t$ in the asymptotic
expansion, Theorem 7.1 requires us to calculate $k=\dim \da1$.
When $\alpha =2$, we obtain $k=\dim \SOnR = n(n-1)/2$. For
$\alpha\neq 2$, we have $k=0$ since $\da1$ is discrete --- 
here we use Conjecture 10.1.3 when $n>2$.

It remains to evaluate the constant $c_1$ in Theorem 7.1
and to verify the assumptions of Theorem
7.1. Since the differential geometries of $\da1$ and $\partial
A_1$ are involved as well as that of $\Lambda_{I(A_1)}$, we need
to calculate the differential and Hessian of $\det$ and $I$.
Lemma 10.0.1 takes care of the former. When dealing with matrices
it is convenient to express $\D I$ and $\D^2I$ on the orthonormal
basis $\Eij$. Thinking of $\Eij$ as an element in the dual of
$\RR^{n^2}\equiv \MnR$, we have $\Eij (M)=\< \Eij ,M\>=M_{i,j}$
for any matrix $M=(M_{i,j})$ in $\MnR$. With this notation, the
gradient and Hessian of $I$ have the following form.

\bigskip

\state{10.1.5. LEMMA.} {\it%
  For any $\alpha\geq 1$ and any $n\times n$ real matrix $x$,
  \vskip 2mm
  \halign{ \qquad\qquad\hfill$#$ & $#$\hfil \cr 
    \D I(x) & = \sum_{1\leq i,j\leq n} \sign (x_{i,j})
	    |x_{i,j}|^{\alpha-1}\Eij \, .  \cr
  \noalign{\vskip 2mm\noindent Moreover, if $\alpha\geq 2$,
           \vskip 2mm}
   \D^2I(x) & = (\alpha -1)\sum_{1\leq i,j\leq n} 
	      |x_{i,j}|^{\alpha -2} \Eij\otimes \Eij\in
	      {\rm M}(n^2,\RR ) \, .\cr
  }
}

\bigskip

\stateit{Proof.} Viewing $\MnR$ as $\RR^{n^2}$, we
have for every $x,h$ in $\MnR$,
$$\displaylines{\qquad
  I(x+\epsilon h) = I(x)+\epsilon \sum_{1\leq i,j\leq n}
	\sign (x_{i,j})|x_{i,j}|^{\alpha -1}h_{i,j}
  \cut
  {} + {\epsilon^2\over 2} \sum_{1\leq i,j\leq n}
	(\alpha -1)|x_{i,j}|^{\alpha-2}h_{i,j}^2+O(\epsilon^2)
  \qquad\cr}
$$
as $\epsilon$ tends to $0$. Since $h_{i,j}=\< \Eij ,h\>$, we obtain
the expression for $\D I$.  Using the polarization formula to express
$\D^2I(x)(h,k)$, we see that
$$
  (h_{i,j}+k_{i,j})^2 - h_{i,j}^2 - k_{i,j}^2
  = 2h_{i,j}k_{i,j}
  = 2\Eij\otimes \Eij (h,k) \, .
$$
This gives the expression for $\D^2I(x)$.\hfill$\qed$

\bigskip

We are equipped to determine the tangent spaces to $\da1$ and
$\Lambda_{I(A_1)}$, from which we will deduce $\det\, G_{A_1}$.

\bigskip

\state{10.1.6. LEMMA.} {\it%
  For any $x$ in $\da1$,
  $$\eqalign{
    T_x\Lambda_{I(A_1)}
      & =\big\{\, h\in\MnR  : \< x^{-1\T},h\> = 0\,\big\} \cr
      & = \{\, x^{-1\T}\,\}^\perp
        = \{\, xh : h\in\MnR \, , \, \trace\, h=0\,\} \, ; \cr
    T_x\da1 
      & = \emptyset\hbox{ if } \alpha\neq 2 \hbox{ (under
	conjecture 10.1.3 for $n\neq 2$).} \, ; \cr
    T_x\da1
      & = \{\, xh : h\hbox{ skewsymmetric}\,\} \hbox{ if }
	\alpha = 2\, .\cr
  }$$
  Consequently, at any $x$ of $\calD_{A_1}\cap \Lambda_{I(A_1)}$,
  \finetune{\hfuzz=10pt} % en fait la formule ne depasse pas
  $$
    T_x\Lambda_{I(A_1)}\ominus T_x\da1 =
    \cases{ \{\, x^{-1\T}\,\}^\perp \quad \hbox{ if $\alpha\neq 2$}
            &\cr
            \noalign{\vskip .05in}
	    \{\, xh : h\in \MnR \, , \, h\hbox{ symmetric }
	    \, , \, \trace (h)=0\,\} &\cr
           \hfill\hbox{ if  $\alpha =2$.}&\cr 
	    \cr}
  $$}
  \finetune{\hfuzz=0pt}

\bigskip

\stateit{Proof.} Since $\Lambda_{I(A_1)}$ is a
level set of $I$, we have for all $x$ in $\calD_{A_1}$,
$$
  T_x\Lambda_{I(A_1)}=\big\{\, \D I(x)\,\big\}^\perp 
  = \{\, x^{-1\T}\,\}^\perp
$$
thanks to Lemmas 10.1.5 and 10.1.1. Since
$$\eqalign{
  \{\, x^{-1\T}\,\}^\perp
  & = \{\, h : \< x^{-1\T},h\> = 0 \,\}
    = \{\, h : \trace (x^{-1}h)=0\,\}  \cr
  & = \{\, xh  : \trace (h)=0\,\}
  \cr}
$$
the expressions for $T_x\Lambda_{I(A_1)}$ follow.

When $\alpha$ is different than $2$, Conjecture 10.1.3 asserts 
that $\da1$ is a finite set, and indeed $T_x\da1=\emptyset$.

For $\alpha$ equal to $2$, the dominating manifold is $\SOnR$. The Lie
algebra of $\SOnR$ is the set of all skewsymmetric matrices
--- see, e.g., Knapp, 1996, \S I.1 --- and the 
expression for $T_x\da1$ follows in this case.

The result on $T_x\Lambda_{I(A_1)}\ominus T_x\da1$ is then clear
since the skewsymmetric matrices are orthogonal to the symmetric
ones.\hfill$\qed$

\bigskip

In order to describe the matrix $G_{A_1}$ involved in Theorem
7.1, recall that the $\ell_p$-norm of a vector
$x\in\RR^{n^2}\equiv \MnR$ is
$$
  |x|_p=\Big( \sum_{1\leq i,j\leq n} |x_{i,j}|^p\Big)^{1/p}
  \, .
$$

\bigskip

\state{10.1.7. LEMMA.} {\it%
  For $x$ in $\calD_{A_1}$, the matrix $G_{A_1}(x)$,
  is obtained in restricting the bilinear form
  $$\displaylines{\qquad
    |x|_{2(\alpha -1)}^{1-\alpha}(\alpha -1)
    \sum_{1\leq i,j\leq n} |x_{i,j}|^{\alpha -2}\Eij\otimes \Eij
    \cut
    {} + |x|_{2(\alpha -1)}^{1-\alpha}
    \lambda^{-1}\sum_{1\leq i,j,k,l\leq n} (x^{-1})_{l,i}
    (x^{-1})_{j,k}\Eij\otimes\Ekl \qquad\cr
  }$$
  to the subspace $T_x\Lambda_{I(A_1)}\ominus T_x\da1$. }

\bigskip

\stateit{Proof.} Given the comment following Theorem 7.1,
it is enough to calculate the second fundamental form of
the hypersurface $\Lambda_{I(A_1)}$ (resp. $\partial A_1$).
Since this hypersurface is the level set of the function $I$
(resp. $\det$), its second fundamental form is the restriction
of $\D^2I/|\D I|$ (resp. $\D^2\det /|\D\det |$) to the tangent space
of $\Lambda_{I(A_1)}$ (resp. $\partial A_1$). Lemma 10.1.5 gives
$$\eqalign{
  {\D^2I(x)\over |\D I(x)|}
  & = { (\alpha -1)\sum_{1\leq i,j\leq n} |x_{i,j}|^{\alpha -2}
      \Eij\otimes\Eij\over 
      \big( \sum_{1\leq i,j\leq n} |x_{i,j}|^{2(\alpha -1)}
      \big)^{1/2} } \cr
  & = { \alpha -1\over |x|_{2(\alpha -1)}^{\alpha -1} }
     \sum_{1\leq i,j\leq n} |x_{i,j}|^{\alpha -2} 
     \Eij\otimes \Eij \, . \cr
  }$$
To calculate the second fundamental form of $\partial A_1$ at a
point $x$ in $\da1$, notice that for $h$ tangent to 
$\partial A_1$, Lemma 10.0.1 implies $\trace (x^{-1}h)=0$. 
Consequently, for $h,k$ in
$T_x\partial A_1$, Lemma 10.0.1 yields $\D^2\det (x)(h,k)=-\trace
(x^{-1}hx^{-1}k)$. Then, Lemma 10.0.2 and (10.0.1) show that 
for $x\in\calD_{A_1}$, the matrix $\D^2\det (x)/|\D\det (x)|$ 
is
$$
  -\sum_{1\leq i,j,l,m\leq n} (x^{-1})_{m,i}(x^{-1})_{j,l}
  E^{i,j}\otimes E^{l,m}\big/ |(x^{-1})^\T |
$$
and the result follows. \hfill$\qed$

\bigskip

\state{10.1.8. LEMMA.} {\it%
  If $\alpha=2$ and $x$ belongs to $\da1 =\SOnR$, then
  $$
    G_{A_1}(x) = {2\over\sqrt n}\Id_{\RR^{(n-1)(n+2)/2}} \, .
  $$}

\bigskip

\stateit{Proof.} For $\alpha =2$ the differential of $I$ is
the identity. If $x$ is in $\da1$, Lemma
10.1.2 forces $|\D I(x)|=|x|=\sqrt{n}$. Furthermore,
since $\calD_{A_1}$ is the special orthogonal group,
$\D\det (x)=(x^{-1})^\T  = x$ on $\calD_{A_1}=\SOnR$. 
Consequently, $|\D\det (x)|=\sqrt{n}$ on $\calD_{A_1}$. 

If $h$, $k$ are in $T_x\Lambda_{I(A_1)}\ominus T_x\da1$, 
and $x$ is in $\SOnR$, Lemma 10.1.6 implies 
$(x^{-1}h)^\T =x^{-1}h$ since $x^{-1}h$ is symmetric as well 
as $\trace (x^{-1}h)=0$. We then infer from Lemma 10.0.1 that
$$\eqalign{
  \D^2\det (x)(h,k)
  & = -\trace (x^{-1}hx^{-1}k)
    = -\trace \big( (x^{-1}h)^\T x^{-1}k\big)
    = -\trace (h^\T k) \cr
  & = -\< h\, ,k\> \, . \cr}
$$
Consequently, the restriction of the bilinear form $\D^2\det$ to
$T_x\Lambda_{I(A_1)}\ominus T_x\da1$ is the identity.
When $x$ is in $\calD_{A_1}$, it follows that
$G_{A_1}(x)(h\, ,k)\allowbreak =2\< h,k\>/\sqrt n$ on
$T_x\Lambda_{I(A_1)}\ominus T_x\da1$. Since
$T_x\Lambda_{I(A_1)}\ominus T_x\da1$ has dimension
$(n-1)(n+2)/2$ thanks to Lemma 10.1.6, the result 
follows.\finetune{\hfill\break\line{\hfill$\qed$}}

\bigskip

In order to obtain an expression for $G_{A_1}$ when $\alpha\neq
2$, we find an explicit orthonormal basis of
$T_x\Lambda_{I(A_1)}\ominus T_x\da1$ for $x$ in $\da1$. 
It is then possible to
express the matrix $G_{A_1}$ in this basis. The construction
goes as follows.

For a real matrix $M={(M_{i,j})}_{1\leq i,j\leq n}$ we denote by
$M_{\sbullet ,j}={(M_{i,j})}_{1\leq i\leq n}$ (resp.
$M_{i,\sbullet}$) the vector in
$\RR^n$ made of its $j$-th column (resp. $i$-th row).

Notice that any $x$ in $\da1$ is also in $\SLnR$, and so is invertible.
For $1\leq i,j\leq n$, $i\neq j$, let $y_j^i\in\RR^n$ be an
orthonormal basis of $\{\, (x^{-1})_{\sbullet ,i}\,\}^\perp$,
the orthogonal subspace in $\RR^n$ of the $i$-th column vector
of $x^{-1}$. Define also 
$e={\big( (x^{-1\T}x^{-1})_{i,i}\big)}_{1\leq i\leq n}\in\RR^n$,
the vector whose coordinates are the diagonal entries of
$x^{-1\T}x^{-1}$. Furthermore, define $y_i^i$, $1\leq i\leq
n-1$ to be an orthonormal basis of $\{\, e\,\}^\perp$ where $\{\,
e\,\}^\perp$ --- in $\RR^n$ --- is equipped with the quadratic form
$\Proj_{\{ e\}^\perp}\diag {\big(( x^{-1\T}x^{-1})_{i,i}\big)}_{1\leq i\leq n}
\big|_{\{\, e\,\}^\perp}$. This quadratic form is the compression to
$\{\, e\,\}^\perp$ of the diagonal matrix obtained by writing the
components of $e$ on its diagonal. We
denote by $y_{k,j}^i$, $1\leq k\leq n$, the components of the
vector $y_j^i$ in $\RR^n$, $1\leq i,j\leq n$. Finally, define
$$\eqalign{
  F^{i,j} & =\sum_{1\leq k\leq n} y_{k,j}^i\Eik\, , \qquad
		i\neq j\, , \, 1\leq i,j\leq n \, , \cr
  F^{i,i} & =\sum_{1\leq l,m\leq n}y_{l,i}^i (x^{-1})_{m,l}
		\Elm \, , \qquad 1\leq i\leq n-1 \, . \cr
  }$$

\bigskip

\state{10.1.9. LEMMA.} {\it%
  In $\MnR\equiv \RR^{n^2}$, the $n^2-1$ vectors $F^{i,j}$, 
  for $i,j$ in $\{\, 1,2,\ldots , n\,\}$ with $(i,j)\neq (n,n)$,
  form an orthonormal basis of $\{\, x^{-1\T}\,\}^\perp$.
  }

\bigskip

\stateit{Proof.} Let us first show that all the
matrices $F^{i,j}$ are orthogonal to $x^{-1\T}$. Indeed, if
$i\neq j$,
$$
  \< F^{i,j},x^{-1\T}\> 
  = \sum_{1\leq k\leq n} y_{k,j}^i (x^{-1})_{k,i}
  = \big\< y_j^i,(x^{-1})_{\sbullet, i}\big\> 
  = 0 \, ,
$$
while for $i=j$,
$$\eqalign{
  \< F^{i,i},x^{-1\T}\>
  & = \sum_{1\leq l,m\leq n} y_{l,i}^i(x^{-1})_{m,l}
    (x^{-1})_{m,l}
    = \sum_{1\leq l\leq n} y_{l,i}^i (x^{-1\T}x^{-1})_{l,l} \cr
  & = \< y_i^i\, ,e\>  = 0 \, . \cr
  }$$
To check that we have an orthonormal basis, we use the identity
$\< E^{i,j},E^{k,l}\> \allowbreak = \delta_{(i,j),(k,l)}$.
If $i\neq j$ and $p\neq q$,
$$
  \< F^{i,j}, F^{p,q}\>
  = \sum_{1\leq k,l\leq n} y_{k,j}^iy_{k,q}^p\delta_{i,p}
  = \< y_j^i\, ,y_q^p\> \delta_{i,p}
  = \delta_{(i,j),(p,q)} \, . 
$$
Next, for $i\neq j$,
$$
  \< F^{i,j},F^{p,p}\> 
  = \sum_{1\leq k\leq n} y_{k,j}^iy_{i,p}^p (x^{-1})_{k,i}
  = y_{i,p}^p \< y_j^i\, ,(x^{-1})_{\sbullet ,i}\>
  = 0 \, .
$$
Finally, for $i,j=1,2,\ldots , n-1$,
$$\eqalign{
  \< F^{i,i},F^{j,j}\>
  & = \sum_{1\leq l,m\leq n} y_{l,i}^iy_{l,j}^j (x^{-1})_{m,l}
    (x^{-1})_{m,l} \cr
  & = \sum_{1\leq l\leq n} y_{l,i}^i y_{l,j}^j 
    (x^{-1\T}x^{-1})_{l,l} \cr
  & = \delta_{i,j} \cr
  }$$
by our choice of $y_i^i$. \hfill$\qed$

\bigskip

Combining Lemmas 10.1.7 and 10.1.9, we can calculate the
$(n^2-1)\times (n^2-1)$-matrix 
${\big( \< G_{A_1}(x)F^{i,j},F^{k,l}\>\big)}$
where $(i,j)$, $(k,l)$ belong to $\{\, 1,\ldots ,n\,\}\times \{\, 1,\ldots
,n\,\}\setminus \{\, (n,n)\,\}$.
This amounts to writing the matrix $G_{A_1}(x)$ in the orthonormal basis
$F^{i,j}$, $1\leq i,j\leq n$, $(i,j)\neq (n,n)$. The explicit
calculation is rather long, and unfortunately does not seem to
simplify much. But the work done so far is all that we need to
implement the approximation numerically.

To conclude the proof of Theorem 10.1.1, it remains to check the
assumptions of Theorem 7.1.

Assumptions (7.1) and (7.2) hold since we assume $\alpha
>1$. Assumption (7.3) is trivial since we assume $\alpha \geq
2$.

Assumption (7.4) is guaranteed by Conjecture 10.1.3 when
$\alpha\neq 2$, while it is trivial for $\alpha =2$.

To check assumption (7.5), use Remark 7.3 and the
calculation of $\D\det$ and $\D I$ made in this section. Indeed,
Lemma 10.1.1 yields
$$
  \< \D I(x),\D\det (x)\>
  = {1\over \lambda}|x|_{2(\alpha -1)}^{2(\alpha -1)} > 0
$$
for any $x$ in $\calD_{A_1}$. This
concludes the proof of Theorem 10.1.4.\hfill$\qed$

\bigskip

When $\alpha =2$, we infer the following corollary.

\bigskip

\state{10.1.10. COROLLARY.} {\it%
  Let ${(X_{i,j})}_{1\leq i,j\leq n}$ be an $n\times n$
  random matrix, with independent and identically coefficients
  all having a standard normal distribution. The distribution
  of  $t^{-1/n}X$
  given $\det X\geq t$ converges weakly* to the uniform
  distribution over $\SOnR$.}

\bigskip

\stateit{Proof.} It is now a straightforward
application of Theorem 7.5.\hfill$\qed$

\bigskip

%\chapternumber={10}
%\chaptername={Random matrices}
\noindent{\fsection 10.2. Random determinants, heavy tails.}
%\chapternumber={10.2}
\section{10.2}{Random determinants, heavy tails}

\bigskip

\noindent We now consider the problem of approximating the tail
probability of the determinant of a random matrix, assuming that its
coefficients are independent and all have a Student like cumulative
distribution function $S_\alpha$.  For this problem, the general
framework proposed so far can be used. The argument is very much like
that used to prove Theorem 8.3.1. For a change, we will give a
probabilistic proof, which is actually inspired by Theorem 5.1,
showing another sort of use of that theorem. This proof will be far
less conceptual, and will give no insights.

The result is as follows.

\bigskip

\state{10.2.1. THEOREM.} {\it%
  Let $X=(X_{i,j})_{1\leq i, j\leq n}$ be a matrix with
  independent and identically distributed coefficients,
  all having a Student like distribution with parameter 
  $\alpha$. Then,
  $$
    P\{\, \det X\geq t\,\}
    \sim {n(2K_{s,\alpha}\alpha^{(\alpha +1)/2})^n\over 2\alpha}
    {(\log t)^{n-1}\over t^\alpha}
    \qquad \hbox{ as } \ttoi \, .
  $$
  }

\bigskip

\stateit{Proof.} We will see why Theorem 5.1 suggests that,
as $t$ tends to infinity,
$$\eqalign{
  P\{\, \det X\geq t\,\}
  & = P\big\{\, \sum_{\sigma\in\frakSigma_n}\epsilon (\sigma )
    \prod_{1\leq i\leq n}X_{i,\sigma (i)}\geq t\,\big\} \cr
  & \sim \sum_{\sigma\in\frakSigma_n}
    P\big\{\, \epsilon (\sigma )\prod_{1\leq i\leq n} 
    X_{i,\sigma (i)}\geq t\,\big\} \cr
  & = n! P\big\{\, \prod_{1\leq i\leq n}X_{1,i}\geq t\,\big\} \,.
    \cr
  }
$$
 
\bigskip

Admitting this relation, our first lemma gives the key estimate.

\bigskip

\state{10.2.2. LEMMA.} {\it%
  Let $X_1,\ldots ,X_n$ be $n$ independent random variables
  with Student-like distribution $S_\alpha$. The product 
  $X_1\ldots X_n$ has upper tail
  $$
    P\{\, X_1\ldots X_n\geq t \,\} \sim
    { (2K_{s,\alpha}\alpha^{\alpha +1\over 2})^n \over
      2\alpha (n-1)!} {(\log t)^{n-1}\over t^\alpha}
    \qquad \hbox{ as } \ttoi \, .
  $$
  Its lower tail is equivalent to its upper tail.
  }

\bigskip

\stateit{Proof.} We proceed by induction. For $n=1$, the result is plain from
the definition of Student-like distributions. Call\ntn{$F_n$} $F_n$ the
cumulative distribution function of the product $X_1\ldots X_n$,
and $c_n$ the constant\ntn{$c_n$}
$$
  c_n 
  = {(2K_{s,\alpha}\alpha^{(\alpha+1)/2})^n\over 2\alpha (n-1)!}
  \, .
$$
Assume that $F_{n-1}$ has the form given in the statement of
the lemma. Then
$$\displaylines{\qquad
  1-F_n(t) = P\{\, X_1\ldots X_{n-1}\geq t/X_n \, ;\, X_n>0\,\}
  \cut
  + P\{\, X_1\ldots X_{n-1}\leq t/n \, ; \, X_n<0\,\} \, .
  \qquad\cr}
$$
Let us evaluate the first probability in the sum. The second
one is either evaluated in the same way, or is obtained from
the first one by changing $X_n$ into $-X_n$. This first probability
can be rewritten as $\int_0^\infty 1-F_{n-1}(t/x) \,\d S_\alpha (x)$.

Let $\delta$ be a positive number. Using the induction
hypothesis, there exists a positive $M$ such that for any $y$
larger than $M$, 
$$
  (1-\delta ) c_{n-1} {(\log y)^{n-2}\over y^\alpha}
  \leq 1-F_{n-1}(y)
  \leq (1+\delta ) c_{n-1} {(\log y)^{n-2}\over y^\alpha} \, .
$$
Moreover, taking $M$ large enough, we also have
$$
  (1-\delta ) {c_1\over y^\alpha}
  \leq 1-S_\alpha (y)
  \leq (1+\delta ) {c_1\over y^\alpha} \, .
$$
Consequently,
$$
  \int_M^{t/M} 1-F_{n-1}(t/x) \,\d S_\alpha (x)
  \leq (1+\delta ) c_{n-1} \int_M^{t/M} 
  {(\log t/x)^{n-2}\over (t/x)^\alpha} \,\d S_\alpha (x) \, .
  \eqno{(10.2.1)}
$$
We integrate by parts, writing
$$\displaylines{
    \int_M^{t/M} x^\alpha (\log t/x)^{n-2}\, \d S_\alpha (x)
    = \Big[ x^\alpha (\log t/x)^{n-2} \big( S_\alpha (x)-1\big)
      \Big]_M^{t/M} 
  \cut
    -\int_M^{t/M} \big( \alpha x^{\alpha-1}(\log t/x)^{n-2}
      - (n-2)x^{\alpha-1}(\log t/x)^{n-3}\big) 
    \big( S_\alpha (x)-1\big)\, \d x \, .
  \cr}
$$
The number $M$ can be taken large enough so that $1/\log z\leq
\delta$ for any $z$ greater than $M$. Then,
$$\displaylines{\quad
    \int_M^{t/M} x^\alpha (\log t/x)^{n-2} \d S_\alpha (x)
  \cut
    \eqalign{ 
      {}\leq{}& O(\log t)^{n-2} + \big(\alpha (1+\delta) c_1+\delta (n-2)\big)
                \int_M^{t/M} {1\over x} 
                \Big( \log {t\over x}\Big)^{n-2} \d x \cr
      {}={}   & O(\log t)^{n-2} + \big(\alpha c_1(1+\delta) +\delta (n-2)\big)
                \int_M^{t/M} {1\over y} (\log y)^{n-2}\d y \cr
      {}={}   & O(\log t)^{n-2} + \big( \alpha c_1(1+\delta)+\delta (n-2)\big)
                {(\log t)^{n-1}\over n-1} \big( 1+ o(1)\big)\cr 
    }%
  \cr}
$$
as $t$ tends to infinity. Therefore, (10.3.1) yields
$$\displaylines{\quad
  \int_M^{t/M}1-F_{n-1}(t/x)\,\d S_\alpha (x)
  \cut
  \leq {(1+\delta )\over t^\alpha} c_{n-1}
  \big( \alpha (1+\delta) c_1+\delta (n-2)\big){(\log t)^{n-1}\over n-1}
  \big( 1+o(1)\big) \, .
  \quad\cr}
$$
We obtain a similar lower bound, replacing $\delta$ by
$-\delta$. 

In the range of integration $x>t/M$, we have
$$
  \int_{t/M}^\infty 1-F_{n-1}(t/x) \,\d S_\alpha (x)
  \leq 1-S_\alpha (t/M) = O(t^{-\alpha}) \qquad
  \hbox{ as } \ttoi \, .
$$
On the other hand, when $x<M$ we have
$$
  \int_0^M \big( 1-F_{n-1}(t/x)\big) \, \d S_\alpha (x)
  \leq 1- F_{n-1} (t/M) 
  = O\Big( {(\log t)^{n-2}\over t^\alpha}\Big)
$$
as $t$ tends to infinity. Since $\delta$ is arbitrary, we proved that
$$
  P\{\, X_1\ldots X_{n-1} \geq t/X_n \, ;\, X_n\geq 0\,\}
  \sim \alpha {c_1c_{n-1}\over n-1} 
  {(\log t)^{n-1}\over t^\alpha} \, ,
$$
as $t$ tends to infinity. Therefore, 
$$
  1-F_n(t)\sim 2\alpha {c_1c_{n-1}\over n-1}
  {(\log t)^{n-1}\over t^\alpha} \, .
$$
Since $2\alpha c_1c_{n-1}=c_n$, the result on the upper tail
follows. The lower tail $F(t)$ as $t$ tends to infinity is
handled in the same way.\hfill$\qed$

\bigskip

The next lemma will allow us to prove that if $\det X\geq t$ and
$t$ is large, it is very unlikely that two different products
$\prod_{1\leq i\leq n} X_{i,\sigma (i)}$,
$\sigma\in\frakSigma_n$, are both of order $t$.

\bigskip

\state{10.2.3. LEMMA.} {\it%
  Let $X_1,\ldots , X_{n+k}$ be $n+k$ independent random
  variables with Student-like distribution function $S_\alpha$.
  For $\alpha >1$, and $k$ positive,
  $$
    P\{\, X_1\ldots X_n\geq t \, ; \, X_{k+1}\ldots X_{k+n}
   \geq t\,\} = o\Big( {(\log t)^{n-1}\over t^\alpha}\Big)
  $$
  as $t$ tends to infinity.
  }

\bigskip

\stateit{Proof.} If $k\geq n$, the result follows from
independence and Lemma 10.2.2. Thus, from now on, assume that
$k<n$. Set $Y=X_1\ldots X_k$, $Z=X_{k+1}\ldots X_n$ and
$U=X_{n+1}\ldots X_{n+k}$. Since these random variables are
independent, we have
$$\displaylines{\qquad
  P\{\, YZ\geq t ; ZU\geq t\,\}
  = \int_{y,u\geq 0}P\Big\{\, Z\geq t \Big( {1\over y}\vee
    {1\over u}\Big)\,\Big\} \, \d F_k(y)\d F_k(u)
  \cut
  +\int_{y,u\leq 0} P\Big\{\, Z\leq t \Big( {1\over y}\wedge
    {1\over u}\Big)\,\Big\} \, \d F_k(y)\d F_k(u) \, .
  \qquad\cr}
$$
Let us evaluate the first integral, the second one being
similar. Using the symmetry in $u$ and $y$, it suffices to prove
that
$$
  \int_{0\leq y\leq u} P\{\, X\geq t/y\,\}\, \d F_k (y)\d F_k(u)
  = o\big( t^{-\alpha}(\log t)^{n-1}\big)
$$
as $t$ tends to infinity.

Let us use the notation $c_k$ as in the proof of Lemma 10.2.2.
Let $\delta$ be an arbitrary positive number. Then, there
exists a positive $M$ such that for any $u>M$
$$\eqalign{
  (1-\delta ) c_{n-k} {(\log u)^{n-k-1}\over u^\alpha }
  & \leq P\{\, Z\geq u\,\} \leq (1+\delta )c_{n-k}
    {(\log u)^{n-k-1}\over u^\alpha} \cr
  (1-\delta )c_k{(\log u)^{k-1}\over u^\alpha}
  & \leq 1-F_k (u)\leq (1+\delta ) c_k 
    {(\log u)^{k-1}\over u^\alpha } \, . \cr
  }
$$
We then have, using Lemma 10.2.2 and the fact that $k$
is strictly less than $n$
$$\displaylines{\qquad
  \int_{\scriptstyle 0<y<u\atop\scriptstyle t/M<u}
  P\{\, Z\geq t/y\,\}\, \d F_k(y)\, \d F_k(u)
  \cut
  \eqalign{
    \leq{}& \int_{t<M<u} \d F_k(u) = O\big(t^{-\alpha}(\log t)^{k+1}\big) \cr
    ={}   &  o\big( t^{-\alpha}(\log t)^{n-1)}\big) \, . \cr}
  \qquad}
$$
Thus, we need to prove that
$$
  \int_{M<y<u<t/M} P\{\, Z\geq t/y\,\} \, \d F_k(y)\d F_k(u)
  = o\big(t^{-\alpha}(\log t)^{n-1}\big)
  \eqno{(10.2.2)}
$$
as $t$ tends to infinity. We first perform the integration in $u$,
obtaining
$$\displaylines{\qquad
  \int_{M<y<u<t/M} P\{\, Z\geq t/y\,\} \, \d F_k(u)\, \d F_k(y)
  \cut
  \eqalign{
   ={}  & \int_{M<y<t/M} P\{\, Z\geq t/y\,\} \big( 1-F_k(y)-1+F_k(t/M)\big)
          \, \d F_k(y) \cr
  \leq{}& \big( 1-F_k(M)\big) \int_{M<y<t/M} P\{\, Z\geq t/y\,\} 
          \, d\, F_k(y)\, . \cr
  }\qquad
  \cr}
$$
We then use the bound on the tail of $Z$ and integrate by parts,
$$\displaylines{\qquad
  \int_{M<y<t/M} P\{\, Z\geq t/y\,\} \, \d F_k(y)
  \cut
  \eqalign{
    \leq{}& (1+\delta ) c_{n-k} \int_{M<y<t/M} \Big( \log {t\over y}\Big)^{n-k-1}
            \Big({y^\alpha\over t^\alpha} \, \d F_k(y) \cr
    \leq{}& (1+\delta ) c_{n-k} \Big[ \Big( \log {t\over y}\Big)^{n-k-1}
            \Big({y^\alpha\over t^\alpha} \big( F_k(y)-1\big)\Big]_M^{t/M} \cr
          & \qquad {}+ (1+\delta ) c_{n-k}\int_{M<y<t/M} 
            {y^{\alpha-1}\over t^\alpha} \Big(\log {t\over y}\Big)^{n-k-2} \times\cr
          & \qquad\qquad
            \Big( \alpha\log {t\over y}+(n-k-1)\Big) \big( F_k(y)-1\big) \, \d y\cr}
  \qquad\cr}
$$
Using the bound on $1-F_k$, we obtain
$$\displaylines{\qquad
  \int_{M<y<t/M} P\{\, Z\geq t/y\,\} \, \d F_k(y)
  \cut
  \eqalign{
   \leq{} & O\big( (\log t)^{n-k-1}/t^\alpha\big)  \cr
          & {} + O(1)\int_{M<y<t/M} {(\log t/y)^{n-k-1}\over t^\alpha} 
            {(\log y)^{k-1}\over y}\, \d y \cr}
  \qquad}
$$
The change of variable $v=(\log y)/\log t$ shows that
$$\displaylines{\qquad
  \int_{M<y<t/M} {(\log t -\log y)^{n-k-1}\over t^\alpha}
  {(\log y)^{k-1}\over y} \, \d y 
  \cut
  \eqalign{
  \sim{} & {(\log t)^{n-2}\over t^\alpha} \int_0^1 (1-v)^{n-k-1} v^{k-1} 
         \, \d v \cr
     ={} & o\Big( {(\log t)^{n-1}\over t^\alpha}\Big) \, .\cr}
  \qquad\cr} 
$$
Consequently, (10.2.2) holds as well as Lemma 10.2.3.\hfill$\qed$

\bigskip

We can now prove Theorem 10.2.1. For a permutation $\sigma$ 
in $\frakSigma_n$, define
$$
  Y_\sigma=\epsilon (\sigma )\prod_{1\leq i\leq n} 
  X_{i,\sigma (i)} \, .
$$
For any fixed positive $\delta$,
$$\displaylines{\quad
  P\{\, \det X\geq t\,\}
  = P\big\{\,\sum_{\sigma\in\frakSigma_n} Y_\sigma\geq t\,\big\}
  \cut
  \eqalign{
  & \geq P\Big\{\, \bigcup_{\sigma\in\frakSigma_n}
    \Big( \{\, Y_\sigma\geq t(1+\delta)\,\} \;
          {\textstyle\bigcap}
          \bigcap_{\scriptstyle\tau\in\frakSigma_n\setminus\{
                   \sigma\}}
   \{\, |Y_\tau|\leq t\delta/n!\,\}\Big)\,\Big\} \cr
  & \geq \sum_{\sigma\in\frakSigma_n}
    P\Big\{\, Y_\sigma\geq t(1+\delta )\,\}\;
          {\textstyle\bigcap}
          \bigcap_{\scriptstyle\tau\in\frakSigma_n\setminus\{
                   \sigma\}} \{\, |Y_\tau|\leq t\delta/n!\,\}
     \,\Big\} \cr
  & \qquad - \sum_{\scriptstyle\sigma_1,\sigma_2\in\frakSigma_n
                   \atop\sigma_1\neq\sigma_2}
    P\{\, Y_{\sigma_1}\geq t(1+\delta ) \, ;\, 
        Y_{\sigma_2}\geq t(1+\delta ) \,\} \cr
  }
  \quad\cr}
$$
From Lemma 10.2.3, we infer that
$$
  \sum_{\scriptstyle\sigma_1,\sigma_2\in\frakSigma_n
                   \atop\sigma_1\neq\sigma_2}
   P\{\, Y_{\sigma_1}\geq t(1+\delta ) \, ;\, 
        Y_{\sigma_2}\geq t(1+\delta ) \,\}
  = o\Big( {(\log t)^{n-1}\over t^\alpha}\Big) \, .
$$
Moreover, if $\tau$ and $\sigma$ are distinct, Lemma 10.2.3 
implies
$$
  P\{\, Y_\sigma \geq t(1+\delta )\hbox{ and }
  |Y_\tau|\geq t/\delta n!\,\} 
  = o\Big( {(\log t)^{n-1}\over t^\alpha}\Big) \, .
$$
Consequently, using Lemma 10.2.2,
$$\displaylines{\qquad
   \sum_{\sigma\in\frakSigma_n}
   P\Big\{\, Y_\sigma\geq t(1+\delta )\,\}\;
          {\textstyle\bigcap}
          \bigcap_{\scriptstyle\tau\in\frakSigma_n\setminus\{
                   \sigma\}} 
         \{\, |Y_\tau|\leq t\delta/n!\,\}
     \,\Big\} 
  \cut
  \eqalign{
  \sim\, & \, n!P\{\, Y_\Id \geq t(1+\delta )\,\} \cr
  \sim\, & \, {n\over 2\alpha}
           (2K_{s,\alpha}\alpha^{\alpha+1\over 2})^n
           {(\log t)^{n-1}\over \big((1+\delta)t\big)^{\alpha+1}}
           \qquad \hbox{ as } \ttoi \, .\cr
  }
  \qquad\cr}
$$
This proves the lower bound
$$
  P\{\, \det\, X\geq t\,\} \geq {n\over 2\alpha} 
  (2K_{s,\alpha}\alpha^{\alpha+1\over 2})^n
  {(\log t)^{n-1}\over t^{\alpha+1}}
  {1+o(1)\over (1+\delta )^{\alpha+1}}
$$
as $t$ tends to infinity.

To obtain a matching upper bound, notice that
$$\displaylines{\quad
  P\{\, \det\, X\geq t\,\} 
  \leq
  P\{\, \exists\, \sigma\in\frakSigma_n \, , \, 
    Y_\sigma\geq t(1-\delta)
  \cut \qquad\quad\hbox{ and } \forall\tau\in\frakSigma_n
    \setminus\{ \sigma\}\, , \, |Y_\sigma|\leq
   t\delta/n! \,\}
%  \hfill\cr\noalign{\vskip 0.05in}\hfill
  \quad\cr\hfill
    + P\{\, \exists\, \tau_1,\tau_2\in\frakSigma_n \, , \, 
    \tau_1\neq\tau_2 \, , \, Y_{\tau_1}\geq \delta t/n! \, ;\,
    Y_{\tau_2}\geq \delta t/n! \,\} 
  \quad\cr}
$$
Applying Lemma 10.2.3 and 10.2.2, we obtain
$$\eqalign{
  P\{\, \det X\geq t\,\}
  & \leq \sum_{\sigma\in\frakSigma_n} 
    P\{\, Y_\sigma\geq t(1-\delta )\,\} 
    +  o\Big( {(\log t)^{n-1}\over t^\alpha}\Big) \cr
  & \sim n{(2K_{s,\alpha}\alpha^{\alpha+1\over 2})^n\over
    2\alpha} {(\log t)^{n-1}\over t^\alpha}
    {1+o(1)\over (1-\delta )^\alpha} \qquad\hbox{ as } 
    \ttoi \, . \cr
  }
$$
Since $\delta$ is arbitrarily small, we proved Theorem
10.2.1.\hfill$\qed$

\bigskip

Let us now show why Theorem 5.1 suggested the
proof of Theorem 10.2.1. Define the set
$$
  A_t=\{\, x\in\MnR \, : \, \det\, x\geq t\,\}
  = t^{1/n} A_1 \, .
$$
Theorem 10.2.1 provides an estimate for the integral
$$
  \int_{A_t} \prod_{1\leq i,j\leq n} s_\alpha (x_{i,j}) 
  \,\d x_{i,j} \, .
$$
The change of variable $Y=\pisa (X)$ leads us to introduce
$$
  B_t=\pisa (A_t) \, .
$$
It allows us to rewrite the integral under consideration as
$$
  \int_{B_t}e^{-I(y)} \,\d y \, ,
$$
where
$$
  I(y)={|y|^2\over 2} + \log (2\pi )^{n^2}
$$
is convex. To minimize $I$ over $B_t$, take a matrix
$x$ in $\partial A_1$ that is in $\SLnR$. 
Then $y=\pisa (t^{1/n}x)$ is on the boundary of $B_t$. 
Furthermore,
$$
  I\big(\pisa (t^{1/n}x)\big) \sim 
  (\log t^{1/n})\sharp\{\, (i,j)\, :\, x_{i,j}\neq 0\,\} \, .
$$
Thus, for $I\big(\pisa (t^{1/n}x)\big)$ to be minimum
asymptotically, $x$ should have as many zero components as
possible, namely $n$. The matrices of $\SLnR$ with $n$
nonvanishing entries form a subgroup which can be described as
follows. Define the matrix
$$
  I_{1,n-1} =\pmatrix{ -1 & 0 \cr
                        0 & \Id_{n-1} \cr} \, .
$$
Let $\DSLnR$ be the subgroup of all diagonal matrices
in $\SLnR$. To a permutation $\sigma$ in $\frakSigma_n$ we
associate the matrix of its permutation representation,
conveniently denoted $\sigma$ as well. Thus, $\sigma
e_i=e_{\sigma (i)}$. Denote by $\frakSigma_{n,+}$ the subgroup of all
even permutation of $n$ elements. Equivalently,
$\frakSigma_{n,+}$ is $\frakSigma_n\bigcap \SLnR$. 
Denote $\frakSigma_{n,-}$ the subset of $\frakSigma_n$ of
all odd permutation matrices. Let $\< I_{1,n-1}\frakSigma_{n,-}\>$
be the subgroup of $\SLnR$ generated by the matrices $I_{1,n-1}\sigma$,
with $\sigma\in\frakSigma_{n,-}$. Then 
$\frakSigma_{n,+}\cup\< I_{1,n-1}\frakSigma_{n,-}\>$ is a group
made of matrices which are, up to the sign of their entries, permutation
matrices, and are of determinant equal to $1$.
This group acts on $\DSLnR$ by 
$$
  (\sigma ,m)\in\big(\frakSigma_{n,+}\cup\< I_{1,n-1}\frakSigma_{n,-}\>\big)
  \times \DSLnR \mapsto \sigma m\in\SLnR \, 
$$ 
Denote by $\big(\frakSigma_{n,+}\cup\< I_{1,n-1}\frakSigma_{n,-}\>\big)\DSLnR$ 
the image of this action. One easily
sees that it is a subgroup of $\SLnR$, made of all the matrices
with exactly $n$ nonvanishing entries. Let $x=\sigma m$ be in
this subgroup. Using Lemma A.1.5 and the fact that $\det\, m =
\prod_{1\leq i\leq d} m_{i,i}=1$,
$$\eqalignno{
  I\big( \pisa (t^{1/\alpha}x)\big)
  & = I\big( \pisa (t^{1/d}m)\big) \cr
  & = \sum_{1\leq i\leq d}\big( \alpha\log (t^{1/n}|m_{i,i}|\big)
    - {1\over 2}\log\log (t^{1/n}|m_{i,i}|) \cr
  & \qquad -2\log (K_{s,\alpha}\alpha^{\alpha/2}2\sqrt\pi )
    \big) + \log (2\pi )^{n^2/2} + o(1) \cr
  & = \alpha\log t -{d\over 2} \log\log t^{1/d}
    - 2d\log (K_{s,\alpha}\alpha^{\alpha/2}2\sqrt \pi ) \cr
  & \qquad + \log (2\pi )^{n^2/2} + o(1) \, . 
  &(10.2.3)\cr
  }
$$
This expression does not depends on $m$. It suggests that
the dominating manifold in our problem should be
$$
  \pisa \Big(t^{1/n}\big(\frakSigma_{n,+}\cup\< I_{1,n-1}\frakSigma_{n,-}\>\big)
             \DSLnR\Big) \, .
$$ 
This set is made up of $n!$ connected components, each component being
$\DSLnR$ composed on the left either by an even permutation, or by
$I_{1,n-1}$ and an odd permutation. As $I$ is invariant under permutations
and composition by $I_{1,n-1}$, all these components should be equally
likely. Since the distribution of the $X_i$'s is asymptotically
symmetric, this suggests the approximation
$$
    P\big\{\, \sum_{\sigma\in\frakSigma_n}\epsilon (\sigma )
    \prod_{1\leq i\leq n}X_{i,\sigma (i)}\geq t \,\big\}
    \sim \sum_{\sigma\in\frakSigma_n}P\big\{\, 
    \prod_{1\leq i\leq n}
    X_{i,\sigma\circ\tau (i)} \geq t\,\big\} \, , 
$$
where $\tau$ is a transposition, depending on $\sigma$, such
that $\frakSigma_n=\{\, \sigma, \sigma\circ\tau \, : \,
\sigma\in \frakSigma_{n,+}\,\}$. The main reason the proof is
complicated using this method is that (10.2.3) is not uniform in $m$. It
is uniform in the range $t^{1/n}m\to\infty$ and $\log
|m_{i,i}|/\log t\to 0$. This is exactly the same problem as the
one we faced in section 8.3, and a similar parameterization can be used.

\bigskip

\setbox1=\vbox{\hbox{\includegraphics{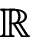}}}

%\noindent{\fsection 10.3. Geometry of the unit ball of 
%  M(n,{\twelvemsbm R}).}
\noindent{\fsection 10.3. Geometry of the unit ball of 
  M({\fsectionmath n},\kern .02in\box1\kern .12in ).}
%\chapternumber={10.3}
\section{10.3}{Geometry of the unit ball of $M(n,\RR)$}
\bigskip

\noindent The purpose of this section is to study some elementary
differential geometric properties of the set $\calS$ of all
real matrices of norm $1$. This will be instrumental in the next
section to obtain results on norm of random matrices.
We will prove --- Propositions
10.3.1 and 10.3.2 --- that this set is a fiber bundle over a Klein
bottle of dimension $2(n-1)$, whose fibers are isomorphic
to the unit ball of $(n-1)\times (n-1)$ real matrices.
We will explicitly calculate various curvatures of this
set.

\bigskip

Recall that the set $\MnR =\RR^{n^2}$ of all $n\times n$ 
matrices with real entries is equipped with the inner
product\ntn{$\< \cdot ,\cdot \>$}
$$
  \langle M,N\rangle = \sum_{1\leq i,j\leq n}M_{i,j}N_{i,j}
  = \trace (MN^\T ) \, .
$$
With this inner product, $\MnR$ is $\RR^{n^2}$ equipped with its
standard inner product. On $\RR^{n^2}$, the Euclidean unit sphere is
a submanifold of dimension $n^2-1$ with constant curvature,
whose geometry is very well understood. However, for
algebraic reasons, it is often more convenient to
equip $\MnR$ with the operator norm\ntn{$\|\cdot \|$}
$$
  \| M\|=\sup\{\, |Mu|\,:\, u\in \Sn1 \,\} \, ,
$$
where $|\cdot |$ is the Euclidean norm in $\RR^n$ and 
$$
 \Sn1 =\{\, x\in\RR^n \,:\, |x|=1\,\}
$$
is the unit sphere centered at the origin. The unit sphere
centered in $(\MnR,\|\cdot \|)$, namely\ntn{$\calS$}
$$
  \calS=\{\, M\in \MnR \,:\, \|M\|=1\,\} \,.
$$
is not as familiar as $S_{n^2-1}$ as far as its
geometry is concerned. We need to understand
what $\calS$ looks like in $\MnR$ identified with $\RR^{n^2}$.

For this purpose, for any $u,v\in\RR^n$, define the 
subspace of matrices\ntn{$H_{u,v}$}
$$
  \Huv=\{\, h\in\MnR \,:\, hu=h^\T v=0\,\} \, .
$$
In what follows, vectors in $\RR^n$ are considered as row vectors,
and so if $u$ belongs to $\RR^n$, then $u^\T $ is a $1\times n$ matrix. We
also use systematically the tensor product notation; if $u$, $v$
are two vectors in $\RR^n$, their tensor product is the matrix
$u\otimes v=vu^\T$. This notation agrees with that used in
section 10.1 when we dealt with vectors in $\RR^{n^2}$.

To understand the geometry of $\calS$, it is convenient
to remove some singular points and define\ntn{$\calS^0$}
$$
  \calS^0 =\{\, M\in\calS \,:\, 1 
  \hbox{ is a simple eigenvalue of } M^\T M\,\} \, .
$$
In $\MnR$, the closure of $\calS^o$ is $\calS$. Proposition 10.3.1
bellow asserts that $\calS^o$
is a smooth submanifold of $\RR^{n^2}$. Moreover, $\calS^0$
is a fiber bundle over
a Klein bottle $\Sn1\otimes\Sn1\equiv S_{n-1}\times S_{n-1}/\{\,\Id,-\Id\,\}$, 
whose fibers are
isomorphic to the unit ball of $H_{e_1,e_1}$ for the operator
norm. So the dimension of the fibers is $(n-1)^2$. We will show
that there are no higher dimensional convex subsets in $\calS^o$
--- this follows from the form of the curvature tensor of $\calS$
given in Theorem 10.3.3. Each fiber is also
is orthogonal to its base point in $\Sn1\times\Sn1/\{\,\Id,-\Id\,\}$.

\bigskip

\state{10.3.1. PROPOSITION} {\it%
  Every matrix $M$ in $\calS$ can be written as $M=u\otimes v+h$
  for some $u,v$ in $\Sn1$, and $h\in \Huv$ with $\|h\|\leq 1$. 
  This decomposition satisfies the following properties:

  \noindent (i) up to the transformation $(u,v)\mapsto (-u,-v)$,
  it is unique if and only if $1$ is a simple eigenvalue of
  $M^\T M$.
  
  \noindent (ii) $\Huv$ is orthogonal to $u\otimes v$ and
  $\dim\Huv =(n-1)^2$ for all $u,v$ in $S_{n-1}$.
  }

\bigskip

\stateit{Proof.} To check that matrices of the form
$M=u\otimes v+h$, with $u,v$ in $\Sn1$ and $\|h\|\leq 1$ 
are of unit norm, notice that the operator norm of such
matrix is at least $1$, since $Mu=v$. On the
other hand, write any vector $x$ of $\RR^n$ as 
$u\langle x,u\rangle +\Proj_{u^\perp}x$ where 
$\Proj_{u^\perp}$ is the projection
onto $\{\,u\,\}^\perp$. Then, apply $M$ to $x$, use that $h$
is a contraction and belongs to $H_{u,v}$ to obtain 
$|Mx|^2\leq |x|^2$, and so $\|M\|\leq 1$.

To prove that all matrices of norm $1$ are of this form, take
$u$ to be a unit eigenvector of $M^\T M$ with eigenvalue $1$.
This vector $u$ is unique up to its sign if and only if $1$
is a simple eigenvalue. Define $v=Mu$ and $h=M-u\otimes v$.
Since
$$
  1=|u|=|M^\T Mu|\leq |Mu|=|v|\leq |u|=1 \, ,
$$
the vector $v$ also belong to $\Sn1$. One easily checks 
that $h$ belongs to $\Huv$. To see why $h$ is a contraction, notice first
that $hu=0$. Moreover, if $w$ is orthogonal to $u$, then
$|hw|=|Mw|\leq |w|$. The uniqueness statement is then clear.

The orthogonality relation (ii) follows from
$\langle u\otimes v,h\rangle = \trace (vu^\T h^\T ) \allowbreak = 0$,
for $h$ belongs to $\Huv$.

To obtain the dimension of $\Huv$, write $R_u$ as an orthogonal
matrix mapping the first vector of the canonical basis of
$\RR^n$, say $e_1$, to $u$. Then $\Huv=R_v^\T H_{e_1,e_1}R_u$.
Hence, $\dim\Huv=\dim H_{e_1,e_1}$. Since the equations
determining $H_{e_1,e_1}$ are
$$
  h_{1,1}=h_{1,2}=\ldots = h_{1,n} = 0
  \hbox{ and } h_{1,1}=h_{2,1}=\ldots = h_{n,1} = 0 \, ,
$$
we have $\dim\Huv=(n-1)^2$ as claimed.\hfill$\qed$

\bigskip

For the unit sphere $\Sn1$ in $\RR^n$, it is an obvious fact that the
tangent space at any point $u$ is just the subspace orthogonal to $u$
in $\RR^n$. So one may wonder if this property has an analogue for the
unit ball $\calS$. Our next proposition shows that this is somewhat
the case and gives an explicit description of the tangent spaces.
This will be useful in calculating the curvature tensor of $\calS^o$.
It also proves that the fibers $\Huv$ are not only orthogonal to
$u\otimes v$ but also to the tangent space 
$T_{u\otimes v}(\Sn1\otimes\Sn1)$. Hence, they point orthogonally to 
the base.

Notice that if $h$ belongs to $\Huv$, the image ${\rm Im} h=h\RR^n$
is included in $\{\,v\,\}^\perp =T_v\Sn1$, 
while ${\rm Im} h^\T \subset\{\,u\,\}^\perp=T_u\Sn1$. For
$u$, $v$ in $\Sn1$ and $h$ in $\Huv$, consider the following 
subspaces of $\MnR$,\ntn{$H_{u,v,h}^1$}\ntn{$H_{u,v,h}^2$}
$$\eqalign{
  H_{u,v,h}^1 &=\{\, a\otimes v-u\otimes (ha) \,:\, a\in T_u\Sn1\,\}\, , \cr
  \noalign{\vskip 0.05in}
  H_{u,v,h}^2 &=\{\, u\otimes b-(h^\T b)\otimes v \,:\, b\in T_v\Sn1\,\} \ . \cr
  }
$$

\bigskip

\state{10.3.2. PROPOSITION.} {\it
  Let $u$, $v$ be in $\Sn1$ and $h$ be in $\Huv$ with $\|h\|< 1$. 
  Then
  
  \noindent (i) $H_{u,v,h}^1\cap H_{u,v,h}^2=\{\,0\,\}\, ;$

  \noindent (ii) $T_{u\otimes v+h}\calS^0 = H_{u,v}\oplus 
  (H_{u,v,h}^1+H_{u,v,h}^2) \, .$

  \noindent The vector $u\otimes v$ is an outward unit normal
  to $\calS^0$ at all points of the form $u\otimes v+h$, with
  $u,v\in\Sn1$ and $h$ a contraction belonging to $\Huv$.
  }

\bigskip

\stateit{Proof.} It is convenient to notice the 
following trivial identity which will be used repeatedly:
for any $a,b,x,y$ in $\RR^n$,
$$
  \langle a\otimes b,x\otimes y\rangle = \trace (ba^\T xy^\T )
  = \langle a,x\rangle \langle b,y\rangle \, .
$$
(i) Let $a$ be in $T_u\Sn1$, and $b$ be in $T_v\Sn1$. Define
$$\eqalign{
  x&=a\otimes v-u\otimes (ha) \in H_{u,v,h}^1 \, , \cr
  y&=u\otimes b-(h^\T b)\otimes v\in H_{u,v,h}^2 \, . \cr
  }
$$
Since $u$ is orthogonal to $a$ and $v$ to $b$, and $h$ is in
$H_{u,v}$,
$$
  |\langle x,y\rangle | 
  = |\langle ha,b\rangle + \langle a,h^\T b\rangle|
  \leq |ha||b|+|a||h^\T b| \, .
$$
Moreover, for the same reasons,
$$
  |x|^2=|a|^2+|ha|^2\, ,\qquad\hbox{ and}\qquad 
  |y|^2=|b|^2+|h^\T b|^2 \, .
$$
Therefore $|\langle x,y\rangle |< |x| |y|$. Thus, $x$ and $y$
cannot be collinear, and $H_{u,v,h}^1\cap H_{u,v,h}^2=
\{\,0\,\}$.

\noindent (ii) The inclusion of $H_{u,v}$ into 
$T_{u\otimes v+h}\calS^0$ is clear: consider the tangent vector
at $0$ of the curve $s\mapsto u\otimes v+(1+s)h\in\calS^0$.

Next, consider two curves $u(s)$, $v(s)$ in $\Sn1$,
with $u(0)=u$, $v(0)=v$, $u'(0)=a$, $v'(0)=b$.
Let $h(s)$ be a curve in $\MnR$ such that 
$h(s)$ is in $H_{u(s),v(s)}$, $h(0)=h$, and $h'(0)=k$. The tangent
vector at $0$ of the curve $u(s)\otimes v(s)+h(s)$ in $\calS^0$
is $a\otimes v + u\otimes b + k$.
Differentiating the relation $h(s)u(s)=h^\T (s)v(s)=0$ at $s=0$
yields
$$
  ku=-ha \qquad\hbox{ and }\qquad k^\T v=-h^\T b \, .
  \eqno{(10.3.1)}
$$
Taking $b=0$, one sees that $k=-u\otimes (ha)$ satisfies
(10.3.1) and so $a\otimes v-u\otimes(ha)$ is in the
tangent space $T_{u\otimes v+h}\calS^0$.
Hence, $H_{u,v,h}^1$ is a subset of $T_{u\otimes v+h}\calS^0$.

Considering $a=0$ and checking that $k=(-h^\T b)\otimes v$
satisfies (10.3.1) yields the inclusion of $H_{u,v,h}^2$ in
$T_{u\otimes v+h}\calS^0$.

The orthogonality of $H_{u,v,h}^1$ and $H_{u,v}$ comes
from the fact that for $a$ in $T_u\Sn1$ and $h$ in $\Huv$,
$$
  \langle a\otimes v -u\otimes (ha) ,h\rangle
  = \trace (av^\T h-u(ha)^\T h)= 0 \, .
$$
Similarly, one proves that $H_{u,v,h}^2$ is orthogonal to $\Huv$.

As a consequence of Proposition 10.3.1, $\calS^0$ is a manifold of
dimension $n^2-1$ and $\dim\Huv=(n-1)^2$. Thus,
$$
  \dim \big(\Huv\oplus (H_{u,v,h}^1+H_{u,v,h}^2)\big)=\dim\calS^0
$$
and we indeed found the whole tangent space to $\calS^0$
--- and not only a subspace.\hfill$\qed$

\bigskip

We can now construct explicitly an orthonormal basis for the tangent
space in which we will express the second fundamental form of the
immersion $\calS^0\subset\RR^{n^2}$, and hence the curvature tensor of
$\calS^0$.

For this purpose, we denote by $e_1^u,\ldots ,e_{n-1}^u$ an orthonormal
basis of $T_u\Sn1$.  Whenever $h$ belongs to $\Huv$, the vector $u$ is
in the kernel of $h^\T h$. Thus, if $\|h\|<1$, the matrix
$\iphth^{-1/2}$ is well defined, and $\{\, u\,\}^\perp$ is an
invariant subspace for this matrix. Consequently, the
vectors\ntn{$a_i$}
$$
  a_i=\iphth^{-1/2}e_i^u\in T_u\Sn1 \, , \qquad 
  1\leq i\leq n-1\, ,
$$
are all in $T_u\Sn1$. They even span $T_uS_{n-1}$, because so do
the $e_i^u$'s and $\| h\| <1$. The matrices\ntn{$f_i$}
$$
  f_i=a_i\otimes v-u\otimes (ha_i) \, , \, 1\leq i\leq n-1\, ,
$$
form an orthonormal basis of $H_{u,v,h}^1$ since an elementary
calculation shows
$$
  \langle f_i,f_j\rangle 
  =\langle a_i,\iphth a_j\rangle
  =\delta_{i,j} \, .
$$

To construct an orthonormal basis of the orthocomplement
$$
  K_{u,v,h}=(H^1_{u,v,h}+H_{u,v,h}^2)\ominus H^1_{u,v,h} \, ,
$$
notice that for $h$ in $\Huv$ with $\|h\|<1$, the subspace
$\{\,v\,\}^\perp$ is invariant under $\iphht^{1/2}$ and
$\imhht^{-1}$. Thus,\ntn{$b_j$}
$$
  b_j=\imhht^{-1}\iphht^{1/2} e_j^v\in T_v\Sn1 \, .
$$
For any $b$ in $T_v\Sn1$, define\ntn{$b^v$}\ntn{$b^u$}
$$\eqalign{
  b^v&=b-2\sum_{1\leq i\leq n-1}\langle b ,\hskip -.6pt ha_i\rangle ha_i
       \in \{\, v\,\}^\perp \, ,\cr
  b^u&=h^\T b-2\sum_{1\leq i\leq n-1}\langle b ,\hskip -.6pt ha_i\rangle a_i
       \in \{\,u\,\}^\perp \, . \cr
  }
$$
The vectors $b_j^u=(b_j)^u$ and $b_j^v=(b_j)^v$ are then 
defined, and so are the matrices\ntn{$g_j$}
$$
  g_j=u\otimes b_j^v-b_j^u\otimes v\in\MnR \, .
$$
Using the bilinearity of the tensor product, we deduce
that $g_j$ belongs to $H_{u,v,h}^1+H_{u,v,h}^2$ since
$u\otimes b_j-(h^\T b_j)\otimes v$ is in $H_{u,v,h}^2$
while $a_i\otimes v-u\otimes (ha_i)$ is in $H_{u,v,h}^1$.

\bigskip

\state{10.3.3. PROPOSITION.} {\it%
  The matrices
  $f_i$, $g_j$, $1\leq i,j\leq n-1$, form an orthonormal
  basis of $H_{u,v,h}^1+H_{u,v,h}^2$.
  }

\bigskip

\stateit{Proof.} It remains for us to prove that the $g_j$'s are
orthonormal, and that they are orthogonal to the $f_i$'s. Since
$b_j^u$ is orthogonal to $u$ and $b_k^v$ to $v$, 
$$
  \< g_j,g_k\> = \< b_j^v,b_k^v\> + \< b_j^u,b_k^u\> \, .
$$
Using the expression of $b_j^u$, $b_k^u$, $b_j^v$ and $b_k^v$, we
write
$$
  \< g_j,g_k\> = \< b_j,Qb_k\> \, ,
$$
where $Q$ is the matrix
$$\displaylines{\qquad
  Q=hh^\T +\Id  + 4\sum_{1\leq i,l\leq n-1} ha_ia_l^\T h^\T
    \< \iphth a_i,a_l\> 
  \cut
  -8\sum_{1\leq i\leq n-1} ha_ia_i^\T h^\T \, .
  \qquad\cr}
$$
Since $\< \iphth a_i,a_l\>=\delta_{i,l}$, we have
$$
  Q= hh^\T +\Id - 4\sum_{1\leq i\leq n-1} ha_ia_i^\T h^\T \, .
$$
Notice that
$$\eqalign{
  \sum_{1\leq i\leq n-1} a_ia_i^\T
  & = \sum_{1\leq i\leq n-1} \iphth^{-1/2} e_i^ue_i^{u\T}
    \iphth^{-1/2} \cr
  \noalign{\vskip 0.05in}
  & = \iphth^{-1/2} \Proj_{u^\perp} \iphth^{-1/2} \, .
  }
$$
Since the image of $h^\T$ is orthogonal to $u$ and
$\{\, u\,\}^\perp$ is invariant under $\iphth^{-1/2}$, 
we obtain
$$
  Q= hh^\T + \Id - 4h \iphth^{-1} h^\T \, .
$$
This expression simplifies further since
$$\displaylines{\qquad
  \imhht \iphht^{-1}\imhht
  \hfill\cr\noalign{\vskip 0.1in}\hfill
  \eqalign{
    =\, & \imhht \Big(\sum_{k\geq 0} (-1)^k (hh^\T )^k 
                      -\sum_{k\geq 0} (-1)^k (hh^\T )^{k+1}\Big)\cr
    =\, & \sum_{k\geq 0}(-1)^k (hh^\T )^k
          -\sum_{k\geq 0}(-1)^k (hh^\T )^{k+1} \cr
        & \qquad\qquad -\sum_{k\geq 0}(-1)^k (hh^\T )^{k+1}
          + \sum_{k\geq 0}(-1)^k (hh^\T )^{k+2} \cr
    =\, & \Id +hh^\T + 4\sum_{k\geq 1} (-1)^k (hh^\T )^k \cr
    =\, & \Id + hh^\T - 4h \iphth^{-1} h^\T \cr
    =\, & Q \cr
  }
  \qquad\cr}
$$
Consequently, replacing $b_j$ by its definition,
$$\eqalign{
  \< g_j,g_k\>
  & = \< \imhht^{-1} \iphht^{1/2}e_j^v \, ,\, \cr
  & \hskip 1in  Q\imhht^{-1}\iphht^{1/2}e_k^v\> \cr
  & = \< e_j^v,e_k^v\> = \delta_{j,k} \, . \cr
  }
$$
To conclude the proof, we calculate
$$\eqalign{
  \< f_i,g_j\>
  & = \< a_i\otimes v - u\otimes (ha_i) \, ,\,
         u\otimes b_j^v-b_j^u\otimes v\> \cr
  & = -\< a_i,b_j^u\> - \< ha_i,b_j^v\> \, . \cr
  }
$$
Since
$$
  b_j^u + h^\T b_j^v = 2h^\T b_j - 2\sum_{1\leq k\leq n-1}
  \< b_j,ha_k\> \iphth a_k \, ,
$$
we deduce that
$$\eqalign{
  \< f_i,g_j\>
  = -2\< a_i,h^\T b_j\> + 2\sum_{1\leq k\leq n} \< b_j,ha_k\>
    \delta_{i,k}
  = 0
  }
$$
as claimed.\hfill$\qed$

\bigskip

Consider an orthonormal basis $h_k$, $1\leq k\leq (n-1)^2$,
of $\Huv$. Furthermore, define the vectors\ntn{$c_j$}
$$
  c_j=-\big[ \Id-2\iphth^{-1}\big] h^\T b_j\in T_u\Sn1 \, .
$$
Quite remarkably, it is possible to explicitely calculate the
curvature tensor of $\calS^0$ through its second fundamental form.

\bigskip

\state{10.3.4. THEOREM.} {\it%
  For $u$, $v$ in $\Sn1$, for $h$ in $\Huv$ with $\|h\|< 1$, 
  the second fundamental form of $\calS$ at $u\otimes v+h$ in
  the orthogonal basis $f_i,g_j,h_k$, $1\leq i,j\leq n-1$,
  $1\leq k\leq (n-1)^2$, is given by the $(n^2-1)\times
  (n^2-1)$ matrix\ntn{$\Pi$}
  $$\Pi =\bordermatrix{
    & {\scriptstyle n-1} & {\scriptstyle n-1} &%
   {\scriptstyle n^2-2n+1} \cr
   {\scriptstyle n-1} & \langle a_i,a_j\rangle &%
      \langle c_i,a_j\rangle & 0 \cr
   {\scriptstyle n-1} & \langle a_i,c_j\rangle &%
      \langle c_i,c_j\rangle +%
      \langle b_i,\imhht b_j\rangle & 0 \cr
    {\scriptstyle n^2-2n+1} & 0 & 0 & 0\cr}
  $$
  }

\bigskip

\stateit{Proof.} Let $N=u\otimes v$ be the
outward unit normal to $\calS^0$ at $u\otimes v+h$ --- see
Proposition 10.3.2. We will denote by $\nabla$ the covariant
derivative on $\calS^0$; that is, for a vector field $X$ defined
on $\calS^0$, for  $u$ a tangent vector vector field and $p$ a
point on $\calS^0$, 
$$
  \nabla_u X(p) = \Proj_{T_p\calS^0} \D X(p)\cdot u \, .
$$
The calculations made in the proof of 
Proposition 10.3.2 show that $f_i$ is the tangent vector at $0$
of a curve $s\mapsto u(s)\otimes v + h(s)$ with $u(0)=u$,
$h(0)=h$ and $u'(0)=a_i$ and $h'(0)=-u\otimes ha_i$. Consequently
$\nabla_{f_i}N=a_i\otimes v$. Moreover,
$$\eqalign{
  \nabla_{g_j}N & =\nabla_{u\otimes b_j-(h^\T b_j\otimes v)}N
 		  - 2\sum_{1\leq i\leq n-1}
                  \langle b_j , ha_i\rangle 
		  \nabla_{u\otimes (ha_i)-a_i\otimes v} N \cr
		& = u\otimes b_j + 2\sum_{1\leq i\leq n-1}
		  \langle b_j,ha_i\rangle a_i\otimes v \, . \cr
  }
$$
Moreover, since $N$ is constant along $h\in\Huv\mapsto
u\otimes v+h$, $\|h\|< 1$, we have $\nabla_{h_k}N=0$.
A routine calculation gives the first entries of the matrix,
namely
$$
 \langle \nabla_{f_i}N,f_j\rangle = \langle a_i,a_j\rangle
    \, .
$$
Next, we have
$$\eqalign{
  \langle \nabla_{g_i}N,f_j\rangle 
  & = -\< h^\T b_i,a_j\> 
    + 2\sum_{1\leq k\leq n-1} \<h^\T b_i,a_k\> \< a_k,a_j\> \cr
  & = \< a_j , 
         -h^\T b_i+2\sum_{1\leq k\leq n-1} a_ka_k^\T h^\T b_i\>
    \, . \cr
  }
$$
Since the image of $h^\T$ is orthogonal to $u$,
$$
  \sum_{1\leq k\leq n-1} a_ka_k^\T h^\T
  = \iphht^{-1} h^\T \, .
$$
This gives the entries in $\< a_i,c_j\>$. 

Finally, we calculate
$$\eqalign{
  \< \nabla_{g_i}N,g_j\>
  & = \big\< u\otimes b_i+2\sum_{1\leq k\leq n-1} \< b_i,ha_k\>
    a_k\otimes v \, ,\, u\otimes b_j^v-b_j^u\otimes v\big\> \cr
  & = \< b_i,b_j^v\> - 2\sum_{1\leq k\leq n-1} \< b_i,ha_k\>
    \< a_k,b_j^u\> \cr
  & = \< b_i,Qb_j\> \, ,\cr
  }
$$
where the matrix $Q$ is
$$
  Q=\Id - 4\sum_{1\leq l\leq n-1} h a_l a_l^\T h^\T
    + 4\sum_{1\leq k,l\leq n-1} h a_k a_l^\T h^\T \< a_k,a_l\>
  \, .  
$$
Again, since the image of $h^\T$ is orthogonal to $u$,
$$
  h\sum_{1\leq l\leq n-1} a_l a_l^\T h^\T
  = h\iphht^{-1} h^\T\, .
$$
Moreover,
$$\displaylines{\qquad
  h\sum_{1\leq k,l\leq n-1} a_ka_l^\T \< a_k,a_l\> h^\T
  \cut
  \eqalign{
    =\, & h\iphht^{-1/2} \sum_{1\leq k,l\leq n-1} e_ke_l^\T
        \< e_k,\iphth^{-1} e_l\>\times \cr
    \noalign{\vskip -.15in}
        & \hskip 2.2in\iphht^{-1/2} h^\T \cr
    =\, & h \iphht^{-2} h^\T \, . \cr
  }
  \qquad\cr}
$$
Consequently,
$$\eqalign{
  Q & = \Id -4 h\iphht^{-1} h^\T + 4 h \iphht^{-2} h^\T \cr
    & = \Id -hh^\T + h\big( \Id -2 \iphht^{-1}\big)^2 h^\T
      \, . \cr
  }
$$
This gives
$$
  \< \nabla_{g_i}N,g_j\> = \< b_i,\imhht b_j\> + \< c_i,c_j\>
$$
as claimed.\hfill$\qed$ 

\bigskip

It follows from Theorem 10.3.4 and elementary results on
immersions that the Riemannian curvature tensor $R$ of $\calS^0$ 
can be calculated explicitly in the basis $f_i,g_j,h_k$.
It is convenient to define $f_{n-1+i}=g_i$ for
$i=1,\ldots ,n-1$ and $f_{2(n-1)+i}=h_i$ for
$i=1,\ldots ,(n-1)^2$. If $X,Y$ are two elements of
$T_{u\otimes v+h}\calS^0$,
then
$$
  \langle R(f_i,f_j)X,Y\rangle =
  X^\T \big( \nabla_{f_i}N(\nabla_{f_j}N)^\T 
  -\nabla_{f_j}N(\nabla_{f_i}N)^\T \big) Y \, .
$$
Thus, $R(f_i,f_j)$ is the compression to the tangent
space of the matrix $\nabla_{f_i}N(\nabla_{f_j}N)^\T 
-\nabla_{f_j}N(\nabla_{f_i}N)^\T$.

\bigskip

As a byproduct of the work done, 
we can prove that there are no flat and nontrivial 
convex subsets in $\calS^0$ besides the unit
balls of the fibers $\Huv$. This result will not be used in the
sequel, but brings more intuition on the shape of the sphere
$\calS$. It is enough to show
that the $2(n-1)\times 2(n-1)$ upper left corner
submatrix of $\Pi$ has no zero eigenvalue. Since $\| h\|< 1$
on our parameterization of $\calS^0$, this follows from the next
result.

\bigskip

\state{10.3.5. PROPOSITION.} {\it%
  The following equality holds,
  $$\displaylines{\qquad
    \det\pmatrix{ \< a_i,a_j\>_{i,j} & \< a_i,c_j\>_{i,j} \cr
                  \< a_i,c_j\>_{i,j} & \<c_i,c_j\>_{i,j}+
                                       \< b_i,\imhht b_j\>_{i,j} 
                                       \cr}
    \hfill\cr\noalign{\vskip 0.1in}\hfill
    = \det \iphht \, \det \iphth \, \det \imhht \, .
    \qquad\cr}
  $$}
\finetune{
%\bigskip
}
\stateit{Proof.} We first calculate a subdeterminant of the
given one. Going back to the definition of the $a_i$'s and using
that $u$ is an eigenvector of $\iphth$ associated to the eigenvalue
$1$, 
$$\eqalign{
  \det (\< a_i,a_j\>)_{i,j}
  & = \det \big( e_i^u,\iphth e_j^u\>\big)_{1\leq i,j\leq n-1} \cr
  & = \det \iphth  \, . \cr
  }
$$
Furthermore, since $v$ is an eigenvalue of $\iphth^{1/2}$ and
$\imhht^{-1}$ associated with the eigenvalue $1$ --- this can be
seen by series expanding and using the fact that $h^\T v=0$ ---
$$\displaylines{\qquad
  \det \big(\< b_i,\imhht b_j\>\big)_{i,j}
  \cut
  \eqalign{
   ={} & \det \big( \< e_i^v, \iphht^{1/2}\imhht^{-1}\iphht^{1/2}
                       e_j^v\>_{i,i}\big) \cr
   ={} & \det \iphht \,\det \imhht \, . \cr
  }
  \qquad\cr}
$$
To conclude the proof, we use the following claim, with
$d_i=\imhht^{1/2}b_i$.

\medskip

\state{Claim.} {\it%
  Let $a_i$, $c_i$, $d_i$, $1\leq i\leq n-1$, be $3(n-1)$ vectors
  in $\RR^n$. Consider the $n\times (n-1)$-matrices  
  $a=(a_1,\ldots ,a_{n-1})$, $c=(c_1,\ldots ,c_{n-1})$ and
  $d=(d_1,\ldots , d_{n-1})$. If the image of $c$ is contained
  in the image of $a$, then
  $$
    \det\bigg( \pmatrix{a^\T\cr c^\T\cr}\pmatrix{a & c\cr} +
               \pmatrix{0 & 0 \cr 0 & d^\T d \cr}\bigg)
    = \det (a^\T a) \,\det (d^\T d) \, . 
  $$
  }

To prove the claim, let $P$ be an orthogonal matrix and $D$ be a
diagonal one such that $d^\T d = PDP^\T$. Writing $M$ for the matrix
whose determinant we want to calculate, we have
$$\eqalign{
  \det \, M
  & = \det\bigg( \pmatrix{ \Id & 0 \cr 0 & P \cr }
                 \pmatrix{ a^\T \cr c^\T \cr }
                 \pmatrix{ a & c \cr}
                 \pmatrix{ \Id & 0 \cr 0 & P^\T \cr }
                 + \pmatrix{ 0 & 0 \cr 0 & D\cr} \bigg) \cr
  \noalign{\vskip 0.05in}
  & = \det \bigg( V + \pmatrix{ 0 & 0 \cr 0 & D\cr }\bigg) \, .
    \cr }
$$
The proof then goes by induction on the dimension of $D$,
noticing that for any real number $\delta$,
$$
  \det \bigg( V + \pmatrix{ 0 & 0 \cr 0 & \delta \cr} \bigg)
  = \delta\, \det\, \widehat V_{m,m} + \det\, V
$$
where $\widehat V_{m,m}$ is the $(m-1)\times (m-1)$ left upper corner
of $V$ and $m$ the dimension of $V$. Consequently, we just need
to prove that $\det\, V=0$. This is clear since the condition
${\rm Im}c\subset {\rm Im}a$ implies that the rank of the matrix
$\pmatrix{ a & c \cr}$ is the dimension of the image of $a$, and
hence the rank of $V$ is at most $\dim{\rm Im} a$. Consequently,
we have
$$
  \det\bigg( V + \pmatrix{ 0 & 0 \cr 0 & D \cr} \bigg)
  = \det\, (a^\T a) \, \det \, D \, .
$$
This proves the claim and concludes the proof of Proposition
10.3.5.\hfill$\qed$

\bigskip

As a consequence of Proposition 10.3.5, the Gauss-Kronecker curvature
of the nonflat part of $\calS^0$ at $u\otimes v+h$ is given by $\det
\iphth\,\det\iphht\det\imhht$. Other curvatures can be calculated as
well, leading to more or less interesting formulas.

\bigskip

%\chapternumber={10}
%\chaptername={Random matrices}
\noindent{\fsection 10.4. Norms of random matrices.}
%\chapternumber={10.4}
\section{10.4}{Norms of random matrices}

\bigskip

\noindent Let us again consider a random matrix
$X={(X_{i,j})}_{1\leq i,j\leq n}$ with independent and
identically distributed coefficients. Its (operator) norm is
$$
  \| X\|=\sup\big\{\, |Xu|  : |u|=1\, , \, u\in\RR^n\,\big\}
  \, .
$$
In this section, we will obtain estimates for
the tail probability
$P\big\{\, \| X\|\geq t\,\big\}$, assuming that the
$X_{i,j}$'s are either symmetric Weibull or Student
like distributed.

In theory, we just need to apply the results of chapter 9. 
Indeed, $\RR^n$ being reflexive, $\| X\|=\sup\big\{\, \< Xu,v\>  : 
u,v\in S_{n-1}\,\big\}$. Since $\< Xu,v\>_{\RR^n}=\< X,u\otimes
v\>_{\RR^{n^2}}$, we see that $\| X\|$ is the supremum of the
linear form $X$ acting on the submanifold $S_{n-1}\otimes S_{n-1}$
of $\RR^{n^2}$.

However, a direct application of the results of chapter 9 in the
case of light tails is not that easy. We will proceed by
using both chapter 7 and ideas from sections 9.1 and 9.2
as well.

Our first result is for light tails.

\bigskip

\state{10.4.1. THEOREM.} {\it%
  Let $X={(X_{i,j})}_{1\leq i,j\leq n}$ be a random matrix with
  independent and identically distributed coefficients, all
  having the symmetric Weibull-like density\ntn{$w_\alpha$}
  $$
    w_\alpha (x) = 
    {\alpha^{1-(1/\alpha )}\over 2\Gamma (1/\alpha )}
    \exp\Big( {-|x|^\alpha\over\alpha}\Big)\, , 
    \qquad x\in\RR \, .
  $$

  \noindent (i) If $\alpha = 2$, then
  $$
    P\big\{\, \| X\|\geq t\,\big\}
    \sim {\sqrt{2\pi}\over 2^{n-1}\Gamma (n/2)^2}
    e^{-t^2/2}t^{2(n-1)-1}
    \qquad\hbox{ as } \ttoi \, .
  $$

  \noindent (ii) If $\alpha >2$, then
  $$\displaylines{
      P\big\{\, \| X\|\geq t\,\big\}
      \sim \bigg( 
      {\alpha^{1-(1/\alpha )}\over 2\Gamma (1/\alpha)}
      \bigg)^{n^2} \times
      \cut
      { (2\pi )^{(n^2-1)/ 2} 2^{2n-1}\over 
        n^{(2-\alpha)(n^2+1)/2} (\alpha -2)^{n-1}
        (\alpha -1)^{(n-1)^2/2} }
      \, { e^{-t^\alpha/(\alpha n^{\alpha -2})} \over
        t^{{\alpha\over 2} (n^2+1)-n^2} } 
    \cr}
  $$
  as $t$ tends to infinity.
  }

\bigskip

\stateit{Proof.}  Define\ntn{$A_t$}
$$
  A_t=\big\{\, x={(x_{i,j})}_{1\leq i,j\leq n}\in\MnR : 
  \|x \|\geq t\,\big\} = tA_1 \, .
$$
We need to evaluate
$$
  \bigg( {\alpha^{1-(1/\alpha )}\over 2\Gamma (1/\alpha )}
  \bigg)^{n^2}
  \int_{tA_1}\exp\bigg( -{1\over\alpha}\sum_{1\leq i,j\leq n}
  |x_{i,j}|^\alpha \bigg) \prod_{1\leq i,j\leq n} \d x_{i,j} \, .
$$
Defines the $\alpha$-homogeneous function\ntn{$I(\cdot )$}
$$
  I(x) = {|x|_\alpha^\alpha\over\alpha } 
  = {1\over\alpha}\sum_{1\leq i,j\leq n} |x_{i,j}|^\alpha \, .
$$
We can apply Theorem 7.1. The first step is to calculate
$I(A_1)$ and $\da1$. To do this, we need a description of
the boundary
$$
  \partial A_1= \big\{\, x\in\MnR  : \| x\|=1\,\big\} \, ,
$$
that is of the sphere of radius $1$ in the space of matrices
endowed with the operator norm. This is provided by Proposition
10.3.1. Let us simply recall here that the
matrices of norm $1$ coincide with all matrices of the form
$u\otimes v+h$, where $u$, $v$ belong to the sphere $\Sn1$ 
and $h$ is an $n\times n$ matrix of operator norm less than $1$,
satisfying $hu=h^Tv=0$. This allows us to find $\da1$.

\bigskip

\state{10.4.2. LEMMA.} {\it%
  The function $I$ is minimum over $\partial A_1$ exactly at
  matrices of the form $u\otimes v$ with

  \noindent (i) $u,v\in\big\{\, n^{-1/2} (\epsilon_1,
  \ldots , \epsilon_n) : \epsilon_i\in\{\, -1,1\,\}\, ,\, 
  1\leq i\leq n\,\big\}$ if $\alpha >2$,

  \noindent (ii) $u,v\in\Sn1$ if $\alpha=2$,

  \noindent (iii) $u,v\in\big\{\, \epsilon e_i  : 
  \epsilon\in\{\, -1,1\,\}\, ,\, 1\leq i\leq n\,\big\}$, if 
  $\alpha <2$.
  }

\bigskip

\stateit{Proof.} As we already mentioned, $\| x\|$ is 
the supremum of the linear form $x\in\RR^{n^2}$ 
acting on $\Sn1\otimes\Sn1$. Proposition 9.1.7 and 
convexity of $A_1^{\rm c}$ implies $I(\partial
A_1)=I_\sbullet (\Sn1\otimes \Sn1)$. Moreover,
Lemma 9.1.9 asserts that
$I_\sbullet (x) = 1/(\alpha |x|_\beta^\alpha)$ where 
$\alpha^{-1}+\beta^{-1}=1$. We can first
calculate the points in $\Sn1\otimes\Sn1$ which minimize
$I_\sbullet$. This is rather easy since
$$
  \alpha I_\sbullet (u\otimes v) 
  = \Big( \sum_{1\leq i,j\leq n} |v_iu_j|^\beta
    \Big)^{-\alpha/\beta} 
  = |v|_\beta^{-\alpha}|u|_\beta^{-\alpha} \, .
$$
Thus, we need to locate the maxima of $|u|_\beta$ on $\Sn1$.

If $\alpha$ is larger than $2$, then $\beta$ is smaller than
$2$. Therefore,
$$
  \Big( {1\over n} \sum_{1\leq i\leq n}|u_i|^\beta 
  \Big)^{1/\beta}
  \leq \Big( {1\over n} \sum_{1\leq i\leq n} u_i^2\Big)^{1/2}
  = {1\over \sqrt n} \, , 
$$
with equality if and only if $|u_i|=1/\sqrt{n}$ for all
$i=1,2,\ldots ,n$. Consequently, 
$$
  \sup\big\{\, |u|_\beta : u\in\Sn1\,\big\}
  = n^{{1\over\beta}-{1\over 2}}
  \, , \qquad \hbox{ if } \alpha >2 \, .
$$

If $\alpha=2$, then $\beta =2$, and $|u|_\beta=1$ over
all $\Sn1$.

Finally, if $\alpha$ is smaller than $2$, then $\beta$ is
larger than $2$. A unit vector $u$
has all its components $u_i$ between $-1$ and $1$. Therefore,
$$
  |u|_\beta
  = \Big( \sum_{1\leq i\leq n} |u_i|^\beta\Big)^{1/\beta}
  \leq \Big( \sum_{1\leq i\leq n}|u_i|^2\Big)^{1/\beta}
  = 1\, ,
$$
with equality if and only if one --- and only one --- of 
the $|u_i|$'s is $1$.

Consequently,
$$
  \alpha I_\sbullet (u\otimes v)
  \geq \cases{  n^{2-\alpha} & if $\alpha >2$, \cr
	       1	     & if $\alpha\leq 2$,\cr
	      }
$$
with equality for\ntn{$u_*$, $v_*$} $(u,v)=(u_*,v_*)$ with $(u_*,v_*)$ exactly
in the following sets,
$$
  \cases{ |u_{*,i}|=|v_{*,i}|=1/\sqrt{n}
		& $i=1,2,\ldots , n$, if $\alpha >2$, \cr
	\noalign{\vskip 2mm}
	  u_*,v_*\in\Sn1
		& if $\alpha =2$, \cr
	\noalign{\vskip 2mm}
	  u_*,v_*\in\big\{\, \epsilon e_i : \epsilon
	  \in\{\, -1,1\,\} \, , \, 1\leq i\leq n\,\big\}
		& if $\alpha <2$. \cr  }
$$
For $u$ and $v$ in $\Sn1$, set
$$
  H_{u,v}=\big\{\, h\in\MnR  : hu=h^Tv = 0 \,\big\} \, .
$$
For $h$ in $H_{u,v}$, we have 
$$
  \< u\otimes v,u\otimes v+h\>
  =\trace (uv^Tvu^T+uv^Th)=1 \, .
$$
Consequently,
$$\eqalign{
  I_\sbullet (u\otimes v)
  & = \inf\big\{\, I(x) : \< x,u\otimes v\> =1\,\big\} \cr
  & \leq \inf\big\{\, I(u\otimes v+h)  : h\in
    H_{u,v}\, , \, \| h\|\leq 1 \,\big\} \, . \cr
  }$$
The inequality
$$
  I_\sbullet (u_*\otimes v_*)
  \leq \inf\big\{\, I(x) : x\in\partial A_1\,\big\} \, .
$$
follows. Observe that $I_\sbullet (u_*\otimes v_*)=I(u_*\otimes v_*)$
for any $\alpha\geq 1$. Since the function 
$h\in H_{u,v}\mapsto I(u\otimes v+h)$ is
convex, as a restriction of a convex function to a convex set, 
the infimum of $I$ over $\partial A_1$ is
achieved only at points $x=u_*\otimes v_*$. On those points
$I_\sbullet$ and $I$ coincide and this concludes the proof.
\hfill$\qed$

\bigskip

It is interesting to realize that the proof of Lemma 10.4.2
relies on the fact that $I_\sbullet$ and $I$ coincide on the
matrices $u_*\otimes v_*$. Geometrically, the matrices $u\otimes
v+h$, $h\in H_{u,v}$ with $\| h\|\leq 1$ frorms a truncated 
cylinder with base $S_{n-1}\otimes S_{n-1}$. What makes the proof
work is that $S_{n-1}$ is the polar reciprocal of its convex hull;
a very special property of the sphere!

\bigskip

Let us now calculate all the terms that come from applying Theorem
7.1. We will then justify that we can indeed apply this theorem in
verifying that its assumptions hold.

Let us first consider the case $\alpha>2$. From Lemma 10.4.2, we deduce
$$
  I(A_1)=n^{2-\alpha}/\alpha \, .
$$
The rescaled dominating manifold
$$\displaylines{\qquad
    \da1 = \Big\{\, u\otimes v : u,v\hbox{ of the form }
    {1\over \sqrt n} (\epsilon_1,\ldots ,\epsilon_n) : 
  \cut
    \epsilon_i\in\{\, -1,1\,\}\, , \, 1\leq i\leq n\,\Big\} \, .
  \qquad\cr}
$$
is of dimension $k=0$. Its Riemannian volume is the 
counting measure
$$
  \calM_{\da1}
  = {1\over 2} \sum_{1\leq i,j\leq n} 
    \sum_{\epsilon_i,\eta_j\in\{ -1,1\}}
    \delta_{{1\over n}(\epsilon_1,\ldots , \epsilon_n)
    \otimes (\eta_1,\ldots ,\eta_n)}
  = \sum\delta_{\epsilon\otimes\eta/n}
$$
where the last sum is over all distinct matrices 
${(\eta_i\epsilon_j/n)}_{1\leq i,j\leq n}$ with coefficients in
$\{\, -1/n,1/n\,\}$. Note that
$\epsilon\otimes\eta$ and $(-\epsilon)\otimes (-\eta)$ are equal,
thus not distinct.

For $x$ in $\da1$, we have
$$
  |\D I(x)|^2
  = \sum_{1\leq i,j\leq n} |x_{i,j}|^{2(\alpha -1)}
  = {n^2\over n^{2(\alpha -1)}}
  = n^{4-2\alpha} \, .
$$
We then need to calculate the curvature term $\det\, G_{A_1}$, 
and hence the fundamental form 
$\Pi_{\Lambda_{I(A_1)},u_*\otimes v_*}$ and 
$\Pi_{\partial A_1,u_*\otimes v_*}$ for $u_*\otimes v_*$ in
$\da1$.

From Theorem 10.3.4 with $h=0$, we deduce that
$$
  \Pi_{\partial A_1}
  =\pmatrix{ \Id_{2(n-1)} & 0 \cr\noalign{\vskip 1mm}%
		0	  & 0 \cr}
  \in {\rm M}({\rm n^2-1},\RR\, ) \, .
$$
On the other hand, the second fundamental form of
$\Lambda_{I(A_1)}$ at $x$ is the restriction to the tangent
space $T_x\Lambda_{I(A_1)}$ of
$$
  {\D^2I(x)\over |\D I(x)|}
  = \Big( \sum_{1\leq i,j\leq n} |x_{i,j}|^{2(\alpha -1)}
  \Big)^{-1/2}\diag 
  {\big(|x_{i,j}|^{\alpha -2}\big)}_{1\leq i,j\leq n} 
  (\alpha -1) \, .
$$
Thus, if $x$ is in $\da1$, we have $|x_{i,j}|=1/n$ and
$$
  {\D^2I(x)\over |\D I(x)|}
  = {1\over n^{2-\alpha }}{\alpha-1\over n^{\alpha -2}}
    \Id_{n^2}
  = (\alpha -1)\Id_{n^2} \, .
$$
Hence, for $x$ in $\da1$,
$$\eqalign{
  \Pi_{\Lambda_{I(A_1)},x} - \Pi_{\partial A_1,x}
  & = (\alpha -1)\Id_{n^2-1} 
    - \pmatrix{ \Id_{2(n-1)} & 0 \cr\noalign{\vskip 1mm} 0 & 0 \cr} \cr
  \noalign{\vskip 3mm}
  & = \pmatrix{ (\alpha -2)\Id_{2(n-1)} & 0 \cr\noalign{\vskip 1mm}
		0 & (\alpha -1)\Id_{n^2-2n+1} \cr} \, .\cr 
  }$$
Therefore, on $\da1$,
$$
  \det\, G_{A_1}(x) = (\alpha -2)^{2(n-1)}
  (\alpha -1)^{n^2-2n+1} \, .
$$
We obtain the constant $c_1$ in Theorem 7.1,
$$\eqalign{
  c_1 
  & = (2\pi )^{(n^2-1)/2} \Big( {1\over n^{2-\alpha}}
    \Big)^{(n^2+1)/2} 
    {1\over (\alpha -2)^{n-1} (\alpha -1)^{(n-1)^2/2}}
    \sharp \da1 \cr
  & = { (2\pi )^{(n^2-1)/2}2^{2n-1} \over
	n^{(2-\alpha )(n^2+1)/2}(\alpha -2)^{n-1}
	(\alpha -1)^{(n-1)^2/2} } \cr
  }$$
Putting all the pieces together,
$$\displaylines{\quad
  P(A_t)
  \sim 
  \Big( {\alpha^{1-(1/\alpha )}\over 2\Gamma(1/\alpha)}
  \Big)^{n^2}
  { e^{-n^{2-\alpha}t^\alpha/\alpha}\over
    t^{(\alpha -2){n^2\over 2}+{\alpha\over 2}} } \times
  \cut
  { (2\pi )^{(n^2-1)/2} 2^{2n-1} \over
	n^{(2-\alpha )(n^2+1)/2}(\alpha -2)^{n-1}
	(\alpha -1)^{(n-1)^2/2} }
  \quad\cr}
$$
as $t$ tends to infinity, which is the result.

Let us now turn to the case $\alpha =2$. From Lemma 10.4.2, we
conclude
$$
  I(A_1)=1/2 \, .
$$
The dominating manifold
$$
  \da1 = \Sn1\otimes \Sn1 
$$
is now of dimension $k=2(n-1)$. Since $S_{n-1}\times S_{n-1}$
is a double covering of $\calD_{A_1}$, 
the Riemannian measure on $\calD_{A_1}$ is half the
product measure on the product of $2$ spheres $S_{n-1}$, 
each having the
Riemannian measure obtained from the Lebesgue measure on $\RR^n$.

On $\da1$, we also have
$$
  |\D I(x)| 
  = \Big( \sum_{1\leq i,j\leq n} |u_iv_j|^2\Big)^{1/2}
  = |u||v|=1 \, .
$$
The computation of the curvature term $\det\, G_{A_1}$ goes the
same way as for $\alpha >2$. Namely, we still have
$$
  \Pi_{\Lambda_{I(A_1)}}
  =\Id_{n^2-1}
  \qquad\hbox{ and } \qquad
  \Pi_{\partial A_1} 
  =\pmatrix{ \Id_{2(n-1)} & 0 \cr 0 & 0 \cr} \, .
$$
In particular,
$$
  \Pi_{\Lambda_{I(A_1)}}-\Pi_{\partial A_1}
  = \pmatrix{ 0 & 0 \cr 0 & \Id_{n^2-1-2(n-1)}\cr }
  = \pmatrix{ 0 & 0 \cr 0 & \Id_{(n-1)^2} \cr}
  \, .
$$
It follows that
$$
  G_{A_1}= \Id_{(n-1)^2}
$$
is of determinant $1$. Therefore, with the notation of Theorem
7.1,
$$
  c_1
  =(2\pi )^{(n^2-2(n-1)-1)/2} {\vol (\Sn1 )^2 \over 2}
  = (2\pi )^{(n-1)^2/2}{\vol (\Sn1 )^2\over 2} \, .
$$
Again, taking all the above estimates into account, we obtain
$$
  P\big\{\, \| X\|\geq t\,\big\}
  \sim {1\over 2} (2\pi )^{((n-1)^2-n^2)/2}\vol (\Sn1 )^2
  e^{-t^2/2}t^{2(n-1)-1}
$$
as $t$ tends to infinity.  This is the result since $\Sn1$ has volume
$2\pi^{n/2}/\Gamma (n/2)$.

It remains to check the assumptions of Theorem 7.1. We already checked
(7.3) and (7.4). Assumption (7.5) is clear as well since the curvature
of $\partial A_1$ is bounded.
\hfill$\qed$

\bigskip

We can now apply Theorem 7.5 to obtain the following result on
conditional distributions. It is worth knowing that $S_1\times
S_1 /\{\, -\Id,\Id\,\}$ is the usual Klein bottle. Hence,
$S_{n-1}\otimes S_{n-1}\equiv S_{n-1}\times S_{n-1}/\{\, -\Id,\Id\,\}$ is a
$2(n-1)$-dimensional Klein bottle.

\bigskip

\state{10.4.3. PROPOSITION.} {\it%
  Let $X={(X_{i,j})}_{1\leq i,j\leq n}$ be a random matrix with
  independent and identically distributed coefficients, all
  having the density 
  $$
    w_\alpha (s) = 
    { \alpha^{1-(1/\alpha)}\over 2\Gamma (1/\alpha )}
    e^{-|s|^\alpha /\alpha} \, , \qquad s\in\RR\, .
  $$
  The conditional 
  distribution of $X/t$ given $\| X\|\geq t$ converges 
  weakly* to a uniform distribution over
  
  \noindent (i) the Klein bottle $\Sn1\otimes\Sn1$ if $\alpha =2$,

  \noindent (ii) the $2^{2n-1}$ matrices of the
  form $u_*\otimes v_*$ for $|u_{*,i}|=|v_{*,i}|=1/\sqrt n$,
  $1\leq i\leq n$, if $\alpha >2$.
  }

\bigskip

\stateit{Proof.} It follows from the calculation of $\da1$ made 
in the proof of Theorem 10.4.1 and Theorem 7.5. \hfill$\qed$

\bigskip

Let us now turn to the problem of estimating the tail
probability of $\| X\|$ when the coefficients $X_{i,j}$ of the
random matrix $X$ are independent and identically distributed
with a Student-like distribution. Given the work done in the previous
sections, this turns to be an easy problem.

\bigskip

\state{10.4.4. THEOREM.} {\it%
  Let $X={(X_{i,j})}_{1\leq i,j\leq n}$ be a random matrix
  with independent coefficients, all having a Student-like
  distribution with parameter $\alpha$. Then,
  $$
    P\big\{\, \| X\|\geq t\,\big\}
    \sim {2n^2K_{s,\alpha}\alpha^{(\alpha -1)/2}\over t^\alpha}
    \qquad\hbox{ as } \ttoi \, .
  $$}

\bigskip

\stateit{Proof.} Notice that the set
$$
  A_t=\big\{\, x\in\MnR  : \| x\|\geq t\,\big\}
  = tA_1
$$
is the complement of a convex set, namely the ball of radius $t$
centered at the origin in $\MnR$ endowed with the operator norm.  In
$\MnR\equiv\RR^{n^2}$, the axial points of this convex set are all the
matrices $\epsilon \Eij$, $\epsilon\in\{\, -1,1\,\}$, $1\leq i,j\leq n$, 
which are of Euclidean norm $1$. There are $2n^2$ such
matrices. Apply Theorem 9.3.1 to obtain the result.\hfill$\qed$

\bigskip

With no extra effort, we can also obtain the following result on
conditional distribution.

\bigskip

\state{10.4.5. PROPOSITION.} {\it%
  Let $X={(X_{i,j})}_{1\leq i,j\leq n}$ be a random matrix
  with independent coefficients, all having a Student-like
  distribution with parameter $\alpha$. The distribution of $X/t$
  given $\| X\|\geq t$ converges weakly* to the uniform
  mixture of the distributions of the matrices $\epsilon Z\Eij$
  with $\epsilon$ in $\{\, -1,1\,\}$ and $Z$ having a Pareto
  distribution,
  $$
    P\{\, Z\geq 1+\lambda\,\} 
    = {1\over (1+\lambda )^\alpha} \, , 
    \qquad \lambda\geq 0 \, .
  $$}

\bigskip

\stateit{Proof.} Apply Corollary 9.3.4. The
axial points of $A_1^{\rm c}$ are the matrices $\epsilon\Eij$,
where $\epsilon$ is in $\{\, -1,1\,\}$. Those matrices are
of unit Euclidean norm in $\RR^{n^2}$. \hfill$\qed$

\bigskip

\centerline{\fsection Notes}
\section{10}{Notes}

\bigskip

\noindent The theory of random matrices has been evolving quite fast
lately. Motivated by applications in physics and in
operator algebras, important progress has been made on the
asymptotic theory as the size of the matrix goes to infinity.
Our fixed size viewpoint is quite different. Amazingly clever
explicit calculations have been made in the Gaussian cases and
some of its variations. A classical reference is Mehta (1991).
Another aspect driven by statistics concerns the Wishart
distribution --- see Johnson and Kotz (1972).

I believe Lemma 10.2.3 is not new, but I have not found a
reference for it. It is very similar to Theorem 2.1 of Rosi\'nski
and Woyczy\'nski (1987), as well as its proof. If we assume that
the $X_i$'s are symmetric, then Lemma 10.2.3 can be deduced from
Rosi\'nski and Woyczy\'nski (1987) in conditioning on the signs
of the $X_i$'s. But here, we assume only asymptotic symmetry of
the tail. Therefore, the signs of $X_i$'s given $|X_i|$ large is
only asymptotically distributed uniformly over $\{\, -1,+1\,\}$.

Theorem 10.1.4 involves the volume of $\SOnR$. It seems to be
calculated in Marinov (1980). But, by ignorance, I have not 
been able to follow his proof. I don't know if Marinov's results
give the volume of $\SOnR$ embedded in $\RR^{n^2}$ of if they
give it up to a proportionality constant.

\vfill\eject

\chapter{11}{Finite sample results for autoregressive processes}
\line{\hfill{\fchapter 11.\quad Finite sample results}}
\vskip .1in
\line{\hfill{\fchapter for autoregressive processes}}

\vskip 54pt

\noindent Autoregressive models are among the simplest and most
widely used models in statistical analysis of time series. Their
classical theory deals mainly with their asymptotic behavior over
a very large time period. In this chapter, we will see that these
models are in fact much more subtle than usually believed. Our
study will build upon results of Chapter 8. Our results will not
exhaust the topic by any mean; they should be considered as an
incentive for further study.

\bigskip

\noindent{\fsection 11.1. Background on autoregressive processes.}
%\chapternumber={11.1}
\section{11.1}{Background on autoregressive processes}
\bigskip

\noindent In order to describe the processes we are interested
in, let us introduce the backward shift\ntn{$B$} $B$ on vectors. For a
vector $u=(u_1,\ldots ,u_n)$ in $\RR^n$, we write $Bu=(0,u_1,\ldots
,u_n)$. Let $\epsilon$ be a mean zero random vector in
$\RR^n$, with
independent and identically distributed components.
We say that the vector $X$ in $\RR^n$ is an autoregressive
process of order $p$ with innovation $\epsilon$ if for some\ntn{$\theta$}
$\theta=(\theta_1,\ldots, \theta_p)$ in $\RR^p$, with
$\theta_p$ not null, it satisfies the equation\ntn{$X$}\ntn{$\epsilon$}
$$
  X=\sum_{1\leq i\leq p} \theta_i B^iX +\epsilon \, .
  \eqno{(11.1.1)}
$$
In a perhaps more explicit form, that means
\vskip\abovedisplayskip
\halign{\kern .7in$#$ &$#$\hfil &$#$ & $#$ & $#$ & $#$ & $#$\cr
  X_1     &=\epsilon_1 & & & & &\cr
  X_2     &=\theta_1X_1       &+&\epsilon_2 & & &\cr
  \kern 0.05in\vdots\hfill& & &  & & & \cr
  X_{p+1} &=\theta_1X_p       &+&\theta_2X_{p-1} & +\;\cdots\; + &
           \theta_pX_1 & + \;\epsilon_{p+1} \cr
  X_{p+2} &=\theta_1X_{p+1}   &+&\theta_2X_p & +\;\cdots\;+ &
            \theta_pX_2 & + \;\epsilon_{p+2} \cr
  \kern 0.05in\vdots\hfill& & & & & & \cr
  X_n     &= \theta_1 X_{n-1} &+&\theta_2X_{n-2} & +\;\cdots\;+ &
            \theta_p X_{n-p} & + \;\epsilon_n \; . \cr
  }
\vskip\belowdisplayskip
For statisticians, the main questions are on estimation and tests
procedures for such models. This means that one observes the vector
$X$, and knows that it is of the form (11.1.1) with the $\epsilon_i$'s
independent and identically distributed. The goal is then to estimate
$\theta$, that is to guess its value based on the knowledge of $X$; or
to perform tests on $\theta$, that is to check if some assumption on
$\theta$ is compatible with the observed value of $X$. How is this
done?

Consider the $n\times p$ matrix $\calX=(BX,\ldots , B^pX)$.
Equation (11.1.1) becomes
$$
  X=\calX \theta +\epsilon \, .
$$
Thus, $X$ is a point in the space spanned by $\calX$ plus a random
vector.  A reasonable guess for $\theta$ is $\theta_{LS}$ such that
$\calX\theta_{LS}$ is the projection of $X$ onto the space spanned by
$BX,\ldots , B^pX$. This is called the least square estimator of
$\theta$. Whenever $\calX$ is of rank $p$, we have
$$
  \theta_{LS}=(\calX^\T\calX)^{-1}\calX^\T X \, .
$$
This can be calculated solely on the observed $X$.

Furthermore, notice that the $(i,j)$-entry of the matrix
$\calX^\T\calX$ is
$$\eqalign{
  \< B^iX,B^jX\> 
  &= \sum_{1+ (i\vee j)\leq k\leq n} X_{k-i}X_{k-j} \cr
  &= \sum_{1+ |j-i|\leq r\leq n-(i\wedge j)}X_rX_{r-|i-j|} \, , \cr
  }
$$
while the $i$-th coordinate of $\calX^\T X$ is $\< B^iX,X\>$.

It is customary to define the empirical autocovariances of order
$k<n$ as\ntn{$\gamma_n(k)$}
$$
  \gamma_n(k) = n^{-1} \sum_{k+1\leq r\leq n} X_rX_{r-k} \, .
$$
Notice that 
$$
  \< B^iX,B^jX\> -n\gamma_n(|i-j|)
  = \sum_{n-(i\wedge j)< r\leq n} X_rX_{r-|i-j|} \, .
$$
Therefore, whenever $X_n=O_P(1)$ as $n$ tends to infinity, 
and $i-j$ is fixed, we have
$$
  \< B^iX,B^jX\> = n\gamma_n(|i-j|)+O_P(1) \qquad
  \hbox{ as } n\to\infty\, .
  \eqno{(11.1.2)}
$$
This explains why the most popular estimator of $\theta$ is not
$\theta_{LS}$ but the following substitute. Define the matrix
$$
  \Gamma_n = \big( \gamma_n(|i-j|)\big)_{1\leq i,j\leq p}
$$
and the vector
$$
  \gamma_n=\big( \gamma_n(i)\big)_{1\leq i\leq p} \, .
$$
If the process $X_n$ is of order $1$ as $n$ tends to infinity, 
then (11.1.2) shows that
$$
  (n^{-1}\calX^\T\calX )^{-1} = \Gamma_n^{-1} + O_P(n^{-1})
$$
while
$$
  \calX^\T X=\gamma_n +O_P(n^{-1}) \, .
$$
Thus, instead of using $\theta_{LS}$, one tends to guess
$\theta$ by
$$
  \hat\theta_n =\Gamma_n^{-1}\gamma_n \, .
$$
Actually, whether we use $\theta_{LS}$ or $\hat\theta_n$ does not
really matter much for our purposes. What we do care about is that $X$ is
a linear function of $\epsilon$ as (11.1.1) shows. The autocovariance
$\gamma_n(k)$ is a quadratic form in $X$, as well as in $\epsilon$, a
classical fact.  It is then plain that tail probabilities of
$\gamma_n(k)$ are relevant to statistics, and that chapter 8 provides
the right estimates.

Maybe in order to fully enjoy the results we are going to prove,
one should know the basics of the classical theory. To make this
text quite selfcontained, let us sketch it. To this end, define
the polynomial
$$
  \Theta (z)=1-\sum_{1\leq i\leq p}\theta_i z^i \, ,
  \qquad z\in\CC \, .
$$
Equation (11.1.1) can be rewritten as
$$
  \Theta (B)X=\epsilon \, .
$$
Denote by $r_1,\ldots , r_p$ the complex roots of $\Theta$. Then
$\Theta (z)=\prod_{1\leq i\leq p}(1-r_i^{-1}z)$. We can write
$1/\Theta (z)$ formally as a series. This is done conveniently
by introducing the vector $r=(r_1,\ldots , r_p)$. Whenever
$s=(s_1,\ldots ,s_p)$ belongs to $\ZZ^p$, 
we write $|s|=s_1+\cdots + s_p$
and $r^s=\prod_{1\leq i\leq p}r_i^{s_i}$. We have
$$\eqalign{
  1/\Theta (z) 
  & = \prod_{1\leq i\leq p} (1-r_i^{-1}z_i)^{-1}
    = \prod_{1\leq i\leq p} \sum_{k\geq 0} r_i^{-k} z_i^k \cr
  & = \sum_{k\geq 0} \Big( \sum_{|s|=k}r^{-s}\Big)z^k \, . \cr
  }
$$
Substituting $B$ for $z$, we can formally define $\Theta
(B)^{-1}$. Writing $\epsilon_i=0$ if $i\leq s$, one can then
check that $X=\Theta (B)^{-1}\epsilon$, that is
$$
  X_n =\sum_{k\geq 0}\Big(\sum_{|s|=r}r^{-s}\Big)\epsilon_{n-k} 
  \, . \eqno{(11.1.3)}
$$
When all the roots $r_i$ are outside the unit circle, there
exists a positive $\eta$ such that 
$$\eqalign{
  \big| \sum_{|s|=k}r^{-s}\big|
  & \leq (1+\eta )^{-k} \sharp \{\, s\, :\, |s|=k\,\} \cr
  & = (1+\eta )^{-k} {k+p-1\choose p-1} \cr
  & \sim (1+\eta )^{-k} {k^p\over (p-1)!}
    \qquad\hbox{ as } k\to\infty \, . \cr
  }
$$
Therefore if the residuals $\epsilon_i$ have a tail which 
decays fast
enough the distribution of $X_n$ converges weakly* to that of
$\sum_{k\geq 0}(\sum_{|s|=k}r^{-s})\epsilon_k$. A simple
condition on the tail of $\epsilon_i$ for this convergence to 
hold is
$$
  \int_1^\infty P\{\, |\epsilon_i|\geq t\,\}
  {\d t\over t} <\infty \, .
$$
A more stringent one is to assume that $\epsilon_i$ is
integrable; in this case $X_n$ converges in $L^1$ as well.

If now some roots are inside the unit disk, we can first assume that
$r_1$ is the unique root with smallest modulus; and therefore $|r_1|$
is less than $1$. Then, (11.1.3) shows that the distribution of
$r_1^nX_n$ converges weakly* to a nondegenerate limit. Since $r_1^n$
converges to $0$ exponentially fast, this amounts to saying that the
process $X_n$ explodes at exponential rate. Some complications occur
if the smallest root is not unique, but $X_n$ still explodes,
essentially at an exponential rate.

As a consequence, the asymptotic behavior of the empirical
autocovariances as $n$ tends to infinity is very different according
to the location of the roots $r_i$ with respect to the unit disk.

In conclusion, the classical theory makes a great deal of the location
of the roots of $\Theta$ with respect to the unit disk.  And it is
essentially all that it cares about, because only the behavior as the
time $n$ goes to infinity is considered. In the following sections, we
will show that it is only a part of the overall behavior of these
processes.

To end this section, let us examine this root question for
autoregressive models of order $1$ and $2$.

For an autoregressive process of order $1$, we write
$X_i=aX_{i-1}+\epsilon_i$. If $|a|<1$, this process is
nonexplosive. This can be represented on the real line as follows.

\setbox1=\vbox to .8in{\ninepointsss
  \hsize=3in
  \kern .5in
  \includegraphics{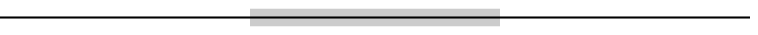}
  \anote{3.05}{.24}{$a$}
  \anote{.93}{.1}{$-1$}
  \anote{1.97}{.1}{$1$}
  \vfill
  \centerline{Shaded region for the nonexplosive domain}
  \hfill
}
\centerline{\box1}

For an autoregressive process of order $2$, we write it as
$X_i=aX_{i-1}+bX_{i-2}+\epsilon_i$. We need to determine
where $a$, $b$ should lie for the polynomial $x^2-ax-b$ to 
have all its roots within the unit disk. If $a^2+4b$
is negative, the roots
are complex, conjugate to each others. They are in the unit
disk if and only if their product is less than $1$, that is if
$-b<1$. If $a^2+4b$ is nonnegative and $a$ is positive, 
the largest root in absolute
value is $(a+\sqrt{a^2+4b})/2$. It is less than $1$ if $a^2+4b
<(2-a)^2 = a^2-4a + 4$, that is $b+a<1$. One can argue similarly
if $a$ is negative, and we obtain the following triangular domain.

\setbox1=\vbox to 2.8in{\ninepointsss
  \hsize=3in
  \kern 2.5in
  \includegraphics{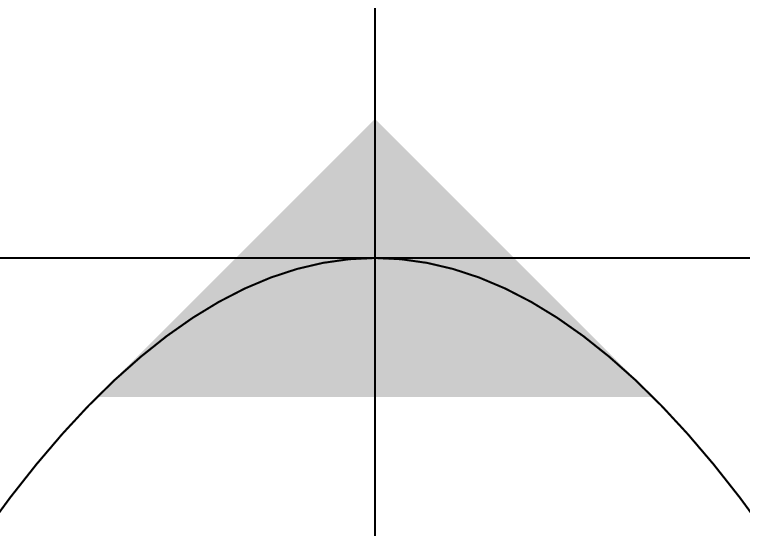}
  \anote{1.4}{2.2}{$b$}
  \anote{3.05}{1.27}{$a$}
  \anote{2.05}{1.35}{$1$}
  \anote{.85}{1.35}{$-1$}
  \anote{3}{.4}{$b=-a^2/4$}
 % \vfill
  \centerline{Shaded region for the nonexplosive domain}
  \vfill\hfill
}
%\vbox{\kern 2.5in}
\centerline{\box1}
%\kern -2,5in
\bigskip

\noindent{\fsection 11.2. Autoregressive process of order 1.}
%\chapternumber={11.2}
\section{11.2}{Autoregressive process of order 1}

\bigskip

\noindent In this section, we investigate the 
tail behavior of the
autocovariances of autoregressive processes of order $1$,\ntn{$a$}
$$\eqalign{
  X_1&=\epsilon_1\cr
  X_n&=aX_{n-1}+\epsilon_n \, , \qquad n\geq 2 \, . \cr
  }
$$
Recall that $\epsilon=(\epsilon_1,\ldots , \epsilon_n)$ and
$X=(X_1,\ldots ,X_n)$. Defining\ntn{$A$}
\setbox2=\hbox{$\mkern-40mu\mddots{19}{23}$}
\wd2=0pt\ht2=0pt\dp2=0pt
\setbox3=\hbox{$\mkern-40mu\mddots{15}{35}$}
\wd3=0pt\ht3=0pt\dp3=0pt
\setbox4=\hbox{$\mkern-23mu\mddots{10}{34}$}
\wd4=0pt\ht4=0pt\dp4=0pt
$$
  A=\pmatrix{ 1         &       &        &    & \cr
	      a         &       &\hbox{\ftitre 0}&    & \cr
              a^2       & \box2 &        &    & \cr
              \mvdots{5}& \box3 &        &    & \cr
	      a^{n-1} &\box4\mkern-5mu\mcdots{5}  
          &\mkern-2mu a^2 &\mkern -2mu a  & 1 \cr} \, ,
$$
we see that $X=A\epsilon$. The matrix of the backward 
shift on $\RR^n$ is\ntn{$B$}
\setbox3=\hbox{$\mkern-20mu\mddots{10}{10}$}
\wd3=0pt\ht3=0pt\dp3=0pt
\setbox4=\hbox{$\mkern-20mu\mddots{7}{12}$}
\wd4=0pt\ht4=0pt\dp4=0pt
\setbox5=\hbox{\raise 5pt\hbox{\kern -13pt\ftitre 0}}
\wd5=0pt\ht5=0pt\dp5=0pt
\setbox6=\hbox{\raise 5pt\hbox{\kern -14pt\ftitre 0}}
\wd6=0pt\ht6=0pt\dp6=0pt
$$
  B=\pmatrix{ 0 &        &   &  \cr
	      1 & \box3 &   &  \box6 \cr
                & \box4 &    &   \cr
	        & \phantom{0}\box5  & \mkern 7mu 1 & \mkern -5mu 0 \cr} \, .
$$
The empirical covariance of order $k$,
$$
  n\gamma_n(k)
  =\sum_{k+1\leq i\leq n} X_iX_{i-k} = \<X,B^kX\>
  = \< A^\T B^kA\epsilon ,\epsilon \>
$$
is a nice quadratic form in $\epsilon$.
If $k\geq n$, then $\gamma_n(k)=0$. We assume from now on that
$k<n$. When $\epsilon$ has a heavy tail, the tail behavior of
$n\gamma_n(k)$ depends on the value of $a$. It is given by the
following result.

\bigskip

\state{11.2.1. THEOREM.} {\it%
  Let $X$ be an autoregressive process of order one,
  with coefficient $a$ and independent
  and identically distributed innovations $\epsilon_i$ having
  a Student-like distribution with parameter $\alpha$. 
  The following expressions
  are equivalent to $P\{\, n\gamma_n(k)\geq t\,\}$ as $t$ 
  tends to infinity.

  \noindent (i) If $a=0$ and $k=0$,
  $$
    K_{s,\alpha}\alpha^{(\alpha-1)/2} 2n\, t^{-\alpha/2} \, .
  $$
  \noindent (ii) If $a=0$ and $k\geq 1$,
  $$
    K_{s,\alpha}\alpha^\alpha  2(n-k)_+ \, t^{-\alpha}\log t \, .
  $$ 
  \noindent (iii) $k$ is even and $a\neq 0$, or $k$ is odd and
  $a>0$,
  $$
    K_{s,\alpha}\alpha^{(\alpha-1)/ 2} 2
    a^{k\alpha/2}\sum_{1\leq i\leq n-k}
    \Big| {1-a^{2i}\over 1-a^2}\Big|^{\alpha /2}\,
     t^{-\alpha/2} \, .
  $$
  \noindent (iv) $k$ is odd and $a<0$,
  $$\displaylines{\qquad
      K_{s,\alpha}^2\alpha^\alpha t^{-\alpha}\log t \times
      \cut
      \sum_{n-k+1\leq i\leq n}\sum_{1\leq j\leq n}
      a^{\alpha |i-k-j|}
      \Big|{ 2-a^{2(n+1)} 
             (a^{2(i\vee (j+k))} - a^{2(j\vee (i+k))})
             \over 1-a^2 }\Big|^\alpha
    \cr}
  $$}
\finetune{
}
It is implicit in the statement that the function $a\mapsto
(1-a^i)/(1-a^2)$ is extended by continuity at $a=1$. Its
value for $a=1$ is $i/2$.

The striking fact is that odd and even autocovariances exhibit
very different decays when $a$ is negative; the former are of order
$t^{-\alpha/2}$, the latter of order $t^{-\alpha}\log t$.

\bigskip

\stateit{Proof.} Define the matrix\ntn{$C$} $C=A^\T B^kA$. We apply the
results of chapter 8. We need to check if the largest diagonal
coefficient of $C$ is positive, zero, or negative. In order to
calculate it, notice that
\vskip\belowdisplayskip
\halign{#\kern .5in&\hfill$#$ & $#$ \hfill\cr
       & A_{i,j}& =\cases{ a^{i-j} & if $i\geq j$,\cr
			  0       & otherwise, \cr} \cr
  \noalign{\vskip 2pt}
  and& & \cr
  \noalign{\vskip 2pt}
       &(B^k)_{i,j} & = \cases{ 1 & if $k+1\leq i\leq n$ and 
                                    $j=i-k$,\cr
				0 & otherwise. \cr}
  \cr}%
\vskip\abovedisplayskip
\noindent Consequently,
$$\eqalign{
  C_{i,i}
  &=\sum_{1\leq j,m\leq n} (A^\T)_{i,m}(B^k)_{m,j}A_{j,i}
   =\sum_{i+k\leq m\leq n} A_{m,i}A_{m-k,i} \cr
\noalign{\vskip .1in}
  &=\cases{ a^k{\displaystyle 1-a^{2(n-i-k+1)}\over 
            \displaystyle 1-a^2} & if
		$a^2\neq 1$ and $n-i-k\geq 0$, \cr
	    a^k(n-i-k+1) & if $a^2=1$ and $n-i-k\geq 0$,\cr
            \noalign{\vskip .08in}
	    0 & if $i\geq n-k+1$. \cr}\cr
  }
$$

Assume $a=0$ and $k=0$. Then $C=\Id$. Statement (i) of Theorem
11.2.1 follows from Theorem 8.2.1.

If $a$ is null and $k$ is nonzero, the matrix $C=B^k$ has 
all its diagonal
elements vanishing. We apply Theorem 8.3.1, calculating
$$
  \sum_{i:C_{i,i}=0}\sum_{1\leq j\leq n}
  |(B^k)_{i,j}+(B^k)_{j,i}|^\alpha
  = 2(n-k)_+ \, .
$$

When $k$ is even and $a$ is nonzero, then $C_{i,i}$ is positive 
for any $1\leq i\leq n-k$. We apply Theorem 8.2.1, calculating
$$
  \sum_{1\leq i\leq n-k}C_{i,i}^{\alpha/2}
  = a^{k\alpha/2}\sum_{1\leq i\leq n-k} 
    \Big({1-a^{2(n-i-k+1)}\over 1-a^2}\Big)^{\alpha/2} \, .
$$
This gives statement (iii), after substituting $n-i-k+1$ for $i$ in
the summation.

Statement (iii), when $k$ is odd and $a$ is positive, follows from 
exactly the same calculation.

Let us now concentrate on $k$ odd and $a$ negative. Then $a^k$ is
negative, and so is $C_{i,i}$ if $1\leq i\leq n-k$, while $C_{i,i}$
vanishes if $i\geq n-k+1$. Therefore, we apply Theorem 8.3.1. We need
to calculate
$$
  \sum_{n-k+1\leq i\leq n} \sum_{1\leq j\leq n}
  |C_{i,j}+C_{j,i}|^\alpha \, .
$$
We have
$$\eqalign{
  C_{i,j}
  & =\sum_{1\leq l,m\leq n} (A^\T)_{i,l} (B^k)_{l,m} A_{m,j}
    =\sum_{i\vee (j+k)\leq l\leq n} a^{l-i}a^{l-k-j} \cr
  & = a^{|i-j-k|}{1-a^{2(n+1-(i\vee (j+k)))}\over 1-a^2} \, .
   \cr
  }
$$
We obtain $C_{j,i}$ by permuting $i$ and $j$. This gives
statement (iv).\hfill$\qed$

\bigskip

How good are these approximations? Looking at the bound in
Theorem 3.1.9 and how we derived Theorem 5.1, we cannot expect
them to be good when we are integrating in a high dimensional
space, that is when $n$ is large.

A plot of the approximations given in Theorem 11.2.1 does not
show much, since all the probabilities go to $0$ as $t$ tends
to infinity. When comparing the tails, it makes more sense 
to look at the relative error. This leads to the plot
$$
  \log P\{\, n\gamma_n(k)\geq t\,\}
$$
as well as the logarithm of the approximation. These should
be approximately in linear relation with $\log t$. Therefore,
the plots below will show the function $t\mapsto P\{\,
\gamma_n(k)\geq t\,\}$ with both axes in logarithmic scale.
Since we do not know a closed formula for $P\{\, n\gamma_n(k)\geq
t\,\}$, we obtained this probability by simulation. We generated
100,000 replicas of $\epsilon$. As $t$ increases, the estimate of the
true probability is based on less and less points; the simulated curve
tends to wiggle as $t$ gets large. The theoretical approximation will
be the smooth curve on the graphs.  The parameters involved are $k$,
$n$, $\alpha$, $a$. We will only consider the autocovariance of order
$1$ in our simulations.  We consider the sample sizes $n=10$, which is
very small, and $n=20$, which is a common order of magnitude in some
applications. We also consider probabilities of interest in
applications, namely between $10^{-1}$ and $10^{-3}$.  Recall that
$5\%$ is about $10^{-1.30}$.

The results are as follows.

Let us first see what happens when the errors have a Cauchy
distribution, corresponding to $\alpha=1$.

When $a=1$, the two plots bellow show that the approximation is
amazingly good.
\bigskip

\setbox1=\hbox to 2.05in{\vbox to 2.5in{\ninepointsss
  \hsize=2.05in
  \kern 2.1in
  \includegraphics{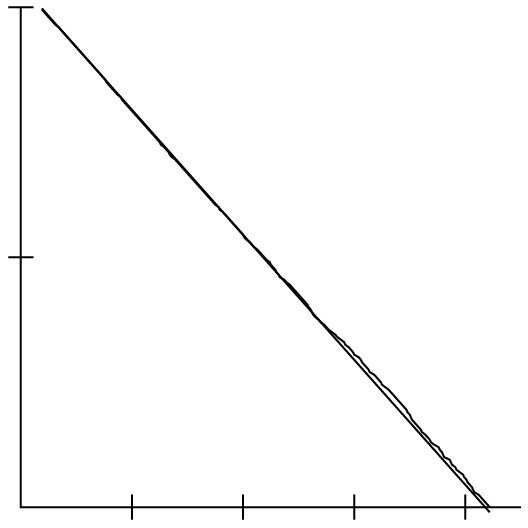}
  \anote{.35}{-.2}{$10^5$}
  \anote{.8}{-.2}{$10^6$}
  \anote{1.25}{-.2}{$10^7$}
  \anote{1.7}{-.2}{$10^8$}
  \anote{-.3}{1}{$10^{-2}$}
  \anote{-.3}{2}{$10^{-1}$}
  \vfill
  \centerline{$n=10$, $a=1$, errors Cauchy}
}}
\setbox2=\hbox to 2.05in{\vbox to 2.5in{\ninepointsss
  \hsize=2.05in
  \kern 2.1in
  \includegraphics{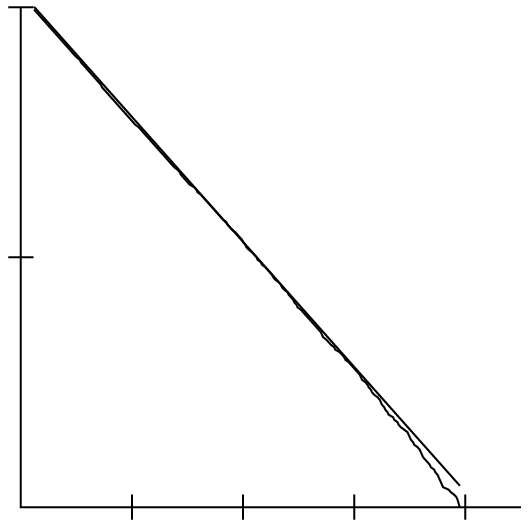}
  \anote{.35}{-.2}{$10^6$}
  \anote{.8}{-.2}{$10^7$}
  \anote{1.25}{-.2}{$10^8$}
  \anote{1.7}{-.2}{$10^9$}
  \anote{-.3}{1}{$10^{-2}$}
  \anote{-.3}{2}{$10^{-1}$}
  \vfill
  \centerline{$n=20$, $a=1$, errors Cauchy}
}}

\line{\kern 0.3in\box1\hfill\box2}
\bigskip

For $a=0.5$ the approximation is also excellent.

\bigskip
\setbox1=\hbox to 2.1in{\vbox to 2.5in{\ninepointsss
  \hsize=2.1in
  \kern 2.1in
  \includegraphics{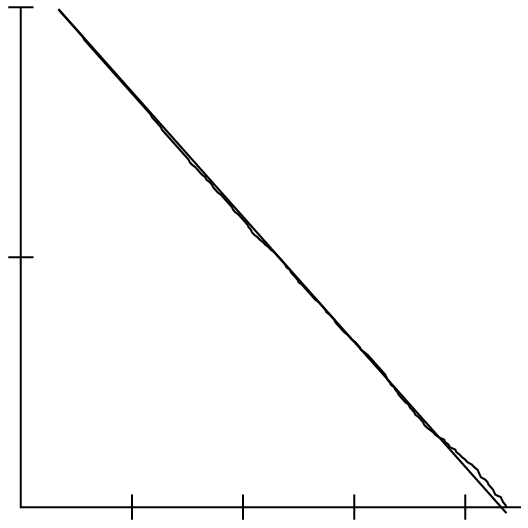}
  \anote{.35}{-.2}{$10^4$}
  \anote{.8}{-.2}{$10^5$}
  \anote{1.25}{-.2}{$10^6$}
  \anote{1.7}{-.2}{$10^7$}
  \anote{-.3}{1}{$10^{-2}$}
  \anote{-.3}{2}{$10^{-1}$}
  \vfill
  \centerline{$n=10$, $a=0.5$, errors Cauchy}
}}
\setbox2=\hbox to 2.05in{\vbox to 2.5in{\ninepointsss
  \hsize=2.05in
  \kern 2.1in
  \includegraphics{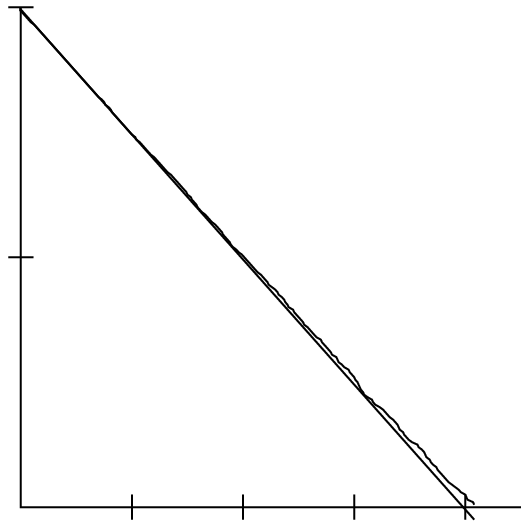}
  \anote{.35}{-.2}{$10^5$}
  \anote{.8}{-.2}{$10^6$}
  \anote{1.25}{-.2}{$10^7$}
  \anote{1.7}{-.2}{$10^8$}
  \anote{-.3}{1}{$10^{-2}$}
  \anote{-.3}{2}{$10^{-1}$}
  \vfill
  \centerline{$n=20$, $a=0.5$, errors Cauchy}
}}
%\wd2=2in

\line{\kern 0.3in\box1\hfill\box2}
\bigskip

In the degenerate case where we need to apply Theorem 8.3.1, 
the coefficient $a$ vanishes. The approximation is not
as good as before. But taking into account that we are
visualizing a relative error, it performs well enough
to be of practical use.

\kern 0.2in

\setbox1=\hbox to 2.05in{\vbox to 2.5in{\ninepointsss
  \hsize=2.05in
  \kern 2.1in
  \includegraphics{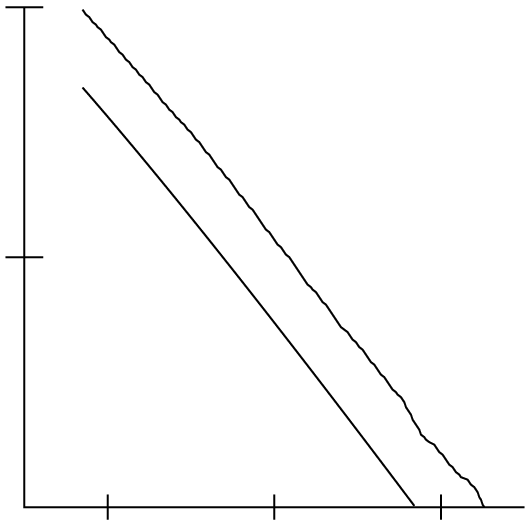}
  \anote{.24}{-.2}{$10^2$}
  \anote{.9}{-.2}{$10^3$}
  \anote{1.56}{-.2}{$10^4$}
  \anote{-.3}{1}{$10^{-2}$}
  \anote{-.3}{2}{$10^{-1}$}
  \vfill
  \centerline{$n=10$, $a=0$, errors Cauchy}
}}
\setbox2=\hbox to 2.05in{\vbox to 2.5in{\ninepointsss
  \hsize=2.05in
  \kern 2.1in
  \includegraphics{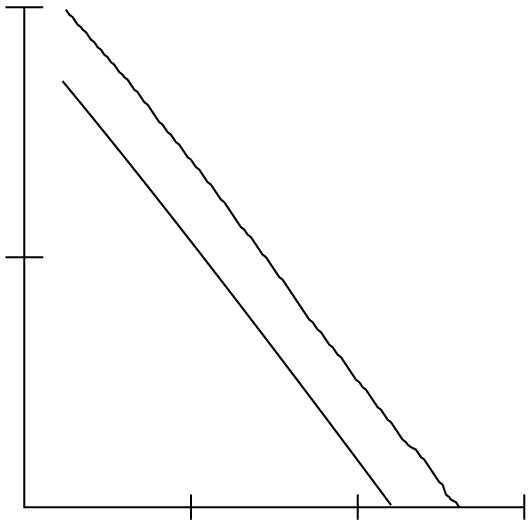}
  \anote{.57}{-.2}{$10^3$}
  \anote{1.23}{-.2}{$10^4$}
%  \anote{1.89}{-.2}{$10^5$}
  \anote{-.3}{1}{$10^{-2}$}
  \anote{-.3}{2}{$10^{-1}$}
  \vfill
  \centerline{$n=20$, $a=0$, errors Cauchy}
}}

\line{\kern 0.3in\box1\hfill\box2}
\bigskip

For $a=-0.5$, the autocovariance tends to be negative since
$\gamma_n(1)/\gamma_n(0)$ is an approximation of $a$. Therefore
we need to go further on the tail to have a good approximation.
For $n=10$, it is still accurate enough to be of 
some practical interest.
One can use the approximation to find critical values at levels
less than $10^{-1.5}\approx 3\%$ say. But as $n$ increases from
$10$ to $20$, the accuracy decreases.

\bigskip
\setbox1=\hbox to 2.05in{\vbox to 2.5in{\ninepointsss
  \hsize=2.05in
  \kern 2.1in
  \includegraphics{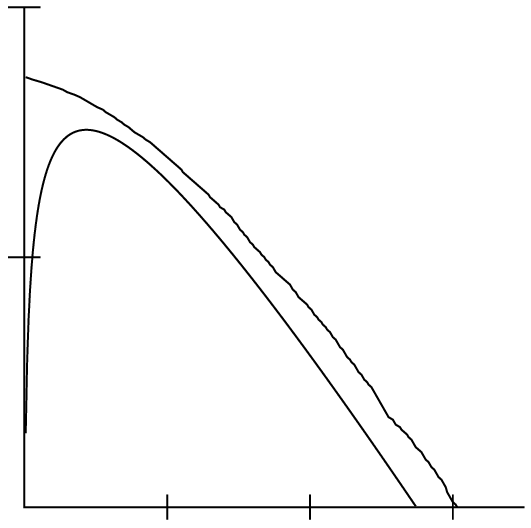}
  \anote{.48}{-.2}{$10^1$}
  \anote{1.05}{-.2}{$10^2$}
  \anote{1.62}{-.2}{$10^3$}
  \anote{-.3}{1}{$10^{-2}$}
  \anote{-.3}{2}{$10^{-1}$}
  \vfill
  \centerline{$n=10$, $a=-0.5$, errors Cauchy}
}}
\setbox2=\hbox to 2.05in{\vbox to 2.5in{\ninepointsss
  \hsize=2.05in
  \kern 2.1in
  \includegraphics{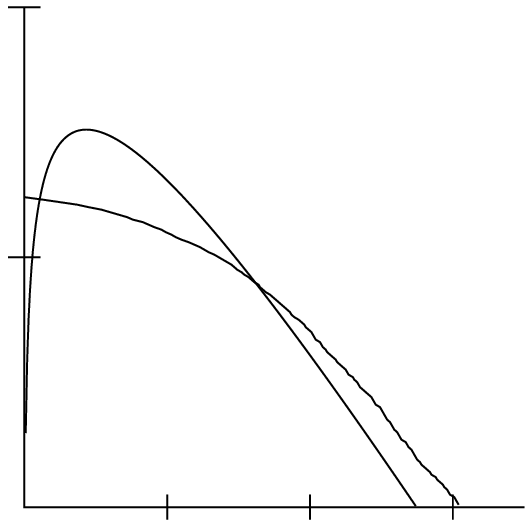}
  \anote{.48}{-.2}{$10^1$}
  \anote{1.05}{-.2}{$10^2$}
  \anote{1.62}{-.2}{$10^3$}
  \anote{-.3}{1}{$10^{-2}$}
  \anote{-.3}{2}{$10^{-1}$}
  \vfill
  \centerline{$n=20$, $a=-0.5$, errors Cauchy}
}}
\line{\kern 0.3in\box1\hfill\box2}

\bigskip

For $a=-1$, we need to go much further in the tail of the
distribution of $\gamma_n(1)$ in order to have a positive
quantile. The approximation is not accurate in the range of
practical interest. As $n$ increases from $10$ to $20$, the
approximation cannot be used in applications.

\setbox1=\hbox to 2.05in{\vbox to 3.1in{\ninepointsss
  \hsize=2.05in
  \kern 2.7in
  \includegraphics{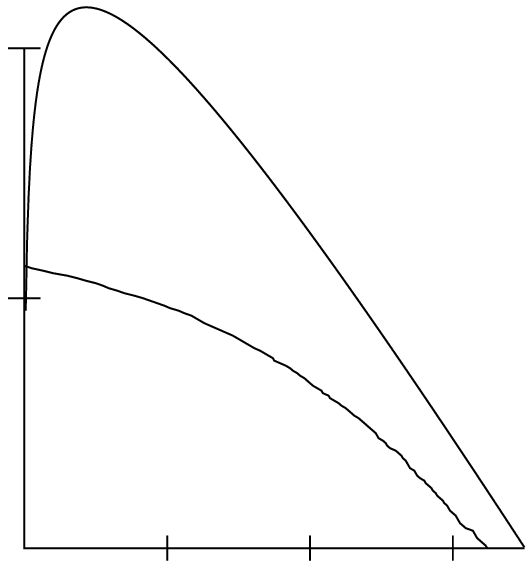}
  \anote{.48}{-.2}{$10^1$}
  \anote{1.05}{-.2}{$10^2$}
  \anote{1.62}{-.2}{$10^3$}
  \anote{-.3}{1}{$10^{-2}$}
  \anote{-.3}{2}{$10^{-1}$}
  \vfill
  \centerline{$n=10$, $a=-1$, errors Cauchy}
}}
\setbox2=\hbox to 2.05in{\vbox to 3.1in{\ninepointsss
  \hsize=2.05in
  \kern 2.7in
  \includegraphics{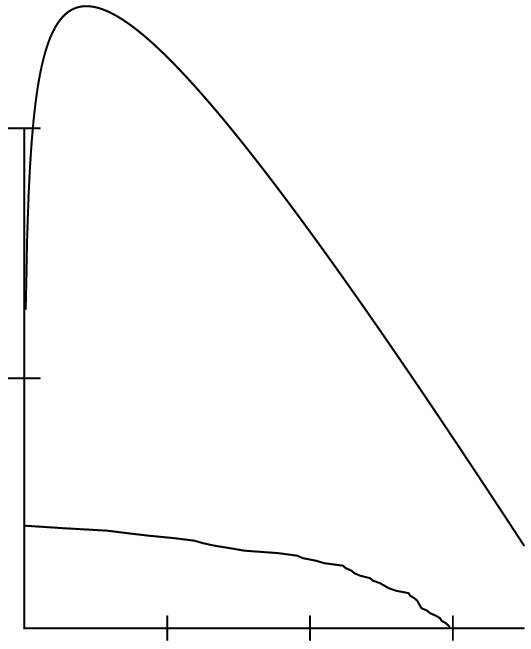}
  \anote{.48}{-.2}{$10^1$}
  \anote{1.05}{-.2}{$10^2$}
  \anote{1.62}{-.2}{$10^3$}
  \anote{-.3}{1}{$10^{-2}$}
  \anote{-.3}{2}{$10^{-1}$}
  \vfill
  \centerline{$n=20$, $a=-1$, errors Cauchy}
}}
\line{\kern 0.3in\box1\hfill\box2}

\bigskip

When $a$ is negative, it makes more sense to approximate the lower tail.
Given what we have done, it is a trivial matter. We state a
result only in the form needed for our discussion.

\bigskip

\state{11.2.2. THEOREM.} {\it%
  Let $X$ be an autoregressive process of order $1$, with
  coefficient $a$ and errors independent and identically
  distributed from a Student-like distribution with
  parameter $\alpha$. If $a$ is negative, then
  $$
    P\{\, n\gamma_n(1)\leq -t \,\}
    \sim K_{s,\alpha}\alpha^{(\alpha-1)/2}2 a^{\alpha/2}
    \sum_{1\leq i\leq n-1} 
    \Big({1-a^{2j}\over 1-a^2}\Big)^{\alpha/2} t^{-\alpha/2}
  $$
  as $t$ tends to infinity.
  }

\bigskip

\stateit{Proof.} With the notation of the proof of Theorem
11.2.1, we need to evaluate the upper tail of $\<-C\epsilon ,
\epsilon\>$. In the proof of Theorem 11.2.1, we shown that when
$k=1$ and $a$ is negative, the matrix $-C$ has its coefficients $C_{n,n}$
vanishing, while $C_{i,i}=a(1-a^{2(n-i)})/(1-a^2)$ for
$i=1,\ldots n-1$. Apply Theorem 8.2.1 to conclude the proof.\hfill$\qed$

\bigskip

There is not much point in reproducing here results on the
approximation of the lower tail. It is enough to say that it
works as expected; that is, the approximation is very sharp when
$a=-1$ or $a=-0.5$. The pictures look identical to those for the
upper tail with positive coefficient $a$.

\bigskip

As $\alpha$ increases, the approximation given in
Theorem 11.2.1 degenerates for
positive values of $a$. For a Student-distribution with
5 degrees of freedom, that is $\alpha=5$,  and $n=10$, their use
starts to be questionable.
\bigskip

\setbox1=\hbox to 2.05in{\vbox to 2.5in{\ninepointsss
  \hsize=2.05in
  \kern 2.1in
  \includegraphics{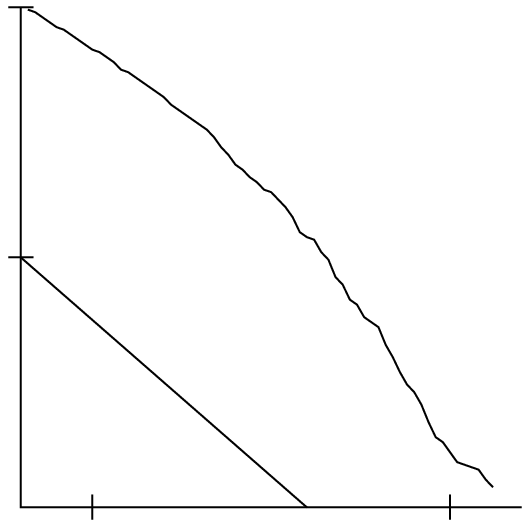}
  \anote{.2}{-.2}{$10^3$}
  \anote{1.6}{-.2}{$10^{3.5}$}
  \anote{-.3}{1}{$10^{-2}$}
  \anote{-.3}{2}{$10^{-1}$}
  \vfill
  \centerline{$n=10$, $a=1$, errors Student(5)}
}}
\setbox2=\hbox to 2.3in{\vbox to 2.5in{\ninepointsss
  \hsize=2.3in
  \kern 2.1in
  \includegraphics{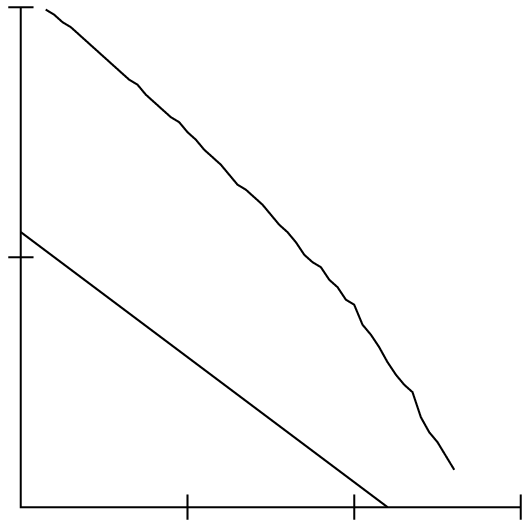}
  \anote{.6}{-.2}{$10^2$}
  \anote{1.25}{-.2}{$10^{2.2}$}
  \anote{1.9}{-.2}{$10^{2.4}$}
  \anote{-.3}{1}{$10^{-2}$}
  \anote{-.3}{2}{$10^{-1}$}
  \vfill
  \centerline{$n=20$, $a=1$, errors Student(5)}
}}
\wd2=155pt

\line{\kern 0.3in\box1\hfill\box2}
%\line{\box1\hfill\box2}

\bigskip

For $a=0$, we appealed to Theorem 8.3.1, and the approximation is not
so good. We would like to point out that this failure can be
seen on a very simple example. Consider two independent random
variables, $X$, $Y$, with density $\alpha / x^{\alpha+1}$ over
$[\, 1,\infty )$. We can calculate explicitly the distribution of
their product,
$$\eqalign{
  P\{\, XY\geq t\,\}
  & = \alpha\int_{x\geq 1} x^{-\alpha-1} P\{\, Y\geq t/x\,\} 
    dx \cr
  & = \alpha\int_{x\geq 1}x^{-\alpha-1}
    \big( (x/t)^\alpha\wedge 1\big) dx \cr
  & = {\alpha\log t\over t^\alpha} + {1\over t^\alpha} \, . \cr
  }
$$
Our approximation picks up the leading term, $\alpha
t^{-\alpha}\log t$ as $t$ tends to infinity. But the growth of
$\log t$ is too slow for the first term to really dominate in the
range where the probability is of order $10^{-1}$ or $10^{-2}$.
One cannot expect a good one-term approximation in such case,
except if $\alpha$ is large.
This suggests that when $a$ is zero, our approximation may improve
when $\alpha$ increases. This is the case. And of course, for
negative $a$, it becomes worse, positive values of the
empirical covariance being less and less likely.

\bigskip

\setbox1=\hbox to 2.2in{\vbox to 2.5in{\ninepointsss
  \hsize=2.2in
  \kern 2.1in
  \includegraphics{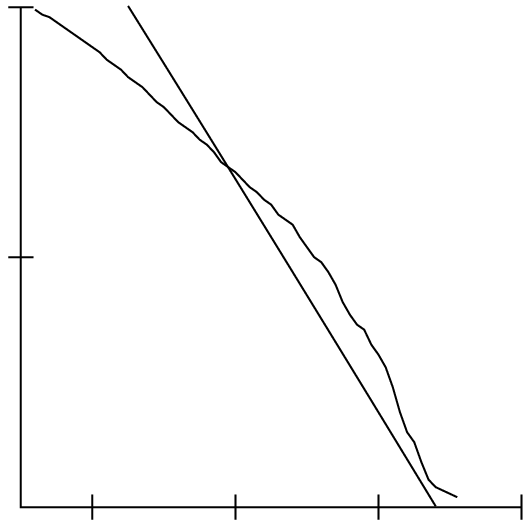}
  \anote{.2}{-.2}{$10^1$}
  \anote{0.8}{-.2}{$10^{1.2}$}
  \anote{1.4}{-.2}{$10^{1.4}$}
  \anote{1.9}{-.2}{$10^{1.6}$}
  \anote{-.3}{1}{$10^{-2}$}
  \anote{-.3}{2}{$10^{-1}$}
  \vfill
  \centerline{$n=20$, $a=0$, errors Student(5)}
}}
\setbox2=\hbox to 2.3in{\vbox to 2.5in{\ninepointsss
  \hsize=2.3in
  \kern 2.1in
  \includegraphics{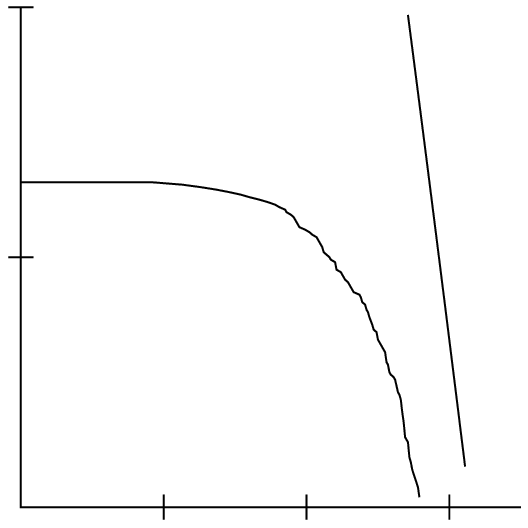}
  \anote{.5}{-.2}{$10^{-1}$}
  \anote{1.1}{-.2}{$1$}
  \anote{1.65}{-.2}{$10$}
  \anote{-.3}{1}{$10^{-2}$}
  \anote{-.3}{2}{$10^{-1}$}
  \vfill
  \centerline{$n=20$, $a=-0.5$, errors Student(5)}
}}
\wd2=150pt

\line{\kern 0.3in\box1\hss\box2}

\bigskip

For $a=0$, the approximation is fairly good, until we have
many moments on the distribution, that is if $\alpha$ is large
enough. Then assuming simply that the errors are normally
distributed may give a better approximation.

\bigskip

\setbox1=\hbox to 2.3in{\vbox to 2.5in{\ninepointsss
  \hsize=2.3in
  \kern 2.1in
  \includegraphics{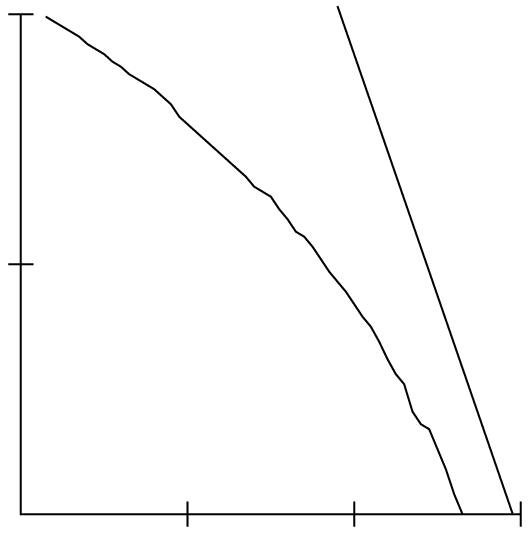}
  \anote{.6}{-.2}{$10^1$}
  \anote{1.3}{-.2}{$10^{1.2}$}
  \anote{1.65}{-.2}{$10^{1.4}$}
  \anote{-.3}{1}{$10^{-2}$}
  \anote{-.3}{2}{$10^{-1}$}
  \vfill
  \centerline{$n=20$, $a=0$, errors Student(10)}
}}

\line{\hfill\box1\hfill}

\bigskip

When our approximation is not so good, one can think of some
alternative techniques. Besides the classical Edgeworth expansion ---
which is poor in term of relative error --- one could also approximate
the Student distribution by a normal one, and then proceed as if
$\epsilon$ were normally distributed. This works well if the normal
distribution has the ``right'' variance. But one has to be aware that
for $\alpha$ small, the right variance {\it is not} that of the
corresponding Student distribution. For instance, for the Student
distribution with 5 degrees of freedom, the normal approximation using
a variance equal to that of the Student distribution poor. One needs a
much larger variance. I tried to approximate the Student distribution
by a normal with the variance such that some quantile of the normal
would be equal to that of the Student. This does not work any better,
in the sense that there is no systematic way to do this kind of
calibration. One should also be aware that the symmetry in the Student
distribution makes our approximations less precise. They would be more
accurate if the errors had a centered Pareto distribution for
instance. Obviously more work is needed to derive a set of
approximations which would cover more or less any regime. It is
doubtful that a single approximation scheme can give satisfactory
results under a very broad class of distributions for the errors and
relatively arbitrary sample size.

\bigskip

Classically, $a$ is estimated by $\hat
a=\gamma_n(1)/\gamma_n(0)$. Let us now consider the test problem

\noindent\qquad $\hbox{H}_0 \, :\,  a\leq a_0$,
         \qquad versus \qquad $\hbox{H}_1 \, :\, a> a_0$.

A possible way to perform this test is to reject the null
hypothesis if $\hat a_n -a_0$ is too large, that is if
$\gamma_n(1)-a_0\gamma_n(0)$ is too large. For a reason which
will be explained in the proof of the next result, it is better 
to use\ntn{$\hat\gamma_n(0)$}
$$
  \hat\gamma_n(0)=n^{-1}\sum_{1\leq i\leq n-1} X_i^2
$$
instead of $\gamma_n(0)$. As $n$ tends to infinity, the classical
theory ensures that it does not make any difference. However, it
does make a difference for our finite sample results. The theory
is much nicer with $\hat\gamma_n (0)$.

The test statistics is again a
quadratic form in $\epsilon$. The following result gives its
tail approximation under the null as well as under the
alternative hypothesis. As one more parameter is involved,
another behavior appears.

\bigskip

\state{11.2.2. THEOREM.} {\it%
  Let $X$ be an autoregressive process of order $1$, with 
  coefficient $a$ and
  independent and identically distributed innovations, all having
  a Student-like distribution with parameter $\alpha$. 
  The tail probability 
  $$
    P\big\{\, n\big(\gamma_n(1)-a_0\hat\gamma_n(0)\big)
    \geq t \,\big\}
  $$
  admits the following equivalent as $t$ tends to infinity.
  
  \noindent (i) If $a>a_0$,
  $$
    K_{s,\alpha}\alpha^{\alpha-1\over 2}
    2(a-a_0)^{\alpha/2}
    \sum_{1\leq i\leq n-1}
    \Big({\displaystyle 1-a^{2(n-i)}\over\displaystyle 1-a^2}
    \Big)^{\alpha/2}\, t^{-\alpha/2} \, .
  $$

  \noindent (ii) If $a_0=a$ and both are nonnegative,
  $$
    K_{s,\alpha}^2\alpha^\alpha 2
    \sum_{1\leq k\leq n-2} (n-k) a^{(k-1)\alpha} \,
    t^{-\alpha}\log t \, ,
  $$
  with the convention that $0^0=1$ when $a=0$.

  \noindent (iii) If $a_0$ is positive and $a<a_0$,
  $$
    c(a_0,a,\alpha,n)\, t^{-\alpha}
  $$ 
  for some function $c(\cdot )$.
  }

\bigskip

In case (iii), we will explain after the proof how to calculate
the function $c(\cdot )$ in a typical case.

\bigskip

\stateit{Proof.} Write $n\big(\gamma_n(1)-a_0\hat\gamma_n(0)\big)
=\epsilon^\T C\epsilon$ with\ntn{$C$}
$$
  C =A^\T BA-a_0A^\T B^\T BA \, .
$$
From the proof of Theorem 11.2.1 we obtain the diagonal terms
$(A^\T BA)_{i,i}$. We calculate
$$
  (A^\T B^\T BA)_{i,i}
  = \sum_{1\leq k\leq n} (BA)_{k,i}^2 
  = \sum_{1\leq k\leq n-1} A_{k,i}^2
  = {1-a^{2(n-i)}\over 1-a^2 } \, .
$$
Consequently,
$$
  C_{i,i}=
  \cases{ (a-a_0){\displaystyle 1-a^{2(n-i)}
                  \over\displaystyle 1-a^2} & if $1\leq i\leq n-1$
						\cr
          0				    & if $i=n$. \cr}
$$

If instead of using $\hat\gamma_n(0)$ we use $\gamma_n(0)$,
the term $(A^\T B^\T BA)_{i,i}$ has denominator $1-a^{2(n-i+1)}$.
The discussion is a bit more involved. The result becomes more
dependent on $n$. It is more complicated to state, but it does
not make much difference as far as the theory goes.

If $a-a_0$ is positive, then all the diagonal coefficients of $C$
but $C_{n,n}$ are positive. We apply Theorem 8.3.1.

If $a-a_0$ is negative, only $C_{n,n}$ is nonnegative. We apply
Theorem 8.3.1. Since we will need it, let us calculate $C_{i,j}$.
First
$$
  (BA)_{k,i}=\cases{ A_{k-1,i} & if $2\leq k\leq n$, \cr 
                     0         & if $k=1$. \cr }
$$
Therefore,
$$
  (BA)_{k,i}=\cases{ a^{k-1-i} & if $2\leq k\leq n$, \cr
                     0         & otherwise. \cr}
$$
Consequently, a little algebra shows that
$$
  (A^\T B^\T BA)_{i,j}
  = \sum_{1\leq k\leq n}(BA)_{k,i}(BA)_{k,j}
  = a^{|i-j|}
    {\displaystyle a^{2n}-a^{2(i\vee j)}\over\displaystyle a^2-1}
  \, .
$$
In particular, for $i<n$,
$$
  C_{n,i}=a^{n-i-1},
$$
while $C_{i,n}=0$ for $1\leq i\leq n$. Therefore
$$
  \sum_{1\leq i\leq n}|C_{i,n}+C_{n,i}|^\alpha
  = \sum_{1\leq i\leq n} a^{\alpha (n-i-1)}
  = \sum_{0\leq i\leq n-2} a^{\alpha^i} \, .
$$

If $a$ and $a_0$ are equal, then all the diagonal coefficients of
$C$ vanish. In this case, for $i$, $j$ distinct,
$$
  C_{i,j}=a^{|i-j-1|}
  {\displaystyle 1-a^{2(n+1-(i\vee (j+1)))}\over 
   \displaystyle 1-a^2}
  - a_0 a^{|i-j|}
  {\displaystyle a^{2n}-a^{2(i\vee j)}\over
   \displaystyle a^2-1} \, .
$$

If $a_0$ is positive and strictly larger than $a$, then all the
diagonal coefficients of $C$ are negative. We need to determine $N(C)$
and apply Theorem 8.2.10.  The statement follows from Theorem
8.2.21.\hfill$\qed$

\bigskip

Some useful information can be deduced from Theorem 11.2.2. A first
qualitative deduction is that if $a_0$ is positive under the null hypothesis,
then the tail probability under consideration has a very different
decay according to the position of $a$ with respect to $a_0$. Thus,
one should probably not use symmetric confidence intervals. It may be
wise to have $a$ somewhere on the right half of the interval.

Next, assume that we want to test with the risk of first type $\eta$
very small.  We can use the approximation under the null hypothesis to
obtain the critical value. Define
$$
  c(a)=2K_{s,\alpha} \alpha^{\alpha-1\over 2}
  (a-a_0)^{\alpha/2}
  \sum_{1\leq i\leq n-1}
  \Big( { 1-a^{2(n-i)}\over 1-a^2 }\Big)^{\alpha/2} \, .
$$
We see that if $0\leq a_0<a$, then $\eta \sim c(a)/t^{\alpha/2}$. 
This gives an approximate critical value 
$t_\eta = \big(c(a)/\eta\big)^{2/\alpha}$. 

The following plot shows the actual risk of the first type when
using the approximate critical value. The ``true'' value of the
risk is obtained by simulation, replicating 100,000 copies of
the vector $\epsilon$. The curves were obtained by linearly
interpolating between the following values for $a$.

\bigskip
{\eightrm\baselineskip=8pt
{
\halign{\hskip 15pt#\hfil\quad &#\hfil\quad &#\hfil\quad &#\hfil\quad &
        #\hfil\quad &#\hfil\quad &#\hfil\quad &#\hfil\quad & 
        #\hfil\quad &#\quad \cr 
0.500 &
0.552 &
0.605 &
0.657 &
0.710 &
0.762 &
0.815 &
0.868 &
0.920 &
0.930 \cr
0.939 &
0.948 &
0.957 &
0.966 &
0.975 &
0.984 &
0.993 &
1.000 &
1.002 &
1.011 \cr
1.019 &
1.028 &
1.037 &
1.046 &
1.055 &
1.064 &
1.073 &
1.082 &
1.091 &
1.100 \cr
1.150 &
1.200 &
1.250 &
1.300 &
1.350 &
1.400 &
1.450 &
1.500 & & \cr
}
}}

\bigskip

\noindent The sample size is $n=20$. The lowest curve is for
$\alpha=1$, the two almost equal curves are for $\alpha=5$ --- the
lower one of the two curves --- and $\alpha=10$. The value $a=1$ is
indicated, as well as the levels 5\% and 10\% .

\setbox1=\hbox to 3.2in{\vbox to 2.3in{\ninepointsss
  \hsize=3in
  \kern 1.9in
  \includegraphics{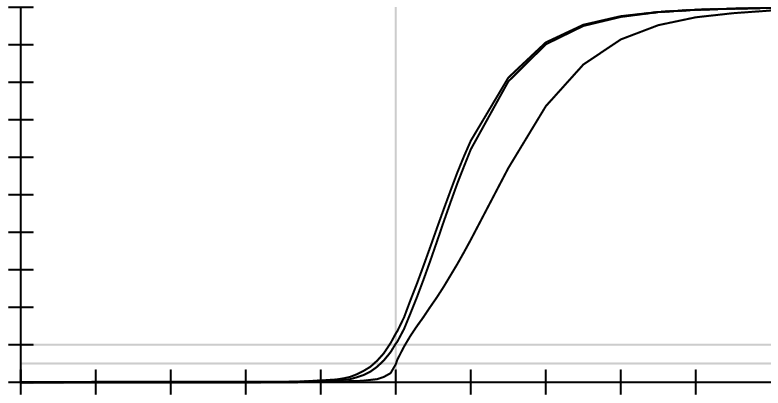}
  \anote{-.25}{.15}{0.1}
  \anote{-.25}{.75}{0.5}
  \anote{-.25}{1.50}{1}
  \anote{.22}{-.15}{0.6}
  \anote{.52}{-.15}{0.7}
  \anote{1.47}{-.15}{1}
  \anote{3.1}{0}{$a$}
  \vfill
  \centerline{$n=20$, $\alpha=1$, $5$, $10$}
}\hfill}
\wd2=155pt

\centerline{\box1}

\medskip

The result is satisfying at first glance. A more careful
examination shows some reason to worry. For $a=1$, the actual risk of
the test is about 10\% when $\alpha=5$ or $10$. This is still
small, but in terms of relative error, this is twice as much as
what we wanted. The fact that the power function grows
moderately fast with $a$ is not a surprise if you plot some
of these processes. On a trajectory of length 20, and with
errors having a Cauchy distribution, an autoregressive process
of order one with $a=1$ looks very similar to one with $a=1.2$
for instance.

\bigskip

\noindent{\fsection 11.3. Autoregressive processes of arbitrary
          order.}
%\chapternumber={11.3}
\section{11.3}{Autoregressive processes of arbitrary
               order}

\bigskip

\noindent In principle, all the results of the previous section 
can be generalized to autoregressive processes of arbitrary order.  As
more parameters are involved in the model, the analysis is
harder. Thus, our goal is rather modest. We will prove that the tail
of the autocovariance of order 1 has different decays according to the
values of the parameters, and the number of observations as
well. Since autoregressive processes of order one are a degenerate
case of those of higher order, the higher order autocovariances would
have an even more complex behavior. For any given value of the
parameters, the results of chapter 8 can be used to do numerical
computations; this is quite easy, and sometimes helpful, but does not
provide further insights.  For autoregressive processes of order 2, we
will obtain some rather precise results, showing the intricacy of
these models.

Let us now consider an autoregressive model of order $p$ as in
(11.1.1). The parameter $\theta=(\theta_1,\ldots , \theta_p)$ is in
$\RR^p$. Our first result shows that $\RR^p$ can be partitioned into
two regions, one where the tail behavior of $n\gamma_n(1)$ is like
$t^{-\alpha/2}$, the other one where it is like $t^{-\alpha}\log t$. The
noticeable fact is that there cannot be other tail behavior, and these
regions are nested when the number of observations varies.

\bigskip

\state{11.3.1. THEOREM.} {\it%
  Consider an autoregressive process of order $p$ with errors
  having a Student-like distribution with parameter $\alpha$.
  There exist nonempty semialgebraic sets\ntn{$R_k$} $R_k$, $k\geq 1$, of
  $\RR^p$, and a positive function $c(\cdot )$ on $\RR^p$, 
  such that
  $$
    P\{\, n\gamma_n(1)\geq t\,\} \sim
    \cases{ c(\theta )t^{-\alpha/2} 
		& if $\theta\in\bigcup_{1\leq k\leq n-1} R_k$, \cr
	    c(\theta )t^{-\alpha}\log t
		& if $\theta\not\in\bigcup_{1\leq k\leq n-1} R_k$.
	  }
  $$
  }

\bigskip

\state{REMARK.} The regions $R_i$ depend of course on the
dimension $p$ of the parameter space. It is remarkable that they
do not depend on $n$, though the quadratic form representing
$n\gamma_n(1)$ does. As $n$ increases, $\bigcup_{1\leq k\leq
n}R_k$ increases, but we will see during the proof that its limit
is a proper nonempty subset of $\RR^p$. The proof also gives
some indications on how to calculate these regions efficiently.

\bigskip

\stateit{Proof of Theorem 11.3.1.} We write $X=A\epsilon$. The
matrix $A$ is more involved than in section 2. Let $e_i$ be the
$i$-th vector from the canonical basis of $\RR^d$, that is
having all its entries vanishing, except the $i$-th one being
$1$. We denote $A_{i,\sbullet}$ the $i$-th row of $A$. Relation
(11.1.1) gives
$$\eqalignno{
  A_{1,\sbullet}   & = e_1 \cr
  A_{2,\sbullet}   & = \theta_1 A_{1,\sbullet}+ e_2 \cr
  \vdots	       & \cr
  A_{p+1,\sbullet} & = \theta_1 A_{p,\sbullet}
    +\theta_2 A_{p-1,\sbullet}+\cdots 
    +\theta_pA_{1,\sbullet} + e_p 
  &(11.3.1)\cr
  }
$$
and for $p+1\leq k\leq n$,
$$
   A_{k,\sbullet} = \theta_1 A_{k-1,\sbullet}
    +\theta_2 A_{k-2,\sbullet}+\cdots 
    +\theta_pA_{k-p,\sbullet} + e_k \, .
  \eqno{(11.3.2)}
$$

Let\ntn{$C$} $C=A^\T BA$ be the matrix of the quadratic form such that
$n\gamma_n(1)=\epsilon^\T C\epsilon$. Since $A_{j,i}=0$ for any
$1\leq j<i\leq n$,
$$
  C_{n,n}=0 \qquad \hbox{ and } \qquad \max_{1\leq i\leq n}
  C_{i,i} \geq 0 \, .
$$
Thus, Theorems 8.2.1 and 8.3.1 show that the only possible tail
behaviors of $n\gamma_n(1)$ are either like $t^{-\alpha/2}$ or like
$t^{-\alpha}\log t$. They show as well the existence of the
nonvanishing function $c(\cdot )$.

Since $A_{i,j}=A_{i+1,j+1}$ for all $i,j\leq n-1$, we have for
$k\geq 1$,
$$\eqalignno{
  C_{n-k,n-k}
  &=\sum_{1\leq i,j\leq n} (A^\T )_{n-k,i} B_{i,j}A_{j,n-k} \cr
  &= \sum_{2\leq j\leq n} A_{j,n-k} A_{j-1,n-k} \cr
  &= \sum_{n-k+1\leq j\leq n} A_{j,n-k} A_{j-1,n-k} \cr
  &= \sum_{1\leq j\leq k} A_{n-k+j,n-k} A_{n-k+j-1,n-k} \cr
  &= \sum_{1\leq j\leq k} A_{j+1,1}A_{j,1} \, .
  &(11.3.3)\cr
  }
$$
Relations (11.3.1) and (11.3.2) show that $A_{j,1}$ is a
polynomial in $a$ and $b$. Hence $C_{n-k,n-k}$ is also a
polynomial in $a$ and $b$. The region
$$
  R_k=\{\, (a,b)\in\RR^2 \, :\, C_{n-k+1,n-k+1}\geq 0\,\}
$$
is then a semialgebraic set, which does not depend on $n$.
Moreover, the largest diagonal coefficient of $C$ is positive if
and only if $(a,b)$ is in $\bigcup_{1\leq k\leq n} R_k$. The
result then follows from Theorems 8.2.1 and 8.3.1.\hfill$\qed$

\bigskip

In practice, to decide in which region we are, we can
numerically compute the diagonal terms of $C$ with formulas
(11.3.1)--(11.3.3). But if one ultimately wants the constant
$c(a,b)$, other elements of $C$ are needed when the tail is like
$t^{-\alpha}\log t$.

The sets $R_k$ do not seem easy to describe in general. For
applications, this may not be so important, since numerical
computation is easy to implement. Understanding their geometry amounts
to understanding the behavior of the roots of inductively defined
polynomials in the $p$ variables $\theta_1,\ldots, \theta_p$.  More
can be said when $p=2$, because polynomials of degree $2$ are well
understood. And this is enough to show how intricate these
autoregressive models are. Thus, from now on, we focus on
autoregressive models of order $2$. We change slightly the notation,
using\ntn{$a$, $b$} $(a,b)$ instead of $(\theta_1,\theta_2)$. Thus,
our model is
$$\eqalign{
  X_1&=\epsilon_1 \, , \cr
  X_2&=aX_1+\epsilon_2 \, ,\cr
  X_k&=aX_{k-1}+bX_{k-2}+\epsilon_k \, , \qquad 3\leq k\leq n
       \, . \cr
  }
$$
We can explicitly write down $R_1$, $R_2$, $R_3$, $R_4$.
Indeed, $C_{n,n}=0$ and thus $R_1=\emptyset$.
We then have
$$
  C_{n-1,n-1}=A_{2,1}A_{1,1}=a \, ,
$$
thus
$$
  R_2=\{\, (a,b) \, : \, a>0\,\} \, .
$$
Since $A_{3,1}A_{2,1}=(a^2+b)a$,
$$
  C_{n-2,n-2}= (a^2+b)a + a = a(a^2+b+1) \, .
$$
Consequently,
$$
  R_3=\{\, (a,b)\, :\, a>0\hbox{ and } b>-1-a^2 \, ;\, 
  \hbox{ or } a<0 \hbox{ and } b<-1-a^2\,\} \, .
$$
But what matters more,
$$
  R_2\cup R_3 =\big( (0,\infty )\times\RR\big) 
  \cup\{\, (a,b)\, :\, a<0 \, , \, b<-1-a^2\,\} \, .
$$
The next figure shows these regions.

We also obtain $A_{4,1}=a(a^2+2b)$, which leads to
$$\eqalign{
  C_{n-3,n-3}
  &= a(a^2+2b)(a^2+b) + a(a^2+b+1) \cr
  &= a\big(2b^2+b(3a^2+1)+a^4+a^2+1\big) \cr
  }
$$
To obtain $R_2\cup R_3\cup R_4$, we need to see when
$C_{n-3,n-3}$ is positive as $a$ is negative. In other words, 
when $a$ is negative and
$$
  2b^2+b(3a^2+1)+1+a^2+a^4 < 0 \, .
$$
This holds if $a<{\displaystyle 1-2\sqrt 2\over\displaystyle 2}$
and
$$
  b_-\leq b\leq b\leq b_+
$$
with
$$
  b_\pm = {-(1+3a^2)\pm\sqrt{a^4-2a^2-7}\over 4} \, .
$$
Thus,
$$\displaylines{\quad
    R_1\cup R_2\cup R_3\cup R_4
    =\big( (0,\infty )\times \RR)\bigcup 
    \{\, (a,b)\, :\, a<0 \, , \, b<-1-a^2\,\}
  \cut
    \bigcup \{\, (a,b) \, : \, a\leq {1-2\sqrt 2\over 2}
    \, ;\, b_-\leq b\leq b_+\,\} \, .
  \quad}
$$

\setbox1=\vbox to 3.2in{\ninepointsss
  \hsize=2.25in
  \kern 3in
  \includegraphics{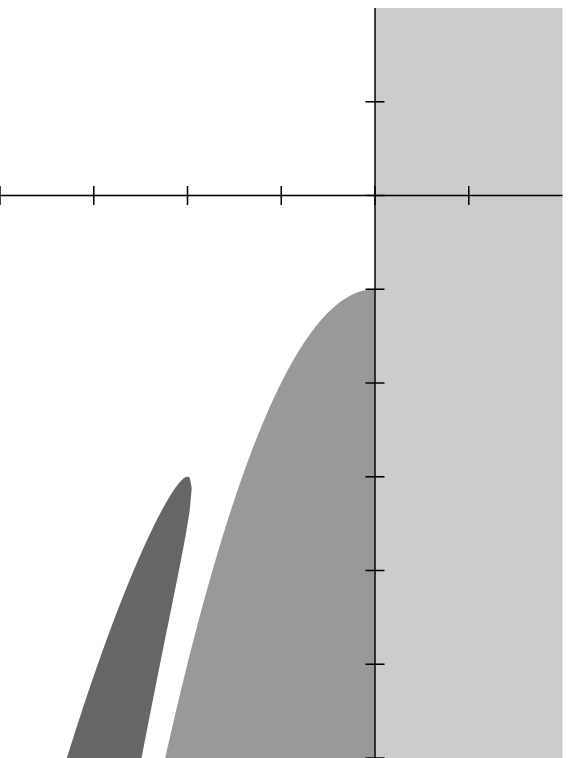}
  \anote{1.38}{2.95}{$b$}
  \anote{2.35}{2.25}{$a$}
  \anote{1.38}{2.6}{1}
  \anote{1.82}{2.1}{1}
  \vfill
  \centerline{Regions $R_2$, $R_2\cup R_3$, $R_2\cup R_3\cup R_4$.}
}
\centerline{\box1}

\bigskip

To calculate $R_5$, we need to solve a cubic equation in $b$.  Closed
form expressions are getting more and more cumbersome, and even
nonexistent. The following picture shows $\bigcup_{1\leq k\leq n} R_k$
for $n=5,6,7$ and $n=10$ in the domain $a$ negative.
\finetune{\vfill\eject}

\setbox1=\hbox{\vbox to 3.4in{\ninepointsss
  \hsize=2in
  \kern 3.2in
  \includegraphics{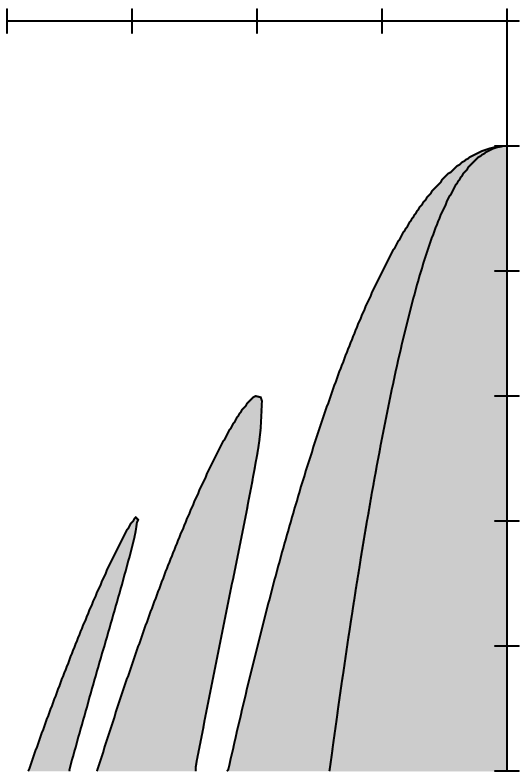}
  \anote{2.1}{3}{$a$}
  \anote{2}{3.1}{$b$}
  \anote{2.1}{2.5}{-1}
  \anote{1.45}{2.8}{-1}
  \vfill
  \centerline{$n=5$}
}}
%\wd1=2in

\setbox2=\hbox{\vbox to 3.4in{\ninepointsss
  \hsize=2in
  \kern 3.2in
  \includegraphics{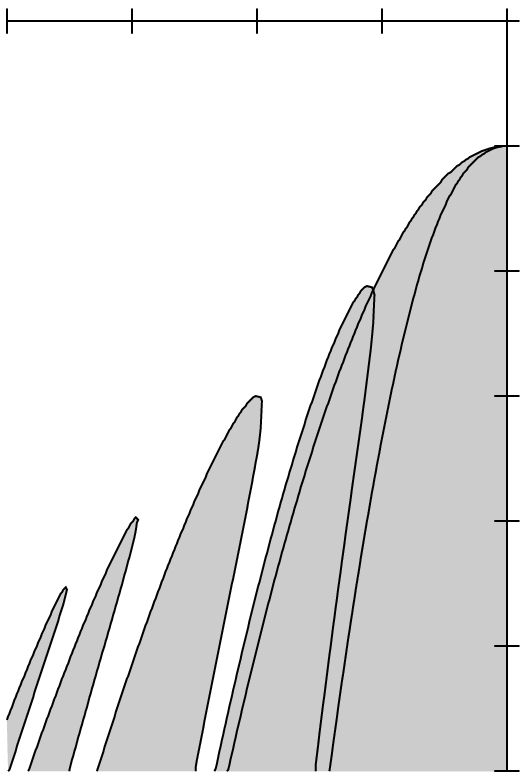}
  \anote{2.1}{3}{$a$}
  \anote{2}{3.1}{$b$}
  \anote{2.1}{2.5}{-1}
  \anote{1.45}{2.8}{-1}
  \vfill
  \centerline{$n=6$}
}}
%\wd2=2in

\line{\box1\hfill\box2}
%\line{\hbox to 2in{\box1}\hfill\hbox to 2in{\box2}}

\setbox1=\hbox{\vbox to 3.4in{\ninepointsss
  \hsize=2in
  \kern 3.2in
  \includegraphics{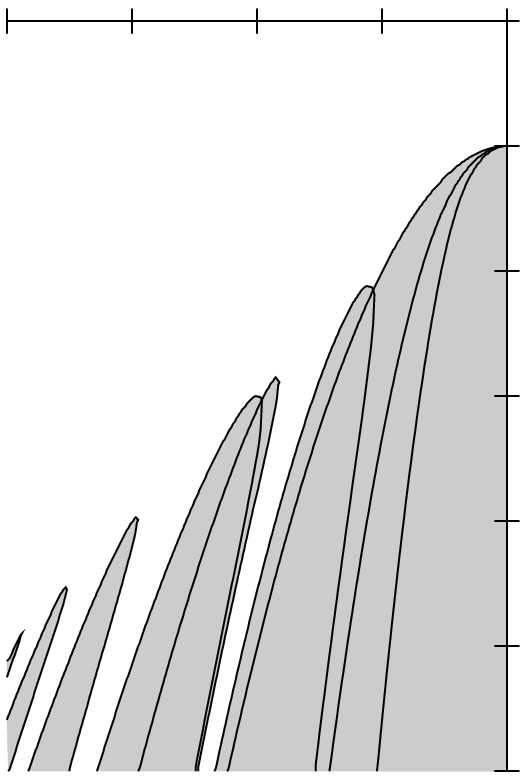}
  \anote{2.1}{3}{$a$}
  \anote{2}{3.1}{$b$}
  \anote{2.1}{2.5}{-1}
  \anote{1.45}{2.8}{-1}
  \vfill
  \centerline{$n=7$}
}}
%\wd2=2in

\setbox2=\hbox{\vbox to 3.4in{\ninepointsss
  \hsize=2in
  \kern 3.2in
  \includegraphics{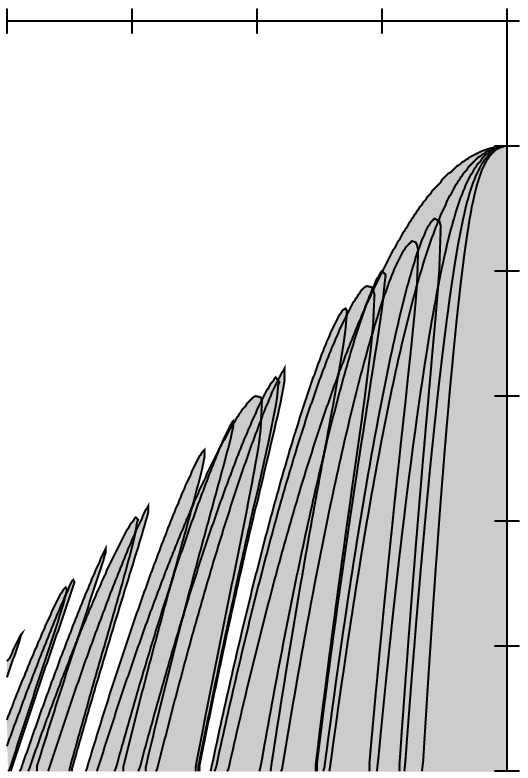}
  \anote{2.1}{3}{$a$}
  \anote{2}{3.1}{$b$}
  \anote{2.1}{2.5}{-1}
  \anote{1.45}{2.8}{-1}
  \vfill
  \centerline{$n=10$}
}}
%\wd2=2in

\line{\box1\hfill\box2}

These pictures can be more or less understood theoretically.
This is the purpose of the next result. Its proof contains even
more information, some of it being important, and we will
discuss further after the proof.

\bigskip

\state{11.3.2. THEOREM.} {\it%
  The closure of $\bigcup_{k\geq 1} R_k$ contains 
  all points $(a,b)$ for which one of the following 
  conditions holds:

  \noindent (i) $a>0$,

  \noindent (ii) $ a\leq 0$ and $b<-a^2-1$,

  \noindent (iii) $a\leq 0$ and $b<\min (-a^2/4 , -a-1)$. 
  }

\bigskip

The region described by the three condition in Theorem 11.3.2
is shaded gray in the following picture. It does not contain
its boundary.

\setbox1=\hbox to 4.2in{\vbox to 3in{\ninepointsss
  \hsize=4.2in
  \kern 3in
  \includegraphics{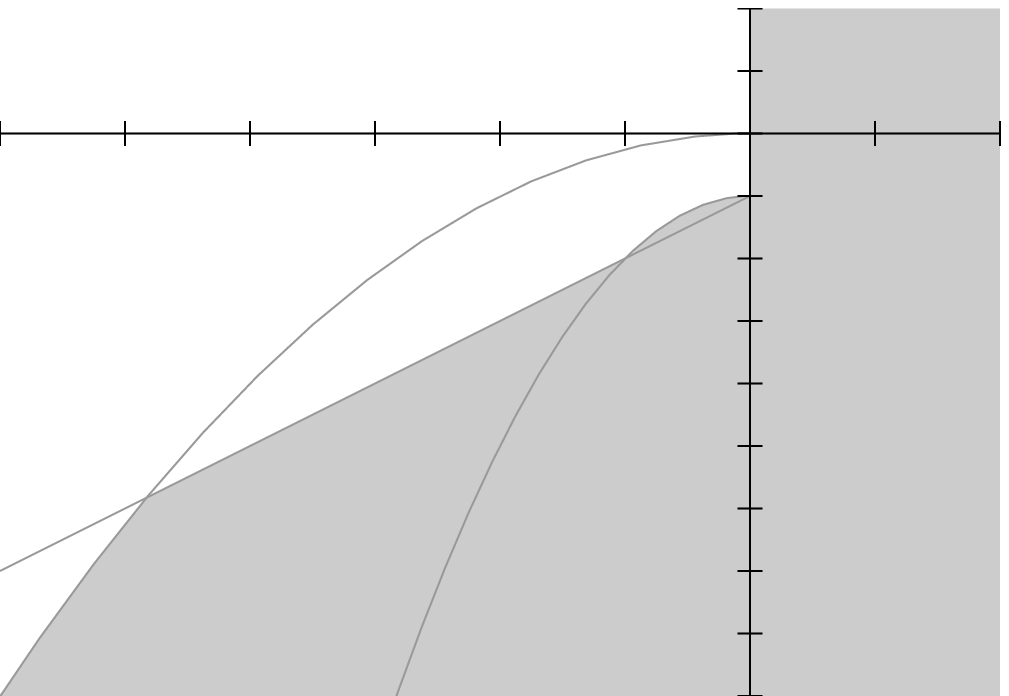}
  \anote{3.48}{2.1}{1}
  \anote{2.87}{2.48}{1}
  \anote{4.1}{2.22}{$a$}
  \anote{2.95}{2.8}{$b$}
  \anote{0}{0.72}{$b=a-1$}
  \anote{0}{-.15}{$b=-a^2/4$}
  \anote{1.5}{-.15}{$b=-a^2-1$}
  \vfill
}\hfill}

\centerline{\box1}

\bigskip

It follows from Theorems 11.3.2 and 11.3.1 that if $(a,b)$ lies
in the gray shaded region, the tail behavior of $P\{\,
n\gamma_n(1)\geq t\,\}$ is typically like $t^{-\alpha /2}$ for 
$n$ large enough. I believe that in the nonshaded domain, 
the tail behavior is like $t^{-\alpha}\log t$; we will
prove this only when $b\geq a^2/4$.

\bigskip

\stateit{Proof of Theorem 11.3.2.} The proof will be done in
examining different regions. We will need several lemmas, and
will actually prove much more than the statement.

\bigskip

\state{11.3.3. LEMMA.} {\it%
  If $a$ is positive, so is the largest diagonal coefficient of
  $C$. If $a$ is nonpositive, then the largest diagonal coefficient
  of $C$ vanishes.
  Consequently, $\bigcup_{k\geq 1}R_k$ contains the region
  $a>0$, but does not intersect the
  region $a\leq 0$ and $b\geq -a^2/4$.
  }

\setbox1=\hbox{\vbox to 3in{\ninepointsss
  \hsize=2in
  \kern 2.8in
  \includegraphics{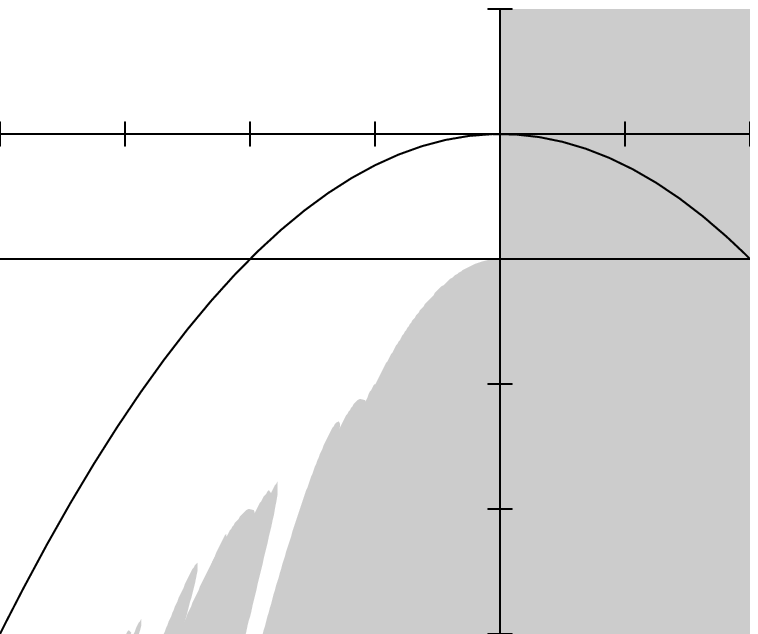}
  \anote{3.1}{2}{$a$}
  \anote{1.95}{2.55}{$b$}
  \anote{1.85}{2.45}{1}
  \anote{2.47}{2.05}{1}
  \anote{1.85}{2.05}{0}
  \anote{.5}{.72}{$b=-a^2/4$}
  \anote{.5}{1.55}{$a=-1$}
  \vfill
%  \centerline{n=10}
}}

\centerline{\box1}

\bigskip

\stateit{Proof.} If $a$ is positive, then $C_{n-1,n-1}=a$ is
positive too. To see what is happens when $a$ is nonpositive, denote
by $u,v$ the roots of the characteristic equation
$x^2-ax-b=0$. Equations (11.3.1)--(11.3.2) yield
$$
  A_{i,j}=ru^{i-j}+sv^{i-j}\, ,\qquad i\geq j \, ,
$$
with initial condition $A_{i,i}=1$ and $A_{i+1,1}=a$. Thus $r$
and $s$ are determined by
$$
  r+s=1\qquad \hbox{ and }\qquad ru+sv=a \, .
$$
If $u$ and $v$ are distinct, that is $b\neq -a^2/4$, 
$$
  A_{j,1}={u^j-v^j\over u-v} \, .
$$
If $b=-a^2/4$, we find
$$
  A_{j,1}=j(a/2)^{j-1} \, .
$$
Consequently, if $b\neq -a^2/4$,
$$
  A_{j+1,1}A_{j,1}={(u^{j+1}-v^{j+1})(u^j-v^j)\over (u-v)^2}
  \, .
  \eqno{(11.3.4)}
$$
Notice that
$$
  (u-v)^2 = (u+v)^2-4uv = a^2+4b \, .
$$

Let us now assume that $a$ is negative and $b$ is positive.  Then
$a^2+4b$ is positive. There is no loss of generality in assuming
$u<0<v$ since the product of the roots, $-b$, is negative. The
inequality
$$
  u^2-v^2=(u+v)(u-v) = -a \sqrt{a^2+4b} \geq 0
$$
forces $|u|\geq |v|$. Consequently, the sign of
$A_{j+1,1}A_{j,1}$ is that of
$$
  \big( (-1)^{j+1}|u|^{j+1} - v^{j+1}\big)
  \big( (-1)^j|u|^j-v^j\big) \, ,
$$
which is negative.
Therefore, the sequence $k\mapsto C_{n-k,n-k}$ is decreasing and
$$
  \max_{1\leq k\leq n} C_{k,k} = C_{n,n} = 0 \, .
$$

If we now assume that $a$ is negative and $-a^2/4<b\leq 0$, the roots
$u,v$ are still real, but both are negative. Thus
$$
  A_{j+1,1}A_{j,1}
  = -{(|u|^{j+1}-|v|^{j+1})(|u|^j-|v|^j)\over a^2+4b}
$$
is nonpositive since the function $x\mapsto x^j$ and $x\mapsto
x^{j+1}$ are increasing on $\RR^+$. The sequence $k\mapsto
C_{n-k,n-k}$ is decreasing and its maximum is $C_{n,n}=0$. This
gives Lemma 11.3.3.\hfill$\qed$

\bigskip

The region left is $a\leq 0$ and $b<-a^2/4$. Our next lemma
covers a part of it.

\bigskip

\state{11.3.4. LEMMA.} {\it%
  If $a\leq 0$ and $b\geq -1$, then the largest diagonal coefficient
  of $C$ vanishes. Therefore, $\bigcup_{k\geq 1}R_k$ does not
  intersect the region $a<0$ and $b\geq -1$.
  }

\bigskip

\stateit{Proof.} If $a^2+4b$ is positive, the result follows from
Lemma 11.3.3. If $a^2+4b$ is negative, equation (11.3.1)--(11.3.3)
gives
$$\eqalign{
  C_{n-k,n-k}
  & = A_{2,1}A_{1,1} + A_{3,1}A_{2,1} 
    + \sum_{3\leq j\leq k} (aA_{j,1}+bA_{j-1,1})A_j \cr
  & = b \, C_{n-k+1,n-k+1} + a\sum_{1\leq j\leq k}A_j^2 \, . \cr
  }
$$
Consequently, for $k\geq 2$,
$$
  C_{n-k,n-k}=b^2 C_{n-k+2,n-k+2} 
  + a(b+1)\sum_{1\leq j\leq k-1} A_j^2 + a A_k^2 \, .
$$
Thus, if $b\geq -1$ and $a\leq 0$,
$$
  C_{n-k,n-k}\leq b^2 C_{n-k+2,n-k+2} \, .
$$
Since $C_{n,n}=0$ and $C_{n-1,n-1}=a<0$ on the given range,
this shows that $C_{n-k,n-k}\leq 0$ for all $k\geq 0$.
\hfill$\qed$

\bigskip

To study the domain $a\leq 0$, $b\leq -1$ and $b<-a^2/4$ is
much more complicated. We assume from now on, and until the end
of the proof of Theorem 11.3.2, that $(a,b)$ is in
this domain. It is then convenient to make a change of
parameterization, setting\ntn{$\phi$}\ntn{$r$}
$$
  a=-2r\cos\phi \, , \qquad 
  b=-r^2 \, ,\qquad
  r\geq 0 \, ,\quad 0\leq\phi\leq \pi/2 \, .
$$
First, this allows us to obtain a closed formula for the
diagonal coefficients of $C$.

\bigskip

\state{13.3.5. LEMMA.} {\it%
  If $a=-2r\cos\phi$ and $b=-r^2$, with $r$ nonnegative and 
  $\phi$ in $[\, 0,\pi/2 \, ]$, then
  $$\displaylines{
    C_{n-k+1,n-k+1} =
      r\,\bigg[ -2(1-r^{2k})\cos\phi\sin^2\phi 
    \cut
      +r^{2(k-1)}(1-r^2)\sin (k\phi )
      \Big( \sin \big((k+1)\phi\big) 
            - r^2\sin\big((k-1)\phi\big)\Big)
      \bigg]\bigg/
     \cr\hfill
      (1-r^2)\big( (1-r^2)^2+4r^2\sin^2\phi\big)\sin^2\phi
    \, .\qquad \cr}
  $$
  }

\bigskip

\stateit{Proof.} Write $u=re^{i\theta}$ and $v=re^{-i\theta}$ for
the roots of the characteristic equation $x^2-ax-b=0$. Setting
$\phi=\pi-\theta$, we obtain
$$\eqalign{
  a &=u+v = 2r\cos \theta=-2r\cos\phi \cr
  b &= -uv = -r^2 \, .\cr
  }
$$
This is the origin of the parameterization. Equations (11.3.3)
and (11.3.4) give for $k\geq 2$,
$$\eqalign{
  C_{n-k+1,n-k+1}
  & = \sum_{1\leq j\leq k-1} 
    { u^{2j+1}-(u+v)(uv)^j+v^{2j+1}\over (u-v)^2} \cr
  & = {1\over (u-v)^2} \Big( u^3 {1-u^{2(k-1)}\over 1-u^2}
    - (u+v)uv{1-(uv)^{k-1}\over 1-uv} \cr
  & \kern 2.3in
    + v^3{1-v^{2(k-1)}\over 1-v^2} \Big) \cr
  }
$$
In this last expression, a part is independent of $k$. It is the
ratio of
$$\displaylines{\qquad
    u^3(1-uv)(1-v^2)-(u+v)uv(1-u^2)(1-v^2)+v^3(1-u^2)(1-uv)
  \cut
    = (u-v)^2(u+v)
    = -4r^2\sin^2\theta\, 2r\cos\theta
  \qquad\qquad\cr
  }
$$
and
$$\displaylines{\qquad
    (u-v)^2(1-u^2)(1-v^2)(1-uv)
  \hfill\cr\noalign{\vskip .05in}\hfill
    \eqalign{=& \, (u-v)^2(1-uv)\big( -(u+v)^2+(1+uv)^2\big) \cr
             =& -(2r\sin\theta )^2 (1-r^2)
                \big((1+r^2)^2-4r^2\cos^2\theta\big) \cr
             =& -4r^2\sin^2\theta \, (1-r^2)
                \big( (1-r^2)^2+4r^2\sin^2\theta\big) \, . \cr}
    \qquad\eqalign{ & \cr & \cr (11.3.5)\cr}
  }
$$
For the part dependent on $k$, we reduce it to the same
denominator, (11.3.5), and obtain the numerator
$$\displaylines{\quad
   (1-uv)\big( -u^{2k+1}-v^{2k+1}+u^2v^2 (u^{2k-1}+v^{2k-1})
     \big)
   \cut
     \hbox to 165pt{$+ (u+v)(uv)^k(1-u^2)(1-v^2)$\hfill}
   \quad\cr\noalign{\vskip .05in}\qquad
     = -(1-r^2)\big( 2r^{2k+1}\cos (2k+1)\theta 
     -r^4 2r^{2k-1}\cos (2k-1)\theta \big)
   \cut\hbox to 165pt{$+ 2r\cos\theta \, r^{2k} 
     \big( (1-r^2)+4r^2\sin^2\theta\big) \, .$\hfil}
   \quad\cr
  }
$$
Adding the part independent of $k$ and that dependent of $k$, we
obtain the numerator
$$\displaylines{
   -8 (r^3-r^{2k+3})\sin^2\theta\cos\theta
   -2r^{2k+1}(1-r^2)
     \big( \cos(2k+1)\theta - r^2\cos(2k-1)\theta
   \cut
     - (1-r^2)\cos\theta\big)
   \cr\hfill
  = - 8(r^3-r^{2k+3})\sin^2\theta\cos\theta 
    + 4r^{2k+1}(1-r^2)\sin k\theta \big( \sin (k+1)\theta
   \cut -r^2\sin^2 (k-1)\theta\big) \, .
   \cr  }
$$
Consequently,
$$\displaylines{
  C_{n-k+1,n-k+1}= r 
  \bigg[
    2(1-r^{2k})\sin^2\theta\cos\theta 
  \cut
    - r^{2(k-1)}(1-r^2)
    \sin k\theta 
    \big( \sin (k+1)\theta - r^2\sin (k-1)\theta\big)
    \bigg]\bigg/
  \cr\hfill
    (1-r^2)\big( (1-r^2)^2+4r^2\sin^2\theta\big) \sin^2\theta
  \, .
  \qquad\cr}
$$
The change of angle $\theta=\pi-\phi$ gives the
result.\hfill$\qed$

\bigskip

In the domain $a\leq 0$ and $b<-1$ with $b^2<-a^2/4$, that is
$r>1$ in our $(r,\phi)$-parameterization, Lemma 11.3.5 shows
that the sign of $C_{n-k+1,n-k+1}$ is that of minus its
numerator. Thus, it has the same sign as
$$\displaylines{\qquad
  2(1-r^{2k})\cos\phi\sin^2\phi 
  \hfill\cr\noalign{\vskip .05in}\hfill
  {}-{}r^{2(k-1)}(1-r^2)
  \sin k\phi \big(\sin (k+1)\phi-r^2\sin(k-1)\phi \big) \, .
  \qquad\cr}
$$
In other words, the sign of $C_{n-k+1,n-k+1}$ is that of
$$
  2\cos \phi\sin^2\phi + r^{2(k-1)}g_k(r^2)
  \eqno{(11.3.6)}
$$
where\ntn{$g_k(\cdot )$}
$$\displaylines{\qquad\qquad
    g_k(s)=
    -\sin (k\phi )(\sin k\phi\cos\phi +\cos k\phi\sin\phi )
  \qquad\hfill\cr\noalign{\vskip .08in}\hfill\qquad
    \eqalign{
      {}+{}&s(-2\cos\phi\sin^2\phi + 2\sin^2k\phi\cos\phi )\cr
      {}-{}&s^2\sin k\phi (\sin k\phi\cos\phi -\cos k\phi\sin\phi )
          \, . \cr
     }\qquad\qquad\cr
  }
$$
We can now explain how to conclude the proof. As $k$ tends to
infinity, the leading term in (11.3.6) is $r^{2(k-1)}g(r^2)$,
since we assume $r>1$. Thus, for large $k$ the diagonal
coefficient $C_{n-k+1,n-k+1}$ is positive whenever $g_k(r^2)$ is
such. Thus, our goal is to determine for which values of
$(s,\phi )$ we can have $g_k(s)$ positive for infinitely many
$k$'s.

The trick is to understand that when $\phi$ is an irrational
multiple of $2\pi$, the sequence $\sin k\phi$ fills $[\, -1,1\, ]$.
Thus, we can consider $\sin k\phi$ almost as a free parameter,
on which we can optimize. This leads us to define the function\ntn{$h(\theta )$}
$$\displaylines{\qquad\qquad
    h(\theta )=
    -\sin \theta\, (\sin\theta\cos\phi +\cos\theta\sin\phi )
  \qquad\hfill\cr\noalign{\vskip .08in}\hfill\qquad
    \eqalign{
      {}+{}&s(-2\cos\phi\sin^2\phi + 2\sin^2\theta\cos\phi )\cr
      {}-{}&s^2\sin\theta (\sin\theta\cos\phi -\cos\theta\sin\phi )
      \, .\cr
    }\qquad\qquad\cr
  }
$$
Formally, this function is obtained by substituting $k\phi$ for $\theta$
in the expression of $g_k$. It is convenient to define
$$
  A=(s-1)^2\cos\phi \, ,\quad
  B=(s^2-1)\sin\phi \, ,\quad
  \hbox{and}\quad
  C=2s\sin^2\phi\cos\phi \, .
$$
Since we implicitly made the change of variable $s=r^2$ after
(11.3.6), the numbers $A$, $B$ and $C$ are all positive --- recall
$0< \phi< \pi/2$ and $r>1$.

The following result will be instrumental.

\bigskip

\state{11.3.6. LEMMA.} {\it%
  The function $h(\theta )$ is maximal at a point $\theta^*$,
  unique modulo $\pi$, and defined by
  $$
    \cos (2\theta^*)=A/\sqrt{A^2+B^2} \, ,\qquad
    \sin (2\theta^*)=B/\sqrt{A^2+B^2} \, .
  $$
  Its maximum value is
  $$
    h(\theta^*)=-{A\over 2}+{\sqrt{A^2+B^2}\over 2} -C \, .
  $$
  }

\bigskip

\stateit{Proof.} We rewrite the function $h(\cdot )$ as 
$$\eqalign{
  h(\theta ) 
  & = -A\sin^2\theta +B\sin\theta\cos\theta -C \cr
  & = -A{1-\cos 2\theta\over 2} + B{\sin 2\theta\over 2} -C 
    \, . \cr
  }
$$
Differentiating with respect to $\theta$, we see that when $h$
is maximum,
$$
  0 = h(\theta^* )
  = -A\sin 2\theta^*+B\cos 2\theta^* \, .
$$
Since neither $A$ nor $B$ vanish, this gives us $\tan\theta^*
=B/A$. Consequently, there exists $\epsilon_1,\epsilon_2$ equal
to either $-1$ or $+1$, such that 
$$
  \cos 2\theta^* = \epsilon_1 A/\sqrt{A^2+B^2}
  \qquad\hbox{ and }\qquad
  \sin 2\theta^* = \epsilon_2 B/\sqrt{A^2+B^2} \, .
$$
At such point, the value of $h(\cdot )$ is
$$
  h(\theta^*)=-{A\over 2}
  +{\epsilon_1A^2+\epsilon_2B^2\over 2\sqrt{A^2+B^2}} -C \, .
$$
It is maximum when $\epsilon_1=\epsilon_2=1$. This determines
$2\theta$ modulo $2\pi$, and therefore, $\theta$ modulo
$\pi$.\hfill$\qed$

\bigskip

We can now determine when the maximum of $h(\cdot )$ is positive.

\bigskip

\state{11.3.7. LEMMA.} {\it%
  In the domain $a<0$ and $b<-1$ with $b<-a^2/4$, the function
  $h(\cdot )$ has a positive supremum if and only if $b<a-1$.
  }

\bigskip

\stateit{Proof.} Using Lemma 11.3.6, the positivity of the
supremum of $h(\cdot )$ is equivalent to
$$
  \sqrt{A^2+B^2} > A+2C \, .
$$
Since $A$, $B$ and $C$ are positive, this is equivalent to
$B^2>4C^2+4AC$. Plugging the expression for $A$, $B$ and $C$ into
this last inequality, we obtain
$$
  (s^2-1)^2 > 16 s^2(1-\cos^2\phi )\cos^2\phi 
  + 8s(s-1)^2\cos^2\phi \, .
$$
Setting $c=\cos^2\phi$, we obtain a quadratic inequality
$$
  16s^2c^2-8s(s^2+1)c+(s^2-1)^2 > 0\, .
  \eqno{(11.3.7)}
$$
The quadratic function of $c$ involved has two positive roots,
$$
  c_-={(s-1)^2\over 4s} \qquad \hbox{ and }\qquad
  c_+={(s+1)^2\over 4s} \, .
$$
Since $c_+>1$ (recall $s=r^2>0)$, inequality (11.3.7) is
equivalent to
$$
  \cos^2\phi < c_- \, .
$$
Going back to the parameterization $a=-2r\cos\phi$ and $b=-r^2$,
that is $a^2=4s\cos^2\phi$ and $b=-s$, we rewrite the above
inequality. After a simplification by $4s$, it gives
$$
  a^2 < (-b-1)^2 \, .
$$
Since $b<-1$ and $a<0$, it is equivalent to $b<a-1$.\hfill$\qed$

\bigskip

We can now state our final lemma.

\bigskip

\state{11.3.8. LEMMA.} {\it%
  Assume that $a<0$ and $b<\min(-a^2/4,-1)$. If $b<a-1$ and
  $\phi$ is an irrational multiple of $2\pi$, then 
  $\limsup_{k\to\infty}g_k(r^2)>0$. On the other hand, if
  $b>a-1$, then there exists a positive $\epsilon$ such that 
  $g_k (r^2)\leq -\epsilon$ for all $k\geq 1$.
  }

\bigskip

\stateit{Proof.} Assume $b<a-1$. Combining Lemmas 11.3.7 and
11.3.8, let $\epsilon$ be a positive number such that $h(\theta )$ is
positive on an $\epsilon$-neighborhood of $\theta^*$. If $\phi$ is an
irrational multiple of $\pi$, the sequence $k\phi$ intersect
$[\,\theta^*-\epsilon ,\theta^*+\epsilon\, ]+2\pi\ZZ$ infinitely
often; this follows from Kronecker's approximation theorem in number
theory --- see, e.g., Hlawka, Schui\ss engeier and Taschner
(1986). Consequently,
$$
  \limsup_{k\to\infty} g_k(r^2)= h(\theta^* ) >0 \, .
$$

If $b>a-1$, then $h(\cdot )$ is a negative function. Since
$g_k(r^2)=h(k\phi )\leq h(\theta^*)$, the result
follows.\hfill$\qed$

\bigskip

To conclude the proof of Theorem 11.3.2, we still assume $a<0$ and
$b<\min(-a^2/4,-1)$. If $b<a-1$, Lemma 11.3.8 shows that whenever
$\phi$ is an irrational multiple of $\pi$, the limit superior of
(11.3.6) is $+\infty$. For such values of $\phi$, the pair $(a,b)$ is
covered by infinitely many regions $R_k$. Thus in the range $b<a-1$,
the only regions not eventually covered by $\bigcup_{k\geq 1}R_k$ are
those for which $\phi$ is a rational multiple of $\pi$. After the
change of parameterization, the complement of this potentially
uncovered set is dense in $b<a-1$.

If $b>a-1$, then (11.3.6) is less than
$$
  2\cos\phi\sin^2\phi + r^{2(k-1)}h(\theta^*) \, ,
  \eqno{(11.3.8)}
$$
which tends to $-\infty$ as $k$ tends to infinity. Therefore,
the pair $(a,b)$ can be covered by at most a finite number of
regions $R_k$. This concludes the proof of Theorem
11.3.2.\hfill$\qed$

\bigskip

Notice that we proved much more than the statement of Theorem
11.3.2. When $a<0$ and $b>-a^2/4$, Lemma 11.3.3 shows that no
region $R_k$ covers $(a,b)$. 

The proof of Lemma 11.3.4 also contains useful information. If
$a<0$ and $b>-a^2/4$, the sequence $k\mapsto C_{n-k,n-k}$ is
decreasing. Thus, the constant $c(a,b)$ in Theorem 11.3.1 and
given by Theorem 8.3.1 is
$$
  c(a,b) = K_{s,\alpha}^2\alpha^\alpha \sum_{1\leq j\leq n}
  |C_{n,j}+C_{j,n}|^\alpha \, .
$$
This expression simplifies further if one notices that
$$
  C_{n,j}=\sum_{2\leq i\leq n} A_{i,n}A_{i-1,j}
  = A_{n,n}A_{n-1,j} = A_{n-1,j} \, .
$$
Thus, $C_{n,j}= A_{n-j-2,1}$ for $1\leq j\leq n-1$. Moreover
$$
  C_{j,n}=0 \qquad \hbox{ for all } 1\leq j\leq n \, .
$$
Consequently,
$$
  c(a,b)
  = K_{s,\alpha}^2\alpha^\alpha
    \sum_{1\leq j\leq n-1} |A_{n-1,j}|^\alpha
  = K_{s,\alpha}^2\alpha^\alpha
    \sum_{1\leq j\leq n-1} |A_{n-j,1}|^\alpha \, . 
$$

In the range $a<0$ and $a-1<b<-a^2/4$, the proof of Lemma 11.3.8 shows
that at most a finite number of regions $R_k$ cover $(a,b)$. Notice
that the bound (11.3.8) shows that the number of such covering regions
is at most the largest $k$ for which (11.3.8) is positive. Ultimately,
this gives an inequality in $k$, $a$ and $b$. Potentially, this could be
used if someone were interested in proving some result for particular
values of $a$ and $b$. However, the pictures above suggest that no
region $R_k$ covers such pair a $(a,b)$; but I don't know how to prove
it.

The proof of Lemma 11.3.9 involves a number theoretic argument
which does not say what happens when $\phi$ is a rational
multiple of $\pi$. The pictures of the regions $R_i$ below 
leave the possibility that some exceptional parabola
$a=-r\cos\phi$, $b=-r^2$ with $\phi\in 2\pi\QQ$ are left
uncovered. Equation (11.3.6) shows that no parabola is left
completely uncovered. Indeed, as $r$ tends to infinity, the
sign of $C_{n-k+1,n-k+1}$ is that of $\sin k\phi\sin(k-1)\phi$.
We claim that the sequence $\sin k\phi\sin(k-1)\phi$ contains
infinitely many positive values whenever $\phi$ is in
$(0,\pi/2)$. Indeed, if $\phi$ is in $(0,\pi/2)$, then $k\phi$
and $(k-1)\phi$ are less than $\pi/2$ apart. When $\phi$ is a
rational multiple of $\pi$ the sequence $(k\phi)_{k\geq 1}$ is
periodic modulo $2\pi$. Consequently, for some $k$, both
$k\phi$ and $(k-1)\phi$ are in $(0,\pi )$ modulo $2\pi$. For
this specific $k$ we have $\sin k\phi\sin (k-1)\phi>0$.

I conjecture that Theorem 11.3.2 is sharp, meaning that the
region described in $(a,b)$ coincides with $\bigcup_{k\geq 1}R_k$.

\bigskip

Combined with our description of $R_1=\{\, (a,b)\, : \,
a>0\,\}$, Lemmas 11.3.4 and 11.3.5 allow us to describe
completely what happens in the stability region.

\bigskip

\state{11.3.10. THEOREM.} {\it%
  For the second order autoregressive process, 
  $X_n=aX_{n-1}+bX_{n-2}+\epsilon_n$,
  assume that the roots of the characteristic equation
  $x^2-ax-b=0$ are inside the unit disk. Then the tail
  of $n\gamma_n(1)$ has the form
  $$
    P\{\, n\gamma_n(1)\geq t\,\} \sim
    \cases{ c(a,b)t^{-\alpha/2}   & if $a>0$ \cr
            \noalign{\vskip .05in}
	    c(a,b)t^{-\alpha} \log t & if $a<0$ \cr} 
    \qquad\hbox{ as } \ttoi \, .
  $$
  }

In particular, in the stability region, the tail of $n\gamma_n(1)$
behaves like $t^{-\alpha/2}$ if $a>0$ and like $t^{-\alpha}\log t$ if
$a<0$. In some sense, this generalizes Theorem 11.2.1 to
autoregressive processes of order $2$. If $b=0$ we recover Theorem
11.2.1.

%%%%%%%%%%%%%%%% DO NOT PRINT %%%%%%%%%%%%%%%%%%%%%%%%%%%%%
\do_not_print{

\state{11.3.4. PROPOSITION.} {\it%
  If $a<0$ and $b\geq -1$ then the largest diagonal coefficient
  of $C$ vanishes.
  }

\bigskip

\stateit{Proof.} If $a^2+4b>0$, the result follows from
Proposition 11.3.2. If $a^2+4b<0$, Lemma 11.3.3 shows that for
$b\geq 1$, (i.e. $r\leq 1$, the sign of $C_{n-k+1,n-k+1}$ is
that of
$$\displaylines{
    f_k= -2\cos\phi\sin\phi^2 (1-r^{2k})
  \cut
    + r^{2(k-1)}(1-r^2)\sin k\phi 
    \big( \sin (k+1)\phi -r^2\sin (k-1)\phi \big) \, . \cr
  }
$$
To prove Proposition 11.3.4, it suffices to prove that $f_k\leq
0$ when $0\leq \phi\leq \pi/2$ and $r\leq 1$. This will be
achieved in using different bounds for $f_k$ according to the
values of $r$ and $\phi$. We start with the trivial case.

If $\sin k\phi \sin (k+1)\phi\leq 0$ and $\sin k\phi \sin
(k-1)\phi\geq 0$, then 
$$
  f_k\leq -2\cos\phi\sin^2\phi (1-r^{2k})\leq 0 \, .
$$
Assume now that $\sin k\phi \sin (k-1)\phi\leq 0$. Since $0\leq
r\leq 1$, we have
$$\eqalign{
  f_k
  &= -2\cos\phi\sin^2\phi (1+r^2+\cdots + r^{2(k-1)})(1-r^2) \cr
  &\qquad + 2r^{2(k-1)}\sin k\phi 
   \big( \sin (k+1)\phi -r^2\sin (k-1)\phi\big) \cr
  &\leq -2\cos\phi\sin^2\phi k r^{2(k-1)}(1-r^2) \cr
  &\qquad + r^{2(k-1)}(1-r^2)\sin k\phi
   \big(\sin (k+1)\phi-\sin (k-1)\phi\big) \cr
  }
$$
Consequently, the inequality $f_k\leq 0$ is implied by
$$
  -2k\cos\phi\sin^2\phi + \sin k\phi
  \big(\sin (k+1)\phi - \sin (k-1)\phi \big) \leq 0 \, .
$$
Expressing $\sin (k+1)\phi$ and $\sin (k-1)\phi$ as functions of
sines and cosines of $\phi$ and $k\phi$, it suffices to have
$$\displaylines{
  -2k\cos\phi\sin^2\phi + 2\sin k\phi\cos k\phi\sin \phi
  \cut
  = \sin\phi\big( \sin (2k\phi )-k\sin 2\phi\big) \leq 0 \, .
  \cr}
$$
Since $\sin 2\phi\geq 0$, one can see by induction on $k$ that
$$
  |\sin (2k\phi)|\leq k\sin 2\phi \, ,
$$
which in turns implies $f_k\leq 0$.

It remains to show that $f_k$ is nonpositive when $\sin
(k-1)\phi\sin k\phi>0$ and $\sin k\phi\sin (k+1)\phi>0$. In this
case,
$$\eqalign{
  f_k
  &\leq -2(1-r^{2k})\cos\phi\sin^2\phi
   + r^{2(k-1)}(1-r^2)\sin k\phi 
   \big(\sin (k+1)\phi + \sin (k-1)\phi\big) \cr
  &= -2(1-r^{2k})\cos\phi\sin^2\phi 
   + r^{2(k-1)}(1-r^2)\sin k\phi 2\sin (k\phi)\cos\phi \cr
  &= -2(1-r^{2(k-1)}) + 2r^{2(k-1)}(1-r^2)\cos\phi 
   (-\sin^2\phi+\sin^2k\phi ) \, . \cr
  }
$$
Consequently, if moreover $\sin^2 k\phi\leq \sin^2\phi$, then
indeed $f_k\leq 0$. 

In the last case, we have  $\sin(k-1)\phi\sin k\phi>0$ and 
$\sin k\phi\sin (k+1)\phi>0$ and $\sin^2k\phi\geq \sin^2\phi$.
We write
$$
  f_k= -2\cos\phi\sin^2\phi + g(r^2) \, ,
$$
where $g$ is the polynomial
$$\eqalign{
   &t^{k-1}\sin k\phi\sin (k+1)\phi \cr
  +&t^k \big(2\cos\phi\sin^\phi -\sin k\phi\sin (k+1)\phi
     -\sin k\phi\sin (k+1)\phi\big) \cr
  +&t^{k+1}\sin k\phi\sin (k-1)\phi \, .\cr
  }
$$
To conclude the proof, we are going to show that $g$ is
nondecreasing. This implies
$$
  f_k\leq -2\cos\phi\sin^2\phi + g(1) = 0 \, .
$$
We calculate
$$\displaylines{
  g'(t)=t^{k-2}\big( (k-1)\sin k\phi\sin (k+1)\phi
  \cut
  \eqalign{+& 2kt\cos\phi (\sin^2\phi-\sin^2k\phi ) \cr
  	   +& (k+1)t^2\sin k\phi\sin (k-1)\phi )\cr}\cr
  }
$$
Thus $g'$ is the product of a monomial and a quadratic term. For
$g'$ to be nonnegative, on $[0,1]$, it suffices to have the
discriminant of the quadratic term negative (recall that we
assume $\sin k\phi\sin (k-1)\phi>0$), that is
$$\eqalign{
  k^2\cos^2\phi (\sin^2\phi-\sin^2k\phi )^2
  &\leq (k+1)(k-1)\sin^2k\phi\sin (k-1)\phi\sin (k+1)\phi \cr
  &=(k^2-1)\sin^2 k\phi (\sin^2k\phi\cos^2\phi \cr
  &\qquad -\cos^2k\phi\sin^2\phi )\, . \cr
  }
$$
We rewrite this inequality as
$$\displaylines{
    (\sin^2k\phi-\sin^2\phi)\cos^2\phi 
    (k^2-\sin^2 k\phi +k^2\sin^2k\phi )
  \cut
    +(k^2-1)\sin^2 k\phi\sin^2\phi (\cos^2\phi-\cos^2k\phi)
    \geq 0 \qquad (11.3.6)\cr
  }
$$
Since we are on the range where $\sin^2k\phi\geq \sin^2\phi$, we
have $\cos^2k\phi\leq \cos^2\phi$. Thus, for (11.3.6) to hold,
it suffices to have
$$
  k^2+(k^2-1)\sin^2\phi \geq 0 \, ,
$$
which is plain. Thus, on the given range of $\phi$ and
$k$, we proved that $g'\geq 0$. This concludes the proof of
Proposition 11.3.4.\hfill$\qed$

} % end do_not_print

%%%%%%%%%%%%%%%% END DO NOT PRINT %%%%%%%%%%%%%%%%%%%%%

\bigskip
\centerline{\fsection Notes}
\section{11}{Notes}

\bigskip

\noindent A basic reference on the classical theory and 
applications of time series is Brockwell and Davis (1987). 

For linear processes with heavy tailed errors, Davis and Resnick
(1986) developed the asymptotic theory as the number of observations
goes to infinity. A slightly different perspective, first
letting $n$ tend to infinity, and then looking at the tail
of the limiting distribution has been investigated in a series
of papers by Mijnheer (1997a, b, c).

When the errors $\epsilon_i$'s have a nearly symmetric
distribution and enough moments, it is tempting to work as if
they were from a Weibull-like distribution, or even from a normal
one. In the later case, the autocovariances are weighted sums of
chi-square random variables. I don't know how to assess the
accuracy of such an approximation.

I somewhat believe that a proper understanding of the regions
$R_i$ in general requires adding some algebraic geometric
tools. Cox, Little and O'Shea (1992, 1998) may be a good
starting point for statisticians interested in pursuing this
research path. I don't know what is the analogue of Theorem 11.3.2
for autoregressive models of order larger than or equal to $3$.

One may think that our choice of heavy tailed errors is the origin
of the complicated tail behavior of the autocovariances. If the
$\epsilon_i$'s have a spherical distribution, we can use the
classical Gaussian trick of diagonalizing the matrix of the
quadratic form. When looking at the tail behavior, one is 
then led to the nontrivial question
of relating the dimension of the largest eigensubspace of the
matrix $C$ to the parameter $\theta$ of the autoregressive model. 
A plot of the spectral gap for the matrix $C$ as a function of
$a$ and $b$ --- thus, for autoregressive of order $2$ --- suggests
an incredibly complex behavior of the tail of the the aucovariance,
even with Gaussian errors.
\vfill\eject

\chapter{12}{Suprema of some stochastic processes}
\line{\hfill{\fchapter 12.\quad Suprema of some}}
\vskip .1in
\line{\hfill{\fchapter stochastic processes}}

\vskip 54pt

\noindent In chapter 9, we studied how the tail of the distribution of the
supremum of a random linear form is related to integration
over some asymptotic sets. The goal in this chapter is to go a
little further, considering some examples which are of
pedagogical interest.

\bigskip

\noindent{\fsection 12.1. Maxima of processes and maxima of
their variances.}
%\chapternumber={12.1}
\section{12.1}{Maximum of processes and maximum of their variance}

\bigskip

\noindent Consider a centered Gaussian process $X(m)$ indexed by
some abstract set\ntn{$M$} $M$. As in chapter 9, write\ntn{$X(M)$}
$$
  X(M)=\sup\{\, X(m) \, :\, m\in M\,\}
$$
for its supremum. Define\ntn{$\sigma^2(M)$}
$$
  \sigma^2(M)=\sup\{\, \Var X(m) \, :\, m\in M\,\}
$$
to be the supremum of its variance. A famous result of Fernique (1970),
Landau and Shepp (1970) asserts that whenever $\sigma^2(M)$ is
finite,
$$
  \lim_{\ttoi} t^{-2}\log P\{\, X(M)\geq t\,\} 
  = -{1\over 2\sigma^2(M)} \, .
  \eqno{(12.1.1)}
$$

This result is often interpreted in saying that, in the logarithmic
scale, the tail of the supremum of the process is driven by the
points of largest variance. A heuristic argument of why this
should be the case is that when the variance is large, the
process tends to fluctuate more. And so, we should expect its
supremum, when large, to be near such a point. As it is, this
heuristic argument could be applied to any process. The aim of
this section is to show that this heuristic is wrong.

We are going to construct some processes whose supremum tail is
driven by the points of smallest variance. This is an
application of ideas developed in chapter 9. The key is to
understand why the maximum variance appears in the Gaussian
case. Proposition 9.2.1 tells us that this happens because
$I_\sbullet$ is proportional to $1/|x|^2$ in the Gaussian case,
as $x$ tends to infinity. That is $I_\sbullet$ 
is inversely proportional to the Euclidean norm 
--- a very specific feature. And the Euclidean norm 
is a monotone function of the variance.
Thus, to build a counterexample to the heuristic, we need to
have a set $M$ and a function $I$ such that $I_\sbullet$ is
minimal on $M$ at points of minimal Euclidean norm.

Let $\alpha$ be in $(1,2)$ and $\theta=\big(\Gamma (1/\alpha )/\Gamma
(3/\alpha )\big)^{\alpha/2}$. The function
$$
  f_\alpha (x) = 
  {\alpha\over 2\sqrt{\Gamma (1/\alpha )\Gamma (3/\alpha)}}
  e^{-\theta |x|^\alpha} \, ,\qquad x\in\RR\, ,
$$
defines a density. It has zero expectation, and unit variance. For
any positive $\beta$, define the unit sphere for the
$\ell_\beta$-norm\ntn{$S_d^{(\beta )}$}
$$
  S_d^{(\beta )} =\{\, x\in\RR^d \, :\, |x|_\beta =1\,\} \, .
$$
If $\beta>1$, define also\ntn{$M_\beta$}
$$
  M_\beta =\{\, \sum_{1\leq i\leq d} \sign (m_i)|m_i|^{\beta-1}e_i \, :\,
  m\in S_d^{(\beta )} \,\} \, .
$$

\bigskip

\state{12.1.1 THEOREM.} {\it%
  Let $X=(X_1,\ldots , X_d)$ be a random vector in $\RR^d$ with
  independent components, all distributed with density
  $f_\alpha$. If $1<\alpha <\beta <2$, then
  $$
    \lim_{\ttoi} t^{-\alpha}\log P\{\, X(M_\beta )\geq t\,\}
    = -{1\over \inf_{m\in M_\beta}\Var X(m)} \, .
  $$
  }

\bigskip

\stateit{Proof.} For $x$ in $\RR^d$, let $I(x)=\sum_{1\leq i\leq
d}\log f(x_i)$. For some constant $c$,
$$
  I(x)=\theta\sum_{1\leq i\leq d}|x_i|^\alpha + c \, .
$$
This is a strictly convex function since $\alpha >1$. It
satisfies (9.2.1). Moreover, Lemma 9.1.9 gives
$$
  I_\sbullet (x)=c+ |x|^{-\alpha}_{\alpha /(\alpha -1)} \, .
$$
Consequently, as $t$ tends to infinity, 
$$\eqalign{
  I_\sbullet (M_\beta /t)
  & \sim t^{-\alpha} \big/
    \sup\Big\{\, \sum_{1\leq i\leq d}|x_i|^{\alpha /(\alpha -1)}
    \, :\, x\in M_\beta \,\Big\}^{\alpha-1} \cr
  & = t^{-\alpha}\big/
    \sup\Big\{\, \sum_{1\leq i\leq d}
    |m_i|^{(\beta-1)\alpha\over\alpha-1} \, :\, 
    m\in S_d^{(\beta)} \,\Big\}^{\alpha-1} \, . \cr
  }
$$
If $m$ belongs to $S_d^{(\beta )}$, then each
component $|m_i|$ is less than or equal to $1$. Then, the
inequality $(\beta -1)/\beta > (\alpha -1)/\alpha$ gives
$$
  \Big( \sum_{1\leq i\leq d}
        |m_i|^{(\beta-1)\alpha\over\alpha-1} 
  \Big)^{\alpha-1}
  \leq
  \Big( \sum_{1\leq i\leq d}
        |m_i|^\beta\Big)^{(1/\beta)\alpha\beta}
  \leq 1\, ,
$$
with equality if and only if $m$ has exactly one coordinate
equal to $1$ or $-1$. Thus Proposition 9.2.1 yields
$$
  \lim_{\ttoi} t^{-\alpha}\log P\{\, X(M)\geq t\,\}
  =-1 \, .
$$
Next the variance of $X(m)$ is $|m|^2$. But if 
$m$ belongs to $M_\beta$, the inequality $|m|^2\geq 1$ holds,
with equality if and only if one of the components of $m$ is $1$
or $-1$, and all the others are zero. This concludes the
proof.\hfill$\qed$

\bigskip

The proof we gave is reasonably short, but a bit mysterious. For
$d=2$, the following picture makes the result obvious if one
keeps in mind Laplace's method. It represents the upper right
quadrant of the plane, and the various sets involved, for
$\alpha=1.2$ and $\beta=1.5$. The set on which the integration
is performed is shaded gray. Its boundary is the polar
reciprocal of $M_\beta$. The set $M_\beta$ is the black line
inside the shaded area. The Euclidean unit sphere is the
white line inside the shaded area. The last black line is the
level set of the density $f_\alpha(x)f_\alpha(y)$.

\medskip

\setbox1=\hbox to 1.2in{\vbox to 1.6in{\ninepointsss
  \hsize=1.05in
  \kern 1.4in
  \includegraphics{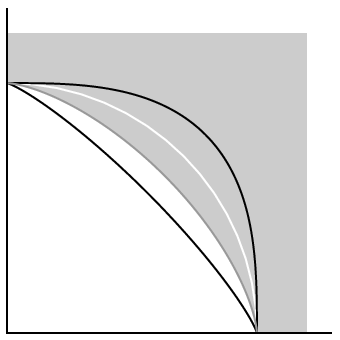}
  \anote{-.1}{-.1}{0}
  \anote{-.1}{1}{1}
  \anote{1}{-.1}{1}
  \vfill
%  \centerline{n=10, a=1, $\scriptstyle\alpha$=1}
}\hfill}
\centerline{\box1}

Notice that for $d=2$ we have a process indexed by the one
dimensional set $M_\beta$. For $d>2$ we have a random field.
More importantly, Theorem 12.1.1 may be a prototype for a
misleading statement! The proof shows that the left hand side
of the statement has no intrinsic connection with the right
hand side! The proof shows that indeed the variance appears by
a pure coincidence. In my opinion, the Gaussian case is not any
different, as Proposition 9.2.1 shows.

\bigskip

\noindent{\fsection 12.2. Asymptotic expansions for the tail of
the supremum of Gaussian processes can be arbitrarily bad.}
%\chapternumber={12.2}
\section{12.2}{ Asymptotic expansions ... can be arbitrarily bad}

\bigskip

\noindent For many reasons, both theoretical and applied, there has
been a large literature devoted to the approximation of the tail
probability of the supremum of Gaussian processes. In view of
(12.1.1), it is quite natural to search for an approximation of the form
$$
  P\{\, X(M)\geq t\,\}
  = t^{\alpha }\exp \Big( {-t^2\over 2\sigma^2(M)}\Big)
  \big( P_n (1/t) + o(t^{-n})\big) \qquad
  \hbox{ as } \ttoi \, 
  \eqno{(12.2.1)}
$$
where $P_n$ is a polynomial of degree $n$. The hope of course is
that for moderate $t$'s, this expansion provides an accurate
approximation.

The aim of this section is to provide an example where such expansion
holds, but, no matter what, provides a poor approximation for fixed
$t$. It is essential to remember that an asymptotic expansion like
(11.2.1), if it exists, is unique --- see, e.g., Olver, 1974, \S
1.7. In particular, it does not depend at all on which method is used
to derive it. Thus, the failure we want to describe is not that of a
particular method. The one proposed in these notes as well as any
other fails, and there is no way arround if one sticks to
approximations of the form (12.2.1).  The basic idea in this section
is to mimic what happened for autoregressive models. Our approximation
was not so good when we appealed to Theorem 8.3.1, because it was
actually quite likely that the largest random variable was not an
$\epsilon_i$ for which $C_{i,i}=0$; even if one conditions by the
appearance of a large deviation, it is quite likely that it is caused
by a large $\epsilon_i$ for which $C_{i,i}=0$.

To build our example we will first give it in a geometric form.
We will discuss afterwards some of its features.

Consider the convext set in $\RR^d$,\ntn{$C$}
$$
  C=\{\, p \, :\, \< p,e_1\>\leq 1\,\}
  \bigcap \{\, p\, : \, |p|\leq 1+\epsilon\,\} \, .
$$
Denote by $M$ the polar reciprocal of $\partial C$. Let $X$ be a random
vector in $\RR^d$, having a centered normal distribution, with
independent components. Then
$$
  X(M)\geq t \hbox{ if and only if } X\not\in tC \, .
$$
Define the polynomials
$$
  P_n(u) = \sum_{0\leq k\leq n} (-1)^k {(2k)!\over 2^{2k}k!} u^k
  \, .
$$

\state{12.2.1 PROPOSITION.} {\it%
  The tail expansion
  $$
    P\{\, X(M)\geq t\,\} = {e^{-t^2/2}\over \sqrt{2\pi}t}
    \big( P_n(1/t)+o(t^{-n})\big) \qquad\hbox{ as } \ttoi
  $$
  holds. However, for any $t\geq 0$, 
  $$
    {\sqrt{2\pi} te^{t^2/2}\over P_n(1/t)}
    P\{\, X(M)\geq t\,\} \geq 2\sqrt{2\pi} 
    e^{-t^2(\epsilon + (\epsilon^2/2))} t^{d-(1/2)} \, .
  $$
  }

\bigskip

Before proving this result, let us see why this provides the
proper example. The first statment lets us hope that
$e^{-t^2/2}P_n(1/t)/\sqrt{2\pi}t$ is a good approximation of the
tail probability of the supremum. The second statment asserts
that it is not the case if $\epsilon$ is small and $t$ is
large, but not too large. The following
plot shows the lower bound for $d=2$, $n=2$ and various values
of $\epsilon$, as a function of the approximation. 
For instance, if $\epsilon =1$, the lower bound
is less than $1$, which means that the approximation may
underestimate. A more interesting value is for $\epsilon=0.5$; when
the approximation is about $10^{-1.5}\approx 3\%$, the lower
bound is about $2$. Thus, the approximation underestimates the
correct probability by a factor at least $2$.
One should keep in mind that for a typical
statistical application, we are interested in $t$'s such that
$P\{\, X(M)\geq t\,\}$ is between $10^{-1}$ and $10^{-2}$.

\setbox1=\hbox to 3.95in{\vbox to 2.7in{\ninepointsss
  \hsize=3.95in
  \kern 2.3in
  \includegraphics{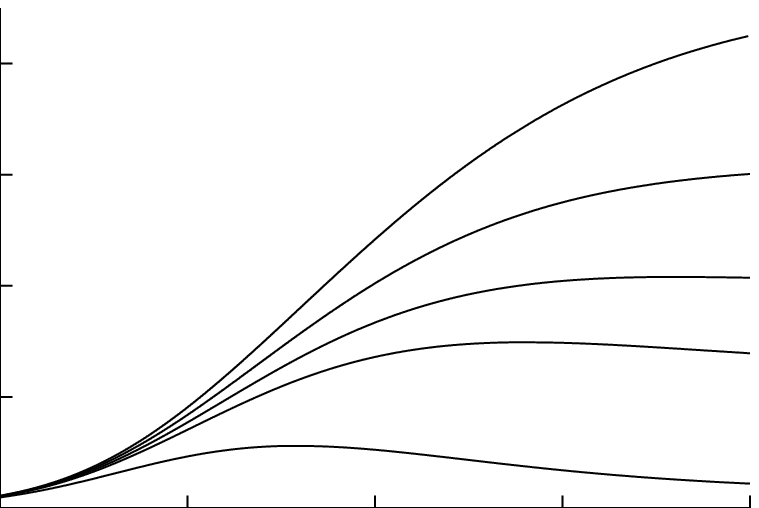}
  \anote{-.1}{-.1}{0}
  \anote{0.65}{-.2}{$10^{-0.5}$}
  \anote{1.4}{-.2}{$10^{-1}$}
  \anote{2.15}{-.2}{$10^{-1.5}$}
  \anote{2.90}{-.2}{$10^{-2}$}
  \anote{-.15}{0.44}{1}
  \anote{-.15}{0.89}{2}
  \anote{-.15}{1.34}{3}
  \anote{-.15}{1.77}{4}
  \anote{3.1}{0.1}{$\epsilon=1$}
  \anote{3.1}{0.6}{$\epsilon=0.6$}
  \anote{3.1}{0.9}{$\epsilon=0.5$}
  \anote{3.1}{1.34}{$\epsilon=0.4$}
  \anote{3.1}{1.86}{$\epsilon=0.3$}
  \anote{3.19}{-0.03}{value of the}
  \anote{3.19}{-0.15}{approximation}
  \vfill
}\hfill}

\centerline{\box1}

We will comment further on this example after we prove our
statement.

\bigskip

\stateit{Proof of Proposition 12.2.1.} We first start with the
obvious bound
$$\eqalignno{
  P\{\, X(M)\geq t\,\}
  & = P\{\, X_1\geq t\hbox{ or } |X|\geq t(1+\epsilon )\,\} \cr
  \noalign{\vskip 0.05in}
  & \cases{\geq P\{\, X_1\geq t\,\} &\cr
           \noalign{\vskip 0.05in}
           \leq P\{\, X_1\geq t\,\} 
           + P\{\, |X|\geq t(1+\epsilon )\,\} \, . &\cr}
  &(12.2.1)\cr
  }
$$
The crude logarithmic estimate of Proposition 2.1 shows that
$$
  P\{\, |X|\geq t (1+\epsilon )\,\} 
  = O(e^{-t^2(1+\epsilon )^2/2} t^{2d} ) 
  \qquad\hbox{ as } \ttoi \, .
$$
On the other hand, the standard asymptotic expansion for the
complementary error function (see, e.g., Olver, 1974, \S 3.1.1)
yields
$$
  P\{\, X_1\geq t\,\} = {e^{-t^2/2}\over \sqrt{2\pi}t}
  \big( P_n(1/t)+o(t^{-n})\big) \qquad \hbox{ as } \ttoi \, .
$$
Thus (12.2.1) provides the asymptotic expansion in the first
assertion of Proposition 12.2.1.

We now use another lower bound, namely
$$
  P\{\, X(M)\geq t\,\} \geq P\{\, |X|\geq t(1+\epsilon )\,\} \, ,
$$
which follows from the equality in (12.2.1). Since $|X|^2$ has a
chi-square distribution with $d$ degrees of freedom, an
integration by parts yields
$$\eqalign{
  P\{\, |X|\geq t(1+\epsilon )\,\}
  & = \int_{t^2(1+\epsilon )^2}^\infty 
    { x^{{d\over 2}-1} e^{-x/2}\over 2^{d/2}\Gamma (d/2) }
    dx \cr
  & \geq 2 t^{d-{1\over 2}}(1+\epsilon )^{d-{1\over 2}}
    e^{-t^2(1+\epsilon )^2/2} \, . \cr
  }
$$
For any $u\geq 0$, the bound $P_n(u)\leq 1$ holds; this comes
from the fact that the asymptotic expansion for the error
function is obtained by integrating by parts, and the
integrations lead to an alternating series --- see, e.g., Olver,
1974, \S 3.1. Thus the second statement of Proposition 12.2.1
follows.\hfill$\qed$

\bigskip

When $d=2$, we can make a very explicit construction of the
process. The polar reciprocal $M$ is just a piece of circle and
a point,
$$
  M=\big\{\, (x,y)\in \RR^2 \, :\, x^2+y^2=(1+\epsilon )^{-2} 
  \, ; x\leq (1+\epsilon )^{-2} \,\big\} \cup \, \{\, (1,0)\,\} \, .
$$
This can be seen by a formal proof, but it is obvious from the
following picture. The domain on which we integrate is shaded.
The unit sphere, or equivalently the level set of the Gaussian
measure, is the dark sphere. The set $M$ is an open arc, in 
black as well, with the point $e_1$ marked.

\setbox1=\hbox to 1.4in{\vbox to 1.7in{\ninepointsss
  \hsize=1.05in
  \kern 1.6in
  \includegraphics{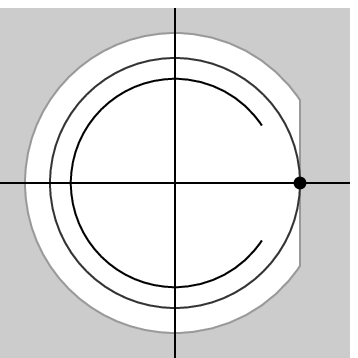}
  \anote{0.6}{0.6}{0}
  \anote{1.2}{0.6}{1}
  \vfill
%  \centerline{n=10, a=1, $\scriptstyle\alpha$=1}
}\hfill}
\centerline{\box1}

We choose $M$ to be reduced --- see definition in section 9.1.
We could as well index the process by a larger set, such as
$$
  (1+\epsilon )^{-1}S_1\cup \,\{\, (1,0)\,\}
$$
or even
$$
  (1+\epsilon )^{-1} S_1\cup \,
  [\, (1+\epsilon )^{-1}e_1 ,e_1\,]
  \, .
$$
This last set can be parameterized as follows. Let
$$
  f(t)=\cases{ {\displaystyle 1\over
                \vbox to 8pt{}\displaystyle 1+\epsilon}
               \big( \cos (2\pi t),
               \sin (2\pi t)\big) & if $0\leq t\leq 1$\cr
             \noalign{\vskip 0.05in}
               \Big( {\displaystyle \epsilon t+1-\epsilon\over
                     \vbox to 8pt{}\displaystyle 1+\epsilon} ,
                     0\Big)       & if $1\leq t\leq 2$. \cr}
$$
The corresponding Gaussian process is
$$
  X(t)=\cases{ {\displaystyle 1\over\vbox to 8pt{}
                \displaystyle 1+\epsilon }
               \big( X_1\cos (2\pi t) +
		X_2\sin (2\pi t)\big) & if $0\leq t\leq 1$ \cr
               \noalign{\vskip 0.05in}
               {\displaystyle \epsilon t+1-\epsilon \over
                \vbox to 8pt{}\displaystyle  1+\epsilon }X_1
                                      & if $1\leq t\leq 2$.\cr
             }
$$
Its variance is
$$
  \Var X(t)=
   \cases{ {\displaystyle 1\over
            \vbox to 8pt{}\displaystyle  (1+\epsilon )^2}
            & if $0\leq t\leq 1$\cr
           \noalign{\vskip 0.05in}
	   \Big({\displaystyle \epsilon t + 1-\epsilon\over
                 \vbox to 8pt{}\displaystyle 1+\epsilon}\Big)^2
            & if $1\leq t\leq 2$.\cr}
$$
Now one can argue that our example is specific; the maximal 
variance is achieved when $t=2$, that is on the boundary of 
the domain. Well, we can always define
$$
  Y(t)=
  \cases{X(t) & if $0\leq t\leq 2$ \cr
         X(4-t) & if $2\leq t\leq 4$. \cr}
$$
This process has maximum variance at $t=2$. And certainly, one
could argue that $\Cov (X(t),X(s))=1$ if $1\leq s,t\leq 2$, and
thus this process is pathological. This argument can be also
ruled out by perturbing each coordinate $Y(t)$ by some very 
tiny multiple of a brownian bridge.

The moral of the story is that the asymptotic expansion should not be
worked out blindly. One should certainly make a careful study of the
covariance of the process and be very cautious when the variance does
not vary much. Notice that the ratio of the maximal variance of the
process to the minimal one is $(1+\epsilon )^2$. For $\epsilon=0.3$,
this is $1.69$, which is not that small. Going back to the lower bound
in Proposition 12.2.1, notice that the polynomial term in the lower
bound is of order $t^{d-(1/2)}$. As the dimension $d$ increases, the
asymptotic expansion gives a worse approximation. It can be arbitrarily
bad by just taking $d$ large enough. Therefore, the constancy of the
variance should be measured with respect to the dimension. It is
therefore very unclear what happens in large dimensions or even in
infinite dimensions. It is also unclear how to assess a priori the
approximating quality of an asymptotic expansion in this context.

\bigskip
\noindent{\fsection 12.3. Maximum of nonindependent Gaussian
random variables.}
%\chapternumber={12.3}
\section{12.3}{Maximum of nonindependent Gaussian
random variables}

\bigskip

\noindent When one wants to simulate numerically a Gaussian process,
there is not much choice other than to discretize it. To what extent
can we obtain an approximation of the distribution of the original
process by that of the corresponding discretization?  There is no
claim that this section brings some new result.  The one we are going
to prove now can be derived from others existing in the
literature. But its derivation may be of pedagogical interest.

Having in mind a discretized process, let $X$ be a Gaussian vector in
$\RR^d$, with mean $0$ and definite positive covariance
matrix\ntn{$\Sigma$} $\Sigma$. Let\ntn{$\sigma^2$}
$\sigma^2=\max_{1\leq i\leq d}\Sigma_{i,i}$ be a largest diagonal
element.

\bigskip

\state{12.3.1. THEOREM.} {\it%
  Let $X$ have a Gaussian distribution, centered, with
  positive definite covariance matrix $\Sigma$. Then
  $$
    P\{\, \max_{1\leq i\leq d} X_i\geq t\,\}
    \sim {\sigma e^{-t^2/2\sigma^2}\over t\sqrt{2\pi}}
    \sharp\{\, i \, :\, \Sigma_{i,i}=\sigma^2\,\} \, .
  $$
  }

\bigskip

\stateit{Proof.} We apply Theorem 7.1. Define
$$
  A_t=tA_1 = t\{\, x\in\RR^d \, :\,
  \max_{1\leq i\leq d} x_i\geq 1\,\} \, .
$$
Set
$$
  I(x)={1\over 2}x^\T\Sigma^{-1}x \, .
$$
This is a convex function, homogenous of degree $\alpha=2$, and
$$
  P\{\, \max_{1\leq i\leq d} X_i\geq t\,\}
  = {1\over (2\pi )^{d/2}(\det \Sigma )^{1/2}}
  \int_{tA_1} e^{-I(x)} \,\d x \, .
$$
To apply Theorem 7.1, we need to minimize $I$ over $A_1$.
Since $\Sigma$ is symmetric, we can diagonalize it and write
$\Sigma= QDQ^\T$ with $D$ diagonal and $Q$ orthogonal. The change
of variable $x=Q^\T y$ shows that
$$\eqalign{
  I(A_1)
  & =\inf\Big\{\, {1\over 2}x^\T\Sigma^{-1} x\, :\, 
    \max_{1\leq i\leq d} x_i \geq 1\,\Big\}  \cr
  & = \inf\Big\{\, {1\over 2} y^\T Dy \, :\, 
    \max_{1\leq i\leq d} \< y,Q^\T e_1\>\geq 1\,\Big\} \, .\cr
  }
$$
Writing the Lagrangian to optimize $y^\T D^{-1}y$ subject to the
constraint $\< y,Q^\T e_i\>=1$ and optimizing over $i$, we obtain
$$
  I(A_1)={1\over 2} \min_{1\leq i\leq d}{1\over e_i^\T De_i}
  = {1\over 2\max_{1\leq i\leq d}\Sigma_{i,i}}
  = {1\over 2\sigma^2} \, .
$$
Moreover, $I(A_1)$ is achieved for the points
$y=DQ^\T e_i/\sigma^2$, or equivalently $x=\Sigma e_i/\sigma^2$.
Thus, the dominating manifold for $A_1$ is
$$
  \calD_{A_1} =
  \Big\{ {\Sigma e_i\over \sigma^2}\, :\, i \hbox{ such that }
  \Sigma_{i,i}=\sigma^2\,\Big\} \, .
$$
It is of dimension $k=0$. Since $\D I=\Sigma^{-1}$, Theorem
7.1 yields
$$
  P\{\,\max_{1\leq i\leq d}X_i\geq t\,\} \sim 
  { e^{-t^2/2\sigma^2}\over t\sqrt{2\pi}(\det\, \Sigma )^{1/2}}
  \sum_{i:\Sigma_{i,i}=\sigma^2} 
  { \sigma^{d+1}\over 
    \big(\det\, G_{A_1}(\Sigma e_i/\sigma^2)\big)^{1/2} }
$$
as $t$ tends to infinity. We need to calculate $G_{A_1}$.  As
mentioned after the statement of Theorem 7.1.1, it is obtained as the
compression of the difference of two second fundamental forms.  The
one for $\partial A_1$ vanishes since $\partial A_1$ is locally a flat
hyperplane. That for the level set of $I$ is $\D^2I/|\D I|$. At
$\Sigma e_1/\sigma^2$, its value is $\sigma^2\Sigma^{-1}$. The tangent
space at $\partial A_1$ at this point is
$$
  \{\,\D I(\Sigma e_i/\sigma^2 )\,\}^\perp 
  = (e_i/\sigma^2)^\perp
  = e_i^\perp \, .
$$
Thus, $G_{A_1}(\Sigma e_i/\sigma^2)$ is the expression of
$\sigma^2\Sigma^{-1}$ to $e_i^\perp$. Therefore
$$
  \det\, G_{A_1}(\Sigma e_i/\sigma^2)
  = \sigma^{2(d-1)}\det\,\<\Sigma^{-1}e_k,e_l\>_{1\leq k,l\leq d
  \atop k,l\neq i} \, .
$$
Thus, it is $\sigma^{2(d-1)}$ times the determinant of the
$(i,i)$-cofactor of $\Sigma^{-1}$, which is
$\Sigma_{i,i}\det\,\Sigma^{-1}$. Consequently,
$$
  {\sigma^{d+1}\over (\det\, \Sigma )^{1/2}\det\, G_{A_1}
  (\Sigma e_i/\sigma^2)^{1/2} } =\sigma \, ,
$$
and this gives the result putting all the estimates
together.\hfill$\qed$

\bigskip

The fact we now want to stress is about discretizing Gaussian
processes to simulate the distribution of their maximum. Since
the maximum of the discretization is less than the maximum of
the original process, this can only give a lower bound. On the
far tail, Theorem 12.3.1 asserts that this lower bound must be
of order $\sigma e^{-t^2/2\sigma^2}/t\sqrt{2\pi}$. This implies two 
things. First, one should include the points of largest variance in
the discretized sequence. This is almost common sense. Second, in
theory, the far tail will be well approximated only if that of the
original process behaves like $\sigma e^{-t^2/2\sigma^2}/t\sqrt{2\pi}$. 
These processes have been characterized by Talagrand (1988).

\bigskip

\noindent{\fsection 12.4. The 
truncated Brownian bridge.}
%\chapternumber={12.4}
\section{12.4}{The truncated Brownian bridge}

\bigskip

\noindent The Brownian bridge\ntn{$B(\cdot )$} $B$ on $[\,0,1\, ]$ 
is a centered Gaussian process with covariance
$$
  {\rm E}B(s)B(t) = st- s\wedge t \, .
$$
It is a classical result that it admits the Karhunen-Lo\`eve
expansion
$$
  B(s) = {\sqrt 2\over \pi}\sum_{k\geq 1} 
  X_k{\sin (k\pi s)\over k} \, , \qquad 0\leq s\leq 1 \, ,
$$
where the $X_k$'s are independent, normally distributed random
variables. It is known --- see, e.g., Billingsley, 1968, \S 11 --- 
that
$$
  P\{\, \sup_{0\leq s\leq 1} B(s)\geq t\,\} = e^{-2t^2} \, .
  \eqno{(12.4.1)}
$$
The aim of this section is to obtain an
approximation for the tail of the supremum of the truncated
series\ntn{$B_d(s)$}
$$
  B_d(s) = {\sqrt 2\over \pi }\sum_{1\leq k\leq d} X_k
  {\sin (k\pi s)\over k} \, .
$$
This example is quite interesting, because we will see that Theorem
5.1 --- or equivalently, Theorem 7.1 --- does not apply. However,
a slight change in the arguments will allow us to obtain the
desired asymptotic equivalence.

Let us first explain why Theorem 7.1 does not apply, and, in
particular, why assumption (7.5) --- or assumption (5.3) if one
uses Theorem 5.1 --- is not satisfied. Recall that $e_1,\ldots ,
e_d$ denotes the canonical basis in $\RR^d$. Define the curve\ntn{$p(s)$}
$$
  p(s)={\sqrt 2\over \pi}\sum_{1\leq k\leq d} 
  {\sin (k\pi s)\over k} e_k \, ,\qquad 0\leq s\leq 1
$$
in $\RR^d$. Introducing the Gaussian vector $Y=(Y_1,\ldots ,
Y_d)$, we have $B_d(s)=\< Y,p(s)\>$. Define
$$
  A_t=\big\{\, y\in\RR^d : \sup_{0\leq s\leq 1}\< y,p(s)\>
  \geq t\,\big\} = tA_1 \, .
$$
Let $I(y)= |y|^2/2$. Then
$$
  P\{\, \sup_{0\leq s\leq 1} B_d(s)\geq t\,\}
  = {1\over (2\pi )^{d/2}} \int_{A_t}e^{-I(y)} \, \d y \, .
$$
To find the dominating manifold, Proposition 9.1.7 combined with Lemma
9.1.9 suggest that we should search for the points $s$ maximizing the
variance of $B_d(s)$,
$$
  \Var B_d(s)= |p(s)|^2
  = {2\over \pi^2} \sum_{1\leq k\leq d} 
    {\sin^2(k\pi s)\over k^2} \, .
$$
Since $\Var B_d(s)= \Var B_d(1-s)$, it suffices to locate the
maximum in $[\, 0,1/2\, ]$. We differentiate the variance, obtaining
$$\eqalign{
  {\d \over \d s} \Var B_d(s)
  & = {2\over\pi^2} \sum_{1\leq k\leq d} 
    {2\sin (k\pi s)\cos (k\pi s) \over k} \cr
  & = {2\over\pi^2} \sum_{1\leq k\leq d} {\sin (2k\pi s)\over k}
    \, . \cr
  }
$$
It follows from Jackson's (1912) theorem --- see, e.g., Andrews,
Askey and Roy, 1999, chapter 7 --- that $(\d / \d s )\Var B_d(s)$
is positive on $(0,1/2)$. Thus $B_d(s)$ has a unique point of
maximal variance for $s=1/2$, no matter what $d$ is. The
maximal value depends on $d$. It is\ntn{$\sigma_d^2$}
$$
  \sigma_d^2=\Var B_d(1/2) 
  = {2\over\pi^2}\sum_{\scriptstyle 1\leq k\leq d \atop 
                       \scriptstyle k\hbox{ \sevenrm odd}} 
    {1\over k^2} \, .
$$
Since the variance is maximum at a unique point, Proposition
9.1.7 suggests that the dominating manifold in our problem is a
single point $p^*=p(1/2)/|p(1/2)|^2$. Its dimension is $k=0$.

We then calculate $I(A_1)=|p^*|^2/2 = 1/(2\sigma_d^2)$.
The level surface $\Lambda_{I(A_1)}$
is the sphere of radius $1/\sigma_d$, centered at the origin. Its
second fundamental form is $\sigma_d\Id_{d-1}$ at every point.

Let us now calculate the second fundamental form of $\partial
A_1$. To parameterize $\partial A_1$, we follow the construction
in section 9.1, with the simplification given in section 9.2 for
the special case of a radial function $I(\cdot )$. The vector
$\tau=p'/|p'|$ is a unit tangent vector to the curve 
$M=p([\,0,1\, ])$. As defined in section 9.1, let $\nu$ be a unit
normal vector to $M$ in $(T_pM+p\RR )\ominus T_pM$, that is
$$
  \nu = {p-\< p,\tau\>\tau\over \sqrt{|p|^2-\<p,\tau\>^2}} \, .
$$
Notice that $|p(s)|^2$ being maximal for $s=1/2$, the vectors
$p(1/2)$ and $\tau (1/2)$ are orthogonal, and  
$\nu(1/2)$ equals $p^*|p(1/2)|$.

Let $X_1,\ldots , X_{d-1}$ be an orthonormal moving frame in
$\RR^d\ominus (T_p\RR+p\RR)$. Using the simplification
pointed out after (9.2.4),
$$
  X(s,u_2,\ldots , u_{d-1}) = {\nu (s)\over \<\nu(s),p(s)\>}
  + \sum_{2\leq j\leq d-1} u_jX_j
  \eqno{(12.4.1)}
$$
defines a local parameterization of $\partial A_1$, or equivalently of
$\partial C_M$, near $p^*$ when $s$ is chosen close to $1/2$ and
$(u_2,\ldots , u_{d-1})$ close to $0$. To calculate the second
fundamental form of $\partial A_1$, recall that $p(s)$ is an outward
normal to $\partial A_1$ at $X(s,u_2,\ldots , u_{d-1})$ thanks to
Lemma 9.1.6. In particular,
$$
  T_X\partial C_M = \{\, p\,\}^\perp 
  = \span \{\, X_2,\ldots , X_{d-1}\,\} \oplus 
    \big( \span (p,\tau )\ominus p\RR\big) \, .
$$
We complete $X_2,\ldots , X_{d-1}$ into an orthonormal basis of
$T_X\partial C_M$ by adding the vector field
$$
  e={\Proj_{p^\perp} \tau\over |\Proj_{p^\perp} \tau| }
  = {1\over |p|} 
    { \tau |p|^2-\< \tau,p\>p\over \sqrt{|p|^2-\<\tau,p\>^2} }
  \, .
$$
We have
$$
  \d p \cdot X_j = 0 \qquad \hbox{ for } j=2,\ldots , d-1 \, ,
$$
expressing the fact that $\partial A_1$ is a ruled surface with
flat generators in the space spanned by $X_2,\ldots  , X_{d-1}$.
Thus, the second fundamental form of $\partial A_1$ vanishes on
$X_2,\ldots , X_{d-1}$. Its matrix in the basis $e,X_2,\ldots ,
X_{d-1}$ is then 
\setbox1=\hbox{\raise -5pt\hbox{\kern 2pt\ftitre 0}}
\wd1=0pt\ht1=0pt\dp1=0pt
\setbox2=\hbox{\raise 5pt\hbox{\kern 5pt\ftitre 0}}
\wd2=0pt\ht2=0pt\dp2=0pt
\setbox3=\hbox{$\mkern-40mu\mddots{5}{4}$}
\wd3=0pt\ht3=3pt\dp3=0pt
$$
  \pmatrix{ {\<\displaystyle\d p\cdot e,e\>\over\displaystyle|p|} 
                  &   & \box1     &  \cr
                  & 0 &        &   \cr
  \noalign{\vskip -.05in}
                  &   & \kern .25in\box3  &   \cr
  \noalign{\vskip .05in}
            \box2 &   &        & \kern -.15in 0 \cr}
  \, .
$$
It remains to calculate $\< \d p\cdot e,e\>$ on the dominating
manifold $p^*=X(1/2,0,\ldots , 0)$. 
The vectors $p(1/2)$ and $\tau(1/2)$
being orthogonal, $e(1/2,0,\ldots , 0)$ and $\tau(1/2)$ are
equal. Since $\tau$ is
a unit tangent vector to the curve $p([\, 0,1\,] )$, we have
$$\eqalign{
  \< \d p(1/2)\cdot e(1/2,0\ldots , 0),e(1/2,0,\ldots, 0)\>
  & = \< \d p\cdot \tau,\tau\> (1/2) \cr
  & =|\tau (1/2)|^2 =1 \, . \cr
  }
$$
In conclusion, on the dominating manifold $M$, and in the basis
$(e,X_2,\ldots , X_{d-1})$, the fundamental form of $\partial
A_1$ and $\Lambda_{I(A_1)}$ are
\setbox1=\hbox{\raise -5pt\hbox{\kern 5pt\ftitre 0}}
\wd1=0pt\ht1=0pt\dp1=0pt
\setbox2=\hbox{\raise 5pt\hbox{\kern 5pt\ftitre 0}}
\wd2=0pt\ht2=0pt\dp2=0pt
\setbox3=\hbox{\raise -5pt\hbox{\kern 5pt\ftitre 0}}
\wd3=0pt\ht3=0pt\dp3=0pt
\setbox4=\hbox{\raise 2pt\hbox{\kern 5pt\ftitre 0}}
\wd4=0pt\ht4=0pt\dp4=0pt
\setbox5=\hbox{$\mkern-40mu\mddots{5}{2}$}
\wd5=0pt\ht5=5pt\dp5=0pt
\setbox6=\hbox{$\mkern-40mu\mddots{5}{2}$}
\wd6=0pt\ht6=5pt\dp6=0pt
$$
  \Pi_{\partial A_1} =\pmatrix{ \sigma &   & \box1  &   \cr
                                       & 0 &        &   \cr
  \noalign{\vskip-.07in}
                                       &   & \kern .23in\box5  &   \cr
  \noalign{\vskip.1in}
                                 \box2 &   &        & \kern -.2in 0 \cr}
  \qquad\hbox{ and } \qquad
  \Pi_{\Lambda_{I(A_1)}} 
  = \pmatrix {\sigma & \box3  &   \cr
         \noalign{\vskip-.07in}
                     & \kern .23in\box6 &        \cr
         \noalign{\vskip.1in}
               \kern -.1in\box4 &        & \kern -.2in\sigma \cr}
 \, .
$$
In particular, $\Pi_{\partial A_1}-\Pi_{\Lambda_{I(A_1)}}$ is
diagonal. Its upper left entry vanishes. Because the dominating
manifold $\calD_{\partial C_M}$ is a point, $\Pi_{\partial
A_1}-\Pi_{\Lambda_{I(A_1)}}=G_{A_1}$; thus $\det\, G_{A_1}$ is null,
and assumption (7.5) does not hold.

The fact that the fundamental forms have the same upper left entry
expresses the fact that along the tangent direction $e(1/2,0,\ldots , 0)$, the
surfaces $\Lambda_{I(A_1)}$ and $\partial C_M$ pull apart very
slowly. The following pictures illustrate this fact when $d=2$ and
$d=3$.

\bigskip

\setbox1=\hbox to 2in{\vbox to 2.45in{\ninepointsss
  \hsize=2in
  \kern 2in
  \includegraphics{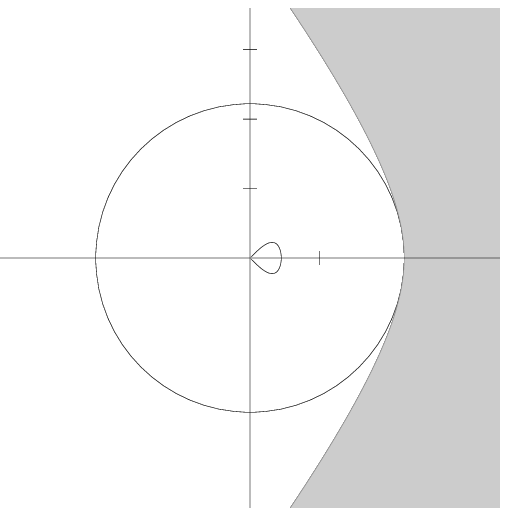}
  \anote{.9}{.9}{0}
  \anote{.85}{1.25}{1}
  %\anote{.85}{1.5}{2}
  %\anote{.85}{1.75}{3}
  \anote{1.25}{.85}{1}
  \anote{1.7}{1.5}{$A_1$}
  \anote{.45}{.4}{$\Lambda_{I(A_1)}$}
  \vfill
  \centerline{$d=2$}
  \kern -.05in
  \centerline{(the little loop next to the 
              origin is $M=p([0,1])$.)}
}\hfill}

\medskip
\centerline{\box1}

\bigskip

\setbox1=\hbox to 2.3in{\vbox to 2.2in{\ninepointsss
  \hsize=2.3in
  \kern 2.1in
  \includegraphics{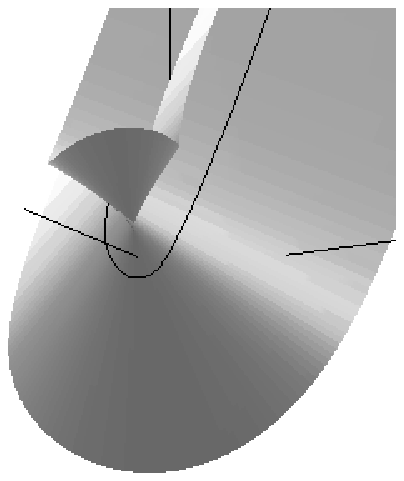}
  \anote{1}{2.05}{$e_3\RR$}
  \anote{.3}{1.2}{$e_1\RR$}
  \anote{2.05}{1.1}{$e_2\RR$}
  \vss
  %\centerline{$d=3$}
  \kern -.05in
}\hfill}

\setbox2=\hbox to 2.3in{\vbox to 2.2in{\ninepointsss
  \hsize=2.3in
  \kern 2.1in
  \includegraphics{Graph/ch12_4_fig2l.eps}
   \anote{1}{2.05}{$e_3\RR$}
  \anote{.3}{1.2}{$e_1\RR$}
  \anote{2.05}{1.1}{$e_2\RR$}
  \vss
  %\centerline{$d=3$}
  \kern -.05in
}\hfill}

\centerline{\leavevmode\kern -.2in\box1\hfill\box2}

\centerline{\ninepointsss $d=3$}

\medskip

\centerline{
\vbox{\hsize=4in\ninepointsss\noindent The boundary of $A_1$
  is a ruled surface. This picture shows the ruled surface obtained
  from the parameterization $X(s,u_2)$, with no
  constraints on $u_2$. The actual set 
  $A_1$ is in the convex hull of the piece containing the origin.
  The curve drawn on the surface is 
  $X(s,0)=\nu(s) /\<\nu(s),p(s)\>$.}
}

\bigskip

\setbox1=\hbox to 2.3in{\vbox to 2.2in{\ninepointsss
  \hsize=2.3in
  \kern 2.1in
  \includegraphics{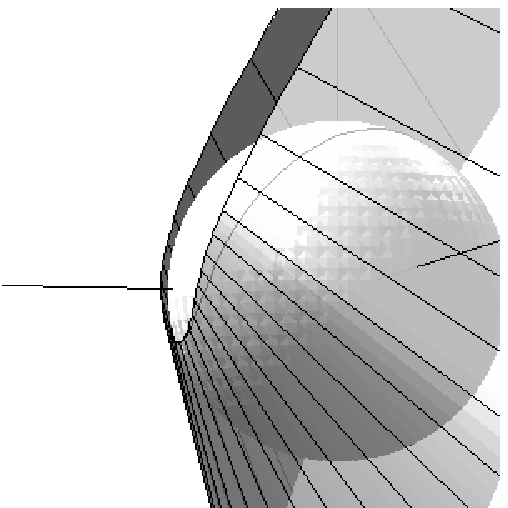}
  \anote{1.2}{2.05}{$e_3\RR$}
  \anote{.3}{0.95}{$e_1\RR$}
  \anote{2.05}{1.1}{$e_2\RR$}
  \vss
  %\centerline{$d=3$}
  \kern -.05in
}\hfill}

\setbox2=\hbox to 2.3in{\vbox to 2.2in{\ninepointsss
  \hsize=2.3in
  \kern 2.1in
  \includegraphics{Graph/ch12_4_fig3l.eps}
   \anote{1.2}{2.05}{$e_3\RR$}
  \anote{.32}{0.95}{$e_1\RR$}
  \anote{2.05}{1.1}{$e_2\RR$}
  \vss
  %\centerline{$d=3$}
  \kern -.05in
}\hfill}

\centerline{\leavevmode\kern -.2in\box1\hfill\box2}
\bigskip

\centerline{
\vbox{\hsize=4in\ninepointsss\noindent 
  In this picture, we can see the sphere $\Lambda_{I(A_1)}$,
  centered and of radius $1/\sigma_3$.
  The transparent surface is $\partial A_1$,
  parametrized by $X(s,u_2)$, but with $u_2$ constrained to be
  negative. Its generators are straight lines, parameterized by
  $u_2$, shown here for $u_2$ negative. 
  Thus, the surface $\partial A_1$ is cut along the
  curve $X(s,0)$. The projection of this curve on the sphere
  $\Lambda_{I(A_1)}$ can be seen. The set $\partial A_1$
  has a unique contact point, $p^*$ with the level set
  $\Lambda_{I(A_1)}$. But, as it can be seen along the curve
  $X(s,0)$, the boundary $\partial A_1$ pulls away very slowly
  from $\Lambda_{I(A_1)}$ near $p^*$.}
}

\bigskip

In the current situation, it is easy to obtain the desired equivalence
for the tail of the supremum, because $I(\cdot )$ is homogeneous. If
this tail were behaved like that of the process at its point of
maximal variance, we would have a tail equivalence in
$$
  P\{\, B_d(1/2)\geq t\,\}
  = 1-\Phi (t/\sigma_d) 
  \sim {1\over \sqrt{2\pi}} {\sigma_d\over t} 
       \exp \Big( -{t^2\over 2\sigma_d^2}\Big) \, .
$$
The tail behavior turns out to be much more surprising.  Even and odd
dimensions $d$ yield different exponents in the polynomial term. The
reason is that for even dimensions, the set $\partial C_M$ pulls away
more slowly along the direction $e(1/2,0,\ldots ,0)$. The contact
between the two surfaces is of order $1$ for odd dimensions, and of
order $3$ for even ones.

\bigskip

\state{12.4.1. THEOREM.} {\it%
  As $t$ tends to infinity, the ratio
  $$
    P\{\, \sup_{0\leq s\leq 1}B_d(s)\geq t\,\} 
    \big/\, \big(1-\Phi (t/\sigma_d)\big)
  $$
  is equivalent to
  \vskip .05in
  
  \noindent (i) $\sqrt{d+1}/\sigma_d$ if $d$ is odd;

  \noindent (ii) $ \sqrt t \, \Gamma (1/4)
    \Big({\displaystyle 6\over\displaystyle d(d+1)}\Big)^{1/4} 
     {\displaystyle \sqrt d\over\displaystyle 2\pi\sqrt 2\sigma_d^2}$
     if $d$ is even and $d\neq 1$;

  \noindent (iii) $1$ if $d=1$.
}

\bigskip

\state{REMARK.} We will see that $\lim_{d\to\infty}\sigma_d^2=1/4$.
Since $d\mapsto \sigma_d^2$ is increasing, we have $\sigma_d^2<1/4$
for any integer $d$. Thus the constant in statement (i) of Theorem 12.4.1
is always larger than $1$. This is in agreement with the failure
of assumption 5.3. In the same spirit, the quantity involved
in statement (ii) is also larger than $1$ when $t$ is large enough.

\bigskip

\stateit{Proof of Theorem 12.4.1.} We first proof statement (iii).
Since 
$$
  B_1(s) = {\sqrt 2\over\pi} X_1\sin (\pi s) \, ,
$$
its supremum is $0$ if $X_1$ is negative, and $\sqrt 2 X_1/\pi$ otherwise.
The result follows. 

From now on, assume that $d\geq 2$. The change of variable $y=tc$ 
allows to write $$\eqalign{
  P\{\, \sup_{0\leq s\leq 1} B_d(s)\geq t\,\}
  & = {1\over (2\pi )^{d/2}} \int_{tA_1}e^{-|y|^2/2} \d y \cr
  & = {t^d\over (2\pi )^{d/2}} \int_{A_1} e^{-t^2|x|^2/2} \d x\, .
    \cr}
$$
Thus, to prove Theorem 12.4.1, we need to estimate the integral
$$
  \int_{A_1} e^{-\lambda |x|^2/2} \d x \, ,
$$
for large values of $\lambda$. It is plain from standard results
on Laplace's method, or from chapter 7, that we can restrict the
integration to the intersection of $A_1$ with an arbitrary
neighborhood of the dominating manifold $p^*$. Making use
of the construction in section 9.1, 
we can write any point $x$ of $A_1$ near $p^*$ as
$$
  x(s,u,v) = {\nu(s)\over \< \nu(s),p(s)\> }
  + \sum_{2\leq j\leq d-1} u_j X_j(s) + v p(s) \, , \qquad
  v\geq 0 \, ,
$$
where $u=(u_2,\ldots , u_{d-1})$ is in a neighborhood of the
origin. Because $X_2,\ldots , X_{d-1}, \nu$ is an orthonormal
basis and $p$ is orthogonal to $X_2,\ldots , X_{d-1}$,
$$
  |x|^2 = {1\over \< \nu,p\>^2} + |u|^2 + v^2|p|^2 + 2 v \, .
$$
To calculate the Jacobian of the change of variable $x\leftrightarrow
(s,u,v)$, is easy. Write
$$\eqalign{
  {\dispar x\over\dispar u_j} &= X_j \, ,\qquad
  {\dispar x\over\dispar v} = p  \, , \cr
  {\dispar x\over\dispar s}   &= {\dispar\over\dispar s}
     \Big( {\nu(s)\over \< \nu(s),p(s)\>}\Big)
     + \sum_{2\leq j\leq d-1} u_j {\dispar \over\dispar s} X_j(s)
     + v p'(s) \, . \cr
  }
$$
Introducing $w=p|p'|^2-\<p,p'\>p'$, we see that $\nu/\<
\nu,p\>=w/\<w,p\>$. Thus,
$$
  \Big( {\nu\over \< \nu,p\>}\Big)'
  = {w'\over \< w,p\>} - w {\<w',p\>+\<w,p'\>\over \< w,p\>^2}\, .
$$
It is straightforward to calculate
$$
  w' = 2p\< p',p''\> - \< p,p''\>p' - \< p,p'\>p'' \, .
$$
Since $|p(s)|^2$ is maximal at $s=1/2$, the tangent vector
$p'(1/2)$ is orthogonal to $p(1/2)$. The specific form of
$p(\cdot )$ in this problem yields that for any $j$, the
derivative $\d^jp/\d s^j$ at $s=1/2$ involves only the vectors
$e_k$ with $k$ odd if $j$ is, and $k$ even otherwise.
Consequently, $p'(1/2)$ and $p''(1/2)$ are orthogonal, and
$$\eqalign{
  & w(1/2) = p|p'|^2(1/2)\, , \qquad w'(1/2) = -\< p,p''\>p'(1/2) \cr
  \noalign{\vskip .05in}
  & {\dispar\over\dispar s}
    \Big( {\nu(s)\over\< \nu(s),p(s)\>}\Big) (1/2)
    = -{\<p,p''\>p'\over |p|^2|p'|^2} (1/2) \, . \cr
  }
$$
Thus, in the orthonormal basis $\tau, X_2,\ldots , X_{d-1}$, 
the Jacobian matrix is 
\setbox1=\hbox{\raise -5pt\hbox{\kern 5pt\ftitre 0}}
\wd1=0pt\ht1=0pt\dp1=0pt
\setbox2=\hbox{\raise 5pt\hbox{\kern 5pt\ftitre 0}}
\wd2=0pt\ht2=0pt\dp2=0pt
$$\pmatrix{ -{\displaystyle\<p,p''\>\over\displaystyle |p|^2|p'|}(1/2)% 
                                               &   &  \box1   &   \cr
                                               & 1 &      &   \cr
                                               &   &\ddots&   \cr
                                        \box2  &   &      & 1 \cr}
  + O(|x-p^*|) 
  \qquad\hbox{ as } x\to p^*\, .
$$

Now, let $\epsilon$ be an arbitrary positive real number. In what
follows, $\eta$ denotes a positive real number, which we will
choose as small as needed. Denote by $D(\eta)$ the domain
$$
  D(\eta) = \{ x(s,u,v)\in \RR^d \, : \, |s-1/2|\leq \eta \, , \,
               |u|\leq\eta \, ,\,
               0\leq v\leq \eta\,\} \, .
$$
From our evaluation of the Jacobian,
if $\eta$ is small enough
$$\displaylines{
    \int_{A_1\cap D(\eta) } e^{-\lambda |x|^2/2}\d x
  \hfill\cr\quad
    \leq (1+\epsilon ) 
    \Big| {\<p,p''\>\over |p|^2|p'|} (1/2)\Big| \, \times
  \hfill\cr\hfill
    \int_{D(\eta )} 
    \exp\bigg( -{\lambda\over 2} 
              \Big( {1\over\< \nu(s),p(s)\>^2} + |u|^2+v^2|p(s)|
                    +2v\Big)
        \bigg) \, \d u\, \d v\, \d s\, .
  }
$$
We first perform the integration in $u$, obtaining the upper bound
$$\displaylines{
    (1+\epsilon ) {|\<p,p''\>|\over |p|^2|p'|}(1/2)
    {(2\pi)^{(d-2)/2}\over \lambda^{(d-2)/2}}\times
  \hfill\cr\hfill
    \int_{|s-1/2|\leq\eta}\int_{0\leq v\leq\eta}
                    \exp \bigg( -{\lambda \over 2}
                    \Big( {1\over\<\nu(s),p(s)\>^2} + v^2|p(s)|^2+2v
                    \Big) \bigg) \, \d s\, \d v 
  \cr}
$$
To perform the integration in $v$, we write
$$\displaylines{\qquad
  \exp\Big( -{\lambda\over 2} \big( v^2|p(s)|^2+2v\big) \Big)
  \cut
  = \Big( \lambda \big( v|p(s)|^2+2\big)\Big)^{-1}
  {\d \over \d v}\exp\Big( {\lambda\over 2} 
  \big( v^2|p(s)|^2+2v\big)\Big)
  \qquad\cr}
$$
and integrate by parts. We obtain
$$
  \int_{0\leq v\leq \eta}
  \exp\Big( -{\lambda\over 2} \big( v^2|p(s)|^2+2v\big) \Big)
  \, \d v
  \sim {1\over\lambda}
$$
as $\lambda$ tends to infinity.
This yields the upper bound
$$
  (1+\epsilon )^2 {|\<p,p''\>|\over |p|^2|p'|} (1/2)
  {(2\pi )^{(d-2)/2}\over \lambda^{d/2}}\int_s
  \exp\Big( -{\lambda\over 2}\<\nu(s),p(s)\>^{-2}\Big) \d s
  \, .
$$
To estimate this last integral boils down to using the classical
Laplace method. We introduce
$$
  \delta = \< \nu,p\>^2 = |p|^2-\< p,\tau\>^2 \, .
$$
We then have
$$\displaylines{
  \int_{|s-1/2|\leq \eta} \exp\Big( -{\lambda\over 2}
  \< \nu(s),p(s)\>^{-2}\Big) \, \d s
  \hfill\cr\noalign{\vskip .1in}\hfill
  = \exp \Big( -{\lambda\over 2} \delta (1/2)\Big)
    \int_{|h|\leq \eta} \exp\Big( -{\lambda\over 2}
    { \delta (1/2)-\delta (1/2+h)
      \over \delta (1/2)\delta (1/2+h) }\Big) \, \d h \, .
  \,\,(12.4.2)\cr}
$$
We then need to obtain a Taylor expansion for $\delta (\cdot )$ 
near $1/2$. To this end, for any integer $m$, we define
$$
  S_m=\sum_{1\leq k\leq d}(-1)^k k^m \, .
$$
Since
$$
  |p(s)|^2 
  = {2\over \pi^2}\sum_{1\leq k\leq d} {\sin^2 (k\pi s)\over k^2}
  = {1\over\pi^2} \sum_{1\leq k\leq d} {1-\cos (2k\pi s)\over k^2}
  \, ,
$$
we easily obtain the derivatives $(\d^m/\d s^m)\big(|p(s)|^2\big)$
at $s=1/2$, for $m=0,1,\ldots , 4$. Using Taylor's formula, 
we infer that
$$
  |p(1/2+h)|^2 = |p(1/2)|^2 + 2h^2S_0 - {2\over 3} \pi^2h^4S_2
  + O(h^5) \, .
$$
We also have
$$
  \<p,p'\>(s) 
  =  {2\over \pi}\sum_{1\leq k\leq d} 
     {\sin (k\pi s)\cos (k\pi s)\over k} 
  = {1\over \pi} \sum_{1\leq k\leq d} {\sin (2k\pi s)\over k}\, .
$$
Hence, a simple calculation of the derivatives at $s=1/2$ and
an application of Taylor's formula give
$$
  \< p,p'\> (1/2+h)
  = 2h S_0 - {4\over 3} \pi^2h^3S_2 + O(h^5) \, .
$$
Finally, 
$$
  |p'(s)|^2
  = 2\sum_{1\leq k\leq d} \cos^2 (k\pi s) 
  = d + \sum_{1\leq k\leq d} \cos (2k\pi s) \, ,
$$
from which we deduce $|p'(1/2)|^2=d+S_0$, and
$$
  |p'(1/2+h)|^2 
  = |p'(1/2)|^2 - 2\pi^2h^2S_2 + O(h^4) \, .
$$
Equipped with these expressions, a little algebra gives
$$\displaylines{
    \delta (1/2+h)=\delta (1/2) 
    + 2S_0\Big( 1-{2S_0\over |p'(1/2)|^2}\Big) h^2
  \cut
    + \pi^2S_2 \bigg( -{2\over 3} - {S_0\over |p'(1/2)|^2}
    \Big( -{16\over 3}+{8S_0\over |p'(1/2)|^2}\Big)\bigg)h^4 
    + O(h^5) \, .
  \cr
  }
$$
The interesting fact is now that $S_0$ equals $-1$ if $d$ is odd, 
and equals $0$ if $d$ is even. Thus, if $d$ is odd,
$$
  \delta (1/2+h)
  =\delta (1/2) - 2\Big( 1+{2\over |p'(1/2)|^2}\Big) h^2 + O(h^4)
  \, ,
$$
while if $d$ is even,
$$
  \delta (1/2+h) = \delta (1/2) - {2\pi^2S_2\over 3} h^4 + O(h^5)
  \, .
$$
Notice that when $d$ is even, $S_2$ is positive since
$$\eqalign{
  \sum_{1\leq k\leq 2m} (-1)^kk^2
  & = \sum_{1\leq k\leq m} (2k)^2-(2k-1)^2
    =\sum_{1\leq k\leq m} (4k-1) \cr
  \noalign{\vskip .05in}
  & = m(2m+1) \, . \cr
  }
$$
Hence, $\delta (s)$ is maximal at $s=1/2$ whatever the parity of
$d$ is.

To conclude the proof of Theorem 12.4.1, assume first that $d$ is
odd. Then, for $\eta$ small enough, (12.4.2) is less than
$$\displaylines{\qquad
  \int_{|h|\leq \eta} \exp \bigg( -\lambda 
  \Big( 1+ {2\over |p'(1/2)|^2} \Big)
  {h^2\over \delta (1/2)^2(1+\epsilon)^2} \bigg) \, \d h
  \cut
  \sim {\displaystyle \sqrt{2\pi} (1+\epsilon)\delta (1/2)
        \over\displaystyle \sqrt{1+2|p'(1/2)|^{-2}} }
  {1\over\sqrt \lambda } \, .
  \qquad\cr}
$$
All the arguments we used to obtain this upper bound can be used
to obtain a lower bound, essentially by changing $\epsilon$ to
$-\epsilon$. Since $\epsilon$ is arbitrary, combining all 
the estimates yields
$$
  P\{\, \sup_{0\leq s\leq 1} B_d(s)\geq t\,\}
  \sim
  \sqrt{d+1\over 2\pi} {1\over t} 
  \exp\Big( -{t^2\over 2\sigma_d^2}\Big)
$$
When $d$ is even, for $\eta$ small enough, the integral 
in (12.4.2) is less than
$$\displaylines{\qquad
    \int_{|h|\leq \eta} \exp \bigg( -{\lambda\over 2} 
    {2\pi^2S_2\over 3} {h^4\over \delta (1/2)^2 (1+\epsilon)^2 } \bigg)
    \, \d h
  \cut
    \sim {\Gamma (1/4)\over \lambda^{1/4}} 
    \Big({3\over S_2}\Big)^{1/4}\sqrt{(1+\epsilon)\delta(1/2)\over\pi}
    \, .
  \qquad\cr}
$$
It then follows that
$$
  P\{\, \sup_{0\leq s\leq 1}B_d(s)\geq t\,\}
  \sim
  {\sqrt d\over 4\pi^{3/2}} {\Gamma(1/4)\over\sigma_d}
  \Big( {6\over d(d+1)}\Big)^{1/4} {1\over \sqrt t}
  \exp\Big( -{t^2\over 2\sigma_d^2}\Big) \, .
$$
This concludes the proof of Theorem 12.4.1.\hfill$\qed$

\bigskip

To conclude this section let us make a remark on $\sigma_d^2$. 
Recall that $\sum_{k\geq 1} (2k+1)^{-2} = \pi^2/8$. Consequently,
$\lim_{d\to\infty} \sigma_d^{-2} = 2$, which is good, given
(12.4.1). We can obtain a more precise estimate of
$\sigma_d^2$. Since
$$
  \int_m^\infty {\d x\over (2x+1)^2}
  \leq \sum_{k\geq m} {1\over (2k+1)^2} 
  \leq \int_m^\infty {\d x\over (2x-1)^2} \, ,
$$
and
$$
  \int_m^\infty {\d x\over (2x+1)^2}
  \leq \int_m^\infty {\d x\over 4x^2}
  \leq \int_m^\infty {\d x\over (2x-1)^2} \, ,
$$
we have
$$\eqalign{
  \Big|\sum_{k\geq m} {1\over (2k+1)^2}-{1\over 4m}\Big|
  & \leq \int_m^\infty {1\over (2x-1)^2} - {1\over (2x+1)^2}
    \, \d x \cr
  & = {1\over 2(4m^2-1)} \, .\cr
  }
$$
Consequently,
$$
  \Big| \sigma_d^2-{1\over 4} - {1\over 2\pi d}\Big|
  \leq {1\over \pi^2 (4d^2-1)} \, .
$$

\bigskip

Though the constant in the exponential term in Theorem 12.4.1,
namely $1/(2\sigma_d^2)$, has the right limit as $d$ tends to
infinity, we cannot take the limit of the polynomial term and
recover (12.4.1). This is caused by the slow convergence
of $B_d$ to the Brownian bridge.

Again, the whole message is to be rather cautious with these
approximations. Strange things may happen and further examination 
is certainly needed if they are to be used in serious applications.

%%%%%%%%%%%%%% DO NOT PRINT %%%%%%%%%%%%%%%%%%%%%
\do_not_print{

We will work almost exclusively with numerical calculations,
though in some instances close formula may be obtained. 

First, we calculate\note{est-ce que la figure est trop haute?}
$$
  \Var B_d(s) = {2\over \pi^2}\sum_{1\leq k\leq d}
  {\sin^2 (k\pi s)\over k^2 } \, .
$$

\setbox1=\hbox to 3.1in{\vbox to 1.8in{\ninepointsss
  \hsize=3.1in
  \kern 1.5in
  \includegraphics{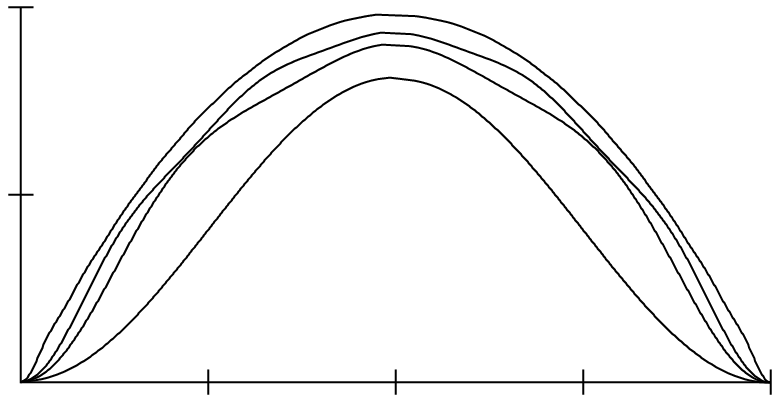}
  \anote{-.1}{0}{0}
  \anote{0.65}{-.15}{0.25}
  \anote{1.45}{-.15}{0.5}
  \anote{2.15}{-.15}{0.75}
  \anote{2.97}{-.15}{1}
  \anote{-.25}{0.75}{0.25}
  \anote{-.2}{1.5}{0.5}
  \vfill
  \centerline{$\Var(B_d(s))$ for $d=1$, $3$, $5$, $20$.}
}\hfill}
\centerline{\box1}

\kern 0.1in
For $d=5$, $10$, $15$, $20$, the maximum of the variance is
found numerically and given in the following table

\centerline{TABLE maximum \`a  d=5 10 15 20}

}
%%%%%%%%%%%%%%%%% END DO NOT PRINT %%%%%%%%%%%%%%%%%%%

\bigskip

\noindent{\fsection 12.5. Polar processes on boundary of convex
sets.}
%\chapternumber={12.5}
\section{12.5}{Polar processes on boundary of convex
sets}

\bigskip

Let $C$ be a bounded convex set in $\RR^d$, with nonempty
interior. The purpose of this section is to construct a simple process
on its boundary. This construction is suggested by Theorem 7.1 and the
results of chapter 9.

We first need to make some remarks on densities proportional to
$e^{-I}$ with $I$ strictly convex and $\alpha$-homogeneous.  Recall
that $\Lambda_c$ denotes the level line $I^{-1}(\{\, c\,\})$. Every
nonzero point $x$ of $\RR^d$ can be written in a unique way as
$x=s\lambda$ for some positive $s$ and $\lambda$ in $\Lambda_1$. In
this $(\lambda ,s)$-coordinate system, the measure $e^{-I}$ can be
rewritten as
$$
  s^{d-1} e^{-s^\alpha }
  |\Proj_{T_\lambda \Lambda_1^\perp}\lambda |
  \, \d s \, \d\calM_{\Lambda_1}(\lambda ) \, .
  \eqno{(13.5.1)}
$$
Conversely, such a measure corresponds to a log-concave and
log-$\alpha$-ho\-mo\-ge\-neous measure on $\RR^d$.

Going back to the convex set $C$ given at the beginning of this
section, assume that it also contains the origin. Let $\Lambda$ be the
polar reciprocal of $\partial C$. For $\lambda$ in $\Lambda$ and
nonnegative $s$, define $I(s\lambda )=s^\alpha$. Then,
$\Lambda=\Lambda_1$ for this specific function $I$. Equip $\RR^d$ with
the log-concave density proportional to (13.5.1), and let $X$ be a
random vector having this density. We can consider the process
$p\in\partial C\mapsto \< X,p\>\in\RR$. We call this process a polar
process on $\partial C$.

The tail distribution of its supremum is given by the following
result.

\bigskip

\state{12.5.1. THEOREM.} {\it%
  For the polar process on $\partial C$ defined above,
  and $1/c=\int_{\RR^d} e^{-I(x)}\, \d x$,
  $$
    P\{\, \sup_{p\in\partial C}\< X,p\> \geq t\,\}
    \sim
    {c\over\alpha} e^{-t^\alpha}t^{d-\alpha}\int_\Lambda
    |\Proj _{(T_\lambda\Lambda )^\perp}\lambda | 
    \, \d\calM_\Lambda (\lambda )
  $$
  as $t$ tends to infinity.
  }

\bigskip

\stateit{Proof.} We apply Theorem 7.1. In view of section 9.1,
the dominating manifold is $\Lambda$, of dimension $k=d-1$.
Moreover, by construction, $\partial A_1=\Lambda_1$. It follows
from Theorem 7.1 that the tail equivalent is of the form
$$
  c e^{-t^\alpha}t^{d-\alpha}
  \int_\Lambda {\d\calM_\Lambda\over |\D I|} \, .
$$
To calculate $\D I$, notice that it is the outward normal to 
$\Lambda$. Its norm is obtained through the identity
$$
  {\d\over \d s} I(\lambda s){\big|}_{s=1}
  = \D I(\lambda )\cdot \lambda 
  = {\d\over \d s}\big(s^\alpha I(\lambda )\big){\big|}_{s=1}
  = \alpha I(\lambda )
  = \alpha \, .
$$
It implies
$$
  |\D I(\lambda )| 
  = \alpha /|\Proj_{(T_\lambda\Lambda )^\perp}\lambda| \, ,
$$
and the result follows.\hfill$\qed$

\bigskip

The correspondence between points $p$ on $\partial C$ and points
$\lambda$ on the polar reciprocal allows us to write
$|\Proj_{(T_\lambda\Lambda )^\perp}\lambda | = 1/|p|$. Thus the
integral in Theorem 12.5.1 can also be rewritten as
$$
  \int_\Lambda {1\over |p(\lambda )|}\d\calM_\Lambda (\lambda ) 
  \, ,
$$
which is a nice formula.

The following pictures are three simulations of a polar process on the
ellipsoid $x^2+2y^2=1$ in $\RR^2$, and $\alpha =1$.  For the pictures
on the left, at every point $p$ of $\partial C$, we draw a segment in
the normal direction to $T_p\partial C$, with length equal to the
value of the process at $p$. The random point $X$ is indicated by the
star
\setbox1=\hbox to 9pt{\includegraphics{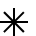}\hfill}\box1 . 
On the right hand side pictures, the process is represented as a curve
held over the ellipsoid, the height of the curve being the value of
the process; these are orthographic projections.

\finetune{\vfill\eject ~ \vskip -.3in}
\setbox1=\hbox to 4in{\vbox{\ninepointsss
  \kern 2.3in
  \includegraphics{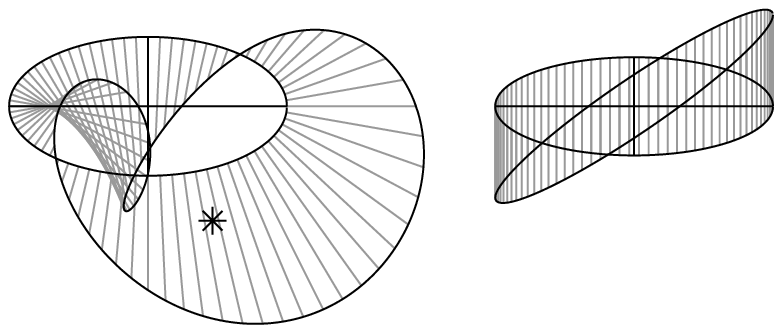}
  \vfill
}\hfill}
\line{\kern -.5in\box1\hfill}

\finetune{\kern -1.2in}

\setbox1=\hbox to 4in{\vbox to 3in{\ninepointsss
  \kern 3in
  \includegraphics{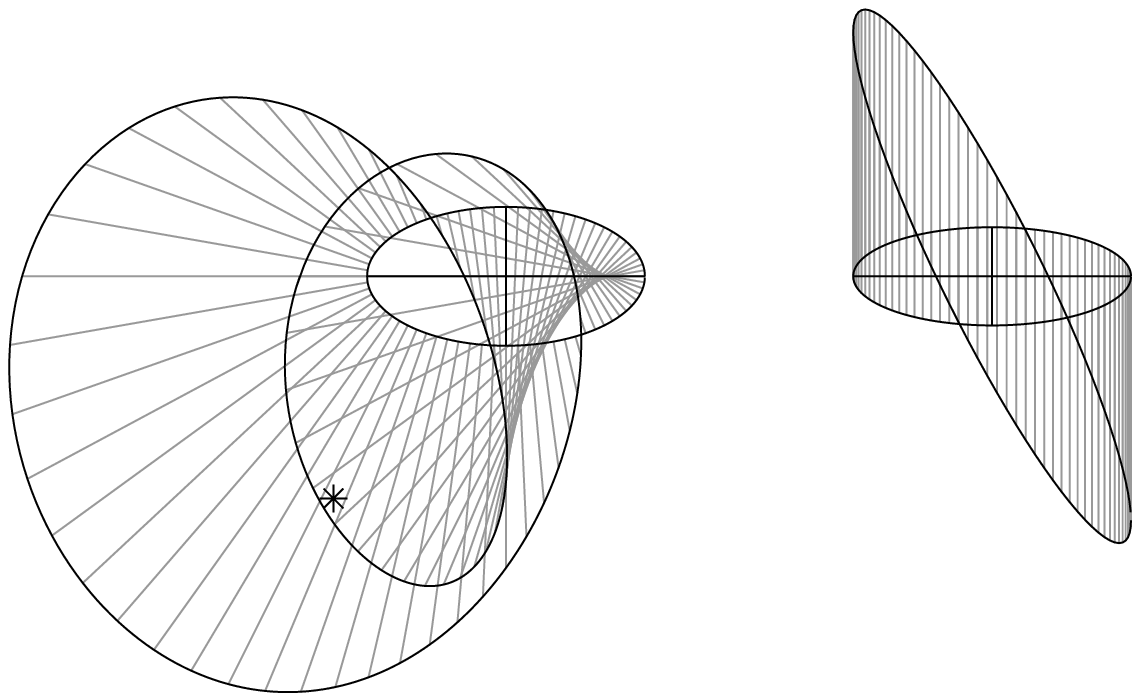}
  \vfill
}\hfill}
\line{\kern -.5in\box1\hfill}

%\finetune{\vfill\eject}
\kern -.3in
\setbox1=\hbox to 4in{\vbox to 4in{\ninepointsss
  \kern 4in
  \includegraphics{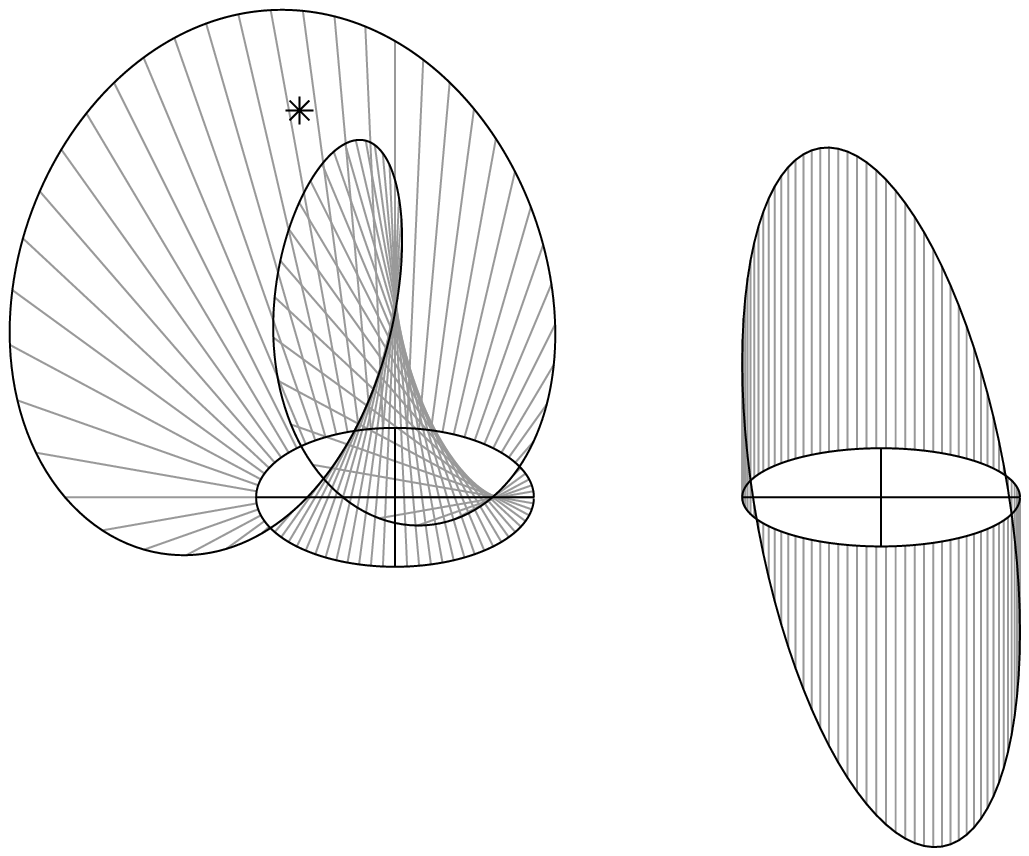}
  \vfill
}\hfill}
\line{\kern -.5in\box1\hfill}

\finetune{\vfill\eject}
\bigskip

\centerline{\fsection Notes}
\section{12}{Notes}

\bigskip

\noindent I do not think the results in this chapter are too serious!
However, I believe the message of caution that some carry is
important. In a different vein, Feynman supposedly said that
when you have a theory, you should show where it works and
where it does not. That may be the point of section 12.2.

Section 12.5 is just a way to generate pretty pictures of
random caustics. I do not know any applications of this
construction.

\vfill
\eject

{
  \def\rightheadline{\ninepointsss Appendix \the\chapternumber\hfill\folio}
  \def\leftheadline{\ninepointsss\folio\hfill Appendix \the\chapternumber}
  ~\vskip 28pt
\chapter{Appendix 1}{Gaussian and Student tails}
{\fchapter
  \line{\hfil Appendix 1.\quad Gaussian and}
  \vskip .1in
  \line{\hfil Student tails}}
\chapternumber={1}
%\vskip 2pt
%\hrule height 2pt depth 2pt 

\vskip 54pt

\noindent This appendix collects a few standard results related to the
tails of the Gaussian and Student distributions. We write
$$
  \underline s_\alpha (x) 
  = K_{s,\alpha} \Big( 1+{x^2\over\alpha}\Big)^{-(\alpha+1)/2}
  \, , \qquad x\in\RR \, , \, \alpha >0
$$
for the Student density (with $\alpha+1$ degrees of freedom), the
constant $K_{s,\alpha}$ ensuring that $s_\alpha$ integrates to $1$
over the real line. We denote by
$$
  \underline S_\alpha (x)=\int_{-\infty}^x s_\alpha (y) \d y
$$
the Student cumulative distribution function. Its tail is given
in the following result.

\bigskip

\state{A.1.1. LEMMA.} {\it%
  We have
  $$
    1-\underline S_\alpha (x) 
    = {K_{s,\alpha}\alpha^{(\alpha -1)/2}\over x^\alpha}
    + O\Big( {1\over x^{\alpha+2}}\Big) 
    \qquad{ as }\ x\to\infty \, .
  $$}

\bigskip

\stateit{Proof.} Notice that as $y$ tends to infinity,
$$\eqalign{
  {1\over 
   \Big(1+{\displaystyle y^2\over\displaystyle\alpha}
        \Big)^{\scriptstyle(\alpha+1)/2}}
  - {1\over\Big({\displaystyle y^2\over\displaystyle\alpha}
                \Big)^{\scriptstyle (\alpha+1)/2}}
  & = { 1-\Big( {\displaystyle\alpha\over\displaystyle y^2}+1
             \Big)^{\scriptstyle (\alpha+1)/2}
      \over \Big( 1+{\displaystyle y^2\over\displaystyle\alpha}
                  \Big)^{\scriptstyle (\alpha+1)/2}} \cr
  & = O\Big( {1\over y^{\alpha+3}}\Big) \, . \cr
  }
$$
Consequently, as $x$ tends to infinity,
$$\eqalignno{
  {1-S_\alpha(x)\over K_{s,\alpha}}
  = \int_x^\infty {\d y\over 
    \Big(1+{\displaystyle y^2\over\displaystyle\alpha}
         \Big)^{\scriptstyle\alpha+1\over \scriptstyle2}}
  & = \int_x^\infty {\alpha^{(\alpha+1)/2}\over y^{\alpha+1}}
    + O\Big( {1\over y^{\alpha+3}}\Big) \d y \cr
  & = {\alpha^{(\alpha-1)/2}\over x^\alpha}
    + O\Big( {1\over x^{\alpha +2}}\Big) &$\qed$\cr
  }$$

\bigskip

Let us now consider a Student-like cumulative distribution 
function $S_\alpha$, i.e., such that
$$
  S_\alpha (-x)\sim 1-S_\alpha (x) 
  \sim {K_{s,\alpha}\alpha^{(\alpha-1)/2}\over x^\alpha}\, ,
$$
as $x$ tends to infinity, and
where the constant $K_{s,\alpha}$ can be any fixed positive number.
We can obtain an asymptotic formula for the high quantiles.

\bigskip

\state{A.1.2. LEMMA.} {\it%
  The following holds,
  $$
    (1-S_\alpha )^\leftarrow (u) 
  \sim { K_{s,\alpha}^{1/\alpha}\alpha^{(\alpha -1)/2\alpha}\over
      u^{1/\alpha} } 
  \qquad\hbox{ as } u\to 0 \, .
  $$}

\bigskip

\stateit{Proof.} Let $x$ tends to infinity and $u$ tends to $0$ 
such that $u=1-S_\alpha (x)$, that is 
$x=(1-S_\alpha )^\leftarrow (u)$. From the Student-like tail,
we infer that
$$
  u\sim {K_{s,\alpha}\alpha^{(\alpha-1)/2}\over x^\alpha} \, ,
$$
which is the result.\hfill$\qed$

\bigskip

We can obtain similar results for the Gaussian distribution with
cumulative distribution function
$$
  \Phi (x) =\int_{-\infty}^x {e^{-y^2/2}\over \sqrt{2\pi}} \d y
  \, .
$$

\bigskip

\state{A.1.3. LEMMA.} {\it%
  We have
  $$
    1-\Phi (x) 
    = {e^{-x^2/2}\over x\sqrt{2\pi}}\big( 1 +o(1)\big)
    \qquad\hbox{ as } x\to\infty \, .
  $$}

\bigskip

\stateit{Proof.} Integrate by parts to obtain
$$\eqalign{
  \sqrt{2\pi}\big( 1-\Phi (x)\big)
  & = \int_x^\infty {1\over y} ye^{-y^2/2} \d y \cr
  & = {e^{-x^2/2}\over x} - \int_x^\infty {e^{-y^2/2}\over y^2}
    \d y \cr
  & = {e^{-x^2/2}\over x} - {e^{-x^2/2}\over 3x^3}
    + \int_x^\infty {e^{-y^2/2}\over 4y^4} \d y \, . \cr
  }$$
Consequently,
$$
  {e^{-x^2/2}\over\sqrt{2\pi}}
  \Big( {1\over x}-{1\over 3x^2}\Big)
  \leq 1-\Phi (x)
  \leq {e^{-x^2/2}\over\sqrt{2\pi}}{1\over x} \, ,
$$
and the result follows.\hfill$\qed$

\bigskip

We can now obtain an asymptotic approximation for high
quantiles. The second statement in the following lemma is
instrumental in the sequel.

\bigskip

\state{A.1.4. LEMMA.} {\it%
  We have
  $$
    \Phi^\leftarrow (1-u) = \sqrt{2\log 1/u}
    -{\log\log 1/u\over 2\sqrt{2\log 1/u}}
    -{\log (2\sqrt{\pi})\over\sqrt{2\log 1/u}} \big(1+o(1)\big)
  $$
  as $u$ tends to $0$, and
  $$
    \Phi^\leftarrow (1-u)^2
    = 2\log 1/u -\log\log 1/u - 2\log (2\sqrt{\pi}) +o(1)
    \qquad\hbox{ as } u\to 0 \, .
  $$}

\bigskip

\stateit{Proof.} As in the proof of Lemma A.1.2, we start with
the equality $u=1-\Phi (x)$, that is $x=(1-\Phi )^\leftarrow
(u)$. We consider $x$ tending to infinity, or equivalently 
$u$ converging to $0$. Lemma A.1.3 implies
$$
  \log u = -{x^2\over 2} -\log x -\log\sqrt{2\pi} + o(1)
  \qquad\hbox{ as } x\to\infty \, , \, u\to 0 \, .
  \eqno{({\rm A}.1.1)}
$$
Consequently, $x=\sqrt{2\log 1/u}(1+x_1)$ with $x_1$ tending to
$0$ as $u$ tends to $0$. Then (A.1.1) yields
$$\eqalignno{
  \log u
  & = \Big(-\log {1\over u}\Big)(1+x_1)^2 
    -{1\over 2} \log\log {1\over u} - \log (2\sqrt{\pi}) \cr
  & \phantom{{}= \log u - 2x_1\log {1\over u} 
             - x_1^2\log {1\over u} } 
    -\log (1+x_1) + o(1) \cr
  & {{} = \log u - 2x_1\log {1\over u} - x_1^2\log {1\over u}}
    -{1\over 2}\log\log {1\over u} \cr
  & \phantom{{}= \log u - 2x_1\log {1\over u} 
             - x_1^2\log {1\over u} } 
    -\log (2\sqrt{\pi})+o(1) \, . 
  &({\rm A}.1.2)\cr
  }$$
We can then calculate
$$
  x_1=-{\log\log 1/u\over 4\log 1/u}(1+x_2)
$$
with $x_2$ tending to $0$ with $u$. But now (A.1.2) implies
$$\eqalign{
  0%
  &{{}={1+x_2\over 2} \log\log 1/u }
    - {1\over 16} {(\log\log 1/u)^2\over\log 1/u}(1+x_2)^2
     \cr
  &\phantom{{}={1+x_2\over 2} \log\log 1/u }
    - {1\over 2}\log\log 1/u  
    - \log (2\sqrt{\pi}) + o(1) \cr
  &  = {x_2\over 2}\log\log {1\over u} -\log (2\sqrt{\pi})
    +o(1) \, , \cr
  }$$
and thus, as $u$ tends to $0$,
$$
  x_2={2\log (2\sqrt{\pi})\over\log\log 1/u}\big( 1+o(1)\big) \, .
$$
Consequently, gathering every piece yields
$$\eqalign{
  x%
  & = \sqrt{2\log 1/u} 
    - {1\over 2\sqrt 2}{\log\log 1/u\over \sqrt{\log 1/u}}
    (1+x_2) \cr
  & = \sqrt{2\log 1/u} 
    -{1\over 2\sqrt{2}}{\log\log 1/u\over\sqrt{\log 1/u}}
    -{1\over\sqrt{2}} {\log (2\sqrt{\pi})\over\sqrt{\log 1/u}}
    \big( 1+o(1)\big) \cr
  }$$
as $u$ tends to $0$, which is the desired expression for
$\Phi^\leftarrow (1-u)$. Square it to obtain that for
$\Phi^\leftarrow(1-u)^2$.\hfill$\qed$

\bigskip

We can now obtain an asymptotic expansion for $\Phi^\leftarrow\circ
S_\alpha$.

\bigskip

\state{A.1.5. LEMMA.} {\it%
  We have
  $$\displaylines{\qquad
    \pisa (x) 
    = \sqrt{2\alpha\log x}
    -{\log\log x\over 2\sqrt{2\alpha\log x}}
    -{ \log(K_{s,\alpha}\alpha^{\alpha/2}2\sqrt{\pi})\over
      \sqrt{2\alpha\log x} }
    \cut
    + o\Big( {1\over\sqrt{\log x}}\Big) \qquad\cr}
  $$
  as $x$ tends to infinity. Consequently,
  $$
    \pisa (x)^2
    = 2\alpha\log x-\log\log x 
    - 2\log (K_{s,\alpha}\alpha^{\alpha/2}2\sqrt{\pi}) + o(1)
  $$
  as $x$ tends to infinity.}

\bigskip

\stateit{Proof.} From Lemma A.1.4, we deduce that
$$\displaylines{\qquad
  \pisa (x)^2
  = -2\log \big(1-S_\alpha (x)\big)
  - \log \big(-\log (1-S_\alpha(x)\big)
  \cut
  -2\log (2\sqrt{\pi}) + o(1)
  \qquad\cr}
$$
as $x$ tends to infinity. But Lemma A.1.1 implies
$$
  -\log \big( 1-S_\alpha (x) \big)
  =\alpha \log x - \log (K_{s,\alpha} \alpha^{(\alpha-1)/2})
  + o(1)  
$$
as $x$ tends to infinity. Consequently,
$$\displaylines{\qquad
    \pisa (x)^2
    = 2\alpha \log x -2\log (K_{s,\alpha} \alpha^{(\alpha-1)/2})
    -\log\log x-\log \alpha
  \cut
    -2\log (2\sqrt{\pi}) +o(1)
  \cr}
$$
as $x$ tends to infinity, which is the second expansion in Lemma A.1.5.
Taking the square root yields the first assertion since
$$\displaylines{\quad
  \pisa (x)
  \hfill\cr\hskip .3in
  = \sqrt{2\alpha\log x} \Big( 1-
     {\log\log x\over 2\alpha\log x} 
     -{ \log (K_{s,\alpha}\alpha^{\alpha /2}2\sqrt{\pi})\over
	\alpha\log x }
    +o\Big( {1\over\log x}\Big)\Big)^{1/2}
  \hfill\cr\hskip .3in
  = \sqrt{2\alpha\log x}\Big(  1-
     {\log\log x\over 4\alpha\log x} 
     -{ \log (K_{s,\alpha}\alpha^{\alpha /2}2\sqrt{\pi})\over
	2\alpha\log x }
     +o\Big( {1\over\log x}\Big)\Big)\, . 
  \hfill\qed\cr
  }
$$

\bigskip

We can also obtain an estimate of $\Saip (x)$.

\bigskip

\state{A.1.6. LEMMA.} {\it%
  We have
  $$
    \log \Saip (x) = {x^2\over 2\alpha} +{1\over\alpha}\log x
    + {1\over\alpha}\log( K_{s,\alpha}\alpha^{\alpha-1\over 2}
    \sqrt{2\pi}) + o(1)
  $$
  as $x$ tends to infinity.}

\bigskip

\stateit{Proof.} Combine Lemma A.1.2 and A.1.3 to
obtain
$$\eqalign{
  \Saip (x)
  & = (1-S_\alpha)^\leftarrow \circ (1-\Phi ) (x) \cr
  & \sim K_{s,\alpha}^{1/\alpha} \alpha^{(\alpha -1)/2\alpha}
    x^{1/\alpha}(2\pi )^{1/2\alpha} e^{x^2/2\alpha}
    \cr
  }$$
as $x$ tends to infinity. The result follows by taking the logarithm.
\hfill$\qed$

\bigskip

In the case of an exact Student distribution, the term $o(1)$ in Lemma
A.1.6 is actually $O(x^{-1/\alpha}e^{-x^2/\alpha})$. In this special
case, the approximation has a terrific accuracy as $x$ tends to
infinity!
\vfill\eject

  ~\vskip 28pt
\chapter{Appendix 2}{Exponential map}
{\fchapter
  \line{\hfil Appendix 2.\quad}
  \line{\hfil Exponential map}}
\chapternumber={2}
\vskip 54pt

\noindent The purpose of this appendix is to state and prove the
following proposition. It gives a bound on the error committed by
linearizing the exponential map over the level set of a function.

\bigskip

\state{A.2.1. PROPOSITION.} {\it%
  Let $I:\RR^d\to\RR$ be a smooth smooth function, and let
  $c$ be a regular value of $I$. Then
  $\Lambda_c=I^{-1}(c)$ is a smooth
  manifold of $\RR^d$. Let
  $$
    M=\sup\Big\{ {|\D^2I(p)|\over |\D I(p)|}\, :\, p\in\RR^d \, , \,
    I(p)=c\Big\} \, .
  $$
  Then, for any $p$ in $\Lambda_c$, any $u$ in $T_p\Lambda_c$ with
  $|u| < 1/4M$,
  $$
    |\exp_p(u)-p-u|\leq M|u|^2 \, .
  $$}

\bigskip

\stateit{Proof.} Let $p$ be a point in $\Lambda_c$ and $v$ a unit
tangent vector to $\Lambda_c$ at $p$. We denote by 
$\gamma (s)=\exp_p(sv)$ the geodesic starting at $p$ in the direction
$u$ on $\Lambda_c$. This parameterization is by arc length.  Recall
that $N=\D I/|\D I|$ is an outward unit normal vector field to the
level set of $I$.  Since
$$
  \d N(p)=\Big( \Id-{1\over 2}NN^\T\Big) {\D^2I\over |\D I|}(p)\, ,
$$
we have $\| \d N(p)\|\leq M$. Consequently,
$$
  \big| N\big( \gamma (s)\big)-N(p)\big|
  =\big|\int_0^s \d N\big(\gamma (r)\big) \cdot\gamma'(r) \d r\big|
  \leq sM \, .
$$
The geodesic $\gamma (\cdot )$ is characterized by the 
parallel transport of its tangent vectors along $\gamma (s)$ 
which can be written as
$$
  \big[ \Id -NN^\T\big(\gamma (s)\big)\big] \gamma''(s)
  = \Proj_{T_{\gamma (s)}\Lambda_c}\gamma''(s) = 0 \, .
$$
Any $w$ in $T_p\Lambda_c$ is orthogonal to $M(p)$. Consequently,
$$
  \< w,\gamma''(s)\>
  = \big\< w,\big[ NN^\T\big(\gamma (s)\big)-NN^\T(p)\big] \cdot
    \gamma''(s)\big\> \, .
$$
By duality,
$$\eqalign{
  |\Proj_{T_p\Lambda_c}\gamma''(s)|
  & \leq \big\| NN^\T\big(\gamma (s)\big)-NN^\T (p)\big\|
    \, | \gamma''(s)| \cr
  & \leq 2\big| N\big(\gamma (s)\big)-NN^\T (p)\big| 
    \, |\gamma''(s)| \cr
  & \leq 2sM|\gamma''(s)| \, .\cr
  }$$
Moreover, since $\big\< N\big(\gamma (s)\big)\, ,\gamma'(s)\big\>
=0$ and $\gamma$ is parametrized by arc length,
$$
  \big|\big\< N\big(\gamma (s)\big)\, , \gamma''(s)\big\> \big|
  = \big| -\big\< \d N\big( \gamma (s)\big)\cdot \gamma'(s)\, ,
    \gamma'(s)\big\>\big|
  \leq M \, .
$$
It follows that
$$\eqalign{
  |\Proj_{(T_p\Lambda_c)^\perp}\gamma''(s)|
  & = |\< N(p)\, ,\gamma''(s)\>| \cr
  & \leq \big| N(p)-N\big(\gamma (s)\big)\big|\, |\gamma''(s)|
    + M\cr
  &\leq sM|\gamma''(s)| + M \, . \cr
  }$$
Consequently, we have the inequality
$$\eqalign{
  |\gamma''(s)|^2
  & = |\Proj_{T_p\Lambda_c}\gamma''(s)|^2 
    + |\Proj_{(T_p\Lambda_c)^\perp}\gamma''(s)|^2 \cr
  \noalign{\vskip 0.05in}
  & \leq 4 s^2M^2|\gamma''(s)|^2 + M^2( s|\gamma''(s)|+1)^2 \, .
    \cr
  }$$
On the range $s\leq 1/(4M)$, this inequality implies
$$
  |\gamma''(s)|
  \leq {|\gamma''(s)|^2\over 4} + {|\gamma''(s)|^2\over 8}
  + 2M^2 \, ,
$$
that is
$$
  |\gamma''(s)|^2\leq 16M^2/5 \leq 4M^2 \, .
$$
Consequently, since $\gamma'(0)=v$,
$$
  |\gamma'(s)-v|
  = \big|\int_0^s\gamma''(t) \d t \big|
  \leq 2Ms \, .
$$
The result follows, since
$$
  |\gamma (s)-p-sv|
  =\big|\int_0^s\big(\gamma'(t)-v\big) \d t\big|
  \leq \int_0^s 2Mt \d t
  = Ms^2
  = M|sv|^2 \, .
  \eqno{\qed}
$$
\vfill\eject

}

{
  \def\rightheadline{\ninepointsss References\hfill\folio}
  \def\leftheadline{\ninepointsss\folio\hfill References}
  ~\vskip 28pt
\chapter{}{References}
{\fchapter\hfill References}
 
\vskip 54pt

\leftskip=\parindent
\parindent=-\parindent

Adler, R.J. (2000). On excursion sets, tube formulas and maxima of 
random fields, {\sl Ann.\ Appl.\ Probab.}, 10, 1--74. 

Andrews, G.E., Askey, R., Roy, R. (1999). {\sl Special Functions},
Cambridge University Press.

Balkema, A.A., Kl\"uppelberg, C., Resnick, S.I. (1999). Limit
laws for exponential families, {\sl Bernoulli}, 5, 951--968. 

Barndorff-Nielsen, O.E. (1978). {\sl Information and Exponential Families
in Statistical Theory}, Wiley.

Barndorff-Nielsen, O.E., Kl\"uppelberg, C. (1999). Tail exactness 
of multivariate saddlepoint approximations, {\sl Scand.\ J. \
Statist.} 26 (1999), 253--264.

Barndorff-Nielsen, O.E. (1990). Approximation of interval probabilities,
{\sl J.\ Roy.\ Statist.\ Soc.}, B, 52, 485--496.

Barbe, Ph., Broniatowski, M. (200?). A sharp Petrov type large deviation
formula, submitted.

Billinsley, P. (1968). {\sl Convergence of Probability Measures}, Wiley.

Bingham, N.A., Goldie, C.M., Teugels, J.L. (1987). {\sl Regular
Variation}, Cambridge University Press.

Bishop, R.L. (1977). Decomposition of cut loci, {\sl Proc.\ Amer.\ Math.
\ Soc.}, 65, 133--136.

Bleistein, N., Handelman, R. (1975). {\sl Asymptotic Expansions of
Integrals}, Holt, Reinhartm Winston, New York. Republished by Dover in 1986.

Bobkov, S., Ledoux, M. (1997). Poincar\'e's inequality and Talagrand's
concentration phenomenon for the exponential distribution, {\sl Probab. \
Theory\ Rel.\ Fields}, 107, 383--400.

Bobokov, S.G., Ledoux, M. (2000). From Brunn-Minkowski to Bracamp-Lieb
and to logarithmic Sobolev inequalities, {\sl Geom.\ Funct.\ Anal.}, 10,
1028--1052.

Bolthausen, E. (1993). Stochastic processes with long range interactions
of the paths, {\sl Doeblin and Modern Probability}, Bauberen, 1991, pp.297--310,
{\sl Contemp. Math.}, 149, Amer. Math. Soc., Providence, R.I.

Breitung, K. (1994). {\sl Asymptotic Approximations for Probability 
Integrals}, {\sl Lecture Notes in Mathematics}, 1592, Springer.

Breitung, K., Hohenbichler, M. (1989). Asymptotic approximations 
for multivariate integrals with an application to
multinormal probabilities, {\sl J. Multivariate\ Anal.}, 30, 80--97. 

Breitung, K., Richter, W.-D. (1996). A geometric approach to an 
asymptotic expansion for large deviation probabilities of Gaussian 
random vectors, {\sl J.\ Multivariate Anal.}, 58, 1--20. 

Broniatowski, M., Fuchs, A. (1995). Tauberian theorems, Chernoff 
inequality, and the tail behavior of finite
convolutions of distribution functions, {\sl Adv.\ Math.}, 116, 12--33. 

Brockwell, P.J., Davis, R.A. (1987). {\sl Time Series: Theory and
Methods}, Springer.

Brown, L.D. (1986). {\sl Fundamentals of Statistical Exponential Families,
with Applications in Statistical Decision Theory}, Institute of
Mathematical Statistics.

Bruce, J.W., Giblin, P.J. (1992). {\sl Curves and Singularities},
2nd ed., Cambridge University Press.

Chavel, I. (1996). {\sl Riemannian Geometry: A Modern Introduction},
Cambridge University Press.

Combet, E. (1982). {\sl Int\'egrales Exponentielles}, {\sl Lecture Notes in
Mathematics}, 937, Springer.

Cox, D., Little, J., O'Shea, D. (1992). {\sl Ideals, Varieties, and
Algorithms}, 2nd ed., Springer.

Cox, D., Little, J., O'Shea, D. (1998). {\sl Using Algebraic Geometry},
Springer.

Csisz\'ar, I. (1984). Sanov property, generalized $I$-projection and 
a conditional limit theorem, {\sl Ann.\ Probab.}, 12,
768--793.

Davis, R., Resnick, S.I. (1986). Limit theory for the sample covariance
and correlation functions of moving averages, {\sl Ann.\ Statist.}, 14, 
533-558.

De Bruijn, N.G. (1958). {\sl Asymptotic Methods in Analysis}, North Holland,
Republished by Dover in 1981.

Diebolt, J., Posse, Ch. (1996). On the density of the maximum of 
smooth Gaussian processes. {\sl Ann.\ Probab.}, 24, 1104--1129. 

Do Carmo, M. P. (1976). {\sl Differential Geometry of Curves and Surfaces},
Prentice Hall.

Do Carmo, M.P. (1992). {\sl Riemannian Geometry}, Birkhauser.

Dembo, A., Zeitouni, O. (1993). {\sl Large Deviations Techniques 
and Applications}, Jones and Bartlett Publishers, Boston,
MA.

Dupuis, P.,  Ellis, R.S. (1997). {\sl A Weak Convergence Approach 
to the Theory of Large Deviations}, Wiley. 

Federer, H. (1969). {\sl Geometric Measure Theory}, Springer.

Feller, W. (1971). {\sl An Introduction to Probability Theory and Its 
  Applications}, 2nd ed., Wiley.

Fernique, X. (1970). Int\'egrabilit\'e des vecteurs gaussiens, {\sl C.R. \
Acad.\ Sci.\ Paris}, 270, 1698-1699.

G\"unther, P. (1960). Einige S\"atze \"uber das Volumenelement eines
Riemannschen Raumes, {\sl Pub.\ Math.\ Deberecen}, 7, 78--93.

Hasanis, T.,  Koutroufiotis, D. (1985). The characteristic 
mapping of a reflector. {\sl J.\ Geom.}, 24, 131--167.

Hlawka, E., Schoi\ss engeir, J., Taschner, R. (1991). {\sl Geometric and 
Analytic Number Theory}, Springer.

Howard, R. (1993). {\sl The Kinematic Formula in Riemannian Homogeneous
Spaces}, {\sl Memoires of the Amer.\ Math.\ Soc.}, 509.

Hwang, C.-R. (1980). Laplace's method revisited: weak convergence 
of probability measures, {\sl Ann.\  Probab.}, 8, 1177--1182. 

Imhoff, J.P. (1961). Computing the distribution of quadratic forms in
normal variables, {\sl Biometrika}, 48, 419--426.

Jackson, D. (1911). Ueler eine trigonometrishe Summe, {\sl Rend.\ Circ. \
Mat.\ Palermo}, 32, 257--262.

Jensen, J.L. (1995). {\sl Saddlepoint Approximations}, Oxford University
Press.

Johnson, N.L., Kotz, S. (1972). {\sl Distributions in Statistics: 
Continuous Multivariate Distributions}, Wiley.

Klingenberg, W. (1959). Contribution to Riemannian geometry in the
large, {\sl Ann. of Math.}, 69, 654--666.

\finetune{\hfuzz=1pt}
Knapp, A.W. (1996). {\sl Lie Groups Beyond an Introduction},
Birkhauser.

\finetune{\hfuzz=0pt}
Landau, H.J., Shepp, L.A. (1970). On the supremum of Gaussian processes,
{\sl Sanky\`a}, A32, 369--378.

LeCam, L. (1986). {\sl Asymptotic Methods in Statistical Decision
Theory}, Springer.

LeCam, L., Yang, G.L. (1990). {\sl Asymptotics in Statistics, Some Basic
Concepts}, Springer.

Ledoux, M. (1996). {\sl Isoperimetry and Gaussian Analysis}, in {\sl
Lectures on Probability Theory and Statistics, Ecole d'\'et\'e de
Probabilit\'es de Saint Flour}, XXIV, 1994, P. Bernard Ed., {\sl Lecture
Notes in Mathematics}, 1648, Springer.

Lieb, E., Loss, M. (1997). {\sl Analysis}, American Mathematical Society.

Marinov, M.S. (1980). Invariant volumes of compact groups, {\sl J. \
Phys.}, A, 13, 3357--3366.

Mehta, M.L. (1991). {\sl Random Matrices}, 2nd edition, Academic Press.

McCleary, J. (1974). {\sl Geometry from a Differentiable Viewpoint},
Cambridge University Press.

Mijnheer, J. (1997a). Asymptotic inference for AR(1) processes with 
(nonnormal) stable errors, {\sl J.\ Math.\ Sci.}, 83, 401--406.

Mijnheer, J. (1997b). Asymptotic inference for AR(1) processes with 
(nonnormal) stable errors, III,  {\sl Comm.\ Statist.\ Stochastic Models}, 
13, 661--672.

Mijnheer, J. (1997c). Asymptotic inference for AR(1) processes with 
(nonnormal) stable errors, IV, A note on the case of a negative
unit root, {\sl J.\ Math.\ Sci.}, 92, 4035--4037.

Morgan, F. (1988). {\sl Geometric Measure Theory, a Beginner's Guide},
Academic Press.

Morgan, F. (1992). {\sl Riemannian Geometry, a Beginner's Guide}, Peters
Ltd.

Olver, F.W.J. (1974). {\sl Asymptotics and Special Functions}, Academic Press.

Piterbarg, V.V. (1996). {\sl Asymptotic methods in the theory of
Gaussian processes and fields, Translated from the Russian 
by V. V. Piterbarg}, Amer.\ Math.\ Soc.

Pollard, D. (1984). {\sl Convergence of Stochastic Processes}, Springer.

Press, W.H., Teukolsky, S.A., Vettering, N.H., Flannery, B.P. (1996). {\sl
Numerical Recipes}, Cambridge University Press.

Rausch, H.E. (1951). A contribution to differential geometry in the large,
{\sl Ann.\ of\ Math.}, 54, 38--55.

Resnick, S.I. (1987). {\sl Extreme Values, Regular Variation, and Point
Processes}, Springer.

Rosi\'nski, J., Woyczy\'nski, W.A. (1987). Multilinear forms in Pareto-like
random variables and product random measures, {\sl Colloquium Mathematicum},
305--313.

Schneider, R. (1993). {\sl Convex Bodies: the Brunn-Minkowski Theory},
Cambridge University Press.

Skorokhod, A.V. (1956). Limit theorem for stochastic processes, 
{\sl Theor.\ Probab.\ Appl.}, 1, 261--290.

Spivak,M. (1970). {\sl A Comprehensive Introduction to Differential
Geometry}, Publish or Perish.

Talagrand, M. (1995). Concentration of measure and isoperimetric
inequalities in product spaces, {\sl Publ.\ Math.\ IHES}, 81, 73--205.

Weyl, H. (1939). On the volume of tubes, Proc.\ Amer.\ Math.\ Soc., 65,
461--472.

\vfill\eject
 % bibliographie
}

{
  \def\rightheadline{\ninepointsss Notations\hfill\folio}
  \def\leftheadline{\ninepointsss\folio\hfill Notations}
  ~\vskip 28pt
\chapter{}{Notation}
{\fchapter\hfill Notation}%

\vskip 54pt

\ninepointsrm

\noindent This index has 3 parts. The first one concerns
general notation, that is used throughout the book, 
eventually with variations in the arguments in each example.

The second part contains the notation introduced in chapters
2--5. This is used almost throughout the book, but in examples
treated in chapters 6--12 is given a specialized meaning.

The third part contains the notation with the concrete meaning
given through chapters 6--12.

Some notation that are used very locally does not appear in the index.

\bigskip

\noindent{\fsection General notation.}

\bigskip

\newdimen\boxdim
\def\notation#1#2{
  \line{\hbox to 1.3in{#1\hfill}\hbox to 3.2in{\vtop{\hsize=3.2in\noindent#2}}}
  }

\notation{$\sharp$}{the cardinal of a set, as in $\sharp\{\, 1, 2, 3\,\}=3$.}
\notation{$\d$}{integration element, as in $\d t$, $\d\calM_M$, $\d x$.}
\notation{$d$}{dimension of the underlying space, $\RR^d$.}
\notation{$\D$, $\D^2$, $\ldots$}{Gradient, Hessian, and higher order differentials.}
\notation{$\exp_p(\cdot )$}{exponential map at $p$ (on a Riemannian manifold).}
\notation{$\calM_M$}{Riemannian measure of a manifold $M$.}
\notation{$\cdot^\perp$}{orthocomplement.}
\notation{$\inj_M(p)$}{radius of injectivity of $p$ in the manifold $M$.}
\notation{$K_M(x,y)$}{sectional curvature of the manifold $M$ along the tangent vector fields $x$ and $y$.}
\kern 5pt
\notation{$\lambda_{\rm min} (M)$, $\lambda_{\rm max} (M)$}{smallest and largest eigenvalue of a matrix $M$.}
\notation{$\partial A$}{boundary of a set.}
\notation{$\partial f(\cdot )/\partial u$}{partial differentiation of a function.}
\notation{$\Pi_{M,p}$}{second fundamental form of the manifold $M$ at $p$.}
\notation{$\RR$, $\RR^d$}{set of real numbers, the Euclidean $d$ dimensional space.}
\notation{${\rm Ricc}$}{Ricci tensor of the level lines of $I$.}
\notation{$S_n$}{unit sphere centered at the origin, of dimension $n$, that is the boundary of the unit ball in $\RR^{n+1}$.}
\kern 5pt
\notation{$S_V(x,r)$}{ball centered at $x$, of radius $r$, in the vector space $V$.}
\notation{$T_pM$}{tangent space of the manifold $M$ at $p$.}
\notation{$N_pM$}{normal space of the manifold $M$ at $p$.}
\notation{$|S|$}{Lebesgue measure of the set $S$.}
\notation{$\omega_n = \pi^n/\Gamma ((n/2)+1)$}{volume of the unit ball of $\RR^n$.}

\finetune{\vfill\eject}

\noindent{\fsection Notation from chapters 2--5.}

\bigskip

{\ninepointsrm

\setbox1=\vtop{\hsize=1.5in\obeylines\parindent=0pt
$A$, 1
$\underline{A}_M$, 53
${\calA}_p(t,v)$, 38
$c_{A,M}^*$, 54
$c_{A,M}$, 54
$\chi_A^F(p)$, 33
$\chi_A^L(p)$, 33
$\calD_A$, 36
$G_A(p)$,54\vss}

\setbox2=\vtop{\hsize=1.5in\obeylines\parindent=0pt
$\Gamma_c$, 19
$I(\cdot )$, 19
$I(A)$, 19
$J\pi_A(\cdot )$, 38
$K_{s,\alpha }$, 2
$K_{w,\alpha }$, 2
$\Lambda_c$, 25
$L(c)$, 19
$\mu_A(B)$, 62\vss}

\setbox3=\vtop{\hsize=1.5in\obeylines\parindent=0pt
$N(x)$, 26
$\omega_A$, 37
$\Pi_{\Lambda_{I(A)},q}^\pi $, 44
$\psi (x,s)=\psi_s(x)$, 26
$\psi_{t*}(\cdot )$, 28
$t_{0,M}(p)$,56
$\tau_A(p)$, 31\vss}

\line{\box1\box2\box3\hfill}
}

\bigskip

\noindent{\fsection Notation from chapters 6--12.}

\bigskip

{ \ninepointsrm

\centerline{\bf Chapter 6}

\medskip

\setbox1=\vtop{\hsize=1.5in\obeylines\parindent=0pt
$K_t$ , p.70
$\nu (p)$ , p.71\vss}
\setbox2=\vtop{\hsize=1.5in\obeylines\parindent=0pt
$f_p(\cdot )$ , p.71
$Q_p$ , p.71\vss}
\setbox3=\vtop{\hsize=1.5in\obeylines\parindent=0pt
$\Pi (\cdot )$ , p.75
$r$ , p.75\vss}
\line{\box1\box2\box3\hfill}

\bigskip

\centerline{\bf Chapter 7}

\medskip

\setbox1=\vtop{\hsize=1.5in\obeylines\parindent=0pt
$A_1$, 79\vss}
\setbox2=\vtop{\hsize=1.5in\obeylines\parindent=0pt
$I(\cdot )$ , 79\vss}
\setbox3=\vtop{\hsize=1.5in\obeylines\parindent=0pt
$\calD_{A_1}$ , 80\vss}
\line{\box1\box2\box3\hfill}

\bigskip

\centerline{\bf Chapter 8}

\medskip

\setbox1=\vtop{\hsize=1.5in\obeylines\parindent=0pt
\quad{\bf\S 1}
$C$, 89
$A_t$, 89
$w_\alpha (\cdot )$, 89
$I(\cdot )$, 89
$H$, 90
\quad{\bf\S 2}
$K_{s,\alpha }$, 93
$A_t$, 93
$J_1$, 93
$\Phi (\cdot )$, 94
$S_\alpha (\cdot )$, 94
$B_t$, 94
$I(\cdot )$, 95
$p(v)$, 95
$p_{\epsilon ,j,t}(v)$, 95
$q(v)$, 96
$q_{\epsilon ,j,t}(v)$, 96
\vss}

\setbox2=\vtop{\hsize=1.5in\obeylines\parindent=0pt
$\gamma_1$, 100
$\calD_{B_t}$, 100
$T$, 108
$\calI$, 108
$V_\calI$, 108
$N(C)$, 108
$J(C)$, 108
$A_t$, 109
$M_\calI$, 109
$\mathaccent "707E {G}(m)$, 109
$\calN_\calI$, 111
$Q(t)$, 111
$p_{\calI,t}(m,v)$, 112
$q_{\calI,t}(m,v)$, 112
$R(t)$, 113
$\gamma_\calI$, 113
$\gamma $, 113
$\calD'_{B_t}$, 116
\vss}

\setbox3=\vtop{\hsize=1.5in\obeylines\parindent=0pt
$\rho_t$, 116
$r_{\calI,t}(m,v)$, 116
$\calD (B,t)$, 116
$q_{\calI,t}(m)$, 116
$r_{\calI,t}(m)$, 116
$p(m)$, 116
$\epsilon_\calC$, 127
$\mu_\calC$, 128
$N_0(C)$, 137
$J_0(C)$, 137
$M_{\epsilon ,\calI}$, 137
$Z_\calI (C)$, 138
$S_\calI (C)$, 139
\quad{\bf\S 3}
$J_*(C)$, 144
$M_\calI$, 144
$\gamma_\calI$, 144
$\gamma $, 145
\vss}

\setbox4=\hbox{\box1\box2\box3\hfill}

\dp4=0pt
\box4
\vfill\eject

\setbox1=\vtop{\hsize=1.5in\obeylines\parindent=0pt
$I(\cdot )$, 145
$R(t)$, 145
$p_\calI(m,v)$, 146
$Q(\cdot )$, 148
\vss}

\setbox2=\vtop{\hsize=1.5in\obeylines\parindent=0pt
$A'_t$, 150
$B'_t$, 150
$\Omega_t$, 154
$B''_t$, 155
\vss}

\setbox3=\vtop{\hsize=1.5in\obeylines\parindent=0pt
$\rho_t$, 157
$q_{\calI,t}(m)$, 157
$r_{\calI,t}(m)$, 157
$\calD_{B''_t}$, 157
\vss}

\line{\box1\box2\box3\hfill}

\bigskip

\centerline{\bf Chapter 9}

\medskip

\setbox1=\vtop{\hsize=1.5in\obeylines\parindent=0pt
$X(M)$, 165
$A_t$, 165
\quad{\bf\S 1}
$C_M$, 166
$M_0$, 167
$N(\cdot )$, 168
$p(\cdot )$, 171
\vss}

\setbox2=\vtop{\hsize=1.5in\obeylines\parindent=0pt
$\nu (p)$, 171
$X_i(p)$, 171
$X(u)$, 171
$\phi _j (p)$, 171
$f_{\scriptscriptstyle \bullet }$, 173
$f^{\scriptscriptstyle \bullet }$, 173
$I_{\scriptscriptstyle \bullet }$, 173
\vss}

\setbox3=\vtop{\hsize=1.5in\obeylines\parindent=0pt
\quad{\bf\S 3}
$A_t$, 186
$B_t$, 186
$\gamma $, 187
$\calD_{B_t}$, 188
\vss}

\line{\box1\box2\box3\hfill}

\bigskip

\centerline{\bf Chapter 10}

\medskip

\setbox1=\vtop{\hsize=1.5in\obeylines\parindent=0pt
${\fam 0\elevenrm M}(n,{\fam \bbfam \elevenbb R})$, 197
$\glnr$, 197
${\fam \frakfam \eleveneufm S}_n$, 198
\quad{\bf\S 1}
$w_\alpha (\cdot )$, 200
$I(\cdot )$, 200
$A_t$, 200
${\fam 0\elevenrm SL}(n,{\fam \bbfam \elevenbb R})$, 201
${\fam 0\elevenrm SO}(n,{\fam \bbfam \elevenbb R})$, 201
\quad{\bf\S 2}
$F_n$, 211
\vss}

\setbox2=\vtop{\hsize=1.5in\obeylines\parindent=0pt
$c_n$, 211
\quad{\bf\S 3}
$\delimiter "426830A  \cdot ,\cdot \delimiter "526930B $, 219
$\delimiter "26B30D \cdot \delimiter "26B30D $, 219
$\calS$, 219
$H_{u,v}$, 219
$\calS^0$, 219
$H_{u,v,h}^1$, 221
$H_{u,v,h}^2$, 221
$a_i$, 222
$f_i$, 223
\vss}

\setbox3=\vtop{\hsize=1.5in\obeylines\parindent=0pt
$b_j$, 223
$b^v$, 223
$b^u$, 223
$g_j$, 223
$c_j$, 225
$\Pi $, 225
\quad{\bf\S 4}
$w_\alpha $, 229
$A_t$, 230
$I(\cdot )$, 230
$u_*$, $v_*$, 231
\vss}

\line{\box1\box2\box3\hfill}

\bigskip

\centerline{\bf Chapter 11}

\medskip

\setbox1=\vtop{\hsize=1.5in\obeylines\parindent=0pt
\quad{\bf\S 1}
$B$, 239
$\theta $, 239
$X$, 239
$\epsilon $, 239
$\gamma _n(k)$, 240
\quad{\bf\S 2}
\vss}
\setbox2=\vtop{\hsize=1.5in\obeylines\parindent=0pt
$a$, 243
$A$, 244
$B$, 244
$C$, 245
$\mathaccent "705E \gamma _n(0)$, 252
$C$, 253
\quad{\bf\S 3}
\vss}
\setbox3=\vtop{\hsize=1.5in\obeylines\parindent=0pt
$R_k$, 256
$C$, 257
$a$, $b$, 258
$\phi $, 264
$r$, 264
$g_k(\cdot )$, 266
$h(\theta )$, 266
\vss}

\line{\box1\box2\box3\hfill}

\bigskip

\centerline{\bf Chapter 12}

\medskip

\setbox1=\vtop{\hsize=1.5in\obeylines\parindent=0pt
\quad{\bf\S 1}
$M$, 273
$X(M)$, 273
\vss}
\setbox2=\vtop{\hsize=1.5in\obeylines\parindent=0pt
$\sigma ^2(M)$, 273
$S_d^{(\beta )}$, 274
$M_\beta $, 274
\vss}
\setbox3=\vtop{\hsize=1.5in\obeylines\parindent=0pt
\quad{\bf\S 2}
$C$, 276\vss}

\line{\box1\box2\box3\hfill}

\finetune{\vfill\eject}

\setbox1=\vtop{\hsize=1.5in\obeylines\parindent=0pt
\quad{\bf\S 3}
$\Sigma $, 280
$\sigma ^2$, 280
\vss}
\setbox2=\vtop{\hsize=1.5in\obeylines\parindent=0pt
\quad{\bf\S 4}
$B(\cdot )$, 282
$B_d(s)$, 282
\vss}
\setbox3=\vtop{\hsize=1.5in\obeylines\parindent=0pt
\line{}
$p(s)$, 283
$\sigma _d^2$, 283
\vss}

\line{\box1\box2\box3\hfill}

}

\vfill\eject
 % index des notations
}

~\vskip 28pt
\chapter{}{}
{\fchapter\hfill Postface}
 
\vskip 54pt\elevenrm
 
\noindent This book is the first of a larger project that I may try to 
complete. A second volume should be devoted to the asymptotic analysis
of multivariate integrals over small wedges and their applications.  A
third one should extend some of the results of the first two volumes
to the infinite dimensional setting, where there are some potentially
amazing applications in the study of stochastic processes.

\bigskip

\line{June 2001\hfill Philippe Barbe}

\vfill\eject

\bye